\documentclass[12pt]{book}

\usepackage{imakeidx}  

\usepackage{amsthm}
\usepackage{amsmath} 
\usepackage{amssymb}

\usepackage{marginnote}  

\usepackage[dvips]{color}
\usepackage{graphicx} 

\usepackage[backref=page]{hyperref}   
\hypersetup{
    colorlinks=true,
    linkcolor=blue,
    filecolor=magenta,
    urlcolor=blue,
		citecolor=red,
    pdfpagemode=FullScreen,
}

\usepackage[toc]{appendix}  

 \usepackage{cancel}

\usepackage{xcolor}

\definecolor{gold}{rgb}{0.85,.66,0}

\usepackage{centernot}  

\newcommand{\sidenote}[1]{\setlength{\marginparwidth}{2cm}%
  \marginnote{\sffamily\raggedright\scriptsize #1}}

\usepackage{dsfont}

\usepackage{enumitem}

\usepackage[singlelinecheck=off,font=footnotesize]{caption}

\usepackage{multicol}  

\usepackage{longtable} 
\usepackage{array}
\newcolumntype{R}{>{\raggedleft\arraybackslash}p{.15\textwidth}}

\usepackage{pdfsync}

\usepackage[latin9]{inputenc}

\usepackage[OT1]{fontenc}
\usepackage[ngerman, english]{babel}

\usepackage{rotating}    

\usepackage{geometry}  
\makeatletter
\DeclareOldFontCommand{\rm}{\normalfont\rmfamily}{\mathrm}
\DeclareOldFontCommand{\sf}{\normalfont\sffamily}{\mathsf}
\DeclareOldFontCommand{\tt}{\normalfont\ttfamily}{\mathtt}
\DeclareOldFontCommand{\bf}{\normalfont\bfseries}{\mathbf}
\DeclareOldFontCommand{\it}{\normalfont\itshape}{\mathit}
\DeclareOldFontCommand{\sl}{\normalfont\slshape}{\@nomath\sl}
\DeclareOldFontCommand{\sc}{\normalfont\scshape}{\@nomath\sc}
\makeatother

\usepackage{csquotes}  
\newcommand\q{\enquote}

\usepackage{physics}   

\newcommand{\ccat}[3]{{#1\, \underset{#3}{\lozenge}\,{#2}}}

\newcommand{\tm}{\times}
\newcommand\sat{\textup{sat}}

\newcommand{\UGS}{\mathrm{UGS}}%

\newcommand{\unit}{\mathds{1}}%
\DeclareMathOperator{\loc}{loc}

\DeclareMathOperator*{\esssup}{ess\,sup}

\newcommand \diag  {\operatorname{diag}}

\newcommand \Limsup {\mathop{\overline{\lim}}}
\newcommand \Liminf {\mathop{\underline{\lim}}}

\newcommand \eps {\varepsilon}

\newcommand \N   {\mathbb{N}}
\newcommand \R   {\mathbb{R}}

\newcommand \C   {\mathbb{C}}
\newcommand \Z   {\mathbb{Z}}
\newcommand \Q   {\mathbb{Q}}

\newcommand \A   {\mathcal{A}}

\newcommand \K   {\mathcal{K}}
\newcommand \Kinf{\mathcal{K_\infty}}

\newcommand \KL  {\mathcal{KL}}

\newcommand \LL  {\mathcal{L}}

\newcommand \PD   {\mathcal{P}}

\newcommand{\Uc}{\ensuremath{\mathcal{U}}}

\newcommand{\Dc}{\ensuremath{\mathcal{D}}}
\newcommand{\Sc}{\ensuremath{\mathcal{S}}}

\newcommand{\vertiii}[1]{{\left\vert\kern-0.25ex\left\vert\kern-0.25ex\left\vert #1 
    \right\vert\kern-0.25ex\right\vert\kern-0.25ex\right\vert}}

\renewcommand{\norm}[1]{\left\lVert #1 \right\rVert} 
\newcommand{\modulus}[1]{\left\lvert #1 \right\rvert} 
\newcommand{\sg}[1]{(#1(t))_{t\geq 0}} 

\newcommand \qrq   {\quad\Rightarrow\quad}

\newcommand \qiq   {\quad\Iff\quad}
\newcommand \srs   {\ \ \Rightarrow\ \ }

\newcommand \Iff   {\Leftrightarrow}

\newcommand \supp{\operatorname{supp}}
\newcommand \id  {\operatorname{id}}

\newcommand \re  {\mathrm{Re}}
\newcommand \im  {\mathrm{Im}}

\newcommand{\normt}[1]{{\left\vert\kern-0.25ex\left\vert\kern-0.25ex\left\vert #1 
		\right\vert\kern-0.25ex\right\vert\kern-0.25ex\right\vert}}

\renewcommand{\ker}{{\rm Ker}\,}
\newcommand{\intt}{{\rm int}\,}

\newcommand{\dist}{{\rm dist}\,}

\newcommand{\Ah}{\hat{A}}

\newcommand{\Rh}{\hat{R}}

\newcommand{\midset}{\;:\;}

\newcommand{\clo}[1]{\overline{#1}}

\newcommand{\soc}[1]{{\color{blue} #1}}           

\newcommand{\todo}[1]{{\color{red}\bf TO DO: #1}} 

\newcommand{\mir}[1]{{\color{red}\bf AM: #1}}     
\newcommand{\amc}[1]{{\color{blue} #1}}           

\newif\ifMath					
\newif\ifEngi					

\newif\ifAndo              
													
\newif\ifExercises					
\newif\ifSolutions          
\newif\ifGerman							
\newif\ifEnglish						

\newif\ifnothabil						

\newif\ifFuture							

\newif\ifConf                    
\newif\ifJournal								 

\newif\ifNOTFORBOOK
\newif\ifFullVersion
\newif\ifExludedDueToSpaceReasons

\usepackage{xifthen}

\newcommand{\einsnorm}[2]{\ensuremath{
    \!\!\;\!\!\!\;
    \left\bracevert\!\!\!\!\!\left\bracevert
    \!
		\ifthenelse{\isempty{#2}}{#1}{#1(#2)}
    \!
      \right\bracevert\!\!\!\!\!\right\bracevert
    \!\!\;\!\!\!\;
  }}

\usepackage{xcolor}
\definecolor{blond}{rgb}{0.98, 0.94, 0.75}
	
\newlength\mytemplen
\newsavebox\mytempbox

\makeatletter
\newcommand\mybluebox{%
    \@ifnextchar[
       {\@mybluebox}%
       {\@mybluebox[0pt]}}

\def\@mybluebox[#1]{%
    \@ifnextchar[
       {\@@mybluebox[#1]}%
       {\@@mybluebox[#1][0pt]}}

\def\@@mybluebox[#1][#2]#3{
    \sbox\mytempbox{#3}%
    \mytemplen\ht\mytempbox
    \advance\mytemplen #1\relax
    \ht\mytempbox\mytemplen
    \mytemplen\dp\mytempbox
    \advance\mytemplen #2\relax
    \dp\mytempbox\mytemplen
    \colorbox{blond}{\hspace{1em}\usebox{\mytempbox}\hspace{1em}}}

\makeatother

\makeatletter
\let\origd=\d
\renewcommand*\d{
  \relax\ifmmode
    \mathrm{d}%
  \else
    \expandafter\origd
  \fi
}\makeatother

\usepackage{mathtools}

\newcommand{\lel}{\left\langle}  
\newcommand{\rir}{\right\rangle} 
\newcommand{\scalp}[2]{ \lel #1, #2 \rir }    

\makeatletter 
\newcommand{\pushright}[1]{\ifmeasuring@#1\else\omit\hfill$\displaystyle#1$\fi\ignorespaces}
\newcommand{\pushleft}[1]{\ifmeasuring@#1\else\omit$\displaystyle#1$\hfill\fi\ignorespaces}
\makeatother 

\newcounter{syscounter}
\newenvironment{sysnum}{\begin{list}{($\Sigma{\arabic{syscounter}}$)}%
{\settowidth{\labelwidth}{($\Sigma4$)}
\settowidth{\leftmargin}{($\Sigma4$)~}%
\usecounter{syscounter}}}
{\end{list}}

\newcounter{WPcounter}

\newtheorem{theorem}{Theorem}[section]
\newtheorem{lemma}[theorem]{Lemma}
\newtheorem{proposition}[theorem]{Proposition}
\newtheorem{corollary}[theorem]{Corollary}
\newtheorem{example}[theorem]{Example}
\newtheorem{ass}[theorem]{Assumption}  
\newtheorem{exercise}{Exercise}[section] 
 
\newtheorem*{solution*}{Solution}

\newtheorem{nnremark}[theorem]{\bf Remark}
\newtheorem{nndefinition}[theorem]{\bf Definition}
\newenvironment{remark}{\begin{nnremark} \rm }{\hfill \hspace*{1pt}\hfill $\lhd$\end{nnremark}}
\newenvironment{definition}{\begin{nndefinition} \rm }{\hfill \hspace*{1pt}\hfill $\lhd$\end{nndefinition}}

\newcommand{\calL}{\mathcal{L}}

\usepackage{pgfplots}
\usepackage{pgf}
\usepackage{tikz} 
\usetikzlibrary{decorations.pathmorphing} 
\usepackage{tikz-3dplot}
\usepgfplotslibrary{fillbetween}
\usetikzlibrary{arrows}
\usetikzlibrary{graphs,decorations.pathmorphing,decorations.markings}
\usetikzlibrary{calc,math}
\usetikzlibrary{patterns}
\usetikzlibrary{shapes}
\usetikzlibrary{tikzmark}

\usetikzlibrary {positioning}
\usetikzlibrary{shadows}
\usetikzlibrary{patterns.meta}
\usetikzlibrary{plotmarks}

\usepackage{subcaption}  

\usepackage[absolute,overlay]{textpos}
\usepackage{scalefnt}   

\tikzdeclarepattern{
  name=hatch,
  parameters={\hatchsize,\hatchangle,\hatchlinewidth},
  bounding box={(-.1pt,-.1pt) and (\hatchsize+.1pt,\hatchsize+.1pt)},
  tile size={(\hatchsize,\hatchsize)},
  tile transformation={rotate=\hatchangle},
  defaults={
    hatch size/.store in=\hatchsize,hatch size=5pt,
    hatch angle/.store in=\hatchangle,hatch angle=0,
    hatch linewidth/.store in=\hatchlinewidth,hatch linewidth=.4pt,
  },
  code={
      \draw[line width=\hatchlinewidth] (0,0) -- (\hatchsize,\hatchsize);
  }
}

\pgfmathsetmacro\weight{1/2}
\pgfmathsetmacro\third{1/3}
\pgfmathsetmacro\twothirds{2/3}

\tikzset{degil/.style={
            decoration={markings,
            mark= at position 0.5 with {
                  \node[transform shape] (tempnode) {$/$};
                  }
              },
              postaction={decorate}
}
}

\usepgfplotslibrary{fillbetween}
\usetikzlibrary{intersections,through}

\makeatletter 
\tikzset{use path/.code=\tikz@addmode{\pgfsyssoftpath@setcurrentpath#1}}
\makeatother

\definecolor{manipulator-color}{RGB}{88,44,44}
\definecolor{manipulator-contour}{rgb}{0.0, 0.18, 0.39}  

\tikzset{>=latex} 

\Andofalse %

\nothabilfalse

\Exercisesfalse

\Solutionsfalse

\Germanfalse

\Englishtrue

\Futurefalse

\newif\ifUNUSED
\UNUSEDfalse

\title{Input-to-state stability of distributed parameter systems}

\author{Andrii Mironchenko}

\makeindex   

\begin{document}


\begin{titlepage}

\bigskip
\bigskip
\bigskip
\bigskip

\begin{center}
  
  \bigskip
\bigskip

\vspace*{2cm}

  \bf {\LARGE Input-to-state stability of distributed parameter systems}

\end{center}
\vspace{0,4cm}
\begin{center}\Large
Andrii Mironchenko \\
  \end{center}
  
  \vspace{3cm}
  
\begin{center}
{\bf Habilitation thesis}\\

  \vspace{4cm}
  
	
	Revised version (\today) of the habilitation thesis\\
	submitted to the Faculty of Computer Science and Mathematics\\
of the University of Passau in June 2022\\
    \end{center}

\begin{figure}[b]
	\centering
		\includegraphics[height=1.6cm]{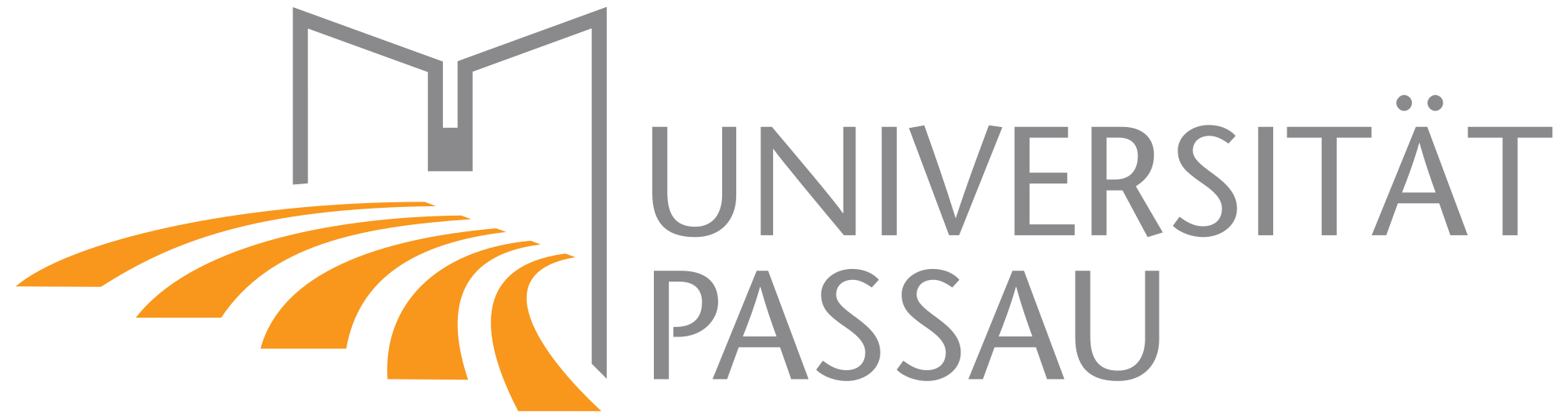}
	\label{fig:unilogo}
\end{figure}

  \end{titlepage}

\newpage
\thispagestyle{empty}
\vspace*{\fill}
{\small 
Date of the habilitation colloquium: 13.12.2022 \\
\bigskip

%
Reviewers: Prof. Dr. Jean-Michel Coron, University Pierre et Marie Curie, France\\
\phantom{GulMentorsss: }Prof. Dr. Marius Tucsnak, University of Bordeaux, France\\
\phantom{GulMentorsss: }Prof. Dr. Fabian Wirth, University of Passau, Germany\\
\phantom{GulMentorsss: }Prof. Dr. Hans Zwart, University of Twente, the Netherlands\\


}

\clearpage

\newpage

\cleardoublepage
\pdfbookmark[0]{Table of Contents}{sumario_label}\tableofcontents

\setcounter{chapter}{0}

\cleardoublepage
\chapter*{\centering \begin{normalsize}Acknowledgements \end{normalsize}}

\addcontentsline{toc}{chapter}{\numberline{}Acknowledgements}

The author thanks to all his coauthors, for many inspiring discussions and ideas shared.
In particular, parts of this thesis were published as joint works with (in alphabetic order)
Sergey Dashkovskiy, 
Jochen Gl\"uck, 
Hiroshi Ito, 
Birgit Jacob,
Iasson Karafyllis, 
Christoph Kawan, 
Miroslav Krstic, 
Navid Noroozi, 
Jonathan Partington, 
Christophe Prieur, 
Fabian Wirth.

I thank to my habilitation mentors Thomas M\"uller-Gronbach, Fabian Wirth, and Hans Zwart for their guidance during the habilitation process. 
I thank to the reviewers Jean-Michel Coron, Marius Tucsnak, Fabian Wirth, and Hans Zwart for their careful evaluation of the submitted version of this work.

I thank to Alexander Kilian, who has proofread this thesis.

Last but not least, I thank to my daughter Darina, whose desire to study the nature around us always inspired me.




\cleardoublepage
\chapter*{Acronyms}
\label{Sec:Acronyms}

\thispagestyle{empty}
\addcontentsline{toc}{chapter}{\numberline{}Acronyms}


\begin{longtable}{R p{.55\textwidth} R} 
0-GAS	  & Global asymptotic stability at zero & Definition~\ref{Stab_Notions_Undisturbed_Systems} \\
0-GATT	  &  Global attractivity at zero & Definition~\ref{Stab_Notions_Undisturbed_Systems}\\
0-UAS	  & Uniform asymptotic stability at zero & Definition~\ref{Stab_Notions_Undisturbed_Systems} \\
0-UGAS	  & Uniform global asymptotic stability at zero & Definition~\ref{Stab_Notions_Undisturbed_Systems} \\
AG	  & Asymptotic gain property & Definition~\ref{def:asymptotic_gain} \\

BECS	&   Bounded energy-convergent state & Definition~\ref{def:BECS}\\
BIC	&   Boundedness-implies-continuation & Definition~\ref{def:BIC}\\
BRS	  &  Boundedness of reachability sets/Bounded reachability sets& Definition~\ref{def:BRS} \\
bUAG	  &   Bounded input uniform asymptotic gain property & Definition~\ref{def:asymptotic_gain}
\\
bULIM	  &   Uniform limit property on bounded sets & Definition~\ref{def:Limit-properties}\\

CEP	  &   Continuity at the equilibrium point& Definition~\ref{def:RobustEquilibrium_Undisturbed}\\

eISS	  & Exponential input-to-state stability & Definition~\ref{Def:eISS}\\


iISS	  &  Integral input-to-state stability &  Definition~\ref{def:iISS} \\
ISpS	  & Input-to-state practical stability  & Definition~\ref{Def:ISpS_wrt_set}\\
ISS	  & Input-to-state stability/Input-to-state stable & Definition~\ref{Def:ISS}\\

LIM	  &  Limit property& Definition~\ref{def:Limit-properties}\\
LISS &  Local input-to-state stability & Definition~\ref{Def:LISS}\\
MAF &  Monotone aggregation function & Definition~\ref{def:Monotone aggregation functions}\\
MBI &  Monotone bounded invertibility & Definition~\ref{def:MBI}\\
MLIM &  Monotone limit property & Definition~\ref{def:MLIM}\\
OBS &  Ordered Banach space & Definition~\ref{def:OBS}\\


ODE, ODEs	  & Ordinary differential equation(s) & \\
PDE, PDEs & Partial differential equation(s) & \\
RFC	  &  Robust forward completeness& Definition~\ref{Def_RFC} \\

sAG	  &  Strong asymptotic gain property& Definition~\ref{def:asymptotic_gain}\\
sISS	  & Strong input-to-state stability & Definition~\ref{Def:sISS}\\
sLIM &  Strong limit property & Definition~\ref{def:Limit-properties}\\


UAG	  &  Uniform asymptotic gain property& Definition~\ref{def:asymptotic_gain}\\

UGAS	  &  Uniform global asymptotic stability & Definition~\ref{def:UGAS}\\
UGATT	  &  Uniform global attractivity& Definition~\ref{def:UniformGlobalAttractivity}\\
UGB	  &  Uniform global boundedness & Definition~\ref{def:ULS-UGS-pGS} \\
UGS	  &  Uniform global stability & Definition~\ref{def:ULS-UGS-pGS} \\
ULIM	  &  Uniform limit property & Definition~\ref{def:Limit-properties}\\
ULS	  &  Uniform local stability & Definition~\ref{def:ULS-UGS-pGS}\\	
WURS	  &  Weak uniform robust stability& Definition~\ref{def:WURS}\\

	  &  & \\
	  &  & \\
	  &  & \\
\end{longtable}

\cleardoublepage
\chapter*{Glossary}
\label{Sec:Lyst_Symb}

\thispagestyle{empty}
\addcontentsline{toc}{chapter}{\numberline{}Glossary}

%
%
  %
%

\begin{center}
\textbf{Sets and numbers}
\end{center}
\vspace{-4mm}
\begin{longtable}{R p{.85\textwidth} } 
	$\N$       & Set of natural numbers: $\N:=\{1,2,\ldots\}$ \\
	$\Z$       & Set of integer numbers \\
	$\R$       & Set of real numbers \\
	$\Z_+$     & Set of nonnegative integer numbers \\
	$\R_+$     & Set of nonnegative real numbers \\
	$\C$       & Set of complex numbers \\
	$U^{\R_+}$ & Set of maps from $\R_+$ to $U$, where $U$ is a certain set \\
  $\re z$    & Real part of a complex number $z$\\
  $\im z$  & Imaginary part of a complex number $z$\\	
	$S^n$  & $:=\underbrace{S \times \ldots \times S}_{n \text{ times}}$ \\
	$\unit$  &   \\ 
  $x^T$  & Transposition of a vector $x \in \R^n$	 \\
	$|\cdot|$  & Euclidean norm in the space $\R^s$, $s \in \N$   \\
	$\|A\|$  & For a matrix $A \in \R^{n \times n}$ we denote $\|A\|:=\sup_{x\neq0}\frac{|Ax|}{|x|}$\\
\end{longtable}

\begin{center}
\textbf{Various notation}
\end{center}
\vspace{-4mm}

\begin{longtable}{R p{.85\textwidth} } 
	$\dist(z,Z)$ & $:=\inf\{\|y-z\|_X: y \in Z \}$. The distance between $z\in X$ and $Z \subset X$ \\
  $B_r(Z)$  & Open ball of radius $r$ around $Z\subset X$, i.e., $\{x\in X:\dist(x,Z)<r\}$\\	
  $B_{r,\; \Uc}$  & Open ball of radius $r$ around $0$ in a normed vector space $\Uc$, i.e., $\{u\in\Uc:\|u\|_{\Uc}<r\}$\\	
  $B_r(x)$  &$:=B_r(\{x\})$\\						
  $B_r$  &$:=B_r(0)$\\	
	$\intt A$ & (Topological) interior of a set $A$ (in a given topology)\\	
	$\clo{A}$ & Closure of a set $A$ (in a given topology)\\	
	$V \subset\subset U$ & $V$ is compactly contained in $U$, that is, $\clo{V}\subset U$\\	
  $\nabla f$  & Gradient of a function $f:\R^n \to \R$ \\
	$f \circ g$  & Composition of maps $f$ and $g$: $f\circ g (s)= f(g(s))$, for $s$ from the domain of definition of $g$   \\
	$\partial G$  &  Boundary of a domain $G$ \\	
	$\mu$  & Lebesgue measure on $\R$ \\
	$I$    & Identity operator (or identity matrix)  \\	
	$\id$  & Identity operator \\	
	$\Limsup$  & $:=\limsup$ (limit superior)\\
  $\Liminf$  & $:=\liminf$ (limit inferior)\\
  $t\to a-0$  & $t$ converges to $a$ from the left\\
  $t\to a+0$  & $t$ converges to $a$ from the right\\
  $t\to -0$  & $t$ converges to $0$ from the left\\
  $t\to +0$  & $t$ converges to $0$ from the right\\
  $x \vee y$  & The supremum (aka join) of vectors $x,y \in\R^n_+$\\	
	$u_t$ & For a function $u:\R_+\to\R^m$, we define $u_t$ by
	\begin{eqnarray*}
u_t(s):=
\begin{cases}
u(s) & \text{, if } s\in[0,t], \\ 
0 & \text{, if } s>t.
\end{cases}	
\end{eqnarray*}

\end{longtable}

\begin{center}
\textbf{Operators, semigroups}
\end{center}
\vspace{-4mm}

\begin{longtable}{R p{.85\textwidth} } 	
	$L(X,U)$  & Space of bounded linear operators from $X$ to $U$ \\
	$L(X)$  & $:=L(X,X)$\\
	$A^*$ & Adjoint operator of an operator $A:D(A)\to X$, $D(A)\subset X$, on a Hilbert space $X$\\
	$s(A)$ & The spectral bound of a closed operator $A$\\
	$\sigma(A)$   & The spectrum of a closed operator $A$\\
	$\rho(A)$     & The resolvent set of a closed operator $A$\\
	$\omega_0(T)$ & The growth bound of a $C_0$-semigroup $T$\\
	$\omega_1(T)$ & The growth bound on $D(A)$ of a $C_0$-semigroup $T$ generated by $A$
\end{longtable}

\begin{center}
\textbf{Sequence spaces}
\end{center}
\vspace{-4mm}

\begin{longtable}{R p{.85\textwidth} } 
	$(x_k)$  & Shorthand notation for a sequence $(x_k)_{k\in\N}$ \\	
	$e_k$  & $k$-th standard basis vector of $\ell_p$, for some $p\in[1,+\infty)$\\	
  $\ell_p$  & Space of the sequences $x=(x_k)$ so that $\|x\|_{\ell_p}<\infty$\\
	$\|x\|_{\ell_p}$ &	$:=\Big(\sum_{k=1}^\infty |x_k|^p\Big)^{1/p}$, for $p\in[1,+\infty)$\\											
	$\|x\|_{\ell_\infty}$ &	$:=\sup_{k\in\N}|x_k|$\\											
	$c_0$  & $:=\{ x=(x_k): x_k\in\R\ \forall k\in\N \; \wedge\; \lim_{k \to \infty} x_k=0\}$ \\
  $c$  & $:=\{ x=(x_k): \exists M>0:\ x_k\in\R\; \wedge\; |x_k|\leq M \ \forall k\in\N\}$\\
  $\ell_{\infty}(I)$ & For a nonempty set $I$ we denote by $\ell_{\infty}(I)$ the Banach space of all functions 
$x:I \rightarrow \R$ with $\|x\|_{\ell_{\infty}(I)} := \sup_{i\in I}|x(i)| < \infty$.\\
  $\ell_{\infty}^+(I)$ & $:= \{ x \in \ell_{\infty}(I) : x(i) \geq 0,\ \forall i \in I \}$.
\end{longtable}

\begin{center}
\textbf{Spaces of continuous functions}
\end{center}
\vspace{-4mm}

\begin{longtable}{R p{.85\textwidth} } 
	$C(X,U)$  & Linear space of continuous functions from $X$ to $U$. \\
	          & We associate with this space the map  $u\mapsto \|u\|_{C(X,U)}:=\sup\limits_{x \in X}\|u(x)\|_U$ \\
	          & That may attain $+\infty$ as a value \\
	$PC_b(\R_+,U)$  & Space of piecewise continuous (right-continuous) functions from $\R_+$ to $U$ with a finite norm
	$\|u\|_{PC_b(\R_+,U)}=\|u\|_{C(\R_+,U)}<\infty$ \\
	$AC(\R_+,U)$  & Space of absolutely continuous functions from $\R_+$ to $U$ with $\|u\|_{C(X,U)}<\infty$ \\
	$C(X)$  &  $:=C(X,X)$\\
  $C_0(\R)$  & $:=\{f \in C(\R): \forall \varepsilon >0 \text{ there is a compact set } K_{\varepsilon} \subset \R: |f(s)|<\varepsilon \; \forall s \in \R \backslash K_{\varepsilon}\}$  \\
  $C_c^k(U)$  & Space of $k \in\N \cup\{\infty\}$ times continuously differentiable functions $f:U \to \R$ with a support, compact in $U$
\end{longtable}

\begin{center}
\textbf{Spaces of measurable functions}
\end{center}
\vspace{-4mm}

Here $U$ denotes an open subset of $\R^n$.

\begin{longtable}{R p{.85\textwidth} } 
  $\|f\|_{L^p(U,\R^m)}$ & $:=\left( \int_U{|f(x)|^p dx} \right)^{\frac{1}{p}}$, for $p\in[1,+\infty)$ \\		
  $\|f\|_{L^\infty(U,\R^m)}$ & $:=\esssup_{x\in U}|f(x)| = \inf_{D \subset U,\ \mu(D)=0} \sup_{x \in U \backslash D} |f(x)|$  \\
	$\|f\|_{\infty}$		& $:=\|f\|_{L^\infty(U,\R^m)}$ \\
  $L^{p}(U,\R^m)$  & The set of Lebesgue measurable maps from $U$ to $\R^m$ with $\|f\|_{L^p(U,\R^m)}<\infty$, for $p\in[1,+\infty]$ \\
  $L^{p}_{\loc}(U,\R^m)$  & $:= \{u:U \to\R^m: u\in L^{p}(V,\R^m) \ \forall V\subset\subset U \}$ , for $p\in[1,+\infty]$.\\
	& The set of Lebesgue measurable and locally essentially bounded functions from $U$ to $\R^m$ \\
  $L^{p}(U)$  & $:=L^{p}(U,\R)$, where $p\in[1,+\infty]$\\
	$L^p(a,b)$  & $:=L^p([a,b]) = L^p([a,b],\R)$, where $p\in[1,+\infty]$ \\
	$\|u\|_{W^{k,p}(U)}$ &:= $\Big(\sum_{|\alpha| \leq k}\int_U{\left| D^\alpha u (x)\right|^p dx} \Big)^{\frac{1}{p}}$, where $k\in\N$, $p\in[1,+\infty)$\\
	$\|u\|_{W^{k,\infty}(U)}$ &:= $\sum_{|\alpha| \leq k}\esssup_{x\in U}{\left| D^\alpha u (x)\right|}$, where $k\in\N$\\
	$W^{k,p}(U)$ & With $k\in\N$, $p\in[1,+\infty]$ is a Sobolev space of functions $u \in L^p(U)$, such that for each multiindex $\alpha$ of order  $\leq k$, the weak derivative $D^\alpha u$ exists and belongs to $L^p(U)$. $W^{k,p}(U)$ is endowed with the norm $\|u\|_{W^{k,p}(U)}$\\
	$W^{k,p}_0(U)$ & The closure of $C_c^\infty(U)$ in the norm of $W^{k,p}(U)$,\quad $k\in\N$, $p\in[1,+\infty]$\\
	$H^k(U)$& $:=W^{k,2}(U)$,\quad $k\in\N$\\
	$H^k_0(U)$& $:=W^{k,2}_0(U)$,\quad $k\in\N$\\
	$H^k(a,b)$  & $:=H^k([a,b]) = H^k([a,b],\R)$,\quad $k\in\N$ \\
	$H^k_0(a,b)$  & $:=H^k_0([a,b]) = H^k_0([a,b],\R)$,\quad $k\in\N$ \\

\end{longtable}

\begin{center}
\textbf{Spaces of vector-valued functions}
\end{center}
\vspace{-4mm}

Let $X$ be a Banach space, and let $I$ be a closed subset of $\R$. We define the following spaces of vector-valued functions

\begin{longtable}{R p{.85\textwidth} } 
$M(I,X)$               &$:= \{f: I \to X: f \text{ is strongly measurable} \}$\\
$L^p(I,X)$             &$:= \Big\{f \in M(I,X): \|f\|_{L^p(I,X)}:=\Big(\int_I\|f(s)\|^p_Xds\Big)^{\frac{1}{p}} < \infty \Big\}$\\
$L^p_{\loc}(\R_+,X)$   &$:= \{f \in L^p([0,t],X)\quad\forall t>0\}$\\
$L^\infty(I,X)$        &$:= \{f \in M(I,X): \|f\|_{L^\infty(I,X)}:=\esssup_{s\in I}\|f(s)\|_X < \infty \}$\\
$L^\infty_{\loc}(I,X)$ &$:= \{f \in L^\infty([0,t],X)\quad\forall t>0\}$
\end{longtable}

\begin{center}
\textbf{Comparison functions}
\end{center}
\vspace{-4mm}

\begin{longtable}{R p{.85\textwidth} } 
  $\PD$  & $:=\left\{\gamma \in C(\R_+):  \gamma(0)=0 \mbox{ and } \gamma(r)>0 \mbox{ for } r>0 \right\}$ \\
	$\K$  & $:=\left\{\gamma \in \PD : \gamma \mbox{ is strictly increasing}  \right\}$\\
	$\K_{\infty}$ & $:=\left\{\gamma\in\K: \gamma\mbox{ is unbounded}\right\}$\\
	$\LL$  & $:=\{\gamma \in C(\R_+): \gamma\mbox{ is strictly decreasing with } \lim\limits_{t\rightarrow\infty}\gamma(t)=0 \}$ \\
	$\KL$  & $:=\left\{\beta \in C(\R_+\times\R_+, \R_+): \beta(\cdot,t)\in{\K},\ \forall t \geq 0, \  \beta(r,\cdot)\in {\LL},\ \forall r > 0\right\}$\\		
\end{longtable}

\cleardoublepage

\chapter*{Preface}

\addcontentsline{toc}{chapter}{\numberline{}Preface}

 Nonlinear distributed parameter systems with both distributed and boundary inputs are used to model a broad range of phenomena, including chemical reactors, fluid and gas dynamics, traffic networks, 
multi-body systems (e.g., robotic arms, flexible elements), 
adaptive optics, 
fluid-structure interactions (e.g., dynamics of aircraft wings), etc.
For many classes of such systems, it is known that small disturbances coming from actuator and observation errors, hidden dynamics, and external disturbances can dramatically reduce the performance, alter the stability, or even destabilize the control system. 
Counteracting to these challenges requires the development of methods for the design of robust controllers and observers for nonlinear distributed parameter systems that ensure the reliability and efficiency of closed-loop systems.
This objective has to be achieved despite the fact that usually infinite-dimensional systems have to be controlled using finitely many (and usually very few) actuators and sensors that can be accessed possibly only at some discrete moments of time and which can frequently be placed merely at the boundary of the spatial domain.

A systematic framework that solves the counterparts of these problems for finite-dimensional nonlinear systems was developed around the notion of \emph{input-to-state stability (ISS)}, introduced by E. Sontag at the dawn of the 1990s. 
ISS combines two different types of stability
behavior: stability in the sense of Lyapunov and input-output stability.
The unified treatment of external and internal stability has
made ISS a central tool in robust stability analysis. ISS plays a vital role in many fields of nonlinear control, including robust stabilization of
nonlinear systems, stabilization via controllers with
saturation, design of robust nonlinear observers, nonlinear
detectability, ISS feedback redesign, stability of nonlinear networked control systems,
supervisory adaptive control, and others.

In this work, we develop an input-to-state stability theory for a broad class of infinite-dimensional systems encompassing partial differential equations (PDEs), time-delay systems, ordinary differential equations (ODEs), ensembles, as well as interconnections of heterogeneous systems consisting of an arbitrary number of finite- and infinite-dimensional components, with both in-domain and boundary couplings.
On the one hand, this provides a solid basis to treat the questions in robust control and observation of infinite-dimensional systems summarized above.
On the other hand, the interest in such a general theory is driven by the desire to create an overarching framework that provides a unifying view of ISS theories for various classes of systems.

To cover such a broad class of distributed parameter systems, we start in Chapter~\ref{chap:Intro} with a definition of a control system as a set of trajectories parametrized by initial conditions and external inputs. Besides its generality, an advantage of this approach is that by tacitly assuming the existence and uniqueness of solutions, we can concentrate on stability issues and prove very general results, such as characterizations of ISS in terms of other stability properties and Lyapunov theorems.
Next, we introduce in Chapter~\ref{chap:Intro} the central notion of input-to-state stability and show how to use Lyapunov functions to verify ISS of several linear and nonlinear hyperbolic and parabolic systems with distributed and boundary inputs. 

In Chapter~\ref{chap:Nonlinear Evolution Equations}, we develop well-posedness theory for linear and nonlinear evolution equations in Banach spaces with distributed and boundary inputs. 
First, we characterize well-posedness and ISS of general linear systems with bounded and admissible input operators based on the methods of the semigroup and admissibility theory. 
We proceed with a thorough well-posedness analysis of nonlinear evolution equations in Banach spaces with Lipschitz continuous nonlinearities and unbounded input operators.
This is an important class of control systems that includes semilinear evolution PDEs with distributed and boundary controls, ODEs, ensembles (infinite ODE networks), important subclasses of delay systems, and switched systems with a finite and infinite number of modes.
We analyze in detail the existence and uniqueness of local and global solutions as well as the regularity of the flow map.

We explain how these results can be applied to the ISS analysis of semilinear boundary control systems. Finally, we analyze semilinear evolution equations governed by analytic semigroups, and derive more precise well-posedness criteria in this case. 
The results shown in Chapter~\ref{chap:Nonlinear Evolution Equations} set a basis for the subsequent stability analysis of this system class.


In Chapter~\ref{chap:Characterizations_ISS_1}, we prove general ISS superposition theorems for input-to-state stability of a broad class of infinite-dimensional systems. These theorems give criteria for ISS in terms of weaker stability properties that significantly improve our understanding of input-to-state stability. Furthermore, they help to derive such essential results as the small-gain theorems in terms of trajectories, non-coercive direct ISS Lyapunov theorems, and integral characterizations of the ISS property. 
Moreover, in this chapter, we introduce the concept of strong input-to-state stability, which is equivalent to ISS in the case of finite-dimensional systems but is weaker than ISS in the infinite-dimensional case.
By means of counterexamples, we explain the differences between finite-dimensional and infinite-dimensional ISS theories.

In Chapter~\ref{chap:Characterizations_ISS_2}, we develop Lyapunov theory for ISS of evolution equations with Lipschitz continuous nonlinearities, whose well-posedness we analyzed in Chapter~\ref{chap:Nonlinear Evolution Equations}.
We introduce the notion of non-coercive ISS Lyapunov functions. We show that the existence of such a function implies ISS, provided that the flow of the system has a certain additional regularity. Here we essentially use the ISS superposition theorems derived in Chapter~\ref{chap:Characterizations_ISS_1}. 
Next, we show that local ISS of general nonlinear systems is equivalent to the uniform local asymptotic stability of an undisturbed system, provided that the nonlinearity satisfies a certain mild regularity requirement.
The developed Lyapunov tools are of central importance for the ISS analysis since, as a rule, the construction of an ISS Lyapunov function is the only realistic way to prove ISS. However, before we start searching for a suitable ISS Lyapunov function, it is important to know that such a function exists. This chapter shows that this is always the case if the nonlinearity is locally Lipschitz continuous on bounded sets in both variables.

The results shown in Chapters~\ref{chap:Characterizations_ISS_1}, \ref{chap:Characterizations_ISS_2} give strong tools for ISS analysis of semilinear PDEs with distributed disturbances. These results naturally extend the classical ISS theory for finite-dimensional systems. However, there are important system classes that are not covered by the systems considered in Chapter~\ref{chap:Characterizations_ISS_2}. Linear boundary control systems are one of such classes that is frequently encountered in the problems of robust control and observer design for linear PDEs with boundary controls that are subject to the actuator and/or measurement disturbances.

	


Despite all advantages of the ISS framework, for some practical systems, input-to-state stability is too restrictive.
ISS excludes systems whose state stays bounded as long as the magnitude of applied inputs and initial states remains below a specific threshold but becomes unbounded when either the input magnitude or the magnitude of an initial state exceeds the threshold. Such behavior is frequently caused by saturation and limitations in actuation and processing rate. 
In Chapter~\ref{chap:Integral input-to-state stability}, we introduce the concept of integral input-to-state stability that captures such nonlinearities. We discuss the available results in the integral ISS theory of linear and bilinear systems and present Lyapunov tools for analysis of iISS for nonlinear parabolic PDEs.

Stability analysis and control of nonlinear systems is a complex task that becomes even more challenging for large-scale systems in view of the nonlinearity of its subsystems, complex topology, and the sheer size of the network.
A prolific approach to study such networks is to overapproximate them by an infinite network and analyze it as a worst-case scenario. 
In Chapter~\ref{chap:Infinite interconnections: Non-Lyapunov methods}, we present a powerful method that helps to analyze such complex networks, provided that their components are ISS and the asymptotic gains describing the influence of the subsystems at each other are known. We present a general framework for the analysis of networks of abstract control systems and derive very general small-gain theorems for stability analysis of infinite networks of infinite-dimensional systems with heterogeneous components. Our stability criteria for the network are based on the properties of monotone nonlinear discrete-time systems induced by so-called gain operators.

In most cases, the ISS of nonlinear systems is verified by constructing an appropriate ISS Lyapunov function.  
It sounds natural to use the information about the ISS Lyapunov functions for subsystems in the formulation of the small-gain theorems instead of using the ISS estimates for subsystems.
In Chapter~\ref{chap:Interconnections}, we show that this is possible. We analyze the stability of finite networks of nonlinear infinite-dimensional systems using the Lyapunov-based small-gain theorems for couplings of ISS and integral ISS systems. We apply our theoretical developments to the analysis of coupled parabolic systems. 


\bigskip

Prior knowledge of ISS theory is not required from a reader. Basic facts about comparison functions and ISS theory of ODEs are collected in Appendices~\ref{chap:Finite-dim_ISS_Theory}, \ref{chap:Comparison functions and principles}. For a detailed treatise of finite-dimensional ISS theory, we refer to \cite{Mir23}.
However, we assume that a reader has basic knowledge of ordinary and partial differential equations as well as of
linear functional analysis and semigroup theory. Some basic definitions and results from the theory of ordered Banach spaces and from PDE theory are summarized in Appendices~\ref{section:ordered-banach-spaces}, \ref{chap:Function_Spaces_Inequalities}.



Please send any corrections and/or comments, no matter how small, to the email:\\
\href{andersmir@gmail.com}{andersmir@gmail.com}
~\\

\begin{flushright}
Andrii Mironchenko\\
Passau, Germany \\
\date{\today}

\end{flushright}

\cleardoublepage
\chapter{Introduction}
\label{chap:Intro}

This chapter defines the notion of input-to-state stability (ISS). Furthermore, we introduce ISS Lyapunov functions that constitute the primary tool for the ISS analysis, at least in the context of nonlinear systems. We show that the existence of an ISS Lyapunov function implies ISS.
Afterward, we show how the ISS can be verified by means of the Lyapunov method for several basic partial differential equations 
(PDE), such as transport equation, reaction-diffusion equation with distributed and boundary controls, Burgers' equation, etc.
In this introductory chapter, we use predominantly quadratic ISS Lyapunov functions and utilize integral inequalities in Sobolev spaces to verify dissipative inequalities for these functions.
Deeper methods for ISS analysis will be developed in the subsequent chapters.

\section{General class of systems}

We start with a general definition of a control system.
\index{control system}
\begin{definition}
\label{Steurungssystem}
Consider the triple $\Sigma=(X,\Uc,\phi)$ consisting of 
\index{state space}
\index{space of input values}
\index{input space}
\begin{enumerate}[label=(\roman*)]  
    \item A normed vector space $(X,\|\cdot\|_X)$, called the \emph{state space}, endowed with the norm $\|\cdot\|_X$.
    \item A normed vector \emph{space of inputs} $\Uc \subset \{u:\R_+ \to U\}$          
endowed with a norm $\|\cdot\|_{\Uc}$, where $U$ is a normed vector \emph{space of input values}.
We assume that the following two axioms hold:
                    
\emph{The axiom of shift invariance}: for all $u \in \Uc$ and all $\tau\geq0$ the time
shift $u(\cdot + \tau)$ belongs to $\Uc$ with \mbox{$\|u\|_\Uc \geq \|u(\cdot + \tau)\|_\Uc$}.

\emph{The axiom of concatenation}: for all $u_1,u_2 \in \Uc$ and for all $t>0$ the \emph{concatenation of $u_1$ and $u_2$ at time $t$}, defined by
\begin{equation}
\ccat{u_1}{u_2}{t}(\tau):=
\begin{cases}
u_1(\tau), & \text{ if } \tau \in [0,t], \\ 
u_2(\tau-t),  & \text{ otherwise},
\end{cases}
\label{eq:Composed_Input}
\end{equation}
belongs to $\Uc$.

    \item A map $\phi:D_{\phi} \to X$, $D_{\phi}\subseteq \R_+ \times X \times \Uc$ (called \emph{transition map}), such that for all $(x,u)\in X \tm \Uc$ it holds that $D_{\phi} \cap \big(\R_+ \times \{(x,u)\}\big) = [0,t_m)\tm \{(x,u)\} \subset D_{\phi}$, for a certain $t_m=t_m(x,u)\in (0,+\infty]$.
		
		The corresponding interval $[0,t_m)$ is called the \emph{maximal domain of definition} of $t\mapsto \phi(t,x,u)$.
		
\end{enumerate}
The triple $\Sigma$ is called a \emph{(control) system}, if the following properties hold:
\index{property!identity}

\begin{sysnum}
    \item\label{axiom:Identity} \emph{The identity property:} for every $(x,u) \in X \times \Uc$
          it holds that $\phi(0, x,u)=x$.
\index{causality}
    \item \emph{Causality:} for every $(t,x,u) \in D_\phi$, for every $\tilde{u} \in \Uc$, such that $u(s) =
          \tilde{u}(s)$ for all $s \in [0,t]$ it holds that $[0,t]\tm \{(x,\tilde{u})\} \subset D_\phi$ and $\phi(t,x,u) = \phi(t,x,\tilde{u})$.
    \item \label{axiom:Continuity} \emph{Continuity:} for each $(x,u) \in X \times \Uc$ the map $t \mapsto \phi(t,x,u)$ is continuous on its maximal domain of definition.
\index{property!cocycle}
        \item \label{axiom:Cocycle} \emph{The cocycle property:} for all
                  $x \in X$, $u \in \Uc$, for all $t,h \geq 0$ so that $[0,t+h]\tm \{(x,u)\} \subset D_{\phi}$, we have
\[
\phi\big(h,\phi(t,x,u),u(t+\cdot)\big)=\phi(t+h,x,u).
\]
\end{sysnum}

\end{definition}

To relate this definition to other concepts available in the literature, consider the following concept:
\begin{definition}
\label{def:strongly continuous nonlinear semigroup}
\index{semigroup!nonlinear}
Let $T(t) : X \to X$, $t \ge 0$, be a family of nonlinear maps.
The family $T:=\{ T(t): t\ge 0\}$ is called a \emph{strongly continuous nonlinear semigroup} if the following holds:
\begin{itemize}
	\item[(i)] For all $x\in X$ it holds that $T(0)(x) =x$.
	\item[(ii)] For all $t_1,t_2\geq 0$ and all $x \in X$ it holds that $T(t_1)\big(T(t_2)(x)\big) = T(t_1+t_2)(x)$.
	\item[(iii)] For each $x\in X$ the map $t\mapsto T(t)(x)$ is continuous.
\end{itemize}
\end{definition}

Take the family $\sg{T}$ as in Definition~\ref{def:strongly continuous nonlinear semigroup}.
Setting $\Uc:=\{0\}$, and defining $\phi(t,x,0):=T(t)(x)$, one can see that $\Sigma:=(X,\{0\},\phi)$ is a control system according to
Definition~\ref{Steurungssystem}.
Indeed, the axioms (i)-(iii) of Definition~\ref{Steurungssystem}, as well as causality, trivially hold. The identity property, continuity of $\phi$ w.r.t.\ time, as well as the cocycle property correspond directly to the axioms (i)-(iii) of
Definition~\ref{def:strongly continuous nonlinear semigroup}.

Abstract linear control systems, considered in Section~\ref{sec:Abstract linear control systems}, are also a special case of the control systems that we consider here.
In a certain sense, Definition~\ref{Steurungssystem} is a direct generalization and a unification of the concepts of strongly continuous nonlinear semigroups with abstract linear control systems.

This class of systems encompasses control systems generated by ordinary
differential equations (ODEs), switched systems, time-delay systems,
evolution PDEs, abstract differential
equations in Banach spaces and many others \cite[Chapter 1]{KaJ11b}.

\section[Forward completeness and BRS]{Forward completeness and boundedness of reachability sets}

\begin{definition}
\label{def:FC_Property} 
\index{forward completeness}
We say that a control system (as introduced in Definition~\ref{Steurungssystem}) is \emph{forward complete (FC)}, if 
$D_\phi = \R_+ \tm X\tm\Uc$, that is for every $(x,u) \in X \times \Uc$ and for all $t \geq 0$ the value $\phi(t,x,u) \in X$ is well-defined.
\end{definition}

Forward completeness alone does not imply, in general, the existence of any uniform bounds on the trajectories emanating from bounded balls that are subject to uniformly bounded inputs (see Example~\ref{0-GAS_but_not_GS}). Systems exhibiting such bounds deserve a special name.
\begin{definition}
\label{def:BRS}
\index{bounded reachability sets}
\index{BRS}
We say that \emph{$\Sigma=(X,\Uc,\phi)$ has bounded reachability sets (BRS)}, if for any $C>0$ and any $\tau>0$ it holds that 
\[
\sup\big\{
\|\phi(t,x,u)\|_X : \|x\|_X\leq C,\ \|u\|_{\Uc} \leq C,\ t \in [0,\tau]\big\} < \infty.
\]
\end{definition}

It is useful to have a quantitative restatement of the BRS property in a comparison-functions-like manner.

We call a function $h: \R_+^3 \to \R_+$ increasing, if $(r_1,r_2,r_3) \leq (R_1,R_2,R_3)$
implies that $h(r_1,r_2,r_3) \leq h(R_1,R_2,R_3)$, where we use the component-wise
partial order on $\R_+^3$. We call $h$ strictly increasing if $(r_1,r_2,r_3)
\leq (R_1,R_2,R_3)$ and $(r_1,r_2,r_3) \neq (R_1,R_2,R_3)$ imply $h(r_1,r_2,r_3) <
h(R_1,R_2,R_3)$.

\begin{lemma}
\label{lem:Boundedness_Reachability_Sets_criterion}
Consider a forward complete system $\Sigma$. The following statements are equivalent:
\begin{enumerate}
    \item[(i)] $\Sigma$ has bounded reachability sets.
    \item[(ii)] There exists a continuous, increasing function $\mu: \R_+^3 \to \R_+$, such that for
all $x\in X, u\in \Uc$ and all $t \geq 0$ we have
 \begin{equation}
    \label{eq:8_ISS}
    \| \phi(t,x,u) \|_X \leq \mu( \|x\|_X,\|u\|_{\Uc},t).
\end{equation}
    \item[(iii)] There exists a continuous function $\mu: \R_+^3 \to \R_+$ such that for
all $x\in X, u\in \Uc$ and all $t \geq 0$ the inequality \eqref{eq:8_ISS} holds.
\end{enumerate}
\end{lemma}

\begin{proof}
(i) $\Rightarrow$ (ii). Define $\tilde\mu: \R_+^3 \to \R_+$ by
\begin{eqnarray}
\label{eq:tildechi-def_ISS_section}
\tilde\mu(C_1,C_2,\tau) := \sup_{\|x\|_X\leq C_1,\: \|u\|_\Uc \leq C_2,\: t \in [0,\tau]} \|\phi(t,x,u)\|_X,
\end{eqnarray}
which is well-defined in view of the item (i). Clearly, $\tilde\mu$ is increasing by definition. In particular, it is
locally integrable.

Define $\hat\mu: (0,+\infty)^3 \to \R_+$ by setting for $C_1,C_2,\tau>0$
\begin{eqnarray*}
\hat\mu(C_1,C_2,\tau) := \frac{1}{C_1C_2\tau} \int_{C_1}^{2C_1}\hspace{-2mm}\int_{C_2}^{2C_2}\hspace{-2mm}\int_{\tau}^{2\tau}\hspace{-2mm} \tilde\mu(r_1,r_2,s) dr_1dr_2ds + C_1C_2\tau.
\end{eqnarray*}
By construction, $\hat\mu$ is strictly increasing and continuous on $(0,+\infty)^3$.

We can enlarge the domain of definition of $\hat\mu$ to all of $\R^3_+$
using monotonicity. To this end we define for $C_2,\tau>0$: $\hat\mu(0,C_2,\tau) := \lim_{C_1\to +0}
\hat\mu(C_1,C_2,\tau)$, for $C_1\geq0,\ \tau>0$: $\hat\mu(C_1,0,\tau) := \lim_{C_2\to +0}
\hat\mu(C_1,C_2,\tau)$ and for $C_1,C_2\geq 0$ we define $\hat\mu(C_1,C_2,0) := \lim_{\tau\to +0}
\hat\mu(C_1,C_2,\tau)$. All these limits are well-defined as $\hat\mu$ is increasing on
$(0,+\infty)^3$, and we obtain that the resulting function is increasing on
$\R_+^3$. Note that the construction does not guarantee that $\hat\mu$ is
continuous. To obtain continuity, choose a continuous strictly
increasing function $\nu: \R_+ \to \R_+$ with 
\[
\nu(r) > \max\{ \hat\mu(0,0,r), \hat\mu(0,r,0), \hat\mu(r,0,0)\},\quad r \geq 0,
\]
and define for $(C_1,C_2,\tau) \geq (0,0,0)$
\begin{eqnarray*} 
\mu(C_1,C_2,\tau) :=
      \max\big\{
  \nu\big(\max\{C_1,C_2,\tau\}\big),\hat\mu(C_1,C_2,\tau) \big\} + C_1C_2\tau.
\end{eqnarray*}
The function $\mu$ is continuous as 
\[
\mu(C_1,C_2,\tau) = \nu\big(\max\{C_1,C_2,\tau\}\big) + C_1C_2\tau
\]
whenever $C_1$, $C_2$ or $\tau$ is small enough.
At the same time we have for $C_1>0$, $C_2>0$ and $\tau >0$ that
\begin{align*}
\mu(C_1,C_2,\tau) &\geq \hat\mu(C_1,C_2,\tau)\\
& \geq
\frac{1}{C_1C_2\tau} \int_{C_1}^{2C_1}\int_{C_2}^{2C_2}\int_{\tau}^{2\tau} 1 dr_1dr_2ds \ \tilde\mu(C_1,C_2,\tau) + C_1C_2\tau \\
&\geq \tilde\mu(C_1,C_2,\tau).
\end{align*}
This implies that (ii) holds with this $\mu$.

(ii) $\Rightarrow$ (iii) is evident.

(iii) $\Rightarrow$ (i) follows due to continuity of $\mu$.
\end{proof}

For a wide class of control systems boundedness of a solution implies the possibility of prolonging it to a larger interval, see \cite[Chapter 1]{KaJ11b}. Next, we formulate this property for abstract systems:
\begin{definition}
\label{def:BIC} 
\index{property!boundedness-implies-continuation}
\index{BIC}
We say that a system $\Sigma$ satisfies the \emph{boundedness-implies-continuation (BIC) property} if for each
$(x,u)\in X \tm \Uc$ such that the maximal existence time $t_m(x,u)$ is finite, and for all $M > 0$, there exists $t \in [0,t_m(x,u))$ with $\|\phi(t,x,u)\|_X>M$.
\end{definition}

\section{Input-to-state stability}

For the formulation of stability properties, we use the standard classes of comparison functions, for which we refer to
Appendix~\ref{chap:Comparison functions and principles}. 
We start with the main concept for this work:
\begin{definition}
\label{Def:ISS}
\index{stability!input-to-state}
\index{input-to-state stability}
\index{ISS}
A system $\Sigma=(X,\Uc,\phi)$ is called \emph{ (uniformly)  input-to-state stable
(ISS)}, if it is forward complete and there exist $\beta \in \KL$ and $\gamma \in \K$
such that for all $ x \in X$, $ u\in \Uc$ and $ t\geq 0$ it holds that
\begin {equation}
\label{iss_sum}
\| \phi(t,x,u) \|_{X} \leq \beta(\| x \|_{X},t) + \gamma( \|u\|_{\Uc}).
\end{equation}
\end{definition}

By $B_{r,\; Z}$  we denote the open ball of radius $r$ around $0$ in a normed vector space $Z$, i.e., $\{u\in Z:\|u\|_{\Uc}<r\}$, and 
 we set for short $B_r:=B_{r,\ X}$.

The local counterpart of the ISS property is
\begin{definition}
\label{Def:LISS}
\index{stability!local input-to-state}
\index{input-to-state stability!local}
\index{LISS}
A system $\Sigma=(X,\Uc,\phi)$ is called \emph{ (uniformly) locally input-to-state stable
(LISS)}, if there exist $\beta \in \KL$, $\gamma \in \K$
and $r>0$ such that for all $ x \in
\clo{B_r}$, $u\in \clo{B_{r,\;\Uc}}$ we have $t_m(x,u)=\infty$, and the inequality \eqref{iss_sum} holds  for all $ t\geq 0$.
\end{definition}

ISS trivially implies the following property, which plays a key role in the stability analysis of the systems without disturbances.
\begin{definition}
\label{def:0-UGAS}
\index{0-UGAS}
\index{stability!uniform global asymptotic at zero}
A system $\Sigma=(X,\Uc,\phi)$ is called \emph{uniformly globally asymptotically stable at zero (0-UGAS)}, if there exists $\beta \in \KL$ such that for all $ x \in X$ we have $t_m(x,0)=+\infty$, and
\begin {equation}
\label{eq:0-UGAS}
\| \phi(t,x,0) \|_{X} \leq \beta(\| x \|_{X},t),\quad t\geq 0.
\end{equation}
\end{definition}

Finally, a stronger than ISS property of exponential ISS is of interest:
\begin{definition}
\label{Def:eISS}
\index{stability!exponential input-to-state}
\index{input-to-state stability!exponential}
\index{eISS}
System $\Sigma=(X,\Uc,\phi)$ is called \emph{ exponentially  input-to-state stable
(eISS)}, if it is forward complete and there exist $M, a>0 $ and $\gamma \in \K$
such that for all $ x \in X$, $ u\in \Uc$ and $ t\geq 0$ it holds that
\begin {equation}
\label{eq:eISS_sum}
\| \phi(t,x,u) \|_{X} \leq Me^{-at}\|x\|_X + \gamma(\|u\|_{\Uc}).
\end{equation}
Furthermore, if in addition $\gamma$ can be chosen as a linear function, then we call $\Sigma$ \emph{exponentially ISS with a linear gain}.
\end{definition}

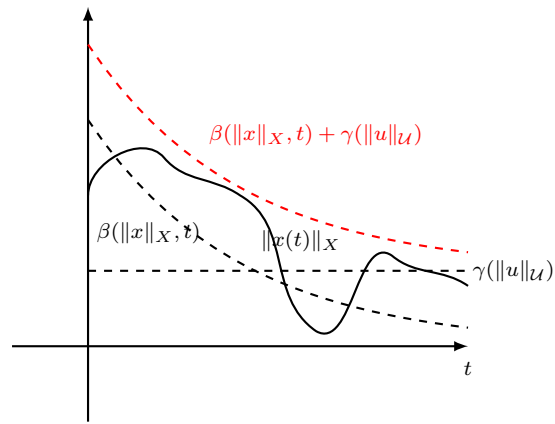
\begin{figure}[ht]
    \centering

\begin{tikzpicture}
\draw[->,thick] (-1,0) --(5,0);
\draw[->,thick] (0,-1) --(0,4.5);
\draw[dashed,thick] (0,1) --(5,1);

\node(gain) at (5.6,1) {\scriptsize $\gamma(\|u\|_{\Uc})$};
\node(trans) at (0.8,1.5) {\scriptsize $\beta(\|x\|_X,t)$};
\node[red] (ISS_sum) at (3,2.8) {\scriptsize $\beta(\|x\|_X,t)+\gamma(\|u\|_{\Uc})$};

\node(x) at (2.8,1.4) {\scriptsize $\|x(t)\|_X$};
\node(t) at (5,-0.3) {\scriptsize $t$};

\draw[-,thick] (0,2) to [out=90,in=130](1,2.5) to [out=-50,in=145](2,2)
to [out=-30,in=145](3,0.2) to [out=-30,in=145](4,1.2)
to [out=-40,in=145](5,0.8);

\draw[dashed, thick,  domain=0:5]  plot (\x, {3*exp(-0.5*\x)});
\draw[dashed, thick, red, domain=0:5]  plot (\x, {1+3*exp(-0.5*\x)});
\end{tikzpicture}
\caption{Typical graph of a norm of a trajectory of an ISS system.}
\label{ISS_graph}
\end{figure}

One of the most common choices for $\Uc$ is the space $\Uc:=PC_b(\R_+,U)$, for a certain normed vector space $U$, with the norm $\| u \|_{\Uc} := \sup\limits_{0 \leq s \leq \infty} \|u(s)\|_U$.
In this case, one can use the alternative definition of the (L)ISS property, which is often used in the literature (see, e.g., \cite{KaJ11}, \cite{HLT08}):
\begin{proposition}
\label{prop:Causality_Consequence}
Let $\Uc:=PC_b(\R_+,U)$, for a certain normed vector space $U$. A system $\Sigma=(X,\Uc,\phi)$ is LISS iff
 there exist $\beta \in \KL$, $\gamma \in \K$
and $r>0$ such that for all $ x \in
\clo{B_r}$, $u\in \clo{B_{r,\;\Uc}}$ we have $t_m(x,u)=\infty$, and the following holds  for all $ t\geq 0$:
\begin {equation}
\label{iss_sum_equiv}
\| \phi(t,x,u) \|_{X} \leq \beta(\| x \|_{X},t) + \gamma( \sup\limits_{0 \leq s \leq t} \|u(s)\|_U).
\end{equation}
\end{proposition}

\begin{proof}
Sufficiency is clear, since $\sup\limits_{0 \leq s \leq t} \|u(s)\|_U \leq \sup\limits_{0 \leq s \leq \infty} \|u(s)\|_U=\|u\|_{\Uc}$.

Now let $\Sigma$ be LISS. Due to causality property of $\Sigma$ the state $\phi(\tau,x,u)$, $\tau \in [0,t]$ of the system $\Sigma$ does not depend on the values of $u(s)$, $s >t$. For arbitrary $t \geq 0$ and $u \in \Uc$ consider another input $\tilde{u} \in \Uc$, defined by
\[
\tilde{u}(\tau):=
\begin{cases}
u(\tau), & \text{ if } \tau \in [0,t], \\ 
u(t),  & \tau >t.
\end{cases}
\]
The inequality \eqref{iss_sum} holds for all admissible inputs, and hence it holds also for $\tilde{u}$. Substituting $\tilde{u}$ into \eqref{iss_sum} and using that
$\|\tilde{u}\|_{\Uc}= \sup\limits_{0 \leq s \leq t} \|u(s)\|_U$, we obtain \eqref{iss_sum_equiv}.
\end{proof}

It is easy to formulate a counterpart of Proposition~\ref{prop:Causality_Consequence} for the (global) ISS property.
A similar result holds for continuous input functions ($\Uc:=C(\R_+,U)$), for the class of strongly measurable and essentially bounded inputs $\Uc:=L^{\infty}(\R_+,U)$ (with $ \mathop{\text{ess} \sup}\limits_{0 \leq s \leq t} \|u(s)\|_U$ instead of $\sup\limits_{0 \leq s \leq t} \|u(s)\|_U$) and many other classes of input functions.

We close the section by showing that trajectories of ISS systems converge to zero for vanishing inputs.

\begin{definition}
\label{def:CIUCS} 
\index{property!convergent input-uniformly convergent state (CIUCS)} 
We say that a forward complete system $\Sigma=(X,\Uc,\phi)$  has a
 \emph{convergent input-uniformly convergent state (CIUCS) property}, 
if for each 
$u\in\Uc$ such that $\lim_{t \to \infty}\|u(t+\cdot)\|_\Uc = 0$, and for any $r>0$, it holds that
\[
\lim_{t \to \infty}\sup_{\|x\|_X\leq r}\|\phi(t,x,u)\|_X = 0.
\] 
\end{definition}

\begin{proposition}
\label{prop:Converging_input_uniformly_converging_state}
Every ISS control system $\Sigma=(X,\Uc,\phi)$ has the CIUCS property. 
\end{proposition}

\begin{proof}
Let $\Sigma=(X,\Uc,\phi)$ be an ISS control system and assume without loss of generality that the corresponding gain $\gamma$ is a $\Kinf$-function.
Pick any $r>0$, any $u\in \Uc$ so that $\lim_{t \to \infty}\|u(\cdot + t)\|_{\Uc} = 0$, and any $\eps>0$. To show the claim of the proposition, we need to show that there is a time $t_\eps=t_\varepsilon(r,u)>0$ so that
\[
\|x\|_X\leq r,\quad t\geq t_\eps \qrq \|\phi(t,x,u)\|_X\leq\eps.
\]

Choose $t_1$ so that $\|u(\cdot + t_1)\|_{\Uc}\leq \gamma^{-1}(\eps)$.
Due to the cocycle property and ISS of $\Sigma$, we have for any $x \in B_r$ that
\begin{eqnarray*}
\|\phi(t+t_1,x,u)\|_X  &=& \big\|\phi\big(t,\phi(t_1,x,u),u(\cdot + t_1)\big)\big\|_X \\
                                             &\leq& \beta\big(\|\phi(t_1,x,u)\|_X,t\big) + \gamma\big(\|u(\cdot+t_1)\|_\Uc\big) \\
                                             &\leq& \beta\big(\|\phi(t_1,x,u)\|_X,t\big) + \eps\\
                                             &\leq& \beta\big(\beta(\|x\|_X,t_1) + \gamma(\|u\|_\Uc),t\big) + \eps\\
                                             &\leq& \beta\big(\beta(r,0) + \gamma(\|u\|_\Uc),t\big) + \eps.
\end{eqnarray*}
Pick any $t_2$ in a way that $\beta\big(\beta(r,0) + \gamma(\|u\|_\Uc),t_2\big) \leq \eps$.
This ensures that $\|\phi(t_2+t_1,x,u)\|_X \leq 2\eps.$
Using the cocycle property once again, we obtain for all $t\geq0$:
\begin{align*}
\|\phi(t+t_2&+t_1,x,u)\|_X \\
&\leq \beta\big(\|\phi(t_2+t_1,x,u)\|_X,t\big) + \gamma\big(\|u(\cdot+t_1+t_2)\|_\Uc\big).
\end{align*}
Due to the axiom of shift invariance we have $\|u(\cdot+t_2+t_1)\|_\Uc \leq \|u(\cdot+t_1)\|_\Uc\leq \gamma^{-1}(\eps)$.
ISS of $\Sigma$ combined with the cocycle property ensures for all $t>0$ that
\begin{eqnarray*}
\|\phi(t+t_2+t_1,x,u)\|_X &\leq& \beta\big(\|\phi(t_2+t_1,x,u)\|_X,t\big) + \eps\\
&\leq& \beta\big(2\eps,0\big) + \eps.
\end{eqnarray*}
Since $\eps>0$ can be chosen arbitrarily small, and since $\beta\big(2\eps,0\big) + \eps \to 0$ as $\eps\to 0$, the claim of the proposition follows.
\end{proof}

\section{Density argument}
\label{sec:Density_argument}

Often it is easier to verify first the ISS of a distributed parameter system on certain dense subspaces of a state space and a space of inputs. For example, Lyapunov functionals may fail to be differentiable on the whole state space, while they are differentiable on certain dense subspaces. We will show in this section that under natural additional assumptions, the ISS of a system over the dense state and input space implies ISS on the whole state and input space.

\begin{definition}
\label{def:Continuous-dependence-systems} 
Let $\Sigma=(X,\Uc,\phi)$ be a forward complete control system. 
We say that the flow
\emph{$\phi$ depends continuously on inputs and on initial states}, if 
 for all $x \in X$, $u \in \Uc$, $T>0$, and all $\eps>0$ there exists $\delta>0$, such that $\forall x' \in X: \|x-x'\|_X< \delta$ and
$\forall u' \in \Uc: \|u-u'\|_{\Uc}< \delta$ it holds that
\[
\|\phi(t,x,u)-\phi(t,x',u')\|_X< \eps, \quad \forall t \in [0,T].
\]
\end{definition}

Now we have the following result
\begin{lemma}
\label{lem:Density_Arg_ISS}
Consider a forward complete control system $\Sigma=(X,\Uc,\phi)$ and let $\phi$ depend continuously on inputs and on initial states.

Let $\hat{X}$, $\hat{\Uc}$ be dense normed vector subspaces of $X$ and $\Uc$, endowed with the norms inherited from $X$ and $\Uc$, respectively, and let $\hat{\Sigma}:=(\hat{X},\hat{\Uc},\phi)$ be the system, generated by the transition map $\phi$ of $\Sigma$ but restricted to the state space $\hat{X}$ and space of admissible inputs $\hat{\Uc}$.

Let $\hat{\Sigma}$ be ISS. Then $\Sigma$ is also ISS with the same $\beta$ and $\gamma$
in the estimate \eqref{iss_sum}.
\end{lemma}

\begin{proof}
Since $\hat{\Sigma}$ is ISS, there exist $\beta \in \KL$ and $\gamma \in \Kinf$, such that $\forall \hat{x} \in \hat{X},\; \forall \hat{u} \in \hat{\Uc}$ and $\forall t\geq 0$ it holds that
\begin {equation}
\label{DichteArg_1}
\| \phi(t,\hat{x},\hat{u}) \|_{X} \leq \beta(\| \hat{x} \|_{X},t) + \gamma( \|\hat{u}\|_{\Uc}).
\end{equation}
Let $\Sigma$ be not ISS with the same $\beta$, $\gamma$. Then there exist $T>0$, $x\in X$, $u \in \Uc$ such that for a certain 
$r=r(T,x,u)>0$ it holds that
\begin {equation}
\label{DichteArg_2}
\| \phi(T,x,u) \|_{X} = \beta(\| x \|_{X},T) + \gamma( \|u\|_{\Uc}) + r.
\end{equation}

From \eqref{DichteArg_1} and \eqref{DichteArg_2}, we obtain
\begin{equation}
\label{DichteArg_3}
\begin{split}
\| \phi(T,x,u) \|_{X} - \| \phi(T,\hat{x},\hat{u}) \|_{X} \geq& \left(\beta(\| x \|_{X},T) - \beta(\| \hat{x} \|_{X},T) \right) \\
& \phantom{aaaa} +  \left(\gamma( \|u\|_{\Uc}) - \gamma( \|\hat{u}\|_{\Uc}) \right) + r.
\end{split}
\end{equation}
Since $\hat{X}$ and $\hat{\Uc}$ are dense in $X$ and $\Uc$, respectively, we can find sequences $(\hat{x}_i) \subset \hat{X}$: $\|x-\hat{x}_i \|_X \to 0$ and $(\hat{u}_i ) \subset \hat{\Uc}$: $\|u-\hat{u}_i \|_{\Uc} \to 0$ as $i\to\infty$.
From \eqref{DichteArg_3}, for all $\eps>0$ there exist $\hat{x}_i$ and $\hat{u}_i$:
\begin {equation}
\label{DichteArg_4}
\| \phi(T,x,u) - \phi(T,\hat{x}_i,\hat{u}_i) \|_{X} \geq
\big| \| \phi(T,x,u) \|_{X} - \| \phi(T,\hat{x}_i,\hat{u}_i) \|_{X}  \big| \geq r - 2\eps.
\end{equation}
This contradicts the assumption of continuous dependence of $\phi$ on initial states and inputs.
Thus, $\Sigma$ is ISS. %
\end{proof}

\section{ISS Lyapunov functions}
\label{sec:ISS_Lyapunov_functions}

For a continuous function $y:\R \to \R$, let
the \emph{right upper Dini derivative} and
the \emph{right lower Dini derivative} be defined by
$D^+y(t):=\Limsup\limits_{h \to +0}\frac{y(t+h)-y(t)}{h}$ and
$D^-y(t):=\Liminf\limits_{h \to +0}\frac{y(t+h)-y(t)}{h}$, respectively.
Furthermore, we use the following notation: $\dot{y}(t):=D^+y(t)$.
In Section~\ref{sec:Dini derivatives}, we give a short overview of the known results for Dini derivatives.

\index{derivative!Lie}
Let $x \in X$ and $V$ be a real-valued function defined in a neighborhood of $x$. The \emph{Lie derivative}
of $V$ at $x$ corresponding to the input $u$ along the trajectories of $\Sigma$ is defined by
\begin{equation}
\label{ISS_LyapAbleitung}
\dot{V}_u(x)=\Limsup \limits_{t \rightarrow +0} {\frac{1}{t}\big(V(\phi(t,x,u))-V(x)\big) }.
\end{equation}
In other words, $\dot{V}_u(x) = D^+V(\phi(t,x,u))\Big|_{t=0}$. 

Lyapunov functions are a powerful tool for investigating ISS and local ISS. 
\begin{definition}
\label{DefLISS_LF}
\index{function!LISS Lyapunov}
\index{function!ISS Lyapunov}
\index{function!eISS Lyapunov}
\index{ISS Lyapunov function}
\index{LISS Lyapunov function}
\index{eISS Lyapunov function}
Let $\Sigma = (X,\Uc,\phi)$ be a time-invariant control system.
A continuous function $V:D \to \R_+$, $D \subset X$, $0 \in \intt(D)=D \backslash \partial D$ is called a \emph{local ISS Lyapunov function (LISS LF) in an implication form for $\Sigma$}, if there exist $\rho_1,\rho_2>0$, $\psi_1,\psi_2 \in \Kinf$, $\chi \in \K$,  and positive definite function $\alpha$, such that $B_{\rho_1} \subset D$, and
\begin{equation}
\label{LyapFunk_1Eig}
\psi_1(\|x\|_X) \leq V(x) \leq \psi_2(\|x\|_X), \quad \forall x \in B_{\rho_1},
\end{equation}
and
$\forall x \in X: \|x\|_X \leq \rho_1, \; \forall u\in \Uc: \|u\|_{\Uc} \leq \rho_2$ it holds:
\begin{equation}
\label{GainImplikation}
 \|x\|_X \geq \chi(\|u\|_{\Uc}) \  \Rightarrow  \ \dot{V}_u(x) \leq -\alpha(V(x)).
\end{equation}
The function $\chi$ is called \emph{ISS Lyapunov gain} for $\Sigma$.

\begin{enumerate}[label=(\roman*)]
		\item If additionally $\alpha(r) = ar$, $a>0$, then $V$ is called a \emph{local exponential ISS (eISS) Lyapunov function in an implication form}.
	\item If in the previous definition $D=X$, $\rho_1=\infty$ and $\rho_2=\infty$, then the function $V$ is called an \emph{ISS Lyapunov function in an implication form}.
\end{enumerate}
\end{definition}

\begin{remark}
\label{rem:Computation-of-dotV} 
Our definition of an ISS Lyapunov function is very general. Yet the computation of $\dot{V}_u(x)$ by definition given in \eqref{ISS_LyapAbleitung} requires the knowledge of the solution for future times, at least on a sufficiently small time interval.

However, one of the main benefits of the Lyapunov theory is the possibility of checking the stability of a system without knowledge of a solution. 
 To retain this feature, for many important classes of systems, it is possible to derive 
useful formulas for the computation of $\dot{V}_u(x)$, which do not involve the knowledge of future dynamics.
For example, for ODE systems \eqref{xdot=f_xu} with a Lipschitz continuous $V$, for almost all $x\in\R^n$ and all continuous $u\in\Uc$ it holds that $\dot{V}_u(x) = \nabla V(x) \cdot f(x,u(0))$, where $\nabla V$ is the gradient of $V$ and $\cdot$ is the standard inner product in $\R^n$.
For time-delay systems, $\dot{V}_u(x)$ can be computed using the so-called Driver derivative.
For many PDE systems, one can give simple formulas for $\dot{V}_u(x)$ for $x$ belonging to a dense subset of the state space $X$.
To obtain the ISS property for all states and inputs, one can use density arguments introduced in Section~\ref{sec:Density_argument}.
\end{remark}

If the input, with respect to which the Lie derivative $\dot{V}_u(x)$ is computed, is clear from the context,  we write simply $\dot{V}(x)$. 

\begin{definition}
\label{def:Equilibrium}
\index{equilibrium}
\index{equilibrium point}
Consider a system $\Sigma=(X,\Uc,\phi)$.
We call $0 \in X$ an \emph{equilibrium point} (of the undisturbed system)  if
$\phi(t,0,0) = 0$ for all $t \geq 0$.
\end{definition}

Now we use the generalized comparison principle (Proposition~\ref{prop:ComparisonPrinciple}) to show the direct Lyapunov theorem.

\index{Direct Lyapunov theorem}
\index{theorem!direct Lyapunov}
\begin{theorem}[Direct Lyapunov theorem]
\label{LyapunovTheorem}
Let $\Sigma = (X,\Uc,\phi)$ be a time-invariant control system satisfying the BIC property.

If $\Sigma$ possesses a (L)ISS Lyapunov function, then it is (L)ISS.
\end{theorem}

\begin{proof}
Let the control system $\Sigma=(X,\Uc,\phi)$ possess a LISS Lyapunov function and $\psi_1,\psi_2,\chi,\alpha$, $\rho_1, \rho_2$ be as in Definition \ref{DefLISS_LF}.

Define $R_1:=\min\{\psi_2^{-1} \circ \psi_1(\rho_1), \psi_1(\rho_1)\}$ and $R_2:=\min\{\rho_2,\chi^{-1}\circ\psi_2^{-1}(R_1)\}$.
Now take an arbitrary control $u \in \Uc$ with $\|u\|_{\Uc} \leq R_2$ and define
\[
I_u:=\{x \in D: V(x) \leq \psi_2 \circ \chi(\|u\|_{\Uc})\}.
\]
Note that as $\|u\|_{\Uc} \leq R_2$, for any $x \in I_u$, it holds that
\[
\psi_1(\|x\|_X)\leq V(x)\leq \psi_2 \circ \chi(\|u\|_{\Uc})\leq \psi_2 \circ \chi(R_2) \leq R_1,
\] 
thus 
$\|x\|_X \leq \psi_1^{-1}(R_1) \leq \rho_1$.
This shows that $I_u \subset B_{\rho_1} \subset D$.

\textbf{Step 1:} We are going to prove that $I_u$ is forward invariant in the following sense:
\begin{eqnarray}
x \in I_u \ \ \wedge \ \  v \in\Uc \ \ \wedge \ \ \|v\|_\Uc \leq \|u\|_\Uc  \qrq \phi(t,x,v) \in I_u \quad \forall t\geq 0.
\label{eq:Forward-Invariance-ISS-Direct-Theorem}
\end{eqnarray}
Assume that this implication does not hold for certain $x\in I_u$ and $v \in\Uc$ satisfying $\|v\|_\Uc \leq \|u\|_\Uc$. There are two alternatives.

First, it may happen that the flow map $\phi(\cdot,x,v)$, corresponding to the initial condition $x$ and input $v$, exists on a finite time interval $[0,t_m(x,v))$ and $\phi(t,x,v) \in I_u$ for $t \in[0,t_m(x,v))$. However, this cannot be true
in view of the BIC property of $\Sigma$.

The second alternative is that there exists $t_1 \in (0,t_m(x,u))$ and $\varepsilon>0$ so that 
$V(\phi(t_1,x,v))> \psi_2\circ\chi(\|u\|_\Uc)  + \varepsilon$.
Define 
\[
t^{*}:=\inf\left\{ t\in [0,t_m(x,u)): V(\phi(t,x,v))> \psi_2\circ\chi(\|u\|_\Uc)  + \varepsilon \right\}.
\]
As $x \in I_u$, and $\phi(\cdot,x,v)$ is continuous, 
\[
\psi_2(\|\phi(t^*,x,v)\|_X) \geq V\left(\phi(t^*,x,v)\right) = \psi_2\circ\chi(\|u\|_\Uc)+\varepsilon \geq \psi_2\circ\chi(\|v\|_\Uc).
\]
Due to \eqref{GainImplikation}, it holds that 
\[
D^+V\left(\phi(t^*,x,v)\right)\leq-\alpha\big(V(\phi(t^*,x,v))\big).
\]
Denote $y(\cdot):=V(\phi(\cdot,x,v))$. Assume that there is $(t_k) \subset (t^*,+\infty)$, such that $t_k \to t^*$, $k\to\infty$ and $y(t_k)\geq y(t^*)$.
Then 
\[
0 \leq \Limsup_{k\to\infty}\frac{y(t_k)-y(t^*)}{t_k-t^*} \leq D^+y(t^*) \leq -\alpha(y(t^*))<0,
\]
a contradiction. Thus, there is some $\delta>0$ such that $y(t) < y(t^*) = \psi_2\circ\chi(\|u\|_\Uc)+\varepsilon$ for all $t\in[t^*,t^*+\delta)$.
This contradicts the definition of $t^*$.
The statement \eqref{eq:Forward-Invariance-ISS-Direct-Theorem} is shown.

Overall, for any $x\in I_u$ we have
\[
\psi_1\left(\|\phi(t,x,u)\|_X\right) \leq V(\phi(t,x,u)) \leq \psi_2\circ \chi(\|u\|_\Uc)\quad \forall t\geq 0,
\]
which implies that for all $t\geq 0$
\begin{eqnarray}
\|\phi(t,x,u)\|_X \leq \psi^{-1}_1 \circ \psi_2\circ \chi(\|u\|_\Uc).
\label{eq:Large_t_estimate}
\end{eqnarray}

\textbf{Step 2:} 
Now take arbitrary $x$: $\|x\|_X \leq R_1$ such that there is $\tau>0$ with $\phi(t,x,u) \notin I_u$ for all $t\in[0,\tau]$. By the axiom of shift invariance (item (ii) of Definition~\ref{Steurungssystem}), and using sandwich bounds \eqref{LyapFunk_1Eig}, we have that 
\[
 \psi_2(\|\phi(t,x,u)\|_X) \geq V(\|\phi(t,x,u)\|_X) \geq \psi_2\circ\chi(\|u(\cdot)\|_{\Uc}) \geq \psi_2\circ \chi(\|u(\cdot+t)\|_{\Uc}),
\]
and using that $\psi_2\in\Kinf$, we have that $\|\phi(t,x,u)\|_X \geq \chi(\|u(\cdot+t)\|_{\Uc})$ for $t\in[0,\tau]$.

As long as $\|\phi(t,x,u)\|_X \leq \rho_1$, using the cocycle property, we have the following differential inequality:
\begin{eqnarray}
D^+V(\phi(t,x,u)) = \dot{V}_{u(t+\cdot)}(\phi(t,x,u)) \leq  - \alpha(V(\phi(t,x,u))).
\label{eq:Direct-LISS-Lyapunov-theorem-decay-large-x}
\end{eqnarray}
Hence, for all such $t$, it holds that 
$V(\phi(t,x,u))\leq V(x)$, which implies that 
\[
\|\phi(t,x,u)\|_X\leq \psi_1^{-1} \circ \psi_2 (\|x\|_X) \leq \psi_1^{-1} \circ \psi_2 (R_1) \leq \rho_1.
\]
Hence, \eqref{eq:Direct-LISS-Lyapunov-theorem-decay-large-x} is valid for all $t\in[0,\tau]$.

Note that the map $t\mapsto V(\phi(t,x,u))$ is continuous, and we can use the generalized comparison principle (Proposition~\ref{prop:ComparisonPrinciple}) 
to infer that $\exists\ \tilde{\beta} \in \KL:\ V(\phi(t,x,u)) \leq \tilde{\beta}(V(x_0),t)$, and consequently:
\begin{equation}
\label{BetaSchaetzung}
\|\phi(t,x,u)\|_X \leq \beta(\|x\|_X,t), \quad t\in[0,\tau],
\end{equation}
where $\beta(r,t)=\psi^{-1}_1 \circ \tilde{\beta}(\psi_2^{-1}(r),t)$, $\forall r,t \geq 0$.

The estimate \eqref{BetaSchaetzung} together with the BIC property ensure that the solution can be prolonged onto an interval that is larger than $[0,\tau]$, and if on this interval the trajectory is still outside of $I_u$, then the estimate 
\eqref{BetaSchaetzung} holds at this larger interval. This indicates that the solution exists at least till the time when it intersects $I_u$.

From the properties of $\KL$ functions it follows, that $\exists t_1$:
\[
t_1:= \inf\{t \geq 0: \phi(t,x,u) \in I_u\}.
\]
By above reasoning, \eqref{BetaSchaetzung} holds on $[0,t_1]$.

For $t \geq t_1$, we have by the cocycle property that 
\[
\phi(t,x,u) = \phi\big(t-t_1,\phi(t_1,x,u),u(t_1+\cdot)\big).
\]
As $\|u(t_1+\cdot)\|_\Uc\leq \|u\|_\Uc$, and $\phi(t_1,x,u) \in I_u$, the property
\eqref{eq:Forward-Invariance-ISS-Direct-Theorem}
ensures that $\phi(t,x,u) \in I_u$ for all $t\geq t_1$, and thus for $t\geq t_1$ the estimate \eqref{eq:Large_t_estimate} is valid. 
We conclude from \eqref{eq:Large_t_estimate} and \eqref{BetaSchaetzung} that
\[
\|\phi(t,x,u)\|_X\leq \beta\left(\|x\|_X,t\right)+\gamma(\|u\|_\Uc) \quad \forall t\ge 0.
\]
Our estimates hold for an arbitrary control $u$: $\|u\|_{\Uc} \leq R_2$. Thus, combining \eqref{BetaSchaetzung} and \eqref{eq:Large_t_estimate}, we obtain the LISS estimate for $\Sigma$ in $(t,x,u) \in \R_+\tm B_{R_1}\tm B_{R_2,\Uc}$.

To prove that the existence of an ISS Lyapunov function ensures the ISS of $\Sigma$, one has to argue as above but with $\rho_1=\rho_2=\infty$.
\end{proof}

\subsection{Exponential ISS Lyapunov functions}

In Definition~\ref{Def:eISS}, we have introduced the concept of exponential ISS, which is stronger than ISS as it presumes exponential convergence rates of a system. To guarantee such a strong type of convergence, eISS Lyapunov functions with special properties are needed.

Implication~\eqref{GainImplikation} with $\alpha(r) = ar$, $r>0$, shows that an eISS Lyapunov function decays exponentially along the trajectory, as long as the state is large enough in comparison to the input.
However, as the functions $\psi_1,\psi_2$ in \eqref{LyapFunk_1Eig} are general nonlinear functions, the trajectory itself does not necessarily possess an exponential convergence rate. Nevertheless, we have the following result:
\begin{theorem}
\label{thm:eISS_LF_with_good_psi_implies_eISS}
Let $\Sigma=(X,\Uc,\phi)$ be a control system and let $V$ be an exponential ISS Lyapunov function with $\psi_1(r) = k_1r^p$ and $\psi_2(r) = k_2 r^p$ for certain $k_1,k_2,p >0$. Then $\Sigma$ is exponentially ISS.
\end{theorem}

\begin{proof}
As $V$ is an exponential ISS Lyapunov function, one can derive using argumentation similar to that in the
proof of Theorem~\ref{LyapunovTheorem}, that $V$ satisfies for a certain $\chi \in\Kinf$ and $a>0$ the following estimate:
\begin{eqnarray*}
V(\phi(t,x,u)) \leq  e^{-at}V(x) + \chi(\|u\|_\Uc),\quad x\in X,\ u\in\Uc,\ t\geq 0.
\end{eqnarray*}
According to the assumptions of the theorem, $V$ satisfies for a certain $p\in \N$ the estimates
\begin{eqnarray*}
k_1 \|x\|_X^p \leq V(x) \leq k_2 \|x\|_X^p,\quad x\in X,
\end{eqnarray*}
and thus
\begin{eqnarray*}
k_1 \|\phi(t,x,u)\|_X^p \leq  e^{-at}k_2 \|x\|_X^p + \chi(\|u\|_\Uc),\quad x\in X,\ u\in\Uc,\ t\geq 0.
\end{eqnarray*}
Hence, for all $x\in X,\ u\in\Uc,\ t\geq 0$ we have
\begin{eqnarray*}
\|\phi(t,x,u)\|_X &\leq& \Big( e^{-at}\frac{k_2}{k_1} \|x\|_X^p + \chi(\|u\|_\Uc)\Big)^{1/p}\\
								  &\leq& \Big(2e^{-at}\frac{k_2}{k_1}\Big)^{1/p} \|x\|_X +\big(2\chi(\|u\|_\Uc)\big)^{1/p},
\end{eqnarray*}
which shows \eqref{eq:eISS_sum} and thus implies exponential ISS of $\Sigma$.
\end{proof}

We formulate the density argument for exponentially ISS systems as a corollary of Lemma~\ref{lem:Density_Arg_ISS}.
\begin{lemma}
\label{lem:Density_Arg_eISS}
Consider a control system $\Sigma=(X,\Uc,\phi)$ and let $\phi$ depend continuously on inputs and on initial states.

Let $\hat{X}$, $\hat{\Uc}$ be dense normed vector subspaces of $X$ and $\Uc$, respectively, and let $\hat{\Sigma}:=(\hat{X},\hat{\Uc},\phi)$ be the system, generated by the transition map $\phi$ of $\Sigma$ but restricted to the state space $\hat{X}$ and space of admissible inputs $\hat{\Uc}$.

Let $\hat{\Sigma}$ be eISS. Then $\Sigma$ is also eISS with the same $M,a$ and $\gamma$
in the estimate \eqref{eq:eISS_sum}.
\end{lemma}

\begin{proof}
The result follows from Lemma~\ref{lem:Density_Arg_ISS}, as it claims that the functions $\beta,\gamma$ from the ISS estimate are inherited by $\hat{\Sigma}$.
\end{proof}

\subsection{Lyapunov functions in a dissipative form}
\label{sec:ISS_LFs_in_a_dissipative_form}

ISS Lyapunov functions in dissipative form are an alternative to ISS Lyapunov functions in implication form.
\begin{definition}
\label{def:LISS-LF-dissipative}
Let $D\subset X$ be open with $0 \in D$.
A continuous function $V:D \to \R_+$ is called a \emph{LISS Lyapunov
  function in a dissipative form} for a system $\Sigma = (X,\Uc,\phi)$,  if there exist $r >0$,
$\psi_1,\psi_2 \in \Kinf$, $\alpha \in \Kinf$ and $\sigma \in \K$
such that $\clo{B_r} \subset D$, and for all $x \in \clo{B_r}$
\begin{equation}
\label{eq:LyapFunk_1Eig_LISS-dissipative}
\psi_1(\|x\|_X) \leq V(x) \leq \psi_2(\|x\|_X),
\end{equation}
and the Dini derivative of $V$ along the trajectories of $\Sigma$ satisfies
\begin{equation}
\label{DissipationIneq}
\dot{V}_u(x) \leq -\alpha(\|x\|_X) + \sigma(\|u\|_{\Uc})
\end{equation}
for all $x \in \clo{B_r}$ and $u\in \clo{B_{r,\Uc}} $.

A function $V:X\to\R_+$ is called an \emph{ISS Lyapunov function in a dissipative form}, if \eqref{eq:LyapFunk_1Eig_LISS-dissipative} and
\eqref{DissipationIneq} hold for all $x\in X, u \in \Uc$.
\end{definition}

We have a simple result.
\begin{proposition}
\label{prop:Dissipative_ISS_LF_implies_implicative_ISS_LF}
Let $\Sigma:=(X,\Uc,\phi)$  be a control system satisfying the BIC property.
If $V$ is an (L)ISS Lyapunov function in a dissipative form, then $V$ is an (L)ISS Lyapunov function in an implication form. 
In particular, $\Sigma$ is (L)ISS.
\end{proposition}

\begin{proof}
Define $\chi(r):=\alpha^{-1}(2\sigma(r))$, $r>0$.
Then $\|x\|_X\geq \chi(\|u\|_\Uc) \Iff \sigma(\|u\|_\Uc) \leq \frac{1}{2}\alpha(\|x\|_X)$, and it holds that
\begin{equation*}
 \|x\|_X \geq \chi(\|u\|_{\Uc}) \  \Rightarrow  \ \dot{V}_u(x) \leq -\frac{1}{2}\alpha(\|x\|_X),
\end{equation*}
and hence $V$ is an ISS Lyapunov function in an implication form. Invocation of Theorem~\ref{LyapunovTheorem} shows that $\Sigma$ is ISS.
\end{proof}

\begin{remark}
\label{rem:Implicative_implies_Dissipative?}
Much harder is the question of whether the existence of an
ISS Lyapunov function in an implication form implies the existence of an ISS Lyapunov function in a dissipative form. 
For the special case of ODE systems, this result has been shown under reasonable regularity assumptions on $V$, and the right-hand side in
\cite[Remark 2.4]{SoW95}. In the infinite-dimensional setting, this result was shown under quite technical assumptions in
\cite[Theorem 3.4]{MiI16}. A fully satisfactory result for infinite-dimensional systems is still missing.
\end{remark}

\begin{remark}
\label{rem:Why_many_Definitions}
ISS Lyapunov functions in an implication form are particularly useful for formulating the small-gain theorems.
At the same time, dissipative ISS Lyapunov functions are particularly useful for studying a so-called integral ISS; see Chapter~\ref{chap:Integral input-to-state stability}. 
The dissipative form of a Lyapunov function resembles the storage function in the theory of dissipative systems \cite{Wil72}, \cite{Wil72b}.
\end{remark}

\subsection{ISS Lyapunov functions and a type of input space}
\label{sec:ISS Lyapunov functions and a type of input space}

For any $u\in\Uc$ and any $\tau\geq 0$ define
\[
u_\tau(s):= \ccat{u}{0}{\tau}(s) = 
\begin{cases}
u(s) &, \text{ if } s\in[0,\tau], \\ 
0 &, \text{ if } s> \tau.
\end{cases}
\]

Let us restate the ISS Lyapunov function concept.
\begin{lemma}
\label{lem:noncoercive_ISS_LF-inf-restatement}
Let $\Sigma:=(X,\Uc,\phi)$  be a control system.
Assume that for any $u\in\Uc$ it holds that $\inf_{\tau > 0}\|u_\tau\|_\Uc \leq \|u\|_\Uc$.

Then a continuous function $V:X \to \R_+$ is an \emph{ISS Lyapunov function in a dissipative form} for the system $\Sigma = (X,\Uc,\phi)$,  if and only if there exist $\psi_1,\psi_2,\alpha \in \Kinf$ and $\sigma \in \K$ such that 
\eqref{eq:LyapFunk_1Eig_LISS-dissipative} holds for all $x \in X$
and the Dini derivative of $V$ along the trajectories of $\Sigma$ for all $x \in X$ and $u\in \Uc$ satisfies
\begin{equation}
\label{DissipationIneq_nc-inf}
\dot{V}_u(x) \leq -\alpha(\|x\|_X)  + \sigma(\inf_{\tau > 0}\|u_\tau\|_{\Uc}).
\end{equation}

\end{lemma}

\begin{proof}
\q{$\Leftarrow$}. Follows from the assumption that $\inf_{\tau > 0}\|u_\tau\|_\Uc \leq \|u\|_\Uc$ for all $u\in\Uc$.

\q{$\Rightarrow$}. Let a continuous function $V:X \to \R_+$ be a non-coercive ISS Lyapunov function.
Pick any $x \in X$ and any $u\in\Uc$.  
Since $\Uc$ is a linear space, $0\in\Uc$ and the axiom of concatenation implies that $u_\tau \in \Uc$ for any $\tau > 0$.
By definition of $\dot{V}_u$, for any $\tau>0$ it holds that $\dot{V}_u(x) = \dot{V}_{u_\tau}(x)$, and thus
\begin{equation*}
\dot{V}_u(x) = \dot{V}_{u_\tau}(x) \leq -\alpha(\|x\|_X)  + \sigma(\|u_\tau\|_{\Uc}).
\end{equation*}
Taking the infimum over $\tau>0$, we obtain \eqref{DissipationIneq_nc-inf}.
\end{proof}

Our definition of ISS Lyapunov functions is defined for any normed vector space $\Uc$, which allows the development of the ISS Lyapunov theory for a very broad class of systems. However, for some input spaces, this definition is far too restrictive, and for other systems, simpler restatements of the ISS Lyapunov function concept may be more useful.
Next, we consider the most important input spaces.

\begin{remark}[Lyapunov functions for systems with $L^p$-inputs]
\label{rem:Lyapunov-functions-for-Lp-spaces} 
For the space \linebreak $\Uc=L^p(\R_+,U)$, for a Banach space $U$ and some $p\in[1,+\infty)$, Definition~\ref{def:LISS-LF-dissipative} 
 is far too restrictive.
Indeed, for any $u\in\Uc = L^p(\R_+,U)$ it holds that 
$\inf_{\tau\ge 0}\|u_\tau\|_{\Uc} = 0$, and thus the inequality \eqref{DissipationIneq_nc-inf}  reduces to
\begin{equation}
\label{DissipationIneq_nc-Lp}
\dot{V}_u(x) \leq -\alpha(\|x\|_X),
\end{equation}
which ensures not only ISS but also ISS with a zero gain (with $\gamma\equiv 0$), also called the uniform global asymptotic stability (see \cite{MiW19a} for more on this notion). 
Coercive and non-coercive ISS Lyapunov theory, which is appropriate for systems with $\Uc=L^p(\R_+,U)$, $p\in[1,+\infty)$, has been developed in \cite{Mir20}.
\end{remark}

\begin{remark}[Piecewise continuous input functions]
\label{rem:Lyapunov-functions-PC-inputs}
If $\Uc=PC_b(\R_+,U)$, for a Banach space $U$, then the dissipation inequality \eqref{DissipationIneq_nc-inf} simplifies to 
\begin{equation}
\label{DissipationIneq_nc-PC_b}
\dot{V}_u(x) \leq -\alpha(\|x\|_X)  + \sigma(\|u(0)\|_{U}),
\end{equation}
which resembles the classical dissipation inequality used in the ISS theory of ODE systems. Similar simplifications occur for ISS Lyapunov functions in implication form.
\end{remark}

\section{Examples}

In this section, we analyze the ISS properties of several basic partial differential equations with distributed as well as boundary inputs. 
To concentrate on the study of ISS, we tacitly assume that all the systems we consider are well-posed with the chosen state space and input space in the sense of mild solutions. 
In each case, we construct an ISS Lyapunov function for the system.

\subsection{Transport equation with a boundary input}
\label{sec:Transport equation with a boundary input}

\index{equation!transport}
Consider the \emph{transport equation}
\begin{eqnarray}
x_t(z,t) + x_z(z,t) =0, \quad (z,t) \in (0,1)\tm(0,+\infty)
\label{eq:TransportEquation}
\end{eqnarray}
subject to the boundary condition
\begin{eqnarray}
x(0,t)=u(t),\quad t>0,
\label{eq:TransportEquation_Dirichlet_BC}
\end{eqnarray}
and the initial condition
\begin{eqnarray}
x(z,0)=x_0(z),\quad z\in(0,1).
\label{eq:TransportEquation_Dirichlet_Init_Cond}
\end{eqnarray}

\begin{proposition}
\label{prop:Transport_Equation_eISS}
Let $X:=L^2(0,1)$ and $\Uc:=PC_b(\R_+,\R)$. Then the system \eqref{eq:TransportEquation}, \eqref{eq:TransportEquation_Dirichlet_BC}, \eqref{eq:TransportEquation_Dirichlet_Init_Cond} is exponentially ISS.
\end{proposition}

\begin{proof}
Pick any $\mu>0$. Choose the following ISS Lyapunov function candidate:
\[
V(x):=\int_0^1 e^{-\mu z}x^2(z)dz.
\]
Clearly, $V$ is coercive as 
\begin{eqnarray}
e^{-\mu} \|x\|^2_X \leq V(x) \leq \|x\|^2_X, \quad x\in X.
\label{eq:Coercive_Estimate_Transport_LF}
\end{eqnarray}
Let us compute the Lie derivative for $V$ along \eqref{eq:TransportEquation} with boundary conditions \eqref{eq:TransportEquation_Dirichlet_BC}. We want to show this first for classical solutions, as this simplifies the computation of the Lie derivative of $V$ considerably.

Assume that $x\in C^1(0,1)$ and $u\in C^1(\R)$. For these initial conditions and boundary inputs, the corresponding classical solution of  \eqref{eq:TransportEquation}, \eqref{eq:TransportEquation_Dirichlet_BC} exists, and we have
\begin{eqnarray*}
\dot{V}_u(x) &=& -\int_0^1e^{-\mu z} 2x(z)x_z(z)dz\\
&=& -e^{-\mu z}x^2(z)\Big|_{z=0}^{z=1} -\mu  \int_0^1e^{-\mu z} x^2(z)dz\\
&=& -e^{-\mu }x^2(1) + x^2(0)  -\mu  \int_0^1e^{-\mu z} x^2(z)dz\\
&\leq& x^2(0)  -\mu  \int_0^1e^{-\mu z} x^2(z)dz\\
&\leq& \|u\|_\Uc^2  -\mu  V(x).
\end{eqnarray*}
This shows that $V$ is an ISS Lyapunov function for \eqref{eq:TransportEquation} with boundary inputs \eqref{eq:TransportEquation_Dirichlet_BC} on the space $C^1(0,1)$, for which the above computations are valid.

As $C^1(0,1)$ is dense in $L^2(0,1)$, we may use the density argument (Lemma~\ref{lem:Density_Arg_eISS}) to show that the system \eqref{eq:TransportEquation}, \eqref{eq:TransportEquation_Dirichlet_BC}
is ISS over the whole state space $L^2(0,1)$.

Now Theorem~\ref{thm:eISS_LF_with_good_psi_implies_eISS} implies that \eqref{eq:TransportEquation}, \eqref{eq:TransportEquation_Dirichlet_BC}, \eqref{eq:TransportEquation_Dirichlet_Init_Cond} is exponentially ISS.
\end{proof}

\begin{remark}
\label{rem:More_General_1d_Systems}
ISS of a much more general class of hyperbolic first-order systems over a 1-dimensional spatial domain is studied through a Lyapunov method in \cite{TPT16}.
\end{remark}

\subsection{Viscous Burgers' equation with Dirichlet boundary conditions and a distributed input}
\label{sec:Burgers_distributed_input}

\index{equation!Burgers'}
Consider a \emph{viscous Burgers' equation} with a distributed input $u$:
\begin{eqnarray}
x_t(z,t) = x_{zz}(z,t) - a x(z,t)x_z(z,t) + b x(z,t) + u(z,t),\quad z\in(0,1),\ t>0,
\label{eq:Viscous_Burgers_distributed_input}
\end{eqnarray}
with $a,b \in\R$, subject to Dirichlet boundary conditions:
\begin{eqnarray}
x(0,t)=x(1,t)=0,\quad t>0.
\label{eq:DirichletConditions_viscous_Burgers_distributed_input}
\end{eqnarray}

For $a=0$, the Burgers' equation \eqref{eq:Viscous_Burgers_distributed_input} reduces to the linear heat equation with Dirichlet boundary conditions.

Burgers' equation is known first of all as a simplification of the Navier-Stokes equations. However, it also has important applications to aerodynamics \cite{Col51}. Furthermore, the inviscid Burgers' equation is the simplest conservation-law system in which the shock-wave phenomenon occurs \cite{Eva10}. In the theory of stochastic PDEs, the Burgers' equation with the addition of stochastic forcing plays an important role and is known as Kardar-Parisi-Zhang equation \cite{KPZ86}. An extensive treatment of Burgers' equation and its generalizations can be found in \cite{Sac87}.

\begin{proposition}
\label{prop:Burgers_Distributed_Input_eISS}
Let $X:=L^2(0,1)$ and $\Uc:= PC_b(\R_+,X)$. Then the system \eqref{eq:Viscous_Burgers_distributed_input}, \eqref{eq:DirichletConditions_viscous_Burgers_distributed_input} is exponentially ISS provided $b<\pi^2$.	
\end{proposition}

\begin{proof}
We are going to show that \eqref{eq:Viscous_Burgers_distributed_input} is ISS with these $X$ and $\Uc$.

Consider the following candidate ISS Lyapunov function:
\begin{eqnarray}
V(x):=\int_0^1 x^2(z) dz = \|x\|^2_{X},
\label{eq:LF_L2_norm_Burgers_distr_input}
\end{eqnarray}
which is evidently coercive.

We choose $x \in C^1(0,1)$ and
$u\in\Uc$ with $v:=u(\cdot,0) \in C(0,1)$.  Then the classical solution $\phi(\cdot,x,u)$ of
\eqref{eq:Viscous_Burgers_distributed_input}, \eqref{eq:DirichletConditions_viscous_Burgers_distributed_input}
exists, and we can differentiate $V$ as follows:

\begin{subequations}
\label{eq:TmpISS_LF_estimate}
\begin{eqnarray}
\dot{V}_u(x) &=& 2 \int_0^1 x(z) \big(x_{zz}(z) - a x(z)x_z(z) + b x(z) + v(z)\big) dz \nonumber\\
&=& 2 \int_0^1 x(z) x_{zz}(z) dz - 2a \int_0^1 x^2(z)x_z(z)dz  \label{eq:TmpISS_LF_estimate_a} \\
&&\qquad\qquad\qquad\qquad\qquad+ 2b \int_0^1 x^2(z) dz +  2\int_0^1 x(z)v(z) dz.
\label{eq:TmpISS_LF_estimate_b}
\end{eqnarray}
\end{subequations}
Partial integration of the first term of \eqref{eq:TmpISS_LF_estimate_a} together with the Dirichlet boundary conditions
\eqref{eq:DirichletConditions_viscous_Burgers_distributed_input} leads to
\begin{eqnarray}
\label{eq:Burgers_key_estimate}
2 \int_0^1 x(z) x_{zz}(z) dz = 2x(z)x_z(z)\big|^{z=1}_{z=0}-2\int_0^1 x^2_z(z)dz = -2\int_0^1 x^2_z(z)dz.
\end{eqnarray}

The only technique which we have used till now is the integration by parts.
The critical point in our argument is now to upperestimate the last term of \eqref{eq:Burgers_key_estimate}
in terms of the Lyapunov function $V(x)$. This can be done utilizing Friedrichs' inequality \eqref{ineq:Friedrichs}:
\begin{eqnarray*}
2 \int_0^1 x(z) x_{zz}(z) dz &\leq & -2 \pi^2  \int_0^1 x^2(z)dz.
\end{eqnarray*}

The second term of \eqref{eq:TmpISS_LF_estimate_a} vanishes in view of Dirichlet boundary conditions:
\begin{eqnarray*}
\int_0^1 x^2(z)x_z(z) dz = \frac{1}{3} x^3(z)\big|^{z=1}_{z=0} = 0.
\end{eqnarray*}
Finally, from \eqref{eq:TmpISS_LF_estimate}, we obtain 
\begin{eqnarray*}
\dot{V}_u(x) &\leq& -2\pi^2 V(x) + 2b V(x)  + 2\int_0^1 x(z)v(z) dz.
\end{eqnarray*}
Utilizing Young's inequality \eqref{thm:Young_Simple}, we have for any $\varepsilon>0$ 
\begin{eqnarray*}
\dot{V}_u(x) &\leq& -2\pi^2 V(x) + 2b V(x)  + \eps V(x) + \frac{1}{\eps}   \int_0^1 v^2(z) dz\nonumber\\
&=& (2b+\eps -2\pi^2) V(x) +  \frac{1}{\eps} \int_0^1 v^2(z) dz.
\label{eq:TmpISS_LF_estimate_final}
\end{eqnarray*}
The last inequality shows that $V$ is an ISS Lyapunov function for the viscous Burgers' equation \eqref{eq:Viscous_Burgers_distributed_input} in the $L^2$ norm provided
\begin{eqnarray}
b<\pi^2.
\label{eq:Burgers_stability_margin}
\end{eqnarray}
Theorem~\ref{thm:eISS_LF_with_good_psi_implies_eISS} now shows that
\eqref{eq:Viscous_Burgers_distributed_input}, \eqref{eq:DirichletConditions_viscous_Burgers_distributed_input}
is eISS with $x \in C^1(0,1)$ and $v \in C(0,1)$.

The density argument (Lemma~\ref{lem:Density_Arg_eISS}) shows that
\eqref{eq:Viscous_Burgers_distributed_input}, \eqref{eq:DirichletConditions_viscous_Burgers_distributed_input}
is ISS on the whole state space $X$ and for all inputs from $\Uc$.
\end{proof}

\begin{remark}
\label{rem:Importance_linear_FuncAn}
Our success in ISS analysis crucially relied on the application of Friedrichs' inequality to the term
\eqref{eq:Burgers_key_estimate}. Using this inequality we have not only proved that
\eqref{eq:Viscous_Burgers_distributed_input}, \eqref{eq:DirichletConditions_viscous_Burgers_distributed_input}
is ISS for some values of parameter $b$, but also obtained a tight estimate \eqref{eq:Burgers_stability_margin} for the values of $b$.
In fact, Friedrichs' inequality is best suited for the nonlinear parabolic initial boundary value problem
\eqref{eq:Viscous_Burgers_distributed_input}, \eqref{eq:DirichletConditions_viscous_Burgers_distributed_input}
since the operator of the highest order in \eqref{eq:Viscous_Burgers_distributed_input}, taken together with the boundary conditions
\eqref{eq:DirichletConditions_viscous_Burgers_distributed_input} is Dirichlet Laplacian,
and Friedrichs' inequality can be obtained by means of a careful spectral analysis of the Dirichlet Laplacian.

This indicates that a precise analysis of ISS for particular PDEs relies heavily on integral inequalities specially tuned for this type of equations. As we will see in the ISS analysis of the Kuramoto-Sivashinskiy equation, this is indeed the case.
\emph{The importance of spectral analysis of differential operators underlines the role of linear functional analysis in the stability analysis of semilinear PDEs.}
\end{remark}

\begin{remark}
\label{rem:Strategy_for_ISS_Analysis}
The strategy that we used to perform ISS analysis in this example will be a sample for the study of more complex equations:
\begin{itemize}
	\item[(i)] Choose an ISS Lyapunov function candidate, say $V$.
	\item[(ii)] Using suitable integral inequalities, show that $V$ is an ISS Lyapunov function on a dense subset of state and input spaces.
	\item[(iii)] Use direct Lyapunov theorem to show ISS on these subspaces of $X$ and $\Uc$.
	\item[(iv)] Use density arguments to show ISS on the whole state and input space.
\end{itemize}
\end{remark}

\subsection{Kuramoto-Sivashinskiy equation with distributed inputs}
\label{sec:KS_equation_with_distributed_input}

\index{equation!Kuramoto-Sivashinskiy}
Consider a \emph{Kuramoto-Sivashinskiy (KS) equation} with distributed inputs:
\begin{eqnarray}
x_t(z,t) = -x_{zzzz}(z,t) - \lambda x_{zz}(z,t) - bx(z,t)x_z(z,t) + u(z,t),
\label{eq:KS_Equation_Dirichlet_BC}
\end{eqnarray}
where $\lambda \in\R$ is an antidiffusion parameter, $b\in\R$ and $u$ is a distributed external input.

Consider \eqref{eq:KS_Equation_Dirichlet_BC} with homogeneous Dirichlet boundary conditions:
\begin{eqnarray}
x(0,t)=x_z(0,t)=x(1,t)=x_z(1,t)=0,\quad t>0.
\label{eq:KS-eq_Dirichlet_Conditions_BI_zero}
\end{eqnarray}

Equation~\eqref{eq:KS_Equation_Dirichlet_BC} was proposed independently by Kuramoto and Tsuzuki \cite{KuT75, Kur78}
as a model for phase turbulence in reaction-diffusion systems and by Sivashinsky \cite{Siv77} as a model for plane flame propagation, describing the influence of both diffusion and thermal conduction of the gas on the stability of a plane flame front.
Furthermore, the KS equation is useful for modeling cellular instabilities \cite{NST85}.

\begin{proposition}
\label{prop:KS_EQ_distr_input_ISS_in_L2_norm}
The system \eqref{eq:KS_Equation_Dirichlet_BC}, \eqref{eq:KS-eq_Dirichlet_Conditions_BI_zero} is exponentially ISS with
$X:=L^2(0,1)$ and $\Uc:= PC_b(\R_+,X)$, $u\in\Uc$ provided that
\begin{eqnarray}
\lambda < 4\pi^2.
\label{eq:KS_eq_lambda_estimate}
\end{eqnarray}
\end{proposition}

\begin{proof}
We apply the strategy described in Remark~\ref{rem:Strategy_for_ISS_Analysis}. 
Consider a candidate ISS Lyapunov function as
\begin{eqnarray}
V(x):=\int_0^1 x^2(z) dz.
\label{eq:KS_eq_LF_L2_norm_Dirichlet_BC}
\end{eqnarray}

For smooth $x$ and $u \in\Uc$ with $u(\cdot,0)=v\in C(0,1)$, we obtain that
\begin{subequations}
\begin{eqnarray}
\dot{V}_u(x) &=& \int_0^1 2x(z) \Big(-x_{zzzz}(z) - \lambda x_{zz}(z) - b x(z)x_z(z) + v(z)\Big) dz \nonumber\\
&=& -2\int_0^1 x(z) x_{zzzz}(z) dz - 2\lambda \int_0^1 x(z)x_{zz}(z)dz  \label{eq:KSeq_TmpISS_LF_estimate_Dirichlet_BC_Dist_Inp_1}\\
&&\qquad\qquad- 2b\int_0^1 x^2(z)x_z(z) dz +  2\int_0^1 x(z)v(z) dz. \label{eq:KSeq_TmpISS_LF_estimate_Dirichlet_BC_Dist_Inp2}
\end{eqnarray}
\label{eq:KSeq_TmpISS_LF_estimate_Dirichlet_BC_Dist_Inp}
\end{subequations}
Partial integration of the first term of \eqref{eq:KSeq_TmpISS_LF_estimate_Dirichlet_BC_Dist_Inp_1} together with the Dirichlet boundary conditions leads to
\begin{eqnarray}
-2\int_0^1 x(z) x_{zzzz}(z) dz &=& -2x(z)x_{zzz}(z)\big|^{z=1}_{z=0}+2\int_0^1 x_z(z)x_{zzz}(z) dz \nonumber\\
&=& -2x(z)x_{zzz}(z)\big|^{z=1}_{z=0} + 2x_z(z)x_{zz}(z)\big|^{z=1}_{z=0} -  2\int_0^1 x^2_{zz}(z) dz\nonumber\\
\big[\text{due to } \eqref{eq:KS-eq_Dirichlet_Conditions_BI_zero}\big] &=& -  2\int_0^1 x^2_{zz}(z) dz.
\label{eq:KSeq_1st_term_Dist_Inp_Dirichlet_BC}
\end{eqnarray}
Partial integration of the second term of \eqref{eq:KSeq_TmpISS_LF_estimate_Dirichlet_BC_Dist_Inp_1} together with the Dirichlet boundary conditions leads to
\begin{eqnarray}
-2\lambda \int_0^1 x(z)x_{zz}(z)dz = 2\lambda\int_0^1 x^2_{z}(z)dz.
\label{eq:KSeq_2nd_term_Dist_Inp_Dirichlet_BC}
\end{eqnarray}
The first term of \eqref{eq:KSeq_TmpISS_LF_estimate_Dirichlet_BC_Dist_Inp2} vanishes in view of Dirichlet boundary conditions:
\begin{eqnarray}
\int_0^1 x^2(z)x_z(z) dz = \frac{1}{3} x^3(z)\big|^{z=1}_{z=0} = 0.
\label{eq:KSeq_3rd_term_Dist_Inp_Dirichlet_BC}
\end{eqnarray}

Combining \eqref{eq:KSeq_1st_term_Dist_Inp_Dirichlet_BC}, \eqref{eq:KSeq_2nd_term_Dist_Inp_Dirichlet_BC}, \eqref{eq:KSeq_3rd_term_Dist_Inp_Dirichlet_BC} we obtain from \eqref{eq:KSeq_TmpISS_LF_estimate_Dirichlet_BC_Dist_Inp} that
\begin{eqnarray*}
\dot{V}_u(x) &\leq& -2 \Big(\int_0^1 x^2_{zz}(z)dz - \lambda\int_0^1 x^2_{z}(z)dz\Big) +  2\int_0^1 x(z)v(z) dz.
\end{eqnarray*}
Now we apply the inequality \eqref{eq:KS_LiK01_inequality} that was specially derived for the Dirichlet problem for the KS equation to the first term on the right-hand side of the above inequality. At the same time, we apply the Cauchy-Bunyakovsky-Schwarz inequality to the second term. For each $\varepsilon>0$, it holds that
\begin{eqnarray}
\dot{V}_u(x) &\leq& -2 \sigma(\lambda)\int_0^1 x^2(z)dz +  \eps \int_0^1 x^2(z)dz + \frac{1}{\eps} \int_0^1 v^2(z) dz \nonumber\\
&=& (\varepsilon-2 \sigma(\lambda)) V(x)  + \frac{1}{\eps} \int_0^1 v^2(z) dz.
\label{eq:Dissipative_Estimate_KS_Equation}
\end{eqnarray}
Here $\sigma$ is defined by \eqref{eq:min-sigman-KS-equation}. 
Let $\lambda<4\pi^2$. According to Lemma~\ref{lem:LowerBound_of_the_Spectrum}, $\sigma(\lambda)>0$ and choosing $\varepsilon$ small enough, we obtain for \eqref{eq:Dissipative_Estimate_KS_Equation} that $V$ is an exponential ISS Lyapunov function in a dissipative form.

Theorem~\ref{thm:eISS_LF_with_good_psi_implies_eISS} implies ISS of \eqref{eq:KS_Equation_Dirichlet_BC}, \eqref{eq:KS-eq_Dirichlet_Conditions_BI_zero} for 
 smooth $x\in L^2(0,1)$ and continuous $u$. The density argument finishes the proof of the proposition.
\end{proof}

One can show that $\lambda=4\pi^2$ is a precise bound for ISS of 
\eqref{eq:KS_Equation_Dirichlet_BC}, \eqref{eq:KS-eq_Dirichlet_Conditions_BI_zero}. For $\lambda>4\pi^2$ the undisturbed linearization of \eqref{eq:KS_Equation_Dirichlet_BC}, \eqref{eq:KS-eq_Dirichlet_Conditions_BI_zero} is unstable, see \cite{LiK01} for more details.

\subsection {Nonlinear reaction-diffusion equation in a Sobolev space}
\label{sec:Reaction-Diffusion_equation_Sobolev_space}

\index{equation!reaction-diffusion}
Consider the following \emph{reaction-diffusion system}:
\begin{subequations}
\begin{align}
x_t(z,t) & = x_{zz}(z,t) - f(x(z,t)) + u(z,t), \quad z \in (0,\pi),\ t>0, \label{eq:Reaction-Diffusion_Dist_Inp_H10_norm-1}\\
x(0,t)   & = x(\pi,t)=0.
\end{align}
\label{eq:Reaction-Diffusion_Dist_Inp_H10_norm}
\end{subequations}
We assume that $f:\R \to\R$ is Lipschitz continuous on bounded sets, monotonically increasing up to infinity, $f(-r)=-f(r)$ for all $r \in \R$ (in particular, $f(0)=0$).

In previous examples, our canonical choice of a state space was $L^2$-space. This choice is reasonable also for \eqref{eq:Reaction-Diffusion_Dist_Inp_H10_norm}. However, if the initial states which we are interested in have more regularity, e.g., belong to a Sobolev space
$H^1_0(0,\pi)$, it is reasonable to analyze ISS w.r.t.\ $H^1_0(0,\pi)$ norm of a solution.
In this section, we show that  ISS of \eqref{eq:Reaction-Diffusion_Dist_Inp_H10_norm} w.r.t.\ state space
$X:=H^1_0(0,\pi)$ and input space $\Uc:=PC_b(\R_+,L^2(0,\pi))$ can be analyzed by means of Lyapunov methods as well.
In what follows we consider the norm in $H^1_0(0,\pi)$ as
\[
\|x\|_{H^1_0(0,\pi)}:=  \Big(\int_0^\pi x_z^2(z) dz\Big)^{\frac{1}{2}},
\]
which is equivalent to the \q{native} $H^1(0,\pi)$-norm due to Friedrichs' inequality.

\begin{proposition}
\label{prop:Reaction-Diffusion_Dist_Inp_H10_norm}
The system \eqref{eq:Reaction-Diffusion_Dist_Inp_H10_norm} is exponentially ISS with
$X:=H^1_0(0,\pi)$ and $\Uc:= PC_b(\R_+,L^2(0,1))$.
\end{proposition}

%
%
%

\begin{proof}
Consider the following ISS Lyapunov function candidate:
\begin{equation}
\label{LF}
V(x):= \int_0^\pi{ \Big( \frac{1}{2} x_z^2(z)+ \int_0^{x(z)}{f(y)dy} \Big) dz},\quad x\in X.
\end{equation}

According to the embedding theorem for Sobolev spaces (see Theorem~\ref{thm5}), every $x\in X = H^1_0(0,\pi)$ belongs actually to 
$C^{0,\frac{1}{2}}(0,\pi)$ (H\"older space with H\"older exponent $\frac{1}{2}$, see Section~\ref{sec:Hoelder spaces} for the definition). Moreover, there exists a constant $c$, which does not depend on $x\in H^1_0(0,\pi)$, such that
\begin{equation}
\label{EmbeddingIneq}
\|x\|_{C^{0,\frac{1}{2}}(0,\pi)} \leq c \|x\|_{H^1_0(0,\pi)},\quad x \in H^1_0(0,\pi).
\end{equation}
In particular, this tells us that if we identify any $x \in X$ with its continuous version, 
then the expression \eqref{LF} makes sense (as $x(z)$ is well-defined).

We are going to prove that $V$ is an ISS Lyapunov function.
To verify the estimates \eqref{LyapFunk_1Eig} for a function $V$, note that $\int_0^{r}{f(y)dy} \geq 0$ for every $r \in \R$ and thus
\begin{equation*}
V(x) \geq \int_0^\pi{\frac{1}{2} x_z^2(z) dz} = \frac{1}{2}  \|x\|^2_{H^1_0(0,\pi)},\quad x\in X.
\end{equation*}

Let us find an estimate from above. We have
\begin{equation*}
V(x)= \int_0^\pi{ \frac{1}{2} x_z^2(z)}\ dz +
\int_0^\pi{ \int_0^{x(z)}{f(y)dy} \ dz},\quad x\in X.
\end{equation*}

Define $\psi:\R_+ \to \R_+$ by 
\[
\psi(r):=\frac{1}{2} r^2 + \sup_{x:\ \|x\|_{H^1_0(0,\pi)} \leq r} \int_0^\pi \int_0^{x(z)}{f(y)dy} dz,\quad r\ge 0.
\]
Inequality \eqref{EmbeddingIneq} and the fact that $\|x\|_{C(0,\pi)} \leq \|x\|_{C^{0,\frac{1}{2}}(0,\pi)}$ for all $x \in C^{0,\frac{1}{2}}(0,\pi)$ imply
\begin{eqnarray*}
\psi(r)&=&\frac{1}{2} r^2 + \sup_{x:\ c\|x\|_{H^1_0(0,\pi)} \leq cr} \int_0^\pi \int_0^{x(z)}{f(y)dy} dz \\
      &\leq& \frac{1}{2} r^2 + \sup_{x:\ \|x\|_{C(0,\pi)} \leq cr} \int_0^\pi \int_0^{x(z)}{f(y)dy} dz \leq \psi_2(r),
\end{eqnarray*}
where $\psi_2(r) := \frac{1}{2} r^2 + \pi \int_0^{cr}{f(y)dy}$, $r\ge 0$.
Since $f$, restricted to positive values of the argument, belongs to $\Kinf$,  $\psi_2$ is also in $\Kinf$.

Finally, for all $x \in H^1_0(0,\pi)$, we have:
\begin{eqnarray}
\label{1Eig_Beispiel}
\frac{1}{2}  \|x\|^2_{H^1_0(0,\pi)} \leq V(x) \leq \psi_2(\|x\|_{H^1_0(0,\pi)}),
\end{eqnarray}
and the property \eqref{LyapFunk_1Eig} is verified.

Let us compute the Lie derivative of $V$ for smooth $x$ and for $u\in \Uc$ with $v:=u(\cdot,0) \in C(0,\pi)$:
\begin{eqnarray*}
\dot{V}_u(x) &=& \int_0^\pi{ x_z(z) x_{zt}(z)+ f\big(x(z)\big) x_t(z) dz} \\
           &=& \left[x_z x_t \right]_{z=0}^{z=\pi} +\int_0^\pi{  - x_{zz}(z)x_t(z)+ f\big(x(z)\big)x_t(z)dz}.
\end{eqnarray*}

From the boundary conditions, it follows that 
$\left[x_z x_t \right]_{z=0}^{z=\pi} = 0$. As $x$ solves \eqref{eq:Reaction-Diffusion_Dist_Inp_H10_norm-1}, we obtain
\[
\dot{V}_u(x)= - \int_0^\pi{ \Big(x_{zz}(z)- f\big(x(z)\big) \Big)^2 dz} + \int_0^\pi{ \Big(x_{zz}(z) - f\big(x(z)\big) \Big)(-v(z)) dz}.
\]
Define
\[
I(x):=\int_0^\pi{ \Big(x_{zz}(z)- f\big(x(z)\big) \Big)^2 dz}.
\]
Using the Cauchy-Schwarz inequality for the second term, we have:
\begin{equation}
\label{LeiAb_V_IS}
\dot{V}_u(x) \leq - I(x) +  \sqrt{I(x)} \ \|v\|_{L^2(0,\pi)}.
\end{equation}
Let us consider $I(x)$:
\begin{eqnarray*}
I(x) &=&\int_0^\pi x_{zz}^2(z) dz -2\int_0^\pi x_{zz}(z) f(x(z)) dz + \int_0^\pi f^2(x(z)) dz \\
     &=& \int_0^\pi x_{zz}^2(z) dz +2\int_0^\pi x_z^2(z) \frac{\partial f}{\partial x}(x(z)) dz
   + \int_0^\pi{f^2(x(z)) dz}  \\
	  &\ge& \int_0^\pi x_{zz}^2(z) dz.
\end{eqnarray*}
For $x \in H^1_0(0,\pi) \cap H^2(0,\pi)$, a version of Friedrichs' inequality \cite[p. 85]{Hen81} implies that
\[
\int_0^\pi{ x_{zz}^2(z) dz} \ge \int_0^\pi{ x_z^2(z) dz}.
\]
Overall, we have:
\begin{equation}
\label{IS_Abschaetzung}
I(x) \geq \|x\|^2_{H^1_0(0,\pi)}.
\end{equation}

Choose the Lyapunov gain as
\[
\chi(r)=a r,\; a>1.
\]
If $\chi(\|v\|_{L^2(0,\pi)}) \le \chi(\|u\|_{\Uc}) \le \|x\|_{H^1_0(0,\pi)}$, we obtain from \eqref{LeiAb_V_IS}, using \eqref{IS_Abschaetzung}:

\begin{equation}
\label{LyapGainAbsch}
\dot{V}_u(x) \leq - I(x) +  \frac{1}{a} \sqrt{I(x)} \|x\|_{H^1_0(0,\pi)} \leq \Big(\frac{1}{a} -  1\Big)I(x) \leq \Big(\frac{1}{a} -  1\Big) \|x\|^2_{H^1_0(0,\pi)}.
\end{equation}
This shows ISS for smooth states and inputs. The density argument finishes the proof.
%
%
%
\end{proof}



\subsection{Lyapunov methods for semilinear parabolic systems with boundary inputs}
\label{sec:Lyapunov methods for semilinear parabolic systems with boundary inputs}

This section shows how the Lyapunov method can be applied to analyze semilinear parabolic systems with \emph{boundary inputs of Neumann type}. 
Consider the Ginzburg-Landau equation:
\begin{eqnarray}
x_t(z,t) = \mu x_{zz}(z,t) + a x(z,t) - x^3(z,t),\quad t>0,\ z\in(0,1),
\label{eq:Semilin-parabolic-Neumann-input}
\end{eqnarray}
where $\mu>0$ is the diffusion coefficient. We investigate \eqref{eq:Semilin-parabolic-Neumann-input} subject to the Neumann boundary input at $z=0$ and to Dirichlet boundary condition at $z=1$:
\begin{subequations}
\label{eq:Neumann-input-parabolic-PDE}
\begin{eqnarray}
x_z(0,t) &=& u(t) ,\quad t>0, \label{eq:Neumann-input-parabolic-PDE-1}\\
x(1,t)&=& 0, \quad t>0. \label{eq:Neumann-input-parabolic-PDE-2}
\end{eqnarray}
\end{subequations}

For consistency with the original paper, we assume that $u \in \Uc:=C^2(\R_+,\R)$, and we understand the solutions in a classical sense, 
but the norm in the state space will be chosen as $L^2(0,1)$-norm
and the norm in the input space $\Uc$ will be chosen as $L^\infty(\R_+)$-norm.
The use of classical solutions calls for adaptation of the ISS concept, and hence we require for the ISS property the validity of the estimate \eqref{iss_sum} only for smooth enough $x$.
The precise definition of the solution concept, state space as well as the proof of well-posedness of the PDE model \eqref{eq:Semilin-parabolic-Neumann-input}, \eqref{eq:Neumann-input-parabolic-PDE}
can be found in \cite{ZhZ18}.

To derive conditions for ISS of 
\eqref{eq:Semilin-parabolic-Neumann-input}, \eqref{eq:Neumann-input-parabolic-PDE}, we use the Lyapunov function candidate
\begin{eqnarray}
V(x):=\|x\|^2_{L^2(0,1)} = \int^1_0 x^2(z) dz.
\label{eq:L2_LF_Linear_diffusion}
\end{eqnarray}

First let us compute the Lie derivative of $V$ for $x$ and $u$ smooth enough:
\begin{eqnarray}
\label{eq:Partial-integration}
\dot{V}_u(x) &=&  2 \int^1_0 x(z) \Big( \mu x_{zz}(z) + a x(z) - x^3(z)\Big) dz \nonumber\\
&=&  2\mu \Big(x_z(z)x(z)\Big)\Big|_{z=0}^{z=1} - 2\mu \int^1_0 \big(x_z(z)\big)^2 dz \nonumber\\
&&\qquad + 2aV(x) - 2 \int^1_0 \big(x^2(z)\big)^2dz.
\end{eqnarray}

Using \eqref{eq:Neumann-input-parabolic-PDE} as well as Jensen's inequality \eqref{ineq:Jensen} for the last term, we obtain 
\begin{eqnarray*}
\dot{V}_u(x)
 &\leq & 2\mu |u(0)| |x(0)|  - 2\mu \int^1_0 \big( x_z(z)\big)^2 dz + 2aV(x) - 2 V^2(x).
\end{eqnarray*}
Using Cauchy's inequality \eqref{ineq:Young_Simple} and afterwards Agmon's inequality \eqref{ineq:Agmon} with boundary condition \eqref{eq:Neumann-input-parabolic-PDE-2},  we obtain for any $\varepsilon>0$:
\begin{eqnarray*}
\hspace{10mm} 2\mu |u(0)| |x(0)| \leq \varepsilon |x(0)|^2 + \frac{\mu^2}{\varepsilon} |u(0)|^2
 \leq \varepsilon \Big(\|x\|^2_{L^2(0,1)} + \|x_z\|^2_{L^2(0,1)}\Big) 
+ \frac{\mu^2}{\varepsilon} |u(0)|^2. 
\end{eqnarray*}
Using this estimate in \eqref{eq:Partial-integration} and rearranging the terms, we obtain
\begin{eqnarray*}
\dot{V}_u(x) 
&\leq& (\varepsilon - 2\mu )\int^1_0 x_z^2(z) dz + (\varepsilon + 2a)V(x) - 2 V^2(x) + \frac{\mu^2}{\varepsilon} |u(0)|^2. 
\end{eqnarray*}
Assuming that $\varepsilon < 2\mu$, we can use Poincare's inequality \eqref{Wirtinger_Variation_Ineq}
 for the first term to obtain that
\begin{eqnarray}
\label{eq:Partial-integration-3}
\dot{V}_u(x) 
&\leq& (\varepsilon - 2\mu )\frac{\pi^2}{4}\int^1_0 x^2(z) dz + (\varepsilon + 2a)V(x) - 2 V^2(x) + \frac{\mu^2}{\varepsilon} |u(0)|^2 \nonumber \\
&=& \Big((\varepsilon - 2\mu )\frac{\pi^2}{4} + \varepsilon + 2a\Big)V(x) - 2 V^2(x) + \frac{\mu^2}{\varepsilon} |u(0)|^2.
\end{eqnarray}
To ensure that the dissipation inequality holds, we have to assume that \linebreak
 $(\varepsilon - 2\mu )\frac{\pi^2}{4} + \varepsilon + 2a \leq 0$. As $\varepsilon>0$ can be chosen arbitrarily small, we obtain the following sufficient condition for ISS of the system 
\eqref{eq:Semilin-parabolic-Neumann-input}, \eqref{eq:Neumann-input-parabolic-PDE}:
\begin{eqnarray}
a< \frac{\mu\pi^2}{4}.
\label{eq:Obtained-criterion}
\end{eqnarray}
The term $- 2 V^2(x)$ in \eqref{eq:Partial-integration-3} shows that outside of the neighborhood of the equilibrium, the convergence rate of the system \eqref{eq:Semilin-parabolic-Neumann-input}, \eqref{eq:Neumann-input-parabolic-PDE}
is faster than exponential.

\subsubsection{Tightness of obtained estimates}
\label{sec:Tightness-of-results} 

Linearizing the system \eqref{eq:Semilin-parabolic-Neumann-input}, \eqref{eq:Neumann-input-parabolic-PDE} for $u\equiv 0$ near the equilibrium, one can see that, for
$a$ and $\mu$ as in \eqref{eq:Obtained-criterion}, the linearized system is exponentially stable, and for 
$a$ and $\mu$ satisfying
\begin{eqnarray}
a> \frac{\mu\pi^2}{4}
\label{eq:Obtained-criterion-unstable}
\end{eqnarray}
the system \eqref{eq:Semilin-parabolic-Neumann-input}, \eqref{eq:Neumann-input-parabolic-PDE} is unstable (even locally near the equilibrium).
This can be done, e.g.,\ by analyzing the spectral properties of the corresponding Laplacian operator, similarly to \cite[Section 1.3]{Hen81}.

We summarize the obtained results in the following proposition:
\begin{proposition}
\label{prop:Semilinear-Parabolic-system-with-Dirichlet-boundary-input} 
Consider the system \eqref{eq:Semilin-parabolic-Neumann-input}, \eqref{eq:Neumann-input-parabolic-PDE}
with $X:=L^2(0,1)$ and $\Uc:=C^2(\R_+,\R)$ (with the supremum norm).

If the condition \eqref{eq:Obtained-criterion} holds, then \eqref{eq:Semilin-parabolic-Neumann-input}, \eqref{eq:Neumann-input-parabolic-PDE} is ISS.

If the condition \eqref{eq:Obtained-criterion-unstable} holds, then 
\eqref{eq:Semilin-parabolic-Neumann-input}-\eqref{eq:Neumann-input-parabolic-PDE} is not uniformly locally asymptotically stable for $u\equiv 0$ (in the sense of Definition~\ref{Stab_Notions_Undisturbed_Systems}(vii)).
\end{proposition}

Thus our analysis (motivated by \cite{ZhZ18}) provides fairly tight results for ISS of the system \eqref{eq:Semilin-parabolic-Neumann-input}, \eqref{eq:Neumann-input-parabolic-PDE}. This indicates that the Lyapunov method combined with the sharp versions of inequalities for the $L^p$ spaces, as Poincare's, Agmon's inequalities, etc., is an efficient method for analyzing ISS of semilinear systems with boundary inputs and for computation of the precise uniform decay rate of the solutions of a system.

\section{Concluding remarks}

\textbf{Control systems.} Definition~\ref{Steurungssystem} of a control system is frequently used within the system-theoretic community at least since the 1960s \cite{KFA69, Wil72}. 
Not all important systems are covered by this definition. In particular,
the input space $C(\R_+,U)$ of continuous $U$-valued functions does not satisfy the axiom of concatenation.
This should not be a big restriction since already piecewise continuous and $L^p$ inputs, which are used in control theory much more frequently than continuous ones, satisfy the axiom of concatenation. 

In a similar spirit, even more general system classes can be considered, containing output maps, time-variant dynamics, the possibility for a solution to jump at certain time instants, systems, which fail to satisfy the cocycle property, etc., see \cite{Kar07a, KaJ11, KaJ11b}.

\textbf{ISS.} The notions of input-to-state stability and an ISS Lyapunov function have been introduced in the context of ODEs in the seminal paper \cite{Son89}.
ISS has quickly become one of the pillars of nonlinear control theory, including robust stabilization, nonlinear observer design, and analysis of large-scale networks, see \cite{KKK95,ArK01,Son08}. 
The dissipative form of a Lyapunov function resembles the storage function in a theory of dissipative systems \cite{Wil72}, \cite{Wil72b}.

The importance of the ISS concept for nonlinear control theory has led to the
development of related notions, refining and/or generalizing ISS in some
sense: integral ISS \cite{Son98, ASW00}, input-to-output stability (IOS)
\cite{JTP94, KPJ08}, local ISS \cite{DaR10}, input-output-to-state
stability (IOSS) \cite{SoW97,KSW01}, input-to-state dynamical stability
(ISDS) \cite{Gru02}, incremental ISS \cite{InS02,Ang09}, to mention a
few.
On the other hand, significant research efforts have been devoted to the
extension of ISS theory to further classes of systems, such as
discrete-time \cite{Ang99}, \cite{SoK03}, \cite{LJH12}, \cite{JiW01,LaN03}, impulsive \cite{HLT08}, \cite{DKM12}, \cite{ChZ09}, \cite{LLX11}, \cite{YSH13}, \cite{DaM13b} switched
\cite{VCL07}, and hybrid \cite{CTG07, CTG08}, \cite{CaT09}, \cite{DaK13}, \cite{LNT14}, \cite{MYL18} systems. See \cite[Chapter 7]{Mir23} for an overview of the results on ISS for different system classes.

For a comprehensive treatment of the ISS theory for finite-dimensional systems, we refer the reader to the recent monograph \cite{Mir23} and to the surveys \cite{Son08,DES11}.

Tremendous progress in the development of the infinite-dimensional ISS theory has brought a number of powerful techniques for the robust stability analysis of distributed parameter systems, including Lyapunov methods \cite{DaM13,ZhZ18}, characterizations of ISS \cite{MiW18b} for nonlinear systems, functional-analytic tools \cite{JNP18,JSZ19,Sch20}, spectral-based methods \cite{KaK16b,KaK19}, monotonicity-based tools \cite{MKK19,ZhZ20b} for linear systems. Most of these tools will be studied in this work. 
We refer to the recent survey~\cite{MiP20} for a detailed literature review.

Proposition~\ref{prop:Converging_input_uniformly_converging_state} is folklore.
Definition~\ref{DefLISS_LF} and Theorem~\ref{LyapunovTheorem} are due to \cite{DaM13}.

\textbf{Examples.} ISS analysis of Burgers' equation with distributed inputs in Section~\ref{sec:Burgers_distributed_input} was never published, though it uses a standard technique. The case of Burgers' equation with a boundary control of Dirichlet or Neumann type is much harder, and there is no fully satisfactory analysis of ISS of such a system, though see \cite{ZhZ19b} for analysis of local ISS of Burgers' equation with boundary inputs using De-Giorgi technique.

Proposition~\ref{prop:KS_EQ_distr_input_ISS_in_L2_norm} is an ISS counterpart of \cite[Theorem 2.1(i)]{LiK01}. 
Stability and boundary control of Kuramoto-Sivashinskiy equation \eqref{eq:KS_Equation_Dirichlet_BC}, \eqref{eq:KS-eq_Dirichlet_Conditions_BI_zero}
has been studied in several papers, see, e.g., \cite{LiK01}. 
Several recent papers on the boundary and (finite-dimensional) in-domain control of Kuramoto-Sivashinskiy equation include \cite{CoL15, KaF18b}.

Example in Section~\ref{sec:Reaction-Diffusion_equation_Sobolev_space} is due to \cite{DaM13}.

The example in Section~\ref{sec:Lyapunov methods for semilinear parabolic systems with boundary inputs} follows the method for analysis of nonlinear reaction-diffusion equations reported in \cite{ZhZ18}. 
The results in \cite{ZhZ18} have been extended in \cite{Sch20} to some classes of abstract semilinear boundary control systems. 
We present the method from \cite{ZhZ18} on a representative example, which shows the essence of the technique.



\cleardoublepage
\chapter[Evolution equations: Well-posedness]{Evolution equations: Well-posedness and properties of the flow}
\label{chap:Nonlinear Evolution Equations}

In Chapter~\ref{chap:Intro}, we have introduced a general class of control systems and considered the method of Lyapunov functions for the analysis of the input-to-state stability for such general systems. 
In this analysis, we have tacitly assumed that our control systems are well-posed. At the same time, usually, well-posedness is not given in advance. The models of the natural and engineering systems are given in terms of the equations of motion, and a crucial step in the analysis of such systems is to prove their well-posedness.
In this chapter, we consider a rather general class of semilinear evolution equations in Banach spaces. This class is broad enough to include important classes of linear and nonlinear evolution equations with nonlinear internal dynamics and distributed and boundary inputs. 
For such systems, we derive sufficient conditions for well-posedness, as well as basic properties of the flow map.

We start with the well-posedness analysis of linear evolution equations with unbounded input operators. This allows characterizing the ISS of linear systems as a combination of well-posedness and exponential stability.
The ISS analysis of such systems is important, particularly for analyzing the robustness of PDE systems w.r.t.\  boundary disturbances and for robust stabilization by means of boundary controllers.

Next, we proceed to the analysis of semilinear evolution equations with unbounded input operators and Lipschitz continuous nonlinearities. The structure of semilinear evolution equations allows combining the \q{linear} methods of admissibility theory with \q{nonlinear} methods, such as fixed point theorems and Lyapunov methods. We consider the case of general $C_0$-semigroups and the special case of analytic semigroups, for which one can achieve stronger results. We discuss semilinear parabolic systems with Dirichlet and Neumann boundary inputs as an example of our developments.

\section{Linear systems}
\label{sec:Linear_systems}

\subsection{Abstract linear control systems}
\label{sec:Abstract linear control systems}

We start by introducing the framework of abstract linear control systems, where we follow the approach in \cite{Wei89b} (see also \cite{TuW09}) rather closely.

Assume that $U$ and $X$ are Banach spaces, $p\in [1,\infty]$ and $\Uc:=L^p(\R_+,U)$, which is a space of strongly measurable functions $u:\R_+\to U$, such that $\|u\|_{L^p(\R_+,U)}:=\Big(\int_{\R_+}\|u(s)\|^p_U ds\Big)^{\frac{1}{p}} < \infty$.

\begin{definition}
\label{def:abstract linear control system}
\index{control system!linear abstract}
An \emph{abstract linear control system} on $X$ and $\Uc$ is a pair $(T,\Phi)$, where $T:= \sg{T} \subset L(X)$ is a strongly continuous semigroup and $\Phi:=\sg{\Phi} \subset L(\Uc,X)$ is so that for all $u,v \in\Uc$ and all $t,\tau\geq 0$ the \emph{composition property} holds:
\begin{eqnarray}
\hspace{8mm}\Phi(t+\tau)(\ccat{u}{v}{\tau}) = T(t)\Phi(\tau)u + \Phi(t)v,\quad
\ccat{u}{v}{\tau}(s):=
\begin{cases}
u(s), & \text{ if } s \in [0,\tau], \\
v(s-\tau),  & \text{ otherwise}.
\end{cases}
\label{eq:Composition-property}
\end{eqnarray}
\end{definition}
We immediately verify some simple properties of the family $\Phi=(\Phi(t))_{t\geq0}$:
\begin{proposition}
\label{prop:Identity-causality-growth}
Let $(T,\Phi)$ be an abstract linear control system. Then:
\begin{enumerate}[label=(\roman*)]
	\item $\Phi(0)=0$.
	\item The \emph{causality} property holds:
\[
\Phi(\tau)u = \Phi(\tau)(\ccat{u}{0}{\tau}) = \Phi(\tau)P(\tau) u,
\]
where we define the projection operator by $P(\tau) u:=\ccat{u}{0}{\tau}$, $u\in\Uc$.
		\item The map $t\mapsto \|\Phi(t)\|$ is nondecreasing.
\end{enumerate}
\end{proposition}

\begin{proof}
(i). Setting $u=v$ and $t=\tau =0$, we see that $\Phi(0)u=2\Phi(0)u$ for all $u$.

(ii). Just take $v=0$ and $t=0$.

(iii). Taking in \eqref{eq:Composition-property} $u=0$ and taking the supremum over all $v$ satisfying $\|v\|_U=1$, we get for any $t,\tau\geq 0$ that
\begin{equation*}
\|\Phi(t)\| \le \|\Phi({t+\tau})\|.
\end{equation*}
\end{proof}

The next result shows the strong continuity of the family $\Phi=(\Phi(t))_{t\geq0}$ for $p<\infty$.
\begin{proposition}
\label{prop:Wei89b-Proposition_2.3}
Let $X$ and $\Uc$ be as in Definition~\ref{def:abstract linear control system} with $p < \infty$ and
let $\Sigma=(T,\Phi)$ be an abstract linear control system
on $X$ and $\Uc$. Then the function
\begin{gather*}
q(t, u) = \Phi(t) u
\end{gather*}
is continuous on the product $\R_+ \times \Uc$. In particular, $\Phi=(\Phi(t))_{t\geq0}$ is a~strongly continuous family of operators.
\end{proposition}

\begin{proof}
Let us show the continuity of $q$ with
respect to $t$. Fix any $u\in\Uc$ and consider the map
\begin{gather*}
f_u(t):=\Phi(t) u,\quad t\geq 0.
\end{gather*}
As $t\mapsto \|\Phi(t)\|$ is nondecreasing, the causality implies for $t\in[0, 1]$ that
\begin{gather*}
\|f_u(t)\|_X = \|\Phi(t) P(t) u\|_X \le \|\Phi(1)\| \cdot \|P(t) u\|_{\Uc}.
\end{gather*}
By the dominated convergence theorem (here we use that $p<\infty$), we see that 
$\|P(t) u\|_{\Uc} \to 0$ as $t \to +0$, hence
\begin{gather*}
\lim_{t\to +0}f_u(t) = 0.
\end{gather*}
The composition property \eqref{eq:Composition-property} and strong continuity of $T$ now imply the right continuity of $f_u$.
To prove the left continuity of $f_u$ at $\tau>0$, consider a~sequence $(\varepsilon_n) \subset [0, \tau]$ such that $\varepsilon_n \to 0$ as $n\to\infty$, and define $u_n(t)=u(\varepsilon_n+t)$. 
Again invoking the dominated convergence theorem for $p<\infty$, we see that $u_n \to u$ as $n\to\infty$.

Clearly,
\begin{gather*}
u = \ccat{u}{u_n}{\varepsilon_n},
\end{gather*}
and using the composition property, we obtain that 
\begin{gather*}
\Phi({\varepsilon_n+(\tau-\varepsilon_n)})u
 = T({\tau-\varepsilon_n})\Phi({\varepsilon_n})u
 + \Phi({\tau-\varepsilon_n}) u_n.
\end{gather*}
This shows
\begin{gather*}
\Phi(\tau)u - \Phi({\tau-\varepsilon_n})u
 = T({\tau-\varepsilon_n})\Phi({\varepsilon_n})u
 + \Phi({\tau-\varepsilon_n})(u_n-u),
\end{gather*}
which yields
\begin{gather*}
\|\Phi(\tau)u - \Phi({\tau-\varepsilon_n})u\|_X
 \le M \cdot \|f_u(\varepsilon_n)\|_X + \|\Phi(\tau)\| \cdot \|u_n-u\|_{\Uc},
\end{gather*}
where $M$ is a~bound for $\|T(t)\|$ on $[0,\tau]$.
Thus, the left continuity of $f_u$ in any $\tau > 0$ is also proved.

Hence, $q$ is continuous in both arguments, and the joint continuity of $q$ follows from the decomposition
\begin{gather*}
\Phi(t) v - \Phi(\tau)u = \Phi(t)(v-u) + (\Phi(t)-\Phi(\tau))u.
\end{gather*}
\end{proof}

The next proposition shows that abstract linear control systems are a subclass of \q{general control systems}: 
\begin{proposition}
\label{prop:Abstract_LCS_are_general_CS}
Let $X$, $\Uc$, $T$ and $\Phi$ be as in Definition~\ref{def:abstract linear control system}. Define
the flow map
\begin{eqnarray}
\phi(t,x,u):=T(t)x + \Phi(t)u,\quad t\geq0,\ x\in X,\ u\in\Uc.
\label{eq:Flow-map}
\end{eqnarray}
Let either $p<\infty$ or $p=\infty$ and $\phi$ be continuous w.r.t.\  time.

Then the triple $\Sigma:=(X,\Uc,\phi)$ is a control system in the sense of Definition~\ref{Steurungssystem}.
\end{proposition}

\begin{proof}
Clearly, the properties (i)--(iv) of Definition~\ref{Steurungssystem} are satisfied for $\Sigma$.
In view of Proposition~\ref{prop:Identity-causality-growth}, $\Phi(0)=0$, which implies the identity axiom for $\Sigma$.

Setting $t=0$ into \eqref{eq:Composition-property}, we obtain that
$\Phi(\tau)(\ccat{u}{v}{\tau}) = \Phi(\tau)u$, for any $u,v\in\Uc$, and thus for all $\tau\geq 0$, $x\in X$ and $u,v\in\Uc$,
it holds that
\[
\phi(\tau,x,\ccat{u}{v}{\tau})=T(\tau)x + \Phi(\tau)(\ccat{u}{v}{\tau}) = T(\tau)x + \Phi(\tau)u =\phi(\tau,x,u).
\]
This shows the causality property of $\Sigma$.

For $p=\infty$, we take the continuity of $\phi$ w.r.t.\  time by assumption, and for $p<\infty$, the continuity axiom follows from Proposition~\ref{prop:Wei89b-Proposition_2.3} and strong continuity of $T$.

Finally, the cocycle property for $\Sigma$ follows by the following computation, valid for all $t,\tau\ge 0$, $x \in X$, and all $u,v\in\Uc$.
\begin{eqnarray*}
\phi(t,\phi(\tau,x,u),v) &=& T(t)\phi(\tau,x,u) + \Phi(t)v \\
&=& T(t)(T(\tau)x + \Phi(\tau)u) + \Phi(t)v \\
&=& T(t+\tau)x + T(t)\Phi(\tau)u + \Phi(t)v \\
&=& T(t+\tau)x + \Phi(t+\tau)(\ccat{u}{v}{\tau}) \\
&=& \phi(t+\tau,x,\ccat{u}{v}{\tau}).
\end{eqnarray*}
\end{proof}

Having a growth rate of the semigroup $T$ of an abstract linear control system $(T,\Phi)$, one can also find the growth rate of $\Phi$.

\begin{proposition}
\label{prop:Proposition_2.5}
Let $X$ and $\Uc$ be as in Definition~\ref{def:abstract linear control system}, and let $(T,\Phi)$ be an abstract linear control system on $X$ and $\Uc$. 
Take $M \ge 1$ and $\omega \in\R$ such that
\begin{gather*}
\|T(t)\| \le M e^{\omega t},\qquad t\geq 0.
\end{gather*}
If $\omega>0$, there is some $L\geq0$ such that
\begin{gather*}
\|\Phi(t)\| \le L e^{\omega t},\qquad t\ge0.
\end{gather*}
If $\omega<0$, then there is $G>0$, such that $\|\Phi(t)\|\leq G$ for all $t\geq 0$.
\end{proposition}

\begin{proof}
From the composition property, we get by induction for any $(u_k)_{k=1}^n\subset \Uc$ that
\begin{equation}
\Phi(n)\left( \dots \left(\ccat{u_1}{u_2}{1}\right)\ccat{}{}{2}
 \dots \ccat{}{}{n-1}u_n\right)
 = T({n-1}) \Phi(1) u_1
 + T({n-2}) \Phi(1) u_2
 + \dots
 + \Phi(1) u_n,
\label{eq:Book_2.4}
\end{equation}
and since $t\mapsto \|\Phi(t)\|$ is nondecreasing, we obtain for $t \in (n-1, n]$:
\begin{gather*}
\|\Phi(t)\| \le \|\Phi(n)\|
 \le \Big(\|T({n-1})\|+\|T({n-2})\|+\dots+\|I\|\Big)
 \cdot \|\Phi(1)\|
 \le M \frac{e^{\omega n}-1}{e^{\omega}-1}\|\Phi(1)\|,
\end{gather*}
so we can take
\begin{gather*}
L = M \frac{e^{\omega}}{e^{\omega}-1}\|\Phi(1)\|.
\end{gather*}
If $\omega<0$, the uniform boundedness of $\sg{\Phi}$ follows easily from \eqref{eq:Book_2.4}.
\end{proof}

The following proposition resolves the question of the eISS of a \emph{well-defined} abstract linear control system to the exponential stability of the semigroup $T$.
\begin{proposition}
\label{prop:ISS_abstract:linear:systems}
Under the assumptions of Proposition~\ref{prop:Abstract_LCS_are_general_CS}, the following statements are equivalent:
\begin{itemize}
    \item[(i)]  $\Sigma$ is eISS with a linear gain.
    \item[(ii)] $T$ is exponentially stable.
\end{itemize}
\end{proposition}

\begin{proof}
%
\noindent(i) $\Rightarrow$ (ii). Clear.


\noindent(ii) $\Rightarrow$ (i). As $T$ is exponentially stable, there are $M\ge 1$, $\lambda>0$ so that $\|T(t)\|\leq M e^{-\lambda t}$ for all $t\geq 0$.
According to Proposition~\ref{prop:Proposition_2.5}, it holds that $\|\Phi(t)\|\leq G$ for a certain $G>0$ and all $t\geq 0$.
Thus, for all $t\geq0$, $x\in X$, $u\in\Uc$, we have
\begin{eqnarray*}
\|\phi(t,x,u)\|_X \leq \|T(t)x\|_X + \|\Phi(t)\| \|u\|_{\Uc}\leq Me^{-\lambda t}\|x\|_X + G \|u\|_{\Uc},
\end{eqnarray*}
which shows eISS of $\Sigma$ with a linear gain $G$.
\end{proof}


\subsection{Representation of abstract control systems}
\label{sec:Representation of abstract control systems}

Assume that $X$ is a Banach space and $A:D(A)\to X$, $D(A)\subset X$ is an infinitesimal generator of a strongly continuous semigroup on $X$ (which already implies that $\rho(A) \neq \emptyset$).
Define the \emph{extrapolation space} $X_{-1}$ as the completion of $X$ with respect to the norm
$ \|x\|_{X_{-1}}:= \|(aI -A)^{-1}x\|_X$ for some $a$ in the resolvent set $\rho(A) \subset\C$ of $A$.

\begin{lemma}
\label{lem:Equivalent-norms-on-X-1} 
The space $X_{-1}$ does not depend on the choice of $a \in\rho(A)$, and different choices of $a \in \rho(A)$ generate equivalent norms on $X_{-1}$.
\end{lemma}

\begin{proof}
An easy exercise.
\end{proof}

Lifting of the state space $X$ to a larger space $X_{-1}$ brings good news: the semigroup  $(T(t))_{t\ge 0}$ extends uniquely to a strongly continuous semigroup  $(T_{-1}(t))_{t\ge 0}$ on $X_{-1}$ whose generator $A_{-1}:X_{-1}\to X_{-1}$ is an extension of $A$ with $D(A_{-1}) = X$, see, e.g.,\ \cite[Section II.5]{EnN00}.

In practice, a system is usually defined by means of differential equations.
The following result shows that every abstract linear control system with $p\neq \infty$ can be represented as a flow map of a certain linear differential equation.

\begin{theorem}
\label{thm:Representation theorem of Weiss}
Let $U$, $X$ be Banach spaces, $p\in[0,+\infty)$ and let $\Uc:=L^p(\R_+,U)$.
Assume that $(T,\Phi)$ is an abstract linear control system on $X$ and $\Uc$.
Then there is a unique operator $B \in L(U, X_{-1})$ such that for any $u \in\Uc$ and any $t\ge 0$ it holds that
\begin{eqnarray}
\Phi(t)u = \int_0^t T_{-1}(t-s) B u(s)ds.
\label{eq:Phi-representation}
\end{eqnarray}
Moreover, for any $x_0\in X$ and $u\in\Uc$ the function $\phi:\R_+\times X \times \Uc \to X$, defined by
\[
\phi(t,x_0,u) = T(t)x_0  + \Phi(t)u
\]
is the (continuous) mild solution of
\[
\dot{x}(t) = A x(t) + B u(t)
\]
with $x(0) = x_0$.
\end{theorem}

According to Theorem~\ref{thm:Representation theorem of Weiss},
every abstract linear control system can be represented as a mild solution of a linear equation. However, not every pair $(A,B)$ where $A$ is a generator of a strongly continuous semigroup and $B \in L(U, X_{-1})$ gives rise to an abstract linear control system.
Theorem~\ref{thm:Representation theorem of Weiss}  tells us that (at least for $p<\infty$) the condition $B\in L(U,X_{-1})$ is a necessary (but not sufficient) condition for that.
Next, we consider the classes of input operators that do generate linear control systems. As we will see, the understanding of such classes of systems is important not only for the well-posedness theory but also for the ISS analysis of linear control systems.
Bounded operators $B$ constitute the simplest class of such operators. They will be discussed in the next section.
Afterward, we turn our attention to classes of so-called admissible operators, which provide us with a much more general and wide-reaching theory.

\subsection{ISS of linear systems with bounded input operators}
\label{sec:Linear_systems_bounded_ops}

Let $X,U$ be Banach spaces. Consider a linear control system
\begin{equation}
\label{eq:Linear_System}
\dot{x}(t)=Ax(t)+ Bu(t), \quad t > 0,
\end{equation}
where $A$ is the generator of a strongly continuous semigroup $T:=(T(t))_{t\ge 0}$ on $X$ and $B\in L(U,X)$.

As we know from \cite[Lemma 10.1.6]{JaZ12}, the (continuous) mild solution of \eqref{eq:Linear_System} corresponding to a certain initial condition $x_0\in X$ and to an external input $u$, is given by 
\begin{eqnarray}
\phi(t,x_0,u):=T(t)x_0 +\int_0^t T(t-s)B u(s)ds,\quad t\ge 0,
\label{eq:Lin_Sys_mild_Solution}
\end{eqnarray}
where the integral is well-defined in the sense of Bochner, see \cite[Chapter 1]{ABH11} for the Bochner integration theory.

Consider for $p\in[1,+\infty]$ the spaces $L^{p}_{\loc}(\R_+,W)$ of \emph{locally Bochner-integrable functions}, consisting of functions $f:\R_+\to W$ so that the restriction $f_I$ of a map $f$ to any compact interval $I \subset \R_+$ belongs to $L^p(I,W)$.

\begin{remark}
\label{rem:L_p-loc-spaces}
Note that the space $L^{1}_{\loc}(\R_+,U)$ is not a normed vector space, as some of its elements have \q{infinite} $L^1$-norm (that is, the $L^1$-norm cannot be introduced for these spaces), and hence according to Definition~\ref{Steurungssystem}, $L^{1}_{\loc}(\R_+,U)$ is not allowed as an input space to a control system. However, ISS of a control system w.r.t.\ inputs in $L^1(\R_+,U)$ implies ISS of the system with $L^{1}_{\loc}(\R_+,U)$, due to causality of the control systems, similarly to Proposition~\ref{prop:Causality_Consequence}.
Therefore, we slightly abuse the use of Definition~\ref{Steurungssystem}, and allow for the input spaces $L^{p}_{\loc}(\R_+,U)$.
\end{remark}

It is reasonable to consider ISS of \eqref{eq:Linear_System} with respect to different choices of input spaces, as defined next:
\begin{definition}
\label{def:Lp-ISS} 
Let $\Uc$ be a normed vector space of admissible inputs. 
We say that the system \eqref{eq:Linear_System} is ISS w.r.t.\  $\Uc$, if $(X,\Uc,\phi)$ with $\phi$ as in \eqref{eq:Lin_Sys_mild_Solution} is a well-posed input-to-state stable control system.

If \eqref{eq:Linear_System} is ISS w.r.t.\  $\Uc=L^p(\R_+,U)$, then we call \eqref{eq:Linear_System} $L^p$-ISS.
\end{definition}

The following lemma shows that \eqref{eq:Linear_System} is a well-posed control system with any $\Uc:=L^{p}_{\loc}(\R_+,U)$, $p\in[1,+\infty]$.
\begin{lemma}
\label{lem:Systems-with-bounded-input-operators-are-well-posed}
Consider the system \eqref{eq:Linear_System}, with $B\in L(U,X)$. Let $\Uc:=L^{p}_{\loc}(\R_+,U)$, $p\in[1,+\infty]$, and let $\Phi$ be defined by
\begin{eqnarray*}
\Phi(t)u = \int_0^t T(t-s) B u(s)ds,\quad t\geq0,\quad u\in\Uc.
\end{eqnarray*}
Then $(T,\Phi)$  is a linear control system in the sense of Definition~\ref{def:abstract linear control system}.

In particular, $(T,\Phi)$ is eISS iff $T$ is exponentially stable. 
\end{lemma}

\begin{proof}
As $B$ is a bounded operator, $\Phi(t)$ is bounded for all $t\geq 0$.
Moreover,
it holds that
\begin{eqnarray*}
\Phi(t+\tau)(\ccat{u}{v}{\tau})
&=& \int_0^{t+\tau} T(t+\tau-s)B (\ccat{u}{v}{\tau})(s) ds \\
&=& \int_0^{\tau} T(t+\tau-s)B u(s) ds + \int_{\tau}^{t+\tau} T(t+\tau-s)B v(s-\tau) ds  \\
&=& T(t)\int_0^{\tau} T(\tau-s)B u(s) ds + \int_{0}^{t} T(t-s)B v(s) ds  \\
&=& T(t)\Phi(\tau)u + \Phi(t)v,
\end{eqnarray*}
which shows the composition property \eqref{eq:Composition-property} and herewith the claim of the lemma.

Application of Proposition~\ref{prop:ISS_abstract:linear:systems} finishes the proof of the lemma.
\end{proof}

\subsection{ISS of linear systems with unbounded input operators}
\label{sec:ISS of linear systems}

Assume that $X$ and $U$ are Banach spaces, $q\in [1,\infty]$, and $\Uc:=L^q(\R_+,U)$.
Let also $A:D(A) \subset X\to X$ be an infinitesimal generator of a strongly continuous semigroup $T:=(T(t))_{t\geq 0}$ on $X$ with a nonempty resolvent set $\rho(A)$.
Consider again the equation \eqref{eq:Linear_System} but now assume that $B$ is an unbounded linear operator from $U$ to $X$.
Unbounded input operators naturally appear in boundary or point control of linear systems, see, e.g., \cite{JaZ12}.

To treat this case, we again consider an extrapolation space $X_{-1}$. 
As we know from the representation Theorem~\ref{thm:Representation theorem of Weiss}, the input operator $B$ must satisfy the condition $B\in L(U,X_{-1})$ in order to give rise to a well-defined control system (at least for $q<\infty$). We assume this for all $q\in[1,\infty]$.

Thus we may consider equation \eqref{eq:Linear_System} on the Banach space $X_{-1}$:
\begin{equation}
\label{eq:Linear_System_extrapolated}
\dot{x}(t)=A_{-1}x(t)+ Bu(t), \quad t > 0,
\end{equation}
and mild solutions of this extrapolated system are given by the variation of constants formula \eqref{eq:Lin_Sys_mild_Solution}, which takes
for every $x \in X$, $u\in L^{1}_{\loc}(\R_+,U)$  and every $t\geq 0$ the form
\begin{eqnarray}
\phi(t,x,u)&:=&T_{-1}(t)x+\int_0^t T_{-1}(t-s)B u(s)ds\nonumber\\
        &=&T(t)x+\int_0^t T_{-1}(t-s)B u(s)ds.
\label{eq:Lifted_Lin_Sys_mild_Solution}
\end{eqnarray}
The last transition is due to the fact that $T_{-1}(t)x = T(t)x$ for all $x\in X$ and $t\geq 0$.

The lifting comes, however, at a price that now the solution $\phi(t,x,u)$ has values in $X_{-1}$.
The formula \eqref{eq:Lifted_Lin_Sys_mild_Solution} defines an $X$-valued function only provided that the value of the integral in \eqref{eq:Lifted_Lin_Sys_mild_Solution} belongs to the state space $X$, despite the fact that what we integrate is in $X_{-1}$.
This motivates the following definition:
\begin{definition}
\label{def:q-admissibility}
\index{operator!$q$-admissible}
The operator $B\in L(U,X_{-1})$ is called a {\em $q$-admissible control operator} for $(T(t))_{t\ge 0}$, where $1\le q\le \infty$, if
there is a $t>0$ so that
\begin{eqnarray}
u\in L^q(\R_+,U) \qrq \int_0^t T_{-1}(t-s)Bu(s)ds\in X.
\label{eq:q-admissibility}
\end{eqnarray}
\end{definition}
Define the operators $\Phi(t): \Uc \to X$ by
\[
\Phi(t)u := \int_0^t T_{-1}(t-s)Bu(s)ds.
\]
Denote also $\Phi:=\{\Phi(t):t\ge0\}$. 

For any $B \in L(U,X_{-1})$, the operators $\Phi(t)$ are well-defined as maps from $\Uc$ to $X_{-1}$ for all $t$. The next result shows that from $q$-admissibility of $B$ it follows that the image of $\Phi(t)$ is in $X$ for all $t \geq 0$ and with $\Uc:=L^q(\R_+,U)$ we have $\Phi(t)\in L(\Uc,X)$ for all $t>0$.

\begin{proposition}
\label{prop:Restatement-admissibility}
Let $X, U$ be Banach spaces and let $q \in [1,\infty]$ be given. Then $B \in L(U,X_{-1})$ is $q$-admissible for $T$ if and only if
for all $t>0$ there is $h_t>0$ so that for all $u \in L^{q}_{loc}(\R_+,U)$ it holds that $\Phi(t)u \in X$ and
\begin{equation}
\label{eq:admissible-operator-norm-estimate}
\left\| \int_0^t T_{-1}(t-s)Bu(s)\,ds\right\|_X \le h_t \|u\|_{L^q([0,t],U)}.
\end{equation}
Furthermore, $(T,\Phi)$ defines an abstract linear control system over $X$ and $\Uc:=L^q(\R_+,U)$.
\end{proposition}

\begin{proof}
Due to $q$-admissibility, take any $t_0$ such that \eqref{eq:q-admissibility} holds with $t:=t_0$.

Take any $\tau\in (0,t_0)$, and any $v \in L^q([0,t_0-\tau],U)$. Define $u:=\ccat{0}{v}{\tau} \in L^q([0,t_0],U)$.
Then $\Phi(t_0)u \in X$ and at the same time
\[
\Phi(t_0)u = \int_{\tau}^{t_0} T_{-1}(t_0-s)Bv(s-\tau)ds = \int_{0}^{t_0-\tau} T_{-1}(t_0 - \tau -s)Bv(s)ds = \Phi(t_0-\tau)v.
\] 
As $\tau \in (0,t_0)$ and $v \in L^q([0,t_0-\tau],U)$ are arbitrary, we see that \eqref{eq:q-admissibility} holds for all $t\in[0,t_0]$.

As in the proof of Lemma~\ref{lem:Systems-with-bounded-input-operators-are-well-posed}, one can show that for all $u_1,u_2 \in \Uc=L^q(\R_+,U)$ and all $t,\tau\geq 0$, it holds that 
\begin{eqnarray}
\Phi(\tau+t)(\ccat{u_1}{u_2}{\tau}) = T_{-1}(t)\Phi(\tau)u_1 + \Phi(t)u_2,
\label{eq:composition property-admissible-operators}
\end{eqnarray}
where the equality is considered in the space $X_{-1}$.
Now take any $t \in [0,t_0]$, and set in \eqref{eq:composition property-admissible-operators} $\tau:=t$. We have that 
\[
\Phi(2t)(\ccat{u_1}{u_2}{t}) = T_{-1}(t)\Phi(t)u_1 + \Phi(t)u_2.
\]
We have shown that $\Phi(t)u_1 \in X$. As $T_{-1}$ coincides with $T$ on $X$, we have that 
$T_{-1}(t)\Phi(t)u_1 = T(t)\Phi(t)u_1 \in X$, and thus $\Phi(2t)(\ccat{u_1}{u_2}{t}) \in X$. 
As $u_1,u_2 \in \Uc$ are arbitrary, we obtain that \eqref{eq:q-admissibility} holds for $t\in[0,2t_0]$.
Iteratively, we see that \eqref{eq:q-admissibility} holds for all $t\geq 0$.

Let us show that $\Phi(t)\in L(\Uc,X)$ for all $t>0$.

As $B\in B(U,X_{-1})$, then for $\mu \in\rho(A)$ it holds that $(\mu I - A)^{-1}B \in L(U,X)$.

Then, as $A$ commutes with the semigroup $T$, and using that closed operators \q{commute} with the Bochner integral
(see \cite[Proposition 1.1.7]{ABH11}), we have that 
\[
\Phi(t)u = \int_0^t T_{-1}(t-s)(\mu I - A)(\mu I - A)^{-1}Bu(s)ds 
= (\mu I - A)\int_0^t T_{-1}(t-s)(\mu I - A)^{-1}Bu(s)ds.
\]
Since $(\mu I - A)^{-1}B$ maps to $X$, and $T_{-1}$ coincides with $T$ on $X$, we have that
\[
\Phi(t)u = (\mu I - A)\int_0^t T(t-s)(\mu I - A)^{-1}Bu(s)ds.
\]
Now 
\[
u \mapsto \int_0^t T(t-s)(\mu I - A)^{-1}Bu(s)ds
\]
is a bounded operator from $\Uc$ to $X$, and thus, $\Phi(t)$ is closed as a product of a closed and a bounded operator.

Since $\Phi(t)$ is defined on the whole $\Uc$, the closed graph theorem ensures that $\Phi(t) \in L(\Uc,X)$.
As $T$ and $T_{-1}$ coincide on $X$, in \eqref{eq:composition property-admissible-operators} we have that 
$T_{-1}(t)\Phi(\tau)u_1 = T(t)\Phi(\tau)u_1$ for all $t,\tau\geq 0$ and all $u_1 \in \Uc$.
Together with \eqref{eq:composition property-admissible-operators}, this shows that $(T,\Phi)$ defines an abstract linear control systems over $X$ and $\Uc$.
\end{proof}

If $B\in L(U,X_{-1})$ is a $q$-admissible operator, then $\phi$ given by \eqref{eq:Lifted_Lin_Sys_mild_Solution} is well-defined on $\R_+\times X\times \Uc$, with $\Uc:=L^q(\R_+,U)$. 
The next natural question is whether $\Sigma:=(X,\Uc,\phi)$ is a control system in the sense of Definition~\ref{Steurungssystem}.
For $q<\infty$, it has an affirmative answer:
\begin{proposition}
\label{prop:q-admissibility-implies-continuity}
Assume that $B\in L(U,X_{-1})$ is a $q$-admissible input operator with $q<\infty$.
Then $(X,\Uc,\phi)$ with $\Uc:=L^q(\R_+,U)$ and $\phi$ given by \eqref{eq:Lifted_Lin_Sys_mild_Solution} is a forward-complete control system in the sense of Definition~\ref{Steurungssystem}, satisfying the CEP and BRS properties.
\end{proposition}

\begin{proof}
By Proposition~\ref{prop:Restatement-admissibility}, $(T,\Phi)$  is an abstract linear control system.
In view of Proposition~\ref{prop:Abstract_LCS_are_general_CS}, $(X,\Uc,\phi)$ is a control system in the sense of Definition~\ref{Steurungssystem}. CEP and BRS properties follow from the fact that $T$ is a strongly continuous semigroup and from the estimate \eqref{eq:admissible-operator-norm-estimate}.
\end{proof}

For $\infty$-admissible control operators, we have the following result:
\begin{proposition}
\label{prop:infty-admissibility-implies-continuity}
Assume that $B\in L(U,X_{-1})$ is an $\infty$-admissible input operator and that $\phi$ is continuous w.r.t.\ time in the norm of $X$.
Then $(X,\Uc,\phi)$ with $\Uc:=L^\infty(\R_+,U)$ and $\phi$ given by \eqref{eq:Lifted_Lin_Sys_mild_Solution} is a forward-complete control system in the sense of Definition~\ref{Steurungssystem}, satisfying the CEP and BRS properties.
\end{proposition}

\begin{proof}
The proof resembles the proof of Proposition~\ref{prop:q-admissibility-implies-continuity}.
\end{proof}

\begin{remark}
\label{rem:Admissibility-and-forward-completeness}
In the language of general control systems, $q$-admissibility means precisely the forward-completeness of control systems for all inputs from $L^q(\R_+,U)$, which is in view of Proposition~\ref{prop:Restatement-admissibility} equivalent to BRS property.
\end{remark}

Although $\phi$ given by \eqref{eq:Lifted_Lin_Sys_mild_Solution} is continuous w.r.t.\ time in the norm of $X_{-1}$ (as a mild solution), this does not imply that $\phi$ is continuous in the $X$-norm.
The question whether $\infty$-admissibility of $B\in L(U,X_{-1})$ implies the continuity of the solution map $\phi$ w.r.t.\ time in the $X$-norm has been open for 30 years (\cite[Problem 2.4]{Wei89b}). Recently, this was positively answered in \cite{JSZ19} for an important subclass of analytic systems on Hilbert spaces:
\begin{theorem}
\label{thm:Contiuity-of-a-map}
Let $A$ generate an exponentially stable analytic semigroup on a Hilbert space $X$, which is similar to a contraction semigroup.
Assume further that $\dim (U) <\infty$.

Then any $B\in L(U,X_{-1})$ is an $\infty$-admissible operator, and the corresponding map $\phi$ is continuous in the norm of $X$ for any $u \in \Uc:=L^\infty(\R_+,U)$. Consequently, $\Sigma:=(X,\Uc,\phi)$ is a control system in the sense of Definition~\ref{Steurungssystem}.
\end{theorem}

\begin{proof}
The result follows from \cite[Theorem 1]{JSZ19}.
\end{proof}

\begin{remark}
\label{rem:Relations-between-admissibility-classes}
Lemma~\ref{prop:Restatement-admissibility} and H\"older's inequality show that $q$-admissibility of $B\in L(U,X_{-1})$ implies $p$-admissibility of $B$ for any $p>q$.
Indeed, from \eqref{eq:admissible-operator-norm-estimate}, we obtain for any $p \in (q,\infty)$:
\begin{eqnarray*}
 \left\| \int_0^t T(s)Bu(s)\,ds\right\|_X &\le& h_t \Big(\int_0^t 1 \cdot |u(s)|^q ds\Big)^{1/q}\\
 &\leq &  h_t \Big(\int_0^t 1^{\frac{p}{p-q}} ds\Big)^{\frac{p-q}{p}} \Big(\int_0^t |u(s)|^p ds\Big)^{1/p}\\
 &\leq &  h_t t^{\frac{p-q}{p}} \Big(\int_0^t |u(s)|^p ds\Big)^{1/p},
\end{eqnarray*}
which is precisely $p$-admissibility of $B$. The case of $p=\infty$ is analogous.

At the same time, in general, $q$-admissibility of $B$ does not imply $p$-admissibility for any $p<q$.
In particular, in \cite[Example 5.2]{JNP18}, an example of an $\infty$-admissible operator was given,
which is not $p$-admissible for any $p<\infty$.
Thus, admissibility gives us a measure for how bad-behaved an unbounded input operator may be, and $\infty$-admissible operators are the \q{worst} type of operators, which still allows for well-posedness of a linear control system \eqref{eq:Linear_System} and at the same time gives
the possibility of obtaining ISS with respect to the supremum norm of the inputs, which is the most classical type of ISS estimates.
Bounded operators are always 1-admissible, and if $X$ is reflexive, then $1$-admissibility of $B$ is equivalent to the fact that $B\in L(U, X)$, see
\cite[Theorem 4.8]{Wei89b}.
\end{remark}

We proceed to stability analysis of linear systems with admissible input operators.
\begin{definition}
\label{def:infinite-time-admissibility}
\index{operator!infinite-time $q$-admissible}
A $q$-admissible operator $B$ is called \emph{infinite-time $q$-admissible operator}, if the constant $h:=h_t$ in Lemma~\ref{prop:Restatement-admissibility} does not depend on $t$.
\end{definition}

\begin{lemma}
\label{lem:GrabowskiLemma}
Let $T$ be an exponentially stable semigroup and let $q \in [1,\infty]$. Then $B$ is $q$-admissible control operator iff it is infinite-time $q$-admissible control operator.
\end{lemma}

\begin{proof}
Let $q < \infty$ (the argument is analogous for $q = \infty$). 
We mimic a proof from \cite[Lemma 1.1]{Gra95}.

Clearly, infinite-time $q$-admissibility implies $q$-admissibility.
Assume that $B$ is $q$-admissible and pick arbitrary $\tau>0$ and $u\in\Uc$.
Consider any $t=k\tau+t'$ for some $k \in \N$ and $t'<\tau$.
{\allowdisplaybreaks
\begin{align*}
 \Big\| \int_0^t & T(s)Bu(s)\,ds\Big\|_X \\
&= \left\| \int_{k\tau}^{k\tau+t'} T(s)Bu(s)\,ds +  \sum_{i=0}^{k-1} \int_{i \tau}^{(i+1) \tau} T(s)Bu(s)\,ds\right\|_X \\
&\leq  \left\| \int_{0}^{t'} T(s+k\tau)Bu(s+k\tau)\,ds \right\|_X +  \left\| \sum_{i=0}^{k-1} \int_{0}^{\tau} T(s+i\tau)Bu(s+i\tau)\,ds\right\|_X \\
&\leq   \|T(k\tau)\| \left\| \int_{0}^{t'} T(s)Bu(s+k\tau)\,ds \right\|_X +  \sum_{i=0}^{k-1} \|T(i\tau)\| \left\|\int_{0}^{\tau} T(s)Bu(s+i\tau)\,ds\right\|_X.
\end{align*}
Using Proposition~\ref{prop:Restatement-admissibility} and assuming without loss of generality that $t\mapsto h_t$ is nondecreasing, we have:
\begin{align*}
 \Big\| \int_0^t & T(s)Bu(s)\,ds\Big\|_X \\
&\leq   Me^{-\lambda k\tau} h_{t'} \Big( \int_0^{t'} \|u(s + k\tau)\|^q_U ds \Big)^{1/q}  +  \sum_{i=0}^{k-1} Me^{-\lambda i\tau} 
h_{\tau} \Big(\int_0^{\tau} \|u(s + i\tau)\|^q_U ds \Big)^{1/q} \\
&\leq   Me^{-\lambda k\tau} h_{\tau} \Big(\int_0^{t} \|u(s)\|^q_U ds \Big)^{1/q} +  \sum_{i=0}^{k-1} Me^{-\lambda i\tau} h_{\tau} \Big(\int_0^{t} \|u(s)\|^q_U ds \Big)^{1/q} \\
&\leq   \sum_{i=0}^{k} Me^{-\lambda i\tau} h_{\tau} \Big(\int_0^{t} \|u(s)\|^q_U ds\Big)^{1/q} \\
&\leq   \sum_{i=0}^{\infty} Me^{-\lambda i\tau} h_{\tau} \Big(\int_0^{t} \|u(s)\|^q_U ds\Big)^{1/q}.
\end{align*}
}
Taking $h:=M h_{\tau} \sum_{i=0}^{\infty} e^{-\lambda i\tau}$, we get the result.
\end{proof}

Now we can characterize ISS of \eqref{eq:Linear_System}:
\begin{theorem}[ISS characterizations for linear systems with unbounded input operators]
\label{thm:ISS-Criterion-lin-sys-with-unbounded-operators}
Let $X$ and $U$ be Banach spaces and $q \in[1,\infty]$. Consider the system \eqref{eq:Linear_System} with $\phi$ defined by \eqref{eq:Lifted_Lin_Sys_mild_Solution}.
The following assertions are equivalent.

\begin{enumerate}[label=(\roman*)]
    \item\label{itm:i-ISS-Crit-Lin-Unbounded} \eqref{eq:Linear_System} is $L^q$-eISS with a linear gain.
    \item\label{itm:ii-ISS-Crit-Lin-Unbounded} \eqref{eq:Linear_System} is $L^q$-ISS.
    \item\label{itm:iii-ISS-Crit-Lin-Unbounded} $T$ is exponentially stable $\ \wedge\ $ $B$ is $q$-admissible 
		$\ \wedge \ $ (in case if $q=\infty$) $\phi$ is continuous w.r.t.\  time in the norm of $X$.
    \item\label{itm:iv-ISS-Crit-Lin-Unbounded} $T$ is exponentially stable  $\ \wedge\ $ $B$ is infinite-time $q$-admissible 
    $\ \wedge \ $ (in case if $q=\infty$) $\phi$ is continuous w.r.t.\  time in the norm of $X$.
\end{enumerate}
\end{theorem}

\begin{proof}
\ref{itm:i-ISS-Crit-Lin-Unbounded} $\Rightarrow$ \ref{itm:ii-ISS-Crit-Lin-Unbounded} $\Rightarrow$ \ref{itm:iii-ISS-Crit-Lin-Unbounded}. Evident.

\ref{itm:iii-ISS-Crit-Lin-Unbounded} $\Rightarrow$ \ref{itm:iv-ISS-Crit-Lin-Unbounded}. Follows by Lemma~\ref{lem:GrabowskiLemma}.

\ref{itm:iv-ISS-Crit-Lin-Unbounded} $\Rightarrow$ \ref{itm:i-ISS-Crit-Lin-Unbounded}. 
Follows from Proposition~\ref{prop:Restatement-admissibility} with $h\equiv const$ and exponential stability of a semigroup $T$.
%
\end{proof}

Similarly in spirit to Theorem~\ref{thm:ISS-criterion-linear-systems-bounded-operators},
Theorem~\ref{thm:ISS-Criterion-lin-sys-with-unbounded-operators} reduces the ISS analysis of linear systems to the     stability analysis of the semigroup and to admissibility analysis of the input operator, which are classical functional-analytic problems with
many powerful tools for its solution, see \cite{JaP04, TuW09}.

\section{Semilinear evolution equations with unbounded input operators}
\label{sec:Semilinear boundary control systems}

Consider infinite-dimensional evolution equations of the form
\begin{subequations}
\label{eq:SEE+admissible} 
\begin{eqnarray}
\dot{x}(t) & = & Ax(t) + B_2f(x(t),u(t)) + Bu(t),\quad t>0,  \label{eq:SEE+admissible-1}\\
x(0)  &=&  x_0, \label{eq:SEE+admissible-2}
\end{eqnarray}
\end{subequations}
where $A: D(A)\subset X \to X$ generates a strongly continuous semigroup $T=\sg{T}$ of boun\-ded linear operators on a Banach space $X$; $U$ is a normed vector space of input values, and $x_0\in X$ is a given initial condition. 
As the input space, we take $\Uc:=L^\infty(\R_+,U)$.

The map $f: X\tm U \to V$ is defined on the whole $X \tm U$ and is mapping to a Banach space $V$. Furthermore, $B:U\to X$, $B_2:V\to X$ are possibly unbounded operators, that belong however to $L(U,X_{-1})$ and $L(V,X_{-1})$ respectively, where $V$ is a Banach space.
Here the extrapolation space $X_{-1}$ is the closure of $X$ in the norm $x \mapsto \|(aI-A)^{-1}x\|_X$, $x \in X$, where $a \in \rho(A)$ (different choices of $a \in\rho(A)$ induce equivalent norms on $X$).

\subsection{Admissible input operators and mild solutions}
\label{sec:Admissible input operators and mild solutions}

Given our analysis of linear theory, we impose
\begin{ass}
\label{ass:Admissibility} 
The operator $B\in L(U,X_{-1})$ is $\infty$-admissible, and the map $(t,u) \mapsto \Phi(t)u$ is continuous on $\R_+\tm L^q(\R_+,U)$. 

In particular, this assumption holds if $B$ is a $q$-admissible operator with $q<\infty$, see \cite[Proposition 2.3]{Wei89b}.
\end{ass}

To define the concept of a mild solution, we also require the following:
\begin{ass}
\label{ass:Integrability} 
Let Assumption~\ref{ass:Admissibility} hold. We assume that for all $u\in L^\infty(\R_+,U)$ and any $x\in C(\R_+,X)$ the map $s\mapsto f\big(x(s),u(s)\big)$ is in $L^\infty_{\loc}(\R_+,V)$, $B_2$ is zero-class $\infty$-admissible, and the map 
\[
t \mapsto \int_0^t T(t-s)B_2f\big(x(s),u(s)\big)ds
\]
is well-defined and continuous on $\R_+$.
\end{ass}

\begin{remark}
\label{rem:Validity-of-integrability-assumption} 
Assumption~\ref{ass:Integrability} holds, in particular, if $B_2 \in L(V,X)$, and
\begin{enumerate}[label=(\roman*)]
	\item $f(x,u)=g(x) + Ru$, $x \in X$, $u\in U$, where $R\in L(U,V)$, and $g$ is continuous on $X$.
Indeed, for a continuous $x$, the map	$s\mapsto g\big(x(s)\big)$ is continuous either, and thus is Riemann integrable. 
The map $s\mapsto T(t-s) B_2Ru(s)$ is Bochner integrable for any $u \in L^1_{\loc}(\R_+,U)$ by \cite[Proposition 1.3.4]{ABH11}, \cite[Lemma 10.1.6]{JaZ12}. 
This ensures that Assumption~\ref{ass:Integrability} holds.

\item If $f$ is continuous on $X \tm U$, and $u$ is piecewise right-continuous, then the map 
 $s\mapsto f\big(x(s),u(s)\big)$ is also piecewise right-continuous, and thus it is Riemann integrable.

\item (ODE systems). Let $X=\R^n$, $U=\R^m$, $A=0$ (and thus $T(t)=\id$ for all $t$), $B_2=\id$, $B=0$, and $f$ be continuous on $X \tm U$. With these assumptions the equations \eqref{eq:SEE+admissible} take form
\begin{align}
\label{eq:ODE}
\dot{x} = f(x,u).
\end{align}
Then for each $u\in L^\infty(\R_+,U)$ and each $x \in C(\R_+,X)$ the map $s\mapsto f(x(s),u(s))$ is Lebesgue integrable, and thus Assumption~\ref{ass:Integrability} holds.

Indeed, as $x$ is a solution of \eqref{eq:ODE} on $[0,\tau)$, $x$ is continuous on $[0,\tau)$. 
By assumptions, $u$ is measurable on $[0,\tau)$, and $f$ is continuous on $\R^n \tm \R^m$. Arguing similarly to  
\cite[Proposition 7]{RoF10} (where it was shown that a composition of a continuous and measurable function defined on a measurable set $E$ is measurable on $E$), we see that the map $q : [0,\tau)\to\R^n$, $q(s):= f(x(s),u(s))$, is a measurable map.
As $u$ is essentially bounded, and $x$ and $f$ map bounded sets into bounded sets, $q$ is essentially bounded on $[0,\tau)$. 
Thus, $q \in L^{\infty}(\R_+,\R^n)$, and thus $q$ is integrable on $[0,\tau)$.
\end{enumerate}
\end{remark}

Next, we define mild solutions of \eqref{eq:SEE+admissible}.
\begin{definition}[Mild solutions]
\label{def:Mild-solution}
\index{solution!mild}
Let Assumptions~\ref{ass:Admissibility}, \ref{ass:Integrability} hold. Let also $\tau>0$ be given. 
A function $x \in C([0,\tau], X)$ is called a \emph{mild solution of \eqref{eq:SEE+admissible} on $[0,\tau]$} corresponding to certain $x_0\in X$ and $u \in L^\infty_{\loc}(\R_+,U)$, if $x$ solves the integral equation
\begin{align}
\label{eq:SEE+admissible_Integral_Form}
x(t)=T(t) x_0 + \int_0^t T(t-s) B_2f\big(x(s),u(s)\big)ds + \int_0^t T_{-1}(t-s) Bu(s)ds. 
\end{align}
Here the integrals are Bochner integrals of $X$-valued maps, and $T_{-1}$ is an extension of the semigroup $T$ to the space $X_{-1}$, see \cite[p.126]{EnN00}.

We say that $x:\R_+\to X$ is a \emph{mild solution of \eqref{eq:SEE+admissible} on $\R_+$} corresponding to 
certain $x_0\in X$ and $u \in L^\infty_{\loc}(\R_+,U)$, if it is a mild solution of \eqref{eq:SEE+admissible} (with $x_0, u$) on $[0,\tau]$ for all $\tau>0$.
\end{definition}

\subsection{Local existence and uniqueness}
\label{sec:Local existence and uniqueness}

Assumptions~\ref{ass:Admissibility}, \ref{ass:Integrability} guarantee that the integral terms in \eqref{eq:SEE+admissible_Integral_Form}
are well-defined. To ensure the existence and uniqueness of mild solutions, we impose further restrictions on $f$.
 
\begin{definition}
We call $f:X \times U \to V$
\begin{enumerate}[label=(\roman*)]
\index{function!Lipschitz continuous}
	\item \emph{Lipschitz continuous  (with respect to the first argument) on bounded subsets of $X$} if for any
$C>0$ there is $L(C)>0$, such that $\forall x,y \in B_C$, $\forall v \in B_{C,U}$ it holds that
\begin{eqnarray}
\|f(y,v)-f(x,v)\|_V \leq L(C) \|y-x\|_X.
\label{eq:Lipschitz}
\end{eqnarray}	
	
  \item  \emph{Lipschitz continuous  (with respect to the first argument) on bounded subsets of $X$, uniformly with respect to the other argument} if for any $C>0$ there is $L(C)>0$, such that \eqref{eq:Lipschitz} holds for all $x,y: \|x\|_X \leq C,\ \|y\|_X \leq C$, and all $v \in U$.
	\item \emph{uniformly globally Lipschitz continuous (with respect to the first argument)} if \eqref{eq:Lipschitz} holds for all $x,y \in X$, and all $u \in \Uc$ with $L = const$ that does not depend on $x,y,u$.
\end{enumerate}
\end{definition}

We omit the indication \q{with respect to the first argument} wherever this is clear
from the context. 

For the well-posedness analysis, we rely on the following assumption on the nonlinearity $f$ in \eqref{eq:SEE+admissible}.

\begin{ass}
\label{Assumption1} 
The nonlinearity $f$ satisfies the following properties:
\begin{itemize}
    \item[(i)] $f: X \times U \to V$ is Lipschitz continuous on bounded subsets of $X$.
    \item[(ii)] $f(x,\cdot)$ is continuous for all $x \in X$.
		\item[(iii)] There exist $\sigma \in \Kinf$ and $c>0$ so that for all $u \in U$ the following holds:
\begin{eqnarray}
\|f(0,u)\|_V \leq \sigma(\|u\|_U) + c.
\label{eq:f0u_estimate}
\end{eqnarray}
\end{itemize}
\end{ass}

Let $S$ be a normed vector space.  
Define the \emph{distance} from $z\in S$ to the set $Z \subset S$ by 
\[
\dist(z,Z):=\inf\{\|y-z\|_S: y \in Z \}.
\]
Further, denote the open ball of radius $r$ around $Z \subset S$ by 
\[
B_{r,S}(Z):=\{y \in X:\dist(y,Z)<r\}.
\]
We set also $B_{r,S}(x):=B_{r,S}(\{x\})$ for $x \in X$, and $B_{r,S}:=B_{r,S}(0)$.
If $S = X$, we write for short $B_{r}(Z):=B_{r,X}(Z)$, $B_{r}(x):=B_{r,X}(x)$, etc.

Finally, for a set $\Sc \subset U$, denote the set of inputs with essential values in $\Sc$ as $\Uc_\Sc$:
\begin{align}
\label{eq:Inputs-constrained}
\Uc_\Sc:=\{u\in\Uc: u(t) \in \Sc, \text{ for a.e. } t\in\R_+\}.
\end{align}

We start with the following sufficient condition for the existence and uniqueness of solutions of a system \eqref{eq:SEE+admissible} with inputs in $\Uc:=L^\infty(\R_+,U)$.
\index{theorem!Picard-Lindel\"of}
\index{Picard-Lindel\"of theorem}
\begin{theorem}[Picard-Lindel\"of theorem]
\label{PicardCauchy}
Let Assumptions~\ref{ass:Admissibility}, \ref{ass:Integrability}, \ref{Assumption1} hold. 
Define $h_0:=\lim_{t\to +0}h_t$, where $h_t$ is defined as in \eqref{eq:admissible-operator-norm-estimate}.

Assume that $(T(t))_{t\ge 0}$ is a strongly continuous semigroup, generated by $A$ and satisfying for certain $M \geq 1$, $\lambda>0$ the estimate
\begin{eqnarray}
\|T(t)\| \leq Me^{\lambda t},\quad t\geq 0.
\label{eq:Picard-Lindeloef-bounds-semigroup}
\end{eqnarray} 
For any compact set $Q \subset X$, any $r>0$, any bounded set $\Sc \subset U$, and any $\delta>0$, there is a time $t_1 = t_1(Q,r,\Sc,\delta)>0$, such that for any $x_0 \in W:=B_r(w)$ for some $w \in Q$, and for any $u\in\Uc_{\Sc}$ there is a unique mild solution of \eqref{eq:SEE+admissible} on $[0,t_1]$, and $\phi([0,t_1],x_0,u) \subset B_{Mr+h_0\|u\|_{L^\infty([0,t_1],U)}+\delta}(w)$.
\end{theorem}

\begin{proof}
First, we show the claim for the case if $Q$ is a single point in $X$, that is, $Q=\{w\}$, for some $\omega\in X$. 
Pick any $C>0$ such that $W := B_r(w) \subset B_C$, and $\Uc_\Sc \subset B_{C,\Uc}$. Pick any $u \in \Uc_\Sc$. Also take any $\delta>0$, and consider the following sets (depending on the parameter $t>0$):
\begin{eqnarray}
\hspace{10mm}Y_{t}:= \big\{x \in C([0,t],X): \sup_{s \in[0,t]} \|x(s) - w\|_X\leq Mr + h_0\|u\|_{L^\infty([0,t],U)} + \delta \big\},
\label{eq:Y_T_Def}
\end{eqnarray}
endowed with the metric 
$\rho_{t}(x,y):=\sup_{s \in [0,t]} \|x(s)-y(s)\|_X$.
As the sets $Y_t$ are closed subsets of the Banach spaces $C([0,t],X)$, for all $t>0$, the space $Y_t$ is a complete metric space.

Pick any $x_0 \in W$. We are going to prove that for small enough $t$, the spaces $Y_{t}$  are invariant under the operator
$\Phi_u$, defined for any $x \in Y_{t}$ and all $\tau \in [0,t]$ by
\begin{align}
\label{eq:Phi-map-PL-Theorem}
\Phi_u(x)(\tau) 
&= T(\tau)x_0 + \int_0^\tau T(\tau-s)B_2f\big(x(s),u(s)\big)ds  + \int_0^\tau T_{-1}(\tau-s) Bu(s)ds.
\end{align}
By Assumptions~\ref{ass:Admissibility}, \ref{ass:Integrability}, the function $\Phi_u(x)$ is continuous for any $x\in Y_t$.

Fix any $t>0$ and pick any $x \in Y_{t}$.
As $x_0 \in W=B_r(w)$, there is $a \in B_r$ such that $x_0=w+a$.
 
Then for any $\tau<t$, it holds that
\begin{align*}
\|\Phi_{u}(x)&(\tau)-w\|_X \\
&\leq   \Big\|T(\tau)x_0 - w\Big\|_X + \Big\|\int_0^\tau T_{-1}(\tau-s) Bu(s)ds\Big\|_X \\
&\qquad\qquad\qquad\qquad\qquad + \Big\|\int_0^\tau T_{-1}(\tau-s) B_2f(x(s),u(s))ds\Big\|_X \\
&\leq  \|T(\tau)(w+a) - w \|_X +  h_\tau \|u\|_{L^\infty([0,\tau],U)} \\
&\qquad\qquad\qquad\qquad\qquad + c_\tau \|f(x(\cdot),u(\cdot))\|_{L^\infty([0,\tau],V)} \\
			&\leq \|T(\tau)w - w\|_X + \|T(\tau)a\|_X + h_\tau \|u\|_{L^\infty([0,\tau],U)}\\
			&\quad+ c_\tau \|f(x(\cdot),u(\cdot)) - f(0,u(\cdot))\|_{L^\infty([0,\tau],V)} + c_\tau \|f(0,u(\cdot))\|_{L^\infty([0,\tau],V)}.
\end{align*}
Now for all $s\in[0,t]$
\begin{align*}
\|x(s)\|_X &\leq \|w\|_X+Mr+ h_0\|u\|_{L^\infty([0,t],U)} +\delta\\
&\leq M(\|w\|_X+r) + h_0 C+\delta\leq (M+h_0)C+\delta =:K.
\end{align*}
In view of Assumption~\ref{Assumption1}(iii), it holds that 
\[
\|f(0,u(s))\|_V \leq \sigma(\|u(s)\|_U) + c,\quad \text{ for a.e. } s \in [0,t].
\]
As $M\geq 1$, it holds that $K>C$, and the Lipschitz continuity of $f$ on bounded balls ensures that there is $L(K)>0$, such that for all $\tau\in[0,t]$ 			
\begin{align*}
\|\Phi_{u}(x)(\tau)-w\|_X &\leq \|T(\tau)w - w\|_X + Me^{\lambda t} r + h_t \|u\|_{L^\infty([0,t],U)} \\
					&\qquad\qquad+ c_\tau \big( L(K) \|x\|_{L^\infty([0,t],X)} + \sigma(\|u\|_{L^\infty([0,t],U)}) + c\big)\\
			&\leq \|T(\tau)w - w\|_X + Me^{\lambda t} r + h_t \|u\|_{L^\infty([0,t],U)} \\
					&\qquad\qquad+ c_t \big( L(K)K + \sigma(C) + c\big).
\end{align*}
Since $T$ is a strongly continuous semigroup, $h_t\to h_0$ as $t\to +0$, and $c_t\to 0$ as $t\to +0$, from this estimate, it is clear that  there exists $t_1$, such that 
\[
\|\Phi_u(x)(t)-w\|_X \leq Mr+ h_0 \|u\|_{L^\infty([0,t_1],U)} +\delta,\quad \text{ for all } t \in [0,t_1].
\]
This means, that $Y_{t}$ is invariant with respect to $\Phi_u$ for all $t \in (0,t_1]$, and $t_1$ does not depend on the choice of $x_0 \in W$.

Now pick any $t>0$, $\tau \in [0, t]$, and any $x, y \in Y_{t}$. Then it holds that
\begin{align*}
\|\Phi_u(x)(\tau) - \Phi_u(y)(\tau)\|_X 
&\leq \Big\|\int_0^\tau T(\tau-s)B_2\big(f(x(s),u(s)) - f(y(s),u(s))\big)ds\Big\|_X \\
&\leq c_\tau \|f(x(\cdot),u(\cdot)) - f(y(\cdot),u(\cdot))\|_{L^\infty([0,\tau],V)} \\
&\leq c_t L(K)\rho_t(x,y) \\
&\leq \frac{1}{2} \rho_t(x,y),
\end{align*}
for $t \leq t_2$, where $t_2>0$ is a small enough real number that does not depend on the choice of $x_0 \in W$.

According to the Banach fixed point theorem, there exists a unique solution of $x(t)=\Phi_u(x)(t)$ on $[0,\min\{t_1,t_2\}]$, which is a  mild solution of \eqref{eq:SEE+admissible}.

\textbf{General compact $Q$.} Till now, we have shown that for any $w\in Q$, any $r>0$, any bounded set $\Sc \subset U$, and any $\delta>0$, there is a time $t_1 = t_1(w,r,\Sc,\delta)>0$ (that we always take the maximal possible), such that for any $x_0 \in W:=B_r(w)$, and for any $u\in\Uc_{\Sc}$ there is a unique solution of \eqref{eq:SEE+admissible} on $[0,t_1]$, and it lies in the ball $B_{Mr+h_0\|u\|_{L^\infty([0,t_1],U)}+\delta}(w)$.

It remains only to show that $t_1$ can be chosen uniformly in $w \in Q$, that is, $\inf_{w \in Q} t_1(w,r,\Sc,\delta)>0$. 
Let this not be so, that is, $\inf_{w \in Q} t_1(w,r,\Sc,\delta)=0$. Then there is a sequence 
$(w_k) \subset Q$, such that the corresponding times $\big(t_1(w_k,r,\Sc,\delta)\big)_{k\in\N}$ monotonically decay to zero.
As $Q$ is compact, there is a converging subsequence of $(w_k)$, converging to some $w^* \in Q$. 
However, $t_1(w^*,r,\Sc,\delta)>0$, which easily leads to a contradiction.
\end{proof}

\begin{corollary}[Picard-Lindel\"of theorem for zero-class admissible $B$ and quasi-contractive semigroups]
\label{cor:PL-for zero-class admissible operators and quasi-contractive semigroups} 
Let Assumptions~\ref{ass:Admissibility}, \ref{ass:Integrability}, \ref{Assumption1} hold.
Let also $B$ be zero-class admissible, and $T$ be a quasi-contractive strongly continuous semigroup. That is, there is $\lambda>0$ such that 
\begin{eqnarray*}
\|T(t)\| \leq e^{\lambda t},\quad t\geq 0.
\end{eqnarray*} 
For any bounded ball $W = B_r(w) \subset X$, any bounded set $\Sc \subset U$, and any $\delta>0$, there is a time $t_1 = t_1(W,\Sc,\delta)>0$, such that for any $x_0 \in W$ and any $u\in\Uc_{\Sc}$ there is a unique solution of \eqref{eq:SEE+admissible} on $[0,t_1]$, and it lies in the ball $B_{r+\delta}(w)$.
\end{corollary}

\begin{proof}
The claim follows directly from Theorem~\ref{PicardCauchy}.
\end{proof}

\begin{remark}
\label{rem:Quasicontractive-SGR-Picard-Lindeloef} 
Without an assumption of quasicontractivity, Corollary~\ref{cor:PL-for zero-class admissible operators and quasi-contractive semigroups} does not hold. 
Consider the special case $f \equiv 0$ and $B\equiv 0$. Then the system \eqref{eq:SEE+admissible} is linear, and for a given $x_0\in X$ the solution of \eqref{eq:SEE+admissible} exists globally and equals $t\mapsto T(t)x_0$. Now take $w:=0$ and pick any $r>0$ and $t_1>0$. Then
\[
\sup_{\tau \in [0,t_1]}\sup_{\|x\|_X\leq r}\|T(\tau)x\|_X= r\sup_{\tau \in [0,t_1]}\|T(\tau)\|.
\]
Since $T$ is merely strongly continuous, the map $t\mapsto \|T(t)\|$ does not have to be continuous at $t=0$, and it may happen that $\lim_{t_1\to 0}\sup_{\tau \in [0,t_1]}\|T(\tau)\| >1$.

Hence, in general, it is not possible to prove that the solution starting at arbitrary $x_0 \in B_r(w)$ will stay in $B_{r+\delta}(w)$ during a sufficiently small and uniform in $x_0 \in B_r(w)$ time. 
\end{remark}

The following example shows that Theorem~\ref{PicardCauchy} does not hold in general if $W$ is a bounded set (and not only a bounded ball over a compact set), even for linear systems governed by contraction semigroups on a Hilbert space.

\begin{example}
\label{examp:PL-difference-to-ODE-case} 
Let $X=\ell_2$, and consider a diagonal semigroup, defined by $T(t)x:=(e^{-kt}x_k)_k$, for all $x=(x_k)_k\in X$ and all $t\geq 0$. This semigroup is strongly continuous and contractive.
Consider a bounded and closed set $W:=\{x \in \ell_2:\ \|x\|_X =1\}$. 
Yet $\|T(t)e_k\|_X = e^{-kt}$, and thus for each $\delta>0$ and for each time $t_1>0$, we can find $k\in\N$, such that 
$\|T(t_1)e_k\|_X<1-\delta$, which means that $T(t_1)e_k \notin B_\delta(W)$.
\end{example}

At the same time, a stronger Picard-Lindel\"of-type theorem can be shown for uniformly continuous semigroups (this encompasses, in particular, the case of infinite ODE systems, also called \q{ensembles}), which fully extends the corresponding result for ODE systems, see \cite[Chapter 1]{Mir23}.

\begin{theorem}[Picard-Lindel\"of theorem for uniformly continuous semigroups]
\label{thm:PicardCauchy-Uniformly-cont-semigroups}
Let Assumptions~\ref{ass:Admissibility}, \ref{ass:Integrability}, \ref{Assumption1} hold, and let $B$ be zero-class 
$\infty$-admissible. 
Let further $T$ be a uniformly continuous semigroup (not necessarily quasicontractive).
For any bounded set $W \subset X$, any bounded set $\Sc \subset U$ and any $\delta>0$, there is a time $\tau = \tau(W,\Sc,\delta)>0$, such that for any $x_0 \in W$, and $u\in\Uc_{\Sc}$ there is a unique solution of \eqref{eq:SEE+admissible} on $[0,\tau]$, and it lies in $B_\delta(W)$.
\end{theorem}

\begin{proof}
Pick any $C>0$ such that $W \subset B_C$, and $\Uc_\Sc \subset B_{C,\Uc}$. 

Take also any $\delta>0$, and consider the following sets (depending on a parameter $t>0$):
\begin{eqnarray}
Y_{t}:= \big\{x \in C([0,t],X): \dist(x(t),W) \leq \delta \ \ \forall t \in [0,t]\big\},
\label{eq:Y_T_Def-uniform}
\end{eqnarray}
endowed with the metric 
$\rho_{t}(x,y):=\sup_{s \in [0,t]} \|x(s)-y(s)\|_X$.
As the sets $Y_t$ are closed subsets of Banach spaces $C([0,t],X)$, $Y_t$ are complete metric spaces for all $t>0$.

Pick any $x_0 \in W$ and any $u \in \Uc_\Sc$. We are going to prove that for small enough $t$, the spaces $Y_{t}$  are invariant under the operator
$\Phi_u$, defined for any $x \in Y_{t}$ and all $\tau \in [0,t]$ by \eqref{eq:Phi-map-PL-Theorem}, which we recall next once again:
\begin{align*}
\Phi_u(x)(\tau) 
&= T(\tau)x_0 + \int_0^\tau T(\tau-s)B_2f\big(x(s),u(s)\big)ds  + \int_0^\tau T_{-1}(\tau-s) Bu(s)ds.
\end{align*}
By Assumptions~\ref{ass:Admissibility}, \ref{ass:Integrability}, the function $\Phi_u(x)$ is continuous.

Fix any $t>0$ and pick any $x \in Y_{t}$.  
Then for any $\tau<t$, it holds that
\begin{align*}
\dist(\Phi_u(x)(\tau),W) &\le \|\Phi_u(x)(\tau)-x_0\|_X \\
&\leq   \Big\|T(\tau)x_0 - x_0\Big\|_X + \Big\|\int_0^\tau T_{-1}(\tau-s) Bu(s)ds\Big\|_X \\
&\qquad\qquad\qquad\qquad\qquad + \Big\|\int_0^\tau T_{-1}(\tau-s) B_2f(x(s),u(s))ds\Big\|_X \\
&\leq   \|T(\tau)-I \|_X \|x_0\|_X +  h_\tau \|u\|_{L^\infty([0,\tau],U)} \\
&\qquad\qquad\qquad\qquad\qquad + c_\tau \|f(x(\cdot),u(\cdot))\|_{L^\infty([0,\tau],V)} \\
			&\leq C \|T(\tau)-I \|_X  + h_\tau \|u\|_{L^\infty([0,\tau],U)}\\
			&\quad+ c_\tau \|f(x(\cdot),u(\cdot)) - f(0,u(\cdot))\|_{L^\infty([0,\tau],V)} + c_\tau \|f(0,u(\cdot))\|_{L^\infty([0,\tau],V)}.
\end{align*}
Now for all $s\in[0,t]$
\begin{align*}
\|x(s)\|_X &\leq C +\delta =: K.
\end{align*}
In view of Assumption~\ref{Assumption1}(iii), it holds that 
\[
\|f(0,u(s))\|_V \leq \sigma(\|u(s)\|_U) + c,\quad \text{ for a.e. } s \in [0,t].
\]
Now Lipschitz continuity of $f$ on bounded balls ensures that there is $L(K)>0$, such that for all $\tau\in[0,t]$ 			
\begin{align*}
\|\Phi_{t}(x)(\tau)-w\|_X 		
			&\leq C \|T(\tau)-I \|_X + h_t \|u\|_{L^\infty([0,t],U)} \\
					&\quad+ c_\tau \big( L(K) \|x\|_{L^\infty([0,t],X)} + \sigma(\|u\|_{L^\infty([0,t],U)}) + c\big)\\
			&\leq C \|T(\tau)-I \|_X + h_t \|u\|_{L^\infty([0,t],U)} \\
					&\quad+ c_t \big( L(K)K + \sigma(C) + c\big).
\end{align*}
Since $T$ is a uniformly continuous semigroup, $h_t\to 0$ as $t\to +0$, and $c_t\to 0$ as $t\to +0$, from this estimate it is clear that  there exists $t_1>0$, depending solely on $C$ and $\delta$, such that 
\[
\dist(\Phi_u(x)(t),W)\le\delta,\quad \text{ for all } t \in [0,t_1].
\]
This means, that $Y_{t}$ is invariant with respect to $\Phi_u$ for all $t \in (0,t_1]$, and $t_1$ does not depend on the choice of $x_0 \in W$.
The rest of the proof is analogous to the proof of Theorem~\ref{PicardCauchy}.
%
%
%
\end{proof}

\begin{remark}
\label{rem:Bounds-for-inputs-in-PL-theorem} 
In the above theorems, the existence time of solutions is uniform with respect to the initial values and inputs with norms bounded by $C$. In general, the existence time of solutions cannot be chosen to be uniform if we do not restrict the norms of the inputs. 
For example, consider the system
\begin{eqnarray*}
\dot{x}(t) =  u(t) x^2(t).
\end{eqnarray*}
The maximal existence time of solutions with $x(0)=1$ and $u(t) \equiv k$ goes to zero as $k \to \infty$.
\end{remark}

\subsection{Well-posedness}
\label{sec:Global-Well-posedness}

\index{solution!extension}
Our next aim is to study the prolongations of solutions and their asymptotic properties.
\begin{definition}
Let $x_1(\cdot)$, $x_2(\cdot)$ be mild solutions of \eqref{eq:SEE+admissible} defined on the intervals $[0,t_1)$ and $[0,t_2)$ respectively, $t_1,t_2>0$.
We call $x_2$ an \emph{extension} of $x_1$ if $t_2>t_1$, and $x_2(t)=x_1(t)$ for all $t \in [0,t_1)$.
\end{definition}

\begin{lemma}
\label{lem:Equality-of-solutions} 
Let Assumptions~\ref{ass:Admissibility}, \ref{ass:Integrability}, \ref{Assumption1} hold.
Take any $x_0 \in X$ and $u\in\Uc$. 
Any two solutions of \eqref{eq:SEE+admissible} coincide in their common domain of existence. 
\end{lemma}

\begin{proof}
Let $\phi_1,\phi_2$ be two mild solutions of \eqref{eq:SEE+admissible} corresponding to $x_0 \in X$ and $u\in\Uc$, defined over $[0,t_1)$
and $[0,t_2)$ respectively. Assume that $t_1\leq t_2$ (the other case is analogous).

By local existence and uniqueness theorem, there is some positive $t_3 \leq t_1$ (we take $t_3$ to be the maximal of such times), such that $\phi_1$ and $\phi_2$ coincide on $[0,t_3)$. If $t_3 = t_1$, then the claim is shown. If $t_3 < t_1$, then by continuity $\phi_1(t_3)=\phi_2(t_3)$.
Now $\psi_1:t\mapsto \phi_1(t_3+t)$ and $\psi_2:t\mapsto \phi_2(t_3+t)$ are two mild solutions for the problem
\[
\dot{x}(t)  = Ax(t) + B_2 f(x(t),u(t+t_3)) + Bu(t_3+\cdot),\quad x(0)= \phi_1(t_3),
\]
By local existence and uniqueness, $\phi_1$ and $\phi_2$ coincide on $[0,t_3+\varepsilon)$ for some $\varepsilon>0$, which contradicts to the maximality of $t_3$.
Hence, $\phi_1$ and $\phi_2$ coincide on $[0,t_1)$.
\end{proof}

\index{solution!maximal}
\index{solution!global}
\begin{definition}
A solution $x(\cdot)$ of \eqref{eq:SEE+admissible} is called 
\begin{enumerate}[label=(\roman*)]
	\item \emph{maximal} if there is no solution of \eqref{eq:SEE+admissible} that extends $x(\cdot)$,
	\item \emph{global} if $x(\cdot)$ is defined on $\R_+$.
\end{enumerate}

\end{definition}

A central property of the system \eqref{eq:SEE+admissible} is
\begin{definition}
\label{def:Well-posedness-finite-dim} 
\index{system!well-posed}
We say that the system \eqref{eq:SEE+admissible} is \emph{well-posed} if for every initial value $x_0 \in X$ and every external input $u \in \Uc$, there exists a unique maximal solution that we denote by $\phi(\cdot,x_0,u):[0,t_m(x_0,u)) \rightarrow X$, where $0 < t_m(x_0,u) \leq \infty$.

We call $t_m(x_0,u)$ the \emph{maximal existence time} of a solution corresponding to $(x_0,u)$.
\end{definition}

\index{flow map}
\index{maximal existence time}

The map $\phi$, defined in Definition~\ref{def:Well-posedness-finite-dim}, and describing the evolution of the system \eqref{eq:SEE+admissible}, 
is called the \emph{flow map}, or just \emph{flow}.
The domain of definition of the flow $\phi$ is 
\[
D_\phi:=\cup_{x_0\in X,\ u\in\Uc}[0,t_m(x_0,u))\times\{(x_0,u)\}.
\] 
We will always deal with maximal solutions. 

We proceed with:
\begin{theorem}[Well-posedness]
\label{thm:Global_well-posedness}
Let Assumptions~\ref{ass:Admissibility}, \ref{ass:Integrability}, \ref{Assumption1} hold. 
Then \eqref{eq:SEE+admissible} is well-posed.
\end{theorem}

\begin{proof}
Pick any $x_0\in X$ and $u\in\Uc$.
Let $\Sc$ be the set of all solutions $\xi$ of \eqref{eq:SEE+admissible}, defined on the corresponding domain of definition $I_\xi=[0,t_{\xi})$.

Define $I:=\cup_{\xi\in\Sc}I_\xi = [0,t^*)$ for some $t^*>0$, which can be either finite or infinite.
For any $t\in I$ there is $\xi\in\Sc$, such that $t \in I_\xi$. 
Define $\xi_{\max}(t):=\xi(t)$.
Because of Lemma~\ref{lem:Equality-of-solutions}, the map $\xi_{\max}$ does not depend on the choice of $\xi$ above and thus is well-defined on $I$.

Furthermore, for each $t\in(0,t^*)$ there is $\xi \in\Sc$, such that $\xi_{\max} = \xi$ on $[0,t+\varepsilon) \subset I$, for a certain $\varepsilon>0$, and thus $\phi(\cdot,x_0,u):=\xi_{\max}$ is a solution of \eqref{eq:SEE+admissible} corresponding to $x_0,u$.
Finally, it is a maximal solution by construction.
\end{proof}

Now we show that well-posed systems \eqref{eq:SEE+admissible} are a special case of general control systems, introduced in Definition~\ref{Steurungssystem}.
\begin{theorem}
\label{thm:SEE-as-control systems} 
Let \eqref{eq:SEE+admissible} be well-posed. Then the triple $(X,\Uc,\phi)$, where $\phi$ is a flow map of \eqref{eq:SEE+admissible}, constitutes a control system in the sense of Definition~\ref{Steurungssystem}.
\end{theorem}

\begin{proof}
By construction, all the axioms in the definition of a control system are fulfilled. 
In particular, continuity holds by the definition of a mild solution.
We only need to check the cocycle property. 

Take  any initial condition $x \in X$, any input $u \in \Uc$, and any $t,\tau\geq 0$, such that 
$[0,t+\tau]\tm \{(x,u)\} \subset D_{\phi}$.
Define an input $v$ by $v(r)=u(r+\tau)$, $r \geq 0$.

Due to \eqref{eq:SEE+admissible_Integral_Form}, we have:
\begin{align*}
\phi(t &+ \tau ,x,u) =  T(t+\tau)x + \int_0^{t+\tau} T_{-1}(t+\tau-s) B_2f(\phi(s,x,u),u(s))ds\\
&\qquad\qquad\qquad\qquad\qquad + \int_0^{t+\tau} T_{-1}(t+\tau-s) Bu(s)ds.
\end{align*}
As $T_{-1}(t)$ is a bounded operator, it can be taken out of the Bochner integral:
\begin{align*}
\phi(t &+ \tau ,x,u) =  T(t)T(\tau)x + T_{-1}(t)\int_0^{\tau} T_{-1}(\tau-s) B_2f(\phi(s,x,u),u(s))ds  \\
&\quad+ T_{-1}(t)\int_0^{\tau} T_{-1}(\tau-s) Bu(s)ds \\
&\quad + \int_{\tau}^{t+\tau} T_{-1}(t+\tau-s) B_2f(\phi(s,x,u),u(s))ds + \int_{\tau}^{t+\tau} T_{-1}(t+\tau-s) Bu(s)ds.
\end{align*}
As $B$ is admissible, we have that $\int_0^{\tau} T_{-1}(\tau-s) Bu(s)ds \in X$. 
Since $T_{-1}(\cdot)$ coincides with $T(\cdot)$ on $X$, we infer
\[
T_{-1}(t)\int_0^{\tau} T_{-1}(\tau-s) Bu(s)ds = T(t)\int_0^{\tau} T_{-1}(\tau-s) Bu(s)ds.
\]
A similar argument applies to the second term of the above sum. Finally,
\begin{align*}
\phi(t + \tau ,x,u) &= T(t) \phi(\tau,x,u) + \int_{0}^{t} T_{-1}(t-s) B_2 f(\phi(s+\tau,x,v),v(s))ds \\
&\qquad\qquad\qquad\qquad\qquad + \int_{0}^{t} T_{-1}(t-s) Bv(s)ds\\
 & =  \phi(t,\phi(\tau,x,u),v),
\end{align*}
and the cocycle property holds.
\end{proof}

We proceed with a general asymptotic property of the flow.
\begin{proposition}
\label{prop:Unboundedness_finite_existence_time}
Let Assumptions~\ref{ass:Admissibility}, \ref{ass:Integrability}, \ref{Assumption1} hold. 
Pick any $x \in X$ and any $u \in \Uc$. 
If $t_m(x,u)$ is finite, then $\|\phi(t,x,u)\|_X \to \infty$ as $t \to t_m(x,u)-0$.
\end{proposition}

\begin{proof}
Pick any $x \in X$ and any $u \in \Uc$, and consider the corresponding maximal solution $\phi(\cdot,x,u)$, defined on $[0,t_m(x,u))$.
Assume that $t_m(x,u)<+\infty$, but at the same time $\Liminf_{t\to t_m(x,u)-0}\|\phi(t,x,u)\|_X <\infty$.
Then there is a sequence $(t_k)$, such that $t_k \to t_m(x,u)$ as $k\to\infty$ and 
$\lim_{k\to\infty}\|\phi(t_k,x,u)\|_X <\infty$. Then also $\sup_{k\in\N}\|\phi(t_k,x,u)\|_X =:C<\infty$.

Let $\tau(C)>0$ be a uniform existence time for the solutions starting in the ball $\clo{B_C}$ subject to inputs of a magnitude not exceeding $\|u\|$, which exists and is positive in view of  Theorem~\ref{PicardCauchy}.
Then the solution of \eqref{eq:SEE+admissible} starting in $\phi(t_k,x,u)$, corresponding to the input $u(\cdot + t_k)$, exists and is unique on $[0,\tau(C)]$ by Theorem~\ref{PicardCauchy}, and by the cocycle property, we have that 
$\phi(\cdot,x,u)$ can be prolonged to $[0,t_k +\tau(C))$, which (since $t_k\to t_m(x,u)$ as $k\to\infty$) contradicts to the maximality of the solution corresponding to $(x,u)$.

Hence $\Liminf_{t\to t_m(x,u)-0}\|\phi(t,x,u)\|_X =\infty$, which implies the claim.
\end{proof}

As a corollary of Proposition~\ref{prop:Unboundedness_finite_existence_time}, we obtain that any uniformly bounded maximal solution of \eqref{eq:SEE+admissible} is a global solution.

\begin{corollary}[BIC property]
\label{cor:BIC-property}
Let Assumptions~\ref{ass:Admissibility}, \ref{ass:Integrability}, \ref{Assumption1} hold. 
The system \eqref{eq:SEE+admissible} satisfies the BIC property.
\end{corollary}

\subsection{Forward completeness and boundedness of reachability sets}
\label{sec:Boundedness of reachability sets}

Local Lipschitz continuity guarantees the existence of local solutions. To ensure the global existence of solutions, stronger requirements on nonlinearity are needed.
\begin{proposition}
\label{prop:Unif-Glob-Lip-and-BRS}
Let Assumptions~\ref{ass:Admissibility}, \ref{ass:Integrability}, \ref{Assumption1} hold.
Let further $f$ be uniformly globally Lipschitz.
Then \eqref{eq:SEE+admissible} is forward complete and has the BRS property.
\end{proposition}

\begin{proof}
By Theorem~\ref{PicardCauchy}, for any $x_0 \in X$ and for any $u \in \Uc$, there exists a mild solution of \eqref{eq:SEE+admissible}, with a maximal existence time $t_m (x_0, u)$, which may be finite or infinite. Let $t_m (x_0, u)$ be finite. 

Let $L>0$ be a uniform global Lipschitz constant for $f$. 
As $\|T(t)\|\leq Me^{\lambda t}$ for some $M\ge 1$, $\lambda\geq 0$ and all $t\geq 0$, for any $t<t_m (x_0, u)$ we have according to the formula \eqref{eq:SEE+admissible_Integral_Form} the following estimates
\begin{align*}
\|\phi(&t,x_0,u)\|_X \leq  \|T(t)\| \|x_0\|_X + \Big\|\int_{0}^{t} T_{-1}(t-s) Bu(s)ds\Big\|_X \\
&\qquad\qquad\qquad\qquad + \Big\|\int_0^t T(t-s)B_2 f(\phi(s,x_0,u),u(s))ds\Big\|_X  \\
			&\leq Me^{\lambda t} \|x_0\|_X + h_t \|u\|_{\Uc} + c_t \|f(\phi(\cdot,x_0,u),u(\cdot))\|_{L^\infty([0,t],V)}\\
			&\leq Me^{\lambda t} \|x_0\|_X + h_t \|u\|_{\Uc} 
			+ c_t \|f(\phi(\cdot,x_0,u),u(\cdot)) - f(0, u(\cdot))\|_{L^\infty([0,t],V)}\\
			&\qquad\qquad\qquad\qquad+ c_t \|f(0,u(\cdot))\|_{L^\infty([0,t],V)}\\
			&\leq Me^{\lambda t} \|x_0\|_X + h_t \|u\|_{\Uc} 
			+ c_t L \|\phi(\cdot,x_0,u)\|_{L^\infty([0,t],X)} + c_t (\sigma(\|u\|_\Uc) + c).
\end{align*}
As $c_t \to 0$ as $t \to +0$, there is some $t_1 \in (0,t_m (x_0, u))$ such that $c_{t_1}L\leq\frac{1}{2}$. Then it holds that 
\begin{align}
\label{eq:BRS-estimate}
\sup_{t \in[0,t_1]}\|\phi(t,x_0,u)\|_X 
			\leq 2\Big(Me^{\lambda t_1} \|x_0\|_X + h_{t_1} \|u\|_{\Uc} + c_{t_1} (\sigma(\|u\|_\Uc) + c)\Big).
\end{align}
Note that $t_1$ does not depend on $x_0$ and $u$. Hence, using cocycle property and with $\phi(t_1,x_0,u)$ instead of $x_0$, we obtain a uniform bound for $\phi(\cdot,x_0,u)$ on $2t_1$, $3t_1$, and so on. Thus, $\phi(\cdot,x_0,u)$ 
is uniformly bounded on $[0,t_m (x_0, u))$, and hence can be prolonged to a larger interval by BIC property, a contradiction to the definition of $t_m (x_0, u)$. Hence, $\Sigma$ is forward complete, and the estimate 
\eqref{eq:BRS-estimate} iterated as above to larger intervals shows the BRS property.
\end{proof}

\subsection{Regularity of the flow map}

We start with a basic result describing the exponential deviation between two trajectories.
\begin{theorem}
\label{thm:Deviation-of-trajectories}
Let Assumptions~\ref{ass:Admissibility}, \ref{ass:Integrability}, \ref{Assumption1} hold.
Take $M \ge 1$, $\lambda\geq 0$  such that $\|T(t)\|\leq Me^{\lambda t}$ for all $t\geq 0$.
Pick any $x_1,x_2 \in X$, any $u\in\Uc$, and let $\phi(\cdot,x_1,u)$ and $\phi(\cdot,x_2,u)$ be defined on a certain common interval $[0,\tau]$.

Then there exists $R=R(x_1,x_2,\tau,u)>0$, such that 
\begin{align}
\label{eq:Exponential-deviation}
\|\phi(t,x_1,u) - \phi(t,x_2,u)\|_X  \leq  2M\|x_1-x_2\|_X e^{Rt},\quad t\in[0,\tau].
\end{align}
\end{theorem}

\begin{proof}
Pick any $x_1,x_2 \in X$, any $u \in \Uc$, and let $\phi_i(t):=\phi(t,x_i,u)$, $i=1,2$ be the corresponding (unique) maximal solutions of \eqref{eq:SEE+admissible} (guaranteed by Theorem~\ref{thm:Global_well-posedness}), defined on $[0,\tau]$, for a certain $\tau>0$.

Set 
\[
K:=\max\big\{\sup_{0 \leq t \leq \tau}\|\phi_1(t)\|_X,\sup_{0 \leq t \leq  \tau}\|\phi_2(t)\|_X,\|u\|_\Uc\big\}<\infty,
\]
 where $K$ is finite due to continuity of trajectories.

Due to  \eqref{eq:SEE+admissible_Integral_Form}, and using Lipschitz continuity of $f$ (see \eqref{eq:Lipschitz}), we have for any $t \in [0,\tau]$:
\begin{align*}
\|\phi_1(t) - \phi_2(t)\|_X &\leq \|T(t)\| \|x_1-x_2\|_X  \\
&\qquad+ \Big\|\int_0^t T(t-s)B_2 \Big(f(\phi_1(s),u(s)) - f(\phi_2(s),u(s))\Big)ds\Big\|_X \\
&\leq Me^{\lambda t} \|x_1-x_2\|_X  +  c_t \|f(\phi_1(\cdot),u) - f(\phi_2(\cdot),u)\|_{L^\infty([0,t],X)}\\
&\leq Me^{\lambda t} \|x_1-x_2\|_X  +  c_t L(K) \|\phi_1(\cdot)-\phi_2(\cdot)\|_{L^\infty([0,t],X)}.
\end{align*}
As $c_t \to 0$ as $t \to +0$, there is some $t_1 \in (0,\tau)$ such that $c_{t_1}L(K)\leq\frac{1}{2}$. 
Note that $t_1$ depends on $\tau$ only (as $K$ does).

 Then, taking the supremum of the previous expression over $[0,t]$, with $t<t_1$, we obtain that 
\begin{align*}
\|\phi_1(t) - \phi_2(t)\|_X 
&\leq 2Me^{\lambda t} \|x_1-x_2\|_X,\quad t\in[0,t_1].
\end{align*}
Take $k\in\N$ such that $kt_1 <\tau$ and $(k+1)t_1>\tau$. 
Then, using the cocycle property, for any $l\in\N$, $l\leq k$ and all $t\in [0,t_1]$ s.t. $lt_1+t<\tau$ we have
\begin{align*}
\|\phi_1(lt_1 + t) - \phi_2(lt_1+t)\|_X 
&\leq (2Me^{\lambda t_1})^{l} 2M e^{\lambda t} \|x_1-x_2\|_X\\
&= 2M e^{l \ln(2M) + \lambda lt_1 + \lambda t} \|x_1-x_2\|_X\\
&= 2M e^{l Rt_1 + \lambda t} \|x_1-x_2\|_X\\
&= 2M e^{R (lt_1 + t)} \|x_1-x_2\|_X,
\end{align*}
where $R>\lambda$. Finally,
\begin{align*}
\|\phi_1(s) - \phi_2(s)\|_X &\leq 2M e^{R s} \|x_1-x_2\|_X,\quad s\in[0,\tau).
\end{align*}
\end{proof}

\begin{definition}
\label{axiom:Lipschitz}
\emph{The flow of a forward complete control system $\Sigma=(X,\Uc,\phi)$, is called Lipschitz continuous on compact intervals (for uniformly bounded inputs)}, if 
for any $\tau>0$ and any $C>0$ there exists $L>0$ so that for any $x_1,x_2 \in \clo{B_C}$,
for all $u \in B_{C,\Uc}$, it holds that 
\begin{eqnarray}
\|\phi(t,x_1,u) - \phi(t,x_2,u) \|_X \leq L \|x_1-x_2\|_X,\quad t \in [0,\tau].
\label{eq:Flow_is_Lipschitz}
\end{eqnarray} 
\end{definition}

Theorem~\ref{thm:Deviation-of-trajectories} estimates the deviation between two trajectories. To have a stronger result, showing the Lipschitz continuity of the flow map $\phi$, we additionally assume the BRS property of \eqref{eq:SEE+admissible}.
\begin{theorem}
\label{thm:Lipschitz-continuity-of-flow}
Let Assumptions~\ref{ass:Admissibility}, \ref{ass:Integrability}, \ref{Assumption1}.
Let further \eqref{eq:SEE+admissible} have BRS property.
 Then \eqref{eq:SEE+admissible} has a flow that is Lipschitz continuous on compact intervals for uniformly bounded inputs.
\end{theorem}

\begin{proof}
Take any $C>0$ and pick any $x_1,x_2 \in B_C$, and any $u \in \Uc$ with
$\|u\|_{\Uc} \leq C$.
Let $\phi_i(\cdot):=\phi(\cdot,x_i,u)$, $i=1,2$ be the corresponding maximal solutions of \eqref{eq:SEE+admissible}. These solutions are global since we assume that \eqref{eq:SEE+admissible} is forward-complete.

As \eqref{eq:SEE+admissible} is BRS, the following quantity is finite for any $\tau>0$:
\[
K:=\sup_{t\in[0,\tau],\ x\in B_C,\ u\in B_{C,\Uc}}\|\phi(t,x,u)\|_X<\infty.
\]
Following the lines of the proof of \eqref{thm:Deviation-of-trajectories}, we obtain the claim.
%
%
%
\end{proof}


To obtain continuity of the flow map with respect to both states and inputs (see Definition~\ref{def:Continuous-dependence-systems}), which is important for the application of the density argument, we impose additional conditions on the nonlinearity $f$.
\begin{theorem}
\label{Continuous_dependence_Thm}
Let Assumptions~\ref{ass:Admissibility}, \ref{ass:Integrability}, \ref{Assumption1} hold.
Let further there exists $q \in \Kinf$ such that for all $C>0$ there is $L(C)>0$: 
for all $x_1,x_2 \in \clo{B_C}$ and all $v_1,v_2 \in \clo{B_{C,U}}$ it holds that
\begin{equation}
\|f(x_1,v_1)-f(x_2,v_2)\|_V \leq L(C) \big( \|x_1-x_2\|_X + q(\|v_1-v_2\|_U) \big).
\label{eq:Lipschitz_On_X_U}
\end{equation}
Further, let \eqref{eq:SEE+admissible} have the BRS property, and $B_2 $ be a bounded operator. 
Then \eqref{eq:SEE+admissible} depends continuously on initial states and inputs.
\end{theorem}

\begin{proof}
Take any $C>0$ and pick any $x_1,x_2 \in \clo{B_C}$ and any $u_1,u_2 \in \clo{B_{C,\Uc}}$, $i=1,2$. Let $\phi_i(\cdot)=\phi(\cdot,x_i,u_i)$, $i=1,2$ be the corresponding global solutions.

Due to  \eqref{eq:SEE+admissible_Integral_Form}, we have:
\begin{align*}
\|\phi_1(t) - \phi_2(t)\|_X &\leq \|T(t)\| \|x_1-x_2\|_X  +h_t \|u_1-u_2\|_{\Uc}\\
&\qquad + \|B_2\| \int_0^t{\|T(t-r)\| \big\|f(\phi_1(r),u_1(r))-f(\phi_2(r),u_2(r))\big\|_Vdr}.
\end{align*}

In view of the boundedness of reachability sets for the system \eqref{eq:SEE+admissible}, we have
\[
K:=\sup_{\|z\|_X\leq C,\ \|u\|_\Uc\leq C,\ t \in [0,\tau]}\|\phi(t,z,u)\|_X < \infty.
\]
As $\|T(t)\|\leq Me^{\lambda t}$ for some $M,\lambda\geq 0$ and all $t\geq 0$, and due to the property 
\eqref{eq:Lipschitz_On_X_U} with $L:=L(K)$ (note that $K\geq C$), we can continue above estimates to obtain
\begin{align*}
\|\phi_1(t) &- \phi_2(t)\|_X \leq Me^{\lambda t} \|x_1-x_2\|_X  + h_t \|u_1-u_2\|_{\Uc} \\
&  + \|B_2\|\int_0^t Me^{\lambda (t-r)}L(K) \Big( \|\phi_1(r)-\phi_2(r)\|_X + q\big(\|u_1(r)-u_2(r)\|_U\big)\Big) dr.
\end{align*}

Now define $z_i(r):=e^{-\lambda r}\phi_i(r)$, $i=1,2$, $r \geq 0$.
Then
\begin{align*}
\|z_1(t) &- z_2(t)\|_X 
\leq M  \|x_1-x_2\|_X  +h_t \|u_1-u_2\|_{\Uc}  + \|B_2\| M L(K)\int_0^t{ \|z_1(r)-z_2(r)\|_X dr} \\
&\qquad\qquad\qquad\qquad + \|B_2\| ML(K)\int_0^t e^{-\lambda r} q\big(\|u_1(r)-u_2(r)\|_U\big) dr\\
&\leq M  \|x_1-x_2\|_X   +h_t \|u_1-u_2\|_{\Uc} + \|B_2\| M L(K)\int_0^t{ \|z_1(r)-z_2(r)\|_X dr} \\
&\qquad\qquad\qquad\qquad + \|B_2\| ML(K)\int_0^t e^{-\lambda r}  dr q\big(\|u_1-u_2\|_\Uc\big)\\
&\leq M  \|x_1-x_2\|_X   +h_t \|u_1-u_2\|_{\Uc} + \|B_2\| M L(K)\int_0^t{ \|z_1(r)-z_2(r)\|_X dr} \\
&\qquad\qquad\qquad\qquad + \frac{\|B_2\| ML(K)}{\lambda} q\big(\|u_1-u_2\|_\Uc\big).
\end{align*}
In view of Gronwall's inequality, we obtain 
\begin{align*}
\|z_1(t) - z_2(t)\|_X  \leq  \Big( M \|x_1-x_2\|_X&  +h_t \|u_1-u_2\|_{\Uc} \\
&+ \frac{\|B_2\| ML(K)}{\lambda} q\big(\|u_1-u_2\|_\Uc\big)\Big) e^{\|B_2\| M L(K) t}
\end{align*}
or in the original variables
\begin{align*}
\|\phi_1(t) - &\phi_2(t)\|_X  \leq \Big( M \|x_1-x_2\|_X\\
  &+h_t \|u_1-u_2\|_{\Uc} + \frac{\|B_2\| ML(K)}{\lambda} q\big(\|u_1-u_2\|_\Uc\big)\Big) e^{ (\|B_2\| M L(K) + \lambda) t}.
\end{align*}
This implies that \eqref{eq:SEE+admissible} depends continuously on inputs and initial states.
\end{proof}

\subsection{Continuity at trivial equilibrium}

Without loss of generality, we restrict our analysis to fixed points of the form $(0,0) \in X \times
\Uc$. 

To describe the behavior of solutions near the equilibrium, the following notion is of importance:
\begin{definition}
\label{def:RobustEquilibrium_Undisturbed}
\index{property!CEP}
Consider a system $\Sigma=(X,\Uc,\phi)$ with equilibrium point $0\in X$.
We say that
\emph{$\phi$ is continuous at the equilibrium} if for every $\eps >0$ and for any $h>0$ there exists a $\delta =
          \delta (\eps,h)>0$, so that $[0,h] \tm B_\delta \tm B_{\delta,\Uc} \subset D_\phi$, and
\begin{eqnarray}
 t\in[0,h],\ \|x\|_X \leq \delta,\ \|u\|_{\Uc} \leq \delta \; \Rightarrow \;  \|\phi(t,x,u)\|_X \leq \eps.
\label{eq:RobEqPoint}
\end{eqnarray}
In this case, we will also say that $\Sigma$ has the \emph{CEP property}.
\end{definition}

CEP property is a \q{local version} of Lyapunov stability and is important, in particular, for the ISS superposition theorems \cite{MiW18b} and for the applications of the non-coercive ISS Lyapunov theory \cite{JMP20}.
\begin{lemma}[Continuity at equilibrium for \eqref{eq:SEE+admissible}]
\label{lem:RobustEquilibriumPoint}
Let Assumptions~\ref{ass:Admissibility}, \ref{ass:Integrability}, \ref{Assumption1} hold, and let $f(0,0)=0$.
Then $\phi$ is continuous at the equilibrium.
\end{lemma}

\begin{proof}
Consider the following auxiliary system
\begin{subequations}
\label{xdot=Ax+f_xu_saturated}
\begin{eqnarray}
\dot{x}(t) & = & Ax(t) + B_2\tilde{f}(x(t),u(t)) + Bu(t),\quad t>0, \\
x(0)  &=&  x_0,
\end{eqnarray}
\end{subequations}
where 
\[
\tilde{f}(x,u):=f\big(\sat(x),\sat_2(u)\big),\quad x\in X,\ u\in U, 
\]
and the saturation function is given for the vectors $z$ in $X$ and in $U$ respectively by
\[
\sat(z):=
 \begin{cases}
           z, &   \|z\|_X\leq 1, \\
           \frac{z}{\|z\|_X}, &   \text{otherwise},
 \end{cases}
\qquad
\sat_2(z):=
 \begin{cases}
           z, &   \|z\|_U\leq 1, \\
           \frac{z}{\|z\|_U}, &   \text{otherwise}. 
 \end{cases} 
\]
As $f$ satisfies Assumption~\ref{Assumption1}, one can show that $\tilde{f}$ is uniformly globally Lipschitz continuous. Hence, \eqref{xdot=Ax+f_xu_saturated} is forward complete and has BRS property by Proposition~\ref{prop:Unif-Glob-Lip-and-BRS}.

We denote the flow of \eqref{xdot=Ax+f_xu_saturated} by $\tilde{\phi} = \tilde{\phi}(t,x,u)$. As $f(x,u)=\tilde{f}(x,u)$ whenever $\|x\|_X\leq 1$ and $\|u\|_U\leq 1$,  it holds also 
\[
\phi(t,x,u) = \tilde{\phi}(t,x,u),
\]
provided that $\|u\|_{\Uc}\leq 1$, $\phi(\cdot,x,u)$ exists on $[0,t]$, and $\|\phi(s,x,u)\|_X\leq 1$ for all $s\in[0,t]$.

Pick any $\eps \in (0,1)$, $\tau \geq 0$, $\delta \in (0,\varepsilon)$, and choose any $x \in B_\delta$, as well as any $u\in B_{\delta,\Uc}$. It holds that
\begin{eqnarray*}
\|\tilde{\phi}(t,x,u) \|_X 
            &\leq& \|\tilde{\phi}(t,x,u) -\tilde{\phi}(t,0,u)\|_X + \|\tilde{\phi}(t,0,u)\|_X.
\end{eqnarray*}
Since \eqref{xdot=Ax+f_xu_saturated} has BRS property, by Theorem~\ref{thm:Lipschitz-continuity-of-flow}, the flow of \eqref{xdot=Ax+f_xu_saturated} is
Lipschitz continuous on compact time intervals. Hence there exists a $L(\tau,\delta)>0$ so that for $t\in[0,\tau]$
\begin{eqnarray}
\label{eq:CEP-estimate-1}
\|\tilde{\phi}(t,x,u) -\tilde{\phi}(t,0,u)\|_X \leq L(\tau,\delta)\|x\|_X \leq L(\tau,\delta)\delta.
\end{eqnarray}

Let us estimate $\|\tilde{\phi}(t,0,u)\|_X$. We have:
\begin{align*}
\|\tilde{\phi}(t,0,u)\|_X 
&\le \Big\|\int_0^t T(t-s)B_2\tilde{f}\big(\tilde{\phi}(s,0,u),u(s)\big)ds\Big\|_X + \Big\|\int_0^t T(t-s)Bu(s)ds\Big\|_X\\
&\leq c_t\esssup_{s\in[0,t]}\big\|\tilde{f}\big(\tilde{\phi}(s,0,u),u(s)\big)\big\|_X + h_t\big\|u\big\|_{L^\infty([0,t],U)}\\
&\leq c_t\esssup_{s\in[0,t]}\big\|\tilde{f}\big(\tilde{\phi}(s,0,u), u(s)\big) - \tilde{f}\big(0, u(s)\big)\big\|_X + c_t\esssup_{s\in[0,t]}\big\| \tilde{f}\big(0, u(s)\big)\big\|_X \\
&\qquad\qquad\qquad + h_t\big\|u\big\|_{L^\infty([0,t],U)}.
\end{align*}
Since $\tilde{f}(0,\cdot)$ is continuous, for any $\eps_2>0$ there exists $\delta_2<\delta$ so that $u(s) \in B_{\delta_2}$ implies that 
$\|\tilde{f}(0,u(s))-\tilde{f}(0,0)\|_X \leq \eps_2$. Since $\tilde{f}(0,0)=0$, for the above $u$ we have
$\|\tilde{f}(0,u(s))\|_X \leq \eps_2$.

As $\tilde{f}$ is uniformly globally Lipschitz, there is $L>0$ such that for the inputs satisfying $\|u\|_\Uc \leq \delta_2$ we have  
\begin{align*}
\|\tilde{\phi}(t,0,u)\|_X \leq c_tL\esssup_{s\in[0,t]}\big\|\tilde{\phi}(s,0,u)\big\|_X + c_t \varepsilon_2 + h_t\delta_2.
\end{align*}
As $c_t\to 0$ for $t\to+0$, there is $t_1>0$, such that $c_{t_1}L\leq \frac{1}{2}$.
 
Then we have that
\begin{align}
\label{eq:CEP-estimate-2}
\|\tilde{\phi}(t,0,u)\|_X \leq 2c_{t_1} \varepsilon_2 + 2h_{t_1}\delta_2,\quad t \leq t_1.
\end{align}

Combining \eqref{eq:CEP-estimate-1} with \eqref{eq:CEP-estimate-2}, we see that whenever $\|x\|_X\leq \delta_2$ and $\|u\|_\Uc \leq \delta_2$, it holds that 
\begin{eqnarray*}
\|\tilde{\phi}(t,x,u) \|_X 
             \leq L(\tau,\delta_2)\delta_2 + 2c_{t_1} \varepsilon_2 + 2h_{t_1}\delta_2,\quad t \leq t_1.
\end{eqnarray*}
Now for any $\varepsilon<1$ we can find $\delta_2<\varepsilon$, such that 
\begin{eqnarray*}
\|\tilde{\phi}(t,x,u) \|_X 
             \leq  \varepsilon,\quad t \leq t_1,\quad \|x\|_X\leq \delta_2,\quad \|u\|_\Uc \leq \delta_2.
\end{eqnarray*}
As $\tilde{\phi}(t,x,u) = \phi(t,x,u)$ whenever $\|\tilde{\phi}(t,x,u)\|_X<1$, we obtain that 
\begin{eqnarray*}
\|\phi(t,x,u) \|_X  \leq  \varepsilon,\quad t \leq t_1,\quad \|x\|_X\leq \delta_2,\quad \|u\|_\Uc \leq \delta_2.
\end{eqnarray*}
Note that $t_1$ depends on $L$ only and does not depend on $\delta_2$. 
Thus, one can find $\delta_3<\delta_2$, such that
\begin{eqnarray*}
\|\phi(t,x,u) \|_X  \leq  \delta_2,\quad t \leq t_1,\ \|x\|_X\leq \delta_3,\ \|u\|_\Uc \leq \delta_3.
\end{eqnarray*}
By cocycle property, we obtain that 
\begin{eqnarray*}
\|\phi(t,x,u) \|_X \leq  \varepsilon,\quad t \leq 2t_1,\ \|x\|_X\leq \delta_3,\ \|u\|_\Uc \leq \delta_3.
\end{eqnarray*}
Iterating this procedure, we obtain that there is some $\omega>0$, such that 
\begin{eqnarray*}
\|\phi(t,x,u) \|_X \leq \varepsilon,\quad t \in [0,\tau],\ \|x\|_X\leq \omega,\ \|u\|_\Uc \leq \omega.
\end{eqnarray*}
This shows the CEP property.
\end{proof}

\section{Boundary control systems}
\label{sec:Boundary_control_systems}

Control systems governed by partial differential equations are defined by PDEs describing the dynamics inside of the spatial domain and boundary conditions, describing the dynamics of the system at the boundary of the domain. Such systems look (at first glance) quite differently from the evolution equations in Banach spaces, studied in Section~\ref{sec:Semilinear boundary control systems}. 
This motivated the development of a theory of abstract boundary control systems that allows for a more straightforward interpretation of PDEs in the language of semigroup theory.

\subsection{Linear boundary control systems}
\label{sec:Linear_Boundary_control_systems}

Let $X$ and $U$ be Banach spaces. Consider a system
\begin{subequations}
\label{eq:BCS}
\begin{align}
\dot{x}(t) &= {\Ah} x(t), \qquad x(0) = x_0, \label{eq:BCS-1}\\
{\Rh}x(t) &= u(t),    \label{eq:BCS-2}
\end{align}
\end{subequations}
where the \emph{formal system operator} $\Ah: D( \Ah ) \subset X \to X$ is a linear operator, the control function $u$ takes values in $U$, and the \emph{boundary operator} $\Rh : D( \Rh ) \subset X \to U$ is linear and satisfies $D(\Ah) \subset D(\Rh)$.

Equations \eqref{eq:BCS} look different from the classic linear infinite-dimensional systems studied previously: 
\begin{equation}
\dot{x}(t) = Ax(t) + Bu(t),
\quad
x(0) = x_0,
\label{eq:BCS-standard-linear-systems}
\end{equation}
where $A$ is the generator of a strongly continuous semigroup, and $B$ is either a bounded or an admissible input operator.

To use for the system \eqref{eq:BCS} the theory we developed for linear systems \eqref{eq:BCS-standard-linear-systems}, we will transform \eqref{eq:BCS} into the form \eqref{eq:BCS-standard-linear-systems}. This can be done only under some additional assumptions.

\begin{definition}
\label{def:BCS}
The system \eqref{eq:BCS} is called a~\emph{linear boundary control system (linear BCS)} if the following conditions hold:
\begin{enumerate}
\item[(i)] The operator $A : D(A) \to X$ with $D(A) = D({\Ah} ) \cap \ker({\Rh})$ defined by
    \begin{equation}
    Ax = {\Ah}x \qquad \text{for} \quad x\in D(A)
        \label{eq:BSC-ass1}
    \end{equation}
    is the infinitesimal generator of a $C_0$-semigroup $(T(t))_{t\geq0}$ on $X$.
\item[(ii)] There is an operator $R \in \mathcal{L}(U,X)$ such that for all $u \in U$ we have $Ru \in D({\Ah})$,
    ${\Ah}R \in \mathcal{L}(U,X)$ and         
        \begin{equation}
    {\Rh}Ru = u, \qquad u\in U.
        \label{eq:BSC-ass2}
    \end{equation}
\end{enumerate}
The operator $R$ in this definition is sometimes called a \emph{lifting operator}. Note that $R$ is not uniquely defined by the properties in item (ii).
\end{definition}

%

Item (i) of Definition~\ref{def:BCS} shows that for $u\equiv 0$, the equations \eqref{eq:BCS} are well-posed.
In particular, as $A$ is the generator of a certain strongly continuous semigroup $T(\cdot)$, for any $x\in D(A)$, it holds that $T(t)x \in D(A)$ and thus $T(t)x \in \ker({\Rh})$ for all $t\geq 0$, which means that \eqref{eq:BCS-2} is satisfied.

Item (ii) of the definition implies, in particular, that the range of the operator ${\Rh}$ equals $U$, and thus the values of inputs are not restricted.

\subsection{Semilinear boundary control systems}
\label{sec:Semilinear-BCS}

Let $(\Ah,\Rh)$ be a linear BCS.
We consider $D(\Ah) \subset X$ as a linear space equipped with the graph norm
\[
\|\cdot\|_{D(\Ah)}:=\|\cdot\|_{X}+\big\|\Ah\cdot\big\|_{X}.
\]

Motivated by \cite{Sch20}, we consider the following class of semilinear boundary control systems.
\begin{definition}
\label{def:N-BCS}
Consider a linear BCS $(\Ah,\Rh)$.
Consider the following system 
\begin{subequations}
\label{eq:NBCS}
\begin{align}
\dot{x}(t)&= \Ah x(t)+f(x(t),w(t)), \ t>0, \label{eq:NBCS-1}\\
{\Rh}x(t) &= u(t),\ t>0,    \label{eq:NBCS-2}\\
x(0) &= x_0,    \label{eq:NBCS-3}
\end{align}
\end{subequations}
with a nonlinearity $f:X\times W\to X$, where $W$ is a Banach space. 
 
The system \eqref{eq:NBCS} we call a \emph{semilinear boundary control system (semilinear BCS)}.
\end{definition}

Following \cite{Sch20}, we define classical solutions to the semilinear BCS \eqref{eq:NBCS}.
\begin{definition}
\label{def:BCS-Classical solutions} 
Let $x_{0}\in D(\Ah)$, $\tau>0$ and $u\in C([0,\tau],U)$.  A function 
\[
x\in C([0,\tau],D(\Ah))\cap C^{1}([0,\tau],X)
\] 
is called a \emph{classical solution} to the semilinear BCS \eqref{eq:NBCS} on $[0,\tau]$ if $x(t)\in X$ for all $t>0$ and the equations \eqref{eq:NBCS} are satisfied pointwise for $t\in (0,\tau]$. 

A function $x:[0,\infty)\to X$ is called \emph{(global) classical solution} to the semilinear BCS \eqref{eq:NBCS}, if $x|_{[0,\tau]}$ is a classical solution on $[0,\tau]$ for every $\tau>0$. 

If $x\in C([0,\tau],D(\Ah))\cap C^{1}((0,\tau],X)$ and $x(t)\in X$ for all $t>0$ and the equations \eqref{eq:NBCS} are satisfied pointwise for $t\in (0,\tau]$, then we say that $x$ is a \emph{classical solution on $(0,\tau]$}.
\end{definition}

Recall the following result (see \cite[Corollary 10.1.4]{JaZ12}):
\begin{proposition}
\label{prop:Useful-formula} 
Let $T$ be a strongly continuous semigroup over a Banach space $X$ with the infinitesimal generator $A$.
Let also $f\in C^1(\R_+, X)$.
Then \linebreak $\int_0^t T(t-r)f(r)dr \in D(A)$ for all $t\geq 0$, and the following holds:
\begin{eqnarray}
A\int_0^t T(t-r)f(r)dr = \int_0^t T(t-r)\dot{f}(r)dr + T(t)f(0) - f(t).
\label{eq:Generator_of_convolution}
\end{eqnarray}
\end{proposition}

The next theorem gives a representation for the (unique) solutions of \eqref{eq:NBCS} \emph{for smooth enough inputs}.
\begin{theorem}
\label{thm:BCS-classical-solution-Representation}
Consider the boundary control system \eqref{eq:BCS} with $f \in C(X\tm W,X)$.
Let $u \in C^2([0,\tau],U)$, $w\in C([0,\tau],W)$ for some $\tau>0$, and let $x_0\in X$ be such that $x_0 - Ru(0) \in D(A)$. Assume that the classical solution of semilinear BCS $\phi(\cdot,x_0,u)$ exists on $[0,\tau]$. Then it can be represented as
\begin{subequations}
\label{eq:BCS-solution-for-smooth-inputs}
\begin{align}
\phi(t,x_0,u) &= T(t)\big(x_0-Ru(0)\big) \nonumber\\
& \qquad + \int_0^t T(t-r)\Big(f(x(r),w(r)) + {\Ah}Ru(r)-R\dot{u}(r)\Big) dr + Ru(t) \label{eq:BCS-solution-for-smooth-inputs-1}\\
&= T(t)x_0 + \int_0^t T(t-r)\Big(f(x(r),w(r))+{\Ah}Ru(r)\Big)dr\nonumber\\
& \qquad - A\int_0^t T(t-r)Ru(r) dr \label{eq:BCS-solution-for-smooth-inputs-2}\\
&= T(t)x_0 + \int_0^t T_{-1}(t-r)\Big(f\big(x(r),w(r)\big)+({\Ah}R - A_{-1}R)u(r)\Big)dr \label{eq:BCS-solution-for-smooth-inputs-3},
\end{align}
\end{subequations}
where $A_{-1}$ and $T_{-1}$ are the extensions of the infinitesimal generator $A$ and of the semigroup $T$ to the extrapolation space $X_{-1}$.
Furthermore, $A_{-1}R \in L(U,X_{-1})$ (and thus ${\Ah}R - A_{-1}R \in L(U,X_{-1})$).
\end{theorem}

The proof is similar to the proof of the linear case (see \cite[Theorem 4.4]{MiP20}, \cite[pp. 93--94]{Sch20}). However, we provide a detailed proof as we strive to be self-contained.

\begin{proof}
\textbf{(a).} 
Pick any $\tau\geq 0$, $u \in C^2([0,\tau],U)$ and any $x_0\in X$ such that $x_0 - Ru(0) \in D(A)$.
Let $x(\cdot):=\phi(\cdot,x_0,u)$ be a classical solution of the semilinear BCS.
Define 
\[
v(t):=x(t) - Ru(t).
\]
We are going to show that this function is a classical solution of the following equation
\begin{eqnarray}
\dot{v} = Av(t) + f\big(v(t) + Ru(t),w(t)\big) - R\dot{u}(t) + \Ah R u(t)
\label{eq:Related-equation}
\end{eqnarray}
on $[0,\tau)$, in the sense that $v\in C([0,\tau),X)$, $v\in C^1((0,\tau),X)$, $v(t) \in D(A)$ for $t\in(0,\tau)$, and \eqref{eq:Related-equation} holds with this $v$.

Firstly, it holds for all $t\ge 0$ that 
\[
\Rh v(t) = \Rh \phi(t,x_0,u) - u = 0.
\]
As $\im(R) \subset D(\Ah)$ and $x$ is a classical solution of semilinear BCS, then $v(t) \in D(\Ah)\cap \ker (\Rh) = D(A)$.

Furthermore,
\begin{eqnarray*}
\dot{v}(t) &=& \Ah x(t)+f(x(t),w(t))  - R\dot{u}(t) \\
&=& \Ah (v(t) + Ru(t)) + f\big(v(t) + Ru(t),w(t)\big)  - R\dot{u}(t)\\
&=& \Ah v(t) + f\big(v(t) + Ru(t),w(t)\big)  - R\dot{u}(t) + \Ah Ru(t).
\end{eqnarray*}
Since $v(t) \in D(A)$, $A = \Ah$ on $D(A)$, and $v \in C^1([0,\tau],X)$, $v$ is a classical solution of \eqref{eq:Related-equation}.

As \eqref{eq:Related-equation} is a classical semilinear evolution equation, $v$ can be represented 
(see \cite[equation (1.2) on page 183]{Paz83}) as 
\begin{eqnarray*}
v(t) &=& T(t)v(0) \\
&& + \int_0^t T(t-r)\big(f\big(v(r) + Ru(r),w(r)\big) + {\Ah}Ru(r)-R\dot{u}(r)\big) dr,
\end{eqnarray*}
and returning to the $x$-variable, we obtain \eqref{eq:BCS-solution-for-smooth-inputs-1}.

\textbf{(b).} By Proposition~\ref{prop:Useful-formula}, it holds that
\begin{align}
 \int_0^t &T(t-r)R\dot{u}(r) dr - Ru(t) =  A\int_0^t T(t-r) Ru(r) dr - T(t)Ru(0).
\label{eq:Corollary-JaZ12}
\end{align}
Substituting this into \eqref{eq:BCS-solution-for-smooth-inputs-1}, we obtain \eqref{eq:BCS-solution-for-smooth-inputs-2}.

\textbf{(c).} Let $A_{-1}$ be the extension of $A$ to the extrapolation space $X_{-1}$, and let $T_{-1}$ be the extrapolated semigroup, generated by $A_{-1}$.
Note that $R\in L(U,X)$, and $D(A_{-1})=X$. Thus, the operator $A_{-1}R$ is well-defined as a linear operator from Banach space $U$ to Banach space $X_{-1}$ with $D(A_{-1}R) = U$. As $A_{-1}$ is the generator of a strongly continuous semigroup, it is closed.
Thus, by \cite[Proposition A.9]{HMM13}, the operator $A_{-1}R$ is closed as a product of a closed and a bounded operator. By closed graph theorem (see, e.g., \cite[Theorem A.3.52]{CuZ20}), $A_{-1}R \in L(U,X_{-1})$.

The map $r \mapsto T(t-r)Ru(r)$ is Bochner integrable in the space $X$ and thus also in $X_{-1}$ (even for any $u \in L^{1}_{\loc}([0,\tau],U)$, see
\cite[Proposition 1.3.4]{ABH11}).
Furthermore, $T(t-r)Ru(r) \in X = D(A_{-1})$ for all $r\in[0,t]$.

Recall that $A_{-1}T_{-1}(s)=T_{-1}(s)A_{-1}$ for all $s\in\R_+$ on $D(A_{-1})$, see, e.g., \cite[Theorem 5.2.2]{JaZ12}.
Consider the map 
\[
w:r \mapsto A_{-1}T_{-1}(t-r)Ru(r) = T_{-1}(t-r)A_{-1}Ru(r).
\]
Since $A_{-1}R \in L(U,X_{-1})$ and $T_{-1}$ is a strongly continuous semigroup on $X_{-1}$, the function $w$ is 
Bochner integrable on $X_{-1}$, by \cite[Proposition 1.3.4]{ABH11}.
Hence by Hille's Lemma (see, e.g., \cite[Theorem 3.7.12]{HiP00}, \cite[Proposition 1.1.7]{ABH11}), we obtain for all $u \in C^2([0,\tau],U)$ that
\begin{align*}
A\int_0^t& T(t-r)Ru(r) dr = A_{-1}\int_0^t T_{-1}(t-r)Ru(r) dr \\
&= \int_0^t A_{-1} T_{-1}(t-r)Ru(r) dr = \int_0^t  T_{-1}(t-r)A_{-1}Ru(r) dr.
\end{align*}
From this, the formula \eqref{eq:BCS-solution-for-smooth-inputs-3} follows.
\end{proof}

An advantage of the representation formula \eqref{eq:BCS-solution-for-smooth-inputs-1} is in the boundedness of the operators $R$ and $\Ah R$ involved in the expression. Its disadvantage is that the derivative of $u$ is employed.
Still, the expression in the right-hand side of \eqref{eq:BCS-solution-for-smooth-inputs-1} makes sense for any $x \in X$ 
and for any $u \in H^1([0,\tau],U)$, and can be called a mild solution of BCS \eqref{eq:BCS}, as is done, e.g., in \cite[p. 146]{JaZ12}.

The formula \eqref{eq:BCS-solution-for-smooth-inputs-3} does not involve any derivatives of inputs, and again is given in terms of a bounded operator ${\Ah}R - A_{-1}R \in L(U,X_{-1})$.
Moreover, if we consider the expression in the right-hand side of \eqref{eq:BCS-solution-for-smooth-inputs-3} in the extrapolation spaces $X_{-1}$, then it makes sense for all $x \in X$ and all $u\in L^{1}_{\loc}(\R_+,U)$, and constitutes a mild solution of 
\begin{eqnarray}
\dot{x}(t) = Ax(t) + f(x(t),w(t)) + Bu(t),
\label{eq:semilinear equation with an admissible operator}
\end{eqnarray}
with
\begin{eqnarray}
B:={\Ah}R-A_{-1}R.
\label{eq:BCS-Input-Operator}
\end{eqnarray}
This motivates us to define the mild solutions of semilinear BCS by means of the formula \eqref{eq:BCS-solution-for-smooth-inputs-3}, as was proposed in \cite{Sch20}.

\begin{remark}
\label{rem:Uniqueness-of-B-boundary-control-systems} 
The operator $B$ is uniquely defined by the boundary control system and does not depend on the choice of the lifting operator $R$, see \cite[Proposition 2.8]{Sch20}.
\end{remark}

\begin{definition}
\label{def:Mild-solution-for-semilin-Boundary-Control-System}
Let $(\Ah,\Rh ,f)$ be a semilinear boundary control system with corresponding $A,R$. Let $x_{0}\in X$,  $\tau>0$, $w\in L_{\loc}^{1}([0,\tau],W)$, and $u\in L_{\loc}^{1}([0,\tau],U)$. A continuous function $x:[0,\tau]\to X$ is called \emph{mild solution} to the semilinear BCS \eqref{eq:NBCS} on $[0,\tau]$ if $x(t)\in X$ for all $t>0$ and $x$ solves 
\begin{align*}
x(t)= T(t)x_{0}+\int_{0}^{t}T_{-1}(t-s)\big(f(x(s),w(s))+Bu(s)\big)ds,
\end{align*}
for all $t\in[0,\tau]$ and where $B=\Ah B_{0}-A_{-1}B_{0}$. A function $x:\R_+ \to X$ is called a \emph{global mild solution} if $x|_{[0,\tau]}$ is a mild solution on $[0,\tau]$ for all $\tau>0$.
\end{definition}

In other words, $x$ is a mild solution of a semilinear BCS \eqref{def:N-BCS}, if $x$ is a mild solution of 
\eqref{eq:semilinear equation with an admissible operator} with $B={\Ah}R - A_{-1}R$.

Thus, semilinear boundary control systems are a special case of semilinear evolution equations studied in Section~\ref{sec:Semilinear boundary control systems}, and we can use our well-posedness and stability theory for semilinear evolution equations to analyze semilinear BCS.

In particular, let us state a criterion for the ISS of linear boundary control systems:
\begin{corollary}
\label{cor:ISS-for-BCS} 
Consider the boundary control system \eqref{eq:BCS} and assume that the corresponding input operator $B$ defined by 
\eqref{eq:BCS-Input-Operator} is $q$-admissible for some $q\in [1,+\infty)$. Assume that $A$ generates an exponentially stable semigroup.
Then for any $p\in[q,+\infty]$
the system $\Sigma:=(X,\Uc,\phi)$ with $\Uc:=L^p(\R_+,U)$ is an ISS control system (with respect to the norm in $\Uc$).
\end{corollary}

\begin{proof}
As $B$ defined by \eqref{eq:BCS-Input-Operator} is $q$-admissible for some $q\in [1,+\infty)$, then it is $p$-admissible for any $p\in[q,+\infty]$ and furthermore, the map $\phi$ is continuous w.r.t.\ time. 
By Proposition~\ref{prop:q-admissibility-implies-continuity}, for all $p\in[q,+\infty]$, the system $\Sigma:=(X,\Uc,\phi)$ with $\Uc:=L^p([0,\infty),U)$, is a forward-complete control system in the sense of Definition~\ref{Steurungssystem}. 
By Theorem~\ref{thm:ISS-Criterion-lin-sys-with-unbounded-operators}, the system $\Sigma$ is ISS.
\end{proof}

\begin{remark}
\label{rem:Weak solutions of boundary control systems}
In this chapter, we consider mild solutions of boundary control systems. 
In some papers, the concepts of weak solutions and strong solutions are used. 
For a detailed discussion of the relationship between all these solution concepts, we refer to \cite[Propositions 2.9, 2.11, Remark 2.10]{Sch20}. 
\end{remark}

\begin{remark}
\label{rem:Checking-adimissibility-of-B} 
The computation of the input operator $B$ using the formula \eqref{eq:BCS-Input-Operator} may be awkward in practice. Other methods for the computation of $B$ can be used, see \cite[Section 10.1]{TuW09}, \cite{EmT00} and \cite[Proposition 2.9]{Sch20}.
Furthermore, in some situations, the admissibility of $B$ and ISS of \eqref{eq:BCS} can be obtained without computation of the operator $B$.
\end{remark}

\section{Semilinear analytic systems}
\label{sec:Semilinear analytic boundary control systems}

\subsection{Recap on fractional powers of operators}
\label{sec:Recap on fractional powers of operators}

In this section, we make a short recap of the fractional powers of sectorial operators, which we need in Section~\ref{sec:Semilinear analytic boundary control systems}. 
A somewhat more detailed account of the properties of fractional powers of sectorial operators can be found, e.g., in \cite[Section 1.4]{Hen81}
and for a detailed treatment, we refer to \cite{Haa06}.

\begin{definition}
\label{def:Gamma function} 
For $z \in\C$ with $\re z >0$ the \emph{gamma function} is defined as
\begin{eqnarray}
\Gamma(z) = \int_0^\infty t^{z-1}e^{-t}dt.
\label{eq:Gamma-function}
\end{eqnarray}
\end{definition}

\begin{definition}
\label{def:Negative-Fractional-powers-operator}
Let $A$ be the generator of an exponentially stable analytic semigroup over a Banach space $X$.
For $\alpha>0$ define
\begin{eqnarray}
(-A)^{-\alpha}:=\frac{1}{\Gamma(\alpha)}\int_0^\infty t^{\alpha-1}T(t)dt.
\label{eq:Fractional-powers-operator}
\end{eqnarray}
\end{definition}
As the semigroup $T$ is analytic, 
the map $t\to T(t)$ is differentiable on $(0,+\infty)$ (see, e.g., \cite[Theorem 4.6, p. 101]{EnN00}), and hence continuous on this interval. As $T$ is also exponentially stable, the integral 
in \eqref{eq:Fractional-powers-operator} converges in the uniform operator topology. 
The following well-known property holds:
\begin{proposition}
\label{prop:Product-formula-fractional-powers}
Let $A$ be the generator of an exponentially stable analytic semigroup over a Banach space $X$. 
For any $\alpha>0$ the operator $(-A)^{-\alpha}$ belongs to $L(X)$, is injective and satisfies
\begin{eqnarray}
(-A)^{-\alpha}(-A)^{-\beta}=(-A)^{-(\alpha+\beta)},\quad \alpha,\beta>0.
\label{eq:Product-formula-fractional-powers}
\end{eqnarray}
\end{proposition}

%

\begin{definition}
\label{def:Positive-Fractional-powers-operator}
Let $A$ be the generator of an exponentially stable analytic semigroup over a Banach space $X$.
For $\alpha>0$ define $(-A)^{\alpha}$ as the inverse of $(-A)^{-\alpha}$ with 
$D((-A)^{\alpha}) = \im((-A)^{-\alpha})$.

By definition, we set $A^0:=I$.
\end{definition}

Let us collect several basic properties of the fractional powers; see \cite[Theorem 6.8, p. 72]{Paz83} for the first three items, and \cite[Exercises, p. 26]{Hen81} for the last one.
\begin{proposition}
\label{prop:Properties-of-positive-fractional-powers} 
Let $A$ be the generator of an exponentially stable analytic semigroup over a Banach space $X$.

\begin{enumerate}[label=(\roman*)]
	\item For $\alpha>0$ the operators $(-A)^\alpha$ are closed and densely defined.
	\item Whenever $\alpha > \beta>0$, we have $D((-A)^\alpha) \subset D((-A)^\beta)$.
	\item For all $\alpha,\beta\in\R$ it holds that 
	\begin{eqnarray}
	(-A)^{\alpha+\beta}x = (-A)^{\alpha}(-A)^{\beta}x,
	\label{eq:produce-rule-fractional}
	\end{eqnarray}
	for every $x \in D((-A)^\gamma)$, where $\gamma=\max\{\alpha,\beta,\alpha+\beta\}$.
	\item $(-A)^\alpha T(t) = T(t)(-A)^\alpha$ on $D((-A)^\alpha)$,\quad  $t\ge 0$.
\end{enumerate}
\end{proposition}
%
%
%
%
%

 Using fractional powers of operators, we can define the spaces that will be important to deal with
nonlinear equations governed by analytic semigroups.
\begin{definition}
\label{def:X-power-alpha}
Let $A$ be the generator of an analytic semigroup over a Banach space $X$.
Let $\alpha\ge 0$. Pick any $\omega>\omega_0(T)$ and define the space $X_\alpha$ and its norm as 
\begin{eqnarray}
X_\alpha:=D((\omega I-A)^\alpha),\quad \|x\|_{X_\alpha}:=\|(\omega I-A)^\alpha x\|_X,\quad x \in X_\alpha.
\label{eq:X-power-alpha}
\end{eqnarray}
\end{definition}
In particular, by definition $X_0 = X$, and $X_1 = D(A)$ (with norms as  in \eqref{eq:X-power-alpha}).

\begin{remark}
\label{rem:well-definiteness-X-power-alpha} 
The choice of different $\omega > \omega_0(T)$ induces the same linear space $X_\alpha$ endowed with an equivalent norm.
\end{remark}

The spaces $X_\alpha$ have a good structure:
\begin{proposition}
\label{prop:X-power-alpha-are-Banach-spaces} 
Let $A$ be the generator of an exponentially stable analytic semigroup over a Banach space $X$.
For any $\alpha\ge 0$, the space $X_\alpha$ is a Banach space. 
Furthermore, for all $\alpha\ge\beta\ge0$, the space $X_\alpha$ is a dense subspace of $X_\beta$, with continuous inclusion.
\end{proposition}

Similarly to the space $X_{-1}$, we introduce
\begin{definition}
\label{def:X_-alpha}
For $\alpha>0$ and $A$ as above, define the spaces $X_{-\alpha}$ as the completion of $X$ with respect to the norm $x \mapsto \|(\omega I-A)^{-\alpha}x\|_X$.
Then we have for any $0<\beta<\alpha<1$ the following chain of continuous inclusions:
\[
X_{-1} \supset X_{-\alpha} \supset X_{-\beta} \supset X = X_0 \supset X_{\beta} \supset X_{\alpha} \supset X_1.
\]
\end{definition}

The following property holds:
\begin{proposition}
\label{prop:Analytic-semigroup-Bound-on-A^alphaT(t)} 
Let $T$ be an analytic semigroup on a Banach space $X$ with a growth bound $\omega_0(T)$ and the generator $A$.
Then for each $\omega,\kappa>\omega_0(T)$ and each $\alpha \in[0,1)$ we have $\im(T(t)) \subset X_\alpha$, and there is $C_\alpha>0$ such that 
\begin{eqnarray}
\|(\omega I-A)^\alpha T(t)\|  \leq \frac{C_\alpha}{t^\alpha}e^{\kappa t},\quad t>0.
\label{eq:Bound-on-A^alphaT(t)}
\end{eqnarray}
Furthermore, the map $t\mapsto (\omega I-A)^\alpha T(t)$ is continuous on $(0,+\infty)$ in the uniform operator topology.
\end{proposition}

\subsection{Semilinear analytic systems and their mild solutions}
\label{sec:Semilinear analytic systems and their mild solutions}

Consider again the system \eqref{eq:SEE+admissible} with $B_2=\id$ that we restate next:
\begin{subequations}
\label{eq:SEE+admissible-analytic} 
\begin{eqnarray}
\dot{x}(t) & = & Ax(t) + f(x(t),u(t)) + Bu(t),\quad t>0,  \label{eq:SEE+admissible-analytic-1}\\
x(0)  &=&  x_0, \label{eq:SEE+admissible-analytic-2}
\end{eqnarray}
\end{subequations}
In Section~\ref{sec:Semilinear boundary control systems}, we have assumed that $f$ is a well-defined map from $X \tm U$ to $X$. 
Although it sounds natural, it is, in fact, a quite restrictive assumption, as already basic nonlinearities, such as polynomial maps, do not satisfy it. 
Indeed, if $f(x) = x^2$, where $x \in X:=L^2(0,1)$, then $f$ maps $X$ to the space $L^1(0,1)$.
However, as $A$ generates an analytic semigroup, the requirements on $f$ can be considerably relaxed. 
Namely, we assume in this section that there is $\alpha \in [0,1]$ such that  $f$ is a well-defined map from $X_\alpha \tm U $ to $X$, where the space $X_\alpha$ is defined in Definition~\ref{eq:X-power-alpha}.

Next, we define mild solutions of \eqref{eq:SEE+admissible-analytic}. Note that the nonlinearity $f$ is defined on $X_\alpha \tm U$, and thus we must require that the mild solution lies in $X_\alpha$ for all positive times. We cannot expect such a nice behavior for general semigroups, but thanks to the smoothing effect of analytic semigroups, this is what we can expect in the analytic case.
\begin{definition}
\label{def:Mild-solution-analytic}
\index{solution!mild}
Let $\tau>0$ be given. 
A function $x \in C([0,\tau], X)$ is called a \emph{mild solution of \eqref{eq:SEE+admissible-analytic} on $[0,\tau]$} corresponding to certain $x_0\in X$ and $u \in L^1_{\loc}(\R_+,U)$, if $x(s) \in X_\alpha$ for $s\in(0,\tau]$, and $x$ solves the integral equation
\begin{align}
\label{eq:SEE+admissible_Integral_Form-analytic}
x(t)=T(t) x_0 + \int_0^t T(t-s) f\big(x(s),u(s)\big)ds + \int_0^t T_{-1}(t-s) Bu(s)ds. 
\end{align}

We say that $x:\R_+\to X$ is a \emph{mild solution of \eqref{eq:SEE+admissible-analytic} on $\R_+$} corresponding to 
certain $x_0\in X$ and $u \in L^1_{\loc}(\R_+,U)$, if it is a mild solution of \eqref{eq:SEE+admissible-analytic} (with $x_0, u$) on $[0,\tau]$ for all $\tau>0$.
\end{definition}

\begin{remark}
Note that if $\alpha = 0$, then $X_\alpha = X_0 = X$, and the concept of a mild solution introduced for general and analytic semigroups coincide.
\end{remark}

\begin{ass}
\label{ass:Regularity-f-analytic} 
Let the following hold:
\begin{enumerate}[label=(\roman*)]
	\item $\alpha \in (0,1)$.
	\item $B\in L(U,X_{-1+\alpha+\varepsilon})$ for sufficiently small $\varepsilon>0$.
	
	\item $f \in C(X_\alpha \tm U, X)$, and $f$ is Lipschitz continuous in the first argument in the following sense: for each $r>0$, there is $L=L(r)>0$ such that for each $x_1,x_2 \in B_{r,X_\alpha}$ and all $u \in B_{r,U}$ it holds that 
\begin{eqnarray}
\|f(x_1,u) - f(x_2,u)\|_X \leq L\|x_1-x_2\|_{X_\alpha}.
\label{eq:Lipschitz-in-x-power-alpha}
\end{eqnarray}
	\item For all $u\in L^\infty(\R_+,U)$ and any $x\in C(\R_+,X)$ with $x((0,+\infty)) \subset X_\alpha$, the map $s\mapsto f\big(x(s),u(s)\big)$ is in $L^w_{\loc}(\R_+,X)$ with $w>\frac{1}{1-\alpha}$.
	\item There is $\sigma\in\Kinf$ such that 
\[
\|f(0,u)\|_X \leq \sigma(\|u\|_U) + c,\quad u\in U.
\]
\end{enumerate}
\end{ass}

%
%

\subsection{Preliminaries for analytic semigroups and admissibility}
\label{sec:Preliminaries for analytic semigroups}

We denote by $\omega_0(T)$ the growth bound of a semigroup $T$. 
We start with the following \q{analytic} version of \cite[Proposition 1.3.4]{ABH11}:
\begin{proposition}
\label{prop:Auxiliary-map-analytic} 
Let $T$ be an analytic semigroup, $\alpha\in[0,1)$, and $g \in L^w_{\loc}(\R_+,X)$ with $w>\frac{1}{1-\alpha}$. Pick any $\omega>\omega_0(T)$. Then the map
\begin{eqnarray}
\xi:t \mapsto \int_0^t (\omega I-A)^\alpha T(t-s)g(s)ds
\label{eq:Continuity-mild-solutions-analytic-extension}
\end{eqnarray}
is well-defined and continuous on $\R_+$.

Furthermore, for any $\kappa>\omega_0(T)$ there is $R=R(\kappa,\alpha)$ such that for any $g \in L^\infty_{\loc}(\R_+,X)$ the following holds:
\begin{align}
\int_0^t \big\|(\omega I-A)^\alpha  T(t-s) g(s)\big\|_X ds  \leq R t^{1-\alpha} e^{\kappa t} \| g\|_{L^\infty([0,t],X)}.
\label{eq:Convolution-analytic-for-Linfty-norm}
\end{align}
\end{proposition}

\begin{proof}
We divide the proof into four parts. 

\textbf{(i). Well-definiteness.} Since $T$ is an analytic semigroup, $T(t)$ maps $X$ to $D(A)$ for any $t>0$. As $X_\alpha \subset D(A)$ for all $\alpha\in[0,1]$, 
the integrand in \eqref{eq:Continuity-mild-solutions-analytic-extension} is in $X$ for a.e. $s\in[0,t)$.
Let us show Bochner integrability of $X$-valued map $s \mapsto (\omega I-A)^\alpha T(t-s)g(s)$ on $[0,t]$.

As $g \in L^{1}_{\loc}(\R_+,X)$, by the criterion of Bochner integrability, $g$ is strongly measurable and $\int_I \|g(s)\|_X ds < \infty$ for any bounded interval $I \subset\R_+$.

Denote by $\chi_{\Omega}$ the characteristic function of the set $\Omega\subset \R_+$. Recall that the map $t\mapsto (\omega I-A)^\alpha T(t)$ is continuous outside of $t=0$ in view of Proposition~\ref{prop:Analytic-semigroup-Bound-on-A^alphaT(t)}. 

If $g(s) = \chi_{\Omega}(s)x$ for some measurable $\Omega \subset \R_+$ and $x \in X$, then the function
\[
s\mapsto (\omega I-A)^\alpha T(t-s)g(s) = (\omega I-A)^\alpha T(t-s)\chi_{\Omega}(s)x
\]
is measurable as a product of a measurable scalar function and a continuous (and thus measurable) vector-valued function.
By linearity, $s\mapsto (\omega I-A)^\alpha T(t-s)g(s)$ is strongly measurable if $g$ is a simple function (see \cite[Section 1.1]{ABH11} for definitions). 

As $g$ is strongly measurable, there is a sequence of simple functions $(g_n)_{n\in\N}$, converging pointwise to $g$ almost everywhere. 
Consider a sequence 
\begin{eqnarray}
\big(s\mapsto (\omega I-A)^\alpha T(t-s)g_n(s)\big)_{n\in\N}
\label{eq:SimpleFun-Seq}
\end{eqnarray}
and take any $s \in [0,t)$ such that $g_n(s) \to g(s)$ as $n\to\infty$.
We have that 
\begin{align*}
\big\|(\omega I-A)^\alpha T(t-s)&g_n(s)- (\omega I-A)^\alpha T(t-s)g(s)\big\|_X \\
&\leq \|(\omega I-A)^\alpha T(t-s)\| \|g_n(s)-g(s)\|_X \to 0,\quad n\to\infty.
\end{align*}
Hence a sequence of strongly measurable functions \eqref{eq:SimpleFun-Seq} converges a.e. to $s\mapsto (\omega I-A)^\alpha T(t-s)g(s)$,
and thus $s\mapsto (\omega I-A)^\alpha T(t-s)g(s)$ is strongly measurable by \cite[Corollary 1.1.2]{ABH11}.

Furthermore, for any $t >0$, using Proposition~\ref{prop:Analytic-semigroup-Bound-on-A^alphaT(t)}, we have that for any $\kappa>\omega_0(T)$ there is $M>0$ such that
\begin{align}
\int_0^t \big\|(\omega I-A)^\alpha  T(t-s) g(s)\big\|_X ds  
&\leq \int_0^t \|(\omega I-A)^\alpha T(t-s)\| \| g(s)\|_X ds \nonumber \\
&\leq M \int_0^t \frac{C_\alpha}{(t-s)^\alpha}e^{\kappa (t-s)} \| g(s)\|_X ds \nonumber\\
&\leq M C_\alpha e^{\kappa t}\int_0^t \frac{1}{(t-s)^\alpha} \| g(s)\|_X ds.
\label{eq:Convolution-boundedness-for-Linfty-norm}
\end{align}

Using H\"older's inequality with $w>\frac{1}{1-\alpha}$ as in the assumptions of the proposition, we obtain
\begin{align}
\int_0^t \|(\omega I-A)^\alpha  T(t-s) g(s)\|_X ds  
&\leq M C_\alpha e^{\kappa t} \Big(\int_0^t \Big(\frac{1}{(t-s)^\alpha}\Big)^p ds\Big)^{\frac{1}{p}}  \Big( \int_0^t \| g(s)\|^w_X ds\Big)^{\frac{1}{w}}\nonumber\\
&\leq \frac{M C_\alpha}{(1-\alpha p)^{1/p}} e^{\kappa t} t^{\frac{1-\alpha p}{p}}  \Big( \int_0^t \| g(s)\|^w_X ds\Big)^{\frac{1}{w}},
\label{eq:Convolution-boundedness}
\end{align}
where $\frac{1}{p} + \frac{1}{w} = 1$, and thus $p$ satisfies $p<\frac{1}{\alpha}$.

Finally, by \cite[Theorem 1.1.4]{ABH11}, 
the map $s \mapsto (\omega I-A)^\alpha  T(t-s)g(s)$ is Bochner integrable on each $[0,t]\subset \R_+$.

\textbf{(ii). Right-continuity.} We would like to show that $\xi$ defined by \eqref{eq:Continuity-mild-solutions-analytic-extension}, is right-continuous on $\R_+$.
To this end, we consider for $h>0$ and $t \ge0$ the expression
\begin{align*}
\xi(t+h) - \xi(t)&= \int_0^{t+h}(\omega I-A)^\alpha{T(t+h-s)g(s)ds} - \int_0^t(\omega I-A)^\alpha{T(t-s)g(s)ds} \\
						&=\int_0^{t}(\omega I-A)^\alpha{\Big(T(t+h-s)-T(t-s)\Big)g(s)ds} \\
						&\qquad\qquad\qquad + \int_t^{t+h}(\omega I-A)^\alpha{T(t+h-s)g(s)ds} \\
						&=(T(h)-I)\int_0^{t}(\omega I-A)^\alpha{T(t-s)g(s)ds} \\
						&\qquad\qquad\qquad+ \int_t^{t+h}{(\omega I-A)^\alpha T(t+h-s)g(s)ds}.
\end{align*}
As shown in the first part of the proof, the integral $\int_0^{t}{(\omega I-A)^\alpha T(t-s)g(s)ds}$ converges in $X$.
Since $T$ is strongly continuous, we have
\[
(T(h)-I)\int_0^{t}{(\omega I-A)^\alpha T(t-s)g(s)ds} \to 0, \quad \text{ whenever } h \to +0.
\]
Furthermore, we estimate for a certain $K$ that does not depend on $h$:
\begin{align*}
\Big\|\int_t^{t+h}{(\omega I-A)^\alpha T(t+h-s)g(s)ds}&\Big\|_X \\
	&\leq \int_t^{t+h} \big\|(\omega I-A)^\alpha T(t+h-s)\big\| \|g(s)\|_X ds\\
	&\leq K\int_t^{t+h} \frac{1}{(t+h-s)^\alpha}  \|g(s)\|_X ds
\end{align*}
and arguing similarly to \eqref{eq:Convolution-boundedness}, and using that $g \in L^w_{\loc}(\R_+,X)$ with $w>\frac{1}{1-\alpha}$,
the last term converges to $0$ as $h \to +0$.
Hence, $\xi(t+h) - \xi(t) \to 0$, as $h \to +0$.

\textbf{(iii). Left-continuity.}
Now let $t>0$ and $h>0$. Consider
\begin{align*}
\xi(t) - \xi(t-h)&= \int_0^{t}{(\omega I-A)^\alpha T(t-s)g(s)ds} - \int_0^{t-h}{(\omega I-A)^\alpha T(t-h-s)g(s)ds} \\
&= \int_0^{t-h}(\omega I-A)^\alpha\big(T(t-s)-T(t-h-s)\big)g(s)ds \\
&\qquad\qquad\qquad + \int_{t-h}^{t}(\omega I-A)^\alpha T(t-s)g(s)ds.
\end{align*}
Then
\begin{align*}
\|\xi(t) &- \xi(t-h)\|_X \le  \Big\|\int_0^{t-h}(\omega I-A)^\alpha\big(T(t-s)-T(t-h-s)\big)g(s)ds\Big\|_X  \\
&\qquad\qquad\qquad\qquad\qquad\qquad +  \Big\|\int_{t-h}^{t}(\omega I-A)^\alpha T(t-s)g(s)ds \Big\|_X \\
  &\le  \int_0^{t}\big\|(\omega I-A)^\alpha \big(T(t-s)-T(t-h-s)\big)g(s)\big\|_X \chi_{[0,t-h]}(s)ds \\
& \qquad\qquad\qquad\qquad\qquad\qquad+ \Big\| \int_{t-h}^{t}(\omega I-A)^\alpha T(t-s)g(s)ds \Big\|_X.
\end{align*}
It holds for all $s \in [0,t-h)$ that
\begin{align*}
\big\|(\omega I-A)^\alpha&\big(T(t-s)-T(t-h-s)\big)g(s)\big\|_X \chi_{[0,t-h]}(s) \\
&\leq \big(\big\|(\omega I-A)^\alpha T(t-s)\big\| + \big\|(\omega I-A)^\alpha T(t-h-s)\big\|\big) \|g(s)\|_X\\
& \leq  K \Big(\frac{1}{(t-s)^\alpha} + \frac{1}{(t-h-s)^\alpha} \Big) \|g(s)\|_X,
\end{align*}
where $K$ depends on $t$, but does not depend on $s,h$.
As in the above arguments, the function on the right-hand side is integrable, and according to the (scalar) Lebesgue theorem on dominated convergence, 
it holds that
\begin{align*}
\lim_{h \to +0}\int_0^{t-h}\big\|&(\omega I-A)^\alpha\big(T(t-s)-T(t-h-s)\big)g(s)\big\|_X ds\\
&=  \int_0^{t} \lim_{h \to +0} \big\|(\omega I-A)^\alpha\big(T(t-s)-T(t-h-s)\big)g(s)\big\|_X ds = 0.
\end{align*}
Finally,
\[
\Big\| \int_{t-h}^{t}(\omega I-A)^\alpha T(t-s)g(s)ds \Big\|_X \to 0, \quad \mbox{ as } h \to +0,
\]
as in the proof of the right-continuity.
This shows that $\|x(t+h)-x(t)\|_X \to 0$ as  $h \to +0$, and hence $x \in C([0,\tau],X)$.

\textbf{(iv). The estimate \eqref{eq:Convolution-analytic-for-Linfty-norm}.} For the last claim, we take $g\in L^\infty_{\loc}(\R_+,X)$ and continue the estimates in \eqref{eq:Convolution-boundedness-for-Linfty-norm} as follows:
\begin{align*}
\int_0^t \big\|(\omega I-A)^\alpha  T(t-s) g(s)\big\|_X ds  
&\leq M C_\alpha e^{\kappa t} \int_0^t \frac{1}{(t-s)^\alpha} ds  \| g\|_{L^\infty([0,t],X)},
\end{align*}
and we obtain \eqref{eq:Convolution-analytic-for-Linfty-norm} with $R=\frac{MC_\alpha}{1-\alpha}$.
\end{proof}

Next, we formulate a sufficient condition for the zero-class admissibility of input operators for analytic systems. 
Part (i) of the following proposition is contained in \cite[Proposition 2.13]{Sch20}. As we would like to show (ii), we give, however, the proof. 

\begin{proposition}
\label{prop:ISS-analytic-systems} 
Assume that $A$ generates an analytic semigroup $T$ and $B\in L(U,X_{-1+\alpha})$ for some $\alpha \in (0,1)$. 
Then:

\begin{enumerate}[label = (\roman*)]
	\item $B$ is zero-class $q$-admissible for $q\in(\frac{1}{\alpha},+\infty]$.
	\item For any $\omega>\omega_0(T)$, any $d \in[0,\alpha)$, and any $\kappa>\omega_0(T)$ there is $R>0$ such that for any $g \in L^\infty_{\loc}(\R_+,X)$ the map
	\begin{align}
	\label{eq:Convolution-map-analytic}
	t\mapsto \int_0^t (\omega I-A)^d  T(t-s) Bg(s) ds
	\end{align}
is continuous in $X$-norm, and the following holds:
\begin{align}
\int_0^t \big\|(\omega I-A)^d  T(t-s) Bg(s)\big\|_X ds  \leq R t^{\alpha-d} e^{\kappa t} \| g\|_{L^\infty([0,t],X)}.
\label{eq:Generalized-admissibility-estimate-analytic}
\end{align}
\end{enumerate}
\end{proposition}

\begin{proof}
\textbf{(i).} Take any $\omega>\omega_0(T)$, and 
consider the corresponding norm on $X_{-1+\alpha}$ as in Definition~\ref{def:X_-alpha}. 
We have
\begin{align*}
\|B\|_{L(U,X_{-1+\alpha})} 
&= \sup_{u\in U:\|u\|_U=1}\|Bu\|_{X_{-1+\alpha}}\\
&= \sup_{u\in U:\|u\|_U=1}\|(\omega I-A)^{-1+\alpha}Bu\|_{X} = \|(\omega I-A)^{-1+\alpha}B\|_{L(U,X)}.
\end{align*}
Thus, the condition $B\in L(U,X_{-1+\alpha})$ is equivalent to $(\omega I-A)^{-1+\alpha}B \in L(U,X)$.
With this in mind, we have
\begin{eqnarray}
T_{-1}(t)B &=& T_{-1}(t)(\omega I-A)^{1-\alpha}(\omega I-A)^{-1+\alpha}B.
\label{eq:Analytic-systems-1}
\end{eqnarray}
Due to \cite[Theorem 2.6.13, p. 74]{Paz83}, on $X_{1-\alpha} = D((\omega I-A)^{1-\alpha})$ it holds that
\[
T_{-1}(t)(\omega I-A)^{1-\alpha} = (\omega I-A)^{1-\alpha}T_{-1}(t).
\]

Now take any $f \in L^q_{\loc}(\R_+,U)$ with $q>\frac{1}{\alpha}$. Representing 
\begin{align*}
\int_0^t T_{-1}(t-s) Bf(s) ds 
= \int_0^t (\omega I-A)^{1-\alpha} T_{-1}(t-s) (\omega I-A)^{\alpha-1} Bf(s) ds
\end{align*}
and applying Proposition~\ref{prop:Auxiliary-map-analytic} and in particular the estimate \eqref{eq:Convolution-boundedness} with $1-\alpha$ instead of $\alpha$, $w:=q$, and with $g:=(\omega I-A)^{-1+\alpha}B f$ we see that $B$ is zero-class $q$-admissible for $q \in (\frac{1}{\alpha},+\infty)$.

\textbf{(ii). }
Similarly, for $d<\alpha$ and $g \in L^\infty_{\loc}(\R_+,U)$, consider the map
\begin{eqnarray*}
s \mapsto (\omega I-A)^d T_{-1}(t-s)Bg(s) &=& (\omega I-A)^{1-\alpha+d}T_{-1}(t-s)(\omega I-A)^{-1+\alpha}Bg(s).
\end{eqnarray*}
By Proposition~\ref{prop:Auxiliary-map-analytic}, this map is Bochner integrable and in view of \eqref{eq:Convolution-analytic-for-Linfty-norm} with $1-\alpha+d$ instead of $\alpha$ and $(\omega I-A)^{-1+\alpha}Bg$ instead of $g$, we see that the map 
	\eqref{eq:Convolution-map-analytic} is continuous and \eqref{eq:Generalized-admissibility-estimate-analytic} holds. 
Setting $d:=0$ in \eqref{eq:Generalized-admissibility-estimate-analytic}, we obtain also zero-class $\infty$-admissibility of $B$.
\end{proof}

\subsection{Local existence and uniqueness}

By Proposition~\ref{prop:ISS-analytic-systems}, the condition $B\in L(U,X_{-1+\alpha+\varepsilon})$ with $\alpha,\varepsilon>0$ implies that $B$ is zero-class $q$-admissible for any $q\in(\frac{1}{\alpha+\varepsilon},+\infty]$. This in turn implies that for such $q$ the map $t\mapsto \int_0^t T_{-1}(t-s) Bu(s)ds$ is continuous for any $u \in L^q(\R_+,U)$, by \cite[Proposition 2.3]{Wei89b}.

By Assumption~\ref{ass:Regularity-f-analytic}(iv), and using Proposition~\ref{prop:Auxiliary-map-analytic}, we see that for any $u\in L^\infty(\R_+,U)$ the map 
\[
t \mapsto \int_0^t T(t-s)f\big(x(s),u(s)\big)ds
\]
is well-defined and continuous.

Hence, if $x\in C(\R_+,X)$ with $x((0,+\infty)) \subset X_\alpha$, then for any  $u\in L^\infty(\R_+,U)$ the right-hand side of 
\eqref{eq:SEE+admissible_Integral_Form-analytic} is a continuous function of time.

Our next result is the local existence and uniqueness theorem for analytic systems with initial states in $X_\alpha$ and the inputs in $\Uc:=L^\infty_{\loc}(\R_+,U)$. Recall the notation $\Uc_\Sc$ from \eqref{eq:Inputs-constrained}.

\begin{theorem}[Picard-Lindel\"of theorem for analytic systems]
\label{PicardCauchy-analytic}
Let Assumption~\ref{ass:Regularity-f-analytic} hold.
Assume that $T$ is an analytic semigroup, satisfying for certain $M \geq 1$, $\lambda>0$ the estimate
\begin{eqnarray*}
\|T(t)\| \leq Me^{\lambda t},\quad t\geq 0.
\end{eqnarray*} 
For any compact set $Q \subset X_\alpha$, any $r>0$, any bounded set $\Sc \subset U$, and any $\delta>0$, there is a time $t_1 = t_1(Q,r,\Sc,\delta)>0$, such that for any $x_0 \in W:=B_{r,X_\alpha}(w)$ for some $w \in Q$, and for any $u\in\Uc_{\Sc}$ there is a unique mild solution of \eqref{eq:SEE+admissible} on $[0,t_1]$, and it lies in the ball $B_{Mr+\delta,X_\alpha}(w)$. 
\end{theorem}

\begin{proof}
First, we show the claim for the case if $Q=\{w\}$ is a single point in $X_\alpha$.

\textbf{(i).} Take any $\omega>\omega_0(T)$, and consider the corresponding space $X_\alpha$. 
Pick any $C>0$ such that $W := B_{r,X_\alpha}(w) \subset B_{C,X_\alpha}$, and $\Uc_\Sc \subset B_{C,\Uc}$. Pick any $u \in \Uc_\Sc$. Take also any $\delta>0$, and consider the following sets (depending on a parameter $t>0$):
\begin{eqnarray*}
Y_{t}:= \big\{y \in C([0,t],X): \|y(s) - (\omega I-A)^\alpha w\|_X\leq Mr + \delta \ \ \forall s \in [0,t]\big\},
\end{eqnarray*}
endowed with the metric 
$\rho_{t}(y_1,y_2):=\sup_{s \in [0,t]} \|y_1(s)-y_2(s)\|_X$.
As the sets $Y_t$ are closed subsets of Banach spaces $C([0,t],X)$, $Y_t$ are complete metric spaces for all $t>0$.

\textbf{(ii).} Pick any $x_0 \in W$. We are going to prove that for small enough $t$, the spaces $Y_{t}$  are invariant under the operator
$\Phi_u$, defined for any $y \in Y_{t}$ and all $\tau \in [0,t]$ by
\begin{align}
\label{eq:FP-Operator}
\Phi_u(y)(\tau) 
&= (\omega I-A)^\alpha T(\tau)x_0  + \int_0^\tau (\omega I-A)^\alpha T_{-1}(\tau-s) Bu(s)ds  \nonumber\\
&\qquad\qquad\qquad + \int_0^\tau (\omega I-A)^\alpha T(\tau-s)f\big((\omega I-A)^{-\alpha} y(s),u(s)\big)ds.
\end{align}

Since $y\in C([0,t],X)$, the map $s \mapsto (\omega I-A)^{-\alpha} y(s)$ is in $C([0,t],X_\alpha)$, as for any $s_1,s_2 \in[0,\tau]$ we have that 
\begin{align*}
\big\|(\omega I-A)^{-\alpha} y(s_1) - (\omega I-A)^{-\alpha} y(s_2)\big\|_{X_\alpha} = \|y(s_1) - y(s_2)\|_X.
\end{align*}

By Assumption~\ref{ass:Regularity-f-analytic}, the map $s\mapsto f\big((\omega I-A)^{-\alpha} y(s),u(s)\big)$ is in $L^w_{\loc}(\R_+,X)$, with $w>\frac{1}{1-\alpha}$.  
Proposition~\ref{prop:Auxiliary-map-analytic} ensures, that the map
\[
\tau\mapsto \int_0^\tau (\omega I-A)^\alpha T(\tau-s)f\big((\omega I-A)^{-\alpha} y(s),u(s)\big)ds
\]
is continuous.

Since $B\in L(U,X_{-1+\alpha+\varepsilon})$, Proposition~\ref{prop:ISS-analytic-systems}(ii) implies that
\[
\tau \mapsto \int_0^\tau (\omega I-A)^\alpha  T_{-1}(\tau-s)  Bu(s)ds
\]
belongs to $C([0,\tau],X)$.

Overall, the function $\Phi_u(y)$ is continuous, and thus $\Phi_u$ maps $C([0,t],X)$ to $C([0,t],X)$.

\textbf{(iii).}
Fix any $t>0$ and pick any $y \in Y_{t}$.
As $x_0 \in W=B_{r,X_\alpha}(w)$, there is $a \in X_\alpha$: $\|a\|_{X_\alpha} <r$ such that $x_0=w+a$.
 
Then for any $\tau<t$, we obtain that
\begin{align*}
\|&\Phi_{t}(y)(\tau)- (\omega I-A)^\alpha w\|_X \nonumber\\
&\leq   \Big\|(\omega I-A)^\alpha T(\tau)x_0 - (\omega I-A)^\alpha w\Big\|_X + \Big\|\int_0^\tau (\omega I-A)^\alpha T_{-1}(\tau-s) Bu(s)ds\Big\|_X \nonumber\\
&\qquad\qquad\qquad + \int_0^\tau \big\|(\omega I-A)^\alpha T(\tau-s)\big\| \big\|f((\omega I-A)^{-\alpha}y(s),u(s))\big\|_Xds. \nonumber
\end{align*}

We substitute $x_0:=w+a$ into the first term on the right-hand side of the above inequality. 
The last term we estimate using \eqref{eq:Bound-on-A^alphaT(t)}.
To estimate the second term, we use that $(\omega I-A)^{\alpha}  B \in L(U,X_{-1+\varepsilon})$. Thus, 
by Proposition~\ref{prop:ISS-analytic-systems}, $(\omega I-A)^{\alpha} B$ is zero-class $\infty$-admissible, and thus there is an increasing continuous function $t \mapsto h_t$ satisfying $h_0=0$, such that:
\begin{align}
\label{eq:technical-estimates-input-operator}
\big\|&\Phi_{t}(y)(\tau)- (\omega I-A)^\alpha w\big\|_X \leq \big\|(\omega I-A)^\alpha T(\tau)w - (\omega I-A)^\alpha w\big\|_X \\
			& \qquad + \|(\omega I-A)^\alpha T(\tau)a\|_X  + h_{\tau} \|u\|_{L^\infty([0,\tau],U)} \nonumber\\
			&\qquad+ \int_0^\tau \frac{C_\alpha e^{\lambda (\tau-s)}}{(\tau-s)^\alpha} \big(\|f(0,u(s))\|_X \nonumber\\
			&\qquad\qquad\qquad\qquad\qquad+ \big\|f((\omega I-A)^{-\alpha}y(s),u(s))-f(0,u(s))\big\|_X\big)ds. \nonumber
\end{align}

To estimate the latter expression, note that 
\begin{itemize}
	\item $\|(\omega I-A)^\alpha a\|_X = \|a\|_{X_\alpha} < r$. 
	\item For all $s\in[0,t]$ we have 
\begin{align*}
\|(\omega I-A)^{-\alpha}y(s)-0\|_{X_\alpha} = \|y(s)\|_X &\leq \|(\omega I-A)^\alpha w\|_X+Mr+\delta\\
&\leq M(\|w\|_{X_\alpha}+r) +\delta\leq MC+\delta :=K.
\end{align*}
	\item In view of Assumption~\ref{ass:Regularity-f-analytic}, it holds that 
\[
\|f(0,u(s))\|_X \leq \sigma(\|u(s)\|_U) + c,\quad \text{ for a.e. } s \in [0,t].
\]
	\item $h$ is a monotonically increasing continuous function.
\end{itemize}
As $M\geq 1$, it holds that $K>C$, and Lipschitz continuity of $f$ on bounded balls ensures that there is $L(K)>0$, such that for all $\tau\in[0,t]$ 			
\begin{align*}
\|&\Phi_{t}(y)(\tau)-(\omega I-A)^\alpha w\|_X \\			
			&\leq \|(\omega I-A)^\alpha T(\tau)w - (\omega I-A)^\alpha w\|_X + \|T(\tau)(\omega I-A)^\alpha a\|_X 
			+ h_{\tau} \|u\|_{L^\infty([0,\tau],U)}\\
			&\quad+ \int_0^\tau \frac{C_\alpha}{(\tau-s)^\alpha} e^{\lambda (\tau-s)} \big(\sigma(\|u(s)\|_U) + c + L(K)\|(\omega I-A)^{-\alpha}y(s)\|_{X_\alpha}\big)ds\\
			&\leq \sup_{\tau\in[0,t]}\|T(\tau)(\omega I-A)^\alpha w - (\omega I-A)^\alpha w\|_X + Me^{\lambda t}r + h_{t} \|u\|_{L^\infty([0,t],U)}\\
			&\quad+ C_\alpha e^{\lambda t} \big(\sigma(C) + c + L(K)K\big) \int_0^t \frac{1}{s^\alpha} ds.
\end{align*}
Since $T$ is a strongly continuous semigroup, and $h_t\to 0$ as $t\to +0$, from this estimate, it is clear that  there exists $t_1$, such that 
\[
\|\Phi_u(y)(t)-w\|_X \leq Mr +\delta,\quad \text{ for all } t \in [0,t_1].
\]
This means, that $Y_{t}$ is invariant with respect to $\Phi_u$ for all $t \in (0,t_1]$, and $t_1$ does not depend on the choice of $x_0 \in W$.

\textbf{(iv).} Now pick any $t>0$, $\tau \in [0, t]$, and any $y_1, y_2 \in Y_{t}$. Then it holds that
\begin{align*}
\|\Phi_u(y_1)&(\tau) - \Phi_u(y_2)(\tau)\|_X \\
&\leq \int_0^\tau \|(\omega I-A)^\alpha T(\tau-s)\|\\
&\qquad\qquad \cdot\big\|f((\omega I-A)^{-\alpha}y_1(s),u(s)) - f((\omega I-A)^{-\alpha}y_2(s),u(s))\big\|_Xds \\
&\leq 
\int_0^t L(K)\frac{C_\alpha}{(\tau-s)^\alpha} e^{\lambda (\tau-s)} \|y_1(s)-y_2(s)\|_X ds \\
&\leq L(K) C_\alpha e^{\lambda t} \int_0^t s^{-\alpha} ds \rho_t(y_1,y_2) \\
&\leq L(K) C_\alpha e^{\lambda t} \frac{t^{1-\alpha}}{1-\alpha} \rho_t(y_1,y_2) \\
&\leq \frac{1}{2} \rho_t(y_1,y_2),
\end{align*}
for $t \leq t_2$, where $t_2>0$ is a small enough real number that does not depend on the choice of $x_0 \in W$.
%

According to Banach fixed point theorem, there exists a unique $y \in Y_t$ that is a fixed point of $\Phi_u$, that is 
\begin{align}
\label{eq:FP-Operator-solution}
y(\tau) 
&= (\omega I-A)^\alpha T(\tau)x_0  + \int_0^\tau (\omega I-A)^\alpha T_{-1}(\tau-s) Bu(s)ds  \nonumber\\
&\qquad\qquad\qquad + \int_0^\tau (\omega I-A)^\alpha T(\tau-s)f\big((\omega I-A)^{-\alpha} y(s),u(s)\big)ds.
\end{align}
on $[0,\min\{t_1,t_2\}]$. 

As $(\omega I-A)^{\alpha}$ is invertible with a bounded inverse, $y$ solves \eqref{eq:FP-Operator-solution} if and only if 
$y$ solves
\begin{align}
\label{eq:FP-Operator-solution-II}
(\omega I-A)^{-\alpha} y(\tau) 
&= T(\tau)x_0  + \int_0^\tau T_{-1}(\tau-s) Bu(s)ds  \nonumber\\
&\qquad\qquad\qquad + \int_0^\tau T(\tau-s)f\big((\omega I-A)^{-\alpha} y(s),u(s)\big)ds.
\end{align}
As $y \in C([0, \min\{t_1,t_2\}],X)$, the map $x:=(\omega I-A)^{-\alpha} y$ is in $C([0, \min\{t_1,t_2\}],X_\alpha)$, and is the unique mild solution of \eqref{eq:SEE+admissible-analytic}. 

\textbf{(v). General compact $Q$.} 
 Similar to the corresponding part of the proof of Theorem~\ref{PicardCauchy}.
~
\end{proof}

\begin{remark}
\label{rem:Existence for initial states in X} 
We have proved our local existence result for initial conditions that are in $X_\alpha$. To ensure local existence and uniqueness 
for the initial states outside of $X_\alpha$, stronger requirements on $f$ have to be imposed, see \cite[Theorems 7.1.5, 7.1.6]{Lun12}.
\end{remark}

Introducing the concepts of maximal solutions and of well-posedness and arguing similar to Sections~\ref{sec:Local existence and uniqueness}, \ref{sec:Global-Well-posedness}, we obtain the following well-posedness theorem.
\begin{theorem}
\label{thm:Well-posedness-analytic}
Let $A$ generate an analytic semigroup, Assumption~\ref{ass:Regularity-f-analytic} hold, and let $\Uc:=L^\infty(\R_+,U)$. 
Then:
\begin{enumerate}[label=(\roman*)]
	\item For each $x\in X_\alpha$ and each $u \in\Uc$, there is a unique maximal solution of \eqref{eq:SEE+admissible-analytic}, defined over the certain maximal interval $[0,t_m(x,u))$. We denote this solution as $\phi(\cdot,x,u)$.
	\item The triple $\Sigma:=(X_\alpha,\Uc,\phi)$ is a well-defined control system in the sense of Definition~\ref{Steurungssystem}.
	\item $\Sigma$ satisfies the BIC property, that is, if for a certain $x \in X_\alpha$ and $u\in\Uc$ we have $t_m(x,u)<\infty$, then 
	 $\|\phi(t,x,u)\|_{X_\alpha} \to \infty$ as $t \to t_m(x,u)-0$.
\end{enumerate}
\end{theorem}

\subsection{Global existence}

Now we are going to derive sufficient conditions for forward completeness and BRS property.
We need the following \q{analytic Gronwall-Bellman lemma}:
%

\begin{proposition}
\label{prop:Gronwall-bellman-analytic-II} 
For any $a,b\geq 0$, any $\alpha,\beta \in[0,1)$, and any $T\in(0,+\infty)$ there is $M=M(b,\alpha,\beta,T)>0$ such that for any integrable function $u:[0,T]\to\R$ satisfying for almost all $t \in[0,T]$ the inequality 
\[
0\leq u(t)\leq at^{-\alpha} + b\int_0^t(t-s)^{-\beta}u(s)ds
\]
it holds for a.e. $t\in[0,T]$ that
\[
0\leq u(t)\leq aMt^{-\alpha}.
\]
\end{proposition}

\begin{proof}
See \cite[p. 6]{Hen81}.
\end{proof}

Motivated by \cite[Section 6.3, Theorem 3.3]{Paz83}, we have the following result guaranteeing the forward completeness and BRS property for semilinear analytic systems.
\begin{theorem}
\label{thm:Global-Existence-Nonlinear-Eq} 
Let $A$ generate an analytic semigroup, Assumption~\ref{ass:Regularity-f-analytic} hold, and let $\Uc:=L^\infty(\R_+,U)$. 
Assume further that there are $L,c>0$ and $\sigma\in\Kinf$ such that 
\begin{eqnarray}
\|f(x,u)\|_X \leq L\|x\|_{X_\alpha} + \sigma(\|u\|_U) + c,\quad x\in X_\alpha,\quad u\in U.
\label{eq:f-linearly-bounded}
\end{eqnarray} 
Then  $\Sigma:=(X_\alpha,\Uc,\phi)$ is a forward complete control system.
\end{theorem}

\begin{proof}
Take any positive $\omega >\omega_0(T)$ and define $X_\alpha$ as in \eqref{eq:X-power-alpha}.

We argue by a contradiction.
Let $\Sigma$ be not forward complete. Then there are $(x_0,u)\in X_\alpha \tm \Uc$ such that $t_m(x_0,u)<\infty$. By Theorem~\ref{thm:Well-posedness-analytic}, we have that $\|\phi(t,x_0,u)\|_{X_\alpha} \to \infty$ as $t \to t_m(x_0,u)-0$.

For $t<t_m(x_0,u)$ denote $x(t):=\phi(t,x_0,u)$. 
As $x(\cdot) \subset X_\alpha$, we can apply $(\omega I-A)^\alpha$ along the trajectory $x(\cdot)$ to obtain
\begin{align*}
(\omega I-A)^\alpha x(t) &= (\omega I-A)^\alpha T(t)x_0 + \int_0^\tau (\omega I-A)^\alpha T_{-1}(\tau-s) Bu(s)ds  \\
			&\qquad	+ \int_{0}^t (\omega I-A)^\alpha T(t-s)f(x(s),u(s))ds.
\end{align*}

We obtain 
\begin{align*}
\|x(t)\|_\alpha 
&= \|(\omega I-A)^\alpha x(t)\|_X  \\
&\leq \|(\omega I-A)^\alpha T(t)x_0\|_X + \Big\|\int_0^\tau (\omega I-A)^\alpha T_{-1}(\tau-s) Bu(s)ds\Big\|_X  \\
			&\qquad\qquad	+ \int_{0}^t \big\|(\omega I-A)^\alpha T(t-s)\big\| \big\|f(x(s),u(s))\big\|_Xds.
\end{align*}
We now estimate the second term as in \eqref{eq:technical-estimates-input-operator}, where $h$ is a continuous increasing function with $h_0=0$. The last term we estimate using \eqref{eq:Bound-on-A^alphaT(t)}. Overall:
\begin{align*}
\|x(t)\|_\alpha 
&\leq Me^{\omega t} \|(\omega I-A)^\alpha x_0\|_X + h_{\tau} \|u\|_{L^\infty([0,t],U)}\\
			&\qquad\qquad	+ \int_{0}^t \frac{C_\alpha}{(t-s)^\alpha}e^{\omega (t-s)} 
\Big(L\|x(s)\|_{X_\alpha} + \sigma(\|u(s)\|_U) + c \Big)ds.
\end{align*}
Defining $z(t):=x(t)e^{-\omega t}$, we obtain from the previous estimate that
\begin{align*}
\|z(t)\|_\alpha 
&\leq M\|(\omega I-A)^\alpha x_0\|_X + \int_{0}^t \frac{C_\alpha}{s^\alpha}ds  \Big(\sigma(\|u\|_{L^\infty(\R_+,U)}) + c \Big)
+ h_{\tau} \|u\|_{L^\infty(\R_+,U)}\\
&\qquad + \int_{0}^t \frac{L C_\alpha}{(t-s)^\alpha}\|z(s)\|_{X_\alpha}ds.  
\end{align*}
Proposition~\ref{prop:Gronwall-bellman-analytic-II} shows that $z$, and hence $x$, is uniformly bounded on $[0,t_m(x_0,u))$, and BIC property (Theorem~\ref{thm:Well-posedness-analytic}(iii)) shows that $t_m(x_0,u)$ is not the finite maximal existence time. A contradiction.
\end{proof}


\subsection{Example: well-posedness of a Burger's equation with a distributed input}

Motivated by \cite[p. 57]{Hen81}, we consider the following semilinear reaction-diffusion equation of Burgers' type on a domain $[0,\pi]$, with distributed input $u$, boundary input $d$ at $z=0$, and homogeneous Dirichlet boundary condition at $\pi$.
\begin{subequations}
\label{eq:Nonlinear-Burgers-Equation}
\begin{align}
x_t &= x_{zz} - xx_z  + f(z,x(z,t))+ u(z,t),\quad z\in (0,\pi),\quad t>0,\\
x&(0,t)=d(t),\quad t>0,\\
x&(\pi,t)=0.
\end{align}
\end{subequations}

Here $f: [0,\pi]\tm \R \to\R$ is  measurable  in  $z$,  locally  Lipschitz  continuous  in  $x$ uniformly  in  $z$, and 
\begin{eqnarray}
|f(z,y)| \leq h(z)g(|y|),\quad \text{for a.e. } z \in[0,\pi], \text{ and all } y \in\R,
\label{eq:f-bound}
\end{eqnarray}
where $h\in L^2(0,\pi)$, and $g$ is continuous, increasing, and both $h$ and $g$ are positive.

We denote $X:=L^2(0,\pi)$. The operator $A:= \frac{d^2}{dz^2}$ with the domain $D(A) = H^2(0,\pi)\cap H^1_0(0,\pi)$ generates an analytic semigroup on $X$.

We assume that the distributed input $u$ belongs to the space $\Uc=L^\infty(\R_+,U)$, with $U:=L^2(0,\pi)$, and the boundary input $d$ belongs to $\Dc:=L^\infty(\R_+,\R)$.

The system \eqref{eq:Nonlinear-Burgers-Equation} can be reformulated as a semilinear evolution equation
\begin{eqnarray}
x_t = Ax + F(x)+u + Bd,
\label{eq:Nonlinear-Burgers-Equation-Evolution-Eq}
\end{eqnarray}
where we slightly abuse the notation and use $x$ as an argument of the evolution equation.

The condition (ii) in Assumption~\ref{ass:Regularity-f-analytic} characterizing the admissibility properties of the boundary input operator $B$ holds in view of \cite[Example 2.16]{Sch20}.
  
The space $X_{\frac{1}{2}}$ corresponding to the operator $A$, is given by (see \cite[Proposition 3.6.1]{TuW09}) 
\[
X_{\frac{1}{2}} = H^1_0(0,\pi),
\]
which is a Banach space with the norm
\[
\|x\|_{\frac{1}{2}}:=\Big|\int_0^\pi |x'(z)|^2 dz\Big|^{\frac{1}{2}}, \quad x \in X_{\frac{1}{2}}.
\]

The nonlinearity $F: X_{\frac{1}{2}} \to X$ in \eqref{eq:Nonlinear-Burgers-Equation-Evolution-Eq} is given by 
\[
F(x)(z) = - x(z)x'(z) + f(z,x(z)).
\]

\begin{proposition}
\label{prop:Burgers-well-posedness} 
For each $x_0 \in X_{\frac{1}{2}}$, each $u \in \Uc=L^\infty_{\loc}(\R_+,U)$, and each boundary input $d \in \Dc=L^\infty_{\loc}(\R_+,\R)$ the system \eqref{eq:Nonlinear-Burgers-Equation} possesses a unique maximal mild solution $\phi(\cdot,x_0,(u,d))$.
The system $\Sigma = (X_{\frac{1}{2}},\Uc\tm\Dc,\phi)$ is a control system satisfying the BIC property.
\end{proposition}

\begin{proof}
We proceed in 3 steps:

\textbf{Step 1: $F$ maps bounded sets of $ X_{\frac{1}{2}}$ to bounded sets of $X$.} 
Since the elements of $X_{\frac{1}{2}}=H^1_0(0,\pi)$ are absolutely continuous functions, using the Cauchy-Schwarz inequality, we obtain that for any $x \in X_{\frac{1}{2}}$ it holds that 
\begin{eqnarray}
\sup_{z\in(0,\pi)}|x(z)| 
&=& \sup_{z\in(0,\pi)}\Big|\int_0^z x'(z) dz\Big| \nonumber \\
&\leq& \sup_{z\in(0,\pi)}\int_0^z |x'(z)| dz \nonumber\\ 
& = & \int_0^\pi |x'(z) |dz  \nonumber\\
&\leq & \Big|\int_0^\pi 1 dz\Big|^{\frac{1}{2}}  \Big|\int_0^\pi |x'(z)|^2 dz\Big|^{\frac{1}{2}}  \nonumber\\
&= &\sqrt{\pi} \|x\|_{\frac{1}{2}}.
\label{eq:sup-norm-estimated-by-Sobolev-norm}
\end{eqnarray}

For any $x \in X_{\frac{1}{2}}$ consider
\begin{eqnarray*}
\|F(x)\|_X^2
&= &\int_0^\pi |F(x)(z)|^2 dz\\
&= &\int_0^\pi |x(z)x'(z) + f(z,x(z))|^2 dz\\
&\leq &\int_0^\pi 2|x(z)x'(z)|^2 + 2|f(z,x(z))|^2 dz.
\end{eqnarray*}
Using \eqref{eq:sup-norm-estimated-by-Sobolev-norm} and \eqref{eq:f-bound}, we continue the estimates as follows:
\begin{eqnarray*}
\|F(x)\|_X^2
&\leq &\int_0^\pi 2|x'(z)|^2 \pi \|x\|^2_{\frac{1}{2}} + 2 |h(z)|^2  |g(|x(z)|)|^2 dz\\
&\leq & 2\pi \|x\|^4_{\frac{1}{2}}  + 2 \int_0^\pi |h(z)|^2 |g(\sqrt{\pi} \|x\|_{\frac{1}{2}})|^2 dz\\
&= & 2\pi \|x\|^4_{\frac{1}{2}}  + 2 \|h\|^2_X |g(\sqrt{\pi} \|x\|_{\frac{1}{2}})|^2.
\end{eqnarray*}
 Taking the square root and using that $\sqrt{a+b}\leq \sqrt{a} + \sqrt{b}$, for all $a,b\geq 0$, we finally obtain
\begin{eqnarray}
\|F(x)\|_X
&\le & \sqrt{2\pi} \|x\|^2_{\frac{1}{2}}  + \sqrt{2} \|h\|_X |g(\sqrt{\pi} \|x\|_{\frac{1}{2}})|.
\label{eq:Estimate-on-F}
\end{eqnarray}
This shows that $F$ is well-defined as a map from $X_{\frac{1}{2}}$ to $X$ and $F$ maps bounded sets of $X_{\frac{1}{2}}$ to bounded sets of $X$.

\textbf{Step 2: $F$ is Lipschitz continuous in $x$.} 
For $x_0 \in X_{\frac{1}{2}}$, there is a neighborhood $V$ of a compact set 
$
\{(z,x_0(z)):z\in[0,\pi]\}
$ 
in $[0,\pi]\tm \R$ and positive constants $L,\theta$ so that for
$(z,x_1) \in V$, $(z,x_2) \in V$ it holds that 
\[
|f(z,x_1) - f(z,x_2)| \leq L |x_1-x_2|.
\]
Thus, there is a neighborhood $U$ of $x_0$ in $X_{\frac{1}{2}}$
such that  $x\in U$ implies that $(z,x(z))\in V$ for a.e. $z\in[0,\pi]$ and for $x_1,x_2 \in U$ it holds that
\begin{align*}
\big\|f(\cdot,x_1(\cdot)) - f(\cdot,x_2(\cdot))\big\|_{X}^2
&= \int_0^\pi\big|f(z,x_1(z)) - f(z,x_2(z))\big|^2dz \\
&\leq L^2  \int_0^\pi|x_1(z)-x_2(z)|^2 dz,
\end{align*}
and using \eqref{eq:sup-norm-estimated-by-Sobolev-norm} we have that
\begin{eqnarray*}
\big\|f(\cdot,x_1(\cdot)) - f(\cdot,x_2(\cdot))\big\|_{X}^2
&\leq& L^2 \pi^2 \|x_1-x_2\|^2_{\frac{1}{2}}.
\end{eqnarray*}
Taking the square root, we have that  
\begin{eqnarray}
\big\|f(\cdot,x_1(\cdot)) - f(\cdot,x_2(\cdot))\big\|_{X}
\leq \pi L \|x_1-x_2\|_{\frac{1}{2}}.
\label{eq:Lipschitz estimate-for-f}
\end{eqnarray}
Finally, for any $x_1,x_2\in X_{\frac{1}{2}}$ it holds that 
\begin{eqnarray*}
\|x_1x_1' - x_2x_2'\|_X &\leq& \|x_1(x_1'-x_2')\|_X + \|(x_1-x_2)x_2'\|_X,
\end{eqnarray*}
and again using \eqref{eq:sup-norm-estimated-by-Sobolev-norm}, we proceed to 
\begin{eqnarray}
\|x_1x_1' - x_2x_2'\|_X 
&\leq& \sqrt{\pi} \|x_1\|_{\frac{1}{2}}  \|x_1-x_2\|_{\frac{1}{2}} + \sqrt{\pi} \|x_1-x_2\|_{\frac{1}{2}} \|x_2\|_{\frac{1}{2}}\nonumber\\
&=& \sqrt{\pi} ( \|x_1\|_{\frac{1}{2}}  + \|x_2\|_{\frac{1}{2}})\|x_1-x_2\|_{\frac{1}{2}}.
\label{eq:Lipschitz estimate-for-cross-term}
\end{eqnarray}
Combining \eqref{eq:Lipschitz estimate-for-f} and \eqref{eq:Lipschitz estimate-for-cross-term}, we obtain the required Lipschitz property for the function $F$.

\textbf{Step 3: Application of general well-posedness theorems.} Finally, Theorem~\ref{PicardCauchy-analytic} shows that the system \eqref{eq:Nonlinear-Burgers-Equation} possesses a unique mild solution for each $x_0 \in X_{\frac{1}{2}}$, each $u \in \Uc=L^\infty_{\loc}(\R_+,U)$, and each boundary input $d \in \Dc=L^\infty_{\loc}(\R_+,\R)$. Theorem~\ref{thm:Well-posedness-analytic} shows that $\Sigma$ is a control system satisfying the BIC property.
\end{proof}

\section{Concluding remarks}
\label{sec:Well-posedness-Concluding remarks}

The framework of linear abstract control systems, discussed in Sections~\ref{sec:Abstract linear control systems}, \ref{sec:Representation of abstract control systems}, and most of the results in these sections, are due to \cite{Wei89b}.

\textbf{ISS of linear systems.}
In our exposition of Section~\ref{sec:ISS of linear systems}, we follow \cite{MiP20}. We adopt the notion of admissibility of input operators from \cite{Wei89b}. For a detailed exposition of the admissibility concepts for input and observation operators, we refer to \cite{TuW09}, \cite{JaP04}. 
Proposition~\ref{prop:Restatement-admissibility} stems from \cite[Proposition 4.2]{Wei89b}; see also \cite[Proposition 4.2.2]{TuW09}. 
Theorem~\ref{thm:ISS-Criterion-lin-sys-with-unbounded-operators} is due to \cite[Proposition 2.10]{JNP18}.
Theorem~\ref{thm:Contiuity-of-a-map} is due to \cite{JSZ19}. 
Lemma~\ref{lem:GrabowskiLemma} is a variation of \cite[Lemma 1.1]{Gra95}.

\textbf{Spectral-based methods and related techniques.} Our development of the ISS theory for linear distributed parameter systems mostly follows that of \cite{JNP18, JSZ19, JMP20, MiP20}. 
There are, however, also alternative approaches. 

In \cite{KaK16b, KaK17a, KaK17b}, the ISS analysis of linear parabolic PDEs with Sturm-Liouville operators over a 1-dimensional spatial domain has been performed using two different methods: (i) the spectral decomposition of the solution, and (ii) the approximation of the solution by means of a  finite-difference scheme. 
We refer to the monograph \cite{KaK19} that provides a systematic description of this method.
A brief description of this method is given in \cite[Section 4.2]{MiP20}.

The spectral method has been used in the papers \cite{LSZ20}, \cite{LhS19} for the study of ISS of Riesz-spectral BCS, which is a modification of the spectral method from  \cite{KaK16b, KaK17a, KaK17b}. 
The essence of the method is a decomposition of $X$ w.r.t.\ the Riesz basis $\{\phi_k:k\in\N\}$.

\textbf{Semilinear evolution equations: Generality.}  In this chapter, we analyzed the well-posedness and properties of the flow for semilinear evolution equations of the form
\begin{subequations}
\label{eq:SEE+admissible-intro} 
\begin{eqnarray}
\dot{x}(t) & = & Ax(t) + B_2f(x(t),u(t)) + Bu(t),\quad t>0,  \label{eq:SEE+admissible-intro1}\\
x(0)  &=&  x_0. \label{eq:SEE+admissible-intro2}
\end{eqnarray}
\end{subequations}
Here $A$ generates a strongly continuous semigroup over a Banach space $X$, the operators $B$ and $B_2$ are admissible, and $f$ is a Lipschitz continuous in the first variable map. 
In our exposition of Sections~\ref{sec:Semilinear boundary control systems},~\ref{sec:Boundary_control_systems},~\ref{sec:Semilinear analytic boundary control systems}, we follow \cite{Mir22c} rather closely, and most of the results are from this reference. 

The class of systems \eqref{eq:SEE+admissible-intro} is rather general:
\begin{itemize}
	\item If $B$ and $B_2$ are bounded operators, \eqref{eq:SEE+admissible-intro} corresponds to the classic semilinear evolution equations covering broad classes of semilinear PDEs with distributed inputs. If $A$ is a bounded operator, such a theory was developed in \cite{DaK74} (for time-variant systems without inputs). In the case of unbounded generators $A$, we refer to \cite{Paz83}, \cite{Hen81}, \cite[Chapter 11]{CuZ20}, \cite{CaH98}, etc. 

	\item If $B_2=0$, and $B$ is an admissible operator, then \eqref{eq:SEE+admissible-intro} reduces to the class of general linear control systems that fully covers linear boundary control systems (see \cite{CuZ20, JaZ12}, \cite{TuW09,TuW14}, \cite{EmT00} for an overview). In particular, this class includes linear evolution PDEs with boundary inputs.
	
	\item Consider a linear system 
\begin{eqnarray}
\dot{x} = Ax + Bu,
\label{eq:linear-open-loop-intro}
\end{eqnarray}
with admissible $B$. Let us apply a feedback controller $u(x)=f(x,d_1)+d_2$ that is subject to additive actuator disturbance $d_2$ and further disturbance input $d_1$. Substituting this controller into \eqref{eq:linear-open-loop-intro}, we arrive 
at systems \eqref{eq:SEE+admissible-intro}, with $B_2=B$. Therefore the systems \eqref{eq:SEE+admissible-intro}  are sometimes called Lipschitz perturbations of linear systems.

	\item In \cite{HuP21}, it was shown that the class of systems \eqref{eq:SEE+admissible-intro} includes 2D Navier-Stokes equations (under certain boundary conditions) with in-domain inputs and disturbances. Furthermore, in \cite{HuP21}, the authors have designed an error feedback controller that guarantees approximate local velocity output tracking for a class of reference outputs.
	Viscous Burgers' equation with nonlinear local terms and boundary inputs of Dirichlet or Neumann type falls into the class \eqref{eq:SEE+admissible-intro} as well.
		
	\item In \cite{Sch20}, it was shown that semilinear boundary control systems with linear boundary operators could be considered a special case of systems \eqref{eq:SEE+admissible-intro}. In this case, it suffices to consider $B_2$ as the identity operator.
Furthermore, in \cite{Sch20}, the well-posedness and input-to-state stability of a class of analytic boundary control systems with nonlinear dynamics and a linear boundary operator were analyzed with the methods of operator theory. 
\end{itemize}

Some highly nonlinear systems (even without inputs) as, e.g., the porous medium equation \cite{Vaz07}, the
nonlinear KdV equation \cite{BaK00} are not covered by \eqref{eq:SEE+admissible-intro}, and should be modeled using methods of nonlinear semigroup theory, which is closely connected to the theory of maximal monotone operators \cite{Bar10}.

\textbf{Semilinear evolution equations: State of the art.} 
The systems \eqref{eq:SEE+admissible-intro} have been studied (up to the assumptions on $f$, and the choice of the space of admissible inputs) in \cite{NaB16} under the assumption that its linearization is an exponentially stable regular linear system in the sense of \cite{TuW09, TuW14, Sta05}.
\cite{NaB16} ensures local well-posedness of regular nonlinear systems assuming the Lipschitz continuity of nonlinearity and invoking regularity of the linearization.
On this basis, the authors show in \cite{NaB16} that  
an error feedback controller designed for robust output regulation of the linearization of a regular nonlinear system
 achieves approximate local output regulation for regular nonlinear systems.

Control of systems \eqref{eq:SEE+admissible-intro} has been studied recently in several papers. In particular, in \cite{NZW19}, the exact controllability of a class of regular nonlinear systems was studied using back-and-forth iterations.
A problem of robust observability was studied for a related class of systems in \cite{JLZ15}.

Stabilization of linear port-Hamiltonian systems by means of nonlinear boundary controllers was studied in \cite{Aug16,RZG17}.
Bounded controls with saturations (a priori limitations of the input signal) have been employed for PDE control in \cite{PTS16,TMP18,MPW21}. 
Recently, several papers appeared that treat nonlinear boundary control systems within the input-to-state stability framework. 
Nonlinear boundary feedback was employed for the ISS stabilization of linear port-Hamiltonian systems in \cite{ScZ21}. 

Several types of infinite-dimensional systems, distinct from \eqref{eq:SEE+admissible-intro}, have been studied as well. 
Time-variant infinite-dimensional semilinear systems are one of such classes. Such systems have been first studied (as far as the author is concerned) for systems without disturbances in \cite{JDP95}.
Recently, in \cite{Sch22}, sufficient conditions for well-posedness and uniform global stability have been obtained for scattering-passive semilinear systems (see \cite[Theorem 3.8]{Sch22}). 

Semilinear systems with outputs are another important extension of \eqref{eq:SEE+admissible-intro}. 
Such systems with globally Lipschitz nonlinearities have been analyzed in \cite[Section 7]{TuW14}, and it was shown that such systems are well-posed and forward complete provided that the Lipschitz constant is small enough. 
In \cite{HCZ19} employing counterexamples, it was shown that semilinear systems with a nonlinear output boundary feedback might fail to be well-posed. Well-posedness of incrementally scattering-passive nonlinear systems with outputs has been analyzed in \cite{SWT22} by applying Crandall-Pazy theorem \cite{CrP69} on the generation of
nonlinear contraction semigroups to a Lax-Phillips nonlinear semigroup representing the system together with its inputs and outputs.

\textbf{Boundary control systems.} For many natural and engineering systems, the interaction of a system with its environment (by controls, inputs, and outputs) occurs at the boundary of the system. Examples of such behavior are given by diffusion equations \cite{AWP12}, vibration of structures \cite{CuZ95}, transport phenomena, etc., with broad applications in robotics \cite{EMJ17}, aerospace engineering \cite{BCC16, PGC13}, and additive
manufacturing \cite{DiK15,HMR15}. Wide classes of port-Hamiltonian systems can be formulated as boundary control systems as well, see \cite{JaZ12,ScJ14}.

The development of the theory of general boundary control systems has been initiated in the pioneering work \cite{Fat68} and was brought further forward by \cite{Sal87}; see \cite{TuW14} for a survey. In the literature, there are several ways how to define a boundary control system; see, e.g., \cite{Fat68}, \cite{CuZ95, JaZ12}, 
\cite{Sal87,TuW09} and \cite{EmT00}. The differences between various methods are discussed in \cite{LGE00}.
In our description of linear boundary control systems, we use the strategy due to \cite{JaZ12}, with some motivation from \cite{EmT00}.

At the same time, the study of nonlinear systems with boundary controls is a much younger subject. For some recent references, we refer to \cite{TuW14,HCZ19,Sch20}. In our exposition of semilinear boundary control systems, we follow \cite{Mir22c}, which is motivated by \cite{Sch20}.


\cleardoublepage
\chapter{ISS superposition theorems}
\label{chap:Characterizations_ISS_1}

ISS superposition theorems play an important role in the finite-dimensional ISS theory. They characterize ISS as a combination of a certain kind of stability and attractivity properties, properly defined for the systems with inputs.
In this chapter, we prove the ISS superposition theorems for infinite-dimensional systems and underline the differences to the finite-dimensional theory. We apply these findings for the derivation of so-called integral characterizations of ISS.
More applications of ISS superposition theorems we will find in Chapter~\ref{sec:Non-coercive ISS Lyapunov functions} where they are used to justify the usage of non-coercive ISS Lyapunov functions, and in Chapter~\ref{chap:Infinite interconnections: Non-Lyapunov methods} where they are applied for the derivation of the small-gain theorems for finite and infinite networks of ISS systems.

\begin{ass}
The control systems considered in this chapter are always forward complete, if not mentioned otherwise.
\end{ass}

\section{Stability properties}
\label{sec:Prelim}

%
%

\subsection{Stability notions for undisturbed systems}

We start with systems without inputs.

\begin{definition}
\label{Stab_Notions_Undisturbed_Systems}
\index{0-ULS}
\index{stability! uniform at zero}
A system $\Sigma=(X,\Uc,\phi)$ is called
\begin{itemize}
\item[(i)] \emph{uniformly stable at zero (0-ULS)}, if there exists \mbox{$\sigma \in \Kinf$} and $r>0$ so that
 \begin{equation}
 \label{eq:1}
\|\phi(t,x,0)\|_X \leq \sigma(\|x\|_X) \quad \forall x \in \clo{B_r}\quad \forall t \geq 0.
 \end{equation}

\index{0-UGS}
\index{stability! uniform global at zero} 
   \item[(ii)] \emph{uniformly globally stable at zero (0-UGS)}, if there exists $
          \sigma \in \Kinf$ so that
          \begin{equation}
              \label{eq:3}
\|\phi(t,x,0)\|_X \leq \sigma(\|x\|_X) \quad\forall x \in X\quad \forall t \geq 0.
          \end{equation}
%

\index{0-GATT}
\index{attractivity!global at zero}
    \item[(iii)] \emph{globally attractive at zero (0-GATT)}, if
\begin{equation}
\label{GATT_LIM}
\lim_{t \to \infty} \left\| \phi(t,x,0) \right\|_X = 0\quad \forall x\in X.
\end{equation}

\index{0-LIM}
\index{limit property at zero}
    \item[(iv)] a system with the \emph{limit property at zero (0-LIM)}, if
\[
\inf_{t \geq 0} \|\phi(t,x,0)\|_X = 0\quad \forall x\in X.
\]

\index{0-UGATT}
\index{attractivity! uniform global at zero}
    \item[(v)] \emph{uniformly globally attractive at zero (0-UGATT)}, if
          for all $\eps, \delta >0$ there is a $\tau=\tau(\eps,\delta) < \infty$ such that
\begin{equation}
\label{UnifGATT}
 t \geq \tau,\quad x \in \clo{B_{\delta}} \qrq \|\phi(t,x,0)\|_X \leq \eps.
\end{equation}

\index{0-GAS}
\index{stability!global asymptotic at zero}
  \item[(vi)] \emph{globally asymptotically stable at zero} (0-GAS), if $\Sigma$ is 0-ULS and 0-GATT.

\index{0-UAS}
\index{stability!uniform asymptotic at zero}

    \item[(vii)] \emph{uniformly asymptotically stable at zero} (0-UAS), if there exists a $ \beta \in \KL$
          and $r>0$, such that for all $x \in \clo{B_r}$
\begin{equation}
\label{UniStabAbschaetzung}
\left\| \phi(t,x,0) \right\|_{X} \leq  \beta(\left\| x \right\|_{X},t)\quad  \forall t\geq 0.
\end{equation}

\end{itemize}
\end{definition}

We stress the difference between the uniform notions 0-UGATT and 0-UGAS and the nonuniform notions 0-GATT and 0-GAS.
For 0-GATT systems, all trajectories converge to the origin,
but their speed of convergence may differ drastically for initial values with the same norm, in contrast to 0-UGATT systems.
The notions of 0-ULS and 0-UGS are uniform in the sense that there exists an upper bound of the norm of trajectories which is equal for initial states with the same norm.

\begin{remark}
\label{0-GAS_strong_stability}
For ODE systems, 0-GAS is equivalent to 0-UGAS, but it is weaker than
0-UGAS in the infinite-dimensional case. For linear systems $\dot{x} =
Ax$, where $A$ generates a strongly continuous semigroup on a Banach space $X$, the
Banach-Steinhaus theorem implies that 0-GAS is equivalent to the strong stability of the associated semigroup $T$ and implies the 0-UGS property.
\end{remark}

For systems that are 0-LIM, trajectories approach the origin arbitrarily closely. Obviously, 0-GATT implies 0-LIM.

\subsection{Stability notions for systems with inputs}

We now consider control systems with inputs.
\begin{definition}
\label{def:ULS-UGS-pGS}
A system $\Sigma=(X,\Uc,\phi)$ is called
\begin{itemize}
\index{ULS}
\index{stability!uniform local}
    \item[(i)] \emph{uniformly locally stable (ULS)}, if there exist $ \sigma \in\Kinf$, $\gamma
          \in \Kinf \cup \{0\}$ and $r>0$ such that for all $ x \in \clo{B_r}$ and all 
					$u \in \clo{B_{r,\Uc}}$ the trajectory $\phi(\cdot,x,u)$ is defined on $\R_+$, and
\begin{equation}
\label{GSAbschaetzung}
\left\| \phi(t,x,u) \right\|_X \leq \sigma(\|x\|_X) + \gamma(\|u\|_{\Uc}) \quad \forall t \geq 0.
\end{equation}

\index{UGS}
\index{stability!uniform global}
  \item[(ii)] \emph{uniformly globally stable (UGS)}, if $\Sigma$ is forward complete and there exist $ \sigma \in\Kinf$, $\gamma
          \in \Kinf \cup \{0\}$ such that for all $ x \in X, u
          \in \Uc$ the estimate \eqref{GSAbschaetzung} holds.

\index{UGB}
\index{uniform global boundedness}					
    \item[(iii)] \emph{uniformly globally bounded (UGB)}, if $\Sigma$ is forward complete and there exist $ \sigma \in\Kinf$, $\gamma
          \in \Kinf \cup \{0\}$ and $c>0$ such that for all $ x \in X$,
          and all $ u \in \Uc$  it holds that
\begin{equation}
\label{pGSAbschaetzung}
\left\| \phi(t,x,u) \right\|_X \leq \sigma(\|x\|_X) +\gamma(\|u\|_{\Uc}) + c \quad \forall t \geq 0.
\end{equation}
\end{itemize}

\end{definition}

We start with a standard reformulation of the $\varepsilon$-$\delta$ formulations of stability in terms of ${\cal K}$-functions.
\begin{lemma}
\label{lem:ULS_restatement}
A system $\Sigma=(X,\Uc,\phi)$ is ULS if and only if for all $\eps>0$ there
exists a $\delta>0$ such that
\begin{equation}
\label{LS_Restatement}
\|x\|_X \leq\delta,\quad \|u\|_{\Uc} \leq \delta,\quad t\geq 0 \qrq \|\phi(t,x,u)\|_X \leq\eps.
\end{equation}
\end{lemma}

\begin{proof}
"$\Rightarrow$"
Let $\Sigma$ be ULS. Let $\sigma,\gamma \in \Kinf$ and $r>0$ be such that
\eqref{GSAbschaetzung} holds for these functions and the neighborhood
specified by $r$.
Let $\eps>0$ be arbitrary and choose
\begin{equation*}
    \delta=\delta(\eps):=\min\left\{\sigma^{-1}\left(\frac{\eps}{2}\right),
\gamma^{-1}\left(\frac{\eps}{2}\right), r\right\}.
\end{equation*}
With this choice, \eqref{LS_Restatement} follows from \eqref{GSAbschaetzung}.

"$\Leftarrow$" Let \eqref{LS_Restatement} hold. For $\varepsilon \geq 0$ define
\begin{align*}
\delta(\eps):=\sup\{ s \geq 0 : \|x\|_X \leq s \ \wedge \  \|u\|_{\Uc} \leq s 
\Rightarrow \sup_{t \geq 0} \|\phi(t,x,u)\|_X \leq\eps  \}.
\end{align*}
Clearly, \eqref{LS_Restatement} implies that $\delta(\cdot)$ is well
defined, increasing and continuous in $0$. Let $\hat{\delta} \in {\cal K}$ be any function with
$\hat{\delta}\leq \delta$ and set $r:= \sup_{s\geq 0}  \hat{\delta}(s) \in
\R_+ \cup {\infty}$. Define $\gamma:=\hat{\delta}^{-1}:[0,r) \to \R$.
Then for $\|x\|_X < r$ and $\|u\|_{\Uc}<r$ we have
\[
\|\phi(t,x,u)\|_X \leq \gamma(\max\{ \|x\|_X, \|u\|_{\Uc} \}) \leq \gamma(\|x\|_X) + \gamma(\|u\|_{\Uc}),
\]
which shows ULS.
\end{proof}

Lemma~\ref{lem:ULS_restatement} can be interpreted as follows:
the system $\Sigma$ is ULS if and only if $\phi$ is continuous at the equilibrium $0$
and the function $\delta$ in Definition~\ref{def:RobustEquilibrium_Undisturbed} is independent of
the time $h$.

The next result shows that UGB is a uniform in time strengthening of the BRS property.
\begin{proposition}
A system $\Sigma=(X,\Uc,\phi)$ is UGB if and only if $\Sigma$ is BRS and a continuous increasing $\mu$ in Lemma~\ref{lem:Boundedness_Reachability_Sets_criterion} can be chosen so that $\mu$ does not depend on $t$.
\end{proposition}

\begin{proof}
If $\Sigma$ is UGB, then $\Sigma$ is BRS and $\mu$ in Lemma~\ref{lem:Boundedness_Reachability_Sets_criterion} can be chosen as
\[
\mu(\|x\|_X,\|u\|_{\Uc}) := \sigma(\|x\|_X) +\gamma(\|u\|_{\Uc}) + c.
\]
Conversely, let Lemma~\ref{lem:Boundedness_Reachability_Sets_criterion} holds with a continuous increasing $\mu=\mu(\|x\|_X,\|u\|_{\Uc})$. Then for any $t\ge0$, $x\in X$, $u\in\Uc$ it holds that
\[
\|\phi(t,x,u)\|_X\leq \max\big\{\mu(\|x\|_X,\|x\|_X),\mu(\|u\|_\Uc,\|u\|_\Uc)\big\}.
\]
Then $\tilde\sigma:r\mapsto \mu(r,r)$ is a continuous increasing function. Define $\sigma(r):=\tilde\sigma(r) - \tilde\sigma(0)$. Then $\sigma\in\Kinf$ and we have for any $t\geq 0$, $x\in X$, $u\in\Uc$:
\[
\|\phi(t,x,u)\|_X\leq \sigma(\|x\|_X) + \sigma(\|u\|_\Uc) + \tilde\sigma(0),
\]
which shows UGB of $\Sigma$.
\end{proof}

Finally, we can characterize uniform global stability.
\begin{lemma}
\label{lem:LS_plus_UGB_equals_GS}
Consider a forward complete system $\Sigma=(X,\Uc,\phi)$. Then
$\Sigma$ is ULS and UGB if and only if it is UGS.
\end{lemma}

\begin{proof}
Assume that $\Sigma$ is ULS and UGB. This means, that there exist $\sigma_1,\gamma_1,\sigma_2,\gamma_2 \in \Kinf$ and $r,c>0$ so that for all $x: \|x\|_X\leq r$ and all $u: \|u\|_{\Uc} \leq r$ it holds that
\begin{equation*}
\left\| \phi(t,x,u) \right\|_X \leq \sigma_1(\|x\|_X) +\gamma_1(\|u\|_{\Uc})
\end{equation*}
and for all $x \in X$ and all $u \in \Uc$ the following estimate holds:
\begin{equation*}
\left\| \phi(t,x,u) \right\|_X \leq \sigma_2(\|x\|_X) +\gamma_2(\|u\|_{\Uc}) + c.
\end{equation*}
Assume without restriction that $\sigma_2 (s) \geq \sigma_1(s)$ and $\gamma_2 (s) \geq \gamma_1(s)$ for all $s \geq 0$.

Pick $k_1,k_2>0$ so that $c=k_1 \sigma_2(r)$ and $c=k_2 \gamma_2(r)$. Then
\[
c \leq c+c \leq  k_1 \sigma_2(\|x\|_X) + k_2 \gamma_2(\|u\|_{\Uc})
\]
for any $x \in X$ and any $u \in \Uc$ so that either $\|x\|_X \geq r$ or $\|u\|_{\Uc} \geq r$.
Thus for all $x \in X$ and all $u \in \Uc$ it holds that
\begin{equation*}
\left\| \phi(t,x,u) \right\|_X \leq (1+k_1)\sigma_2(\|x\|_X) + (1+k_2)\gamma_2(\|u\|_{\Uc}).
\end{equation*}
This shows uniform global stability of $\Sigma$.
\end{proof}

\subsection{Attractivity properties for systems with inputs}

We define the attractivity properties for systems with inputs.
\begin{definition}
\label{def:asymptotic_gain}
\index{AG}
\index{property!asymptotic gain}
A forward complete system $\Sigma=(X,\Uc,\phi)$ has the
\begin{itemize}
    \item[(i)] \emph{asymptotic gain (AG) property}, if there is a $
          \gamma \in \Kinf  \cup \{0\}$ such that for all $\eps >0$, for
          all $x \in X$ and for all $u \in \Uc$ there exists a
          $\tau=\tau(\eps,x,u) < \infty$ such that
\begin{equation}
\label{AG_Absch}
t \geq \tau\ \qrq \|\phi(t,x,u)\|_X \leq \eps + \gamma(\|u\|_{\Uc}).
\end{equation}

\index{sAG}
\index{property!strong asymptotic gain}
  \item[(ii)] \emph{strong asymptotic gain (sAG) property}, if there is a $
    \gamma \in \Kinf \cup \{0\}$ such that for all $x \in X$ and for all
    $\eps >0$ there exists a $\tau=\tau(\eps,x) < \infty $ such that
\begin{equation}
\label{sAG_Absch}
u \in \Uc,\quad t \geq \tau \qrq \|\phi(t,x,u)\|_X \leq \eps + \gamma(\|u\|_{\Uc}).
\end{equation}

\index{UAG}
\index{property! uniform asymptotic gain}
    \item[(iii)] \emph{uniform asymptotic gain (UAG) property}, if there
          exists a
          $ \gamma \in \Kinf \cup \{0\}$ such that for all $ \eps, r
          >0$ there is a $ \tau=\tau(\eps,r) < \infty$ such
          that for all $u \in \Uc$ and all $x \in B_{r}$
\begin{equation}
\label{UAG_Absch}
u \in \Uc,\quad x \in B_{r},\quad t \geq \tau \qrq \|\phi(t,x,u)\|_X \leq \eps + \gamma(\|u\|_{\Uc}).
\end{equation}

\index{bUAG}
\index{property!bounded input uniform asymptotic gain}
\item[(iv)] \emph{bounded input uniform asymptotic gain (bUAG) property}, if there
          exists a
          $ \gamma \in \Kinf \cup \{0\}$ such that for all $ \eps, r
          >0$ there is a $ \tau=\tau(\eps,r) < \infty$ such
          that
\begin{equation}
\label{eq:bUAG_Absch}
\|u\|_{\Uc}\leq r,\quad x \in B_{r},\quad t \geq \tau \qrq \|\phi(t,x,u)\|_X \leq \eps + \gamma(\|u\|_{\Uc}).
\end{equation}
\end{itemize}

\end{definition}

All three properties AG, sAG, and UAG imply that all trajectories converge to the ball of radius $\gamma(\|u\|_{\Uc})$ around the origin as $t \to \infty$.
The difference between AG, sAG, and UAG is in the kind of dependence of
$\tau$ on the states and inputs.
In UAG systems, this time depends (besides $\eps$) only on the norm of the state. In sAG systems, it depends on the state $x$ (and may vary for different states with the same norm), but it does not depend on $u$. In AG systems, $\tau$ depends on $x$ and $u$.
For systems without inputs, the AG and sAG properties are reduced to
0-GATT, and the UAG property becomes 0-UGATT.

Clearly, UAG property implies bUAG property. It is not completely clear, whether the converse holds without further assumptions, as, e.g., BRS property.
However, for a system satisfying the bUAG property, it is not generally possible to verify UAG (and even sAG) property without an increase of the gain, as explained in the following example.
\begin{example}
\label{examp:AG_UAG}
    Consider the system
        \begin{equation}
        \label{Ex:AG_UAG_Difference}
        \dot{x}(t)=-\frac{1}{1+|u(t)|} x(t),
        \end{equation}
where $u\in \Uc:=L^{\infty}(\R_+,\R)$ is a globally bounded input and $x(t)\in\R$ is the state.

As follows from the analysis in \cite[Example 2.46]{Mir23}, this system is forward complete, has BRS property, has bUAG property with a zero gain, has UAG property with arbitrary small linear gain, and is ISS with arbitrary small linear gain. However, it is not sAG with a zero gain.
\end{example}

Next, we define properties similar to AG, sAG, and UAG, which formalize
reachability of the $\eps$-neighborhood of the ball $B_{\gamma(\|u\|_{\Uc})}$ by trajectories of $\Sigma$.
\begin{definition}
\label{def:Limit-properties}
We say that $\Sigma=(X,\Uc,\phi)$ has the
\begin{itemize}
\index{LIM}
\index{property!limit}
    \item[(i)] \emph{limit property (LIM)} if there exists
          $\gamma\in\Kinf\cup\{0\}$ such that
for all $x\in X$, $u \in \Uc$ and $\eps>0$ there is a $t=t(x,u,\eps)$:
\[
\|\phi(t,x,u)\|_X \leq \eps + \gamma(\|u\|_{\Uc}).
\]
\index{sLIM}
\index{property!strong limit}
  \item[(ii)] \emph{strong limit property (sLIM)}, if there exists $\gamma\in\Kinf\cup\{0\}$ so that for every $\eps>0$
and for every $x\in X$ there exists $\tau = \tau(\eps,x)$ such that
\begin{eqnarray}
u\in\Uc \qrq \exists t\leq\tau:\  \|\phi(t,x,u)\|_X \leq \eps + \gamma(\|u\|_{\Uc}).
\label{eq:sLIM_ISS_section}
\end{eqnarray}

\index{ULIM}
\index{property!uniform limit}
  \item[(iii)] \emph{uniform limit property (ULIM)}, if there exists
    $\gamma\in\Kinf\cup\{0\}$ so that for every $\eps>0$ and for every $r>0$ there
    exists a $\tau = \tau(\eps,r)$ such that
\begin{eqnarray}
\|x\|_X \leq r \ \wedge \ u\in\Uc \qrq \exists t\leq\tau:\  \|\phi(t,x,u)\|_X \leq \eps + \gamma(\|u\|_{\Uc}).
\label{eq:ULIM_ISS_section}
\end{eqnarray}

\index{bULIM}
\index{property!uniform limit on bounded sets}
	\item[(iv)] \emph{uniform limit property on bounded sets (bULIM)}, if there exists
    $\gamma\in\Kinf\cup\{0\}$ so that for every $\eps>0$ and for every $r>0$ there is a $\tau = \tau(\eps,r)$ such that
\begin{eqnarray}
\|x\|_X \leq r \ \wedge \ \|u\|_\Uc \leq r \qrq \exists t\leq\tau:\  \|\phi(t,x,u)\|_X \leq \eps + \gamma(\|u\|_{\Uc}).
\label{eq:bULIM_ISS_section}
\end{eqnarray}
\end{itemize}
\end{definition}

\begin{remark}
\label{rem:LIMAG}
It is easy to see that AG is equivalent to the existence of a
$\gamma\in\Kinf$ for which
\[
x\in X \ \wedge \ u\in\Uc \qrq \mathop{\overline{\lim}}_{t\to\infty}\|\phi(t,x,u)\|_X\leq\gamma(\|u\|_{\Uc}),
\]
and LIM is equivalent to the existence of a $\gamma\in\Kinf$ so that
\[
x\in X \ \wedge \ u\in\Uc \qrq \inf_{t\geq0}\|\phi(t,x,u)\|_X\leq\gamma(\|u\|_{\Uc}).
\]
In particular, AG implies LIM. On the other hand, it is easy to see that LIM and UGS together imply AG.
\end{remark}

\begin{remark}
For systems without inputs, the notions of sLIM and LIM coincide and for infinite-dimensional systems, they are
strictly weaker than the ULIM, even for linear systems generated by $C_0$-semigroups (in this case ULIM is equivalent to exponential stability of the semigroup, and LIM is at least not stronger than strong stability of a semigroup), see \cite[Proposition 5.1]{Mir17a}.
\end{remark}

\begin{remark}
\label{rem:ODEs_LIM_and_ULIM}
 For ordinary differential equations (under natural regularity assumptions), the notions of LIM, sLIM, and ULIM coincide; see Proposition~\ref{prop:ULIM_equals_LIM_in_finite_dimensions}.
\end{remark}

\section{ISS superposition theorem}
\label{sec:Motivation_MainResult}

\begin{figure*}[tbh]
\centering
\begin{tikzpicture}[>=implies,thick]
\small

\node (ISS) at (2.7,0.4) {ISS};
\node (bUAG) at (0.8,0.4) {bUAG};

 
\draw [rounded corners] (-6.55,1.2) rectangle (3.1,0.1);
\node (Thm5) at (-1.5,0.9) {\footnotesize Thm.~\ref{thm:UAG_equals_ULIM_plus_LS}};


 \node (Thm8) at (-3.2,-1) {\footnotesize Thm.~\ref{wISS_equals_sAG_GS}};


\node (Rem) at (-4,-2.35) {\footnotesize Rem.~\ref{rem:LIMAG}};

\draw [rounded corners] (-6.55,-1.8) rectangle (-0.2,-0.7);

\draw[thick,double equal sign distance,<->] (ISS) to (bUAG);

\node (bULIM_UGS) at (-2,0.4) {bULIM$\,\wedge\,$UGS};
\node (bULIM_ULS) at (-5.2,0.4) {bULIM$\,\wedge\,$ULS};

\draw[thick,double equal sign distance,<->] (bULIM_UGS) to (bUAG);
\draw[thick,double equal sign distance,<->] (bULIM_ULS) to (bULIM_UGS);

\node (sISS) at (-0.6,-1.5) {sISS};
\node (sAG_UGS) at (-2.7,-1.5) {sAG$\,\wedge\,$UGS};
\node (sLIM_UGS) at (-5.4,-1.5) {sLIM$\,\wedge\,$UGS};
\draw[<->,double equal sign distance] (sISS) to (sAG_UGS);

\draw[<->,double equal sign distance] (sLIM_UGS) to (sAG_UGS);

\draw[thick,double equal sign distance,->] (-3.3,0) to (-3.3,-0.6);
\draw[red,double equal sign distance,->,degil]  (-2.95,-0.6) to (-2.95,0);
\node[red] (notImp_ix) at (-2.45,-0.3) {\footnotesize(ix)};

\draw[red,double equal sign distance,->,degil]  (0.7,-2.45) to (0,-2.45);
\node[red] (vi) at (1.1,-2.45) {\footnotesize(vi)};

\draw[red,double equal sign distance,<-,degil]  (0.7,-2.95) to (0,-2.95);
\node[red] (vii) at (1.1,-2.95) {\footnotesize(vii)};

\draw[red,double equal sign distance,->,degil]  (0.7,-1.3) to (0,-1.3);
\node[red] (iv) at (1.1,-1.3) {\footnotesize(iv)};

\draw[red,double equal sign distance,<-,degil]  (0.7,-1.7) to (0,-1.7);
\node[red] (v) at (1.1,-1.7) {\footnotesize(v)};



\node (AG_UGS) at (-2.7,-2.75) {AG$\,\wedge\,$UGS};
\node (LIM_UGS) at (-5.4,-2.75) {LIM$\,\wedge\,$UGS};

\node (AG_ULS) at (-2.7,-3.9) {AG$\,\wedge\,$ULS};

\node (AG_0UGAS) at (2.7,-1.5) {AG$\,\wedge\,$0-UGAS};
\node (AG_0UAS) at (2.7,-2.75) {AG$\,\wedge\,$0-UAS};
\node (AG_0LS) at (2.7,-3.9) {AG$\,\wedge\,$0-ULS};
\node (AG_0GAS) at (5.5,-3.9) {AG$\,\wedge\,$0-GAS};


\draw[double equal sign distance,->]  ($(AG_UGS)+(-0.2,-0.3)$) to ($(AG_ULS)+(-0.2,0.3)$);
\draw[red,double equal sign distance,<-,degil]  ($(AG_UGS)+(0.2,-0.3)$) to ($(AG_ULS)+(0.2,0.3)$);
\node[red] (viii) at ($\weight*(AG_ULS)+\weight*(AG_UGS) + (0.75,0)$) {\footnotesize(viii)};

\draw[double equal sign distance,->]  ($(AG_0UGAS)+(-0.2,-0.3)$) to ($(AG_0UAS)+(-0.2,0.3)$);
\draw[red,double equal sign distance,<-,degil]  ($(AG_0UGAS)+(0.2,-0.3)$) to ($(AG_0UAS)+(0.2,0.3)$);
\node[red] (ii) at ($\weight*(AG_0UAS)+\weight*(AG_0UGAS) + (0.6,0)$) {\footnotesize(ii)};

\draw[double equal sign distance,->]  ($(AG_0UAS)+(-0.2,-0.3)$) to ($(AG_0LS)+(-0.2,0.3)$);
\draw[red,double equal sign distance,<-,degil]  ($(AG_0UAS)+(0.2,-0.3)$) to ($(AG_0LS)+(0.2,0.3)$);
\node[red] (i) at ($\weight*(AG_0LS)+\weight*(AG_0UAS) + (0.6,0)$) {\footnotesize(i)};

%

\coordinate (ISS_1) at ($(ISS) + (-0.5,0)$);
\coordinate (ISS_2) at ($(ISS) + (0,-0.4)$);
\coordinate (AG_0UGAS_1) at ($(AG_0UGAS) + (-0.5,0)$);
\coordinate (coord1) at ($(AG_0UGAS_1) + (0,0.3)$);
\coordinate (coord2) at ($(ISS_1) + (0,-0.4)$);

\draw[red,double equal sign distance,->,degil]  (coord1) to (coord2);
\node[red] (iii) at ($\weight*(ISS_1) + \weight*(AG_0UGAS_1) + (-0.45,-0.15)$) {\footnotesize(iii)};

\draw[thick,double equal sign distance,->] ($(sAG_UGS)+(-0.2,-0.37)$) to ($(AG_UGS)+(-0.2,0.27)$);

  \coordinate (A2) at ($(sAG_UGS)+(0.2,-0.37)$);
  \coordinate (B2) at ($(AG_UGS)+(0.2,0.27)$);
\node[red] (Quest) at (-2.5,-2.23) {\footnotesize ???};
\draw[double equal sign distance,-] (B2) to ($(B2)+(0,0.09)$);
\draw[double equal sign distance,->] ($(B2)+(0,0.42)$) to (A2);

\draw[double equal sign distance,->]  ($(AG_UGS)+(-0.2,-0.3)$) to ($(AG_ULS)+(-0.2,0.3)$);

\draw[thick,double equal sign distance,<->] (AG_UGS) to (LIM_UGS);

\draw[thick,double equal sign distance,->] (ISS_2) to (AG_0UGAS);

\draw[thick,double equal sign distance,<->]  (AG_0LS) to (AG_0GAS);

\draw[thick,double equal sign distance,->] ($(AG_ULS)+(0.9,-0.12)$)  to ($(AG_0LS)+(-1.1,-0.12)$);

  \coordinate (A1) at ($(AG_ULS)+(0.9,0.12)$);
  \coordinate (B1) at ($(AG_0LS)+(-1.1,0.12)$);
\path (B1) -- node (Q) {\footnotesize\color{red}{???}} (A1);
    \draw[thick,double equal sign distance,->] (B1) -- (Q) -- (A1);


\end{tikzpicture}

\caption[caption]{Relations between stability properties of
  infinite-dimensional systems, which have bounded reachability sets and whose flow is continuous at trivial equilibrium:
\begin{itemize}
\setlength{\itemindent}{5mm}
   \item Black arrows show implications or equivalences which hold for general control systems
in infinite dimensions.
\item {\color{red}Red arrows (with the negation sign)}
are implications which do not hold, and for which the counterexamples exist, see Remark~\ref{rem:Nonimplications}.
\item Black arrows with question marks inside mean that it is not known right now (as far as the author is concerned), whether these converse implications hold or not.
\end{itemize}
}
\label{fig:ISS_Equiv}
\end{figure*}

The central result in this chapter is the following theorem:
\index{theorem!ISS Superposition}
\begin{theorem}[ISS superposition theorem]
\label{thm:MainResult_Characterization_ISS}
Let $\Sigma=(X,\Uc,\phi)$ be a forward complete system satisfying
the BRS and the CEP property.
Then the relations depicted in Figure~\ref{fig:ISS_Equiv} hold.

Black arrows show implications or equivalences which hold for the class of
infinite-dimensional systems defined in Definition~\ref{Steurungssystem};
the red arrows (with the negation sign) are  implications
which do not hold due to the counterexamples presented in this chapter;
the arrows with question marks indicate that it is not known to us whether these implications hold.
\end{theorem}

\begin{proof}
The result follows from Theorems~\ref{thm:UAG_equals_ULIM_plus_LS},
\ref{wISS_equals_sAG_GS}. The counterexamples for the red arrows are
discussed in Section~\ref{sec:Counterexamples}, see Remark~\ref{rem:Nonimplications}, where an overview of the counterexamples is provided.
\end{proof}

\begin{remark}
\label{rem:ISS_for_ODEs}
For ODEs with a sufficiently regular right-hand side, all the combinations in Figure~\ref{fig:ISS_Equiv} are
equivalent,
as for ODE systems AG\,$\wedge$\,0-GAS is equivalent to ISS by Proposition~\ref{Characterizations_ODEs}.
In contrast, for infinite-dimensional systems, these notions are divided into several groups, which are not equivalent to each other.
\end{remark}

The proof of Theorem~\ref{thm:MainResult_Characterization_ISS} will be given in several steps.

\begin{enumerate}[label=(\roman*)]
	\item First, in Theorem~\ref{thm:UAG_equals_ULIM_plus_LS}, we characterize the (uniform) ISS property as a combination of ULIM, BRS, and ULS properties. This implies the equivalences in the upper level in Figure~\ref{fig:ISS_Equiv}.
	\item In Section~\ref{sec:WeakISS}, we introduce the concept of strong input-to-state stability (sISS).
For nonlinear ODE systems, this is equivalent to ISS, see Proposition~\ref{prop:ULIM_equals_LIM_in_finite_dimensions}, and for linear infinite-dimensional systems without inputs, sISS is equivalent to strong stability of the associated semigroup $T$ (which justifies the term "strong" for this notion).
We show in Theorem~\ref{wISS_equals_sAG_GS} that
\begin{center}
sISS \quad $\Iff$\quad sAG $\wedge$ UGS \quad $\Iff$\quad sLIM $\wedge$ UGS.
\end{center}
	\item On the other hand, ISS implies the combination AG $\wedge$ 0-UGAS, which is very different from sISS: sISS does not imply the existence of a uniform
convergence rate for the undisturbed system (and thus, it does not imply
0-UGAS). At the same time, AG $\wedge$ 0-UGAS does not ensure the uniform global stability property for a system with inputs.
	\item Below the level of sISS and AG $\wedge$ 0-UGAS, there are further levels with even weaker properties.
The counterexamples, discussing delicate properties of
infinite-dimensional systems and giving the necessary counterexamples for Figure~\ref{fig:ISS_Equiv}, are discussed in detail in Section~\ref{sec:Counterexamples}.
\end{enumerate}

Our main result in this section is the following characterization of ISS for a general class of control systems.
\begin{theorem}[ISS superposition theorem]
\label{thm:UAG_equals_ULIM_plus_LS}
Let $\Sigma=(X,\Uc,\phi)$ be a forward complete control system. The following statements are equivalent:
\begin{enumerate}[label=(\roman*)]
    \item\label{itm:ISS-Characterization-bounded-properties-1} $\Sigma$ is ISS.
    \item\label{itm:ISS-Characterization-bounded-properties-2} $\Sigma$ is UAG and UGS.
    \item\label{itm:ISS-Characterization-bounded-properties-3} $\Sigma$ is bUAG and UGS.
	  \item\label{itm:ISS-Characterization-bounded-properties-4} $\Sigma$ is bUAG, CEP and BRS.
    \item\label{itm:ISS-Characterization-bounded-properties-5} $\Sigma$ is bULIM, ULS and BRS.
    \item\label{itm:ISS-Characterization-bounded-properties-6} $\Sigma$ is bULIM and UGS.
\end{enumerate}
\end{theorem}

\begin{proof}
\ref{itm:ISS-Characterization-bounded-properties-1} $\Rightarrow$ \ref{itm:ISS-Characterization-bounded-properties-2}.
Estimating $\beta(\|x\|_X,t) \leq \beta(\|x\|_X,0)$ for all $x\in X$ and all $t\geq 0$, we see that ISS implies UGS. The UAG property follows by Lemma~\ref{ISS_implies_UAG}.

\ref{itm:ISS-Characterization-bounded-properties-2} $\Rightarrow$ \ref{itm:ISS-Characterization-bounded-properties-3} $\Rightarrow$ \ref{itm:ISS-Characterization-bounded-properties-4}. Clear.

\ref{itm:ISS-Characterization-bounded-properties-4} $\Rightarrow$ \ref{itm:ISS-Characterization-bounded-properties-5}.
A combination bUAG $\wedge$ CEP implies ULS by Lemma~\ref{UAG-ULS}.

\ref{itm:ISS-Characterization-bounded-properties-5} $\Rightarrow$ \ref{itm:ISS-Characterization-bounded-properties-6}.
By Proposition~\ref{prop:ULIM_plus_mildRFC_implies_UGB}, bULIM $\wedge$ BRS property implies UGB property.
UGB $\wedge$ ULS properties imply UGS by Lemma~\ref{lem:LS_plus_UGB_equals_GS}.

\ref{itm:ISS-Characterization-bounded-properties-6} $\Rightarrow$ \ref{itm:ISS-Characterization-bounded-properties-1}.
bULIM $\wedge$ UGS implies bUAG by Lemma~\ref{lem:ULIM_plus_GS_implies_UAG}.
bUAG $\wedge$ UGS implies UAG by Lemma~\ref{lem:UGS_und_bUAG_imply_UAG}.
UAG $\wedge$ UGS implies ISS by Lemma~\ref{lem:UAG_implies_ISS}.
\end{proof}

Next, we show a series of auxiliary results, which we used in the proof of Theorem~\ref{thm:UAG_equals_ULIM_plus_LS}.

\begin{lemma}
\label{ISS_implies_UAG}
If a system $\Sigma=(X,\Uc,\phi)$ is ISS, then it is UAG.
\end{lemma}

\begin{proof}
Let $\Sigma$ be ISS with the corresponding $\beta\in\KL$ and $\gamma\in\Kinf$.
Take arbitrary $\eps, r >0$.  Define $\tau=\tau(\eps,r)$ as
the solution of the equation $\beta(r,\tau)=\eps$ (if it exists,
then it is unique because of the monotonicity of $\beta$ in the second
argument, if it does not exist, we set $\tau(\eps,r):=0$). Then for
all $t \geq \tau$, all $x\in X$ with $\|x\|_X \leq r$ and all $u \in \Uc$ we have
\begin{eqnarray*}
\|\phi(t,x,u)\|_X &\leq& \beta(\|x\|_X,t) + \gamma(\|u\|_{\Uc})  \\
&\leq& \beta(\|x\|_X,\tau) + \gamma(\|u\|_{\Uc}) \\
&\leq& \eps + \gamma(\|u\|_{\Uc}),
\end{eqnarray*}
and the implication \eqref{UAG_Absch} holds.
\end{proof}

\begin{lemma}
\label{UAG-ULS}
If a system $\Sigma=(X,\Uc,\phi)$ is bUAG and CEP, then it is ULS.
\end{lemma}

\begin{proof}
    We will show that \eqref{LS_Restatement} holds so that the claim
    follows from Lemma~\ref{lem:ULS_restatement}.
		Let $\tau$ and $\gamma$ be the
    functions from the bUAG definition. Let $\eps>0$ and $\tau:=\tau(\eps/2,1)$.
    Pick any $\delta_1\in (0,1]$ so that
    $\gamma(\delta_1)<\eps/2$. Then for all $x \in X$ with $\|x\|_X \leq 1$
    and all $u \in \Uc$ with $\|u\|_{\Uc} \leq\delta_1$ we have
\begin{equation}\label{eq:2}
\sup_{t\geq \tau}\|\phi(t,x,u)\|_X \leq\frac{\eps}{2}+\gamma(\|u\|_{\Uc})<\eps.
\end{equation}

Since $\Sigma$ is CEP, there is some $\delta_2=\delta_2(\eps,\tau)>0$ so that
\[
\|\eta\|_X \leq\delta_2 \ \wedge \  \|u\|_{\Uc} \leq \delta_2 \quad \Rightarrow \quad \sup_{t\in[0,\tau]}\|\phi(t,\eta,u)\|_X \leq\eps.
\]
Together with \eqref{eq:2}, this proves  \eqref{LS_Restatement} with $\delta:=\min\lbrace 1,\delta_1,\delta_2 \rbrace$.
\end{proof}

We proceed with
\begin{proposition}
\label{prop:ULIM_plus_mildRFC_implies_UGB}
Assume that a system $\Sigma=(X,\Uc,\phi)$ is BRS and has the bULIM property. Then $\Sigma$ is UGB.
\end{proposition}

\begin{proof}
Pick any $r>0$ and set $\eps:=\frac{r}{2}$.
Since $\Sigma$ has the bULIM property, there exists a $\tau=\tau(r)$ (more precisely, $\tau=\tau(\frac{r}{2}, \max\{r,\gamma^{-1}(\tfrac{r}{4})\})$ from bULIM property), so that
\begin{eqnarray}
x \in\clo{B_r},\, \|u\|_{\Uc} \leq \gamma^{-1}(\tfrac{r}{4})
\ \Rightarrow \ \exists t\leq \tau:\ \|\phi(t,x,u)\|_X \leq \frac{r}{2} +\gamma(\|u\|_{\Uc})  \leq \frac{3r}{4}.
\label{eq:ULIM_pGS_ISS_Section2}
\end{eqnarray}
Without loss of generality, we can assume that $\tau$ is increasing in $r$. In particular, it is locally integrable.
Defining $\bar\tau(r):=\frac{1}{r}\int_r^{2r}\tau(s)ds$ for $r>0$, we see that $\bar\tau(r) \geq \tau(r)$ and $\bar\tau$ is continuous.
For any $r_2>r_1>0$ via the change of variables $s= \frac{r_2}{r_1}w$, we have also that
\begin{eqnarray*}
\bar\tau(r_2) &=& \frac{1}{r_2}\int_{r_2}^{2r_2}\tau(s)ds = \frac{1}{r_2} \int_{r_1}^{2r_1}\tau\Big(\frac{r_2}{r_1}w\Big)\frac{r_2}{r_1}dw \\
&>& \frac{1}{r_1}\int_{r_1}^{2r_1}\tau(w)dw = \bar\tau(r_1),
\end{eqnarray*}
and thus $\bar\tau$ is increasing. We define further $\bar\tau(0):=\lim_{r\to+0}\bar\tau(r)$ (the limit exists as $\bar\tau$ is increasing).

Since $\Sigma$ is BRS, Lemma~\ref{lem:Boundedness_Reachability_Sets_criterion} implies that there exists a continuous, increasing function $\mu: \R_+^3 \to \R_+$, such that for
all $x\in X, u\in \Uc$ and all $t \geq 0$ estimate \eqref{eq:8_ISS} holds. This implies that
\begin{eqnarray}
x\in \clo{B_r},\, \|u\|_{\Uc} \leq \gamma^{-1}(\tfrac{r}{4}),\, t \leq \bar\tau(r) \ \Rightarrow\ \|\phi(t,x,u)\|_X \leq \tilde\sigma(r),
\label{eq:ULIM_pGS_ISS_Section3}
\end{eqnarray}
where $\tilde\sigma:r\mapsto \mu\big(r,\gamma^{-1}(\tfrac{r}{4}),\bar\tau(r)\big)$, $r\geq 0$
is continuous and increasing, since $\mu,\gamma,\bar\tau$ are continuous increasing functions.

Using identity axiom ($\Sigma$\ref{axiom:Identity}), we see from \eqref{eq:ULIM_pGS_ISS_Section3} that $\tilde\sigma(r)\geq \|\phi(0,x,u)\|_X = \|x\|_X$ for $x \in \clo{B_r}$, and thus $\tilde\sigma(r) \geq r > \frac{3r}{4}$ for any $r>0$.

Assume that there exist  $x\in \clo{B_r}$, $u\in\Uc$ with $\|u\|_{\Uc} \leq \gamma^{-1}(\tfrac{r}{4})$ and $t\geq0$ so that $\|\phi(t,x,u)\|_X  > \tilde \sigma(r)$. Define
\[
t_m:=\sup\{s \in[0,t]: \|\phi(s,x,u)\|_X \leq r\}\geq 0.
\]
The quantity $t_m$ is well-defined, since $\|\phi(0,x,u)\|_X = \|x\|_X
\leq r$ due to the identity property ($\Sigma$\ref{axiom:Identity}).

In view of the cocycle property ($\Sigma$\ref{axiom:Cocycle}), it holds that
\[
\phi(t,x,u) = \phi\big(t-t_m,\phi(t_m,x,u),u(\cdot + t_m)\big),
\]
and $u(\cdot + t_m) \in\Uc$, since $\Sigma$ satisfies the axiom of shift invariance.
Assume that $t-t_m \leq\tau(r)$. Since $\|\phi(t_m,x,u)\|_X\leq r$,
\eqref{eq:ULIM_pGS_ISS_Section3} implies that $\|\phi(t,x,u)\|_X \leq \tilde \sigma(r)$ for all $t \in [t_m,t]$.
Otherwise, if $t-t_m >\tau(r)$, then due to \eqref{eq:ULIM_pGS_ISS_Section2} there exists $t^* < \tau(r)$, so that
\[
\big\|\phi\big(t^*,\phi(t_m,x,u),u(\cdot + t_m)\big)\big\|_X = \|\phi(t^*+t_m,x,u)\|_X \leq \frac{3r}{4},
\]
which contradicts the definition of $t_m$, since $t_m+t^*<t$.
Hence
\begin{eqnarray}
x \in\clo{B_r},\, \|u\|_{\Uc} \leq \gamma^{-1}(\tfrac{r}{4}),\, t\geq
0 \ \Rightarrow \
  \|\phi(t,x,u)\|_X \leq \tilde\sigma(r).
\label{eq:ULIM_pGS_ISS_Section4}
\end{eqnarray}
Denote $\sigma(r):=\tilde\sigma(r) - \tilde\sigma(0)$, for any $r\geq 0$. Clearly, $\sigma\in\Kinf$.

For each $x\in X$, $u\in\Uc$ define $r:=\max\{\|x\|_X,4\gamma(\|u\|_{\Uc})\}$.
Then \eqref{eq:ULIM_pGS_ISS_Section4} immediately shows for all $x\in X,\ u\in\Uc,\ t\geq 0$ that
\begin{eqnarray*}
\|\phi(t,x,u)\|_X &\leq& \sigma\big(\max\{\|x\|_X,4\gamma(\|u\|_{\Uc})\}\big) + \tilde\sigma(0) \\
&\leq& \sigma(\|x\|_X) + \sigma\big(4\gamma(\|u\|_{\Uc})\big) + \tilde\sigma(0),
\end{eqnarray*}
which shows UGB of $\Sigma$.
\end{proof}

\begin{lemma}
\label{lem:ULIM_plus_GS_implies_UAG}
If a system $\Sigma=(X,\Uc,\phi)$ is bULIM and UGS, then $\Sigma$ is bUAG.
\end{lemma}

\begin{proof}
Without loss of generality, assume that $\gamma$ in the definitions of ULIM
and UGS is the same (otherwise, take the maximum of the two).

Pick any $\eps>0$ and any $r>0$. By the uniform limit property, there exists $\gamma\in\Kinf$, independent of $\eps$ and $r$,
and $\tau=\tau(\eps,r)$ so that \eqref{eq:bULIM_ISS_section} holds.

In view of the cocycle property, we have from the UGS property that for the above $x,u,t$ and any $s\geq 0$
\begin{eqnarray*}
\|\phi(s+t,x,u)\|_X &=& \big\|\phi(s,\phi(t,x,u),u(s+\cdot))\big\|_X \\
            &\leq& \sigma\big(\|\phi(t,x,u)\|_X\big) + \gamma(\|u\|_{\Uc})\\
            &\leq& \sigma\big(\eps + \gamma(\|u\|_{\Uc})\big) + \gamma(\|u\|_{\Uc}).
\end{eqnarray*}
Now let $\tilde\eps := \sigma(2 \eps)>0$.
Using the evident inequality $\sigma(a+b)\leq \sigma(2a) + \sigma(2b)$, valid for any $a,b\geq 0$, we proceed to
\begin{eqnarray*}
\|\phi(s+t,x,u)\|_X \leq \tilde\eps + \tilde\gamma(\|u\|_{\Uc}),
\end{eqnarray*}
where $\tilde\gamma(r) = \sigma(2\gamma(r)) + \gamma(r)$.

Overall, for any $\tilde\eps>0$ and any $r>0$ there exists $\tau=\tau(\eps,r) = \tau(\frac{1}{2}\sigma^{-1}(\tilde\eps),r)$,
so that
\begin{eqnarray*}
\|x\|_X\leq r\ \wedge\ \|u\|_{\Uc}\leq r \ \wedge\ t\geq \tau \qrq \|\phi(t,x,u)\|_X \leq \tilde\eps + \tilde\gamma(\|u\|_{\Uc}),
\end{eqnarray*}
which shows the bUAG property of $\Sigma$.
\end{proof}

The following lemma shows how the bUAG property can be \q{upgraded} to the UAG property.
\begin{lemma}
\label{lem:UGS_und_bUAG_imply_UAG}
If a system $\Sigma=(X,\Uc,\phi)$ is a UGS and bUAG control system, then $\Sigma$ is UAG.
\end{lemma}

\begin{proof}
Pick any $\eps>0$, $r>0$, and let $\tau$ and $\gamma$ be as in the formulation of the bUAG property.
Let $x \in B_r$ and let $u\in\Uc$ be arbitrary.
If $\|u\|_\Uc\leq r$, then \eqref{UAG_Absch} is the desired implication.

Let $\|u\|_\Uc > r$. Hence it holds that $\|u\|_\Uc > \|x\|_X$.
Due to uniform global stability of $\Sigma$, it holds for all $t, x, u$ that
\[
\|\phi(t,x,u)\|_X \leq \sigma(\|x\|_X) + \gamma(\|u\|_{\Uc}),
\]
where we assume that $\gamma$ is the same as in the definition of a bUAG property (otherwise pick the maximum of both).
For $\|u\|_\Uc > \|x\|_X$ we obtain that
\[
\|\phi(t,x,u)\|_X \leq \sigma(\|u\|_\Uc) + \gamma(\|u\|_{\Uc}),
\]
and thus for all $x \in X$, $u\in \Uc$ it holds that
\begin{equation*}
t \geq \tau \quad \Rightarrow \quad \|\phi(t,x,u)\|_X \leq \eps + \gamma(\|u\|_{\Uc}) +\sigma(\|u\|_\Uc),
\end{equation*}
which shows the UAG property with the asymptotic gain $\gamma + \sigma$.
\end{proof}

The final technical lemma of this section is
\begin{lemma}
If a system $\Sigma=(X,\Uc,\phi)$ is UAG and UGS, then $\Sigma$ is ISS.
\label{lem:UAG_implies_ISS}
\end{lemma}

\begin{proof}
Assume that $\Sigma$ is UAG and UGS and that $\gamma$ in \eqref{GSAbschaetzung} and \eqref{UAG_Absch} are the same (otherwise, pick $\gamma$ as a maximum of both of them).
Fix arbitrary $r \in \R_+$. We are going to
construct a function $\beta \in \KL$ so that \eqref{iss_sum} holds.

From global stability, it follows that there exist $\gamma,\sigma \in \Kinf$ such
that for all $ t\geq 0$, all $x \in \clo{B_r}$ and all $ u \in \Uc$ we have
\begin{equation}
    \label{lem7:help}
\|\phi(t,x,u)\|_X \leq \sigma(r) + \gamma(\|u\|_{\Uc}).
\end{equation}
Define $\eps_n:= 2^{-n}  \sigma(r)$, for $n \in \N$. The UAG property implies that there exists a sequence of times
$\tau_n:=\tau(\eps_n,r)$, which we may without loss of generality assume
to be strictly increasing, such that for all $x \in \clo{B_r}$ and all $u \in \Uc$
\[
\|\phi(t,x,u)\|_X \leq \eps_n + \gamma(\|u\|_{\Uc})\quad  \forall t \geq \tau_n.
\]
From \eqref{lem7:help}, we see that we may set $\tau_0 := 0$.
Define $\omega(r,\tau_n):=\eps_{n-1}$, for $n \in \N$ and $\omega(r,0):=2\eps_0=2\sigma(r)$.

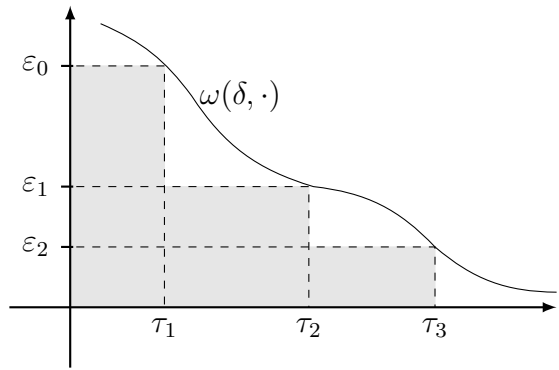
\begin{figure}[ht]
\centering
\begin{tikzpicture}[scale=0.8]
\filldraw[gray!20] (0,0) rectangle (1.55,4);
\filldraw[gray!20] (1.55,0) rectangle (3.925,2);
\filldraw[gray!20] (3.925,0) rectangle (6,1);

\draw[thick,->] (-1,0) -- (8,0);
\draw[thick,->] (0,-1) -- (0,5);
\draw[thick] (0.05,4) -- (-0.15,4) node[left] {$\eps_0$};
\draw[thick] (0.05,2) -- (-0.15,2) node[left] {$\eps_1$};
\draw[thick] (0.05,1) -- (-0.15,1) node[left] {$\eps_2$};
\draw (0.5,4.7) to[bend left=15] (2,3.5) to[bend right=20] (4,2) to[bend left=20] (6,1) to[bend right=20] (8,0.25);

\draw[dashed] (0.1,4) -- (1.55,4) -- (1.55,0) node[below] {$\tau_1$};
\draw[dashed] (0.1,2) -- (3.925,2) -- (3.925,0) node[below] {$\tau_2$};
\draw[dashed] (0.1,1) -- (6,1) -- (6,0) node[below] {$\tau_3$};

\node (Omega) at (2.8,3.5) {$\omega(\delta,\cdot)$};
\end{tikzpicture}
\caption{Construction of the function $\omega(\delta,\cdot)$.}
\label{Omega-Konstruktion}
\end{figure}

Now extend the definition of $\omega$ to a function $\omega(r,\cdot) \in \LL$.
We obtain for $t \in (\tau_n,\tau_{n+1})$, $n=0,1,\ldots$ that
\[
\|\phi(t,x,u)\|_X \leq \eps_n + \gamma(\|u\|_{\Uc})< \omega(r,t) + \gamma(\|u\|_{\Uc}),\quad x\in B_r,\quad u\in\Uc.
\]
Doing this for all $r \in \R_+$, we obtain the definition of the function $\omega$.

Now define $\hat \beta(r,t):=\sup_{0 \leq s \leq r}\omega(s,t) \geq
\omega(r,t)$ for $(r,t) \in \R_+^2$. From this definition, it follows that,
for each $t\geq 0$, $\hat\beta(\cdot,t)$ is
increasing in $r$ and $\hat\beta(r,\cdot)$ is decreasing in $t$ for each $r>0$ as
every $\omega(r,\cdot) \in \LL$.
Moreover, for each fixed $t\geq0$, $\hat \beta(r,t) \leq \sup_{0 \leq s \leq r}\omega(s,0)=2\sigma(r)$. This implies that $\hat\beta$ is continuous in the first argument at $r=0$ for any fixed $t\geq0$.

By \cite[Proposition 9]{MiW19b}, $\hat\beta$ can be upper bounded by certain $\tilde{\beta}\in \KL$, and the estimate
\eqref{iss_sum} is satisfied with such a $\beta$.
\end{proof}

\section{Integral characterization of ISS}
\label{sec:Relations between ISS and norm-to-integral ISS}

In this section, we characterize ISS by notions involving integrals of the trajectories, which is interesting in its own right. Still, it will also be instrumental for establishing non-coercive ISS Lyapunov theorems in Section~\ref{sec:Non-coercive ISS Lyapunov functions}.

\begin{definition}
\label{def:norm-i_ISS} 
\index{input-to-state stability!norm-to-integral}
We call a system $\Sigma=(X,\Uc,\phi)$ \emph{norm-to-integral ISS} if there are $\alpha, \psi, \sigma\in\Kinf$ so that
for all $x\in X$, $u\in\Uc$ and $t\geq 0$ it holds that
\begin{equation}
\label{eq:ni-ISS}
\int_0^t \!\alpha\big(\|\phi(s,x,u)\|_X\big) ds \leq \psi(\|x\|_X) + t \sigma(\|u\|_{\Uc}).
\end{equation}
\end{definition}

For the special case of $\Uc:=L^\infty(\R_+,U)$, we introduce one more stability notion.
\begin{definition}
\label{def:i-i_ISS} 
Consider a forward complete control system $\Sigma:=(X,\Uc,\phi)$ with the input space $\Uc:=L^\infty(\R_+,U)$, where $U$ is any normed vector space.

\index{input-to-state stability!integral-to-integral}
We call $\Sigma$ \emph{integral-to-integral ISS} if there are $\alpha, \psi, \sigma\in\Kinf$ so that
for all $x\in X$, $u\in\Uc$ and $t\geq 0$ it holds that
\begin{equation}
\label{eq:ii-ISS}
\int_0^t \!\alpha(\|\phi(s,x,u)\|_X) ds \leq \psi(\|x\|_X) +\! \int_0^t\! \sigma(\|u(s)\|_{U})ds.
\end{equation}
\end{definition}

We start with
\begin{proposition}
\label{prop:ISS-implies-norm-i-ISS} 
If a forward complete control system is ISS, then it is norm-to-integral ISS.
\end{proposition}

\begin{proof}
Let $\Sigma=(X,\Uc,\phi)$ be a forward complete ISS control system with a transient $\beta\in\KL$.
By Sontag's $\KL$-lemma (Corollary~\ref{Sontags_KL_Lemma-2}), there are $\xi_1, \xi_2 \in\Kinf$ so that $\beta(r,t) \leq \xi_1^{-1}(e^{-t}\xi_2(r))$ for all $r,t\in\R_+$.
ISS of $\Sigma$ now implies that there is $\gamma\in\Kinf$ such that the following holds:
\begin {equation*}
\| \phi(t,x,u) \|_{X} \leq \xi_1^{-1}(e^{-t}\xi_2(\| x \|_{X})) + \gamma( \|u\|_{\Uc}),\quad t\geq 0, \ x\in X, \ u\in\Uc.
\end{equation*}
Define $\bar{\xi}(r):=\xi_1(\frac{1}{2}r)$, $r\in\R_+$. Using the inequality $\bar{\xi}(a+b)\leq \bar{\xi}(2a) + \bar{\xi}(2b)$, which is valid for all $a,b\in\R_+$, we obtain that for all $x\in X$, $u\in\Uc$ and $t\geq 0$ it holds that
\begin {equation}
\label{eq:ISS-implies-iiISS-auxiliary}
\bar{\xi}(\| \phi(t,x,u) \|_{X}) \leq e^{-t}\xi_2(\| x \|_{X}) + \xi_1(\gamma( \|u\|_{\Uc})).
\end{equation}
Integrating \eqref{eq:ISS-implies-iiISS-auxiliary}, we see that 
\begin {equation*}
\int_0^t\bar{\xi}(\| \phi(s,x,u) \|_{X}) ds \leq \xi_2(\| x \|_{X}) + t\xi_1 \circ \gamma( \|u\|_{\Uc}),\quad t\geq 0, \ x\in X, \ u\in\Uc.
\end{equation*}
This shows norm-to-integral ISS of $\Sigma$.
\end{proof}

Next, we show that norm-to-integral ISS implies ISS for a class of forward-complete control systems satisfying the  CEP and BRS properties. We start with

\begin{proposition}
\label{prop:ncISS_implies_ULIM} 
Let $\Sigma=(X,\Uc,\phi)$ be a forward complete control system. If $\Sigma$ is  norm-to-integral ISS, then $\Sigma$ is ULIM.
\end{proposition}

\begin{proof}
As $\Sigma$ is norm-to-integral ISS, the estimate \eqref{eq:ni-ISS} holds. 
Now define $\gamma(r):=\alpha^{-1} \big(2\sigma(r)\big)$, $r\in\R_+$ and
$\tau(r,\eps):=2(\psi(r)+1)(\alpha(\eps))^{-1}$ for any $r,\eps>0$.

Assume that $\Sigma$ is not ULIM with these $\gamma$ and $\tau$.
Then there are some $\eps>0$, $r>0$, $x\in \clo{B_r}$ and $u\in\Uc$ so that $\|\phi(t,x,u)\|_X >\eps + \gamma(\|u\|_\Uc)$
for all $t\in[0,\tau(r,\eps)]$.

Using Proposition~\ref{prop:Kfun_lower_estimates}, we have for these $\eps,x,u$ and all $t\in[0,\tau(r,\eps)]$ that:
\begin{align*}
\int_0^t \alpha(\|\phi(s,x,u)\|_X) ds &\geq \int_0^t \alpha\big(\eps + \gamma(\|u\|_\Uc)\big) ds \\
& \geq \int_0^t \frac{1}{2}\alpha(\eps) + \sigma(\|u\|_\Uc) ds \\
& =  \frac{t}{2}\alpha(\eps) + t\sigma(\|u\|_\Uc) .
\end{align*}
In particular, for $t:=\tau(r,\eps)$ we obtain that
\begin{eqnarray}
\label{ncISS_Eq2}
\int_0^{\tau(r,\eps)} \alpha(\|\phi(s,x,u)\|_X) ds \geq\psi(r)+1 + \tau(r,\eps)\sigma(\|u\|_\Uc).
\end{eqnarray}
Combining estimates \eqref{eq:ni-ISS} and \eqref{ncISS_Eq2}, we see that
\begin{eqnarray*}
\psi(r)+1 \leq \psi(\|x\|_X) \leq\psi(r),
\end{eqnarray*}
a contradiction. This shows that  $\Sigma$ is ULIM.
\end{proof}

Next, we provide a sufficient condition for the  ULS property. 
\begin{proposition}
\label{prop:ncISS_plus_CEP_implies_ULS} 
Let $\Sigma=(X,\Uc,\phi)$ be a forward complete control system satisfying the CEP property. If $\Sigma$ is norm-to-integral ISS, then $\Sigma$ is ULS.
\end{proposition}

\begin{proof}
Let $\Sigma$ be norm-to-integral ISS with the corresponding $\alpha,\psi,\sigma$ as in Definition~\ref{def:norm-i_ISS}.

By Lemma~\ref{lem:ULS_restatement}, $\Sigma$ is ULS if and only if for any $\eps>0$ there
is a $\delta>0$ such that 
\begin{equation*}
\|x\|_X \leq\delta\ \wedge \ \|u\|_{\Uc} \leq \delta\ \wedge \  t\geq 0 \;\; \Rightarrow\;\;
\|\phi(t,x,u)\|_X \leq\eps.
\end{equation*}
Seeking a contradiction, assume that $\Sigma$ is not ULS. Then there is $\varepsilon>0$ such that for any $\delta>0$ there are
 $x \in B_\delta$, $u\in B_{\delta,\Uc}$ and $t\geq 0$ such that $\|\phi(t,x,u)\|_X = \varepsilon$.
Then there are sequences $( x_k )$ in $X$, $( u_k)$ in $\Uc$, and $(t_k) \subset \R_+$ such that $x_k \to 0$ as $k \to \infty$, $u_k\to 0$ as $k\to\infty$ and
\begin{equation*}
    \| \phi(t_k,x_k,u_k) \|_X = \varepsilon,\quad k \in \N.
\end{equation*}

Since $\Sigma$ is CEP, for the above $\eps$ there is a $\delta_1 = \delta_1(\eps,1)$ so that
\begin{eqnarray}
\|x\|_X \leq \delta_1 \ \wedge \ \|u\|_\Uc \leq \delta_1\ \wedge \ t\in [0,1] \srs \|\phi(t,x,u)\|_X < \eps.
\label{eq:ncLF_CEP_def_tilde_delta}
\end{eqnarray}
Define for $k\in\N$ the following time sequence:
\[
t^{1}_k := \sup\{t\in[0,t_k]: \|\phi(t,x_k,u_k)\|_X \leq \delta_1\},
\]
if the supremum is taken over a nonempty set, and define 
$t^{1}_k:= 0$ otherwise.

Again as $\Sigma$ is CEP, for the above $\delta_1$ there is a $\delta_2>0$ so that
\begin{eqnarray}
\|x\|_X\leq \delta_2\ \wedge \ \|u\|_\Uc \leq \delta_2\ \wedge \ t\in [0,1] \srs \|\phi(t,x,u)\|_X < \delta_1.
\label{eq:ncLF_CEP_def_bar_delta}
\end{eqnarray}
Without loss of generality, we assume that $\delta_2$ is chosen small enough so that
\begin{eqnarray}
\alpha(\delta_1) > \psi(\delta_2).
\label{eq:Relation_tilde_delta_and_bar_delta}
\end{eqnarray}
We now define
\[
t^{2}_k := \sup\{t\in[0,t_k]: \|\phi(t,x_k,u_k)\|_X \leq \delta_2\},
\]
if the supremum is taken over a nonempty set, and set $t^{2}_k:= 0$ otherwise.

Since $u_k\to 0$ and $x_k\to 0$ as $k\to\infty$, there is $K>0$ so that $\|u_k\|_\Uc \leq \delta_2$ and 
$\|x_k\|_X \leq \delta_2$ for $k\geq K$.

From now on, we always assume that $k\geq K$.

Using \eqref{eq:ncLF_CEP_def_tilde_delta}, \eqref{eq:ncLF_CEP_def_bar_delta} and the cocycle property,
it is not hard to show that for $k\geq K$ it must hold that $t_k\geq 2$,
as otherwise, we arrive at a contradiction to 
$\|\phi(t_k,x_k,u_k)\|_X=\eps$.

Assume that $t_k - t^{1}_k <1$. This implies (since $t_k\geq 2$), that
$t^{1}_k>0$. By the cocycle property, we have  
\begin{eqnarray*}
\|\phi(t_k,x_k,u_k)\|_X= \|\phi(t_k - t^{1}_k,\phi(t^{1}_k,x_k,u_k),u_k(\cdot + t^{1}_k)\|_X.
\end{eqnarray*}
The axiom of shift invariance justifies the inequalities
\[
\|u_k(\cdot + t^{1}_k)\|_\Uc \leq \|u_k\|_\Uc \leq \delta_2\leq \delta_1.
\]
Since $\|\phi(t^{1}_k,x_k,u_k)\|_X = \delta_1$, and $t_k - t^{1}_k <1$, we have by \eqref{eq:ncLF_CEP_def_tilde_delta} that 
$\|\phi(t_k,x_k,u_k)\|_X < \eps$, a contradiction.
Hence $t_k - t^{1}_{k}\geq 1$ for all $k\geq K$. 

Analogously, we obtain that  $t^{1}_k - t^{2}_k\geq 1$ and $t_k - t^{2}_k\geq 2$.

Define 
\[
x^{2}_k:=\phi(t^{2}_k,x_k,u_k),\quad u^{2}_k := u_k(\cdot + t^{2}_k)
\]
and
\[
x^{1}_k:=\phi(t^{1}_k,x_k,u_k),\quad u^{1}_k := u_k(\cdot + t^{1}_k).
\]
Due to the axiom of shift invariance $u^{1}_k, u^{2}_k \in\Uc$ and 
\[
\|u^{1}_k\|_\Uc \leq \|u^{2}_k\|_\Uc \leq \|u_k\|_\Uc \leq \delta_2.
\]
Also by the definition of $t^{2}_k$ we have $\|x^{2}_k\|_X =\delta_2$.

Applying \eqref{eq:ni-ISS}, we obtain for $t:=t_k - t^{2}_k$ that 
\begin{align}
\label{nc_LS_proof_Eq1}
\int_0^{t_k - t^{2}_k} \alpha(\|\phi(s,x^{2}_k,u^{2}_k)\|_X) ds &\leq \psi(\|x^{2}_k\|_X) + (t_k - t^{2}_k)\sigma(\|u^{2}_k\|_{\Uc}) \nonumber\\
&\leq \psi(\delta_2) + (t_k - t^{2}_k)\sigma(\|u_k\|_{\Uc}).
\end{align}
On the other hand, by changing the integration variable and using the cocycle property, we obtain that
\begin{align*}
\int_0^{t_k - t^{2}_k}\hspace{-4mm} \alpha(\|\phi(s,x^{2}_k,u^{2}_k)&\|_X) ds
= \int_0^{t^{1}_k - t^{2}_k} \hspace{-4mm}\alpha(\|\phi(s,x^{2}_k,u^{2}_k)\|_X) ds +
\int_{t^{1}_k - t^{2}_k}^{t_k - t^{2}_k} \hspace{-4mm}\alpha(\|\phi(s,x^{2}_k,u^{2}_k)\|_X) ds\\
&= \int_0^{t^{1}_k - t^{2}_k} \hspace{-4mm}\alpha(\|\phi(s,x^{2}_k,u^{2}_k)\|_X) ds +
\int_{0}^{t_k - t^{1}_k} \hspace{-4mm}\alpha(\|\phi(s+ t^{1}_k - t^{2}_k,x^{2}_k,u^{2}_k)\|_X) ds\\
&= \int_0^{t^{1}_k - t^{2}_k} \hspace{-4mm}\alpha(\|\phi(s,x^{2}_k,u^{2}_k)\|_X) ds +
\int_{0}^{t_k - t^{1}_k} \hspace{-4mm}\alpha(\|\phi(s,x^{1}_k,u^{1}_k)\|_X) ds.
\end{align*}
For the last transition we have used that for every $s\in[0,t_k - t^{1}_k]$
\begin{align*}
\phi(s+ t^{1}_k - t^{2}_k,x^{2}_k,u^{2}_k)
&= \phi\big(s+ t^{1}_k - t^{2}_k,\phi(t^{2}_k,x_k,u_k),u_k(\cdot + t^{2}_k)\big)\\
&= \phi\big(s,\phi(t^{1}_k - t^{2}_k,\phi(t^{2}_k,x_k,u_k),u_k(\cdot + t^{2}_k)),u_k(\cdot + t^{1}_k)\big)\\
&= \phi\big(s,\phi(t^{1}_k,x_k,u_k),u_k(\cdot + t^{1}_k)\big)\\
&= \phi(s,x^{1}_k,u^{1}_k).
\end{align*}
By definition of $t^{2}_k$ and $t^{1}_k$, we have that 
\begin{align*}
\|\phi(s,x^{2}_k,u^{2}_k)\|_X  
&= \big\|\phi\big(s,\phi(t^{2}_k,x_k,u_k),u_k(\cdot + t^{2}_k)\big)\big\|_X \\
&= \|\phi(s+t^{2}_k,x_k,u_k)\|_X 
\geq \delta_2,\qquad s\in[0,t_k-t^{2}_k],
\end{align*}
and, analogously,
\[
\|\phi(s,x^{1}_k,u^{1}_k)\|_X \geq \delta_1,\quad s\in[0,t_k-t^{1}_k].
\]
Continuing the above estimates and using that $t_k - t^{1}_k\geq 1$ and $\alpha(\delta_1)>\alpha(\delta_2)$, we arrive at
\begin{align*}
\int_0^{t_k - t^{2}_k} \alpha(\|\phi(s,x^{2}_k,u^{2}_k)\|_X) ds 
&\geq (t^{1}_k - t^{2}_k)\alpha(\delta_2) + (t_k - t^{1}_k)\alpha(\delta_1)\\
&\geq (t^{1}_k - t^{2}_k)\alpha(\delta_2) + (t_k - t^{1}_k-1)\alpha(\delta_1) + \alpha(\delta_1)\\
&\geq (t_k - t^{2}_k - 1)\alpha(\delta_2) + \alpha(\delta_1).
\end{align*}
Since $t_k - t^{2}_k\geq 2$ and in view of \eqref{eq:Relation_tilde_delta_and_bar_delta}, we derive
\begin{eqnarray}
\label{nc_LS_proof_Eq2}
\int_0^{t_k - t^{2}_k} \alpha(\|\phi(s,x^{2}_k,u^{2}_k)\|_X) ds >
\frac{t_k - t^{2}_k}{2}\alpha(\delta_2) + \psi(\delta_2).
\end{eqnarray}
Combining the inequalities \eqref{nc_LS_proof_Eq1} and \eqref{nc_LS_proof_Eq2}, we obtain
\begin{eqnarray*}
\frac{t_k - t^{2}_k}{2}\alpha(\delta_2)  < (t_k - t^{2}_k)\sigma(\|u_k\|_{\Uc}).
\end{eqnarray*}
This leads to
\begin{eqnarray*}
\frac{1}{2}\alpha(\delta_2)  < \sigma(\|u_k\|_{\Uc}), \quad k\geq K.
\end{eqnarray*}
Finally, since $\lim_{k\to\infty}\|u_k\|_\Uc= 0$, letting $k\to\infty$ we come to a contradiction.
This shows that $\Sigma$ is ULS.
\end{proof}

Finally, we characterize ISS in terms of norm-to-integral ISS.
\begin{theorem}
\label{thm:ncISS_LF_sufficient_condition_NEW}
Let $\Sigma=(X,\Uc,\phi)$ be a forward complete control system. Then $\Sigma$ is ISS if and only if $\Sigma$ is norm-to-integral ISS and has CEP and BRS properties. 
\end{theorem}

\begin{proof}
\q{$\Rightarrow$}. Clearly, ISS implies CEP and BRS properties. By Proposition~\ref{prop:ISS-implies-norm-i-ISS}, ISS implies  norm-to-integral ISS.

\q{$\Leftarrow$}.
Propositions~\ref{prop:ncISS_implies_ULIM} and \ref{prop:ncISS_plus_CEP_implies_ULS} imply that $\Sigma$ is ULIM and ULS.
Since $\Sigma$ is BRS, Theorem~\ref{thm:UAG_equals_ULIM_plus_LS} shows that $\Sigma$ is ISS.
\end{proof}

\subsection{Remark on input-to-state practical stability}

In some cases, it is impossible (as in quantized control)
or too costly to construct feedback that results in an ISS closed-loop system.
For these applications, one defines the following relaxation of the ISS property:
\begin{definition}
\label{Def:ISpS_wrt_set}
\index{ISpS}
\index{input-to-state stability!practical}
A forward complete control system $\Sigma=(X,\Uc,\phi)$ is called \emph{(uniformly) input-to-state practically stable (ISpS)}, if there exist $\beta \in \KL$, $\gamma \in \Kinf$ and $c>0$
such that for all $ x \in X$, $ u\in \Uc$ and $ t\geq 0$ the following holds:
\begin {equation}
\label{isps_sum}
\| \phi(t,x,u) \|_X \leq \beta(\| x \|_X,t) + \gamma( \|u\|_{\Uc}) + c.
\end{equation}
\end{definition}

In Theorem~\ref{thm:ncISS_LF_sufficient_condition_NEW}, we require that the CEP property holds. If the CEP property is not available, we can still infer input-to-state practical stability of $\Sigma$, using the main result in \cite{Mir19a}.
\begin{theorem}
\label{thm:ncISpS_LF_sufficient_condition_NEW}
Let $\Sigma=(X,\Uc,\phi)$ be a forward complete control system with BRS. If $\Sigma$ is norm-to-integral ISS, then $\Sigma$ is ISpS.
\end{theorem}

\begin{proof}
Proposition~\ref{prop:ncISS_implies_ULIM} implies that $\Sigma$ is ULIM.
Since $\Sigma$ is also BRS, \cite[Theorem III.1]{Mir19a} shows that $\Sigma$ is ISpS.
\end{proof}

\section{Strong input-to-state stability}
\label{sec:WeakISS}

In this section, we introduce and study a notion of strong input-to-state stability, which is
weaker than ISS and is equivalent to the strong stability of linear undisturbed systems (which explains the terminology). 

\begin{definition}
\label{Def:sISS}
\index{sISS}
\index{input-to-state stability!strong}
A forward complete system $\Sigma=(X,\Uc,\phi)$ is called \emph{strongly input-to-state stable
(sISS)}, if there exist $\gamma \in \K$, $\sigma \in \Kinf$ and $\beta: X \times \R_+ \to \R_+$, so that

\begin{enumerate}[label=(\roman*)]
    \item  $\beta(x,\cdot) \in \LL$ for all $x \in X$, $x \neq 0$.
  \item $\beta(x,t) \leq \sigma(\|x\|_X)$ for all $x \in X$ and all $t \geq 0$.
    \item For all $ x \in X$, all $u\in \Uc$, and all $t\geq 0$ it holds that
\begin{equation}
\label{wISS_sum}
\| \phi(t,x,u) \|_{X} \leq \beta(x,t) + \gamma( \|u\|_{\Uc}).
\end{equation}
\end{enumerate}
\end{definition}

Clearly, ISS implies sISS, but the converse implication doesn't hold for infinite-dimensional systems in general, as shown in the following example. 
For the systems without inputs, there are numerous examples of linear PDE systems which are strongly stable but not exponentially stable, see \cite{Oos00, AvL98}, etc.

%
%

\begin{example}[sISS does not imply LISS, even for linear systems]
\label{examp:sISS_not_ISS}
Consider the following system
	\begin{eqnarray}
	\dot{x}_k = -\frac{1}{k}x_k + \frac{1}{k^3}u,\quad k\in\N.
	\label{eq:simple_sISS_linear_system}
	\end{eqnarray}
	Here $u \in L^{\infty}(\R_+,\R)$, $x_k(r)\in\R$ and $X:=\ell_1$.

This system is not LISS, as it is not exponentially stable for $u\equiv 0$.

The solution of the $k$-th subsystem of \eqref{eq:simple_sISS_linear_system} is given by
\begin{eqnarray*}
\phi_k(t,x_k,u) &=& e^{-\frac{1}{k}t}x_k + \frac{1}{k^3} \int_0^t e^{-\frac{1}{k}(t-s)}u(s)ds.
\end{eqnarray*}
Hence
\begin{eqnarray*}
|\phi_k(t,x_k,u)| &\leq& e^{-\frac{1}{k}t}|x_k| + \frac{1}{k^3} \|u\|_{\infty} \int_0^t e^{-\frac{1}{k}(t-s)}ds\\
									&=& e^{-\frac{1}{k}t}|x_k| + \frac{1}{k^3} k (1-e^{-\frac{1}{k}t}) \|u\|_{\infty}\\
									&=& e^{-\frac{1}{k}t}|x_k| + \frac{1}{k^2} \|u\|_{\infty}.
\end{eqnarray*}
Let $x=(x_k) \in \ell_1$. Then
\begin{eqnarray*}
\|\phi(t,x,u)\|_X &\leq& \sum_{k=1}^\infty e^{-\frac{1}{k}t}|x_k| + \sum_{k=1}^\infty \frac{1}{k^2} \|u\|_{\infty}\\
									&=& \sum_{k=1}^\infty e^{-\frac{1}{k}t}|x_k| + \frac{\pi^2}{6} \|u\|_{\infty}.
\end{eqnarray*}
Define $\beta(x,t):= \sum_{k=1}^\infty e^{-\frac{1}{k}t}|x_k|$. 

Then $\beta(x,\cdot)\in\LL$ for any $x\in X$ and $\beta(x,t)\leq \sum_{k=1}^\infty |x_k| = \|x\|_X$ for all $t,x$. 
Hence, \eqref{eq:simple_sISS_linear_system} is sISS.
\end{example}

Strong ISS can be characterized as follows.
\begin{theorem}[sISS superposition theorem]
\label{wISS_equals_sAG_GS}
Let $\Sigma=(X,\Uc,\phi)$ be a forward complete control system.
The following statements are equivalent.
\begin{enumerate}[label=(\roman*)]
  \item $\Sigma$ is sISS.
  \item $\Sigma$ is sAG and UGS.
  \item $\Sigma$ is sLIM and UGS.
\end{enumerate}
\end{theorem}

\begin{proof}
(i) $\Rightarrow$ (ii).
Let $\Sigma$ be sISS with the corresponding
$\beta:X \times \R_+ \to \R_+$ and $\sigma \in \Kinf$ and $\gamma \in
\mathcal{K}$. By definition, $\Sigma$ is UGS characterized by $\sigma$ and
$\gamma$.

Fix any $x \in X$ and any $\eps >0$. Define $\tau=\tau(\eps,x)$ as the
solution of the equation $\beta(x,\tau)=\eps$ (if this solution exists,
then it is unique because of the monotonicity of $\beta$ in the second argument; if it does not exist, we set $\tau(\eps,x)=0$). Then for all $t \geq \tau$ and all $u \in \Uc$
\begin{eqnarray*}
\|\phi(t,x,u)\|_X &\leq& \beta(x,t) + \gamma(\|u\|_{\Uc})\\
  &\leq& \beta(x,\tau) + \gamma(\|u\|_{\Uc}) \\
    &\leq& \eps + \gamma(\|u\|_{\Uc}),
\end{eqnarray*}
and the estimate \eqref{sAG_Absch} holds.
Thus, sISS implies sAG.

(ii) $\Rightarrow$ (iii). This is clear.

(iii) $\Rightarrow$ (ii). Can be proved along the lines of Lemma~\ref{lem:ULIM_plus_GS_implies_UAG}.

(ii) $\Rightarrow$ (i). Assume that $\Sigma$ is UGS and sAG. We are going to construct $\beta:X \times \R_+ \to \R_+$ with the properties as in Definition~\ref{Def:sISS}, so that \eqref{wISS_sum} holds.

Uniform global stability of $\Sigma$ implies that there exist $\sigma,\gamma \in\Kinf$ so that for all $t \geq 0$, all $u \in \Uc$, and all $x \in X$ it holds that
\[
\|\phi(t,x,u)\|_X \leq \sigma(\|x\|_X) + \gamma(\|u\|_{\Uc}).
\]
Define $\eps_n(x):= 2^{-n} \sigma(\|x\|_X)$, for each $x \in X$, $n \in \N$. Due to
 sAG, there exists a sequence of times $\tau_n = \tau_n(x):=\tau(\eps_n(x),x)$ that we
 assume without loss of generality to be strictly increasing in $n$, such that
 \[
 \|\phi(t,x,u)\|_X \leq \eps_n + \gamma(\|u\|_{\Uc}) \quad \forall u \in \Uc,\ \forall t \geq \tau_n.
 \]
 Define $\beta(x,\tau_n):=\eps_{n-1}$, for $n \in \N$.

 Now extend the function $\beta(x,\cdot)$ for $t \in \R_+ \backslash
 \{\tau_n:\ n \in \N\}$ so that $\beta(x,\cdot) \in \LL$ and $\beta(x,t) \leq
 2 \sigma(\|x\|_X)$ for all $t \geq 0$ (this can be done by choosing the values of $\beta(x,t)$ sufficiently small for $t \in [0,\tau_1)$).

 The function $\beta$ satisfies the estimate \eqref{wISS_sum}, because for
 all $t \in (\tau_n,\tau_{n+1})$ it holds that $ \|\phi(t,x,u)\|_X \leq \eps_n + \gamma(\|u\|_{\Uc})< \beta(x,t) + \gamma(\|u\|_{\Uc})$.
 Performing this procedure for all $x \in X$, we obtain the definition of the function $\beta$.
 This shows $\Sigma$ is sISS.
\end{proof}

For sISS systems, we have the following counterpart of Proposition~\ref{prop:Converging_input_uniformly_converging_state}:
\begin{proposition}
\label{prop:Converging_input_converging_state}
Let $\Sigma=(X,\Uc,\phi)$ be an sISS control system. If $u\in\Uc$ and $\lim_{t \to \infty}\|u(\cdot + t)\|_{\Uc} = 0$, then $\lim_{t \to \infty}\|\phi(t,x,u)\|_X = 0$ for any $x \in X$.
\end{proposition}

\begin{proof}
Let $\Sigma=(X,\Uc,\phi)$ be an sISS control system and assume without loss of generality that the corresponding gain $\gamma$ is a $\Kinf$-function.
Pick any $x \in X$ and any $u\in \Uc$ so that $\lim_{t \to \infty}\|u(\cdot + t)\|_{\Uc} = 0$.

Pick any $\eps>0$. We are going to show that there is a time $t_\eps>0$ so that $\|\phi(t,x,u)\|_X\leq\eps$ for all
$t\geq t_\eps$.

Choose $t_1$ so that $\|u(\cdot + t_1)\|_{\Uc}\leq \gamma^{-1}(\eps)$.
Due to the cocycle property and sISS of $\Sigma$, we have
\begin{eqnarray*}
\|\phi(t+t_1,x,u)\|_X  &=& \big\|\phi\big(t,\phi(t_1,x,u),u(\cdot + t_1)\big)\big\|_X \\
                                             &\leq& \beta\big(\phi(t_1,x,u),t\big) + \gamma\big(\|u(\cdot+t_1)\|_\Uc\big) \\
                                             &\leq& \beta\big(\phi(t_1,x,u),t\big) + \eps.
\end{eqnarray*}
Now pick any $t_2$ in a way that $\beta\big(\phi(t_1,x,u),t_2\big) \leq \eps$.
This ensures that $\|\phi(t_2+t_1,x,u)\|_X \leq 2\eps.$
Using the cocycle property once again, we obtain for all $t\geq0$:
\begin{align*}
\|\phi(t+t_2&+t_1,x,u)\|_X \\
&\leq \beta\big(\phi(t_2+t_1,x,u),t\big) + \gamma\big(\|u(\cdot+t_1+t_2)\|_\Uc\big).
\end{align*}
Due to the axiom of shift invariance, we have $\|u(\cdot+t_2+t_1)\|_\Uc \leq \|u(\cdot+t_1)\|_\Uc\leq \gamma^{-1}(\eps)$.
Strong ISS of $\Sigma$ ensures for all $t>0$ that
\begin{eqnarray*}
\|\phi(t+t_2+t_1,x,u)\|_X &\leq& \sigma\big(\|\phi(t_2+t_1,x,u)\|_X\big) + \eps\\
&\leq&\sigma(2\eps) + \eps.
\end{eqnarray*}
Since $\eps>0$ can be chosen arbitrarily small, and since $\sigma(2\eps)+\eps\to 0$ as $\eps\to 0$, the claim of the proposition follows.
\end{proof}

\section{Counterexamples}
\label{sec:Counterexamples}

In this section, we construct:
\begin{itemize}
  \item two nonlinear systems $\Sc^1$, $\Sc^3$ without inputs,
    \item two nonlinear systems $\Sc^2$, $\Sc^4$ with inputs,
\end{itemize}
providing counterexamples that show that the following implications are false (note that the axioms ($\Sigma$\ref{axiom:Identity})--($\Sigma$\ref{axiom:Cocycle}) are fulfilled for all $\Sc^i$):
\begin{itemize}
    \item[$\Sc^1$:] FC\,$\wedge$\,0-GAS\,$\wedge$\,0-UAS\ \ $\not\Rightarrow$\ \ BRS.
    \item[$\Sc^2$:] FC\,$\wedge$\,0-UGAS\,$\wedge$\,AG\,$\wedge$\,LISS\ \ $\not\Rightarrow$\ \ BRS.
    \item[$\Sc^3$:] FC\,$\wedge$\,BRS\,$\wedge$\,0-GAS\,$\wedge$\,0-UAS\ \ $\not\Rightarrow$\ \ 0-UGS.
    \item[$\Sc^4$:] FC\,$\wedge$\,BRS\,$\wedge$\,0-UGAS\,$\wedge$\,AG\,$\wedge$\,LISS\ \ $\not\Rightarrow$\ \ UGS.
\end{itemize}

System $\Sc^1$ shows that already for undisturbed systems, nonuniform global attractivity does not ensure that the solution map $\phi(t,\cdot)$ maps bounded balls into bounded balls. And even if it does, then global stability still cannot be guaranteed, as clarified by system $\Sc^3$. This shows that the difference between nonuniform attractivity and stability is much more significant for nonlinear infinite-dimensional systems than for ODEs.

\begin{remark}
\label{rem:Nonimplications}
Before we proceed to the constructions of the systems $\Sc^i$, let us show how they justify the negated implications depicted in Figure~\ref{fig:ISS_Equiv}.
\begin{itemize}
    \item[(i)] Follows from Example~\ref{examp:sISS_not_ISS}.
    \item[(ii)] Follows by construction of $\Sc^3$.
    \item[(iii)] Follows by construction of $\Sc^4$.
    \item[(iv)] Follows by construction of $\Sc^4$.
    \item[(v)] Follows from Example~\ref{examp:sISS_not_ISS}.
    \item[(vi)] Follows by construction of $\Sc^3$ and/or $\Sc^4$.
    \item[(vii)] Follows from Example~\ref{examp:sISS_not_ISS}.
    \item[(viii)] Follows by construction of $\Sc^4$.
    \item[(ix)] Follows from Example~\ref{examp:sISS_not_ISS}.
\end{itemize}
In addition, in Example~\ref{ex:0-UGAS_sAG_AG_zero_gain_LISS_not_ISS}, we show by example that 0-UGAS $\wedge$ sAG $\wedge$ AG with zero gain $\wedge$ UGS with zero gain $\wedge$ LISS with zero gain do not imply ISS (and even do not imply ULIM). Hence, the properties of the \q{second} level (sISS and AG $\wedge$ 0-UGAS) are not only different from each other (in the sense that they do not imply each other), but also even taken together, they do not imply ISS.

Finally, systems $\Sc^1$ and $\Sc^2$ show that the systems with
global nonuniform attractivity properties together with very strong
properties near the equilibrium may not even be BRS.
\end{remark}

\begin{example}[FC\,$\wedge$\,0-GAS\,$\wedge$\,0-UAS\ \ $\not\Rightarrow$\ \ BRS]
\label{0-GAS_but_not_GS}
According to Remark~\ref{0-GAS_strong_stability}, 0-GAS implies 0-UGS for \emph{linear}
infinite-dimensional systems. Now, we show that for nonlinear systems, 0-GAS does not even imply the BRS of the undisturbed system.
Consider the infinite-dimensional system $\Sc^1$ defined by
\begin{eqnarray}
\Sc^1: \left
\{\begin{array}{l}
\Sc^1_k:
\left
\{\begin{array}{l}
\dot{x}_k = -x_k + x_k^2 y_k - \frac{1}{k^2} x_k^3,\\
\dot{y}_k = - y_k.
\end{array}
\right.
\\
k \in\N,
\end{array}
\right.
\label{eq:GAS_not_GS}
\end{eqnarray}
with the state space
\begin{eqnarray}
X:=\ell_2= \left\{ (z_k)_{k=1}^{\infty}: \sum_{k=1}^{\infty} |z_k|^2 <\infty , \quad z_k = (x_k,y_k) \in \R^2 \right\}.
\label{eq:l2_space}
\end{eqnarray}
 We show that $\Sc^1$ is forward complete, 0-GAS, and 0-UAS but does not have bounded reachability sets.
Furthermore, the nonlinearity that defines $\Sc^1$ is Lipschitz continuous on bounded balls w.r.t.\  the state $(x,y)$.

Before we give detailed proof of these facts, let us give an informal explanation of this phenomenon.
If we formally place $0$ into the definition of $\Sc^1_k$ instead of the term $- \frac{1}{k^2} x_k^3$, then the state of $\Sc^1_k$ (for each $k$) will exhibit a finite
escape time, provided $y_k(0)$ is chosen large enough.
The term  $- \frac{1}{k^2} x_k^3$ prevents the solutions of $\Sc^1_k$ from growing to infinity: the solution then looks like a pike, that is stopped by the damping $- \frac{1}{k^2} x_k^3$, and converges to 0 since $y_k(t) \to 0$ as $t\to\infty$.
However, the larger is $k$, the higher will be the corresponding pikes, and hence there is no uniform bound for the solution of $\Sc^1$ starting from a bounded ball.
Now, we proceed to the rigorous proof.

\textbf{Step 1: Right-hand side of \eqref{eq:GAS_not_GS} is Lipschitz continuous on bounded balls.} 
Denote the right-hand side of \eqref{eq:GAS_not_GS} by $f$, pick any $R>0$, any $(x,y) \in B_R$,	 any $(\tilde{x}, \tilde{y}) \in B_R$, and consider 
\begin{align*}
\|f(x,y)-f(\tilde{x},&\tilde{y})\|_X^2 
=
\sum_{k=1}^\infty \Big|-x_k + x_k^2 y_k - \frac{1}{k^2} x_k^3 -\big(-\tilde{x}_k + \tilde{x}_k^2 \tilde{y}_k - \frac{1}{k^2} \tilde{x}_k^3\big)\Big|^2 + |y_k - \tilde{y}_k|^2\\
&\le
\sum_{k=1}^\infty\Big| |x_k-\tilde{x}_k| + x_k^2 |y_k - \tilde{y}_k| + |x_k^2 - \tilde{x}_k^2| |\tilde{y}_k|  + \frac{1}{k^2} |x_k^3 - \tilde{x}_k^3|\Big|^2 + |y_k - \tilde{y}_k|^2.
\end{align*}
As $|x_k|\leq R$, $|\tilde{x}_k|\leq R$, $|y_k|\leq R$, we see that for all $k\in\N$ we have that
\[
|x_k^3 - \tilde{x}_k^3| = |x_k - \tilde{x}_k| |x_k^2 + x_k\tilde{x}_k + \tilde{x}_k^2 |\leq 3R^2 |x_k - \tilde{x}_k|
\]
and 
\[
|x_k^2 - \tilde{x}_k^2| |\tilde{y}_k| \leq  2 R^2 |x_k - \tilde{x}_k|.
\]
Thus, we continue the above computations as follows:
\begin{align*}
\|f(x,y)-&f(\tilde{x},\tilde{y})\|_X^2 \\
&\le
\sum_{k=1}^\infty\Big| |x_k-\tilde{x}_k| + R^2 |y_k - \tilde{y}_k| + 2R^2|x_k - \tilde{x}_k|  + 3R^2 |x_k - \tilde{x}_k|\Big|^2 + |y_k - \tilde{y}_k|^2\\
&\le
C(R) \sum_{k=1}^\infty |x_k-\tilde{x}_k|^2 + |y_k - \tilde{y}_k|^2\\
&\le
C(R) \|(x,y) - (\tilde{x},\tilde{y})\|^2_X,
\end{align*}
for a certain $C=C(R)$. This shows Lipschitz continuity of $f$ on bounded balls.

\textbf{Step 2: $\Sc^1$ is 0-UAS.} Indeed, for $r<1$ the Lyapunov function $V(z) = \|z\|^2_{\ell_2}=\sum_{k=1}^{\infty} (x_k^2+y_k^2 )$ satisfies
for all $z_k$ with $|x_k| \leq r$ and $|y_k| \leq r$, $k \in \N$, the estimate
\begin{equation}
\label{eq:Sigma1-Vdot}
\begin{aligned}
\dot{V}(z) &= 2 \sum_{k=1}^{\infty} \Big(-x^2_k + x_k^3 y_k - \frac{1}{k^2} x_k^4- y_k^2\Big) \\
&\leq 2 \sum_{k=1}^{\infty} \Big(-x^2_k + |x_k|^3 |y_k| - \frac{1}{k^2} x_k^4- y_k^2\Big)\\
&\leq 2 \sum_{k=1}^{\infty}  \Big( (r^2 - 1)x^2_k - y_k^2\Big) \\
&\leq 2(r^2 - 1)V(z).
\end{aligned}
\end{equation}
We see that $V$ is an exponential local Lyapunov function for the system \eqref{eq:GAS_not_GS}
and thus \eqref{eq:GAS_not_GS} is locally uniformly asymptotically
stable. Indeed, it is not hard to show that the domain of attraction
contains $\{ z \in \ell_2 : |x_k| < r, |y_k| < r, \, \forall k \}$.

\textbf{Step 3: $\Sc^1$ is forward complete and globally attractive.}
First, we point out that every $\Sc^1_k$ is 0-GAS (and hence 0-UGAS, since $\Sc^1_k$ is finite-dimensional). This follows from the fact that any $\Sc^1_k$ is a cascade interconnection of an ISS $x_k$-system (with $y_k$ as an input) and a globally asymptotically stable $y_k$-system, see \cite{Son89}.

Furthermore, for any $z(0) \in \ell_2$ there exists a finite $N >0$ such that
$|z_k(0)| \leq \tfrac{1}{2}$ for all $k \geq N$. Decompose the norm of $z(t)$ as follows:
\[
\|z(t)\|_{\ell_2} = \sum_{k=1}^{N-1} |z_k(t)|^2 + \sum_{k=N}^{\infty} |z_k(t)|^2.
\]
According to the previous arguments, $\sum_{k=1}^{N-1} |z_k(t)|^2 \to 0$ as $t \to 0$ since all $\Sc^1_k$ are 0-UGAS for $k=1,\ldots,N-1$.

Since $|z_k(0)| \leq \tfrac{1}{2}$ for all $k \geq N$, we can apply the computations as in \eqref{eq:Sigma1-Vdot}
in order to obtain (for $r:=\frac{1}{2}$) that
\[
\frac{d}{dt}\Big(\sum_{k=N}^{\infty} |z_k(t)|^2\Big) \leq -\frac{3}{2} \sum_{k=N}^{\infty} |z_k(t)|^2.
\]
Hence
$\sum_{k=N}^{\infty} |z_k(t)|^2$ decays monotonically and exponentially to $0$ as $t \to \infty$.
Overall, $\|z(t)\|_{\ell_2} \to 0$ as $t \to \infty$ which shows that $\Sc^1$ is forward complete, 0-GAS and 0-UAS.

\textbf{Step 3: $\Sc^1$ is not BRS.}
To prove this, it is enough to show that there exists an $r>0$ and $\tau>0$ so that for any $M>0$, there exist $z \in \ell_2$ and $t \in [0,\tau]$ so that $\|z\|_{\ell_2} =r$ and
$\|\phi(t,z,0)\|_{\ell_2} > M$.

Let us consider $\Sc^1_k$.
For $y_k \geq 1$ and for $x_k \in [0,k]$ it holds that
\begin{equation}
\dot{x}_k \geq -2x_k + x_k^2.
\label{eq:Counterex_tmp}
\end{equation}
Pick an initial state $x_k(0)=c>0$ (which is independent of $k$) so that the solution of $\dot{x}_k = -2x_k + x_k^2$ blows up to infinity in time $t^*=1$. Now pick $y_k(0)=e$ (Euler's constant) for all $k=1,2,\ldots$. For this initial condition, we obtain $y_k(t) = e^{1-t} \geq 1$ for $t \in [0,1]$.
And consequently for $z_k(0)=(c,e)^T$ there exists a time $\tau_k\in
(0,1)$ such that $x_k(\tau_k) = k$ for the solution of  $\Sc^1_k$.

Now consider an initial state $z(0)$ for $\Sc^1$, where
$z_k(0)=(c,e)^T$ and $z_j(0)= (0,0)^T$ for $j\neq k$.
For this initial state, we have that $\|z(t)\|_{\ell_2} = |z_k(t)|$ and
\[
\sup_{t \geq 0}\|z(t)\|_{\ell_2}=\sup_{t \geq 0}|z_k(t)| \geq |x_k(\tau_k)| \geq  k.
\]
As $k\in\N$ was arbitrary, the system $\Sc^1$ is not BRS.
\hfill \qedsymbol
\end{example}


\begin{example}[FC\,$\wedge$\,0-UGAS\,$\wedge$\,AG\,$\wedge$\,LISS\ \ $\not\Rightarrow$\ \ BRS]
\label{0-UGAS_AG_LISS_but_not_GS}
In this modification of Example~\ref{0-GAS_but_not_GS},	 it is demonstrated that 0-UGAS $\wedge$ AG $\wedge$ LISS does not imply BRS.
Let $\Sc^2$ be defined by
\begin{eqnarray*}
\Sc^2: \left
\{\begin{array}{l}
\Sc^2_k:
\left
\{\begin{array}{l}
\dot{x}_k = -x_k + x_k^2 y_k |u|- \frac{1}{k^2} x_k^3,\\
\dot{y}_k = - y_k.
\end{array}
\right.
\\
k=1,2,\ldots,
\end{array}
\right.
\end{eqnarray*}
And let the state space of $\Sc^2$ be $\ell_2$ (see \eqref{eq:l2_space}) and its input space be $\Uc:=PC_b(\R_+,\R)$.

Evidently, this system is 0-UGAS. Also, it is clear that $\Sc^2$ is not BRS, since for $u \equiv 1$, we obtain exactly the system from Example~\ref{0-GAS_but_not_GS}, which is not BRS. The proof that this system is forward complete, LISS and AG with zero gain mimics the argument we exploited to show 0-GATT of Example~\ref{0-GAS_but_not_GS}. Thus, we omit it.
\hfill \qedsymbol
\end{example}

\begin{example}[FC\,$\wedge$\,BRS\,$\wedge$\,0-GAS\,$\wedge$\,0-UAS\ \ $\not\Rightarrow$\ \ 0-UGS]
\label{ex:FC_0GAS_BRS_not_GS}
We construct a counterexample in 3 steps.

\textbf{Step 1.} Let us revisit Example~\ref{0-GAS_but_not_GS} and find
useful estimates from above for the dynamics of the subsystems
$\Sc^1_k$.

\textbf{Step 1.1.}
We first note that for initial conditions $z^0_k = (x^0_k,y^0_k)$ with $x^0_ky^0_k
\leq 0$ we have for the corresponding solution $z_k(\cdot)$ of $\Sc^1_k$ that $x_k(t)y_k(t) \leq 0$. Indeed, if $y_k(0)<0$, then $y_k(t) <0$ for all $t>0$ for which the trajectory is defined. As $x_k(0)\geq0$, we see from the equation for $\Sc^1_k$ that $\dot{x}_k(t)<0$ as long as $x_k(t)>0$. However, if $x_k(\tau)=0$ for a certain $\tau>0$, then $x_k(t)=0$ also for $t\geq \tau$ in view of uniqueness of solutions. Hence, $x_k$ cannot cross zero, and thus $x_k(t)y_k(t)\leq 0$ for all time of existence of solutions. The case when $y_k(0)\geq 0$ can be treated similarly.

Now define 
\begin{align}
\label{eq:Lyapunov-function-V-k-aux}
V_k(z_k):=x_k^2+y_k^2,\quad z_k = (x_k,y_k) \in\R^2.
\end{align}
We have for any $z_k\in\R^2$ that 
\[
\dot{V}_k(z_k) = 2 \Big(-x^2_k + x_k^2 x_k y_k - \frac{1}{k^2} x_k^4- y_k^2\Big) \leq 0.
\]
This shows that $x^0_ky^0_k\leq 0$  implies $|z_k(t)| \leq |z^0_k|$ for all $t\geq0$.

\textbf{Step 1.2.} It is easy to check that for each $k\in\N$ and each
$z_k(0)=(x_k(0),y_k(0)) \in \R^2$ with $y_k(0)x_k(0) >0$
the solution of $\Sc^1_k$ for the
initial condition $z_k(0)$ (see \eqref{eq:GAS_not_GS}) can be estimated
in norm by
\begin{eqnarray}
\label{eq:x_k-aux-estimate-counterex}
|x_k(t)| \leq |\hat{x}_k(t)|
\end{eqnarray}
where $\hat{x}_k(t)$ is the first component of the solution of the system
\begin{eqnarray}
\hat{\Sc}^1_k:
\left
\{\begin{array}{l}
\dot{\hat{x}}_k(t) = \hat{x}_k^2(t) y_k(0), \\
\dot{y}_k = - y_k.
\end{array}
\right.
\label{eq:Examp3_HilfSys}
\end{eqnarray}
with initial condition $\hat{z}_k(0) = (\hat{x}_k(0),y_k(0)) =
(x_k(0),y_k(0))$.

Indeed, taking $V_k$ as defined by \eqref{eq:Lyapunov-function-V-k-aux}, and differentiating it along the solutions of $S_k^1$ and of $\hat{S}_k^1$, see see that $\dot{V}_k(z_k(t)) \leq \dot{V}_k(\hat{z}_k(t))$. 
This implies that $|z_k(t)| \leq |\hat{z}_k(t)|$, and since the $y_k$-dynamics in both systems are the same, we arrive at
\eqref{eq:x_k-aux-estimate-counterex}.

\textbf{Step 1.3.} 
This solution of the $\hat{x}_k$-subsystem of \eqref{eq:Examp3_HilfSys} reads as
\[
\hat{x}_k(t) = \frac{x_k(0)}{1-ty_k(0)x_k(0)},
\]
and this solution exists for $t\in[0,\frac{1}{y_k(0)x_k(0)})$.

Now pick any $R>0$ and assume that $z_k(0) = (x_k(0),y_k(0)) \in B_R$.
Since
\[
\frac{1}{y_k(0)x_k(0)}\geq \frac{2}{y^2_k(0) + x^2_k(0)} \geq \frac{2}{R^2},
\]
the solutions of \eqref{eq:Examp3_HilfSys} for any initial condition $z_k(0)\in B_R$ exist at least on the time interval $[0,2R^{-2})$.
Moreover, for every such solution for
$t\in [0,(2y_k(0)x_k(0))^{-1})$ (and in particular for $t\in [0,R^{-2}]$) it holds that
\[
|\hat{x}_k(t)| \leq 2|x_k(0)|.
\]
Overall, for each $R>0$, all $k\in\N$, all $z_k(0)=(x_k(0),y_k(0)) \in
B_R$ and all
$t\in [0,R^{-2}]$
 the solution of $\Sc^1_k$
corresponding to the initial condition $z_k(0)$ satisfies
\begin{eqnarray}
|z_k(t)| =\sqrt{x_k^2(t) + y_k^2(t)} \leq \sqrt{\hat{x}_k^2(t) + y_k^2(t)} \leq 2 |z_k(0)|.
\label{eq:Examp3_local_estimate_Sigma_Sys}
\end{eqnarray}

\textbf{Step 2.} Now we are going to modify the system $\Sc^1$
by using time transformations.
 Define
$\tilde{x}_k(t):=x_k(\frac{t}{k})$, $\tilde{y}_k(t):=y_k(\frac{t}{k})$,
for any $t\geq 0$ and any $k\geq 1$. In other words, we make the time of the $k$-th mode $k$
times slower than the time of $\Sc^1_k$. This new system we denote by
$\tilde\Sc^1$. The equations defining $\tilde\Sc^1$
are
\begin{eqnarray}
\tilde\Sc^1: \left
\{\begin{array}{l}
\tilde\Sc^1_k:
\left
\{\begin{array}{l}
\dot{\tilde{x}}_k = \frac{1}{k}\big(-\tilde{x}_k + \tilde{x}_k^2 \tilde{y}_k - \frac{1}{k^2} \tilde{x}_k^3\big),\\
\dot{\tilde{y}}_k = - \frac{1}{k} \tilde{y}_k.
\end{array}
\right.
\\
k=1,2,\ldots
\end{array}
\right.
\label{eq:GAS_BRS_not_GS}
\end{eqnarray}
Again the state space of $\tilde\Sc^1$ is $\ell_2$, see \eqref{eq:l2_space}.

We have seen that $\Sc^1$ fails to satisfy the BRS property since, with the increase of $k$, the solutions of subsystems $\Sc^1_k$ at a given time $t$ have larger spikes.
A non-uniform change of clocks in $\Sc^1$, performed above, makes such behavior impossible.
At the same time, $\tilde\Sc^1$ still is not 0-UGS. Next, we show detailed proofs of these facts.

From the computation in \eqref{eq:Sigma1-Vdot} it is easy to obtain that
for the dynamics of $\tilde\Sc^1$ we have
$\dot V(z) \leq 0$ if $\|z\|_{\ell_2}\leq 1$. It
follows that for all $z(0)\in \ell_2$ with $\|z(0)\|_{\ell_2}\leq 1$ we have
\[
\|z(t)\|_X\leq \|z(0)\|_X,
\]
and therefore $\tilde\Sc^1$ is 0-ULS.

Forward completeness and global attractivity of $\tilde\Sc^1$ can be shown along the lines of Example~\ref{0-GAS_but_not_GS}.
This shows that $\tilde\Sc^1$ is 0-GAS.

Let us prove that $\tilde\Sc^1$ is BRS.
Pick any $R>0$, any time $\tau>0$ and any $z\in \ell_2$: $\|z\|_{\ell_2} \leq R$. In view of \eqref{eq:Examp3_local_estimate_Sigma_Sys}, we have for any $k\in\N$:
\begin{equation}
|\tilde{z}_k(t)| \leq 2 |z_k(0)|\,\quad
\forall \ t\in [0,kR^{-2}].
\label{eq:Examp3_est}
\end{equation}

Hence there is a $N=N(R,\tau)$ so that the estimate \eqref{eq:Examp3_est}
holds for all $z\in B_R$, all $k\geq N$ and for all $t\in[0,\tau]$.
Thus, for all $z\in B_R$ and all $t\in[0,\tau]$ we have
\begin{eqnarray*}
\|z(t)\|^2_X &=& \sum_{k=1}^{N-1}|z_k(t)|^2 + \sum_{k=N}^{\infty}|z_k(t)|^2 \\
&\leq& \sum_{k=1}^{N-1}|z_k(t)|^2 + 4\sum_{k=N}^{\infty}|z_k(0)|^2 \\
&\leq& \sum_{k=1}^{N-1}|z_k(t)|^2 + 4R^2.
\end{eqnarray*}
Since $N$ is finite and depends on $R$ and $\tau$ only, and since every $z_k$-subsystem is GAS, it is clear that
\[
\sup \{ \|z(t)\|_{\ell_2}  \midset \|z(0)\|_{\ell_2}\leq R,\ t\in[0,\tau]\} <\infty,
\]
so that $\tilde\Sc^1$ is BRS.

To show that $\tilde\Sc^1$ is not 0-UGS, recall the construction in Example~\ref{0-GAS_but_not_GS}.

Let $c>0$ be such that the solution of $\dot{x}_k = -2x_k + x_k^2$ starting at $x_k(0)=c$, blows up to infinity in time $t^*=1$. 
Consider an initial state $z(0)$ for $\tilde\Sc^1$, where
$\tilde{z}_k(0)=(c,e)^T$ and $\tilde{z}_j(0)= (0,0)^T$ for $j\neq k$.
For this initial state we have that $\|\tilde{z}(t)\|_{\ell_2} = |\tilde{z}_k(t)| \geq |\tilde{x}_k(t)| = |x_k(\frac{t}{k})|$.
And hence
\[
\sup_{t \geq 0}\|\tilde{z}(t)\|_{\ell_2} \geq |\tilde{x}_k(k\tau_k)| = |x_k(\tau_k)|  \geq  k.
\]
As $k\in\N$ was arbitrary, this shows that $\tilde\Sc^1$ is not 0-UGS.

\textbf{Step 3.} Let $c$ be as in Step 2. Let $a:=\min\{c,\tfrac{1}{2}\}$
and choose a smooth function $\xi:\R\to\R$ with
\[
\xi(s):=
\begin{cases}
s & \text{, if } |s| \leq \frac{a}{4}, \\
0&\text{, if } |s| > \frac{a}{2}, \\
\in [-\frac{a}{2},-\frac{a}{4}]\cup [\frac{a}{4},-\frac{a}{2}]  & \text{, otherwise. }
\end{cases}
\]
Now consider the modification of $\tilde\Sc^1$, which we denote $\Sc^3$.
\begin{eqnarray}
\Sc^3: \left
\{\begin{array}{l}
\Sc^3_k:
\left
\{\begin{array}{l}
\dot{\tilde{x}}_k = -\xi(\tilde{x}_k) +  \frac{1}{k}\big(-\tilde{x}_k + \tilde{x}_k^2 \tilde{y}_k - \frac{1}{k^2} \tilde{x}_k^3\big),\\
\dot{\tilde{y}}_k = -\xi(\tilde{y}_k) - \frac{1}{k} \tilde{y}_k.
\end{array}
\right.
\\
k=1,2,\ldots
\end{array}
\right.
\end{eqnarray}
The additional dynamics generated by $\xi$ improve the stability
properties of $\Sc^3$ compared to $\tilde\Sc^1$.
In particular, since $\tilde\Sc^1$ is forward complete, 0-GAS, BRS,
$\Sc^3$ also has these properties.
Moreover, in a neighborhood of the origin, the dynamics of $\Sc^3_k$ are dominated by the terms $-\xi(\tilde{x}_k) = -\tilde{x}_k$ and $-\xi(\tilde{y}_k) = -\tilde{y}_k$, which renders $\Sc^3$ 0-UAS.
This can be justified, e.g., by a Lyapunov argument, as in Example~\ref{0-GAS_but_not_GS}.

Now, since $\xi(s)=0$ for $s > \frac{a}{2}$, the argument, used to show that $\tilde\Sc^1$ is not 0-UGS, shows that $\Sc^3$ is again not 0-UGS.
\hfill \qedsymbol
\end{example}

\begin{example}[FC\,$\wedge$\,BRS\,$\wedge$\,0-UGAS\,$\wedge$\,AG\,$\wedge$\,LISS\ \ $\not\Rightarrow$\ \ UGS]
\label{examp:FC_BRS_0UGAS_AG_LISS_not_UGS}
Consider the system $\Sc^4$ with the state space $\ell_2$ (see \eqref{eq:l2_space}) and its input space be $\Uc:=PC_b(\R_+,\R)$.
\begin{eqnarray}
\Sc^4: \left
\{\begin{array}{l}
\Sc^4_k:
\left
\{\begin{array}{l}
\dot{\tilde{x}}_k = -\xi(\tilde{x}_k) +  \frac{1}{k}\big(-\tilde{x}_k + \tilde{x}_k^2 \tilde{y}_k |u| - \frac{1}{k^2} \tilde{x}_k^3\big),\\
\dot{\tilde{y}}_k = -\xi(\tilde{y}_k) - \frac{1}{k} \tilde{y}_k.
\end{array}
\right.
\\
k=1,2,\ldots
\end{array}
\right.
\label{eq:GAS_BRS_AG_not_UGS}
\end{eqnarray}
Since this example is a combination of the two previous ones, we omit all details and only mention that
$\Sc^4$ is forward complete, BRS, 0-UGAS, LISS, and AG with zero gain, but at the same time, $\Sc^4$ is not UGS.
\hfill \qedsymbol
\end{example}

\section{Systems without inputs}
\label{No_disturbances}

In this section, we classify the stability notions for abstract systems
$\Sigma=(X,\{0\},\phi)$
without inputs. This simplified picture can help to understand the general case, and at the same time, it is interesting in its own right.

 
\begin{figure*}[tbh]
\centering
\begin{tikzpicture}[>=implies,thick]


\node (UGAS) at (4.3,5.5) {0-UGAS};
\node (UGATT) at (0.3,5.5) {0-UGATT};
\node (ULIM_ULS) at (-3.5,5.5) {0-ULIM\,$\wedge$\,0-ULS};
\node (LF) at (7.2,5.5) {$\exists$ coercive LF};
\node (ncLF) at (10.2,5.5) {$\exists$ nc LF};

\node (UAS_GATT) at (4.3,4.5) {0-GATT\,$\wedge$\,0-UAS};
\node (GAS_GS) at  (0.3,4.5) {0-GATT\,$\wedge$\,0-UGS};
\node (LIM_GS) at (-3.5,4.5) {0-LIM\,$\wedge$\,0-UGS};

\node (GAS) at (0.3,3.5) {0-GAS};
\node (LIM_LS) at (-3.5,3.5) {0-LIM\,$\wedge$\,0-ULS};
\node (GATTLS) at (4.3,3.5) {0-GATT\,$\wedge$\,0-ULS};


\draw[->,degil,double equal sign distance] (2,4.25) to (2.6,4.25);
\draw[<-,degil,double equal sign distance] (2,4.75) to (2.6,4.75);


\node (S-1) at (4.85,5) {{\scriptsize(1)(2)}};
\draw[->,double equal sign distance] (UGAS) to (UAS_GATT);
\node (S1) at (4.75,4) {\scriptsize(1)};
\draw[->,double equal sign distance] (UAS_GATT) to (GATTLS);

\node (S0) at (-0.15,5) {\scriptsize(1)};
\draw[->,double equal sign distance] (UGATT) to (GAS_GS);

\node (S1.2) at (-0.25,4) {{\scriptsize(1)(2)}};
\draw[->,double equal sign distance] (GAS_GS) to (GAS);


\draw[<->,double equal sign distance](ULIM_ULS) to (UGATT);
\draw[<->,double equal sign distance] (UGAS) to (UGATT);

\draw[<->,double equal sign distance] (LF) to (ncLF);
\draw[<->,double equal sign distance]  (UGAS) to (LF);

\draw[<->,double equal sign distance] (GAS_GS) to (LIM_GS);

\draw[<->,double equal sign distance] (GAS) to (GATTLS);
\draw[<->,double equal sign distance] (GAS) to (LIM_LS);

\end{tikzpicture}

\caption{Characterizations of 0-UGAS for systems, satisfying BRS and CEP properties. Here \q{nc LF} means non-coercive Lyapunov function. Implications marked by (1) resp. (2)
  become equivalences for
  {\scriptsize (1)} ODE systems, see e.g. \cite[Proposition 2.5]{LSW96} and
{\scriptsize (2)} linear systems (as a consequence of the Banach-Steinhaus theorem).}
\label{UGAS_Equiv}
\end{figure*}

\begin{lemma}
\label{0LIM_0LS_0GAS}
$\Sigma$ is 0-LIM and 0-ULS if and only if $\Sigma$ is 0-GAS.
\end{lemma}

\begin{proof}
It is clear that 0-GAS implies 0-LIM and 0-ULS. Hence, we only prove the
converse direction.

Pick any $\eps_1>0$. Since $\Sigma$ is 0-ULS, there is a
$\delta_1=\delta_1(\eps_1)>0$ so that  $\|x\|_X \leq \delta_1$ implies
$\|\phi(t,x,0)\|_X \leq \eps_1$ for all $t \geq 0$.

Pick any $x \in X$. Since $\Sigma$ is 0-LIM, there exists a
$T_1=T_1(x)>0$  such that  $\|\phi(T_1,x,0)\|_X \leq \delta_1$.
By the semigroup property, $\phi(t+T_1,x,0)=\phi(t,\phi(T_1,x,0),0)$ and consequently $\|\phi(t+T_1,x,0)\|_X \leq \eps_1$ for all $t \geq 0$.

Pick a sequence $(\eps_i)$ with $\eps_i \to 0$ as $i \to \infty$. According to the above argument, there exists a sequence of times  $T_i=T_i(x)$ such that $\|\phi(t,x,0)\|_X \leq \eps_i$ for all $t \geq T_i$, and thus $\|\phi(t,x,0)\|_X \to 0$ as $t \to \infty$. This shows that $\Sigma$ is 0-GATT, and since we assumed that $\Sigma$ is 0-ULS,
$\Sigma$ is also 0-GAS.
\end{proof}

Now we can state the main result of this section.
\begin{proposition}
\label{Main_Prop_Undisturbed_Systems}
For the system $\Sigma$ without inputs, the relations depicted in Figure~\ref{UGAS_Equiv} hold.
\end{proposition}

\begin{proof}
The equivalences on the uniform level follow directly from the equivalence between UAG and ISS, as well as from
Theorem~\ref{t:ISSLyapunovtheorem}.
By definition, 0-GAS is equivalent to 0-GATT $\wedge$ 0-ULS, and it is equivalent to 0-LIM $\wedge$ 0-ULS according to Lemma~\ref{0LIM_0LS_0GAS}.

The implications (2) follow since 0-UAS $\Leftrightarrow$ 0-UGAS and 0-ULS $\Leftrightarrow$ 0-UGS for linear systems.
Finally, (1) is well-known.

The observation that 0-UAS $\wedge$ 0-GATT is not implied by and does not imply 0-GAS $\wedge$ 0-UGS follows from Example~\ref{0-GAS_but_not_GS} and since the strong stability of strongly continuous semigroups is weaker than exponential stability.
\end{proof}

\section{Concluding remarks}

%

The notions of an asymptotic gain, uniform asymptotic gain, and limit property are due to \cite{SoW95, SoW96}.
The notions of strong asymptotic gain, strong limit property, and uniform limit property have been introduced in \cite{MiW18b}.
Uniform limit property with zero gain (uniform weak attractivity) has been introduced in \cite{MiW19a} and thoroughly studied in \cite{Mir17a, MiW19b}. It is motivated by the concept of weak attractivity, which is classical in the dynamical systems theory \cite{Bha66}, \cite{BhS02}.

Theorem~\ref{thm:MainResult_Characterization_ISS} has been first proved in \cite{MiW18b}. 
The ISS superposition theorem (Theorem~\ref{thm:UAG_equals_ULIM_plus_LS}) was shown in \cite[Theorem 5]{MiW18b} and is based on the classical characterizations of ISS for ODE systems \cite{SoW96}. 
The bUAG property was formally introduced in \cite{Mir19b}, but the idea of using this kind of attractivity was also used previously, e.g., in \cite{Tee98}. Using the bUAG property, in \cite{Mir19b}, the ISS characterization from \cite[Theorem 5]{MiW18b} was somewhat strengthened in \cite[Theorem 6.7]{Mir21}. We show in Theorem~\ref{thm:UAG_equals_ULIM_plus_LS} the ISS superposition theorem in this latter formulation.
 
Integral characterizations of ISS shown in Section~\ref{sec:Relations between ISS and norm-to-integral ISS} are due to \cite{JMP20}. For ODEs, the corresponding results have been shown in \cite{Son98}. The superposition theorem for the strong ISS property (Theorem~\ref{wISS_equals_sAG_GS}) is due to \cite[Theorem 12]{MiW18b}. Strong ISS of linear infinite-dimensional systems was analyzed in \cite{NaS18}.

Counterexamples are taken mostly from \cite{Mir16,MiW18b}.
All our counterexamples are given for \q{infinite ODEs}, also called \q{ensembles}. 
Thus, it is possible that the ISS superposition Theorem~\ref{thm:UAG_equals_ULIM_plus_LS} can be strengthened for other particular classes of infinite-dimensional systems. For time-delay systems, some progress has been achieved in \cite{MiW17e}, \cite{KPC22}. See also \cite{CKP22} for a survey on ISS of time-delay systems.

Characterizations of ISS greatly simplify the proofs of other important results, such as small-gain theorems for ODEs \cite{DRW07} and hybrid systems \cite{CaT09, DaK13}, Lyapunov-Razumikhin theory for time-delay systems \cite{Tee98}, \cite{DKM12},
 non-coercive ISS Lyapunov theorems \cite{MiW18b, JMP20}, relations between ISS and nonlinear $L^2 \to L^2$ stability \cite{Son98}, to name a few examples.

Input-to-state practical stability (ISpS) has been introduced in \cite{JTP94, JMW96}, and extends an earlier concept of practical asymptotic stability of dynamical systems \cite{LLM90}. 
The ISpS superposition theorem (Theorem~\ref{thm:ncISpS_LF_sufficient_condition_NEW}) for infinite-dimensional systems has been shown in \cite{Mir19a}. Some of the characterizations derived in \cite{Mir19a} were previously shown in the special case of ODE systems in \cite{SoW96}, and some other ones are new already for the ODE case. 
ISpS is extremely useful for the stabilization of stochastic control systems \cite{ZhX13}, control under quantization errors \cite{ShL12, JLR09}, sample-data control \cite{NKK15}, the study of interconnections of nonlinear systems by means of small-gain theorems \cite{JTP94,JMW96}, etc.

As we have already seen, the relationships between the stability notions depend on the type of dynamics (linear, nonlinear), the dimension of the state space (finite or infinite), on the type of the state space (Banach space, Hilbert space, normed vector space, etc.), and ultimately on the particular class of systems.

For example, already the fact that the state space is infinite-dimensional implies that the closed bounded balls are not compact, which breaks a great deal of the proofs which have been devised in the ISS theory for ODEs.
In particular, the equivalence between forward completeness and BRS property, as well as between LIM and ULIM properties, fails (in general) for infinite-dimensional spaces with Lipschitz right-hand sides.
At the same time, despite the infinite-dimensionality of the state space, linear systems over Banach spaces with admissible operators still enjoy the BRS property, and attractivity (strong stability of a semigroup for zero controls) implies Lyapunov stability.
However, linear infinite-dimensional systems over normed vector spaces, which are not Banach spaces, already lose some of these nice properties \cite{JaW15}.


\cleardoublepage
\chapter{Lyapunov criteria}
\label{chap:Characterizations_ISS_2}

Lyapunov functions is a key tool for the stability analysis of nonlinear systems. 
In Section~\ref{sec:ISS_Lyapunov_functions}, we have already introduced the concept of a coercive ISS Lyapunov function, shown that the existence of such a function implies ISS of a system, and used this concept to verify ISS of several important (non)linear PDE systems.
In this chapter, we undertake an in-depth analysis of Lyapunov methods. 
Using ISS superposition theorems, we show that in many cases, the coercivity of ISS Lyapunov functions can be dropped, and already existence of a non-coercive ISS Lyapunov function is sufficient for ISS of a nonlinear abstract control system.
We show that for evolution equations with nonlinearities that are Lipschitz both with respect to states and inputs, ISS is equivalent to the existence of an ISS Lyapunov function.
Finally, we give an explicit construction of ISS Lyapunov functions for linear systems with bounded input operators and (under certain restrictions) with unbounded input operators.

\section{Non-coercive ISS Lyapunov functions}
\label{sec:Non-coercive ISS Lyapunov functions}

In Theorem~\ref{LyapunovTheorem},  we have seen that the existence of an ISS Lyapunov function implies ISS.
However, the construction of ISS Lyapunov functions for infinite-dimensional systems, especially nonlinear ones, is a challenging task. Already for undisturbed linear systems over Hilbert spaces, \q{natural} Lyapunov function candidates constructed via solutions of Lyapunov equations are of the form $V(x):= \lel Px,x \rir$, where $\lel\cdot ,\cdot \rir$ is a scalar product in $X$ and $P$ is a linear bounded positive operator whose spectrum may contain $0$. Such $V$ fail to satisfy the bound from below in \eqref{LyapFunk_1Eig}, and possess only the weaker property $V(x)>0$ for $x\neq 0$.
In this section, we show how non-coercive ISS Lyapunov functions can be used to show ISS.

\begin{definition}
\label{def:noncoercive_ISS_LF}
\index{Lyapunov function!non-coercive ISS}
Consider a control system $\Sigma = (X,\Uc,\phi)$.
A continuous function $V:X \to \R_+$ is called a \emph{non-coercive ISS Lyapunov function} for $\Sigma$, if there exist $\psi_2,\alpha \in \Kinf$ and $\sigma \in \K$ such that
\begin{equation}
\label{LyapFunk_1Eig_nc_ISS}
0 < V(x) \leq \psi_2(\|x\|_X), \quad \forall x \in X \setminus \{0\},
\end {equation}
and the Dini derivative of $V$ along the trajectories of $\Sigma$ for all $x \in X$ and $u\in \Uc$ satisfies
\begin{equation}
\label{DissipationIneq_nc}
\dot{V}_u(x) \leq -\alpha(\|x\|_X)  + \sigma(\|u\|_{\Uc}).
\end{equation}

Moreover, if \eqref{DissipationIneq_nc} holds just for $u=0$, we call $V$ a \emph{non-coercive Lyapunov function} for the undisturbed system $\Sigma$.
\end{definition}

Note that continuity of $V$ and the estimate \eqref{LyapFunk_1Eig_nc_ISS} imply that $V(0)=0$.

The next proposition shows that the norm-to-integral ISS property arises naturally in the theory of ISS Lyapunov functions:
\begin{proposition}
\label{prop:ncLF_implies_norm-to-integral-ISS} 
Let $\Sigma=(X,\Uc,\phi)$ be a forward complete control system. Assume
that there exists a non-coercive ISS Lyapunov function for $\Sigma$. Then $\Sigma$ is norm-to-integral ISS.
\end{proposition}

\begin{proof}
Assume that $V$ is a non-coercive ISS Lyapunov function for $\Sigma$ with
the corresponding $\psi_2,\alpha,\sigma$.
Pick any $u\in\Uc$ and any $x\in X$. As we assume forward completeness of $\Sigma$, the trajectory $\phi(\cdot,x,u)$ exists for all $t\geq0$ and due to \eqref{DissipationIneq_nc}, we have for any $t>0$ that:
\begin{equation}
\label{DissipationIneq_nc-specified}
\dot{V}_{u(t+\cdot)}\big(\phi(t,x,u)\big) \leq -\alpha(\|\phi(t,x,u)\|_X)  + \sigma(\|u(t+\cdot)\|_{\Uc}).
\end{equation}
By definition of $\dot{V}$, and using the cocycle property for $\Sigma$, we have that
\begin{align*}
\dot{V}_{u(t+\cdot)}\big(\phi(t,x,u)\big)
&=\mathop{\overline{\lim}} \limits_{h \rightarrow +0} {\frac{1}{h}\Big(V\big(\phi(h,\phi(t,x,u),u(t+\cdot))\big)-V\big(\phi(t,x,u)\big)\Big) }\\
&=\mathop{\overline{\lim}} \limits_{h \rightarrow +0} {\frac{1}{h}\Big(V\big(\phi(t+h,x,u)\big)-V\big(\phi(t,x,u)\big)\Big) }.
\end{align*}
Defining $y(t):=V\big(\phi(t,x,u)\big)$, we see that 
\begin{eqnarray}
\dot{V}_{u(t+\cdot)}\big(\phi(t,x,u)\big) = D^+ y(t),
\label{eq:Lie-der-and-Dini-derivative}
\end{eqnarray}
and $y(0)= V(x)$ due to the identity axiom of the system $\Sigma$.

In view of the continuity axiom of $\Sigma$, for fixed $x,u$ the map $\phi(\cdot,x,u)$ is continuous, and thus $t\mapsto -\alpha(\|\phi(t,x,u)\|_X)$ is continuous as well.

For $t\geq 0$, define $G(t):=\int_0^t \alpha(\|\phi(s,x,u)\|_X) ds$ and
$b(t):=\sigma(\|u(t+\cdot)\|_{\Uc})$.
Note that by the axiom of shift invariance, $b$ is nonincreasing.
As $G$ is continuously differentiable, we can rewrite the inequality \eqref{DissipationIneq_nc-specified} as
\begin{eqnarray}
D^+y(t) \leq - \frac{d}{dt}G(t) + b(t). 
\label{eq:Shorthand-notation-dissip-ineq}
\end{eqnarray}
Pick any $r>0$ and define $b(s)=b(0)$ for $s \in [-r,0]$. As $b$ is a
nonincreasing function on $[-r,\infty)$, it holds for any $t\geq 0$ that $b(t) \leq
\lim_{h\to +0}b(t-r+h)$, and by the final inequality in
Lemma~\ref{lem:Integrals-of-monotone-functions} applied to $b(\cdot - r)$,
we obtain for all $t\ge 0$
\[
b(t) \leq D_+ \int_0^{t} b(s-r) ds = - D^+ \Big(-\int_0^{t} b(s-r) ds\Big).
\]
Thus, \eqref{eq:Shorthand-notation-dissip-ineq} implies that 
\begin{eqnarray}
D^+y(t) + \frac{d}{dt}G(t) + D^+ \Big(-\int_0^{t} b(s-r) ds\Big) \leq 0.
\label{eq:Shorthand-notation-dissip-ineq-2}
\end{eqnarray}
Due to 
\[
D^+(f_1(t)+f_2(t))\leq D^+(f_1(t)) + D^+(f_2(t)),
\]
which holds for any functions $f_1, f_2$ on the real line, this implies that 
\begin{eqnarray*}
D^+\Big(y(t) + G(t) -\int_0^{t} b(s-r) ds\Big) \leq 0.
\end{eqnarray*}
It follows from 
Proposition~\ref{prop:fundamental theorem of calculus-Dini-derivative}
 that $t\mapsto y(t) + G(t) -\int_0^{t} b(s-r)
ds$ is nonincreasing. As $G(0) = 0$, it follows that for all $r>0$ and all $t\ge 0$
\begin{eqnarray*}
y(t) + G(t) -\int_0^{t} b(s-r) ds \leq y(0) = V(x). 
\end{eqnarray*}
As $b$ is bounded, we may pass to the limit $r \to 0$ and obtain
\begin{eqnarray*}
y(t) + G(t) -\int_0^{t} b(s) ds \leq y(0) = V(x). 
\end{eqnarray*}
Now $y(t)\geq 0$ for all $t\in\R_+$, and so for all $t\ge0$, $x \in X$ and $u\in\Uc$
\begin{align}
\label{eq:old-int-to-int-ISS}
\int_0^t \alpha(\|\phi(s,x,u)\|_X) ds 
&\leq \psi_2(\|x\|_X) + \int_0^{t} \sigma(\|u(s+\cdot)\|_{\Uc}) ds\\
&\leq \psi_2(\|x\|_X) + t \sigma(\|u\|_{\Uc}).\nonumber
\end{align}
This completes the proof.
\end{proof}

\begin{remark}
\label{rem:ncLF-and-integral-to-integral-ISS} 
Similarly to Proposition~\ref{prop:ncLF_implies_norm-to-integral-ISS}, one can show that 
for a forward-complete system $\Sigma=(X, PC_b(\R_+,U), \phi)$ the existence of a non-coercive ISS Lyapunov function implies integral-to-integral ISS.
\end{remark}

We can now state our main result on non-coercive ISS Lyapunov functions.

\begin{theorem}{}
    \label{t:ISSLyapunovtheorem}
Let $\Sigma=(X,\Uc,\phi)$ be a forward complete control system with the input space 
$\Uc:=L^\infty(\R_+,U)$, which is CEP and BRS.    
If there exists a (non-coercive) ISS Lyapunov function for $\Sigma$, then $\Sigma$ is ISS.    
\end{theorem}

\begin{proof}
By Proposition~\ref{prop:ncLF_implies_norm-to-integral-ISS}, $\Sigma$ is norm-to-integral ISS. The application of Theorem~\ref{thm:ncISS_LF_sufficient_condition_NEW} finishes the proof.
\end{proof}

\begin{remark}
\label{rem:Why-Linfty} 
Actually, Theorem~\ref{t:ISSLyapunovtheorem} holds without the restriction of the input space to be $\Uc:=L^\infty(\R_+,U)$ (the proof remains the same in the general case). However, as discussed in Section~\ref{sec:ISS Lyapunov functions and a type of input space}, 
the concept of an ISS Lyapunov function that we have employed is too restrictive for the systems with $\Uc:=L^p(\R_+,U)$, $p\in[1,+\infty)$. To underline this fact, we have added an extra assumption $\Uc:=L^\infty(\R_+,U)$ to the formulation of Theorem~\ref{t:ISSLyapunovtheorem}.
\end{remark}

For evolution equations \eqref{InfiniteDim} satisfying Assumption~\ref{Assumption1}, Theorem~\ref{t:ISSLyapunovtheorem} can be strengthened:
\begin{theorem}
\label{thm:ncISS_LF_sufficient_condition}
Let \eqref{InfiniteDim} satisfy Assumption~\ref{Assumption1} and have BRS property.
If there exists a non-coercive ISS Lyapunov function for \eqref{InfiniteDim}, then \eqref{InfiniteDim} is ISS.
\end{theorem}

\begin{proof}
Lemma~\ref{lem:RobustEquilibriumPoint} ensures that \eqref{InfiniteDim} has CEP property. 
Theorem~\ref{t:ISSLyapunovtheorem} shows the claim.
\end{proof}

\section{Alternative definitions of ISS Lyapunov functions}

In this section, let $\Uc:=L^\infty(\R_+,U)$. In this case, our definition of an ISS Lyapunov function seems to be the most natural. 
For forward complete finite-dimensional systems with $f$ as in the statement of Corollary~\ref{cor:ni-ISS-and-ISS-ODEs}, existence of a non-coercive ISS Lyapunov function implies not only norm-to-integral ISS, but also integral-to-integral ISS, which follows from Proposition~\ref{prop:ncLF_implies_norm-to-integral-ISS} and 
 Corollary~\ref{cor:ni-ISS-and-ISS-ODEs}.
However, for infinite-dimensional systems, we were able to show only a somewhat weaker property (which is still stronger than norm-to-integral ISS but weaker than integral-to-integral ISS), see \eqref{eq:old-int-to-int-ISS}. Thus, a question remains whether the existence of an ISS Lyapunov function (coercive or non-coercive) as defined in Definition~\ref{def:noncoercive_ISS_LF} implies integral-to-integral ISS for forward complete systems.
Although we do not have an answer to this problem, in the following proposition, we show that if an ISS Lyapunov function satisfies a somewhat stronger dissipative estimate, then integral-to-integral ISS can be verified.

\begin{proposition}
\label{prop:ncLF_implies_ii_ISS} 
Let $\Sigma=(X,\Uc,\phi)$ be a forward complete control system with the input space $\Uc:=L^\infty(\R_+,U)$. 
Assume that there is a continuous function $V:X \to \R_+$, $\psi_2,\alpha \in \Kinf$ and $\sigma \in \K$ such that
\eqref{LyapFunk_1Eig_nc_ISS} holds and the Dini derivative of $V$ along the trajectories of $\Sigma$ for all $x \in X$ and $u\in \Uc$ satisfies
\begin{equation}
\label{DissipationIneq_nc-Linfty-stronger-form}
\dot{V}_u(x) \leq -\alpha(\|x\|_X)  + \Big[D_{+,\tau}\Big(\int_0^\tau\sigma(\|u(s)\|_{U})ds\Big)\Big]_{\tau=0},
\end{equation}
where $D_{+,\tau}$ means that the lower right-hand Dini derivative is taken with respect to the argument $\tau$.

Then $\Sigma$ is integral-to-integral ISS.
\end{proposition}

Before we prove this proposition, note that the estimate \eqref{DissipationIneq_nc-Linfty-stronger-form} implies \eqref{DissipationIneq_nc}, since for any $x \in X$ and $u\in \Uc$ it holds that 
\begin{eqnarray*}
\Big[D_{+,\tau}\Big(\int_0^\tau\sigma(\|u(s)\|_{U})ds\Big)\Big]_{\tau=0}
\leq
\Big[D_{+,\tau}\Big(\int_0^\tau\sigma(\|u\|_{\Uc})ds\Big)\Big]_{\tau=0}
= \sigma(\|u\|_{\Uc}).
\end{eqnarray*} 
Thus, function $V$ as in Proposition~\ref{prop:ncLF_implies_ii_ISS} is an ISS Lyapunov function for $\Sigma$, with a (potentially) stronger dissipative estimate.

\begin{proof}
Assume that $V$ is a non-coercive ISS Lyapunov function for $\Sigma$ with
the corresponding $\psi_2,\alpha,\sigma$.
Pick any $u\in\Uc$ and any $x\in X$. As we assume forward completeness of $\Sigma$, the trajectory $\phi(\cdot,x,u)$ exists for all $t\geq0$ and due to \eqref{DissipationIneq_nc}, we have for any $t>0$ that:
\begin{equation}
\label{DissipationIneq_nc-specified-i-to-i-ISS}
\dot{V}_{u(t+\cdot)}\big(\phi(t,x,u)\big) \leq -\alpha(\|\phi(t,x,u)\|_X) + \Big[D_{+,\tau}\Big(\int_0^\tau\sigma(\|u(t+s)\|_{U})ds\Big)\Big]_{\tau=0}.
\end{equation}
The last term can be rewritten in a simpler form:
\begin{align*}
\Big[D_{+,\tau}\Big(\int_0^\tau\sigma(\|u(t&+s)\|_{U})ds\Big)\Big]_{\tau=0}
=\Big[D_{+,\tau}\Big(\int_t^{t+\tau}\sigma(\|u(s)\|_{U})ds\Big)\Big]_{\tau=0}\\
&=\Big[D_{+,\tau}\Big(\int_0^{t+\tau}\sigma(\|u(s)\|_{U})ds\Big)\Big]_{\tau=0}
=D_{+,t}\Big(\int_0^{t}\sigma(\|u(s)\|_{U})ds\Big).
\end{align*}
Arguing as in Proposition~\ref{prop:ncLF_implies_norm-to-integral-ISS}, we define 
$y(t):=V\big(\phi(t,x,u)\big)$ and verify that the equality  
$\dot{V}_{u(t+\cdot)}\big(\phi(t,x,u)\big) = D^+ y(t)$ 
holds and $y(0)= V(x)$.

In view of the continuity axiom of $\Sigma$, for fixed $x,u$ the map $\phi(\cdot,x,u)$ is continuous, and thus $t\mapsto -\alpha(\|\phi(t,x,u)\|_X)$ is continuous as well. Hence, we can rewrite the inequality \eqref{DissipationIneq_nc-specified-i-to-i-ISS} as
\begin{eqnarray*}
D^+y(t) \leq - \frac{d}{dt}\int_0^t \alpha(\|\phi(s,x,u)\|_X) ds + D_+\Big(\int_0^t\sigma(\|u(s)\|_{U})ds\Big). 
\end{eqnarray*}
As $-D_+\big(\int_0^t\sigma(\|u(s)\|_{U})ds\big) = D^+\big(-\int_0^t\sigma(\|u(s)\|_{U})ds\big)$, we proceed to
\begin{eqnarray*}
D^+y(t) + \frac{d}{dt}\int_0^t \alpha(\|\phi(s,x,u)\|_X) ds + D^+ \Big(-\int_0^t\sigma(\|u(s)\|_{U})ds\Big) \leq 0.
\end{eqnarray*}
As in the proof of Proposition~\ref{prop:ncLF_implies_norm-to-integral-ISS}, by subadditivity of $D^+$, we obtain that
\begin{eqnarray*}
D^+\Big(y(t) + \int_0^t \alpha(\|\phi(s,x,u)\|_X) ds -\int_0^t\sigma(\|u(s)\|_{U})ds\Big) \leq 0.
\end{eqnarray*}
It follows from 
Proposition~\ref{prop:fundamental theorem of calculus-Dini-derivative} that $t\mapsto y(t) + \int_0^t \alpha(\|\phi(s,x,u)\|_X) ds -\int_0^t\sigma(\|u(s)\|_{U})ds$ 
is nonincreasing. 
For $t\geq 0$, define $G(t):=\int_0^t \alpha(\|\phi(s,x,u)\|_X) ds$.
As $G(0) = 0$, for all $r>0$ we have
\begin{eqnarray*}
y(t) + G(t) -\int_0^t\sigma(\|u(s)\|_{U})ds \leq y(0) = V(x)\leq  \psi_2(\|x\|_X). 
\end{eqnarray*}
Now $y(t)\geq 0$ for all $t\in\R_+$, and so
\begin{align*}
\int_0^t \alpha(\|\phi(s,x,u)\|_X) ds 
\leq \psi_2(\|x\|_X) + \int_0^t\sigma(\|u(s)\|_{U})ds.
\end{align*}
This completes the proof.
\end{proof}

\section{Characterization of local input-to-state stability}
\label{sec:Main_Result_AND_Structure}

Consider infinite-dimensional evolution equations of the form
\begin{subequations}
\label{InfiniteDim} 
\begin{eqnarray}
\dot{x}(t) & = & Ax(t) + f(x(t),u(t)),\quad t>0,  \label{InfiniteDim-1}\\
x(0)  &=&  x_0, \label{InfiniteDim-2}
\end{eqnarray}
\end{subequations}
where $A: D(A)\subset X \to X$ generates a strongly continuous semigroup $\sg{T}$ of boun\-ded linear operators on a Banach space $X$; $U$ is a normed vector space of input values, $x_0\in X$ is a given initial condition, and $f:X\times U \to X$.
As the space of admissible inputs, we consider the space $\Uc:=PC_b(\R_+,U)$ of globally bounded, piecewise continuous functions from $\R_+$ to $U$ that are right-continuous.

Clearly, systems \eqref{InfiniteDim} are a special case of a more general class of systems \eqref{eq:SEE+admissible} studied previously.

In the rest of the chapter, we analyze nonlinear systems of the form \eqref{InfiniteDim} that 
satisfy Assumption~\ref{Assumption1} with $V=X$. As we have shown in Chapter~\ref{chap:Nonlinear Evolution Equations}, such systems give rise to well-posed control systems $(X,\Uc,\phi)$, where  $\phi(\cdot,x,u)$ is the unique maximal (mild) solution of \eqref{InfiniteDim} associated with 
an initial condition $x \in X$ at $t=0$, and input $u \in \Uc$.

In this section, we derive a characterization of LISS property for the systems of the form \eqref{InfiniteDim}.
In addition to Assumption~\ref{Assumption1}, we suppose that $f(0,0)=0$, i.e., $x \equiv 0$ is an equilibrium point of \eqref{InfiniteDim}.

We rely upon the converse Lyapunov theorem for undisturbed systems \eqref{InfiniteDim}. To this end, define:
\begin{definition}
\label{def:0-UAS}
A continuous function $V:D \to \R_+$, $0 \in \intt(D) \subset X$,  is called a \emph{0-UAS Lyapunov function},  if there exist $r>0$, $\psi_1,\psi_2 \in \Kinf$ and $\alpha \in \Kinf$ such that $\{x \in X: \|x\|_X \leq r \} \in D$  and
$\forall x \in X: \|x\|_X \leq r$ the following estimate holds
\begin{equation}
\label{LyapFunk_1Eig_UAS}
\psi_1(\|x\|_X) \leq V(x) \leq \psi_2(\|x\|_X),
\end{equation}
and the Lie derivative of $V$ along the trajectories of the system \eqref{InfiniteDim} with  $u = 0$ satisfies
\begin{equation}
\label{DissipationIneq_0_UAS}
\dot{V}_0(x) \leq -\alpha(\|x\|_X)
\end{equation}
for all $x \in X: \|x\|_{X} \leq r$.
\end{definition}

The classical converse Lyapunov theorem for undisturbed systems \eqref{InfiniteDim} is formulated as follows, see, e.g., \cite[Theorem 4.2.1]{Hen81}:
\begin{proposition}
\label{UAS_ConverseLT}
Let $f(\cdot,0)$ be Lipschitz continuous in a certain neighborhood of $x=0$.
If \eqref{InfiniteDim} is 0-UAS, then there exists a Lipschitz continuous 0-UAS Lyapunov function for \eqref{InfiniteDim}.
\end{proposition}

\begin{remark}
Note that in  \cite[Theorem 4.2.1]{Hen81} analyticity of a semigroup is assumed. However, this assumption has not been used in the proof of \cite[Theorem 4.2.1]{Hen81} and was made just because \cite{Hen81} is devoted to analytic semigroups. Thus, \cite[Theorem 4.2.1]{Hen81} is still true for merely strongly continuous semigroups.
\end{remark}

The next result (that is closely related to the known fact about the robustness of the 0-UAS property \cite[Corollary 4.2.3]{Hen81}) shows that a Lipschitz continuous 0-UAS Lyapunov function for \eqref{InfiniteDim} is, under a certain assumption on the nonlinearity $f$, also a LISS Lyapunov function for \eqref{InfiniteDim}. 
\begin{proposition}
\label{Converse_LISS_Lyapunov_Theorem}
Let Assumption~\ref{Assumption1} hold, $f(0,0)=0$, and let there exist $\sigma \in \K$ and $\rho >0$ so that for all $v \in U$: $\|v\|_U \leq \rho$ and all $x \in X$: $\|x\|_X  \leq \rho$ we have
\begin{eqnarray}
\|f(x,v)-f(x,0)\|_X \leq \sigma(\|v\|_U).
\label{eq:Estimate_f_concerning_u}
\end{eqnarray}
Let $V$ be a Lipschitz continuous 0-UAS Lyapunov function for \eqref{InfiniteDim}. Then $V$ is also a LISS Lyapunov function
for \eqref{InfiniteDim}.
\end{proposition}

\begin{proof}
Let $V:D \to \R_+$, $D \subset X$, with $0 \in \intt(D)$ be a Lipschitz continuous (0-UAS) Lyapunov function for \eqref{InfiniteDim}, which satisfies \eqref{DissipationIneq_0_UAS} for all $x \in X$: $\|x\|_X \leq r$.

Let Assumption~\ref{Assumption1} hold and pick $\rho$ given in the formulation of the proposition so that $B_\rho \subset D$.

Lemma~\ref{lem:RobustEquilibriumPoint} implies that the flow of \eqref{InfiniteDim} is continuous at the equilibrium. This implies that there exist $r_1 \in (0,r), r_2 \in (0,\rho)$ and $t^*>0$ so that for all 
$x\in X: \|x\|_X \leq r_1$ and all $u \in \Uc$: $\|u\|_{\Uc} \leq r_2$ the solution $\phi(\cdot,x,u)$ exists on $[0,t^*]$ and $\phi([0,t^*],x,u) \subset B_r \subset D$.

We are going to prove that $V$ is a LISS Lyapunov function for \eqref{InfiniteDim}.
To this end, we derive a dissipative estimate for $\dot{V}_{u}(x)$ for all $x,u$: $\|x\|_X \leq r_1$ and $\|u\|_{\Uc} \leq r_2$. We have:
\begin{eqnarray*}
\hspace{-8mm}
\dot{V}_{u}(x) &=&\mathop{\overline{\lim}} \limits_{t \rightarrow +0} {\frac{1}{t}\Big(V\big(\phi(t,x,u)\big)-V(x)\Big)} \\
				 &=&\mathop{\overline{\lim}} \limits_{t \rightarrow +0} \frac{1}{t}\Big(V\big(\phi(t,x,0)\big)-V(x) \\
				&& \qquad\qquad+ V\big(\phi(t,x,u)\big)-V\big(\phi(t,x,0)\big)\Big) \\
				 &\leq& \dot{V}_{0}(x) + \mathop{\overline{\lim}} \limits_{t \rightarrow +0} {\frac{1}{t}\Big(V\big(\phi(t,x,u)\big)-V\big(\phi(t,x,0)\big)\Big)}
\end{eqnarray*}
Since $V$ is a 0-UAS Lyapunov function for \eqref{InfiniteDim}, due to \eqref{DissipationIneq_0_UAS}, it holds for a certain $\alpha \in \Kinf$ that
\begin{eqnarray*}
\dot{V}_{u}(x) &\leq& -\alpha(\|x\|_X) +\mathop{\overline{\lim}} \limits_{t \rightarrow +0} {\frac{1}{t} \Big|V\big(\phi(t,x,u)\big)-V\big(\phi(t,x,0)\big)\Big|}.
\end{eqnarray*}
Since $\phi(t,x,u) \in D$ for all $x,u$: $\|x\|_X \leq r_1$ and $\|u\|_{\Uc} \leq r_2$, and since $V$ is Lipschitz continuous on $D$, there exists $L>0$ so that 
\begin{eqnarray}
\hspace{-8mm}
\dot{V}_{u}(x) \leq -\alpha(\|x\|_X) + L \mathop{\overline{\lim}} \limits_{t \rightarrow +0} {\frac{1}{t} \|\phi(t,x,u)-\phi(t,x,0)\|_X}.
\label{LISS_Converse_Proof_Estimate}
\end{eqnarray}
The variation of constants formula implies: 
\begin{align*}
\|&\phi(t,x,u)-\phi(t,x,0)\|_X \\
&= \Big\| \int_0^t T(t-s)\big( f(\phi(s,x,u),u(s)) - f(\phi(s,x,0),0) \big)ds \Big\|_X \\
&\leq  \int_0^t \big\|T(t-s)\big\|\big\| f\big(\phi(s,x,u),u(s)\big) - f\big(\phi(s,x,0),0\big) \big\|_X ds.
\end{align*}
Denote $M:=\sup_{0\leq s \leq t^*} \|T(s)\|< \infty$. With this notation, we proceed to
\begin{align*}
\|\phi(t,x,u)-&\phi(t,x,0)\|_X  \\
&\leq  \int_0^t M  \Big( \big\| f\big(\phi(s,x,u),u(s)\big) - f\big(\phi(s,x,0),u(s)\big)\big\|_X \\
&\qquad\qquad+ \big\|f\big(\phi(s,x,0),u(s)\big) - f\big(\phi(s,x,0),0\big) \big\|_X \Big)ds.
\end{align*}
Recalling that $\|\phi(t,x,u)\|_X \leq \rho$ for all $t \in [0,t^*]$, using the inequality \eqref{eq:Estimate_f_concerning_u}  and due to Lipschitz continuity of $f$ w.r.t.\  the first argument, there exists $L_2>0$:
\begin{align*}
\|\phi(t,x,u)-&\phi(t,x,0)\|_X \\
&\leq  M  L_2 \int_0^t \hspace{-1mm} \big\|\phi(s,x,u)  {-} \phi(s,x,0)\big\|_X ds {+} Mt \sigma \big(\hspace{-1mm}\sup_{0 \leq s \leq t}\hspace{-1mm}\|u(s)\|_U\big) \\
&\leq  M L_2t \sup_{0\leq s \leq t} \hspace{-1mm} \|\phi(s,x,u) {-} \phi(s,x,0)\|_X {+} M t \sigma \big(\hspace{-1mm}\sup_{0 \leq s \leq t}\hspace{-1mm}\|u(s)\|_U\big) .
\end{align*}
The right-hand side of the above inequality is nondecreasing in $t$ and consequently
\begin{align*}
\mathop{\sup}_{0\leq s \leq t} \|\phi(s,x,u)-&\phi(s,x,0)\|_X  \\
&\leq  M L_2t \sup_{0\leq s \leq t} \hspace{-1mm} \|\phi(s,x,u)  {-} \phi(s,x,0)\|_X {+} Mt \sigma \big(\hspace{-1mm}\sup_{0 \leq s \leq t}\hspace{-1mm}\|u(s)\|_U\big) .
\end{align*}
Pick $t< t^*$ small enough so that $1 - M L_2 t >0$. Then
\begin{align*}
\sup_{0\leq s \leq t}\hspace{-1mm} \|\phi(s,x,u){-}\phi(s,x,0)\|_X \leq  \frac{M t}{1{-}M L_2 t} \sigma \big(\hspace{-1mm}\sup_{0 \leq s \leq t}\|u(s)\|_U\big) .
\end{align*}
Using this estimate in \eqref{LISS_Converse_Proof_Estimate}, taking the limit $t \to +0$, and recalling that $u$ is a right-continuous function, we obtain for all $x:\, \|x\|_X \leq r_1$ and all $u \in \Uc$: $\|u\|_{\Uc} \leq r_2$ the LISS estimate
\begin{eqnarray}
\dot{V}_{u}(x) \leq -\alpha(\|x\|_X) + LM \sigma \big(\|u(0)\|_U\big).
\label{LISS_Converse_Proof_Estimate_Final}
\end{eqnarray}
This shows that $V$ is a LISS Lyapunov function for \eqref{InfiniteDim}.
\end{proof}

Using the previous proposition, we obtain our main result:
\begin{theorem}
\label{Characterization_LISS}
Let Assumption~\ref{Assumption1} hold (with $V=X$), $f(0,0)=0$, and let there exist $\sigma \in \K$ and $\rho >0$ so that for all $v \in U$: $\|v\|_U \leq \rho$ and all $x \in X$: $\|x\|_X  \leq \rho$ we have
\begin{eqnarray*}
\|f(x,v)-f(x,0)\|_X \leq \sigma(\|v\|_U).
\end{eqnarray*}
Then for the system \eqref{InfiniteDim} the following properties are equivalent:
\begin{enumerate}
	\item[(i)] 0-UAS.
	\item[(ii)] Existence of a Lipschitz continuous 0-UAS Lyapunov function.
	\item[(iii)] Existence of a Lipschitz continuous LISS Lyapunov function.
	\item[(iv)] LISS.
\end{enumerate}
\end{theorem}

\begin{proof}
The claim of the theorem is a consequence of Propositions~\ref{UAS_ConverseLT},~\ref{Converse_LISS_Lyapunov_Theorem}, Theorem~\ref{LyapunovTheorem}, and of the obvious fact that LISS implies 0-UAS.
\end{proof}

\begin{remark}
\label{rem:Attraction_regions_LISS_systems} 
The proof of Theorem~\ref{Characterization_LISS} is based on converse Lyapunov theorems for the 0-UAS property, and thus, it does not give a precise estimate for the region in $X \times \Uc$ for which the LISS estimate \eqref{iss_sum} is valid (for given $\beta,\gamma$).
The estimation of the region in which the LISS estimate is valid is a hard problem, which has not been touched so far for infinite-dimensional systems. 
For ODE systems, this problem has been approached by means of a numerical construction of 
continuous piecewise affine LISS Lyapunov functions using the linear programming method \cite{LBG15}.
\end{remark}

Next, we show by means of an example that the additional assumption in Theorem~\ref{Characterization_LISS} cannot be dropped in infinite dimensions.

\begin{example}[0-UGAS\,$\wedge$\,sAG\,$\wedge$\ AG with zero gain\,$\wedge$\ UGS with zero gain\ \ $\not\Rightarrow$\ \ LISS]
\label{ex:0-UGAS_sAG_AG_zero_gain_not_LISS}
Consider a system $\Sigma$ with the state space $X=l_1:=\{ (x_k): \sum_{k=1}^{\infty} |x_k| <\infty  \}$
and with the input space $\Uc:=PC_b(\R_+,\R)$.

Let the dynamics of the $k$-th mode of $\Sigma$ be given by 
\begin{eqnarray}
\dot{x}_k(t) = -\frac{1}{1+k|u(t)|}x_k(t).
\label{eq:CounterEx_AG_UGAS_no_ISS}
\end{eqnarray}
We use the notation $\phi_k(t,x_k,u)$ for the state of the \mbox{$k$-th} mode of \eqref{eq:CounterEx_AG_UGAS_no_ISS}.
Then $\phi(t,x,u) = (\phi_k(t,x_k,u))$ (we indicate here that the dynamics of different modes are independent of each other).

Clearly, $\Sigma$ is 0-UGAS, since for $u \equiv 0$ its dynamics are given by $\dot{x}=-x$.
At the same time, the inequality $\|\phi(t,x,u)\|_X \leq \|x\|_X$ holds for all $t \geq 0$, $x \in X$ and $u \in \Uc$, and thus $\Sigma$ is UGS with zero gain.

Next, we show step-by-step that:
\begin{enumerate}
	\item[(i)] $\Sigma$ satisfies Assumption~\ref{Assumption1} (and thus belongs to the class of systems which we consider).
	\item[(ii)] $\Sigma$ is AG with zero gain.
	\item[(iii)] $\Sigma$ is sAG with arbitrarily small (but necessarily nonzero) linear gain. In other words, one must pay for uniformity, but this payment can be made arbitrarily small.
	\item[(iv)] Finally, $\Sigma$ is not LISS. 
\end{enumerate}

\textbf{(i):}
The flow of $\Sigma$ is continuous at the equilibrium point since as we mentioned above the system $\Sigma$ is UGS with zero gain. Next, we show that Assumption~\ref{Assumption1} also holds.

Note that $f$ is globally Lipschitz, since for any $x,y \in X$ and any $v \in U$
\begin{align*}
\|f(x,v){-}f(y,v)\|_X=\sum_{k=1}^{\infty} \frac{1}{1+k|v|}|x_k - y_k| &\leq \sum_{k=1}^{\infty} |x_k - y_k|\\
& =  \|x-y\|_X.
\end{align*}
Now pick any $x \in X$ and let us show that $f(x,\cdot)$ is continuous at any $v \neq 0$.
Consider
\begin{align*}
\|f(x,v)-f(x,v_1)\|_X 
&=\sum_{k=1}^{\infty} \Big| \frac{1}{1+k|v_1|} - \frac{1}{1+k|v|}\Big| |x_k| \\
&=\sum_{k=1}^{\infty} \frac{k\big| |v_1|-|v|\big|}{(1+k|v_1|)(1+k|v|)} |x_k| \\
&\leq |v_1-v| \sup_{k \in \N} \frac{k}{(1+k|v_1|)(1+k|v|)}  \sum_{k=1}^{\infty} |x_k|.
\end{align*}
Since $v \neq 0$ and we consider $v_1$ that are close to $v$, we assume that $v_1 \neq 0$ as well.
Then
\begin{align*}
\sup_{k \in \N} \frac{k}{(1+k|v_1|)(1+k|v|)} &\leq \sup_{k \in \N} \frac{k}{1+k(|v_1|{+}|v|)}
=
\lim_{k \to \infty} \frac{k}{1+k(|v_1|{+}|v|)}
=
\frac{1}{|v_1|{+}|v|}.
\end{align*}
Overall, we get that 
\begin{eqnarray*}
\|f(x,v)-f(x,v_1)\|_X \leq  \frac{|v_1-v|}{|v_1|+|v|}  \|x\|_X.
\end{eqnarray*}
Let $|v_1-v| \leq \delta$ for some $\delta \in (0,\tfrac{|v|}{2})$.
Then $|v_1| \geq |v|/2$ and
\begin{eqnarray*}
\|f(x,v)-f(x,v_1)\|_X \leq  \frac{\delta}{3/2|v|}  \|x\|_X.
\end{eqnarray*}
Now for any $\eps >0$ pick $\delta < \tfrac{|v|}{2}$ so that $\frac{\delta}{3/2|v|}  \|x\|_X < \eps$ (which is always possible since $v \neq 0$).
This shows the continuity of $f(x,\cdot)$ for $v \neq 0$.
Note that the choice of $\delta$ depends on $\|x\|_X$ but does not depend on $x$ itself.

Next, we show that $f(x,\cdot)$ is continuous at zero as well, but it is no more uniform w.r.t.\  the first argument of $f$.
Pick again any $x \in X$.
For any $\eps>0$ there exists $N=N(\eps,x)>0$ so that  $\sum_{k=N+1}^{\infty} |x_k| < \frac{\eps}{2}$. We have for this $x$:
\begin{align*}
\|f(x,0)-f(x,v_1)\|_X 
&=  \sum_{k=1}^{\infty} \frac{k|v_1|}{1+k|v_1|} |x_k| \\
&=  \sum_{k=1}^{N} \frac{k|v_1|}{1+k|v_1|} |x_k| +  \sum_{k=N+1}^{\infty}\hspace{-1mm} \frac{k|v_1|}{1+k|v_1|} |x_k| \\
&\leq  \sum_{k=1}^{N} \frac{N|v_1|}{1+N|v_1|} |x_k| +  \sum_{k=N+1}^{\infty}  |x_k| \\
&<  \frac{N|v_1|}{1+N|v_1|} \|x\|_X + \frac{\eps}{2} \\
&\leq  N|v_1| \|x\|_X + \frac{\eps}{2}.
\end{align*}
Pick $\delta := \frac{\eps}{2\|x\|_X N}$.
Then for all $v_1 \in U$: $|v_1| \leq \delta$ it follows that $\|f(x,0)-f(x,v_1)\|_X < \eps$.
This shows the continuity of $f(x,\cdot)$ at zero for any given $x \in X$. However, note that $\delta$ depends on $x$ (since $N$ does) and not only on $\|x\|_X$.
Overall, we have shown that Assumption~\ref{Assumption1} is fulfilled for $\Sigma$.

\textbf{(ii): $\Sigma$ is AG with zero gain.} For any $x \in X$, for any $\eps >0$ there exists $N = N(x,\eps) \in \N$ so that $\sum_{k=N+1}^{\infty} |x_k(0)| <\frac{\eps}{2}$.

The norm of the state of $\Sigma$ at time $t$ can be estimated as follows:
\begin{eqnarray*}
\hspace{-8mm}
\|\phi(t,x,u)\|_X    
&=&  \sum_{k=1}^{N} |\phi_k(t,x_k,u)| + \sum_{k=N+1}^{\infty} |\phi_k(t,x_k,u)| \\
&\leq& \sum_{k=1}^{N} |\phi_k(t,x_k,u)| + \sum_{k=N+1}^{\infty} |\phi_k(0,x_k,u)| \\
&\leq& \sum_{k=1}^{N} |\phi_k(t,x_k,u)| +  \frac{\eps}{2}.
\end{eqnarray*}
Now we estimate the state of the $k$-th mode of our system for all $k = 1,\ldots,N$:
\begin{eqnarray}
\label{eq:CounterEx_kth_mode_estimate}
|\phi_k(t,x_k,u)| &=& e^{-\int_0^t \frac{1}{1+k|u(s)|} ds} |x_k(0)| \nonumber\\
 &\leq& e^{-\int_0^t \frac{1}{1+k\|u\|_{\Uc}} ds} |x_k(0)| \nonumber\\
 &=& e^{- \frac{1}{1+k\|u\|_{\Uc}}t} |x_k(0)| \\
 &\leq& e^{- \frac{1}{1+N\|u\|_{\Uc}}t} |x_k(0)|, \nonumber
\end{eqnarray}
which holds for any $u \in \Uc$ and any $x_k(0) \in \R$. Using this estimate, we proceed to 
\begin{eqnarray*}
\|\phi(t,x,u)\|_X  \leq \sum_{k=1}^{N} e^{- \frac{1}{1+N\|u\|_{\Uc}}t} |x_k(0)| +  \frac{\eps}{2}.
\end{eqnarray*}
Clearly, for any $u \in \Uc$ there exists $\tau_a=\tau_a(x,\eps,u)$ so that $\sum_{k=1}^{N} e^{- \frac{1}{1+N\|u\|_{\Uc}}t} |x_k(0)| \leq \frac{\eps}{2}$ for $t \geq \tau_a$.

Overall, we see that for any $x \in X$, for any $t \geq 0$ and for all $u \in \Uc$ there exists $\tau_a=\tau_a(x,\eps,u)$, so that for all $t \geq \tau_a$ it holds that $\|\phi(t,x,u)\|_X \leq \eps.$
This shows that $\Sigma$ is AG with $\gamma \equiv 0$.
However, the time $\tau_a$ depends on $u$, and thus the above argument doesn't tell us whether the system is strongly AG.

\textbf{(iii):} Next, we show that $\Sigma$ is sAG, but we should pay for this by adding a linear gain. However, this linear gain can be made arbitrarily small. 

In \eqref{eq:CounterEx_kth_mode_estimate}, we have obtained the following estimate for the state of the $k$-th mode of $\Sigma$:
\begin{eqnarray*}
|\phi_k(t,x_k,u)| \leq e^{- \frac{1}{1+k\|u\|_{\Uc}}t} |x_k(0)|.
\end{eqnarray*}
This expression can be further estimated as:
\begin{equation}
\begin{split}
|\phi_k(t,x_k,u)| \leq& \ e^{- \frac{1}{1+k \max\{\|u\|_{\Uc},\frac{2^k}{r}|x_k(0)|\}}t} \cdot \max\{|x_k(0)|,\frac{r}{2^k}\|u\|_{\Uc}\}.
\end{split}
\label{CounterEx_TMP_Eq1}
\end{equation}
For $|x_k(0)| \geq \frac{r}{2^k} \|u\|_{\Uc}$, we obtain 
\begin{eqnarray}
|\phi_k(t,x_k,u)| \leq e^{- \frac{1}{1+k \frac{2^k}{r}|x_k(0)|}t} |x_k(0)|.
\label{CounterEx_TMP_Eq2}
\end{eqnarray}
For $|x_k(0)| \leq \frac{r}{2^k} \|u\|_{\Uc}$, the inequality \eqref{CounterEx_TMP_Eq1} implies
\begin{align}
|\phi_k(t,x_k,u)| \leq e^{- \frac{1}{1+k \|u\|_{\Uc}}t} \frac{r}{2^k}\|u\|_{\Uc} \leq \frac{r}{2^k}\|u\|_{\Uc}.
\label{CounterEx_TMP_Eq3}
\end{align}
Overall, for any $x_k(0) \in \R$ and any $u \in \Uc$, we obtain from \eqref{CounterEx_TMP_Eq2} and \eqref{CounterEx_TMP_Eq3}:
\begin{align}
|\phi_k(t,x_k,u)| \leq e^{- \frac{1}{1+k \frac{2^k}{r}|x_k(0)|}t} |x_k(0)| +  \frac{r}{2^k}\|u\|_{\Uc}.
\label{CounterEx_TMP_Eq4}
\end{align}
Having this estimate for the state of the $k$-th mode, we proceed to the estimate for the whole state of $\Sigma$:
\begin{align}
\hspace{-4mm}
\|\phi(t,x,u)\|_X  &= \sum_{k=1}^{\infty} |\phi_k(t,x_k,u)| \nonumber\\
&\leq \sum_{k=1}^{\infty} e^{- \frac{1}{1+k \frac{2^k}{r}|x_k(0)|}t} |x_k(0)| +  \sum_{k=1}^{\infty} \frac{r}{2^k}\|u\|_{\Uc}\nonumber \\
&= \sum_{k=1}^{\infty} e^{- \frac{1}{1+k \frac{2^k}{r}|x_k(0)|}t} |x_k(0)| +  r \|u\|_{\Uc}.
\label{CounterEx_TMP_Eq5}
\end{align}
This estimate is true for all $t\geq 0$, all $x \in X$, all $u \in \Uc$, and for any $r >0$.

Now, we apply the trick used above in the proof that the system is AG. 
For any $x \in X$, for any $\eps >0$ there exists $N = N(x,\eps) \in \N$ so that $\sum_{k=N+1}^{\infty} |x_k(0)| <\frac{\eps}{2}$.
Using this fact, we continue estimates from \eqref{CounterEx_TMP_Eq5}:
\begin{align*}
\|\phi(t,x,u)\|_X  &\leq \sum_{k=1}^{N} e^{- \frac{1}{1+k \frac{2^k}{r}|x_k(0)|}t} |x_k(0)| \\
&\qquad\qquad + \sum_{k=N+1}^{\infty} e^{- \frac{1}{1+k \frac{2^k}{r}|x_k(0)|}t} |x_k(0)| +  r \|u\|_{\Uc}. \\
 &\leq \sum_{k=1}^{N} e^{- \frac{1}{1+k \frac{2^k}{r}|x_k(0)|}t} |x_k(0)| 
+ \frac{\eps}{2} +  r \|u\|_{\Uc}.
\end{align*}
For above $\eps$ and $x$, we can find sufficiently large time $\tau_a=\tau_a(\eps,x)$ (clearly, $\tau_a$ depends also on $N$, but $N$ itself depends only on $x$ and $\eps$), so that for all $u \in \Uc$ and all $t \geq \tau_a$
\[
\sum_{k=1}^{N} e^{- \frac{1}{1+k \frac{2^k}{r}|x_k(0)|}t} |x_k(0)| \leq  \frac{\eps}{2}.
\]

Overall, we obtain that for any $r >0$, for all $\eps >0$ and for any $x \in X$ there exists $\tau_a= \tau_a(t,x)$ such that 
for all $u \in \Uc$ and all $t \geq \tau_a$
\begin{eqnarray*}
\|\phi(t,x,u)\|_X  &\leq& \eps +  r \|u\|_{\Uc}.
\end{eqnarray*}
This shows that $\Sigma$ satisfies the strong AG property with the gain $\gamma (s) = r s$.

To finish the proof of (iii), we show that $\Sigma$ is not sAG for the gain $\gamma \equiv 0$.

Indeed, pick $\eps:=\tfrac{1}{2}$, $x= e_1= (1,0,\ldots,0,\ldots)^T$, and consider a certain constant input $u \equiv c$. The corresponding solution will have only one nonzero component - the first one. The norm of the state equals
\[
\|\phi(t,e_1,u)\|_X = |\phi_1(t,1,u)|= e^{-\frac{1}{1+c}t}. 
\]
If $\Sigma$ were sAG with the zero gain, then it would exist a time $\tau_a$, which does not depend on $u$, so that 
for all $c$ and all $t \geq \tau_a$ it holds that $e^{-\frac{1}{1+c}t} \leq \frac{1}{2}$.
But this is false since $e^{-\frac{1}{1+c}\tau_a}$ monotonically increases to 1 as long as $c \to \infty$.
Thus, $\Sigma$ is not sAG with zero gain.

\textbf{(iv):} We argue by contradiction, that $\Sigma$ is not LISS.
Assume that $\Sigma$ is LISS and hence there exist $r>0$, $\beta \in \KL$ and $\gamma \in \Kinf$ so that 
the inequality \eqref{iss_sum} holds for all $x \in X$: $\|x\|_X \leq r$ and for all $u \in \Uc$: $\|u\|_{\Uc} \leq r$.

Pick a constant input $u(t) = \eps$ for all $t\geq 0$, where $\eps>0$ is chosen so that 
\[
\max\{\eps,3\gamma(\eps)\} \leq r.
\]
This is always possible since $\gamma \in \Kinf$.
Denote $c:=\gamma(\eps)$. 

The LISS property of $\Sigma$ implies that for all $x \in X$ with the norm $\| x \|_{X} = 3c \leq r$ and for all $t \geq 0$ it holds that 
\begin {equation}
\label{CounterEx_ISS_ineq2}
\| \phi(t,x,u) \|_{X} \leq \beta(3c,t) + c.
\end{equation}
Since $\beta \in \KL$, there exists $t^*$ (which depends on $c$, but does not depend on $x$) so that $\beta(3c,t^*) = c$\footnote{Note that any $\beta \in \KL$ satisfying the LISS estimate automatically satisfies $\beta(r,0) \geq r$ for all $r>0$ (consider $t=0$ and $u\equiv 0$ in the LISS estimate). } and thus for all $x \in X$ with $\| x \|_{X} = 3c$ it must hold
\begin {equation}
\label{CounterEx_ISS_ineq3}
\| \phi(t^*,x,u) \|_{X} \leq 2c.
\end{equation}
Now consider the initial states of the form $x_{<k>}~=~3c~\cdot~e_k$, where $e_k$ is the $k$-th standard basis vector of $X=l_1$. Certainly, $\|x_{<k>}\|_X= 3c$ for all $k \in \N$.
Then 
\begin{eqnarray*}
\|\phi(t^*,x_{<k>},u)\|_X = |\phi_k(t^*,3c,u)| = e^{- \frac{1}{1+k\eps}t^*} 3c.
\end{eqnarray*}
Now, for any $t^{*}>0$ there exists $k$ so that $e^{- \frac{1}{1+k\eps}t^*} 3c > 2c$,
which contradicts to \eqref{CounterEx_ISS_ineq3}.
This shows that $\Sigma$ is not LISS. 
\end{example}

\begin{remark}
It is possible to check directly that the second assumption of Theorem~\ref{Characterization_LISS} does not hold.
\end{remark}

In Example~\ref{ex:0-UGAS_sAG_AG_zero_gain_not_LISS}, we have shown that a system which is 0-UGAS, sAG, AG with zero gain and UGS with zero gain, does not have to be LISS.
Next, we modify this example to show that \emph{if the system in addition to the above list of properties is LISS, this still does not guarantee ISS}.

\begin{example}[0-UGAS\,$\wedge$\,sAG\,$\wedge$\ AG with zero gain\,$\wedge$\ UGS with zero gain\,$\wedge$\ LISS\ \ $\not\Rightarrow$\ \ ISS]
\label{ex:0-UGAS_sAG_AG_zero_gain_LISS_not_ISS}
Consider a system $\Sigma$ with the state space $X:=l_1$ and input space $\Uc:=PC_b(\R_+,\R)$. Let also the dynamics of the $k$-th mode of $\Sigma$ be given by 
\begin{eqnarray}
\dot{x}_k(t) = -\frac{1}{1+|u(t)|^k}x_k(t).
\label{eq:CounterEx_AG_UGAS_LISS_no_ISS}
\end{eqnarray}
We continue to use the notation $\phi_k(t,x_k,u)$ for the state of the $k$-th mode \eqref{eq:CounterEx_AG_UGAS_LISS_no_ISS}.

As in Example~\ref{ex:0-UGAS_sAG_AG_zero_gain_not_LISS}, one can prove that this system satisfies Assumption~\ref{Assumption1}, is 0-UGAS, UGS with zero gain, AG with zero gain, and sAG with a nonzero gain (we skip this proof since it is completely analogous). Moreover, for $u$: $\|u\|_{\Uc} \leq 1$, for all $x \in X$ and all $t\geq 0$, it holds that
\begin{eqnarray}
\|\phi(t,x,u)\|_X \leq e^{-\frac{1}{2}t}\|x\|_X.
\label{eq:CounterEx_AG_UGAS_LISS_no_ISS_state_estimate}
\end{eqnarray}
Thus, $\Sigma$ is LISS with zero gain and with $r=1$.

The proof that $\Sigma$ is not ISS goes along the lines of Example~\ref{ex:0-UGAS_sAG_AG_zero_gain_not_LISS}, with the change that the norm of the considered inputs should be larger than 1. 
\end{example}

\section{Converse ISS Lyapunov theorems for semilinear systems}
\label{sec:input-state-stab}

In this section, we prove a converse ISS Lyapunov theorem for systems of the form \eqref{InfiniteDim}.

We follow the method developed in \cite{SoW95} for systems described by ODEs.
Consider the problem of global stabilization of
\eqref{InfiniteDim} by means of feedback laws that are subject to
multiplicative disturbances with a magnitude bounded by $1$. Let $\varphi:X \to \R_+$ be Lipschitz continuous on bounded sets and consider inputs
\begin{eqnarray}
u(t):=d(t)\varphi(x(t)), \quad t\geq 0,
\label{eq:Multiplicative_Feedbacks}
\end{eqnarray}
where $d \in \Dc$, where 
\begin{eqnarray}
\Dc := \{d:\R_+ \to D, \text{ piecewise continuous}\},\qquad D :=\{d \in U: \|d\|_U \leq 1\}.
\label{eq:Dcal-set}
\end{eqnarray}
Applying this feedback law to \eqref{InfiniteDim}, we obtain the closed-loop system 
\begin{eqnarray}
\label{eq:Modified_InfDimSys_With_Disturbances}
\dot{x}(t)&=&   Ax(t)+f(x(t),d(t)\varphi(x(t))) \nonumber\\
          &=:&  Ax(t) + g(x(t),d(t)).
\end{eqnarray}

\subsection{Basic properties of the closed-loop system}

In general, Assumption~\ref{Assumption1} does not guarantee Lipschitz continuity of $g$. To
ensure this, we assume that $f$ is Lipschitz continuous also with respect to $u$.
\begin{definition}
\label{def:Bi-Lipschitz-continuity} 
We call $f:X \tm U \to X$ Lipschitz continuous on bounded balls of $X \tm U$, if 
for each $R>0$ there is $L>0$ such that 
\begin{align}
x_1,x_2 \in \clo{B_R},\quad  &u_1,u_2 \in \clo{B_{R,U}} \nonumber\\
&\qrq \|f(x_1,u_1) - f(x_2,u_2)\|_X \leq L(\|x_1-x_2\|_X + \|u_1-u_2\|_U).
\label{eq:Bilinear-Lipschitz-continuity}
\end{align}
\end{definition}

It is easy to see that if $f$ is Lipschitz continuous on bounded balls of $X \tm U$, then Assumption~\ref{Assumption1} holds, and in particular, \eqref{InfiniteDim} is well-posed.

The next lemma shows that $g$ in \eqref{eq:Modified_InfDimSys_With_Disturbances} is Lipschitz continuous.

\begin{lemma}
\label{lem:Regularity_of_g}
Let $f$ be Lipschitz continuous on bounded balls of $X \tm U$.
Then $g$ (defined by \eqref{eq:Modified_InfDimSys_With_Disturbances}) is Lipschitz continuous on bounded balls of $X$, uniformly with respect to the second argument, i.e., 
$\forall C>0 \; \exists L_g(C)>0$, such that $\forall x,y \in \clo{B_C}$ and  $\forall d \in D$, it holds that
\begin{eqnarray*}
\|g(x,d)-g(y,d)\|_X \leq L_g(C) \|x-y\|_X.
\end{eqnarray*}
\end{lemma}

\begin{proof}
Pick arbitrary $C>0$, $x,y \in \clo{B_C}$, and $d \in D$. Then
\begin{align*}
\|g(x,d)-g(y,d)\|_X &= \|f(x,d\varphi(x)) - f(y,d\varphi(y))\|_X \\
&= \|f(x,d\varphi(x)) - f(y,d\varphi(x)) + f(y,d\varphi(x)) - f(y,d\varphi(y))\|_X \\
&\leq \|f(x,d\varphi(x)) - f(y,d\varphi(x))\|_X 
 + \|f(y,d\varphi(x)) - f(y,d\varphi(y))\|_X.
\end{align*}
Since $\varphi$ is Lipschitz continuous, it is bounded on $\clo{B_C}$
by a bound $R$. Without loss of generality, let $R>C$. As $f$ is Lipschitz continuous on bounded balls of $X \tm U$, and as
$\|d\|_{U}\leq 1$, we can upper bound the first summand by
$L_f^1(R)\|x-y\|_X$ and the second by $L_f^2(R) |\varphi(x)-\varphi(y)|$. The
claim now follows from the local Lipschitz continuity of $\varphi$.
\end{proof}

In particular, Lemma~\ref{lem:Regularity_of_g} shows that the system \eqref{eq:Modified_InfDimSys_With_Disturbances} is well-posed. Denote the maximal mild solution of \eqref{eq:Modified_InfDimSys_With_Disturbances} at time $t$, starting at
$x \in X$ and with disturbance $d \in \Dc$ by $\phi_{\varphi}(t,x,d)$.
On its interval of existence, $\phi_{\varphi}(t,x,d)$ coincides with
the solution of \eqref{InfiniteDim} for the input $u(t)=d(t)\varphi(\phi_{\varphi}(t,x_0,d))$.
\footnote{Forward completeness of \eqref{InfiniteDim} does not imply forward completeness of \eqref{eq:Modified_InfDimSys_With_Disturbances}. For  example, consider $\dot{x} = -x + u$, $u(t) = d\cdot x^2(t)$ for $d >0$.}

\begin{remark}
Lipschitz continuous feedbacks need not lead to Lipschitz continuous $g$ if $f$ is not Lipschitz with respect to inputs. Consider, e.g.,
$\dot{x}(t) = (u(t))^{1/3}$ and $u(t) := x(t)$.
\end{remark}

Next, we enlist several properties of the closed-loop system \eqref{eq:Modified_InfDimSys_With_Disturbances} that will be important for us in the sequel.
\begin{definition}
\label{Def_RFC}
\index{RFC}
\index{forward completeness!robust}
The system \eqref{eq:Modified_InfDimSys_With_Disturbances} is called \emph{robustly forward complete (RFC)} 
if for any $C>0$ and any $\tau>0$ it holds that 
\[
\sup_{\|x\|_X\leq C,\ d\in \Dc,\ t \in [0,\tau]}\|\phi_{\varphi}(t,x,d)\|_X < \infty.
\]
\end{definition}

\begin{definition}
\label{def:Lipschitz-continuity-of-the-flow-modified-sys}
Assume that \eqref{eq:Modified_InfDimSys_With_Disturbances}  is forward complete.
We say that the \emph{flow of \eqref{eq:Modified_InfDimSys_With_Disturbances} is Lipschitz continuous on compact intervals uniformly in $\Dc$}, if 
for any $\tau>0$ and any $R>0$ there exists $L>0$ so that for any $x,y \in
\clo{B_R}$, for all $t \in [0,\tau]$ and for all $d \in \Dc$ it holds that 
\begin{eqnarray}
\|\phi_{\varphi}(t,x,d) - \phi_{\varphi}(t,y,d) \|_X \leq L \|x-y\|_X.
\label{eq:Flow_is_Lipschitz-modified-system}
\end{eqnarray}    
\end{definition}

\begin{definition}
\label{def:UGAS}
\index{UGAS}
\index{stability!uniform global asymptotic}
The system \eqref{eq:Modified_InfDimSys_With_Disturbances} is called \emph{uniformly globally asymptotically
stable (UGAS)} if \eqref{eq:Modified_InfDimSys_With_Disturbances} is forward complete and there exists certain $\beta \in \KL$ such that 
\begin{equation}
\label{UGAS_wrt_D_estimate}
 d\in \Dc,\quad x \in X,\quad t \geq 0  \qrq  \|\phi_{\varphi}(t,x,d)\|_X \leq \beta(\|x\|_X,t).
\end{equation}
\end{definition}

UGAS can be characterized with the help of uniform global attractivity.
\begin{definition}
\label{def:UniformGlobalAttractivity}
\index{UGATT}
\index{attractivity!uniform global}
The system \eqref{eq:Modified_InfDimSys_With_Disturbances} is called \emph{uniformly globally attractive (UGATT)}, if \eqref{eq:Modified_InfDimSys_With_Disturbances} is forward complete and for any $r,\eps >0$ there exists $\tau=\tau(r,\eps)$ such that for all $d \in \Dc$ it holds that 
\begin{equation}
\|x\|_X \leq r, \quad t \geq \tau(r,\eps) \qrq \|\phi_{\varphi}(t,x,d)\|_X \leq \eps.
\label{eq:UAG_with_zero_gain}
\end{equation}
\end{definition}

\begin{definition}
\label{def:UniformGlobalStability}
\index{stability!global Lyapunov}
The system \eqref{eq:Modified_InfDimSys_With_Disturbances} is called \emph{globally Lyapunov stable}, if \eqref{eq:Modified_InfDimSys_With_Disturbances} is forward complete and there exists $\sigma \in \Kinf$ so that
\begin{equation}
d\in \Dc,\quad x \in X,\quad t \geq 0 \qrq \|\phi_{\varphi}(t,x,d)\|_X \leq \sigma(\|x\|_X).
\label{eq:GS}
\end{equation}
\end{definition}

The following characterization of UGAS is analogous to Lemma~\ref{lem:UAG_implies_ISS}, see also \cite[Theorem 2.2]{KaJ11b}.
\begin{proposition}
\label{prop:UGAS_Characterization}
The system \eqref{eq:Modified_InfDimSys_With_Disturbances} is UGAS if and only if \eqref{eq:Modified_InfDimSys_With_Disturbances} is  UGATT and globally Lyapunov stable.
\end{proposition}

Coercive Lyapunov functions for the UGAS property are defined as follows:
\begin{definition}
\label{def:UGAS_LF_With_Disturbances}
A continuous function $V:X \to \R_+$ is called a \emph{Lyapunov function} for \eqref{eq:Modified_InfDimSys_With_Disturbances},  if there exist
$\psi_1,\psi_2 \in \Kinf$ and $\alpha \in \Kinf$
such that 
\begin{equation}
\label{LyapFunk_1Eig_UGAS}
\psi_1(\|x\|_X) \leq V(x) \leq \psi_2(\|x\|_X) \quad \forall x \in X
\end{equation}
holds 
and the Dini derivative of $V$ along the trajectories of system
\eqref{eq:Modified_InfDimSys_With_Disturbances} satisfies for all $x \in
X$, and all $d \in \Dc$ that
\begin{equation}
\label{DissipationIneq_UGAS_With_Disturbances}
\dot{V}_d(x) \leq -\alpha(\|x\|_X).
\end{equation}
\end{definition}

The following converse Lyapunov theorem will be crucial for our developments \cite[Section 3.4]{KaJ11b}:
\begin{theorem}
\label{LipschitzConverseLyapunovTheorem-1}
Let \eqref{eq:Modified_InfDimSys_With_Disturbances} be UGAS and let its
flow be Lip\-schitz continuous on compact intervals uniformly in $\Dc$. Then
\eqref{eq:Modified_InfDimSys_With_Disturbances} admits a Lyapunov
function, which is Lipschitz continuous on bounded balls.
\end{theorem}

We will need the following property, which formalizes the robustness of \eqref{InfiniteDim} with respect to the feedback \eqref{eq:Multiplicative_Feedbacks}.
\begin{definition}
\label{def:WURS}
\index{WURS}	
\index{stability!weak uniform robust asymptotic}
The system \eqref{InfiniteDim} is called \emph{weakly uniformly robustly asymptotically
stable (WURS)}, if there exists $\varphi:X \to \R_+$, which is Lipschitz continuous on bounded sets, and $\psi \in \Kinf$ such that 
\[
\varphi(x) \geq \psi(\|x\|_X) \quad \forall x\in X,
\]
and \eqref{eq:Modified_InfDimSys_With_Disturbances} is uniformly globally asymptotically stable with that $\varphi$.
\end{definition}

The next proposition shows how the WURS property of system \eqref{InfiniteDim} reflects the regularity of the solutions of \eqref{eq:Modified_InfDimSys_With_Disturbances}.

\begin{proposition}
\label{prop:Flow_of_g}
Consider a forward complete system \eqref{InfiniteDim}.
Assume that
\begin{enumerate}
    \item[(i)] $f$ is Lipschitz continuous on bounded balls of $X \tm U$.
    \item[(ii)] \eqref{InfiniteDim} is WURS.
\end{enumerate}
Then for any $\varphi$ satisfying the conditions of
Definition~\ref{def:WURS}, the closed-loop system
\eqref{eq:Modified_InfDimSys_With_Disturbances} has a flow that is Lipschitz continuous on compact intervals uniformly in $\Dc$.
\end{proposition}

\begin{proof}
Since \eqref{InfiniteDim} is WURS and $\varphi$ is a stabilizing feedback
as required in Definition~\ref{def:WURS}, system
\eqref{eq:Modified_InfDimSys_With_Disturbances} is forward complete and
UGAS. Let $\beta\in{\cal KL}$ be a bound as in \eqref{UGAS_wrt_D_estimate}. Then \eqref{eq:Modified_InfDimSys_With_Disturbances} 
is RFC, as for any $C>0$ and any $\tau>0$ 
\[
\sup_{\|x\|_X\leq C,\ d\in \Dc,\ t \in [0,\tau]}\|\phi_{\varphi}(t,x,d)\|_X \leq \beta(C,0) < \infty.
\]
Assumption (i) together with Lemma~\ref{lem:Regularity_of_g} imply that
$g$ is Lipschitz continuous on bounded sets
uniformly in the second argument. 
Now, Assumptions~\ref{ass:Admissibility}, \ref{ass:Integrability} trivially hold for the system \eqref{eq:Modified_InfDimSys_With_Disturbances}. As $f$ is Lipschitz continuous in both arguments, and $\varphi$ is Lipschitz continuous as well, one can see that $g$ is continuous in both variables. Furthermore, the system \eqref{eq:Modified_InfDimSys_With_Disturbances} is UGAS, and thus it holds that $g(0,d)=0$ for all $d$. This shows that Assumption~\ref{Assumption1} holds for \eqref{eq:Modified_InfDimSys_With_Disturbances} as well. 
Now, Theorem~\ref{thm:Lipschitz-continuity-of-flow} applied to the system \eqref{eq:Modified_InfDimSys_With_Disturbances}
implies the claim.
\end{proof}

\subsection{Main result}

\begin{figure}[tbh]
\centering
\begin{tikzpicture}[>=implies,thick]
\node (ISS) at (1,5) {(1) is ISS};
\node (WURS) at   (3.5,3.5) {(1) is WURS};
\node (ISSLF) at (6,5) {$\exists$ ISS-LF for (1)};

\node (Lem4) at (1.2,4) {\scriptsize Lemma~\ref{lem:ISS_implies_WURS}};
\node (Lem_and_Prop) at (7.5,4) {\scriptsize Lemma~\ref{lem:WURS_implies_ISS_LF}, based on Theorem~\ref{LipschitzConverseLyapunovTheorem-1}};
\node (Prop1) at (3.4,5.3) {\scriptsize Theorem~\ref{LyapunovTheorem}};

\draw[thick,double equal sign distance,->] (ISS) to (WURS);
\draw[thick,double equal sign distance,->] (WURS) to (ISSLF);
\draw[thick,double equal sign distance,->] (ISSLF) to (ISS);

\end{tikzpicture}
\caption{Converse ISS Lyapunov theorem (Theorem~\ref{ISS_Converse_Lyapunov_Theorem}).}
\label{ISS_CLT}
\end{figure}
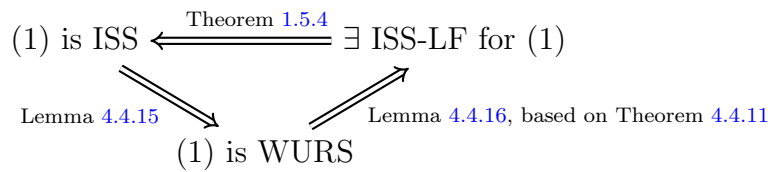

In this section, we  prove that for system \eqref{InfiniteDim} with
$f$ being Lipschitz continuous on bounded balls of $X \tm U$ the notions depicted in Figure~\ref{ISS_CLT} are equivalent.

\begin{remark}
Lipschitz continuity of $f$ on bounded balls of $X \tm U$ was needed already in the ODE case \cite{SoW95} to ensure the Lipschitz continuity of the flow of the auxiliary system \eqref{eq:Modified_InfDimSys_With_Disturbances}.
In the infinite-dimensional case, without Lipschitz continuity of $f$ in $u$, the nonlinearity $g$ in \eqref{eq:Modified_InfDimSys_With_Disturbances} is merely continuous, which does not guarantee even the existence of solutions of \eqref{eq:Modified_InfDimSys_With_Disturbances} in general, as Peano theorem does not hold for infinite-dimensional systems \cite{God75, HaJ10}.
\end{remark}

First, we show in Lemma~\ref{lem:ISS_implies_WURS} that ISS implies WURS. 
Next, we apply Theorem~\ref{LipschitzConverseLyapunovTheorem-1} to prove that WURS of \eqref{InfiniteDim} implies the existence of a Lipschitz continuous coercive ISS Lyapunov function for \eqref{InfiniteDim}.
Finally, the direct Lyapunov theorem
(Theorem~\ref{LyapunovTheorem}) completes the proof.

\begin{lemma}
\label{lem:ISS_implies_WURS}
If \eqref{InfiniteDim} is ISS, then \eqref{InfiniteDim} is WURS.
\end{lemma}

\begin{proof}
The proof goes along the lines of \cite[Lemma 2.12]{SoW95}. 
Let \eqref{InfiniteDim} be ISS. To prove that \eqref{InfiniteDim} is WURS, we are going to use Proposition~\ref{prop:UGAS_Characterization}. 

Since \eqref{InfiniteDim} is ISS, there exist $\beta \in \KL$ and $\gamma \in \Kinf$ so that \eqref{iss_sum} holds for any $t \geq 0$, $x \in X$, $u \in \Uc$.
Define $\alpha(r):=\beta(r,0)$, for $r \in \R_+$. Substituting $u \equiv 0$ and $t=0$ into \eqref{iss_sum} we see that $\alpha(r) \geq r$ for all $r \in \R_+$.

Pick any $\sigma \in \Kinf$ so that $\sigma (r) \leq \gamma^{-1}\big(\frac{1}{4} \alpha^{-1}(\frac{2}{3}r)\big)$ for all $r\geq 0$.
We may choose Lipschitz continuous on bounded sets maps $\varphi:X\to \R_+$ and $\psi \in \Kinf$ such
that  $\psi(\|x\|_X) \leq \varphi(x) \leq \sigma(\|x\|_X)$
(just pick a Lipschitz continuous on bounded sets $\psi \in \Kinf$ and set
$\varphi(x):=\psi(\|x\|_X)$ for all $x \in X$, which guarantees that
$\varphi$ is Lipschitz continuous on bounded sets).

We are going to show that for all $x \in X$, all $t \geq 0$ and all $d \in \Dc$ it holds that
\begin{eqnarray}
\gamma\Big(\big\|d(t) \varphi(\phi_{\varphi}(t,x,d))\big\|_U\Big) \leq \frac{\|x\|_X}{2}. 
\label{eq:ISS_implies_WURS_Estimate}
\end{eqnarray}
First, we show that \eqref{eq:ISS_implies_WURS_Estimate} holds for all
times $t\geq 0$ small enough. 
Since $\alpha^{-1}(r) \leq r$ for all $r >0$, we have
\begin{align*}
\gamma\Big(\big\|d(t) \varphi(\phi_{\varphi}(t,x,d))\big\|_U\Big)
& \leq \gamma\big(\sigma(\|\phi_{\varphi}(t,x,d)\|_X)\big) \\
&\leq \frac{1}{4} \alpha^{-1}\Big(\frac{2}{3}\|\phi_{\varphi}(t,x,d)\|_X\Big) 
\leq \frac{1}{6} \|\phi_{\varphi}(t,x,d)\|_X.
\end{align*}
For any $d \in \Dc$ and any $x \in X$ the latter expression can be made smaller than 
$\frac{1}{2} \|x\|_X$ by choosing $t$ small enough, since $\phi_{\varphi}$
is continuous in $t$.

Now pick any $d \in \Dc$, $x \in X$ and define $t^* = t^*(x,d)$ by
\begin{equation*}
t^*:= \inf \left\{t \geq 0: \gamma\Big(\|d(t)\|_U \big|\varphi(\phi_{\varphi}(t,x,d))\big|\Big) > \frac{\|x\|_X}{2}\right\}.
\end{equation*}
By the first step, we know $t^*>0$.
Assume that $t^*< \infty$ (otherwise, our claim is true).
Then \eqref{eq:ISS_implies_WURS_Estimate} holds for all $t \in [0,t^*)$.
Thus, for all $t \in [0,t^*)$ it holds that 
\begin{align*}
\|\phi_{\varphi}(t,x,d)\|_X &\leq \beta(\|x\|_X,t) + \frac{\|x\|_X}{2} \\  
&\leq \beta(\|x\|_X,0) + \frac{\alpha(\|x\|_X)}{2}  
 =  \frac{3}{2} \alpha(\|x\|_X).
\end{align*}
Using this estimate, we find that
\begin{align*}
\gamma\Big(\|d(t^*)\|_U \big|\varphi(\phi_{\varphi}(t^*,x,d))\big|\Big) 
&\leq 
\frac{1}{4} \alpha^{-1}\Big(\frac{2}{3}\|\phi_{\varphi}(t^*,x,d)\|_X\Big)  \\
&\leq
\frac{1}{4} \alpha^{-1}(\|x\|_X)  \leq \frac{1}{4} \|x\|_X.
\end{align*}
But this contradicts the definition of $t^*$.
Thus, $t^* = +\infty$.

Now we see that for any $x \in X$, any $d \in \Dc$ and all $t \geq 0$, we have 
\begin{equation}
\label{eq:Helps_for_RFC}
\|\phi_{\varphi}(t,x,d)\|_X \leq \beta(\|x\|_X,t) + \frac{\|x\|_X}{2},
\end{equation}
which shows global Lyapunov stability of \eqref{eq:Modified_InfDimSys_With_Disturbances}.

Since $\beta \in \KL$, there exists a $t_1=t_1(\|x\|_X)$ so that $\beta(\|x\|_X,t_1) \leq \frac{\|x\|_X}{4}$ and consequently
\begin{equation*}
 d\in \Dc,\; x \in X,\; t \geq t_1 \quad \Rightarrow \quad \|\phi_{\varphi}(t,x,d)\|_X \leq \frac{3}{4} \|x\|_X.
\end{equation*}

By induction, we obtain that there exists a strictly increasing sequence of times $(t_k)_{k}$, which depends on the norm of $\|x\|_X$ but is independent of $x$ and $d$ so that
\begin{equation*}
\|\phi_{\varphi}(t,x,d)\|_X \leq \Big(\frac{3}{4} \Big)^{k} \|x\|_X,
\end{equation*}
for all $x \in X$, any $d \in \Dc$ and all $t \geq t_k$.

This means that for all $\eps>0$ and for all $\delta >0$ there exist a time $\tau=\tau(\delta)$ so that 
for all $x \in X$ with $\|x\|_X \leq \delta$, for all $d \in \Dc$ and for all $t \geq \tau$ we have
\[
\|\phi_{\varphi}(t,x,d)\|_X \leq \eps.
\]
This shows uniform global attractivity of \eqref{eq:Modified_InfDimSys_With_Disturbances}.

Now Proposition~\ref{prop:UGAS_Characterization} shows that \eqref{eq:Modified_InfDimSys_With_Disturbances} is UGAS and thus 
\eqref{InfiniteDim} is WURS.
\end{proof}

\begin{lemma}
\label{lem:WURS_implies_ISS_LF}
If \eqref{InfiniteDim} is WURS, and $f$ is Lipschitz continuous on bounded balls of $X \tm U$, then there exists an
ISS Lyapunov function for \eqref{InfiniteDim}, which is Lipschitz continuous on bounded sets. 
\end{lemma}

\begin{proof}
Let \eqref{InfiniteDim} be WURS. This means that 
\eqref{eq:Modified_InfDimSys_With_Disturbances} is UGAS over $\Dc$ for
suitable $\varphi, \psi$ chosen in accordance with Definition~\ref{def:WURS}.
Proposition~\ref{prop:Flow_of_g} and Theorem~\ref{LipschitzConverseLyapunovTheorem-1} imply that
there exists a Lipschitz continuous on bounded sets Lyapunov function $V:X\to\R_+$, satisfying \eqref{LyapFunk_1Eig_UGAS}
for certain $\psi_1,\psi_2 \in \Kinf$ and whose Lie derivative along the solutions of \eqref{eq:Modified_InfDimSys_With_Disturbances} satisfies for all $x \in X$ and for all $d \in \Dc$ the estimate
\begin{eqnarray}
\hspace{-8mm}\dot{V}_{\varphi,d}(x) := \mathop{\overline{\lim}} \limits_{t \rightarrow +0} {\frac{1}{t}\big(V(\phi_\varphi(t,x,d))-V(x)\big) }\leq -\alpha(V(x)).
\label{eq:WURS_LF_Ineq}
\end{eqnarray}
Since \eqref{eq:Modified_InfDimSys_With_Disturbances} is UGAS, $\phi_\varphi(\cdot,x,d)$ is well-defined on $\R_+$ for any $d\in\Dc$ and any $x\in X$.
Now pick any $x\in X$. 
We claim that for all $u\in\Uc$ so that $\|u\|_\Uc \leq \varphi(x)$ it holds that 
\begin{eqnarray}
\dot{V}_u(x)\leq -\alpha(V(x)).
\label{eq:WURS_LF_Ineq_2}
\end{eqnarray}
From \eqref{eq:WURS_LF_Ineq_2} it follows that $V$ is an ISS Lyapunov function for \eqref{InfiniteDim} in implication form with Lyapunov gain $\chi:=\psi^{-1}$.

The following technical argument shows \eqref{eq:WURS_LF_Ineq_2}. Let
$x\in X$ and $u\in\Uc$ be so that $\|u\|_\Uc \leq \varphi(x)$.
If $x=0$, then $u\equiv 0$, and $\dot{V}_0(0)=\dot{V}_{\varphi,0}(0)\leq -\alpha(V(0)) = 0$, and \eqref{eq:WURS_LF_Ineq_2} holds.

Now let $x\neq 0$. Then there is $d \in\Dc$ so that $u(t) = d(t)\varphi(x)$
for $t\geq0$. Define further $u_2(t) := d(t)\varphi(\phi_{\varphi}(t,x,d))$,
$x(t):=\phi(t,x,u)$, $x_2(t):=\phi(t,x,u_2)$, for $t\in\R_+$.

Let us compute the Lie derivative of $V$:
{\allowdisplaybreaks
\begin{align*}
\dot{V}_u(x) &= \mathop{\overline{\lim}} \limits_{t \rightarrow +0} {\frac{1}{t}\big(V(x(t))-V(x)\big) }\\
 &= \mathop{\overline{\lim}} \limits_{t \rightarrow +0} {\frac{1}{t}\Big(V(x(t))-V(x_2(t)) + V(x_2(t)) -V(x)\Big) }\\
&\leq  \mathop{\overline{\lim}} \limits_{t \rightarrow +0}
{\frac{V(x(t))-V(x_2(t))} {t}}+ \mathop{\overline{\lim}} \limits_{t \rightarrow +0} {\frac{V(x_2(t)) -V(x)}{t} }\\
&\leq  \mathop{\overline{\lim}} \limits_{t \rightarrow +0}
{\frac{|V(x(t))-V(x_2(t))|}{t}} -\alpha(V(x)),
\end{align*}
}
where we have used \eqref{eq:WURS_LF_Ineq} and that for $t\in\R_+$ it holds that
\[
x_2(t) = \phi\big(t,x,d(\cdot)\varphi(\phi_{\varphi}(\cdot,x,d))\big) = \phi_\varphi(t,x,d).
\]
As $V$ is Lipschitz continuous on bounded sets, there is a constant $L>0$ so that 
\begin{eqnarray}
\dot{V}_u(x) \leq \mathop{\overline{\lim}} \limits_{t \rightarrow +0} {\frac{L}{t}\big\|x(t) - x_2(t)\big\|_X} -\alpha(V(x)).
\label{eq:Computation_Lie_derivative_u}
\end{eqnarray}
As $f$ is Lipschitz continuous on bounded balls of $X \tm U$, and since there is $M>0$ so that $\|T(r)\|\leq M$ for $r\in[0,1]$, we obtain for $t\in[0,1]$:
{\allowdisplaybreaks
\begin{align*}
\big\|x(t) - x_2(t)\big\|_X
&=  \Big\|\int_0^t T(t-s) \big(f(x(s),u(s)) - f(x_2(s),u_2(s)) \big) ds\Big\|_X \\%
&\leq  \int_0^t \big\|T(t-s)\big\| \big\|f(x(s),u(s)) - f(x_2(s),u_2(s)) \big\|_X ds \\%
&\leq  M \int_0^t  L_f \left(\|x(s) - x_2(s)\|_X +  \|u(s)-u_2(s))\|_U\right) ds\\
&\leq  M L_f \left(t \sup_{s\in[0,t]}\big\|u(s) {-} u_2(s) \big\|_U +  \int_0^t \big\|x(s) {-} x_2(s)\big\|_X ds\right),  
\end{align*}
} 
for a certain constant $L_f > 0$.
Gronwall's inequality \eqref{lem:Gronwall}  implies that
\begin{align*}
\big\|&x(t) - x_2(t)\big\|_X \leq  M L_f t \sup_{s\in[0,t]}\big\|u(s) - u_2(s) \big\|_U e^{ML_f t},\quad t\in[0,1].
\end{align*}
Substituting this into \eqref{eq:Computation_Lie_derivative_u}, we obtain that
\begin{eqnarray}
\hspace{-8mm} \dot{V}_u(x) &\leq& \mathop{\overline{\lim}} \limits_{t \rightarrow +0}\big( L M L_f \sup_{s\in[0,t]}\big\|u(s) - u_2(s) \big\|_U e^{ML_f t}\big) -\alpha(V(x)).
\label{eq:Computation_Lie_derivative_u_2}
\end{eqnarray}
Furthermore, as $\|d(s)\|_U\leq 1$ for all $s\in[0,t]$, we have
\begin{align*}
\sup_{s\in[0,t]} \big\|u(s) - u_2(s) \big\|_U &= \sup_{s\in[0,t]}\big\|d(s)\big(\varphi(x) - \varphi(\phi_{\varphi}(s,x,d))\big) \big\|_U \\
&\leq  \sup_{s\in[0,t]}|\varphi(x) - \varphi(\phi_{\varphi}(s,x,d)) |.
\end{align*}
Since $\phi_{\varphi}(s,x,d) \to x$, as $s \to 0$, and since $\varphi$ is Lipschitz continuous,
we obtain
\[
\sup_{s\in[0,t]} \big\|u(s) - u_2(s) \big\|_U \to 0,\quad t\to 0.
\]
Substituting this into \eqref{eq:Computation_Lie_derivative_u_2}, we see that \eqref{eq:WURS_LF_Ineq_2} holds 
for all $x \in X$ and all $u \in \Uc$ with $\|u\|_{\Uc} \leq \psi(\|x\|_X)$.
\end{proof}

We conclude our investigation with the following characterization of the ISS property:
\begin{theorem}
\label{ISS_Converse_Lyapunov_Theorem}
Let $f$ be Lipschitz continuous on bounded balls of $X \tm U$. Then the following statements are equivalent:
\begin{enumerate} 
    \item[(i)] \eqref{InfiniteDim} is ISS.
    \item[(ii)] \eqref{InfiniteDim} is WURS.
    \item[(iii)] There exists a Lipschitz continuous on bounded balls coercive ISS Lyapunov function for
          \eqref{InfiniteDim}.
\end{enumerate}
\end{theorem}

\begin{proof}
The claim follows from Theorem~\ref{LyapunovTheorem} and Lemmas~\ref{lem:ISS_implies_WURS} and \ref{lem:WURS_implies_ISS_LF}. 
\end{proof}

\section{Lyapunov characterization for ISS of linear systems with bounded input operators}
\label{sec:Lyapunov_Theorem_systems_bounded_operators}

Our converse ISS Lyapunov theorem (Theorem~\ref{ISS_Converse_Lyapunov_Theorem}) is based upon the converse UGAS Lyapunov theorem
(Theorem~\ref{LipschitzConverseLyapunovTheorem-1}), which is non-constructive in the sense that the Lyapunov function that is provided in the proof of this result is constructed based on the knowledge of the flow, which is usually unknown.

Next, we give exhaustive criteria for ISS of linear systems
\begin{equation}
\label{eq:Linear_System-ISS-LF}
\dot{x}(t)=Ax(t)+ Bu(t), \quad t > 0,
\end{equation}
where $A$ generates a $C_0$-semigroup $T$ over a Banach space $X$, and $B\in L(U,X)$.
In particular, we provide an explicit expression for the ISS Lyapunov function for the ISS system \eqref{eq:Linear_System-ISS-LF} 
based on the knowledge of the semigroup $T$ generated by $A$.

\begin{theorem}[ISS characterizations for linear systems with bounded input operators]
\label{thm:ISS-criterion-linear-systems-bounded-operators}
The following statements are equivalent for the system \eqref{eq:Linear_System-ISS-LF} with $B\in L(U,X)$:
\begin{enumerate}[label=(\roman*)]
    \item\label{itm:0-ISS-LinSysCrit} \eqref{eq:Linear_System-ISS-LF} is ISS w.r.t.\  $PC_b(\R_+,U)$.
    \item\label{itm:i-ISS-LinSysCrit} \eqref{eq:Linear_System-ISS-LF} is eISS w.r.t.\  $L^p(\R_+,U)$ for all $p\in[1,+\infty]$.
    \item\label{itm:ii-ISS-LinSysCrit} \eqref{eq:Linear_System-ISS-LF} is eISS w.r.t.\  $L^p(\R_+,U)$ for some $p\in[1,+\infty)$.
    \item\label{itm:iv-ISS-LinSysCrit} \eqref{eq:Linear_System-ISS-LF} is 0-UGAS.
    \item\label{itm:v-ISS-LinSysCrit} The semigroup $T$ is exponentially stable.
    \item\label{itm:vi-ISS-LinSysCrit} The function $V:X\to \R_+$ defined by
    \begin{eqnarray}
V(x)=\int_0^{\infty} \|T(t) x\|_X^2 dt
\label{eq:LF_LinSys_Banach-formulation}
\end{eqnarray}
is a non-coercive ISS Lyapunov function for \eqref{eq:Linear_System-ISS-LF} with $\Uc:=PC_b(\R_+,U)$.
    \item\label{itm:vii-ISS-LinSysCrit} The function $V^\gamma:X\to \R_+$ defined by
\begin{equation}
    \label{eq:eqnorm}
    V^\gamma (x) := \max_{s\geq 0} \| e^{\gamma s} T(s) x\|_X,
\end{equation}
is a uniformly globally Lipschitz continuous coercive ISS Lyapunov function for \eqref{eq:Linear_System-ISS-LF} with $\Uc:=PC_b(\R_+,U)$.
Here $\gamma\in (0,\lambda)$, and $\lambda>0$ is so that $\|T(t)\|\leq Me^{-\lambda t}$ for some $M>0$ and all $t\geq 0$.
    \item\label{itm:viii-ISS-LinSysCrit} There is a non-coercive ISS Lyapunov function for \eqref{eq:Linear_System-ISS-LF} with $\Uc:=PC_b(\R_+,U)$.
\end{enumerate}
\end{theorem}

\begin{proof}
\ref{itm:0-ISS-LinSysCrit} $\srs$ \ref{itm:iv-ISS-LinSysCrit}. Clear.

\ref{itm:i-ISS-LinSysCrit} $\srs$ \ref{itm:ii-ISS-LinSysCrit} $\srs$ \ref{itm:iv-ISS-LinSysCrit}. Clear.

\ref{itm:iv-ISS-LinSysCrit} $\srs$ \ref{itm:v-ISS-LinSysCrit}. Let \eqref{eq:Linear_System-ISS-LF} be 0-UGAS and pick $u\equiv 0$. 
Then there is $\beta\in\KL$ such that $\|T(t) x\|_X \leq \beta(1,t)$ for
all $t \geq 0$ and for all $x$ with  $\|x\|_X=1$. Since $\beta \in \KL$,
there exists a $t^*>0$ such that $\|T({t^*})\|<1$ and thus $T$ is an exponentially
stable semigroup.

\ref{itm:v-ISS-LinSysCrit} $\srs$ \ref{itm:i-ISS-LinSysCrit}. Take $M,\lambda>0$ such that $\|T(t)\| \leq Me^{-\lambda t}$.
Pick also any $x_0 \in X$ and $u\in\Uc$ and estimate the norm of the solution $x(t)=\phi(t,x_0,u)$ of \eqref{eq:Linear_System-ISS-LF}:
\begin{eqnarray*}
\|\phi(t,x_0,u)\|_X  &= & \Big\|T(t)x_0 + \int_0^t{T(t-r)B u(r)dr}\Big\|_X \\
						&\leq & \|T(t)\|\|x_0\|_X + \int_0^t{\|T(t-r)\|\|B\| \|u(r)\|_U dr}, \\
					&	\leq & M e^{-\lambda t} \|x_0\|_X + M \|B\| \int_0^t{e^{-\lambda (t-r)} \|u(r)\|_U dr}.
\end{eqnarray*}
Now, estimating $e^{-\lambda (t-r)}  \leq 1$, $r \leq t$, we obtain that \eqref{eq:Linear_System-ISS-LF} is eISS w.r.t.\  $L^1(\R_+,U)$. 
From here, we also easily obtain eISS w.r.t.\  $L^\infty(\R_+,U)$.

To prove the claim for $p\in(1,+\infty)$, pick any $q > 1$ so that $\frac{1}{p}+ \frac{1}{q} = 1$. 
We continue the above estimates, using H\"older's inequality:
\begin{align*}
\|\phi(t,x_0,u)\|_X  &\leq  M e^{-\lambda t} \|x_0\|_X  + M \|B\| \int_0^t{e^{-\frac{\lambda}{2} (t-r)} e^{-\frac{\lambda}{2} (t-r)} \|u(r)\|_U dr}  \\
&\leq  M e^{-\lambda t} \|x_0\|_X + M \|B\|  \Big(\int_0^t{e^{- \frac{q \lambda}{2} (t-r)}dr} \Big)^{\frac{1}{q}} \cdot \Big(\int_0^t{e^{-\frac{p \lambda}{2} (t-r)}\|u(r)\|^p_U dr} \Big)^{\frac{1}{p}} \\
 &\leq  M e^{-\lambda t} \|x_0\|_X + M \|B\|  \Big( \frac{2}{q\lambda} \Big)^{\frac{1}{q}}  \Big(\int_0^t{\|u(r)\|^p_U dr} \Big)^{\frac{1}{p}} \\
& =  M e^{-\lambda t} \|x_0\|_X + w \|u\|_{L^p(\R_+,U)},
\end{align*}
where $w:=M \|B\|  \Big( \frac{2}{q\lambda} \Big)^{\frac{1}{q}}$.
This means that \eqref{eq:Linear_System-ISS-LF} is eISS w.r.t.\  $L^p(\R_+,U)$ also for $p>1$.

\ref{itm:v-ISS-LinSysCrit} $\srs$ \ref{itm:vi-ISS-LinSysCrit}. Follows by Proposition~\ref{ConverseLyapunovTheorem_LinearSystems}.

\ref{itm:v-ISS-LinSysCrit} $\srs$ \ref{itm:vii-ISS-LinSysCrit}. Follows by Proposition~\ref{p:Lyapunovfunction2}.

\ref{itm:vi-ISS-LinSysCrit}, \ref{itm:vii-ISS-LinSysCrit} $\srs$ \ref{itm:viii-ISS-LinSysCrit}. Clear.

\ref{itm:viii-ISS-LinSysCrit} $\srs$ \ref{itm:0-ISS-LinSysCrit}. Follows by Theorem~\ref{t:ISSLyapunovtheorem} (CEP an BRS properties for a linear system \eqref{eq:Linear_System-ISS-LF} with $\Uc:=PC_b(\R_+,U)$ and $B\in L(U,X)$ are always fulfilled).
\end{proof}

\begin{remark}
\label{rem:Coercivity-of-quadratic-LFs-in-special-cases}
The Lyapunov function defined by \eqref{eq:LF_LinSys_Banach-formulation} is, in general, non-coercive. However, it is coercive for some special classes of systems.
\end{remark}

Theorem~\ref{thm:ISS-criterion-linear-systems-bounded-operators} shows that
the question about ISS of a control system \eqref{eq:Linear_System-ISS-LF} with $B\in L(U, X)$ can be reduced to the analysis of the exponential stability of a semigroup $T$, generated by $A$, which is a classical problem in semigroup theory.
On the other hand, for each input-to-state stable linear system with a bounded input operator, there is a quadratic non-coercive as well as a homogeneous coercive ISS Lyapunov function.
The Lyapunov function \eqref{eq:LF_LinSys_Banach-formulation} is widely used for linear systems without inputs, see, e.g., \cite[Theorem 5.1.3]{CuZ95}, and in a Hilbert space setting, it is obtained via the solution of an operator Lyapunov equation, see Theorem~\ref{prop:Converse_ISS_Lyapunov_theorem_lin_Systems_with_bounded_input_operators}. The contribution of Theorem~\ref{thm:ISS-criterion-linear-systems-bounded-operators}\ref{itm:vi-ISS-LinSysCrit} is that this function is also a non-coercive ISS Lyapunov function for \eqref{eq:Linear_System-ISS-LF} with $B\in L(U, X)$. The fact that such natural Lyapunov functions, as \eqref{eq:LF_LinSys_Banach-formulation} are non-coercive, is one of the reasons motivating the study of non-coercive Lyapunov functions.

\begin{remark}\label{REM:1}
For finite-dimensional linear systems \eqref{eq:Linear_System-ISS-LF}, 0-UGAS (and thus ISS) is equivalent to the \emph{strong stability of a semigroup $T$}, i.e., to the property that for all $x\in X$ it holds that $\phi(t,x,0)=T(t)x \to 0$ as $t\to\infty$.
For infinite-dimensional linear systems \eqref{eq:Linear_System-ISS-LF}, strong stability of a semigroup $T$ is much weaker than exponential stability, even for bounded $A$. Moreover, in this case, \eqref{eq:Linear_System-ISS-LF} may have unbounded solutions in the presence of a bounded and exponentially converging to zero input, see, e.g., \cite[p. 247]{MaP11}. See Section~\ref{sec:WeakISS} for the ISS-like nonlinear counterpart of a strong stability concept for systems with inputs.
\end{remark}

We continue with several results that we used in the proof of Theorem~\ref{thm:ISS-criterion-linear-systems-bounded-operators}, where we assume $\Uc:=PC_b(\R_+,U)$.

We start with a technical lemma; its proof is straightforward and is omitted.
\begin{lemma}
\label{AuxiliaryEquality}
Let $B \in L(U,X)$ and let $T$ be a $C_0$-semigroup. Then for any $u\in \Uc:=PC_b(\R_+,U)$ it holds that
\begin{eqnarray}
\lim\limits_{h \rightarrow +0} \frac{1}{h} \int_0^h{T(h-s) B u(s)ds} = B u(0).
\label{eq:Convolution_Bu0}
\end{eqnarray}
\end{lemma}

\begin{proposition}
\label{ConverseLyapunovTheorem_LinearSystems}
If $T$ is an exponentially stable semigroup, then $V:X \to \R_+$, defined as
\begin{eqnarray}
V(x)=\int_0^{\infty} \|T(t) x\|_X^2 dt, \quad x\in X,
\label{eq:LF_LinSys_Banach}
\end{eqnarray}
is a non-coercive ISS Lyapunov function for \eqref{eq:Linear_System-ISS-LF}, which is
Lipschitz continuous on bounded sets.
 Moreover, for all $ x \in X$, $ u \in \Uc:=PC_b(\R_+,U)$ and $ \eps>0$, it holds that 
\begin{equation}
\label{Final_Lyapunov_Inequality}
\dot{V}_u(x) \leq -\|x\|_X^2 + \frac{\eps M^2}{2\lambda} \|x\|_X^2 + \frac{M^2}{2\lambda \eps}  \|B\|^2 \|u(0)\|_U^2,
\end{equation}
where $M, \lambda>0$ are so that 
\begin{equation}
\|T(t) \| \leq M e^{-\lambda t}.
\label{eq:SemigroupEstimate}
\end{equation}
\end{proposition}

\begin{proof}
As $T$ is exponentially stable, there exist $M,\lambda>0$ such that \eqref{eq:SemigroupEstimate} holds.

Consider $V:X \to \R_+$ as defined in \eqref{eq:LF_LinSys_Banach}.
We have
\begin{eqnarray}
V(x) \leq \int_0^{\infty}\hspace{-1.5mm} \|T(t) \|^2\|x\|_X^2 dt \leq \frac{M^2}{2\lambda} \|x\|_X^2.
\label{eq:LF_Upper_estimate}
\end{eqnarray}

Let $V(x)=0$. Then $\|T(t) x\|_X \equiv 0$ a.e. on $[0,\infty)$. Strong continuity of $T$ implies that $x = 0$, and thus 
\eqref{LyapFunk_1Eig_nc_ISS} holds.

Next, we estimate the Dini derivative of $V$: 
{\allowdisplaybreaks
\begin{align*}
\dot{V}_u(x)& = \mathop{\overline{\lim}} \limits_{h \rightarrow +0} {\frac{1}{h}(V(\phi(h,x,u))-V(x)) } \\
    &=  \mathop{\overline{\lim}} \limits_{h \rightarrow +0} {\frac{1}{h}\Big(\int_0^{\infty} \|T(t) \phi(h,x,u)\|_X^2 dt - \int_0^{\infty} \|T(t) x\|_X^2 dt \Big) } \\
    &=  \mathop{\overline{\lim}} \limits_{h \rightarrow +0}
\frac{1}{h}\Big(\int_0^{\infty} \Big\|T(t) \Big(T(h)x + \int_0^h{T({h-s}) B
  u(s)ds}\Big) \Big\|_X^2 dt \\
    & \qquad -\int_0^{\infty} \|T(t) x\|_X^2 dt \Big)  \\
    &=  \mathop{\overline{\lim}} \limits_{h \rightarrow +0} \frac{1}{h}\Big(\int_0^{\infty} \Big\|T({t+h}) x 
		+ T(t) \int_0^h{T(h-s) B u(s)ds}\Big\|_X^2 dt \\
		    & \qquad - \int_0^{\infty} \|T(t) x\|_X^2 dt \Big) \\  
   &\leq  \mathop{\overline{\lim}} \limits_{h \rightarrow +0} \frac{1}{h}\Big(\int_0^{\infty} \Big(\Big\|T(t+h) x\Big\|_X \\
    &\qquad + \Big\|T(t) \int_0^h\hspace{-1.5mm}{T(h-s) B u(s)ds}\Big\|_X \Big)^2 dt - \int_0^{\infty} \hspace{-1.5mm}\|T(t) x\|_X^2 dt \Big)  \\
   &=  I_1 + I_2,
\end{align*}
}
where $I_1:= \mathop{\overline{\lim}} \limits_{h \rightarrow +0} \frac{1}{h}\Big(\int_0^{\infty} \|T(t+h) x\|_X^2 dt - \int_0^{\infty} \|T(t) x\|_X^2 dt \Big)$,
\begin{align*}
\text{and}\ \ I_2  &:= \mathop{\overline{\lim}} \limits_{h \rightarrow +0} \frac{1}{h}  \int_0^{\infty} \Big( 2\big\|T(t+h) x\big\|_X \Big\|T(t) \int_0^h\hspace{-1.5mm}{T(h-s) B u(s)ds}\Big\|_X \\
 &\quad\quad\quad\quad + \Big\|T(t) \int_0^h{T(h-s) B u(s)ds}\Big\|_X^2 \Big) dt.
\end{align*}
Let us compute $I_1$:
\begin{align*}
I_1 
&= \mathop{\overline{\lim}} \limits_{h \rightarrow +0} \frac{1}{h}\Big(   \int_h^{\infty} \|T(t) x\|_X^2 dt - \int_0^{\infty} \|T(t) x\|_X^2 dt \Big) \\
    &= \mathop{\overline{\lim}} \limits_{h \rightarrow +0} - \frac{1}{h}   \int_0^h \|T(t) x\|_X^2 dt 
    = -\|x\|_X^2.
\end{align*}

Now we proceed with $I_2$:
\begin{align*}
I_2 &\leq \mathop{\overline{\lim}} \limits_{h \rightarrow +0}  \int_0^{\infty}\hspace{-1.5mm} 2\big\|T(t+h) x\big\|_X \Big\|T(t) \frac{1}{h} \int_0^h\hspace{-1.5mm}{T(h-s) B u(s)ds}\Big\|_X dt \\
      & \qquad\qquad+ \mathop{\overline{\lim}} \limits_{h \rightarrow +0} \int_0^{\infty}\hspace{-1mm} \frac{1}{h} \Big\|T(t) \int_0^h{T(h-s) B u(s)ds}\Big\|_X^2 dt.
\end{align*}
The limit of the second term equals zero since
\begin{align*}
\mathop{\overline{\lim}} \limits_{h \rightarrow +0} \int_0^{\infty}& \frac{1}{h} \Big\|T(t)   \int_0^h{T(h-s) B u(s)ds}\Big\|_X^2 dt \\
&\leq \mathop{\overline{\lim}} \limits_{h \rightarrow +0} \int_0^{\infty} \frac{1}{h} M^4 e^{-2\lambda t}\|B\| \|u\|_{\cal U}
h^2 dt 
= \ 0.
\end{align*}
To compute the limit of the first term, note that 
\begin{align*}
2\big\|T(t+h) x\big\|_X \Big\|T(t) \frac{1}{h} & \int_0^h{T(h-s) B u(s)ds}\Big\|_X \\
&\leq 2 M \|x\|_X \|T(t) \| M  \|B\| \sup_{r \in [0,h]}\|u(r)\|_{U} \\
&\leq 2 M^3 \|x\|_X  \|B\| \|u\|_{\cal U} e^{-\lambda t}.
\end{align*}
Applying the dominated convergence theorem, Lemma~\ref{AuxiliaryEquality} and Young's inequality, we obtain for any $\eps >0$ that
\begin{eqnarray}
\hspace{-8mm}
I_2 & = & \int_0^{\infty} 2\|T(t) x\|_X \|T(t)  B u(0)\|_X dt \\
    & \leq & \int_0^{\infty} \eps \|T(t) x\|_X^2 + \frac{1}{\eps} \|T(t)  B u(0)\|_X^2 dt \nonumber\\
   & \leq & \int_0^{\infty} \eps \|T(t) \|^2 dt \|x\|_X^2 + \frac{1}{\eps} \int_0^{\infty} \|T(t) \|^2 \|B u(0)\|_X^2 dt \nonumber\\ 
  & \leq &  \frac{\eps M^2}{2\lambda} \|x\|_X^2 + \frac{M^2}{2\lambda \eps}  \|B\|^2 \|u(0)\|_U^2. \nonumber
\label{eq:I_2-estimate-LFs-linear systems}
\end{eqnarray}
Overall, we obtain that for all $x \in X$, all $u \in {\cal U}$, and all $\eps>0$, the inequality \eqref{Final_Lyapunov_Inequality} holds. Considering $\eps < \frac{2\lambda}{M^2}$ this shows that $V$ is a non-coercive ISS Lyapunov function (in dissipative form) for \eqref{eq:Linear_System-ISS-LF}.

It remains to show the Lipschitz continuity of $V$ on bounded sets.
Pick arbitrary $r>0$ and any $x,y \in \clo{B_r}$. It holds that
\begin{align*}
|V(x) - V(y)| &= \Big| \int_0^{\infty} \|T(t) x\|^2_X - \|T(t) y\|^2_X dt \Big| \\
&\leq   \int_0^{\infty} \Big|\|T(t) x\|_X - \|T(t) y\|_X\Big| \big(\|T(t) x\|_X + \|T(t) y\|_X\big) dt \\
&\leq   \int_0^{\infty} \|T(t) x - T(t) y\|_X \big(\|T(t) x\|_X + \|T(t) y\|_X\big) dt \\
&\leq   \int_0^{\infty} Me^{-\lambda t} \|x - y\|_X Me^{-\lambda t} (\|x\|_X + \|y\|_X) dt \\
&\leq   \frac{M^2 r}{\lambda} \|x-y\|_X,
\end{align*}
which shows the Lipschitz continuity of $V$ on bounded sets.
\end{proof}

\begin{remark}
The ISS Lyapunov function
$V$ defined in \eqref{eq:LF_LinSys_Banach} is not coercive in general. 
Noncoercivity of $V$ defined by \eqref{eq:LF_LinSys_Banach} implies that the system
\[\dot{x}=Ax, \quad y=x\]
is not exactly observable on $[0,\infty)$ (even though we can measure the full state!), see \cite[Corollary 4.1.14]{CuZ95}.
The reason for this is that for any given exponential decay rate, there are
states decaying at a faster rate, and thus part of the information about
the state disappears \q{infinitely fast}.
\end{remark}

Below we provide another construction of ISS Lyapunov functions for the system \eqref{eq:Linear_System-ISS-LF} with a bounded input operator.
It is based on a standard construction in the analysis of $C_0$-semigroups, see, e.g., \cite[eq. (5.14)]{Paz83}.


\begin{proposition}
    \label{p:Lyapunovfunction2}
    Let $T$ be exponentially stable with $M\ge 1$, $\lambda>0$ such that
\eqref{eq:SemigroupEstimate} holds and let $0<\gamma < \lambda$.
Then $V^\gamma :X \to \R_+$, defined by
\begin{equation}
    \label{eq:eqnorm-recall}
    V^\gamma (x) := \max_{s\geq 0} \| e^{\gamma s} T(s) x\|_X,
\end{equation}
is a coercive globally Lipschitz continuous ISS Lyapunov function for \eqref{eq:Linear_System-ISS-LF}.

Furthermore, it defines an equivalent norm on $X$, and for
any $u\in {\cal U}$, $x\in X$, we
have the dissipation inequality
\begin{equation}
    \label{eq:dissipation2}
    \dot{V}^\gamma_u(x) \leq - \gamma \ V^\gamma(x) + V^\gamma (Bu(0)) \,.
\end{equation}
\end{proposition}

\begin{proof}
As $T(0)=I$, and in view of \eqref{eq:SemigroupEstimate}, we have for any $x \in X$ that 
\[
\|x\|_X \leq V^\gamma(x) \leq M \max_{s\geq 0}e^{\gamma s}e^{-\lambda s}\|x\|_X = M\|x\|_X.
\]
It is easy to check that $V^\gamma$ is a norm, and the previous estimate indicates that it is equivalent to the original norm on $X$. In particular, $V^\gamma$ is coercive.

We continue by showing Lipschitz continuity of $V^\gamma$. Pick any $x,y \in X$ and assume that
$V^\gamma(x) > V^\gamma(y)$.
Then
\begin{align*}
V^\gamma(x) - V^\gamma(y) &= \max_{s\geq 0} \| e^{\gamma s} T(s) x\|_X - \max_{s\geq 0} \| e^{\gamma s} T(s) y\|_X \\
& \leq \max_{s\geq 0} \big(\| e^{\gamma s} T(s) x\|_X - \| e^{\gamma s} T(s) y\|_X\big) \\        
& \leq \max_{s\geq 0} \big| \| e^{\gamma s} T(s) x\|_X - \| e^{\gamma s} T(s) y\|_X \big|\\
& \leq \max_{s\geq 0} \| e^{\gamma s} T(s)(x-y)\|_X 
 \leq M\|x-y\|_X.        
\end{align*}
The case $V^\gamma(y) > V^\gamma(x)$ can be treated analogously.
This shows that $V^\gamma$ is globally Lipschitz continuous. 

Further, note that for all $x \in X$ and all $t\geq 0$ we have that 
\begin{align}
    \label{eq:eqnorm2}
    V^\gamma(T(t) x) 
 = \max_{s\geq 0} \| e^{\gamma s} T(s) T(t)  x\|_X =
    e^{-\gamma t} \max_{s\geq 0} \| e^{\gamma (s+t)} T(s+t) x\|_X {\leq}
    e^{-\gamma t} V^\gamma(x) \,.
\end{align}
To obtain an infinitesimal estimate, we compute, using the
    triangle inequality ($V^\gamma$ is a norm), the estimate
    \eqref{eq:eqnorm2}, and Lemma~\ref{AuxiliaryEquality},
{\allowdisplaybreaks		
    \begin{align*}
    \dot{V}^\gamma_u&(x) = \mathop{\overline{\lim}} \limits_{h \rightarrow +0} {\frac{1}{h}(V^\gamma(\phi(h,x,u))-V^\gamma(x)) } \\
    &=  \mathop{\overline{\lim}} \limits_{h \rightarrow +0}
    \frac{1}{h}\Big( V^\gamma \Big(T(h) x + \int_0^h{T(h-s) B u(s)ds}\Big) - V^\gamma(x) \Big)  \\
    &\leq \mathop{\overline{\lim}} \limits_{h \rightarrow +0}
    \frac{1}{h}\Big( V^\gamma \big(T(h) x\big) + V^\gamma
    \Big(\int_0^h{T(h-s) B u(s)ds}\Big) - V^\gamma(x) \Big)\\
    &\leq  \mathop{\overline{\lim}} \limits_{h \rightarrow +0}
    \frac{1}{h}\Big( (e^{-\gamma h} -1 ) V^\gamma(x) + V^\gamma
    \Big(\int_0^h{T(h-s) B u(s)ds}\Big)  \Big)    \\
 &=  - \gamma \ V^\gamma(x) + V^\gamma (Bu(0)) \,.
    \end{align*}
}		
In the last inequality, we have used the continuity of $V^\gamma$.
This shows that $V^\gamma$ is an ISS Lyapunov function (in dissipative form) and that
\eqref{eq:dissipation2} holds. 
\end{proof}


\section[ISS Lyapunov functions for linear systems]{Construction of ISS Lyapunov functions for linear systems with unbounded input operators}
\label{sec:Lyapunov_Theorem_systems_unbounded_operators}

In this section, we investigate the applicability of Lyapunov methods to the analysis of linear systems with unbounded input operators.
We start with good news:
\begin{proposition}
\label{prop:for-linear-systems-ncISSLF-implies-ISS}
Assume that $X,U$ are Banach spaces, $q\in [1,+\infty]$ and $B$ is a $q$-admissible operator.
If $q=\infty$, assume additionally that the mild solution $\phi$ of \eqref{eq:Linear_System-ISS-LF} is continuous w.r.t.\ time in the norm of $X$.
Then the existence of a non-coercive ISS Lyapunov function for \eqref{eq:Linear_System-ISS-LF} implies ISS of \eqref{eq:Linear_System-ISS-LF}.
\end{proposition}

\begin{proof}
Under the assumptions of this proposition, \eqref{eq:Linear_System-ISS-LF} is a well-defined control system, satisfying the CEP and BRS properties, see Propositions~\ref{prop:q-admissibility-implies-continuity}, \ref{prop:infty-admissibility-implies-continuity}. 
The claim follows from Theorem~\ref{t:ISSLyapunovtheorem}.
\end{proof}

Proposition~\ref{prop:for-linear-systems-ncISSLF-implies-ISS} enables us to use non-coercive ISS Lyapunov functions as in Theorem~\ref{thm:ISS-criterion-linear-systems-bounded-operators}\ref{itm:vi-ISS-LinSysCrit}.
However, we have \q{to compensate} the influence of unbounded input operators, which poses further requirements on ISS Lyapunov functions.

\subsection{Non-coercive converse Lyapunov theorems for general linear systems}

In the remainder of this section, we specialize to the case of complex Hilbert spaces $X$ endowed with the scalar 
product $\scalp{\cdot}{\cdot}_X$ and to the input spaces
$\Uc:=L^\infty(\R_+,U)$, where $U$ is a Banach space.
The adjoint operator of an operator $P\in L(X)$ we denote by $P^*$, and we call $P\in L(X)$ self-adjoint if $P=P^*$.

We have constructed in Theorem~\ref{thm:ISS-criterion-linear-systems-bounded-operators}\ref{itm:vi-ISS-LinSysCrit} 
a non-coercive ISS Lyapunov function for ISS linear systems with bounded input operators $B\in L(U,X)$.
A Hilbert space version of this result is as follows:
\begin{proposition}
\label{prop:Converse_ISS_Lyapunov_theorem_lin_Systems_with_bounded_input_operators} 
Let $X$ be a complex Hilbert space, $\Uc:=L^\infty(\R_+,U)$, where $U$ is a Banach space, $A$ generate a strongly continuous semigroup on $X$ and let $B\in L(U,X)$.
The system \eqref{eq:Linear_System-ISS-LF} is ISS if and only if there is a self-adjoint operator $P \in L(X)$ so that $\scalp{Px}{x}_X>0$ for $x\neq 0$
and $P$ solves the Lyapunov equation

\begin{eqnarray}
\scalp{Px}{Ax}_X +\scalp{Ax}{Px}_X = -\|x\|^2_X,\quad x\in D(A).
\label{eq:Gn_Lyap_Equation}
\end{eqnarray}
Furthermore, $V:X\to\R_+$ defined by
\begin{align}
V(x) =\scalp{Px}{x}_X
\label{eq:QuadraticLF}
\end{align}
is a non-coercive  ISS Lyapunov function for \eqref{eq:Linear_System-ISS-LF}.
\end{proposition}
It is well-known that \eqref{eq:QuadraticLF} is another representation of the function \eqref{eq:LF_LinSys_Banach-formulation}, see the proof of \cite[Theorem 4.1.3]{CuZ20}.

Proposition~\ref{prop:Converse_ISS_Lyapunov_theorem_lin_Systems_with_bounded_input_operators} 
is not valid for systems with merely admissible operators $B$.
In contrast to systems with bounded input operators, we need further assumptions that relate the operators $P$ and $A$.
The following proposition provides sufficient conditions of this type.
\begin{theorem}
\label{thm:Gen_ISS_LF_Construction}
Let $A$ be  the generator of a $C_0$-semigroup $(T(t))_{t\ge 0}$ on a complex Hilbert space $X$ 
 and let $\Uc:=L^\infty(\R_+,U)$, where $U$ is a Banach space.

Assume that there is an operator $P\in L(X)$ satisfying the following conditions:
\begin{itemize}
	\item[(i)] \label{item:Gen_Converse_Lyap_Theorem_1} $P$ satisfies
\begin{eqnarray}
\re\scalp{Px}{x}_X > 0,\quad x\in X\backslash\{0\}.
\label{eq:Positivity}
\end{eqnarray}
\item[(ii)] \label{item:Gen_Converse_Lyap_Theorem_3}  $\im(P) \subset D(A^*)$.
\item[(iii)] \label{item:Gen_Converse_Lyap_Theorem_4} $PA$ has an
  extension to a bounded operator on $X$, that is, $PA\in L(X)$. (We
  also denote this
  extension by $PA$.)
\item[(iv)] \label{item:Gen_Converse_Lyap_Theorem_2} $P$ satisfies the Lyapunov inequality 
\begin{eqnarray}
\re\scalp{(PA+A^*P)x}{x}_X \leq  -\scalp{x}{x}_X,\quad x\in D(A).
\label{eq:LyapIneq}
\end{eqnarray}
\end{itemize}
Then for any $\infty$-admissible input operator $B\in L(U,X_{-1})$ the function
\begin{equation}
\label{Lyap}
 V(x) := \re\, \langle Px,x\rangle_X
 \end{equation}
is a non-coercive ISS Lyapunov function for \eqref{eq:Linear_System-ISS-LF}. It 
satisfies for each $\eps>0$, all $x_0 \in X$, and all $u \in \Uc$ the dissipation inequality
\begin{align}
\label{ex:Dissipative_Inequality_linear_system}
 \dot{V}_u(x_0) \le  (\varepsilon-1) \|x_0\|_X^2 + c(\eps) \|u\|^2_\infty,
\end{align}
where 
\begin{align}
\label{eq:ceps-def}
c(\eps)&:= \frac{1}{4\varepsilon}\big(\|A^\ast P\|_{L(X)} +\|PA\|_{L(X)}\big)^2 \|A^{-1}_{-1}B\|_{L(U,X)}^2M^2\\
&\qquad\qquad\qquad+ M\|A^\ast P\|_{L(X)} \| A_{-1}^{-1}B\|_{L(U,X)}
\kappa(0),
\nonumber
\end{align}
and
$\kappa(0)=\lim_{t \searrow 0}\kappa(t)$, where  $\kappa(t)>0$ is the smallest constant satisfying
\begin{equation}\label{kappa} \left\|\int_0^t T_{-1}(t-s)Bu(s)\, ds\right\|_X \le \kappa(t) \|u\|_{\infty}, \end{equation}
for every $u\in L^\infty([0,t),U)$. (The existence of the constants $\kappa(t)$ is implied by the $\infty$-admissibility of $B$.)

In particular, the existence of a non-coercive ISS Lyapunov function \eqref{Lyap} implies that \eqref{eq:Linear_System-ISS-LF} is ISS for any $\infty$-admissible $B$.
\end{theorem}

\begin{remark}
\label{rem:Why-Real-parts} 
Note that we have to take the real parts of the expressions in \eqref{Lyap} and \eqref{eq:LyapIneq}, as we deal with complex Hilbert spaces and we do not assume that $P$ is a positive operator on $X$.
\end{remark}

\begin{remark}
\label{rem:Classical_Form_Lyapunov_Equation} 
If, in addition to the assumptions of Theorem
\ref{thm:Gen_ISS_LF_Construction}, the operator $P$ is self-adjoint, then equation \eqref{eq:Gn_Lyap_Equation} is equivalent to \eqref{eq:LyapIneq}. 
\end{remark}

\begin{proof}
Due to Proposition~\ref{prop:for-linear-systems-ncISSLF-implies-ISS}, if $V$ is a non-coercive ISS Lyapunov function for \eqref{eq:Linear_System-ISS-LF}, then \eqref{eq:Linear_System-ISS-LF} is ISS.

By the assumptions
\[
0 < V(x) \leq \|P\|_{L(X)}\|x\|^2_X, \quad x\in X\backslash \{0\},
\]
and thus \eqref{LyapFunk_1Eig_nc_ISS} holds.
It remains to show the dissipation inequality \eqref{DissipationIneq_nc} for $V$.

The operator $A:D(A)\subset X\rightarrow X$ is densely defined as an infinitesimal generator of a $C_0$-semigroup, and hence $A^*$ is well-defined and also the generator of a $C_0$-semigroup, see \cite[Corollary 10.6]{Paz83}. In particular, this implies that $A^*$ is a closed operator. 
Since $P \in L(X)$, the operator $S:=A^*P$ with the domain $D(S):=\{x\in X: Px \in D(A^*)\}$ is a closed operator, see, e.g., \cite[Proposition A.9]{HMM13}.
However, by our assumptions $\im(P) \subset D(A^*)$, which implies $D(S) = X$,
and thus $S = A^*P \in L(X)$ by the closed graph theorem. In particular, the term $\|A^\ast P\|_{L(X)}$ in 
\eqref{eq:ceps-def} makes sense.

For $x_0\in X$ and $u\in \Uc$, we have
{\allowdisplaybreaks
\begin{align}
V(\phi(t,x_0,u))&-V(x_0) \label{eq:MainEquality}\\
&= \re\, \left\langle P\left(T(t)x_0+\int_0^t T_{-1}(t-s)Bu(s)ds\right),\right.\nonumber\\
&\qquad  \left. T(t)x_0+\int_0^t T_{-1}(t-s)Bu(s)ds\right\rangle_X -\re\,\langle Px_0,x_0\rangle_X\nonumber\\
&= \re\,\langle PT(t)x_0,x_0\rangle_X -\re\,\langle Px_0,x_0\rangle_X \nonumber\\
&+ \re\,\langle PT(t)x_0,T(t)x_0\rangle_X - \re\,\langle PT(t)x_0,x_0\rangle_X \label{eq:line2}\\
&+ \re\, \Big\langle PT(t)x_0 ,\int_0^t T_{-1}(t-s)Bu(s)ds\Big\rangle_X \label{eq:line3}\\
&\qquad + \re\,\left\langle P\int_0^t T_{-1}(t-s)Bu(s)ds ,T(t)x_0\right\rangle_X  \label{eq:line4}\\
&+ \re\, \left\langle P\int_0^t T_{-1}(t-s)Bu(s)ds ,\int_0^t T_{-1}(t-s)Bu(s)ds\right\rangle_X.   \label{eq:line5}
\end{align}
}
The terms in the line \eqref{eq:line2} of the previous expression can be transformed into:
\begin{align}
\re\,\langle PT(t)x_0,T(t)&x_0\rangle_X - \re\,\langle PT(t)x_0,x_0\rangle_X = \re\,\langle PT(t)x_0,T(t)x_0 - x_0\rangle_X. \label{eq:line2_Estimate}
\end{align}

Applying \cite[Theorem 4.1.3]{CuZ20} 
to the operator $\tfrac{1}{2}(P+P^*)$, we see that the conditions
 (i), (ii), and (iv) imply that $A$ generates an exponentially stable
 semigroup. Hence, the spectral bound of $A$ is negative, and thus $0\in\rho(A)$ and consequently $A^{-1}\in L(X)$  exists. Further, the exponential stability of the $C_0$-semigroup $(T(t))_{t\ge 0}$ implies 
\[ \|T(t)\|_{L(X)} \le Me^{-\omega t},\quad t\ge 0,\]
for some constants $M$, $\omega>0$. Thanks to $\rho(A)=\rho(A_{-1})$ the operator $A_{-1}^{-1}$ exists  as well. 

By \cite[Theorem II.5.5]{EnN00} 
the map $A:D(A) \to X$ can be continuously extended to the linear isometry
 $A_{-1}$ which maps $(X,\|\cdot\|_X)$ onto $(X_{-1},\|\cdot\|_{X_{-1}})$.
Hence $A_{-1}^{-1}$,
mapping $(X_{-1},\|\cdot\|_{X_{-1}})$ onto $(X,\|\cdot\|_X)$, is again a
  linear isometry, and thus a bounded operator.
As $B\in L(U,X_{-1})$, we have that $A^{-1}_{-1}B \in L(U, X)$. In particular, 
$T_{-1}(t-s)A_{-1}^{-1}Bu(s) = T(t-s)A_{-1}^{-1}Bu(s) \in X$ for all $t\geq 0$ and all $s \in[0,t]$.
Due to the fact, that $A_{-1}^{-1}$ and $T_{-1}(t-s)$ commute, we obtain
{\allowdisplaybreaks
\begin{align}
\Big\| A_{-1}^{-1} \int_0^t T_{-1}(t-s)Bu(s)ds \Big\|_X
&=   \Big\| \int_0^t A_{-1}^{-1}T_{-1}(t-s)Bu(s)ds \Big\|_X\nonumber\\
&=  \Big\| \int_0^t T_{-1}(t-s)A_{-1}^{-1}Bu(s)ds \Big\|_X.\nonumber\\
&\leq    \int_0^t \|T(t-s)\|_{L(X)}\|A_{-1}^{-1}B\|_{L(U,X)}\|u(s)\|_U ds \nonumber\\
&\leq   \int_0^t M ds \|A_{-1}^{-1}B\|_{L(U,X)}\|u\|_\infty \nonumber\\
&\leq   M t \|A_{-1}^{-1}B\|_{L(U,X)}\|u\|_\infty.  \label{eq:est}
\end{align}}

Since $\im(P)\subset D(A^*)$,  we estimate the expression in
\eqref{eq:line3} using the Cauchy-Schwarz inequality and \eqref{eq:est}:
{\allowdisplaybreaks
\begin{align}
 \re\, \Big\langle P T(t)x_0 ,\int_0^t T_{-1}(t-s)&Bu(s)ds\Big\rangle_X\nonumber\\
 &=\re\, \Big\langle P T(t)x_0 ,A A^{-1}\int_0^t  T_{-1}(t-s)Bu(s)ds\Big\rangle_X\nonumber\\
&=\re\,\Big\langle  A^\ast P T(t)x_0 , A_{-1}^{-1} \int_0^t T_{-1}(t-s)Bu(s)ds\Big\rangle_X\nonumber\\
&\leq \|A^\ast P T(t)x_0\|_X \cdot \Big\| A_{-1}^{-1} \int_0^t T_{-1}(t-s)Bu(s)ds \Big\|_X\nonumber\\
&\leq  \|A^\ast P\|_{L(X)} \|T(t)x_0\|_X \cdot  M t \|A_{-1}^{-1}B\|_{L(U,X)}\|u\|_\infty.  \label{eq:line3_Estimate}
\end{align}
}
To upperestimate the expression \eqref{eq:line4}, we use again \eqref{eq:est} to obtain 
\begin{align}
\re\,\Big\langle P\int_0^t T_{-1}(t-s)Bu(s)ds,& T(t)x_0\Big\rangle_X\nonumber\\
&= \re\,\left\langle PAA^{-1}\int_0^t T_{-1}(t-s)Bu(s)ds ,T(t)x_0\right\rangle_X\nonumber\\
&\leq \|PA\|_{L(X)} \cdot  M t \|A_{-1}^{-1}B\|_{L(U,X)}\|u\|_\infty \|T(t)x_0\|_X.  \label{eq:line4_Estimate_b}
\end{align}
Finally, we estimate the expression \eqref{eq:line5} using \eqref{eq:est} and $\kappa$ as defined in \eqref{kappa}
\begin{align}
\re\,\Big\langle P\int_0^t T_{-1}&(t-s)Bu(s)ds ,\int_0^t T_{-1}(t-s)Bu(s)ds\Big\rangle_X\nonumber\\
& =  \re\,\left\langle P\int_0^t \!T_{-1}(t-s)Bu(s)ds ,AA^{-1}_{-1}\int_0^t \!T_{-1}(t-s)Bu(s)ds\!\!\right\rangle_X\nonumber\\
& =  \re\,\left\langle A^*P\!\int_0^t \!T_{-1}(t-s)Bu(s)ds ,A^{-1}_{-1}\!\int_0^t \!T_{-1}(t-s)Bu(s)ds\!\!\right\rangle_X\nonumber\\
&\leq  \|A^*P\|_{L(X)} \kappa(t)\|u\|_\infty \cdot Mt \|A_{-1}^{-1}B\|_{L(U,X)}\|u\|_\infty.   \label{eq:line5_Estimate}
\end{align}
Substituting \eqref{eq:line2_Estimate}, \eqref{eq:line3_Estimate}, 
\eqref{eq:line4_Estimate_b}, and \eqref{eq:line5_Estimate} into \eqref{eq:MainEquality}, we obtain:
\begin{align*}
V(\phi(t,x_0,u))-V(x_0) \le \,& \re\, \langle PT(t)x_0 - Px_0, x_0\rangle_X + \re\,\langle PT(t)x_0,T(t)x_0 - x_0\rangle_X\\
&+\|A^\ast P\|_{L(X)} \|T(t)x_0\|_X \cdot  M t \|A_{-1}^{-1}B\|_{L(U,X)}\|u\|_\infty\\
&+ \|PA\|_{L(X)} \cdot  M t \|A_{-1}^{-1}B\|_{L(U,X)}\|u\|_\infty \|T(t)x_0\|_X\\
& + \|A^*P\|_{L(X)} \kappa(t)\|u\|_\infty \cdot Mt \|A_{-1}^{-1}B\|_{L(U,X)}\|u\|_\infty.
\end{align*}
For $x_0\in X$, we have $A^{-1}x_0\in D(A)$ and we obtain from the definition of the generator $A$ that
\begin{align*}
\mathop{\overline{\lim}}_{t\searrow 0} \re\, \frac{1}{t}\langle PT(t)x_0 - Px_0, x_0\rangle_X
& =\mathop{\overline{\lim}}_{t\searrow 0} \re\, \frac{1}{t}\langle PA[ T(t)A^{-1}x_0 - A^{-1}x_0], x_0\rangle_X\\
&= \re\scalp{PAx_0}{x_0}_X
\end{align*}
and similarly
\begin{align*}
\mathop{\overline{\lim}}_{t\searrow 0} \frac{1}{t}\re\,\langle PT(t)x_0,T(t)x_0 - x_0\rangle_X 
& =\mathop{\overline{\lim}}_{t\searrow 0} \re\, \frac{1}{t}\langle A^*PT(t)x_0 , T(t)A^{-1}x_0 - A^{-1}x_0\rangle_X\\
&= \re\scalp{A^* Px_0}{x_0}_X.
\end{align*}
This implies for every $\varepsilon>0$ that (recall the definition of $\kappa(0)$ before \eqref{kappa})
{\allowdisplaybreaks
\begin{align*}
 \dot{V}_u(x_0) =\mathop{\overline{\lim}}_{t\searrow 0} \frac{1}{t}&(V(\phi(t,x_0,u))-V(x_0))\\
&\le  \re\scalp{PAx_0}{x_0}_X + \re\,\scalp{A^*Px_0}{x_0}_X \\
&\qquad  + \|A^\ast P\|_{L(X)} \|x_0\|_X  \|A_{-1}^{-1}B\|_{L(U,X)} M  \|u\|_\infty\\
&\qquad+  \|PA\|_{L(X)}\|  A_{-1}^{-1}B\|_{L(U,X)}  M  \|u\|_\infty  \|x_0\|_X\\
&\qquad +   \|A^\ast P\|_{L(X)} \|  A_{-1}^{-1}B\|_{L(U,X)}  M  \kappa(0)\|u\|^2_\infty\\
&=    \re\scalp{PAx_0}{x_0}_X + \re\,\scalp{A^*Px_0}{x_0}_X \\
& \qquad +  \| x_0\|_X (\|A^\ast P\|_{L(X)} + \|PA\|_{L(X)})  \|A_{-1}^{-1}B\|_{L(U,X)} M  \|u\|_\infty\\
& \qquad+  \|A^\ast P\|_{L(X)} \|  A_{-1}^{-1}B\|_{L(U,X)}  M  \kappa(0)\|u\|^2_\infty.
\end{align*}
Using Young's inequality and the estimate \eqref{eq:LyapIneq}, we arrive at
\begin{align*}
 \dot{V}_u(x_0) &\le  - \|x_0\|_X^2 + \varepsilon  \|x_0\|_X^2  + \frac{(\|A^\ast P\|_{L(X)} +\|PA\|_{L(X)})^2 \|A^{-1}_{-1} B\|_{L(U,X)}^2M^2 }{4\varepsilon} \|u\|^2_\infty  \\
&\qquad + \|A^\ast P\|_{L(X)}\|  A_{-1}^{-1}B\|_{L(U,X)}  M  \kappa(0)\|u\|^2_\infty,
\end{align*}
which shows the dissipation inequality \eqref{ex:Dissipative_Inequality_linear_system}, and thus also
\eqref{DissipationIneq_nc}.
}
\end{proof}

\begin{remark}
\label{rem:Reformulation_of_Linear_Theorem} 
Theorem~\ref{thm:Gen_ISS_LF_Construction} has been formulated as a direct
Lyapunov theorem. However, the following reformulation as a partial converse result is also possible. Assume that \eqref{eq:Linear_System-ISS-LF} is ISS, and the solution $P$ of the Lyapunov equation \eqref{eq:Gn_Lyap_Equation} satisfies $\im(P) \subset D(A^*)$ and $PA$ is bounded. Then \eqref{Lyap} is an ISS Lyapunov 
function for \eqref{eq:Linear_System-ISS-LF}. 
\end{remark}

It is of virtue to compare the ISS Lyapunov theorem for bounded input operators (Proposition~\ref{prop:Converse_ISS_Lyapunov_theorem_lin_Systems_with_bounded_input_operators}) and ISS Lyapunov theorem for admissible input operators (Theorem~\ref{thm:Gen_ISS_LF_Construction}). The ISS Lyapunov function candidate considered in both these results is the same. What differs is the assumptions and the set of input operators, for which this function is indeed an ISS Lyapunov function. Proposition~\ref{prop:Converse_ISS_Lyapunov_theorem_lin_Systems_with_bounded_input_operators} states that if the semigroup generated by $A$ is exponentially stable, then there is an operator $P$, which satisfies the assumptions (i) and (iv) of Theorem~\ref{thm:Gen_ISS_LF_Construction} and the condition $P=P^*$. Furthermore, \eqref{eq:QuadraticLF} is an ISS Lyapunov function for \eqref{eq:Linear_System-ISS-LF} for any bounded input operator. 
\emph{Thus, the additional key assumptions which we impose to tackle the unboundedness of an input operator are assumptions (ii) and (iii).} 
We note that with these assumptions, \eqref{Lyap} is an ISS Lyapunov function for any $\infty$-admissible operator $B$.

In the next subsections, we show the applicability of Theorem~\ref{thm:Gen_ISS_LF_Construction} for some important special cases.
We start with sufficient conditions, guaranteeing that Theorem~\ref{thm:Gen_ISS_LF_Construction} can be applied with $P=-A^{-1}$.
Then we show that these sufficient conditions are fulfilled for broad classes of systems generated by subnormal operators. Finally, we proceed to diagonal semigroups (whose generators are self-adjoint operators), and finally, we give a construction of a non-coercive ISS Lyapunov function for the heat equation with Dirichlet boundary inputs.

\subsection{A special case: $P=-A^{-1}$}
\label{sec:choosing_A_inverse}

In this section, we give sufficient conditions for the applicability of Theorem~\ref{thm:Gen_ISS_LF_Construction} with $P:=-A^{-1}$.

\begin{proposition}
\label{prop:Conv_ISS_LF_Theorem_LinOp}
Let $A$ be  the generator of an exponentially stable $C_0$-semigroup $(T(t))_{t\ge 0}$ on a (complex) Hilbert space $X$
 and let $\Uc:=L^\infty(\R_+,U)$, where $U$ is a Banach space.

Further, assume that
\begin{itemize}
	\item[(a)] $D(A)\subset D(A^\ast)$.
	\item[(b)] There is a certain $\delta\in(0,1)$ such that for every $x \in X$ we have 
\begin{equation}\label{eqn:a2}
 \re\, \langle A^*A^{-1}x, x\rangle_X + \delta \|x\|_X^2\ge 0.
\end{equation}
	\item[(c)] $  \re\,\langle Ax,x\rangle_X <0$ holds for  every $x\in D(A)\backslash\{0\}$.
\end{itemize}
Then
\begin{equation}
\label{Lyap-Apower--1}
 V(x) := - \re\, \langle A^{-1}x,x\rangle_X
 \end{equation}
is an ISS Lyapunov function for  \eqref{eq:Linear_System-ISS-LF} for any $\infty$-admissible operator $B\in L(U,X_{-1})$.
\end{proposition}

\begin{proof}
As $A$ generates an exponentially stable semigroup, $0 \in\rho(A)$ and thus $P:=-A^{-1} \in L(X)$. We show step by step 
that this choice of $P$ satisfies all the requirements (i)--(iv) of Theorem~\ref{thm:Gen_ISS_LF_Construction}.\\
(i). For any $x\in X\backslash\{0\}$ there is $y \in D(A)\backslash\{0\}$ so that $x = Ay$. By the assumptions of the proposition, it holds that
\[
V(x) = -\re\scalp{y}{Ay}_X = -\re\scalp{Ay}{y}_X >0.
\]
(ii). We have $\im(P)=\im(A^{-1}) = D(A) \subset D(A^*)$, which holds by our assumptions.\\
(iii). Trivial as $PA = -I$.\\
(iv).  
By assumptions, there is a $\delta<1$ so that
\begin{eqnarray*}
\re\scalp{(PA+A^*P)x}{x}_X = \re\scalp{(-I-A^*A^{-1})x}{x}_X
												&=& -\scalp{x}{x}_X - \re\scalp{A^*A^{-1} x}{x}_X\\
												&\leq& -(1-\delta)\scalp{x}{x}_X,
\end{eqnarray*}
and thus, $P$ satisfies the Lyapunov inequality up to a scaling coefficient (and $\tilde{P}:=\frac{1}{1-\delta}P$ satisfies precisely \eqref{eq:LyapIneq}).

Hence, all assumptions of Theorem~\ref{thm:Gen_ISS_LF_Construction} are satisfied, and an application of 
Theorem~\ref{thm:Gen_ISS_LF_Construction} shows the claim.
\end{proof}

Proposition~\ref{prop:Conv_ISS_LF_Theorem_LinOp} can be used to obtain converse non-coercive ISS Lyapunov theorems for linear systems with unbounded input operators and subnormal generators, see \cite[Corollary 6.7]{JMP20}. Next, we present an important special case of this result for systems governed by self-adjoint operators:
\begin{corollary}
\label{cor:Self-adjoint-A}
Let $X$ be a complex Hilbert space and let $(A,D(A))$ be a negative definite self-adjoint operator on $X$, that is
$\scalp{Ax}{x}_X <0$ for all $x\in D(A)\backslash\{0\}$ and $D(A) = D(A^*)$ with $Ax = A^*x$ for all $x\in D(A)$.
Then 
\[
V(x) := - \langle A^{-1}x,x\rangle_X
\]
is a non-coercive ISS Lyapunov function for \eqref{eq:Linear_System-ISS-LF} with any $\infty$-admissible operator $B$.
\end{corollary}

\begin{proof}
Take $P:=-A^{-1}$. Since $A = A^*$, for all $x \in X$ we have $\scalp{A^{-1}x}{x}_X\in\R$, as
\[
\scalp{A^{-1}x}{x}_X = \scalp{A^{-1}x}{A^*A^{-1}x}_X = \scalp{AA^{-1}x}{A^{-1}x}_X = \scalp{x}{A^{-1}x}_X.
\]
All assumptions of Proposition~\ref{prop:Conv_ISS_LF_Theorem_LinOp} are fulfilled, and its application shows the claim.
\end{proof}

\subsection{ISS Lyapunov functions for input-to-state stable diagonal systems}
\label{sec:Example}

Consider a linear system \eqref{eq:Linear_System-ISS-LF} with the state space
\[
X=\ell_2(\N):=\big\{ x=(x_k): \|x\|_{X}=\Big(\sum_{k=1}^{\infty} |x_k|^2 \Big)^{1/2}< \infty  \big\}
\] 
endowed in the usual way with the scalar product $\lel  \cdot,\cdot \rir_{\ell_2}$.
Let $U:=\R$ and $\Uc:=L^\infty(\R_+,U)$.

Consider the operator $A:X\to X$, defined by $Ae_k=- \lambda_k e_k$, where $e_k$ is the $k$-th unit
vector of $\ell_2(\N)$ and $\lambda_k \in \R$ with $\lambda_k<\lambda_{k+1}$ for all $k$, $\lambda_1>\varepsilon>0$ and $\lambda_k \to \infty$  as $k \to \infty$.

\index{system!diagonal}
The operator $A$ can be represented using the spectral decomposition
\begin{eqnarray}\label{a1}
Ax:=\sum_{k=1}^\infty - \lambda_k \lel  x,e_k \rir_{\ell_2} e_k, \quad x \in D(A),
\label{eq:Generator_diagonal_semigroups}
\end{eqnarray}
with
\begin{eqnarray}\label{a2}
D(A) = \{x \in \ell_2(\N): \sum_{k=1}^\infty - \lambda_k \lel  x,e_k \rir_{\ell_2} e_k \text{ converges}\}.
\label{eq:Domain_generator_diagonal_semigroups}
\end{eqnarray}

We have the following result:
\begin{proposition}
\label{prop:Converse_Lyapunov_Theorem_Diagonal_System} 
Let $A$ be given by \eqref{a1}, \eqref{a2} and  $B \in L(\mathbb C^m,X_{-1})$. Then \eqref{eq:Linear_System-ISS-LF} is $L_\infty$-ISS and 
\begin{eqnarray}
V(x):=\sum_{k=1}^\infty \frac{1}{\lambda_k} \lel  x,e_k \rir_{\ell_2}^2
\label{eq:ISS_LF_diagonal_semigroups}
\end{eqnarray}
is a non-coercive ISS Lyapunov function for \eqref{eq:Linear_System-ISS-LF}.
\end{proposition}

\begin{proof}
By assumptions, the operator  $A$ is self-adjoint with  $\sigma(A) \subset (-\infty,0)$. Thus the assumptions of Corollary~\ref{cor:Self-adjoint-A} hold. 
Furthermore, the inverse of $A$ is given by
\begin{eqnarray}
A^{-1}x:=\sum_{k=1}^\infty - \frac{1}{\lambda_k} \lel  x,e_k \rir_{\ell_2} e_k,
\label{eq:Inverse_Generator_diagonal_semigroups}
\end{eqnarray}
and thus, the Lyapunov function \eqref{Lyap-Apower--1} has the form \eqref{eq:ISS_LF_diagonal_semigroups}.
\end{proof}


%
%

It is easy to see that the ISS Lyapunov function $V$ given by \eqref{eq:ISS_LF_diagonal_semigroups} is not coercive since $\lambda_k \to \infty$ as $k \to \infty$. 

\subsection{ISS Lyapunov functions for the heat equation with Dirichlet boundary input}
\label{sec:ISS Lyapunov functions for a heat equation with Dirichlet boundary input}

It is well-known that the classical heat equation with Dirichlet boundary inputs is ISS with $\Uc:=L^\infty(\R_+,U)$, which has been verified by means of several different methods \cite{JNP18, KaK16b, MKK19}. However, no constructions for ISS Lyapunov functions have been proposed. 
In the next example, we show that by using Theorem~\ref{thm:Gen_ISS_LF_Construction}
 one can construct a non-coercive ISS Lyapunov function for this system.
\begin{example}
\label{ex1}
Consider the heat equation on $[0,1]$ with a Dirichlet boundary input at the point $1$:
\begin{subequations}
\label{eq:Linear-heat-equation}
\begin{align}
x_t(z,t)& =a x_{zz}(z,t), \quad z\in(0,1),~ t>0,\\
  x(0,t)&=0,  \quad x(1,t)=u(t), \quad t>0,\\
x(z,0)&=x_{0}(z),
\end{align}
\end{subequations}
where $a>0$. 
We choose
 $X=L^{2}(0,1)$, $U=\mathbb C$ and  $\Uc:=L^\infty(\R_+,\C)$.

The system \eqref{eq:Linear-heat-equation} can be rewritten in the form \eqref{eq:Linear_System-ISS-LF} with $A$ defined by
\begin{align*}
Af={}af'', \quad f\in D(A):={}\left\{f\in H^{2}(0,1) : f(0)= f(1)=0\right\}.
\end{align*}
Here $H^2(0,1)$ denotes the Sobolev space of functions $f \in L^2(0,1)$, which have weak derivatives of order $\leq 2$, all of which belong to $L^2(0,1)$.
It is well-known that $A$ is a self-adjoint operator on $X$ generating an exponentially stable analytic $C_0$-semigroup on $X$. 

By Theorem~\ref{thm:Contiuity-of-a-map}, $B\in L(U, X_{-1})$ is $\infty$-admissible, for every $x_0\in X$ and $u\in L^\infty(0,\infty)$ the corresponding mild solution is continuous with respect to time and  $\kappa(0)=0$.
In \cite{JNP18}, the following ISS estimates have been shown:
\begin{align*}
\|\phi(t,x_0,u)\|_{L^2(0,1)}&\le {\rm{e}}^{-a\pi^{2} t}\|x_0\|_{L^2(0,1)} +\frac{1}{\sqrt{3}} \|u\|_{L^\infty(0,t)},\\
\|\phi(t,x_0,u)\|_{L^2(0,1)}&\le {\rm{e}}^{-a\pi^{2} t}\|x_0\|_{L^2(0,1)} +c \left( \int_0^t |u(s)|^pds \right)^{1/p},
\end{align*}
for every $x_0\in X$, $u\in \Uc$, $p> 2$ and some constant $c=c(p)>0$. 
The direct application of Corollary~\ref{cor:Self-adjoint-A} shows that 
\begin{eqnarray*}
 V(x) &=& -  \langle A^{-1}x,x\rangle_X = \int_0^1 \left(\int_z^1  (z -\tau) x(\tau) d\tau \right)\overline{x(z)} dz 
 \end{eqnarray*}
is a non-coercive ISS Lyapunov function for \eqref{eq:Linear-heat-equation}, and thus \eqref{eq:Linear-heat-equation} is ISS.
\end{example}

\section{ISS superposition theorems for Lipschitz continuous systems}

For systems \eqref{InfiniteDim} satisfying an additional property which was needed in Theorem~\ref{Characterization_LISS} to show the equivalence of LISS and 0-UAS, one can strengthen the ISS superposition theorems shown in Theorem~\ref{thm:UAG_equals_ULIM_plus_LS}:
\begin{theorem}
\label{thm:MainResult_Characterization_ISS_EQ_Banach_Spaces}
Let \eqref{InfiniteDim} be forward complete, satisfy Assumption~\ref{Assumption1}, $f(0,0)=0$, and let there exist $\sigma \in \K$ and $\rho >0$ so that for all $v \in U$: $\|v\|_U \leq \rho$ and all $x \in X$: $\|x\|_X  \leq \rho$ we have
\begin{eqnarray*}
\|f(x,v)-f(x,0)\|_X \leq \sigma(\|v\|_U).
\end{eqnarray*}
The following statements are equivalent:
\begin{itemize}
    \item[(i)] \eqref{InfiniteDim}   is ISS.
    \item[(ii)] \eqref{InfiniteDim}  is bUAG and BRS.
    \item[(iii)] \eqref{InfiniteDim} is bULIM, ULS, and BRS.
    \item[(iv)] \eqref{InfiniteDim}  is bULIM and UGS.
    \item[(v)] \eqref{InfiniteDim}   is bULIM, 0-ULS, and BRS.
\end{itemize}
\end{theorem}

\begin{proof}
    By Lemma~\ref{lem:RobustEquilibriumPoint}, our assumptions imply the CEP property of \eqref{InfiniteDim}. In conjunction with the ISS superposition theorem for general forward complete control systems (Theorem~\ref{thm:UAG_equals_ULIM_plus_LS}), this shows the equivalence of  (i)--(iv).
Clearly, (iii) implies (v). 

(v) $\Rightarrow$ (iii). By the BRS property, the value
  \begin{equation*}
      \tilde{\beta}(r,t) := \sup \{ \|\phi(t,x,0)\|_X \midset x \in B_r \} 
  \end{equation*}
is finite for all $(r,t) \in \R_+^2$. The function $\tilde{\beta}$ is
nondecreasing in $r$. 
By 0-ULS of \eqref{InfiniteDim},  $\tilde\beta$ is continuous in the first argument at $0$.

Also, for fixed $r\geq 0$, we claim that $\lim_{t\to\infty}
 \tilde{\beta}(r,t) = 0$ by bULIM and 0-ULS. To see this, let $\sigma$ be
 the function characterizing 0-ULS. Given $\varepsilon>0$, we may by bULIM
 choose a $\tau >0$ such that for all $x\in B_r$ there is a $t\leq \tau$
 with 
 \begin{equation*}
     \|\phi(t,x,0)\|_X \leq \sigma^{-1}(\varepsilon). 
 \end{equation*}
 By 0-ULS and the cocycle property, it follows that $\|\phi(t,x,0)\|_X \leq
 \varepsilon$ for all $t\geq \tau$ and all $x \in B_r$ so that we have the desired
 convergence. We now have that for all $x\in X$ and all $t\geq 0$
\begin{equation*}
     \|\phi(t,x,0)\|_X \leq \max \{ \tilde{\beta}(\|x\|_X,t+s) \midset s \geq 0 \}    + \|x\|_X e^{-t}. 
 \end{equation*}
This upper bound is a well-defined function of $(\|x\|_X,t)$, continuous w.r.t.\  the first argument at $\|x\|_X=0$, strictly increasing in $\|x\|_X$ and
strictly decreasing to 0 in $t$. 
It is easy to see that there is a $\beta \in \mathcal{KL}$ so that
\begin{equation*}
     \|\phi(t,x,0)\|_X \leq \beta(\|x\|_X,t),
 \end{equation*}
and thus \eqref{InfiniteDim} is 0-UGAS.

Hence, Theorem~\ref{Characterization_LISS} implies LISS (and, in particular, ULS) of \eqref{InfiniteDim}.
\end{proof}

Finally, we summarize the obtained criteria for well-posed forward complete systems \eqref{InfiniteDim} with Lipschitz continuous nonlinearities in Figure~\ref{fig:ISS_Equiv_Evolution_Equations}.

\begin{figure*}[tbh]
\centering
\begin{tikzpicture}[>=implies,thick]
\small

\node (ISS) at (2.7,0.4) {ISS};
\node (bUAG) at (0.8,0.4) {bUAG};

\node (ISSLF) at (5.5,0.4) {\color{blue} $\exists$ ISS-LF};
 
\draw [rounded corners] (-6.55,1.2) rectangle (3.1,0.1);
\node (Thm5) at (-1.5,0.9) {\footnotesize Thm.~\ref{thm:UAG_equals_ULIM_plus_LS}};

\node (Thm6) at (-7.2,-0.1) {\footnotesize \color{blue} Thm.~\ref{thm:MainResult_Characterization_ISS_EQ_Banach_Spaces}};
\node (Thm3) at (4.2,0.9) {\footnotesize \color{blue} Thm.~\ref{ISS_Converse_Lyapunov_Theorem}
};

 \node (Thm8) at (-3.2,-1) {\footnotesize Thm.~\ref{wISS_equals_sAG_GS}};

 \node (Thm2) at (4.6,-2.38) {\footnotesize \color{blue} Thm.~\ref{Characterization_LISS}};

\node (Rem) at (-4,-2.35) {\footnotesize Rem.~\ref{rem:LIMAG}};

\draw [rounded corners] (-6.55,-1.8) rectangle (-0.2,-0.7);

\draw[thick,double equal sign distance,<->] (ISS) to (bUAG);
\draw[thick,double equal sign distance, blue,dashed,<->] (ISSLF) to (ISS);

\node (bULIM_UGS) at (-2,0.4) {bULIM$\,\wedge\,$UGS};
\node (bULIM_ULS) at (-5.2,0.4) {bULIM$\,\wedge\,$ULS};
\node (bULIM_0ULS)at (-8.7,0.4) {\color{blue} bULIM$\,\wedge\,$0-ULS};

\draw[thick,double equal sign distance,<->] (bULIM_UGS) to (bUAG);
\draw[thick,double equal sign distance,<->] (bULIM_ULS) to (bULIM_UGS);
\draw[thick,double equal sign distance,blue,dashed,<->] (bULIM_0ULS) to (bULIM_ULS);

\node (sISS) at (-0.6,-1.5) {sISS};
\node (sAG_UGS) at (-2.7,-1.5) {sAG$\,\wedge\,$UGS};
\node (sLIM_UGS) at (-5.4,-1.5) {sLIM$\,\wedge\,$UGS};
\draw[<->,double equal sign distance] (sISS) to (sAG_UGS);

\draw[<->,double equal sign distance] (sLIM_UGS) to (sAG_UGS);

\draw[thick,double equal sign distance,->] (-3.3,0) to (-3.3,-0.6);
\draw[red,double equal sign distance,->,degil]  (-2.95,-0.6) to (-2.95,0);
\node[red] (notImp_ix) at (-2.45,-0.3) {\footnotesize(ix)};

\draw[red,double equal sign distance,->,degil]  (0.7,-2.45) to (0,-2.45);
\node[red] (vi) at (1.1,-2.45) {\footnotesize(vi)};

\draw[red,double equal sign distance,<-,degil]  (0.7,-2.95) to (0,-2.95);
\node[red] (vii) at (1.1,-2.95) {\footnotesize(vii)};

\draw[red,double equal sign distance,->,degil]  (0.7,-1.3) to (0,-1.3);
\node[red] (iv) at (1.1,-1.3) {\footnotesize(iv)};

\draw[red,double equal sign distance,<-,degil]  (0.7,-1.7) to (0,-1.7);
\node[red] (v) at (1.1,-1.7) {\footnotesize(v)};



\node (AG_UGS) at (-2.7,-2.75) {AG$\,\wedge\,$UGS};
\node (LIM_UGS) at (-5.4,-2.75) {LIM$\,\wedge\,$UGS};

\node (AG_ULS) at (-2.7,-3.9) {AG$\,\wedge\,$ULS};

\node (AG_0UGAS) at (2.7,-1.5) {AG$\,\wedge\,$0-UGAS};
\node (AG_LISS) at (5.5,-2.75) {\color{blue}AG$\,\wedge\,$LISS};
\node (AG_0UAS) at (2.7,-2.75) {AG$\,\wedge\,$0-UAS};
\node (AG_0LS) at (2.7,-3.9) {AG$\,\wedge\,$0-ULS};
\node (AG_0GAS) at (5.5,-3.9) {AG$\,\wedge\,$0-GAS};


\draw[double equal sign distance,->]  ($(AG_UGS)+(-0.2,-0.3)$) to ($(AG_ULS)+(-0.2,0.3)$);
\draw[red,double equal sign distance,<-,degil]  ($(AG_UGS)+(0.2,-0.3)$) to ($(AG_ULS)+(0.2,0.3)$);
\node[red] (viii) at ($\weight*(AG_ULS)+\weight*(AG_UGS) + (0.75,0)$) {\footnotesize(viii)};

\draw[double equal sign distance,->]  ($(AG_0UGAS)+(-0.2,-0.3)$) to ($(AG_0UAS)+(-0.2,0.3)$);
\draw[red,double equal sign distance,<-,degil]  ($(AG_0UGAS)+(0.2,-0.3)$) to ($(AG_0UAS)+(0.2,0.3)$);
\node[red] (ii) at ($\weight*(AG_0UAS)+\weight*(AG_0UGAS) + (0.6,0)$) {\footnotesize(ii)};

\draw[double equal sign distance,->]  ($(AG_0UAS)+(-0.2,-0.3)$) to ($(AG_0LS)+(-0.2,0.3)$);
\draw[red,double equal sign distance,<-,degil]  ($(AG_0UAS)+(0.2,-0.3)$) to ($(AG_0LS)+(0.2,0.3)$);
\node[red] (i) at ($\weight*(AG_0LS)+\weight*(AG_0UAS) + (0.6,0)$) {\footnotesize(i)};

%

\coordinate (ISS_1) at ($(ISS) + (-0.5,0)$);
\coordinate (ISS_2) at ($(ISS) + (0,-0.4)$);
\coordinate (AG_0UGAS_1) at ($(AG_0UGAS) + (-0.5,0)$);
\coordinate (coord1) at ($(AG_0UGAS_1) + (0,0.3)$);
\coordinate (coord2) at ($(ISS_1) + (0,-0.4)$);

\draw[red,double equal sign distance,->,degil]  (coord1) to (coord2);
\node[red] (iii) at ($\weight*(ISS_1) + \weight*(AG_0UGAS_1) + (-0.45,-0.15)$) {\footnotesize(iii)};

\draw[thick,double equal sign distance,->] ($(sAG_UGS)+(-0.2,-0.37)$) to ($(AG_UGS)+(-0.2,0.27)$);

  \coordinate (A2) at ($(sAG_UGS)+(0.2,-0.37)$);
  \coordinate (B2) at ($(AG_UGS)+(0.2,0.27)$);
\node[red] (Quest) at (-2.5,-2.23) {\footnotesize ???};
\draw[double equal sign distance,-] (B2) to ($(B2)+(0,0.09)$);
\draw[double equal sign distance,->] ($(B2)+(0,0.42)$) to (A2);

\draw[double equal sign distance,->]  ($(AG_UGS)+(-0.2,-0.3)$) to ($(AG_ULS)+(-0.2,0.3)$);

\draw[thick,double equal sign distance,<->] (AG_UGS) to (LIM_UGS);

\draw[thick,double equal sign distance,->] (ISS_2) to (AG_0UGAS);
\draw[thick,double equal sign distance,blue,dashed,<->]  (AG_LISS) to (AG_0UAS);

\draw[thick,double equal sign distance,<->]  (AG_0LS) to (AG_0GAS);

\draw[thick,double equal sign distance,->] ($(AG_ULS)+(0.9,-0.12)$)  to ($(AG_0LS)+(-1.1,-0.12)$);

  \coordinate (A1) at ($(AG_ULS)+(0.9,0.12)$);
  \coordinate (B1) at ($(AG_0LS)+(-1.1,0.12)$);
\path (B1) -- node (Q) {\footnotesize\color{red}{???}} (A1);
    \draw[thick,double equal sign distance,->] (B1) -- (Q) -- (A1);


\end{tikzpicture}
\caption[caption]{Relations between stability properties of
  infinite-dimensional systems, which have a robust equilibrium point and
  bounded reachability sets:
\begin{itemize}
\setlength{\itemindent}{5mm}
   \item Black arrows show implications or equivalences which hold for general control systems
in infinite dimensions.
\item {\color{red}Red arrows (with the negation sign)}
are implications which do not hold, in view of available counterexamples, see Remark~\ref{rem:Nonimplications}.
\item {\color{blue}Blue dashed equivalences} are proved only for systems of the form \eqref{InfiniteDim} and under certain additional conditions.
\item Black arrows with question marks inside mean that it is not known right now (as far as the author is concerned), whether the converse implications hold or not.
\end{itemize}
}
\label{fig:ISS_Equiv_Evolution_Equations}
\end{figure*}

\section{Concluding remarks}		

The use of non-coercive Lyapunov functions for linear systems without disturbances is quite classic. 
Lyapunov functions, constructed by using operator Lyapunov equation  \cite[Theorem 4.1.3]{CuZ20}, are in general non-coercive. 
Certainly, it has a special structure (in fact, this Lyapunov function is precisely the functional used in Datko-Pazy Theorem with $p=2$). However, Littman's and Rolewicz's extensions of Datko-Pazy result \cite{Lit89}, \cite[Theorem 3.2.2]{Nee96} indicate that also other types of non-coercive Lyapunov functions are certificates for the exponential stability of a semigroup. More recently, non-coercive converse Lyapunov theorems for infinite-dimensional switched systems were reported in \cite{HaS11}.

For nonlinear systems, non-coercive Lyapunov functions have been introduced in \cite{MiW19a}, where they were used for the analysis of the uniform global asymptotic stability. The existence of non-coercive Lyapunov functions has been related to the integral-type stability notions for nonlinear systems in \cite{MiW19b}, which refines some results in \cite{MiW19a}.
The results on non-coercive ISS Lyapunov functions discussed in Section~\ref{sec:Non-coercive ISS Lyapunov functions} are due to \cite{JMP20} and constitute extensions of the results from \cite{MiW19a} to systems with inputs.

Theorem~\ref{Characterization_LISS} (characterization of LISS) has been proved in \cite[Theorem 4]{Mir16}. The proof technique has been motivated by \cite[Corollary  4.2.3]{Hen81}. 
Specialized to ODE systems, this result was reported in \cite[Theorem 1.4]{LJH14}. A slightly weaker result claiming that for ODE systems, 0-GAS implies LISS, has been shown in \cite[Lemma I.1]{SoW96}.
Equivalence between 0-UAS and LISS is reminiscent of the equivalence between asymptotic and so-called total stability (see \cite[Theorems 56.3, 56.4]{Hah67}). 

The converse Lyapunov Theorem~\ref{ISS_Converse_Lyapunov_Theorem} was shown in \cite[Theorem 5]{MiW17c}. For ODEs, a stronger result was shown in  \cite[Theorem 1]{SoW95}, where the existence of a smooth ISS Lyapunov function for an ISS system has been proved. 
Our proof strategy follows that of \cite{SoW95}. A key step in the proof is the invocation of the converse Lyapunov theorem for the auxiliary system with disturbances \eqref{eq:Modified_InfDimSys_With_Disturbances}. For ODE systems, a smooth converse Lyapunov theorem from \cite{LSW96} is used in \cite{SoW95}. The author is not aware of such a smooth converse result for infinite-dimensional systems. Thus, we rely on Lipschitz continuous converse Lyapunov 
Theorem~\ref{LipschitzConverseLyapunovTheorem-1}, due to \cite[Section 3.4]{KaJ11b}.

Theorem~\ref{thm:MainResult_Characterization_ISS_EQ_Banach_Spaces} (ISS superposition theorems for semilinear evolution equations) is due to \cite[Theorem 6]{MiW18b}.


\textbf{Lyapunov characterization for ISS of linear systems with bounded input operators.}
Theorem~\ref{thm:ISS-criterion-linear-systems-bounded-operators} was stated in \cite{MiP20}. Some of its equivalences are due to \cite[Proposition 4]{MiI16}. The other ones, based on Propositions~\ref{ConverseLyapunovTheorem_LinearSystems}, \ref{p:Lyapunovfunction2}, are due to \cite[Propositions 6,7]{MiW17c}.

\textbf{Converse ISS Lyapunov theorems for linear systems.}
The results in Section~\ref{sec:Lyapunov_Theorem_systems_unbounded_operators} are due to \cite{JMP20}.
Despite these positive results, our present knowledge of Lyapunov methods for linear systems with unbounded operators is insufficient, and many questions are open.
The very first question is whether ISS of \eqref{eq:Linear_System-ISS-LF} implies the existence of a (coercive or non-coercive) ISS Lyapunov function.

In Section~\ref{sec:Transport equation with a boundary input}, we have constructed a coercive quadratic ISS Lyapunov function for the transport equation with a boundary input, and in \cite{TPT18}, it was shown that a similar method achieves 
a coercive ISS Lyapunov function also for linear hyperbolic systems of conservation laws. 
In Section~\ref{sec:Lyapunov methods for semilinear parabolic systems with boundary inputs}, we have used the coercive quadratic ISS Lyapunov function to show ISS of the (linear or nonlinear) heat equation with a Neumann boundary input.

At the same time, in spite of a lot of efforts, no coercive ISS Lyapunov functions have been proposed for the heat equation with a Dirichlet boundary control \eqref{eq:Linear-heat-equation}, even though a non-coercive one for this system exists, see Section~\ref{sec:ISS Lyapunov functions for a heat equation with Dirichlet boundary input}.

The admissibility degree of the Neumann input operator in the heat equation is better than the admissibility degree of the Dirichlet input operator for the same equation. Thus, one may conjecture the reason for the possible failure of the existence of a coercive quadratic ISS Lyapunov function is a bad admissibility degree of the input operator. 
On the other hand, till now, almost exclusively quadratic ISS Lyapunov functions have been used for the ISS analysis. Thus, possibly, one has to search for other classes of ISS Lyapunov functions. 

Also, it is not known whether there are linear ISS control systems for which there is no coercive ISS Lyapunov function, but non-coercive ISS Lyapunov functions do exist.
As ISS can be characterized as the exponential stability of an underlying semigroup combined with the admissibility of the input operator, a Lyapunov characterization of the admissibility could probably be helpful. Such a characterization exists for the case of 2-admissibility, see
\cite[Theorem 3.1]{HaW97}.



\cleardoublepage

\chapter{Integral input-to-state stability}
\label{chap:Integral input-to-state stability}

Despite all advantages of the ISS framework, for some practical systems, input-to-state stability is too restrictive.
This is because ISS excludes systems whose state stays bounded as long as the magnitude of applied inputs and of initial states remains below a specific threshold but becomes unbounded when either the input magnitude or the magnitude of an initial state exceeds the threshold. Such behavior is frequently caused by saturation and limitations in actuation and processing rate. 
The idea of integral input-to-state stability (iISS) is to capture such nonlinearities \cite{Son98, ASW00}.

\section{Basic properties of iISS systems}

In this section, we assume that $\Sigma$ is a forward complete control system in the sense of Definition~\ref{Steurungssystem} with $\Uc$ which is a linear subspace of the space $L^{1}_{\loc}(\R_+,U)$ of locally Bochner integrable $U$-valued functions.

\begin{definition}
\label{def:iISS}
\index{iISS}
\index{input-to-state stability!integral}
A forward complete system $\Sigma=(X,\Uc,\phi)$ is called \emph{integral input-to-state stable (iISS)} if there exist $\theta \in \Kinf$, $\mu \in \K$ and $\beta \in \KL$ such that for all $(t, x, u) \in \R_+ \tm X \tm \Uc$ it holds that
\begin{equation}
\label{iISS_Estimate}
\|\phi(t,x,u)\|_X \leq \beta(\|x\|_X,t) + \theta\left(\int_0^t \mu(\|u(s)\|_U)ds\right).
\end{equation}
\end{definition}

If the integral in the right-hand side of \eqref{iISS_Estimate} diverges, we consider it as equal to $+\infty$.

\begin{remark}{(iISS and forward completeness)}
\label{rem:iISS-and-Forward-completeness} 
We defined iISS only for forward complete systems.
Alternatively, one could require in the iISS definition that the estimate \eqref{iISS_Estimate} holds only for $(t, x, u) \in D_\phi \subset \R_+ \tm X \tm \Uc$.
If $\Sigma$ satisfies the BIC property (see Definition~\ref{def:BIC}), then an easy argument shows that $\Sigma$ is forward complete.
The same argument could be applied to the ISS property as well.

Finally, note that iISS is stronger than forward completeness together with 0-UGAS, even for nonlinear ODE systems, see an example in \cite[Section~V]{ASW00}.
\end{remark}

In the ISS theory of ODEs and delay systems, as well as for many other classes of evolution equations, $L^\infty$ space or the space $PC_b$ of piecewise continuous functions is a natural choice of the input space. For such choices of the input space, the ISS estimate \eqref{iss_sum} provides an upper bound on the system's response with respect to the maximal magnitude of the applied input. In contrast, \emph{integral ISS gives an upper bound of the system's response with respect to a kind of energy fed into the system, described by the integral on the right-hand side of \eqref{iISS_Estimate}.}

\begin{definition}
\label{def:BECS} 
\index{BECS}
\index{property!bounded energy-convergent state}
We say that a forward complete system $\Sigma = (X,\Uc,\phi)$ has the \emph{bounded energy-convergent state (BECS)} property if there is $\xi\in\K$ such that 
\begin{eqnarray}
\int_0^\infty \xi(\|u(s)\|_U)ds <\infty \qrq \forall r>0\ \ 
\lim_{t \to \infty}\sup_{\|x\|_X\leq r}\|\phi(t,x,u)\|_X = 0.
\label{eq:BECS} 
\end{eqnarray}
\end{definition}

Similarly to Proposition~\ref{prop:Converging_input_uniformly_converging_state}, we have the following basic properties of iISS systems, belonging to the folklore of the ISS theory:
\begin{proposition}
\label{prop:iISS-implies-bounded energy convergent state} 
Let $\Sigma$ be an iISS control system. Then $\Sigma$ is 0-UGAS and satisfies the BECS property with $\xi:=\mu$.
\end{proposition}

\begin{proof}
Let a forward complete system $\Sigma=(X,\Uc,\phi)$ be iISS with corresponding functions $\theta\in\Kinf$ and $\mu\in\K$. Pick any $u\in \Uc$ so that $\int_0^\infty \mu(\|u(s)\|_U)ds <\infty$.

To show the claim of the proposition, we need to show that for all  $\eps>0$ and $r>0$ there is a time $t^*=t^*(\eps,r)>0$ so that
\[
\|x\|_X \leq r\quad \wedge \quad  t \geq t^* \qrq \|\phi(t,x,u)\|_X \leq \eps.
\]
Pick any $r,\varepsilon>0$ and choose $t_1>0$ so that $\int_{t_1}^\infty \mu(\|u(s)\|_U)ds \leq \theta^{-1}(\frac{\eps}{2})$. 

Note that $u(\cdot +t_1)\in\Uc$ for all $t_1\geq 0$ in view of the axiom of shift-invariance.
Due to the semigroup property and integral ISS of $\Sigma$, we have for all $t\ge 0$ and all $x \in B_r$ that
\begin{eqnarray*}
\|\phi(t+t_1,x,u)\|_X  &=& \big\|\phi\big(t,\phi(t_1,x,u),u(\cdot + t_1)\big)\big\|_X \\
       &\leq& \beta\big(\|\phi(t_1,x,u)\|_X,t\big) + \theta\Big(\int_0^t \mu(\|u(s+t_1)\|_U)ds\Big) \\
       &\leq& \beta\Big(\beta(\|x\|_X,t_1) + \theta\Big(\int_0^{t_1} \mu(\|u(s)\|_U)ds\Big),t\Big)\\
			&&\qquad\qquad+ \theta\Big(\int_{t_1}^{t+t_1} \mu(\|u(s)\|_U)ds\Big)\\
       &\leq& \beta\Big(\beta(r,0) + \theta\Big(\int_0^{\infty} \mu(\|u(s)\|_U)ds\Big),t\Big) + \frac{\eps}{2}.
\end{eqnarray*}
Pick any $t_2$ in a way that 
\[
\beta\big(\beta(r,0) + \theta\Big(\int_0^{\infty} \mu(\|u(s)\|_U)ds\Big),t_2\big) \leq \frac{\eps}{2}.
\]
This ensures that for all $t\geq 0$ and all $x \in B_r$
\begin{eqnarray*}
\|\phi(t+t_2+t_1,x,u)\|_X 
&\leq& \beta\Big(\beta(r,0) + \theta\Big(\int_0^{\infty} \mu(\|u(s)\|_U)ds\Big),t+t_2\Big) + \frac{\eps}{2}\\
&\leq& \beta\Big(\beta(r,0) + \theta\Big(\int_0^{\infty} \mu(\|u(s)\|_U)ds\Big),t_2\Big) + \frac{\eps}{2}\\
&\leq& \eps.
\end{eqnarray*}
Since $\eps>0$ and $r>0$ are arbitrary, the claim of the proposition follows. 
\end{proof}

A \q{density argument} that we have shown in Lemma~\ref{lem:Density_Arg_ISS} for ISS, can also be shown for iISS:
\begin{proposition}
\label{lem:Density_Arg_iISS}
Consider a control system $\Sigma=(X,\Uc,\phi)$ and let $\phi$ depend continuously on inputs and on initial states.

Let $\hat{X}$, $\hat{\Uc}$ be dense normed vector subspaces of $X$ and $\Uc$, endowed with inherited norms from $X$ and $\Uc$ respectively, and let $\hat{\Sigma}:=(\hat{X},\hat{\Uc},\phi)$ be the system, generated by the transition map $\phi$, which is the same as in $\Sigma$ but restricted to the state space $\hat{X}$ and space of admissible inputs $\hat{\Uc}$.

Let $\hat{\Sigma}$ be iISS. Then $\Sigma$ is also iISS with the same $\beta, \theta$ and
$\mu$ in the estimate \eqref{iISS_Estimate}.
\end{proposition}
For the system $\hat{\Sigma}$ whose state space $\hat{X}$ is a normed vector space, we define the notion of iISS in the same way as it is done in Definition~\ref{def:iISS} for systems on a Banach space $X$.

\begin{proof}
Since $\hat{\Sigma}$ is iISS, there exist $\beta \in \KL$ and $\mu\in\K$, $\theta \in \Kinf$, such that $\forall \hat{x} \in \hat{X},\; \forall \hat{u} \in \hat{\Uc}$ and $\forall t\geq 0$ it holds that
\begin {equation}
\label{DichteArg_1_iISS}
\| \phi(t,\hat{x},\hat{u}) \|_{X} \leq \beta(\| \hat{x} \|_{X},t) + \theta \Big(\int_0^t \mu(\|\hat{u}(s)\|_{U})ds\Big).
\end{equation}
Let $\Sigma$ be not iISS with the same $\beta,\mu,\theta$. Then there exist $t^*>0$, $x\in X$, $u \in \Uc$:
\begin {equation}
\label{DichteArg_2_iISS}
\| \phi(t^*,x,u) \|_{X} = \beta(\| x \|_{X},t^*) + \theta \Big(\int_0^{t^*} \mu(\|u(s)\|_U)ds\Big) + r,
\end{equation}
where $r=r(t^*,x,u)>0$.
From \eqref{DichteArg_1_iISS} and \eqref{DichteArg_2_iISS}, we obtain
\begin{align}
\| \phi(t^*,x,u) \|_{X} - & \| \phi(t^*,\hat{x},\hat{u}) \|_{X} \geq
\beta(\| x \|_{X},t^*) - \beta(\| \hat{x} \|_{X},t^*)
\nonumber \\
&\quad + \theta \Big(\int_0^{t^*}\hspace{-1.2ex}\mu(\|u(s)\|_U)ds\Big) - \theta \Big(\int_0^{t^*}\hspace{-1.2ex}\mu(\|\hat{u}(s)\|_{U})ds\Big) + r.
\label{DichteArg_3_iISS}
\end{align}
Since $\hat{X}$ and $\hat{\Uc}$ are dense in $X$ and $\Uc$, respectively,
and since the map 
\[
u \mapsto \theta \Big(\int_0^{t^*} \mu(\|u(s)\|_U)ds\Big)
\]
is continuous, we can find sequences $(\hat{x}_i) \subset \hat{X}$: $\|x-\hat{x}_i \|_X \to 0$ and $(\hat{u}_i) \subset \hat{\Uc}$: $\|u-\hat{u}_i \|_{\Uc} \to 0$.
From \eqref{DichteArg_3_iISS} it follows that for each arbitrary $\eps>0$,
there exist $\hat{x}_i$ and $\hat{u}_i$ such that
\begin{equation*}
\| \phi(t^*,x,u) - \phi(t^*,\hat{x}_i,\hat{u}_i) \|_{X}  \geq
\left| \| \phi(t^*,x,u) \|_{X} - \| \phi(t^*,\hat{x}_i,\hat{u}_i) \|_{X}  \right| \\
\geq r - 2\eps.
\end{equation*}
This contradicts the assumption of continuous dependence of $\Sigma$ on initial states and inputs.
Thus, $\Sigma$ is iISS with the same $\theta$, $\beta$ and $\mu$
in \eqref{iISS_Estimate}.
\end{proof}

For linear systems with bounded input operators, the notions of ISS and integral ISS are equivalent; see Theorem~\ref{thm:ISS-criterion-linear-systems-bounded-operators}. For nonlinear systems as well as for linear systems with unbounded input operators, the situation is more complex.

\begin{definition}\label{def:iISSV}
\index{Lyapunov function!iISS}
Consider a control system $\Sigma=(X,\Uc,\phi)$ with the input space $\Uc:=PC_b(\R_+,U)$.
A continuous function $V:X \to \R_+$ is called a \emph{non-coercive iISS Lyapunov function} for $\Sigma$, if there exist
$\psi_2 \in \Kinf$, $\alpha \in \PD$ and $\sigma \in \K$ 
such that 
\begin{equation}
\label{LyapFunk_1Eig_iISS_noncoercive}
0 < V(x) \leq \psi_2(\|x\|_X), \quad \forall x \in X\backslash\{0\},
\end{equation}
and the Lie derivative of $V$ along the trajectories of the system \eqref{InfiniteDim} satisfies 
\begin{equation}
\label{eq:iISS-DissipationIneq}
\dot{V}_u(x) \leq -\alpha(\|x\|_X) + \sigma(\|u(0)\|_U),\quad x\in X,\ u\in\Uc.
\end{equation}
If additionally there is a $\psi_1\in\Kinf$ so that 
\begin{equation}
\label{LyapFunk_1Eig_iISS_coercive}
\psi_1(\|x\|_X) \leq V(x) \leq \psi_2(\|x\|_X), \quad x \in X,
\end {equation}
then $V$ is called a \emph{(coercive) iISS Lyapunov function} for $\Sigma$.
\end{definition}

Note that in the definition of an iISS Lyapunov function, the decay rate $\alpha$ is a $\PD$-function, and in the definition of an ISS Lyapunov function in a dissipative form $\alpha \in\Kinf$.

The direct Lyapunov theorem for the iISS property reads as follows (and is proved similarly to the corresponding result for ODE systems \cite[p. 1088]{ASW00}):
\begin{proposition}
\label{PropSufiISS}
Consider a control system $\Sigma=(X,\Uc,\phi)$ with $\Uc:=PC_b(\R_+,U)$, satisfying the BIC property.
If there is a coercive iISS Lyapunov function for $\Sigma$, then $\Sigma$ is iISS.
\end{proposition}

%
\begin{proof}
As $\alpha \in\PD$, by Proposition~\ref{prop:Positive-definiteness-criterion}, there are $\omega\in\Kinf$ and $\sigma\in\LL$ such that $\alpha(r)\geq \omega(r)\sigma(r)$ for all $r\geq 0$.
Pick any $\tilde{\alpha}\in\PD$ satisfying
\begin{eqnarray}
\tilde{\alpha}(r) \geq \omega(\psi_2^{-1}(r)) \cdot \sigma(\psi_1^{-1}(r)),\quad r\in\R_+.
\label{eq:tilde-rho-iISS}
\end{eqnarray} 
By \eqref{eq:iISS-DissipationIneq} and \eqref{LyapFunk_1Eig_iISS_coercive}, we have for any $x\in X$, $u\in \Uc$ that 
\begin{eqnarray}
\label{eq:iISS_decrease-iISS-direct-LF-theorem-proof}
\dot{V}_u(x) 
	&\leq& -\alpha(\|x\|_X) + \sigma(\|u(0)\|_U)  \nonumber\\
	&\leq& -\omega(\|x\|_X)\sigma(\|x\|_X) + \sigma(\|u(0)\|_U) \nonumber\\
	&\leq& - \omega\circ\psi_2^{-1}(V(x)) \cdot \sigma\circ\psi_1^{-1}(V(x)) + \sigma(\|u(0)\|_U) \nonumber\\
	&\leq& - \tilde{\alpha}(V(x)) + \sigma(\|u(0)\|_U).
\end{eqnarray}
Now take any input $u \in\Uc$, any initial state $x \in X$, and consider the corresponding solution $\phi(\cdot,x,u)$ that is defined on $[0,t_m(x,u))$.

As $V$ is continuous and $\phi$ is continuous w.r.t. $t$, the map
$y:t \mapsto V(\phi(t,x,u))$ is continuous, and for all $t \in [0,t_m(x,u))$, it holds by 
\eqref{eq:iISS_decrease-iISS-direct-LF-theorem-proof} that 
\begin{align*}
D^+ V(\phi(t,x,u)) &\leq - \tilde{\alpha}\big(V(\phi(t,x,u))\big) + \sigma(\|u(t)\|_U).
\end{align*}

By Proposition~\ref{prop:comparison-principle-with-inputs}, there is $\beta\in\KL$ (that does not depend on $x$ and $u$) such that
\begin{eqnarray*}
V(\phi(t,x,u)) \leq \beta(V(x),t) + \int_0^t 2 \sigma(\|u(s)\|_U) ds,\quad t\in [0,t_m(x,u)),
\end{eqnarray*}
and thus by \eqref{LyapFunk_1Eig_iISS_coercive}
\begin{eqnarray*}
\psi_1(\|\phi(t,x,u)\|_X) \leq \beta(\psi_2(\|x\|_X),t) + \int_0^t 2 \sigma(\|u(s)\|_U) ds,\quad t\in [0,t_m(x,u)).
\end{eqnarray*}
Using that $\psi_1^{-1}(a+b)\leq \psi_1^{-1}(2a) + \psi_1^{-1}(2b)$ for all $a,b\geq 0$, we finally obtain for
$t\in [0,t_m(x,u))$ that 
\begin{eqnarray}
\label{eq:Direct-iISS-Lyapunov-Theorem}
\|\phi(t,x,u)\|_X 
&\leq& \psi_1^{-1}\Big( \beta(\psi_2(\|x\|_X),t) + \int_0^t 2 \sigma(\|u(s)\|_U) ds\Big)\nonumber\\
&\leq& \psi_1^{-1}\Big(2 \beta\big(\psi_2(\|x\|_X),t\big)\Big) + \psi_1^{-1}\Big(4\int_0^t \sigma(\|u(s)\|_U) ds\Big).
\end{eqnarray}
As $\Sigma$ satisfies the BIC property, this estimate ensures that $t_m(x,u)=\infty$.
This shows iISS.
\end{proof}

\begin{remark}
\label{rem:Why-PC-functions-in-iISS} 
Please note that for general input spaces $\Uc$, allowed by Definition~\ref{Steurungssystem}, the term $u(0)$ in \eqref{eq:iISS-DissipationIneq} may have no sense, e.g.,\ this is the case for $\Uc=L^\infty(\R_+,\R)$.
To allow for control systems with general input spaces, the definition of an iISS Lyapunov function has to be changed. For example, one could mimic the definition of the ISS Lyapunov function and require in the definition of an iISS Lyapunov function instead of \eqref{eq:iISS-DissipationIneq} a validity of the inequality
\begin{equation}
\label{DissipationIneq-iISS-2}
\dot{V}_u(x) \leq -\alpha(\|x\|_X) + \sigma(\|u\|_\Uc).
\end{equation}
Then, using the same argument as in the proof of Proposition~\ref{PropSufiISS}, it is possible to infer from the existence of such a Lyapunov function a certain variation of the iISS property, which is, however, possibly weaker than iISS.
\end{remark}

\begin{definition}
\label{def:ISS-wrt-small-inputs}
\index{input-to-state stability!with respect to small inputs}
The system $\Sigma=(X,\Uc,\phi)$ is called \emph{input-to-state stable
(ISS) with respect to small inputs}, if there exist $\beta \in \KL$, $\gamma \in \Kinf$ and $R>0$
such that for all $ x \in X$, $ u\in \Uc$: $\|u\|_\Uc\leq R$ and $ t\geq 0$ it holds that
\begin {equation}
\label{iss-small-inputs-sum}
\| \phi(t,x,u) \|_{X} \leq \beta(\| x \|_{X},t) + \gamma( \|u\|_{\Uc}).
\end{equation}
\end{definition}

Although integral ISS systems share several nice properties with ISS systems, there are also important distinctions.
\begin{example}
\label{examp:iISS-not-strong-iISS} 
Let $X=\R$ and $\Uc= PC_b(\R_+,\R)$. Consider the system
\begin{eqnarray}
\dot{x} = -\frac{x}{1+x^2} + u.
\label{eq:iISS-not-strong-iISS}
\end{eqnarray}
This system is iISS. However, it does not have the convergent input-convergent state property (compare to Proposition~\ref{prop:Converging_input_uniformly_converging_state}), and constant inputs of arbitrarily small magnitude induce unbounded trajectories, provided that the initial state is chosen sufficiently large.
\end{example}

The obstructions demonstrated in Example~\ref{examp:iISS-not-strong-iISS} motivate the following strengthening of the iISS property, which is well-studied for ODE systems \cite{CAI14, CAI14b}:
\begin{definition}
\label{def:strong-iISS} 
\index{input-to-state stability!strong integral}
The system $\Sigma=(X,\Uc,\phi)$ is called \emph{strongly integral input-to-state stable (strongly iISS)} if $\Sigma$ is iISS and ISS with respect to small inputs.
\end{definition}

\section{Bilinear systems}
\label{sec:Bilinear_systems}

One of the simplest classes of nonlinear control systems are bilinear systems which form a bridge between the linear and the nonlinear theories and are important in numerous applications such as biochemical reactions, quantum-mechanical processes \cite{PaY08,BDK74}, reaction-diffusion-convection processes controlled by means of catalysts \cite{Kha03}, etc.

The following example indicates that most bilinear systems are not ISS:
\begin{example}
\label{examp:1dim_bilinear_system} 
A simple example of a strongly integral ISS system that is not ISS is given by a one-dimensional bilinear system
\begin{eqnarray}
\dot{x} = -x + xu,
\label{eq:SimpleBilinSys}
\end{eqnarray}
with $x(t) \in\R$ and $u\in PC_b(\R_+,\R)$.
Clearly, for $u\equiv 2$, the solution is unbounded, which is not possible for an ISS system. At the same time, \eqref{eq:SimpleBilinSys} is iISS. This can be shown by using an iISS Lyapunov function $V(x) = \ln(1+x^2)$, $x\in\R$.
Trivially, \eqref{eq:SimpleBilinSys} is also ISS with respect to small inputs, and thus \eqref{eq:SimpleBilinSys} is strongly iISS.
\end{example}

At the same time, all bilinear finite-dimensional 0-UGAS systems are iISS \cite{Son98}, and even strongly iISS \cite[Corollary 2]{CAI14}.
These results have been extended in \cite{MiI16, MiW15} for generalized bilinear distributed parameter systems with a bounded bilinear term. Here, we present these results and put them into the perspective of the strong iISS property.

Consider the special case of systems \eqref{InfiniteDim} of the form
\begin{equation}
\label{BiLinSys}
\begin{array}{l}
{\dot{x}(t)=Ax(t)+ Bu(t) + C(x(t),u(t)),} \\
x(0)=x_0,
\end{array}
\end{equation}
where $A$ generates a strongly continuous semigroup on $X$,  $B \in L(U,X)$, and $C: X \times U \to X$ satisfies the Assumption~\ref{Assumption1}, and, furthermore, 
\begin{align}
\exists \xi \in \K:   \quad \|C(x,u)\|_X \leq \|x\|_X \xi(\|u\|_U) \quad \forall x\in X,\ u\in U.
\label{eq:BilinOperator}
\end{align}
As we did for systems \eqref{InfiniteDim}, we assume that inputs belong to the space $\Uc:=PC_b(\R_+,U)$.

Next, we present a criterion for strong integral input-to-state stability of \eqref{BiLinSys}. 

\begin{proposition}
\label{ConverseLyapunovTheorem_BilinearSystems}
Let \eqref{BiLinSys} satisfy the assumption \eqref{eq:BilinOperator}.
The following statements are equivalent:
\begin{itemize}
    \item[(i)]   \eqref{BiLinSys} is strongly iISS.
    \item[(ii)]   \eqref{BiLinSys} is iISS.
    \item[(iii)]  \eqref{BiLinSys} is 0-UGAS.
    \item[(iv)]  $A$ generates an exponentially stable $C_0$-semigroup.
\end{itemize}
\end{proposition}

\begin{proof}
(i) $\Rightarrow$ (ii) $\Rightarrow$ (iii). These implications are clear from the definitions.

(iii) $\Rightarrow$ (iv). In view of \eqref{eq:BilinOperator}, $C(x,0)=0$. Hence the input-free system \eqref{BiLinSys} is linear, and the implication follows by Theorem~\ref{thm:ISS-criterion-linear-systems-bounded-operators}.

(iv) $\Rightarrow$ (i). {\bf Step 1: Integral ISS.}
Take $M,\lambda>0$ such that $\|T(t)\|\leq Me^{-\lambda t}$ for all $t\geq 0$. Pick any $u\in\Uc$, and $x_0$ and any $t$ such that $x(\cdot):=\phi(\cdot,x_0,u)$ is well-defined on $[0,t]$. The mild solution of \eqref{BiLinSys} has the form:
\begin{equation*}
\label{iISS_Syst_IntEq}
x(t)=T(t)x_0 + \int_0^t{T(t-r) \big( Bu(r)+ C(x(r),u(r))\big)dr} .
\end{equation*}
From $B \in L(U,X)$ and inequality \eqref{eq:BilinOperator} we have
\begin{align*}
\label{iISS_Syst_IntEq_Estim1}
\|x(t)\|_X & \leq  \|T(t)\|\|x_0\|_X + \int_0^t \|T(t-r)\| \big( \|B\|\|u(r)\|_U+ \|C(x(r),u(r))\|_X\big)dr \\
& \leq  Me^{-\lambda t} \|x_0\|_X + \int_0^t Me^{-\lambda (t-r)} \big( \|B\|\|u(r)\|_U + \|x(r)\|_X \xi(\|u(r)\|_U)\big)dr .
\end{align*}
We multiply both sides of the inequality by $e^{\lambda t}$ and define $z(t): = x(t) e^{\lambda t}$. From $\lambda>0$, we obtain
\begin{align*}
\|z(t)\|_X &\leq M \Big( \|z(0)\|_X + \|B\| \int_0^t e^{\lambda r} \|u(r)\|_U dr\Big)  + \int_0^t{M \|z(r)\|_X \xi(\|u(r)\|_U)dr}.
\end{align*}
Since $q: t \mapsto M\big( \|z(0)\|_X + \|B\| \int_0^t e^{\lambda r} \|u(r)\|_U dr \big)$ is a nondecreasing function, Gronwall's inequality (Lemma~\ref{lem:Gronwall}) yields
\begin{align*}
\|z(t)\|_X &\leq M\Big( \|z(0)\|_X + \|B\| \int_0^t e^{\lambda r} \|u(r)\|_U dr \Big)  e^{\int_0^t{M \xi(\|u(r)\|_U)dr}}.
\end{align*}
Returning to the original variables and using $\lambda>0$, we have
\begin{align*}
\|x(t)\|_X
&\leq M\Big( e^{-\lambda t}\|x_0\|_X + \|B\| \int_0^t e^{-\lambda (t-r)} \|u(r)\|_U dr\Big)  e^{\int_0^t{M \xi(\|u(r)\|_U)dr}} \\
&\leq M\Big( e^{-\lambda t}\|x_0\|_X + \|B\| \int_0^t \|u(r)\|_U dr\Big)  e^{\int_0^t{M \xi(\|u(r)\|_U)dr}}.
\end{align*}
Applying to both sides the function $\alpha \in \Kinf$ defined for all $r \geq 0$ as $\alpha(r)=\ln(1+r)$ results in
\begin{align*}
\alpha(\|x(t)\|_X ) &\leq \ln\Big(1+ M\big( e^{-\lambda t}\|x_0\|_X +\|B\| \int_0^t \|u(r)\|_U dr\big) e^{\int_0^t{M \xi(\|u(r)\|_U)dr}}\Big).
\end{align*}
Since for all $a,b \in \R_+$ it holds that
\[
\ln(1+ae^b) \leq \ln((1+a)e^b)=\ln(1+a) + b ,
\]
we obtain
\begin{align*}
\alpha(\|x(t)\|_X ) \leq \ln\Big( 1+ M\big( e^{-\lambda t}\|x_0\|_X  &+\|B\| \int_0^t \|u(r)\|_U dr\big) \Big) \\
 &  + \int_0^t{M \xi(\|u(r)\|_U)dr}.
\end{align*}
Moreover, for all $a,b \in \R_+$ it holds that
\[
\ln(1+a + b) \leq \ln((1+a)(1+b))=\ln(1+a) + \ln(1+b),
\]
which implies
\begin{align*}
\alpha(\|x(t)\|_X ) &\leq \ln\big(1+ M e^{-\lambda t}\|x_0\|_X\big)\\
&\qquad + \ln\Big(1+ M\|B\| \int_0^t \|u(r)\|_U dr\Big)  + \int_0^t{M \xi(\|u(r)\|_U)dr}.
\end{align*}
Since $\beta:(s,t) \mapsto \ln(1+ M e^{-\lambda t}s)$ is a $\KL$-function, with the help of
$\alpha^{-1}(a+b)\leq \alpha^{-1}(2a)+\alpha^{-1}(2b)$ holding for
any $a, b\in\R_+$, we obtain the estimate of the form 
\eqref{iISS_Estimate} on the domain of existence of $\phi(\cdot,x_0,u)$. 
This implies that $\phi(\cdot,x_0,u)$ is uniformly bounded on its domain of existence. As Assumption~\ref{Assumption1} holds, \eqref{BiLinSys} has the BIC property by Corollary~\ref{cor:BIC-property}, and hence $\phi(\cdot,x_0,u)$ is defined globally, \eqref{BiLinSys} is forward complete and iISS.

{\bf Step 2: ISS w.r.t.\  small inputs.}
Pick any $\gamma \in (0,\lambda)$. Let us show that
\[
V^\gamma(x):=\max_{r\geq 0}\|e^{\gamma r}T(r)x\|_X,\quad x \in X,
\]
is a coercive ISS Lyapunov function for \eqref{BiLinSys} subject to the input space $\clo{B_{R,\Uc}}=\{u\in\Uc:\|u\|_\Uc \leq R\}$ for $R>0$ small enough.

First note, that $V^\gamma(x)\geq \|x\|_X$ for any $x\in X$ and that $V^\gamma$ is globally Lipschitz continuous in view of Proposition~\ref{p:Lyapunovfunction2}.

To obtain an infinitesimal estimate, we compute, using the triangle inequality ($V^\gamma$ is a norm), we have the following estimate for all $u$:
    \begin{align*}
    \dot{V}^\gamma_u(x) &= \mathop{\overline{\lim}} \limits_{h \rightarrow +0} {\frac{1}{h}\big(V^\gamma(\phi(h,x,u))-V^\gamma(x)\big) } \\
    &=  \mathop{\overline{\lim}} \limits_{h \rightarrow +0}
    \frac{1}{h}\Big( V^\gamma \Big(T(h) x + \int_0^h{T(h-s) B u(s)ds} \\
		&\qquad\qquad\qquad\qquad\qquad+ \int_0^h{T(h-s) C\big(\phi(s,x,u),u(s)\big)ds}\Big) - V^\gamma(x) \Big)  \\
    &\leq  \mathop{\overline{\lim}} \limits_{h \rightarrow +0}
    \frac{1}{h}\Big( V^\gamma \big(T(h) x\big) + V^\gamma\Big(\int_0^h{T(h-s) B u(s)ds}\Big) \\
		&\qquad\qquad\qquad\qquad\qquad+ V^\gamma\Big(\int_0^hT(h-s) C\big(\phi(s,x,u),u(s)\big)\Big)- V^\gamma(x) \Big).
\end{align*}

Let $u\in\Uc$ be so that $\|u\|_\Uc\leq R$ for a certain $R>0$ that will be specified later.
It holds that:
\begin{align*}
\mathop{\overline{\lim}} \limits_{h \rightarrow +0}&\frac{1}{h} V^\gamma\Big(\int_0^hT(h-s) C\big(\phi(s,x,u),u(s)\big)ds\Big)\\
&\leq
\mathop{\overline{\lim}} \limits_{h \rightarrow +0}\frac{1}{h} \max_{r\geq 0}e^{\gamma r}\|T(r)\|\int_0^h \|T(h-s)\| \| C\big(\phi(s,x,u),u(s)\big)\|_Xds \\
&\leq
\mathop{\overline{\lim}} \limits_{h \rightarrow +0}\frac{1}{h} \max_{r\geq 0}Me^{(\gamma - \lambda) r}\int_0^h Me^{-\lambda(h-s)} \|\phi(s,x,u)\|_X \xi(\|u(s)\|_U)ds \\
&\leq \mathop{\overline{\lim}} \limits_{h \rightarrow +0}\frac{1}{h} M^2\int_0^h \|\phi(s,x,u)\|_X \xi(\|u(s)\|_U)ds \\
& =  M^2 \|\phi(0,x,u)\|_X \xi(\|u(0)\|_U) \leq M^2 \xi(R) \|x\|_X  \leq M^2 \xi(R) V^\gamma(x).
\end{align*}
Here we have used the continuity of $\phi$ with respect to time as well as the piecewise continuity of $u$.

With this estimate and arguing as in Proposition~\ref{p:Lyapunovfunction2}, we obtain that
    \begin{align*}
    \dot{V}^\gamma_u&(x) \leq  - \gamma \ V^\gamma(x) + V^\gamma (Bu(0)) + M^2 \xi(R) V^\gamma(x).
    \end{align*}
Choosing $R>0$ so that  $M^2\xi(R) <\gamma$, we see that $V^\gamma$ is an ISS Lyapunov function for
\eqref{BiLinSys} for inputs in $\clo{B_{R,\Uc}}$, and thus \eqref{BiLinSys} is ISS for inputs in $\clo{B_{R,\Uc}}$ by Theorem~\ref{LyapunovTheorem}.
%
\end{proof}

Assuming a stronger requirement on $C$, we can also show that iISS implies the existence of a non-coercive iISS Lyapunov function of a certain form. Note that right now, there are no results proving that the existence of a non-coercive iISS Lyapunov function (possibly under several further restrictions) implies the iISS of a control system. However, the following result may motivate further research in this area.

\begin{proposition}
\label{prop:iISS-noncoercive-LFs-II}
Let \eqref{BiLinSys} satisfy the assumption
\begin{align}
\exists \xi \in \K:\quad   x\in X,\ u\in U,\ t\geq 0 \qrq \|T(t)C(x,u)\|_X \leq \|T(t)x\|_X \xi(\|u\|_U).
\label{eq:BilinOperator-strengthened}
\end{align}

The following statements are equivalent:
\begin{itemize}
    \item[(i)-(iii)]   \eqref{BiLinSys} is strongly iISS/iISS/0-UGAS.
    \item[(iv)]  $A$ generates an exponentially stable semigroup.
    \item[(v)] The function $W:X \to\R_+$, defined by
    \begin{eqnarray}
W(x)=\ln \Big(1 + \int_0^{\infty} \|T(t)x\|_X^2 dt \Big),\quad x\in X,
\label{eq:LF_BiLinSys_Banach}
\end{eqnarray}
is a non-coercive iISS Lyapunov function for \eqref{BiLinSys}.
\end{itemize}
\end{proposition}

\begin{proof}
The equivalence of the assertions (i)--(iv) is shown in Proposition~\ref{ConverseLyapunovTheorem_BilinearSystems}.
 
(iv) $\Rightarrow$ (v). 
As $T$ is an exponentially stable semigroup, there exist $M,\lambda>0$ so that
$\|T(t)\|\leq Me^{-\lambda t}$, for all $t\geq 0$.
Consider the Lyapunov function for the linear system (with $C\equiv 0$), introduced in \eqref{eq:LF_LinSys_Banach-formulation}
\[
V(x):=\int_0^{\infty} \|T(t)x\|_X^2 dt,\quad x\in X.
\]
The chain rule for Dini derivatives (Lemma~\ref{Lemma2}) implies
\begin{eqnarray}
\dot{W}_u(x) = \frac{1}{1+V(x)}\dot{V}_u(x).
\label{eq:Ln_LF_ChainRule}
\end{eqnarray}
To compute $\dot{V}_u(x)$, we perform the same derivations as in
Proposition~\ref{ConverseLyapunovTheorem_LinearSystems}, and obtain the estimate \eqref{eq:I_2-estimate-LFs-linear systems}
with the difference, that instead of $Bu(t)$ we need
to consider $Bu(t)+C(x(t),u(t))$.We have
\begin{eqnarray*}
\dot{V}_u(x) 
  & = &  -  \|x\|_X^2 +\int_0^{\infty} 2\|T(t) x\|_X \big\|T(t) \big( B u(0)+C(x,u(0))\big)\big\|_X dt \\
  & \leq &  -  \|x\|_X^2 +\int_0^{\infty} 2\|T(t) x\|_X \big(\|T(t)B u(0)\|_X + \|T(t)C(x,u(0))\|_X\big) dt \\
  & \leq &  -  \|x\|_X^2 +\int_0^{\infty} 2\|T(t) x\|_X\|T(t)B u(0)\|_X dt\\
	&&\qquad\qquad\qquad +\int_0^{\infty} 2\|T(t) x\|_X \|T(t)C(x,u(0))\|_X dt.
\end{eqnarray*}
Using Young's inequality for the second term, and \eqref{eq:BilinOperator-strengthened} for the last one, we obtain for any $\eps>0$ the estimate:
\begin{eqnarray*}
\dot{V}_u(x) 
	& \leq &  -  \|x\|_X^2 + \int_0^{\infty} \eps \|T(t) x\|_X^2 + \frac{1}{\eps} \|T(t)  B u(0)\|_X^2 dt\\
	&&\qquad\qquad\qquad  + \int_0^{\infty} 2\|T(t) x\|_X \|T(t)x\|_X \xi(\|u(0)\|_U) dt \\	
	& \leq &  -  \|x\|_X^2 + \frac{\eps M^2}{2\lambda} \|x\|_X^2 + \frac{M^2}{2\lambda \eps}  \|B\|^2 \|u(0)\|_U^2
	+2V(x) \xi(\|u(0)\|_U).
\end{eqnarray*}

Due to \eqref{eq:Ln_LF_ChainRule}, we obtain
\begin{eqnarray*}
\dot{W}_u(x) &\leq & \Big(\frac{\eps M^2}{2\lambda} - 1 \Big) \frac{\|x\|_X^2}{1+V(x)}
+  \frac{M^2}{2\lambda \eps} \frac{\|B\|^2 \|u(0)\|_U^2}{1+V(x)}  \\
&&  +  \frac{2 V(x)}{1+V(x)} \xi(\|u(0)\|_U) \\
						&\leq & \Big(\frac{\eps M^2}{2\lambda} - 1 \Big) \frac{\|x\|_X^2}{1+\frac{M^2}{2\lambda} \|x\|_X^2} +  \frac{M^2}{2\lambda \eps} \|B\|^2 \|u(0)\|_U^2  +  2 \xi(\|u(0)\|_U).
\end{eqnarray*}
Considering $\eps \in (0,\tfrac{2\lambda}{M^2})$, we see that $W$ is a non-coercive iISS Lyapunov function for \eqref{BiLinSys}.

(v) $\Rightarrow$ (iv). This implication follows by Datko-Pazy Theorem (see \cite[Lemma 8.1.2]{JaZ12}, \cite[Lemma 5.1.2]{CuZ95}).
\end{proof}

\begin{remark}
Condition \eqref{eq:BilinOperator-strengthened} holds, e.g., for $C(x,u)=xBu$, where $B \in L(U,\R)$.
\end{remark}

\section{Examples}
\label{sec:Examples-iISS}

\subsection{ISS and iISS}
\label{sec:Example-MiI16}

In this concluding section, we illustrate our findings on an example of an iISS parabolic system.
Let $c>0$ and $L>0$.
Consider the following reaction-diffusion system
\begin{subequations}
\label{eq:iISS-example-nonlinear}
\begin{eqnarray}
x_t(z,t) &=& c x_{zz}(z,t) + \frac{x(z,t)}{1+|z-1|x(z,t)^2}u(z,t), \\
x(0,t) &=& x(L,t)=0,
\end{eqnarray}
\end{subequations}
on the region $(z,t) \in (0,L) \times (0,\infty)$ of
the $\R$-valued functions $x(z,t)$ and $u(z,t)$.

We choose the state space as $X:=L^2(0,L)$ and the input space as $\Uc:=PC_b(\R_+,U)$, with $U:=C(0,L)$.

It is easy to see that the above system is a generalized bilinear system since its nonlinearity satisfies \eqref{eq:BilinOperator}. This system is 0-UGAS. Therefore, it is iISS for any $L>0$. Below we explicitly construct an iISS Lyapunov function for this system. Afterward, we prove the ISS of this system for $L<1$.

Define
\[
W(x):=\int_0^L{x^2(z)dz} = \|x\|^2_{L^2(0,L)}.
\]
Take any smooth $x \in X$, any $u\in \Uc$ and denote $v:=u(\cdot,0)$.
Since $1+|z-1|x(z)^2\geq 1$, we obtain
\begin{align*}
\dot{W}_u(x) &= 2 \int_0^L{x(z)\Big( c x_{zz}(z) + \frac{x(z)}{1+|z-1|x(z)^2}v(z)\Big)}dz \\
&\leq -2c  \int_0^L \left( x_{z}(z) \right)^2 dz + 2 \int_0^L{x^2(z)|v(z)|dz}.
\end{align*}
Applying Friedrichs' inequality (Proposition~\ref{theorem:Friedrichs}) to the first term, we proceed to:
\begin{align*}
\dot{W}_u(x) \leq  -2c \left( \frac{\pi}{L} \right)^2 W(x) + 2 W(x) \|v\|_{C(0,L)}.
\end{align*}
Choosing
\begin{align}
V(x):=\ln\left(1+W(x)\right)
\label{eq:LF_1st_Sys}
\end{align}
yields
\begin{align}
\dot{V}_u(x) &\leq
-2c \!\left( \frac{\pi}{L} \right)^2\!\frac{ W(x)}{1\!+\!W(x)}
\!+2 \frac{ W(x)}{1\!+\!W(x)} \|v\|_{C(0,L)}    \nonumber\\
&\leq   -2c \!\left( \frac{\pi}{L} \right)^2\!\frac{\|x\|_{L^2(0,L)}^2}{1\!+\!\|x\|_{L^2(0,L)}^2}
+   2 \|v\|_{C(0,L)},
\label{Examp_Conclusing_iISS_estim} \\
&=
-\alpha(\|x\|_{L^2(0,L)}) + \sigma(\|v\|_{C(0,L)}) ,
\nonumber
\end{align}
where
\begin{align}
\alpha(s)=
2c \!\left( \frac{\pi}{L} \right)^2\!\frac{s^2}{1\!+\!s^2}
, \quad
\sigma(s)=2s .
\label{Examp_Conclusing_iISS_estimb}
\end{align}
Thus, Proposition \ref{PropSufiISS}
establishes iISS of \eqref{eq:iISS-example-nonlinear} irrespective of a value of $L>0$.

The above derivations hold for smooth enough functions $x$, and ISS for all $x \in L^2(0,L)$ follows due to Proposition~\ref{lem:Density_Arg_iISS}.

Interestingly, when $L<1$, the system \eqref{eq:iISS-example-nonlinear} is ISS for
the input space $\Uc:=PC_b(\R_+,U)$ with both $U=C(0,L)$ and $U=L^2(0,L)$.
To verify this, we first note that
\begin{align}
\sup_{s\in\R}
\left|\frac{s}{1+|z-1|s^2}\right|=
\frac{1}{2\sqrt{1-z}}
\end{align}
holds for $z<1$.
Assume $L<1$.  Using the same Lyapunov function $W$ we obtain for smooth $x$ and $v:=u(\cdot,0) \in C(0,L)$ that
\begin{align}
\dot{W}_u(x) &\leq 
2 \int_0^L x(z)c x_{zz}(z)dz
+ 2 \int_0^L \frac{1}{2\sqrt{1-z}}|x(z)v(z)|dz
\nonumber \\
&\leq
-2c \left( \frac{\pi}{L} \right)^2 \|x\|_{L^2(0,L)}^2
+\frac{1}{\sqrt{1-L}}\|x\|_{L^2(0,L)}\|v\|_{L^2(0,L)}
\nonumber \\
&\leq
-\left(2c \left( \frac{\pi}{L} \right)^2-w\right) \|x\|_{L^2(0,L)}^2
+\frac{1}{4(1-L)w}\|v\|_{L^2(0,L)}^2
\label{Example_Estim}
\end{align}
for $0<w<2c({\pi}/{L})^2$.  Recall that
$\|v\|_{L^2(0,L)}^2\leq L\|v\|_{C(0,L)}^2$ for $v \in C(0,L)$.
Thus, by virtue of Proposition \ref{PropSufiISS},
system \eqref{eq:iISS-example-nonlinear} is ISS whenever $L<1$. We stress that the coefficient of $\|u\|_{L^2(0,L)}^2$
in \eqref{Example_Estim} goes to $\infty$ as $L$ tends to $1$ from below. Hence, the ISS estimate \eqref{Example_Estim} is valid only if $L<1$.

For the choice of the space of input values $U=L^p(0,L)$ with $p \geq 1$, the case of
$L \geq 1$ does not allow us to have an ISS estimate like
\eqref{Example_Estim}.
In fact, if $L \geq 1$ and $U=L^p(0,L)$ for any $p \geq 1$, the right-hand side of \eqref{eq:iISS-example-nonlinear} system is undefined.

To see this, take $u:z \mapsto |z-1|^{-\tfrac{1}{2p}} \in L^p(0,L)$ and
$x:z \mapsto |z-1|^{-\tfrac{1}{2}+\tfrac{1}{2p}} \in L^2(0,L)$.
Then for any fixed $t$ we have $f_t(x,u):z \mapsto \frac{x(z,t)}{1+|z-1|x(z,t)^2}u(z,t) \notin L^2(0,L)$.
Thus, according to our formulation of \eqref{InfiniteDim}, the system \eqref{eq:iISS-example-nonlinear} is not well-defined for $U=L^p(0,L)$  for any real $p \geq 1$.

For the choice of the input space $U=C(0,L)$, we expect that the
system \eqref{eq:iISS-example-nonlinear} is not ISS for $L \geq 1$,
but we have not proved it at this time.
The blow-up of the $\sigma$-term in \eqref{Example_Estim} corresponding to
the dissipation inequality \eqref{eq:iISS-DissipationIneq} for $V=W$ suggests the
absence of ISS for the system \eqref{eq:iISS-example-nonlinear} in the case $L \geq 1$.
It is worth noticing that the iISS estimate
\eqref{Examp_Conclusing_iISS_estim} is valid for all $L>0$,
that is, for all $L>0$ we do not have a blowup of $\sigma$ in \eqref{Examp_Conclusing_iISS_estimb}, and $\alpha$ does not
become a zero function. In fact, one can recall the idea
demonstrated by Proposition~\ref{PropSufiISS} with
Definition~\ref{def:iISSV}.
An iISS Lyapunov function characterizes the absence of ISS by
only allowing the decay rate $\alpha$ in the dissipation
inequality \eqref{eq:iISS-DissipationIneq} to satisfy
$\Liminf_{s\to\infty}\alpha(s)<\lim_{s\to\infty}\sigma(s)$.
The iISS Lyapunov function yields ISS when
$\Liminf_{s\to\infty}\alpha(s)\geq \lim_{s\to\infty}\sigma(s)$, which is not the case in \eqref{Examp_Conclusing_iISS_estimb}.
Being able to uniformly characterize iISS irrespectively of whether
systems are ISS or not
should be advantageous in many applications. For instance,
ISS of subsystems is not necessary for the stability of
their interconnections, and there are examples of
UGAS interconnections involving iISS systems which are not ISS
\cite{Ito06,ItJ09,AnA07,KaJ12}.

\subsection{ISS Lyapunov functions for a class of nonlinear parabolic systems: Sobolev state space}
\label{sec:ISS_Parabolic}

The choice of an $L^p$ space both as a state space and as a space of input values is standard in ISS theory for parabolic systems \cite{MaP11,DaM13}. However, as we have seen in Section~\ref{sec:Example}, for iISS systems, such a choice is often not possible. Since our final goal is to consider interconnections of iISS and ISS systems, we need to have constructions of ISS Lyapunov functions with Sobolev state spaces. This section provides one such construction.

Consider the system 
\begin{equation}
\label{Nonlinear_Parabolic_W_12q}
x_t(z,t) = c x_{zz}(z,t) + f\big(x(z,t),x_z(z,t)\big)+ u(z,t),\quad z \in (0,L),\quad t>0,
\end{equation}
with the Dirichlet boundary conditions
\begin{eqnarray}
x(0,t) = x(L,t) =0, \quad t \geq 0.
\label{eq:BoundaryConditions_W_12q}
\end{eqnarray}
The next theorem gives a sufficient condition for ISS 
of \eqref{Nonlinear_Parabolic_W_12q} with respect to the state space $X:=W^{1,2q}_0(0,L)$, $q \in \N$, and the space $U$ of 
input values by the construction of a Lyapunov function.

\begin{theorem}\label{theorem:parabolic_ISS_W_12q_norm}
Suppose there exists a convex continuous function $\eta: \R_+\to \R$, so that 
\begin{eqnarray}
\int_0^L \Big( y_z(z) \Big)^{2q-2} y_{zz}(z)f\left(y(z),y_{z}(z)\right) dz \geq 
\int_0^L \eta\left(\left(y_{z}(z)\right)^{2q}\right) dz
\label{eq:Assumption_on_f}
\end{eqnarray}
holds for all twice continuously differentiable $y \in W^{1,2q}_0(0,L)$. Also let there exist $\epsilon >0$ so that 
\begin{eqnarray}
\hat{\alpha}(s):=
\frac{\pi^2}{q^2L^2}(c-\epsilon)s+L\eta\left(\frac{s}{L}\right)
\geq 0 \quad \forall s\in\R_+ . 
\label{eq:Assumption_on_a}
\end{eqnarray}
Then 
\begin{align}
V(y):=\int_0^L \Big( y_z(z) \Big)^{2q} dz = \|y\|_{W^{1,2q}_0(0,L)}^{2q},\quad y \in W^{1,2q}_0(0,L),
\label{eq:LF_W_12q_norm}
\end{align}
is an ISS Lyapunov function of 
\eqref{Nonlinear_Parabolic_W_12q}, \eqref{eq:BoundaryConditions_W_12q} with respect to the space $U=L^{2q}(0,L)$ of input values 
and $U=W_0^{1,2q}(0,L)\cap W^{2,2q}(0,L)$ as well. 
\end{theorem}

\begin{proof}
Assume for a while that the state $y \in X$ is smooth, and $u \in \Uc$ is such that $v:=u(\cdot,0)\in C(0,L)$. 
Along the solution of 
\eqref{Nonlinear_Parabolic_W_12q}, \eqref{eq:BoundaryConditions_W_12q}, 
the function $V$ given as in \eqref{eq:LF_W_12q_norm} satisfies 
\begin{align*}
\dot{V}_u(y) &=  \; 2q \int_0^L \Big( y_z(z) \Big)^{2q-1} y_{zt}(z) dz \\
					 &=  -2q (2q-1) \int_0^L y_{t} (z) \Big( y_z(z) \Big)^{2q-2} y_{zz}(z) dz +2q \Big( y_z(z) \Big)^{2q-1} y_{t}(z) \Big|_{z=0}^L.
\end{align*}
Due to \eqref{eq:BoundaryConditions_W_12q}, we have $y_{t}(z) \big|_{z=0}^L=0$ and consequently 					
\begin{align*}
\dot{V}_u(y) = -2q (2q-1) \int_0^L \Big( y_z(z) \Big)^{2q-2} y_{zz}(z) 
				\cdot \Big(c y_{zz}(z) + f(y(z),y_{z}(z))+ v(z) \Big) dz.				
\end{align*}
Next we utilize \eqref{eq:Assumption_on_f} to obtain
\begin{align}
\tfrac{1}{2q (2q-1)} \dot{V}_u(y) 
&\leq -c \int_0^L \Big( y_z(z) \Big)^{2q-2} \Big(y_{zz}(z)\Big)^2 dz 
-\int_0^L \eta\left(\left(y_{z}(z)\right)^{2q}\right) dz
\nonumber\\
				    & \quad 
- \int_0^L \Big( y_z(z) \Big)^{2q-2} y_{zz}(z) v(z)  dz.
\label{BasicEstimate_For_Other_U}
\end{align}
For any $\omega>0$, using Young's inequality $y_{zz}(z) v(z) \leq \frac{\omega}{2} \Big(y_{zz}(z)\Big)^2 + \tfrac{1}{2\omega}v^2(z)$, we get
\begin{align}
\tfrac{1}{2q (2q-1)} \dot{V}_u(y) &\leq \left(\frac{\omega}{2}-c\right) \int_0^L \Big( y_z(z) \Big)^{2q-2} \Big(y_{zz}(z)\Big)^2 dz 
-\int_0^L \eta\left(\left(y_{z}(z)\right)^{2q}\right) dz
\nonumber\\
				    & \quad 
+ \frac{1}{2\omega}\int_0^L \Big( y_z(z) \Big)^{2q-2} v^2(z) dz.
\label{eq:W_12q_Temporary_Estimate}
\end{align}
It is easy to see that
\[
\int_0^L \Big( y_z(z) \Big)^{2q-2} \Big(y_{zz}(z)\Big)^2 dz
= \frac{1}{q^2} \int_0^L \Big( \frac{\partial}{\partial z} \Big(\Big( y_z(z) \Big)^{q} \Big)\Big)^2 dz.
\]
Due to Friedrichs' inequality  \eqref{ineq:Friedrichs} we proceed to
\begin{equation}
\int_0^L \Big( y_z(z) \Big)^{2q-2} \Big(y_{zz}(z)\Big)^2 dz
\geq \frac{\pi^2}{q^2L^2} \int_0^L \Big( y_z(z) \Big)^{2q}dz  = \frac{\pi^2}{q^2L^2}V(y).
\label{eq:LinTermEstimate}
\end{equation}
Define 
\begin{eqnarray}
\xi(s):=\frac{1}{L}\hat{\alpha}(Ls) = \frac{\pi^2}{q^2L^2}(c-\epsilon)s+\eta\left(s\right)
, \quad  \forall s\in\R_+ . 
\label{eq:Assumption_xi}
\end{eqnarray}
Here we assume that $\epsilon \in (0,c)$.
The convexity of $\eta$ implies the convexity of $\xi$. 
Due to the definition \eqref{eq:LF_W_12q_norm} of $V$, 
Jensen's inequality \eqref{ineq:Jensen} yields 
\begin{eqnarray}
\int_0^L\xi\left(\left(y_{z}(z)\right)^{2q}\right) dz \geq \hat{\alpha}(V(y)).
\label{eq:xihatal}
\end{eqnarray}
We continue the estimates of $\dot{V}$ in the 
cases of $q=1$ and $q>1$ separately.

In the case of $q=1$, inequality \eqref{eq:W_12q_Temporary_Estimate} implies for $\omega < 2c$
\begin{align}
\dot{V}_u(y) &\leq 
2\Big(\frac{\omega}{2}-c\Big)\frac{\pi^2}{L^2}V(y) 
-2\int_0^L \eta\left(\left|y_{z}(z)\right|^{\!2}\right)\!dz
+ \frac{1}{\omega} \|v\|^{2}_{L^{2}(0,L)}
\nonumber \\
& = 
2\Big(\frac{\omega}{2}-c\Big)\frac{\pi^2}{L^2}V(y) 
-2\int_0^L \xi\left(\left|y_{z}(z)\right|^{\!2}\right)\!dz
\nonumber \\
& \qquad \qquad \qquad \qquad
+2\left(c-\epsilon\right)\frac{\pi^2}{L^2} \int_0^L \left|y_z(z)\right|^{\!2} dz
+ \frac{1}{\omega} \|v\|^{2}_{L^{2}(0,L)}
\nonumber \\
& \leq 
2\Big(\frac{\omega}{2}{-}\epsilon\Big)\frac{\pi^2}{L^2}V(y) 
-2\hat{\alpha}(V(y))
+ \frac{1}{\omega} \|v\|^{2}_{L^{2}(0,L)}.
\label{eq:W_12q_Final_q_equals_2}
\end{align}
Here, the last inequality uses \eqref{eq:xihatal}. 
Recall that $\epsilon>0$ and pick any $\omega\in(0,2\epsilon)$. The above inequality shows the needed estimate for $\dot{V}(y)$ provided $y$ and $u$ are smooth enough.

Next, assume that $q>1$. We apply Young's inequality \eqref{ineq:Young} to the last term in \eqref{eq:W_12q_Temporary_Estimate} with any $\omega_2>0$ as follows
\begin{equation*}
\Big( y_z(z) \Big)^{2q-2} v^2(z) 
\leq
\frac{1}{q\omega_2} v^{2q}(z) + \omega^{\frac{1}{q-1}}_2\frac{q-1}{q} 
\Big( y_z(z) \Big)^{2q}.
\end{equation*}
Putting this expression into \eqref{eq:W_12q_Temporary_Estimate} and using \eqref{eq:Assumption_xi}, we obtain finally
\begin{align}
\tfrac{1}{2q (2q-1)} \dot{V}_u(y) 
&\leq \Big( \Big(\frac{\omega}{2}-c\Big)\frac{\pi^2}{q^2L^2} +\omega^{\frac{1}{q-1}}_2 \frac{q-1}{2\omega q}\Big)V(y) 
\nonumber \\
& \qquad \qquad 
-\int_0^L \hspace{-1ex}\eta\left(\left(y_{z}(z)\right)^{2q}\right) dz
+ \frac{1}{2\omega \omega_2 q} \|v\|^{2q}_{L^{2q}(0,L)} 
\nonumber \\
&\leq \Big( \Big(\frac{\omega}{2}-\epsilon\Big)\frac{\pi^2}{q^2L^2} +\omega^{\frac{1}{q-1}}_2 \frac{q-1}{2\omega q}\Big)V(y) 
\nonumber \\
& \qquad \qquad  - \hat{\alpha}(V(y)) + \frac{1}{2\omega \omega_2 q} \|v\|^{2q}_{L^{2q}(0,L)} .
\label{eq:W_12q_Final_Estimate}
\end{align}
It is easy to see that by choosing $\omega>0$ and $\omega_2>0$ small enough, we can ensure  
\[
\Big(\frac{\omega}{2}-\epsilon\Big)\frac{\pi^2}{q^2 L^2}+\omega^{\frac{1}{q-1}}_2 \frac{q-1}{2\omega q}<0.
\] 
In the above proof, we have used several times integration by parts as well as partial derivatives of $y$ and $u$. Thus, the derivations are justified if the functions are smooth enough. However, the right-hand side in the inequalities \eqref{eq:W_12q_Final_q_equals_2} and \eqref{eq:W_12q_Final_Estimate} is continuous in the norms in $X$ and $U$. The left-hand side $\dot{V}_u(y)$ is also continuous in the norm of $X$.

Therefore, \eqref{eq:W_12q_Final_q_equals_2} and \eqref{eq:W_12q_Final_Estimate} hold for all $y \in X$ and all $u \in U$. 
Due to \eqref{eq:Assumption_on_a}, 
the function $V$ is an ISS Lyapunov function of \eqref{Nonlinear_Parabolic_W_12q} with respect to the input space $U=L^{2q}(0,L)$, $q \geq 1$. 
Application of Friedrichs' inequality to $\|u\|^{2q}_{L^{2q}(0,L)}$ 
proves that $V$ is an ISS Lyapunov function
with respect to $U=W_0^{1,2q}(0,L)\cap W^{2,2q}(0,L)$. 
\end{proof}

According to \eqref{DissipationIneq}, iISS allows the decay rate of $V$ to be much slower 
for large magnitudes of state variables (since $\alpha \in \PD$ can be bounded) than ISS can allow. 
This indicates that significantly different constructions for iISS 
Lyapunov functions are needed. The next section is devoted to this question.


\subsection{{\rm i}ISS of a class of nonlinear parabolic systems: $L^p$ state space}
\label{sec:parablic}

Consider the system 
\begin{equation}
\label{GekoppelteNonLinSyst_ODE-PDE}
x_t(z,t) = c x_{zz}(z,t) + f(x(z,t),u(z,t)) ,\quad t>0,
\end{equation}
defined on the spatial domain $(0,L)$ with 
\begin{eqnarray}
x(0,t) x_z(0,t) = x(L,t) x_z(L,t)=0, 
\quad t \geq 0 
\label{eq:BoundaryConditions_1eq}
\end{eqnarray}
which represents boundary conditions of Dirichlet, Neumann or mixed type. 
The state space for \eqref{GekoppelteNonLinSyst_ODE-PDE} we choose as $X:=L^{2q}(0,L)$ for some $q \in \N$ and input space we take as $U:=L^{\infty}(0,L)$ and $H^1_0(0,L)$.

Define the following ODE associated with 
\eqref{GekoppelteNonLinSyst_ODE-PDE} given by
\begin{eqnarray}
\dot{w}(t)=f(w(t),r(t)) , \quad w(t),\ r(t)\in\R.
\label{eq:ODEsys}
\end{eqnarray}
The next theorem provides the construction of an iISS Lyapunov function for a class of nonlinear systems of the form \eqref{GekoppelteNonLinSyst_ODE-PDE}.
\begin{theorem}
\label{theorem:parabplic_iISS}
Suppose that $W: w \mapsto w^{2q}$ satisfies for some $\alpha\in\Kinf\cup\{0\}$, $\sigma \in \K$ and all $w \in\R$ and $r\in\R$ that
\begin{eqnarray}
\label{eq:LF_for_ODEsysi}
\dot{W}_r(w): = 2q w^{2q-1} f(w,r) \leq - \alpha(W(w)) + W(w)\sigma(|r|).
\end{eqnarray}
Let any of the following conditions hold:
\begin{enumerate}[label=(\roman*)]
\item\label{item:iISS_boundary} 
$x(0,t)=0$ for all $t\geq 0$ or $x(L,t)=0$ for all $t\geq 0$.
\item\label{item:iISS_alpha} 
$\alpha$ is convex and $\Kinf$. 
\end{enumerate}
Let $X:=L^{2q}(0,L)$. Then an iISS Lyapunov function of \eqref{GekoppelteNonLinSyst_ODE-PDE} with 
\eqref{eq:BoundaryConditions_1eq} with respect to the input space $\Uc:=PC_b(\R_+,U)$ with either $U=L^{\infty}(0,L)$ or $U=H^1_0(0,L)$, is given by
\begin{align}
V(x)=\ln(1+Z(x)) , 
\label{eq:LF_for_PDEi}
\end{align}
where $Z$ is defined as 
\begin{align}
Z(x)=\int_0^L W(x(z)) dz = \|x\|_{L^{2q}(0,L)}^{2q}.
\label{eq:LF_for_PDE}
\end{align}
Furthermore, if $\alpha$ is convex and satisfies 
\begin{align}
\label{eq:parabolic_alporder}
\Liminf_{s\to\infty}\frac{\alpha(s)}{s}=\infty, 
\end{align}
then $V$ given above is 
an ISS Lyapunov system of \eqref{GekoppelteNonLinSyst_ODE-PDE} with 
\eqref{eq:BoundaryConditions_1eq} with respect to 
$U=L^{\infty}(0,L)$ as well as $U=H^1_0(0,L)$.
\end{theorem}
\begin{proof}
Consider $Z$ given by \eqref{eq:LF_for_PDE} and let $U=L^{\infty}(0,L)$.
Pick any smooth $y \in X$ and any $u\in\Uc$ with $v:=u(\cdot,0) \in C(0,L)$.
Using \eqref{eq:LF_for_ODEsysi}, we have
\begin{align}
\dot{Z}_u(y) &= 2q \int_0^L y^{2q-1}(z) 
\cdot \left(c y_{zz}(z) + f(y(z),v(z)) \right)dz 
\nonumber \\
&\leq -2q (2q-1)c \int_0^L y^{2q-2}(z) \left(y_{z}(z) \right)^2 dz 
\nonumber \\
& \quad\quad\quad\quad\quad+ \!\int_0^L \hspace{-1.2ex}\!\Big(-\alpha\big(W(y(z))\big) \!+\! W(y(z)) \sigma(|v(z)|)\!\Big) dz  
\nonumber \\
&\leq - \frac{2(2q-1)c}{q} \int_0^L \left(\frac{\partial}{\partial z}(y^q)(z) \right)^2 dz 
-\int_0^L \hspace{-.2ex}\alpha\big(W(y(z))\big)dz 
 + Z(y) \sigma(\|v\|_{U}).
\label{eq:parabolic_iISS_V}
\end{align}
In the last estimate, we have used boundary conditions \eqref{eq:BoundaryConditions_1eq}. 

First, suppose that item \ref{item:iISS_boundary} holds. 
Since $y \in W^{2q,1}(0,L)$, we have $y^q \in L^2(0,L)$ and $\frac{d}{dz}(y^q)=qy^{q-1}\frac{dy}{dz} \in L^2(0,L)$ due to H\"older's inequality (as $\frac{dy}{dz} \in L^{2q}(0,L)$). 
Overall, we have $y^q \in W^{2,1}(0,L)$. Applying Poincare's inequality \eqref{Wirtinger_Variation_Ineq} to the first term in 
\eqref{eq:parabolic_iISS_V}, we obtain 
\begin{align*}
\dot{Z}_u(y) \leq - \frac{2(2q-1)c}{q} \frac{\pi^2}{4L^2}Z(y) + Z(y) \sigma(\|v\|_{U}) 
\end{align*}
with the help of $y(0)=0$ or $y(L)=0$. Defining $V$ as in \eqref{eq:LF_for_PDEi} results in 
\begin{align}
\dot{V}_u(y) &\leq - \frac{2(2q-1)c}{q} \frac{\pi^2}{4L^2} \frac{\|y\|_X^{2q}}{1+\|y\|_X^{2q}} + \sigma(\|v\|_{U}) . 
\label{Zderivative}
\end{align}
Recall that $X=L^{2q}(0,L)$.
Hence, $V$ is 
an iISS Lyapunov function of \eqref{GekoppelteNonLinSyst_ODE-PDE} 
with boundary conditions \eqref{eq:BoundaryConditions_1eq} 
for the space $U=L^\infty(0,\infty)$. 

Next, assume that item \ref{item:iISS_alpha} is satisfied. 
Due to the convexity of $\alpha$, Jensen's inequality in 
\eqref{eq:parabolic_iISS_V} allows us to obtain 
\begin{align*}
\dot{Z}_u(y) \leq 
-L\alpha\left(\frac{1}{L}Z(y(z))\right)
+   Z(y) \sigma(\|v\|_{U}) .
\end{align*}
For $V$ in \eqref{eq:LF_for_PDEi} we have 
\begin{align}
\dot{V}_u(y) \leq
-\frac{L\alpha\left(\frac{1}{L}\|y\|_X^{2q}\right)}{1+\|y\|_X^{2q}}
+ \sigma(\|v\|_{U}) .
\label{Zderivative2}
\end{align}
Thus, the function $V$ is 
an iISS Lyapunov function of \eqref{GekoppelteNonLinSyst_ODE-PDE} 
with boundary conditions \eqref{eq:BoundaryConditions_1eq} 
for the space of input values $U=L^\infty(0,L)$. 
Since \eqref{eq:parabolic_alporder} implies 
\[
\Liminf_{s\to\infty} \frac{\alpha(s)}{1+Ls}=\infty,
\] 
the above inequality guarantees that $V$ is 
an ISS Lyapunov function in the case of \eqref{eq:parabolic_alporder}. 

Finally, 
to deal with the space $H^1_{0}(0,\pi)$ for the input values, 
we recall Agmon's inequality \eqref{ineq:Agmon}, which implies for $v \in H^1_0(0,\pi)$
\begin{align*}
\|v\|^2_{L^{\infty}(0,L)} \leq 
\|v\|^2_{L^{2}(0,L)} + 
\|v_z\|^2_{L^{2}(0,L)}.
\end{align*}
This inequality yields 
\begin{align}
\|v\|^2_{L^{\infty}(0,L)} \leq 
\left(\frac{L^2}{\pi^2}+1\right) \|v\|^2_{H^1_{0}(0,L)}
\label{eq:agmonfried0}
\end{align}
with the help of Friedrichs' inequality \eqref{ineq:Friedrichs}. 
Substitution of \eqref{eq:agmonfried0} into 
\eqref{Zderivative} and \eqref{Zderivative2}, proves that $V$ is 
an iISS Lyapunov function of 
\eqref{GekoppelteNonLinSyst_ODE-PDE}-\eqref{eq:BoundaryConditions_1eq} 
with respect to $U=H^1_{0}(0,\pi)$ under either item \ref{item:iISS_boundary} or item \ref{item:iISS_alpha}.
\end{proof}

\begin{remark}
We want to draw the reader's attention to the choice of the input space. First, we have proved the iISS of the system \eqref{eq:bilinear1} for the input space $L^{\infty}(0,L)$. For many applications, this choice of input space is reasonable and sufficient. However, when considering interconnections of control systems, the input to one system is the state of another system. Thus, having $L^{\infty}(0,L)$ as an input space of the first subsystem automatically means that it is the state space of another subsystem, which complicates the proof of its ISS since the constructions of Lyapunov functions for this choice of state space are hard to find (e.g., how to differentiate such Lyapunov functions?), if possible. As we have seen in Section~\ref{sec:ISS_Parabolic}, this is not the case if we choose $H^1_0(0,L)$ as the state space. This underlines the role of Agmon's inequality in our constructions, which made possible the transition from the space $L^{\infty}(0,L)$ to $H^1_0(0,L)$ in the previous theorem. 
\end{remark}

Note that the term $W(w)\sigma(|r|)$ in \eqref{eq:LF_for_ODEsysi} allows analyzing PDEs \eqref{GekoppelteNonLinSyst_ODE-PDE} with bilinear or generalized bilinear terms which do not possess the ISS property. 


\subsection{Integral ISS of a class of nonlinear parabolic systems: Sobolev state space}
\label{sec:parablicsobolev}

Instead of the $L^p$ state space that we used for characterizing iISS
in Section \ref{sec:parablic}, for a class of parabolic systems, 
this section demonstrates that iISS can be established 
with Sobolev state spaces. 
We consider 
\begin{align}
x_t(z,t) = c x_{zz}(z,t) + f\big(x(z,t),x_z(z,t)\big)
+ x_z(z,t)u(z,t)
\label{eq:bilinear1}
\end{align}
defined for $(z,t) \in (0,L) \times (0,\infty)$ with the Dirichlet boundary conditions
\begin{eqnarray}
x(0,t) = x(L,t) =0, \quad \forall t \geq 0 .
\label{eq:bilinear1boundary}
\end{eqnarray}
We take $X:=W^{1,2q}_0(0,L)$, $q \in \N$. 
Modifying Theorem \ref{theorem:parabolic_ISS_W_12q_norm}, 
one can verify the following. 

\begin{theorem}\label{theorem:bilinear1}
Suppose that 
\eqref{eq:Assumption_on_f} holds for all $y \in X$ with 
some convex continuous function $\eta: \R_+\to \R$ 
and some $\epsilon\in\R_+$ such that 
\eqref{eq:Assumption_on_a} holds. 
If $\epsilon>0$, then 
the function $V$ given by 
\begin{align}
&
V(y)=\ln(1+Z(y)) , 
\label{eq:bilinear1V}
\\
&
Z(y)=\int_0^L \Big( y_z(z) \Big)^{2q} dz = \|y\|_{W^{1,2q}_0(0,L)}^{2q} 
\label{eq:bilinear1Z}
\end{align}
is an iISS Lyapunov function of 
\eqref{eq:bilinear1}-\eqref{eq:bilinear1boundary} with respect to the space $U=L^\infty(0,L)$ of input values and $U=H^1_0(0,L)$ as well. 
\end{theorem}
\begin{proof}
As in the proof of Theorem \ref{theorem:parabolic_ISS_W_12q_norm}, 
in the case $q=1$, along solutions of 
\eqref{eq:bilinear1}-\eqref{eq:bilinear1boundary},  
the function $Z$ given by \eqref{eq:bilinear1Z} satisfies for $\hat{\alpha}$ given by \eqref{eq:Assumption_on_a} (and depending on $\eta$ from \eqref{eq:Assumption_on_f}) the following estimates:
\begin{align}
\dot{Z}_u(y) 
&\leq 
2\Big(\frac{\omega}{2}-\epsilon\Big)\frac{\pi^2}{L^2}Z(y) 
-2\hat{\alpha}(Z(y))
+ \frac{1}{\omega}\int_0^L \Big( y_z(z) \Big)^{\!2} v^2(z) dz
\nonumber \\
&\leq  
2\Big(\frac{\omega}{2}-\epsilon\Big)\frac{\pi^2}{L^2}Z(y) 
-2\hat{\alpha}(Z(y))
+ \frac{Z(y)}{\omega} \|v\|^{2}_{L^\infty(0,L)}
\label{eq:bilinear1vdot1}
\end{align}
for any $\omega>0$. 
Due to \eqref{eq:bilinear1V}, we have 
\begin{align}
\dot{V}_u(y) 
&\leq  
2\Big(\frac{\omega}{2}-\epsilon\Big)\frac{\pi^2}{L^2}
\frac{\|y\|_X^{2}}{1+\|y\|_X^{2}}
-\frac{2\hat{\alpha}(\|y\|_X^2)}{1+\|y\|_X^{2}}
+ \frac{1}{\omega} \|v\|^{2}_{L^\infty(0,L)}.
\label{eq:bilinear1vdot2}
\end{align}
Pick $\omega\in(0,2\epsilon)$. Then $\epsilon >0$ and property 
\eqref{eq:bilinear1vdot2} imply that 
$V$ is an iISS Lyapunov function with respect to the space 
$  
U=L^\infty(0,L)$ of input values.

Next, consider $q>1$. Again, following the argument used in 
the proof of Theorem \ref{theorem:parabolic_ISS_W_12q_norm}, 
we obtain
\begin{align}
\tfrac{1}{2q (2q-1)} \dot{Z}_u(y) 
&\leq 
\Big( \big(\frac{\omega}{2}-\epsilon\big)\frac{\pi^2}{q^2L^2} +\omega^{\frac{1}{q-1}}_2 \frac{q-1}{2\omega q}\Big)Z(y) 
\nonumber \\
& \quad\qquad \qquad  - \hat{\alpha}(Z(y)) + \frac{Z(y)}{2\omega \omega_2 q} 
\|v\|^{2q}_{L^\infty(0,L)} .
\label{eq:bilinear1vdot3}
\end{align}
for any $\omega, \omega_2>0$. 
From  \eqref{eq:bilinear1V} it follows that 
\begin{align}
\tfrac{1}{2q (2q-1)} \dot{V}_u(y)
&\leq 
\Big( \Big(\frac{\omega}{2}-\epsilon\Big)\frac{\pi^2}{q^2L^2} +\omega^{\frac{1}{q-1}}_2 \frac{q-1}{2\omega q}\Big)
\frac{\|y\|_X^{2q}}{1+\|y\|_X^{2q}}
\nonumber \\
& \quad\qquad \qquad  - 
\frac{\hat{\alpha}(\|y\|_X^{2q})}{1+\|y\|_X^{2q}}
+ \frac{1}{2\omega \omega_2 q} 
\|v\|^{2q}_{L^\infty(0,L)}.
\label{eq:bilinear1vdot4}
\end{align}
This inequality with sufficiently small $\omega, \omega_2>0$ 
implies that 
$V$ is an iISS Lyapunov function with respect to the space 
$U=L^\infty(0,L)$ of input values. 

To deal with the space $H^1_{0}(0,\pi)$ for the input values, 
substitute \eqref{eq:agmonfried0} into 
\eqref{eq:bilinear1vdot2} and \eqref{eq:bilinear1vdot4}. 
\end{proof}


\section[Integral ISS of linear systems]{Integral ISS of linear systems with unbounded input operators}
\label{sec:iISS of linear systems}

The ISS estimate is defined in terms of the norms in $X$ and $\Uc$, which is one of the reasons for the rather elegant characterization of
ISS in terms of the exponential stability of a semigroup and admissibility of the input operator.
In contrast to this, integral ISS is defined in terms of the integration of an input function with a nonlinear scaling, see \eqref{iISS_Estimate}, which does not necessarily produce a norm. This makes the characterization of iISS more involved.

Recall that for linear systems with bounded input operators, ISS and iISS are equivalent notions in view of Theorem~\ref{thm:ISS-criterion-linear-systems-bounded-operators}. This remains true for linear systems with unbounded operators if $\Uc$ is an $L^p$-space with $p<\infty$.

\begin{proposition}
\label{prop:iISS_equals_ISS_p_less_infty}
Let $X$ and $U$ be Banach spaces.
If the system \eqref{eq:Linear_System} is $L^p$-ISS for some $p\in[1,\infty)$, then \eqref{eq:Linear_System} is $L^q$-ISS for all $q \in [p,\infty]$ and iISS.
\end{proposition}

\begin{proof}
As \eqref{eq:Linear_System} is $L^p$-ISS, $B$ is $p$-admissible and $T$ is exponentially stable.
In view of Remark~\ref{rem:Relations-between-admissibility-classes}, $B$ is $q$-admissible for all $q \in[p,\infty]$, which implies $L^q$-ISS for all such $q$.

Furthermore, for $\Uc:=L^p(\R_+,U)$ with $p\in[1,\infty)$ the ISS estimate \eqref{iss_sum} is automatically an iISS estimate with
$\theta(r) = \gamma(r^{1/p})$ and $\mu(r) = r^p$, $r\geq 0$.
\end{proof}

The case $p=\infty$ needs special care. First of all, integral ISS implies ISS:
\begin{proposition}
\label{prop:iISS_implies_ISS}
Let $X$ and $U$ be Banach spaces and let $\Uc:=L^\infty(\R_+,U)$.
If \eqref{eq:Linear_System} is iISS, then $T$ is an exponentially stable semigroup, and $B$ is an infinite-time $\infty$-admissible operator.
If additionally $\phi$ is continuous w.r.t.\ time in the norm of $X$, then \eqref{eq:Linear_System} is $L^\infty$-ISS.
\end{proposition}

\begin{proof}
Integral ISS estimate for $u\equiv 0$ implies exponential stability of the semigroup $T$.
Furthermore, integral ISS implies well-posedness of the system \eqref{eq:Linear_System} for $\Uc:=L^\infty(\R_+,U)$, which is equivalent to $\infty$-admissibility of $B$.
According to Lemma~\ref{lem:GrabowskiLemma}, $B$ is infinite-time $\infty$-admissible.
As $\phi$ is continuous w.r.t.\ time in the norm of $X$, Theorem~\ref{thm:ISS-Criterion-lin-sys-with-unbounded-operators}
shows $L^\infty$-ISS of \eqref{eq:Linear_System}.
\end{proof}

\begin{remark}
\label{rem:Terminology_of_JNP18_JSZ19}
In this work, continuity of (mild) solutions is a part of the definition of the system, which is in accordance with the well-posedness concepts for linear infinite-dimensional control systems in Section~\ref{sec:Abstract linear control systems}.
However, in the terminology of \cite{JNP18, JSZ19}, the continuity of solutions is not a part of the definition of ISS. Hence, in the terminology of \cite{JNP18, JSZ19},
$L^\infty$-ISS of \eqref{eq:Linear_System} is equivalent to the properties that $T$ is an exponentially stable semigroup and $B$ is infinite-time $\infty$-admissible operator, and Proposition~\ref{prop:iISS_implies_ISS} (now without an assumption of the continuity of $\phi$) reads as follows:
\begin{center}
\eqref{eq:Linear_System} is iISS $\qrq$ \eqref{eq:Linear_System} is $L^\infty$-ISS.
\end{center}
\end{remark}

In view of Proposition~\ref{prop:iISS_implies_ISS}, we know that integral ISS implies $L^\infty$-ISS. It is an open problem whether the converse implication holds. 

By definition of a coercive ISS Lyapunov function for general infinite-dimensional systems, it follows that any such function is automatically a
coercive iISS Lyapunov function, and thus a system possessing a coercive ISS Lyapunov function is both ISS and iISS.
Hence, if there are ISS control systems that are not iISS (i.e., if the answer to the above problem is negative), then these systems do not possess a coercive ISS Lyapunov function.

A partial positive answer for this open problem was obtained in \cite{JSZ19} for parabolic systems on Hilbert spaces employing holomorphic functional calculus.
\begin{proposition}{\cite[Theorem 2]{JSZ19}}
\label{prop:iISS-equals-ISS-analytic-case}
Assume that $A$ generates an exponentially stable, analytic semigroup on a Hilbert
space $X$, which is similar to a contraction semigroup.
Let $U$ be so that $\dim (U) <\infty$, $\Uc:=L^\infty(\R_+,U)$.

The following statements are equivalent:
\begin{itemize}
    \item[(i)]   $B\in L(U,X_{-1})$.
    \item[(ii)]  \eqref{eq:Linear_System} is $L^\infty$-ISS.
    \item[(iii)] \eqref{eq:Linear_System} is iISS.
\end{itemize}
\end{proposition}

\begin{proof}
(iii) $\Rightarrow$ (ii). Integral ISS of \eqref{eq:Linear_System} w.r.t.\ $\Uc$ implies that $B$ is an $\infty$-admissible operator, and in particular, $B\in L(U,X_{-1})$.
By Theorem~\ref{thm:Contiuity-of-a-map}, the triple $\Sigma:=(X,\Uc,\phi)$ is a control system, in particular, $\phi$ is continuous w.r.t.\ time in the $X$-norm. Now Proposition~\ref{prop:iISS_implies_ISS} shows $L^\infty$-ISS of \eqref{eq:Linear_System}.

(ii) $\Rightarrow$ (i). Clear.

(i) $\Rightarrow$ (iii). See \cite[Theorem 2]{JSZ19}.
\end{proof}

A graphical representation of the relationship between ISS and iISS of linear systems is depicted in Table~\ref{table:overview}.

\begin{remark}
\label{rem:iISS-and-Orlicz}
A notable result \cite[Theorem 3.1]{JNP18} shows that iISS of \eqref{eq:Linear_System} is equivalent to the existence of an Orlicz space $W$ so that \eqref{eq:Linear_System} is ISS with $\Uc:=W$. 
\end{remark}

\begin{table}
\center
\def\arraystretch{1.5}%
\begin{tabular}{c|c|c|c}
& \parbox[t]{2.7cm}{\centering Eq.~\eqref{eq:Linear_System},\\
$B$ bounded}&\parbox[t]{2.7cm}{\centering Eq.~\eqref{eq:Linear_System}, \\
$B$ unbounded\vspace{0.2cm} } &\parbox[t]{3.7cm}{\centering Eq.~\eqref{InfiniteDim},\\ $f$ Lipschitz in $(x,u)$}\\
\hline
$\mathrm{dim}\, X<\infty$ &\multicolumn{2}{|c|}{ISS $\iff$ iISS}&ISS $\Longrightarrow \atop{\centernot\Longleftarrow}$ iISS\\
\hline
$\mathrm{dim}\, X=\infty$&ISS $\iff$ iISS& ISS $\impliedby\atop \left(\stackrel{?}{\Longrightarrow}\right)$ iISS& not known \\
\hline
\end{tabular}
\medskip

\caption{(Taken from \cite{JNP18}). Relations between ISS and iISS with respect to $L^{\infty}$ and under the assumption that $\phi$ is continuous in the $X$-norm, in various settings.}
\label{table:overview}
\end{table}

\section{Diagonal systems}
\label{sec:diagonal_systems}

Now we study ISS and iISS of a well-studied (see, e.g., \cite[Sections 2.6, 5.3]{TuW09}) class of linear parabolic diagonal systems.
%
%

For $q\in[1,\infty)$ denote by $\ell_q$ the Banach space of sequences $(a_k)_{k}$, $a_k\in\C$, so that $\sum_{k=1}^\infty |a_k|^q<\infty$.
\begin{definition}
\label{def:q-Riesz-basis}
\index{$q$-Riesz basis}
Let $X$ be a separable Hilbert space.
We say that a sequence of vectors $(\phi_k)$ is a \emph{$q$-Riesz basis} of $X$
if $(\phi_k)$ is a Schauder basis of $X$ and
for certain constants $c_1,c_2>0$ and all sequences $(a_k)_{k} \in \ell_{q}$ we have
\[
c_1\sum_{k=1}^\infty |a_k|^q \le \left\| \sum_{k=1}^\infty a_k \phi_k \right\|^q_X\le  c_2 \sum_{k=1}^\infty |a_k|^q.
\]
\end{definition}

Assume that $U=\mathbb C$ and $q\in[1,\infty)$. Further, let $X$ be a separable Hilbert space and let $X$ possess a $q$-Riesz basis consisting of eigenvectors $(\phi_k)$ of $A$ with eigenvalues $(\lambda_k)$.

We further assume that the sequence $(\lambda_k)_{k}$ lies in $\mathbb C$ with $\sup_{k\in\N}\re( \lambda_{k})<0$ and that there exists a constant $R>0$ such that $|\im\, \lambda_{k}|\le R|\re\, \lambda_{k}|$, $k\in \N$.

Hence the linear operator $A \colon D(A)\subset X \rightarrow X$, defined by
\[ A\phi_k = \lambda_k \phi_k, \quad k\in\N, \qquad D(A)= \Big\{\sum_{k=1}^\infty x_k\phi_k: (x_k)_{k}\in \ell_q \ \wedge \ \sum_{k=1}^\infty |x_k \lambda_k|^q <\infty\Big\},\]
generates an analytic semigroup $(T(t))_{t\ge 0}$ on $X$
(see \cite[Theorem 1.3.4]{Hen81} and note that $-A$ is a sectorial operator).
This semigroup is given by $T(t)\phi_k= {\rm{e}}^{t\lambda_k}\phi_k$, for all $t\geq 0$ and $k\in\N$.

In \cite[Theorem 4.1]{JNP18} the following result has been shown:
\begin{proposition}
\label{prop:Diagonal-Systems}
Let $\dim(U)<\infty$ and $\Uc:=L^\infty(\R_+,U)$. Assume that the operator $A$ possesses a $q$-Riesz basis of $X$ that consists  of eigenvectors
$(\phi_k)_{k}$  with eigenvalues $(\lambda_k)_{k}$ lying in a sector in the open left half-plane $\mathbb C_-$ with   $\sup_{k\in\N}\re( \lambda_{k})<0$ and $B\in L(U,X_{-1})$.

Then the system \eqref{eq:Linear_System} is an $L^\infty$-iISS and $L^\infty$-ISS control system in the sense of Definition~\ref{Steurungssystem}.
\end{proposition}

\begin{proof}
Assume first that $U=\C$. By \cite[Theorem 4.1]{JNP18}, it follows that $B$ is an $\infty$-admissible operator, and
the mild solution $\phi$ of \eqref{eq:Linear_System} satisfies the ISS and iISS estimates. The latter fact implies that $B$ is a so-called
zero-class $\infty$-admissible operator, i.e., in \eqref{eq:q-admissibility} it holds that $h_t \to 0$ as $t\to 0$.
Finally, \cite[Proposition 2.5]{JNP18} shows that $\phi$ is continuous.
This shows for $U=\C$ the triple $\Sigma:=(X,\Uc,\phi)$ is a control system in the sense of Definition~\ref{Steurungssystem} and
furthermore $\Sigma$ is ISS.

The case of general finite-dimensional $U$ can be reduced to the one-dimensional case \cite[Proposition 4]{JSZ19}.
\end{proof}

Proposition~\ref{prop:Diagonal-Systems} gives one more (in addition to Proposition~\ref{prop:iISS-equals-ISS-analytic-case})
 class of systems for which $L^\infty$-iISS and $L^\infty$-ISS are equivalent concepts.

%
%

A simple criterion guaranteeing that $B$ belongs to $L(\C,X_{-1})$ can be found in \cite[p. 882]{JNP18}.
In particular, Proposition~\ref{prop:Diagonal-Systems} shows ISS of the one-dimensional heat equation with a Dirichlet boundary condition, see \cite{JNP18}.

As an example, we consider a class of Riesz-spectral boundary control systems studied in \cite{LSZ20}, \cite{LhS19}.

\begin{definition}{(see \cite[Definition 2.3.1]{CuZ95})}
\label{def: Riesz-spectral operator}
\index{operator!Riesz-spectral}
Let $X$ be a Hilbert space and $A:~D(A) \subset X \rightarrow X$ be a linear, closed operator. For $n \in \N$, let $\lambda_n$ be the eigenvalues of $A$ and $\phi_n \in D(A)$ the corresponding eigenvectors. $A$ is called a \emph{Riesz-spectral operator} if:
\begin{enumerate}[label=(\roman*)]
    \item $\left\{ \phi_n: n \in \mathbb{N} \right\}$ is a 2-Riesz basis.
    \item The closure of $\{ \lambda_n: n \in \mathbb{N} \}$ is totally disconnected, i.e.,\ 
    no two points $\lambda,\mu \in \clo{\{\lambda_n:n\in\N\}}$ can be connected by a segment entirely lying in $\clo{\{\lambda_n:n\in\N\}}$.
\end{enumerate}
\end{definition}

By \cite[Theorem 2.3.5]{CuZ95}, the spectrum $\sigma(A)$ of a Riesz-spectral operator $A$ is given by $\sigma(A):=\clo{\{\lambda_n:n\in\N\}}$,
and the growth bound of a semigroup $T$ generated by $A$ can be computed as 
\begin{eqnarray}
\omega_0:= \sup_{i\in\N} \re\lambda_i <0,
\label{eq:Growth-bound}
\end{eqnarray}
that is, $A$ satisfies the \emph{spectrum determined growth assumption}.

As an application of Proposition~\ref{prop:Diagonal-Systems}, we obtain
\begin{proposition}
\label{prop:Riesz-spectral-sys} 
Let $X$ be a separable Hilbert space and consider BCS \eqref{eq:BCS}. Assume that the operator $A$ is a Riesz-spectral operator that generates an exponentially stable analytic semigroup. 
Furthermore, assume that $\dim (U) <\infty$ and $\Uc:=L^\infty(\R_+,U)$.
Then $B$ given by \eqref{eq:BCS-Input-Operator} is an admissible operator and the system
\eqref{eq:BCS} is ISS with these $X$ and $\Uc$.
\end{proposition}

\begin{proof}
Follows from Proposition~\ref{prop:Diagonal-Systems}.
\end{proof}

We note that in the papers \cite{LSZ20}, \cite{LhS19}, a different method has been employed for the study of ISS of Riesz-spectral BCS, which is a modification of the spectral method, briefly explained in Section~\ref{sec:Well-posedness-Concluding remarks}.
The essence of the method is a decomposition of $X$ w.r.t.\ the Riesz basis $\{\phi_k:k\in\N\}$.

\section{Concluding remarks}

The concepts of integral input-to-state stability (integral ISS) and iISS Lyapunov functions have been introduced in \cite{Son98}. 
Since then, they have been studied intensively. 
Characterizations of iISS property have been developed in a series of papers \cite{Son98,ASW00,ASW00b,AIS04b}.
The notion of the strong integral ISS has been introduced in \cite{CAI14} and further studied in \cite{CAI14b}.

We refer to \cite{Mir23} and \cite{Son08} for an overview of the available results in the integral ISS theory of finite-dimensional systems. 
For infinite-dimensional systems, the iISS theory is much less developed, and many questions remain open \cite[Section 10.1]{MiP20}.

Proposition~\ref{PropSufiISS} is due to \cite{MiI16}. The proof is analogous to the corresponding result for ODE systems \cite[p. 1088]{ASW00}.

In Propositions~\ref{ConverseLyapunovTheorem_BilinearSystems}, \ref{prop:iISS-noncoercive-LFs-II}, we summarized several results due to \cite[Theorem 4.2]{MiI16}, \cite[Proposition 5]{MiW15}, \cite[Proposition 3.34]{MiP20}. 
Examples in Section~\ref{sec:Examples-iISS} are taken from \cite{MiI15b,MiI16}. 

Results in Sections~\ref{sec:iISS of linear systems}, \ref{sec:diagonal_systems} are due to \cite{JNP18}.
Proposition~\ref{prop:iISS_implies_ISS} is due to \cite[Proposition 2.10]{JNP18} and Proposition~\ref{prop:iISS-equals-ISS-analytic-case} is due to \cite[Theorem 2]{JSZ19}.

Bilinear systems with admissible control operators have been studied in \cite{HJS22}, and the results have been applied to the controlled Fokker-Planck equation. Nevertheless, the iISS theory of bilinear systems with unbounded input operators remains much less developed than the ISS theory of linear systems with admissible input operators.



\cleardoublepage
\chapter{Infinite interconnections: Non-Lyapunov methods}
\label{chap:Infinite interconnections: Non-Lyapunov methods}


In real-world systems, the size of networks and the number of connections between their components may be unknown and time-varying, which poses new challenges for stability analysis and control design of such networks. A promising approach to this problem is to over-approximate the network by an infinite network and to perform the stability analysis and control design for this infinite over-approximation \cite{CIZ09,JoB05b,DaD03}. This approach has received significant attention during the last two decades. In particular, a large body of literature is devoted to spatially invariant systems and/or linear systems consisting of an infinite number of components interconnected by the same pattern \cite{BPD02,BaV05,BeJ17,CIZ09}.

Recently, a number of works appeared devoted to the stability and control of nonlinear infinite networks of ordinary differential equations that are not necessarily spatially invariant \cite{DMS19a,DaP20,KMS21,MNK21, MKG20}. Small-gain analysis of infinite networks is especially challenging since the gain operator, collecting the information about the internal gains, acts on an infinite-dimensional space, in contrast to couplings of finitely many systems of arbitrary nature. This calls for a careful choice of the infinite-dimensional state space of the overall network and motivates the use of the theory of positive operators on ordered Banach spaces for the small-gain analysis.

In this chapter, we provide rather general nonlinear trajectory-based ISS small-gain theorems for continuous-time infinite networks, whose components may be infinite-dimensional systems of different types.

\section{Infinite interconnections}
\label{sec:Infinite interconnections}

In this section, we introduce (feedback) interconnections of an arbitrary number of control systems indexed by some nonempty set $I$. For each $i \in I$, let $(X_i,\|\cdot\|_{X_i})$ be a normed vector space which will serve as the state space of a control system $\Sigma_i$ (in the sense of Definition~\ref{Steurungssystem}). Before we can specify the space of inputs for $\Sigma_i$, we first have to construct the overall state space. In the following, we use the sequence notation $(x_i)_{i \in I}$ for functions with domain $I$. The overall state space is then defined as%
\begin{equation*}
 	X := \Bigl\{ (x_i)_{i \in I} 
 	{\in \prod_{i \in I} X_i : } 
 	\, \sup_{i \in I} \|x_i\|_{X_i} < \infty \Bigr\}.%
\end{equation*}
It is a vector space with respect to pointwise addition and scalar multiplication, 
and we can turn it into a normed space in the following way:

\begin{proposition}
The state space $X$ is a normed space with respect to the norm
	\begin{equation*}
 		\|x\|_X := \sup_{i\in I}\|x_i\|_{X_i}.%
	\end{equation*}
	If all of the spaces $(X_i,\|\cdot\|_{X_i})$ are Banach spaces, then so is $(X,\|\cdot\|_X)$.
\end{proposition}

\begin{proof}
Let $(x^n)_{n \in \N}$ be a Cauchy sequence in $(X,\|\cdot\|_X)$. Then, for every $\eps>0$, there exists $N \in \N$ so that for all $n,m \geq N$ and $i\in I$ we have $\|x^n_i - x^m_i\|_{X_i} \leq \eps$. In particular, $(x^n_i)_{n\in\N}$ is Cauchy in $(X_i,\|\cdot\|_{X_i})$. Let $x_i \in X_i$ be its limit and put $x := (x_i)_{i \in I}$. We have%
\begin{equation*}
  \sup_{i \in I} \|x_i\|_{X_i} \leq \sup_{i \in I} \left(\|x_i - x_i^N\|_{X_i} + \|x_i^N\|_{X_i}\right) \leq \eps + \|x^N\|_X < \infty%
\end{equation*}
implying $x \in X$, and $x^n$ converges to $x$. %
\end{proof}

We also define for each $i \in I$ the normed vector space $X_{\neq i}$ by the same construction as above, but for the restricted index set $I \setminus \{i\}$. Then $X_{\neq i}$ can be identified with the closed linear subspace $\{ (x_j)_{j \in I} \in X : x_i = 0 \}$ of $X$.%

Now consider for each $i \in I$ a control system (in the sense of Definition~\ref{Steurungssystem}) of the form%
\begin{equation}
\label{eq:Sigma-i-inf-interconnections}
  \Sigma_i = (X_i,PC_b(\R_+,X_{\neq i}) \tm \Uc,\bar{\phi}_i),%
\end{equation}
where $PC_b(\R_+,X_{\neq i})$ is the space of all globally bounded piecewise continuous functions $w:\R_+ \rightarrow X_{\neq i}$, with the norm $\|w\|_{\infty} = \sup_{t \geq 0}\|w(t)\|_{X_{\neq i}}$. The norm on $PC_b(\R_+,X_{\neq i}) \tm \Uc$ is defined by%
\begin{equation}\label{eq_product_input_norm}
  \|(w,u)\|_{PC_b(\R_+,X_{\neq i}) \tm \Uc} := \max\left\{ \|w\|_{\infty}, \|u\|_{\Uc} \right\}.%
\end{equation}
Here we assume that $\Uc \subset U^{\R_+}$ for some vector space $U$, and $\Uc$ satisfies the axioms of shift invariance and concatenation. Then, by the definition of $PC_b(\R_+,X_{\neq i})$ and the norm \eqref{eq_product_input_norm}, these axioms are also satisfied for the product space $PC_b(\R_+,X_{\neq i}) \tm \Uc$.%

\begin{definition}
\label{def_interconnection}
\index{feedback interconnection}
Given the control systems $\Sigma_i$ ($i \in I$) as defined in \eqref{eq:Sigma-i-inf-interconnections}, assume that there is a map $\phi:D_\phi \to X$, defined on $D_\phi \subset \R_+ \tm X \tm \Uc$, such that:
\begin{enumerate}
	\item[(i)] For each $x \in X$ and each $u \in \Uc$ there is $\varepsilon>0$ such that $[0,\varepsilon] \tm \{(x,u)\} \subset D_\phi$.
	\item[(ii)] Furthermore, 
the components $\phi_i$ of the transition map $\phi:D_{\phi} \rightarrow X$ satisfy%
\begin{equation*}
  \phi_i(t,x,u) = \bar{\phi}_i(t,x_i,(\phi_{\neq i},u)) \quad \forall(t,x,u) \in D_{\phi},%
\end{equation*}
where $\phi_{\neq i}(\cdot) = (\phi_j(\cdot,x,u))_{j \in I \setminus \{i\}}$ for all $i \in I$.\footnote{By the causality axiom, we can assume that $\phi_{\neq i}$ is globally bounded, since $\bar{\phi}_i(t,x_i,(\phi_{\neq i},u))$ does not depend on the values $\phi_{\neq i}(s)$ with $s > t$, and on the compact interval $[0,t]$, $\phi_{\neq i}$ is bounded because it is continuous.}

We also assume that $\phi$ is maximal in the sense that no other map $\tilde\phi:\tilde{D}_\phi \to X$ with $\tilde{D}_\phi \supset D_\phi$ exists, which satisfies all of the above properties, and coincides with $\phi$ on $D_\phi$.
\end{enumerate}

If the map $\phi$ is unique with the above properties, and if $\Sigma = (X,\Uc,\phi)$ is a control system satisfying the BIC property, then $\Sigma$ is called the \emph{(feedback) interconnection} of the systems $\Sigma_i$.

We then call $X_{\neq i}$ the space of \emph{internal input values}, $PC_b(\R_+,X_{\neq i})$ the space of \emph{internal inputs}, and $\Uc$ the space of \emph{external inputs} of the system $\Sigma_i$. Moreover, we call $\Sigma_i$ the \emph{$i$-th subsystem} of $\Sigma$.
\end{definition}

The stability properties introduced above are defined in terms of the {norm} of the whole input, and this is not suitable for the consideration of coupled systems, as we are interested not only in the collective influence of all inputs on a subsystem but in the influence of particular subsystems on a given subsystem. The following definition provides the needed flexibility.

Let $Z$ be a Banach space. For $g \in PC(\R_+,Z)$ and $t>0$ we denote 
$\|g\|_{[0,t]} := \sup_{s\in[0,t]}|g(s)|$.

\begin{definition}\label{def_subsys_iss_semimax}
\index{input-to-state stability!semi-maximum formulation}
Given the spaces $(X_j,\|\cdot\|_{X_j})$, $j\in I$, and the system $\Sigma_i$ for a fixed $i \in I$, we say that $\Sigma_i$ is \emph{input-to-state stable (ISS) (in semi-maximum formulation)} if $\Sigma_i$ is forward complete and there are $\gamma_{ij},\gamma_j \in \K \cup \{0\}$ for all $j \in I$, and $\beta_i \in \KL$ such that for all initial states $x_i \in X_i$, all internal inputs $w_{\neq i} = (w_j)_{j\in I \setminus \{i\}} \in PC_b(\R_+,X_{\neq i})$, all external inputs $u \in \Uc$ and $t \geq 0$:%
\begin{equation*}
  \|\bar{\phi}_i(t,x_i,(w_{\neq i},u))\|_{X_i} \leq \beta_i(\|x_i\|_{X_i},t) + \sup_{j \in I}\gamma_{ij}(\|w_j\|_{[0,t]}) + \gamma_i(\|u\|_{\Uc}).%
\end{equation*}
Here we assume that the functions $\gamma_{ij}$ satisfy $\sup_{j \in I}\gamma_{ij}(r) < \infty$ for every $r \geq 0$ (implying that the sum on the right-hand side
is finite) and $\gamma_{ii} = 0$.
\end{definition}

The functions $\gamma_{ij}$ and $\gamma_i$ in this definition are called \emph{(nonlinear) gains}. 

Next we write $\ell_{\infty}(I)$ for the Banach space of all functions 
$x:I \rightarrow \R$ with $\|x\|_{\ell_{\infty}(I)} := \sup_{i\in I}|x(i)| < \infty$. 
Moreover, we make $\ell_{\infty}(I)$ to an ordered Banach space by introducing the order via the positive cone $\ell_{\infty}^+(I)$:
\[
\ell_{\infty}^+(I) := \{ x \in \ell_{\infty}(I) : x(i) \geq 0,\ \forall i \in I \}.
\]
We refer to Appendix~\ref{section:ordered-banach-spaces} for the basics of the theory of ordered Banach spaces, upon which we base ourselves in this chapter.

Assuming that all systems $\Sigma_i$, $i\in I$, are ISS in semi-maximum formulation, we can define a nonlinear monotone operator $\Gamma_{\otimes}:\ell_{\infty}^+(I) \rightarrow \ell_{\infty}^+(I)$ from the gains $\gamma_{ij}$ by
\begin{equation}
\label{eq:Gain-operator-semimax}
  \Gamma_{\otimes}(s) := \bigl(\sup_{j\in I}\gamma_{ij}(s_j)\bigr)_{i\in I},\quad s = (s_i)_{i\in I} \in \ell_{\infty}^+(I).%
\end{equation}
This operator and similar operators, composed of individual gains, are frequently called \emph{gain operators}.

In general, $\Gamma_{\otimes}$ is not well-defined. It is easy to see that the following assumption is equivalent to $\Gamma_{\otimes}$ being well-defined.

\begin{ass}\label{ass_gammamax_welldef}
For every $r>0$, we have $\sup_{i,j \in I}\gamma_{ij}(r) < \infty$.
\end{ass}

\begin{lemma}
Assumption \ref{ass_gammamax_welldef} is equivalent to the existence of $\zeta\in\Kinf$ and $a\geq0$ such that $\sup_{i,j\in I}\gamma_{ij}(r) \leq a + \zeta(r)$ for all $r\geq 0$.
\end{lemma}

\begin{proof}
Obviously, the implication \q{$\Leftarrow$} holds. Conversely, define $\xi: \R_+ \to \R_+$ by
\begin{equation*}
  \xi(r) := \sup_{i,j\in I}\gamma_{ij}(r).
\end{equation*}
As a supremum of continuous increasing functions, $\xi$ is lower semicontinuous and nondecreasing on its domain of definition. As $\xi(r)$ is finite for every $r\geq 0$ by assumption, define
\begin{equation*}
\tilde{\xi}(r):=
\begin{cases}
0 & \mbox{if } r=0,\\
\xi(r)-a & \mbox{if } r>0,%
\end{cases}
\end{equation*}
where $a:=\lim_{r\to+0}\xi(r)\geq 0$ (the limit exists as $\xi$ is nondecreasing). By construction, $\tilde{\xi}$ is nondecreasing, continuous at $0$ and satisfies $\tilde{\xi}(0)=0$. Hence, $\tilde{\xi}$ can be upper bounded by a certain $\zeta\in \Kinf$ (this follows from a more general result in \cite[Proposition~9]{MiW19b}). Overall, $\sup_{i,j\in I} \gamma_{ij}(r) \leq a + \zeta(r)$ for all $r\geq 0$. 
\end{proof}

Also, observe that $\Gamma_{\otimes}$, if well-defined, is a monotone operator. That is, for all $s^1,s^2 \in \ell_{\infty}^+(I)$ the following implication holds
\begin{equation*}
  s^1 \leq s^2 \quad\Rightarrow\quad \Gamma_{\otimes}(s^1) \leq \Gamma_{\otimes}(s^2).
\end{equation*}

\begin{remark}\label{rem_lg_sg_operator}
If all gains $\gamma_{ij}$ are linear, then $\Gamma_{\otimes}$ satisfies the following two properties:%
\begin{enumerate}[label=(\roman*)]
\item $\Gamma_{\otimes}$ is a homogeneous operator of degree one, i.e., $\Gamma_{\otimes}(as) = a\Gamma_{\otimes}(s)$ for all $a \geq 0$ and $s \in \ell_{\infty}^+(I)$.
\item $\Gamma_{\otimes}$ is subadditive, i.e., $\Gamma_{\otimes}(s^1 + s^2) \leq \Gamma_{\otimes}(s^1) + \Gamma_{\otimes}(s^2)$ for all $s^1,s^2 \in \ell_{\infty}^+(I)$.
\end{enumerate}
\end{remark}

Finally, we provide a criterion for continuity of $\Gamma_{\otimes}$ 
\begin{proposition}
Assume that the family $\{\gamma_{ij}\}_{(i,j) \in I^2}$ is pointwise equicontinuous, i.e., for every $r \in \R_+$ and $\eps>0$ there exists $\delta>0$ such that $|\gamma_{ij}(r) - \gamma_{ij}(s)| \leq \eps$ whenever $(i,j) \in I^2$ and $|r - s| \leq \delta$. Then $\Gamma_{\otimes}$ is well-defined and continuous.%
\end{proposition}

\begin{proof}
First, we show that $\Gamma_{\otimes}$ is well-defined. Fixing some $r>0$, the family $\{\gamma_{ij}\}_{(i,j)\in I^2}$ is uniformly equicontinuous on the compact interval $[0,r]$, which follows by a compactness argument. Hence, we can find $\delta>0$ so that $|s_1 - s_2| \leq \delta$ with $s_1,s_2 \in [0,r]$ implies $|\gamma_{ij}(s_1) - \gamma_{ij}(s_2)| \leq 1$ for all $i,j$. We can assume that $\delta$ is of the form $r/n$ for an integer $n$. Then%
\begin{equation*}
  \gamma_{ij}(r) = \sum_{k=0}^{n-1} \Bigl[\gamma_{ij}\Bigl(\frac{k+1}{n}r\Bigr) - \gamma_{ij}\Bigl(\frac{k}{n}r\Bigr)\Bigr] \leq n < \infty%
\end{equation*}
for all $(i,j) \in I^2$. Hence, $\Gamma_{\otimes}$ is well-defined.

Now we prove continuity. Choose any $\varepsilon>0$, fix some $s^0 \in \ell_{\infty}^+(I)$ and let $s \in \ell_{\infty}^+(I)$ so that $\|s - s^0\|_{\ell_{\infty}(I)} \leq \delta$ for some $\delta>0$ to be determined. By the required equicontinuity, we can choose $\delta$ small enough so that $|\gamma_{ij}(s^0_j) - \gamma_{ij}(s_j)| \leq \eps$ for all $(i,j)$ as $|s^0_j - s_j| \leq \|s^0 - s\|_{\ell_{\infty}(I)} \leq \delta$. This also implies%
\begin{equation*}
  \|\Gamma_{\otimes}(s^0) - \Gamma_{\otimes}(s)\|_{\ell_{\infty}(I)} = \sup_{i \in I}\Bigl|\sup_{j \in I} \gamma_{ij}(s^0_j) - \sup_{j \in I} \gamma_{ij}(s_j)\Bigr| \leq \eps.%
\end{equation*}
In the last inequality, we use the estimate
\begin{equation*}
  \sup_{j \in I} \gamma_{ij}(s^0_j) - \sup_{j \in I} \gamma_{ij}(s_j) \leq \sup_{j \in I} (\gamma_{ij}(s_j) + \eps) - \sup_{j \in I} \gamma_{ij}(s_j) = \eps,%
\end{equation*}
and the analogous estimate in the other direction. 
\end{proof}

Another formulation of ISS for the systems $\Sigma_i$ is as follows. In this formulation, we need to assume that $I$ is countable.
\begin{definition}\label{def_subsys_iss_sum}
\index{input-to-state stability!summation formulation}
Assume that $I$ is a nonempty countable set. Given the spaces $(X_j,\|\cdot\|_{X_j})$, $j\in I$, and the system $\Sigma_i$ for a fixed $i \in I$, we say that $\Sigma_i$ is \emph{input-to-state stable (ISS) (in summation formulation)} if $\Sigma_i$ is forward complete and there are $\gamma_{ij},\gamma_j \in \K \cup \{0\}$ for all $j \in I$, and $\beta_i \in \KL$ such that for all initial states $x_i \in X_i$, all internal inputs $w_{\neq i} = (w_j)_{j\in I \setminus \{i\}} \in PC_b(\R_+,X_{\neq i})$, all external inputs $u \in \Uc$ and $t \geq 0$:%
\begin{equation*}
  \|\bar{\phi}_i(t,x_i,(w_{\neq i},u))\|_{X_i} \leq \beta_i(\|x_i\|_{X_i},t) + \sum_{j \in I}\gamma_{ij}(\|w_j\|_{[0,t]}) + \gamma_i(\|u\|_{\Uc}).%
\end{equation*}
Here we assume that the functions $\gamma_{ij}$ are such that $\sum_{j \in I}\gamma_{ij}(r) < \infty$ for every $r \geq 0$ (implying that the sum on the right-hand side
is finite) and $\gamma_{ii} = 0$.
\end{definition}

\begin{remark}
\label{rem:Why-many-formulations} 
If a network has finitely many components, ISS in summation formulation and ISS in semi-maximum formulation are equivalent concepts. Nevertheless, even for finite networks, the gains in semi-maximum formulation and the gains in summation formulation are distinct, and for some systems, one formulation is better than the other one in the sense that it produces tighter (and thus smaller) gains. This motivates the interest in analyzing both formulations. We illustrate this by examples in Sections~\ref{examp: a linear spatially invariant system}, \ref{examp: a nonlinear spatially invariant system}. In fact, also more general formulations of input-to-state stability for networks are studied in the literature, using the formalism of monotone aggregation functions, see Section~\ref{sec:ISS-SGT-Lyapunov-form}.
\end{remark}

Assuming that all systems $\Sigma_i$, $i\in I$, are ISS, we can define a nonlinear monotone operator $\Gamma_{\boxplus}:\ell_{\infty}^+(I) \rightarrow \ell_{\infty}^+(I)$ from the gains $\gamma_{ij}$ as follows:%
\begin{equation*}
  \Gamma_{\boxplus}(s) := \Bigl(\sum_{j\in I}\gamma_{ij}(s_j)\Bigr)_{i \in I},\quad s = (s_i)_{i\in I} \in \ell_{\infty}^+(I).%
\end{equation*}
Again, $\Gamma_{\boxplus}$ might not be well-defined. Hence, we need to make an appropriate assumption.%

\begin{ass}\label{ass_gammasum_welldef}
For every $r > 0$, we have%
\begin{equation*}
  \sup_{i \in I}\sum_{j \in I} \gamma_{ij}(r) < \infty.%
\end{equation*}
\end{ass}

\begin{remark}
Assume that all the gains $\gamma_{ij}$, $(i,j) \in I^2$, are linear functions. Then the gain operator $\Gamma_{\boxplus}$ can be regarded as a linear operator on $\ell_{\infty}(I)$ and Assumption \ref{ass_gammasum_welldef} {is equivalent to $\Gamma_{\boxplus}$ being a bounded} linear operator on $\ell_{\infty}(I)$.
\end{remark}

\begin{proposition}
Assume that the operator $\Gamma_{\boxplus}$ is well-defined. A sufficient criterion for continuity of $\Gamma_{\boxplus}$ is that each $\gamma_{ij}$ is a $C^1$-function and%
\begin{equation*}
  \sup_{i \in I} \sum_{j \in I} \sup_{0 < s \leq r} \gamma_{ij}'(s) < \infty,\quad r > 0.%
\end{equation*}
\end{proposition}

\begin{proof}
Fix $s^0 = (s^0_j)_{j\in I} \in \ell_{\infty}^+(I)$ and $\eps>0$. Let $s \in \ell_{\infty}^+(I)$ with $\|s - s^0\|_{\ell_{\infty}(I)} = \sup_{i \in I}|s_i - s^0_i| \leq \delta$ for some $\delta > 0$, to be determined later. Then%
\begin{equation*}
  \|\Gamma_{\boxplus}(s^0) - \Gamma_{\boxplus}(s)\|_{\ell_{\infty}(I)} = \sup_{i \in I}\Bigl| \sum_{j \in I} (\gamma_{ij}(s^0_j) - \gamma_{ij}(s_j)) \Bigr|.
\end{equation*}
Using the assumption that each $\gamma_{ij}$ is a $C^1$-function and writing $s^0_{\max} := \|s^0\|_{\ell_{\infty}(I)}$, we can estimate this by%
\begin{align*}
  & \|\Gamma_{\boxplus}(s) - \Gamma_{\boxplus}(s^0)\|_{\ell_{\infty}(I)} \leq \sup_{i \in I} \sum_{j \in I} |\gamma_{ij}(s_j) - \gamma_{ij}(s^0_j)| \\
	                                   &\leq \sup_{i \in I} \sum_{j \in I} \sup_{r \in [s^0_j - \delta,s^0_j + \delta]}|\gamma_{ij}'(r)| |s_j - s^0_j|
	                                   \leq \delta \sup_{i \in I} \sum_{j \in I} \sup_{r \leq s^0_{\max}+\delta} \gamma_{ij}'(r).%
\end{align*}
By assumption, the last supremum is finite, which implies that $\delta$ can be chosen small enough so that the whole expression is less than $\eps$. 
\end{proof}

We also need versions of UGS for the systems $\Sigma_i$.

\begin{definition}\label{def_subsys_ugs_semimax}
Given the spaces $(X_j,\|\cdot\|_{X_j})$, $j\in I$, and the system $\Sigma_i$ for a fixed $i \in I$, we say that $\Sigma_i$ is \emph{uniformly globally stable (UGS) (in semi-maximum formulation)} if $\Sigma_i$ is forward complete and there are $\gamma_{ij},\gamma_j \in \K \cup \{0\}$ for all $j \in I$, and $\sigma_i \in \Kinf$ such that for all initial states $x_i \in X_i$, all internal inputs $w_{\neq i} = (w_j)_{j\in I \setminus \{i\}} \in PC_b(\R_+,X_{\neq i})$, all external inputs $u \in \Uc$ and $t \geq 0$:%
\begin{equation*}
  \|\bar{\phi}_i(t,x_i,(w_{\neq i},u))\|_{X_i} \leq \sigma_i(\|x_i\|_{X_i}) + \sup_{j \in I}\gamma_{ij}(\|w_j\|_{[0,t]}) + \gamma_i(\|u\|_{\Uc}).%
\end{equation*}
Here we assume that the functions $\gamma_{ij}$ are such that $\sup_{j \in I}\gamma_{ij}(r) < \infty$ for every $r \geq 0$ (implying that the sum on the right-hand side is finite) and $\gamma_{ii} = 0$.
\end{definition}

\begin{definition}\label{def_subsys_ugs_sum}
Let $I$ be a countable index set. Given the spaces $(X_j,\|\cdot\|_{X_j})$, $j\in I$, and the system $\Sigma_i$ for a fixed $i \in I$, we say that $\Sigma_i$ is \emph{uniformly globally stable (UGS) (in summation formulation)} if $\Sigma_i$ is forward complete and there are $\gamma_{ij},\gamma_j \in \K \cup \{0\}$ for all $j \in I$, and $\sigma_i \in \Kinf$ such that for all initial states $x_i \in X_i$, all internal inputs $w_{\neq i} = (w_j)_{j\in I \setminus \{i\}} \in PC_b(\R_+,X_{\neq i})$, all external inputs $u \in \Uc$ and $t \geq 0$:%
\begin{equation*}
  \|\bar{\phi}_i(t,x_i,(w_{\neq i},u))\|_{X_i} \leq \sigma_i(\|x_i\|_{X_i}) + \sum_{j \in I}\gamma_{ij}(\|w_j\|_{[0,t]}) + \gamma_i(\|u\|_{\Uc}).%
\end{equation*}
Here we assume that the functions $\gamma_{ij}$ are such that $\sum_{j \in I}\gamma_{ij}(r) < \infty$ for every $r \geq 0$ (implying that the sum on the right-hand side is finite) and $\gamma_{ii} = 0$.
\end{definition}

\section{Stability of discrete-time systems}
\label{sec:Stability of discrete-time systems}

In this section, we study the stability properties of the system
\begin{equation}
	\label{eq_monotone_system}
 	x(k+1) \leq A(x(k)) + u(k),\quad k \in \Z_+.
\end{equation}
Here, $(X,X^+)$ is an ordered Banach space, $A: X^+ \rightarrow X^+$ is a nonlinear operator on the cone $X^+$, and the input $u$ is an element of $\ell_{\infty}(\Z_+,X^+)$, where the latter space is defined as%
\begin{equation*}
  \ell_{\infty}(\Z_+,X^+) := \{u = (u(k))_{k\in\Z_+} : u(k) \in X^+,\ \|u\|_{\infty} := \sup_{k\in\Z_+}\|u(k)\|_X < \infty\}.%
\end{equation*}
A \emph{solution} of the equation~\eqref{eq_monotone_system} is a mapping $x: \Z_+ \to X^+$ that satisfies~\eqref{eq_monotone_system}.
We call a mapping $x: \Z_+ \to X^+$ \emph{decreasing} if $x(k+1) \le x(k)$ for all $k \in \Z_+$.

As we will see, for the small-gain analysis of infinite interconnections, the properties of the gain operator and the discrete-time system \eqref{eq_monotone_system} induced by the gain operator are essential. So we now relate the stability of the system \eqref{eq_monotone_system} to the properties of the operator $A$.

\begin{definition}
\label{def:MLIM}
\index{MLIM}
\index{property!monotone limit}
The system \eqref{eq_monotone_system} has the \emph{monotone limit property (MLIM)} if there is $\xi \in \Kinf$ such that for every $\eps>0$, every constant input $u(\cdot) :\equiv w \in X^+$ and every decreasing solution $x: \Z_+ \to X^+$ of \eqref{eq_monotone_system},
there exists $N = N(\eps,u,x(\cdot)) \in \Z_+$ with%
\begin{equation*}
  \|x(N)\|_X \leq \eps + \xi(\|w\|_X).%
\end{equation*}
\end{definition}

\begin{definition}
\label{def:MBI}
\index{MBI}
\index{property!monotone bounded invertibility}
Let $(X,X^+)$ be an ordered Banach space, and let $A:X^+ \rightarrow X^+$ be a nonlinear operator. We say that $\id - A$ has the \emph{monotone bounded invertibility (MBI) property} if there exists $\xi \in \Kinf$ such that for all $v,w \in X^+$ the following implication holds:
\begin{equation*}
  (\id - A)(v) \leq w \quad \Rightarrow \quad \|v\|_X \leq \xi(\|w\|_X).%
\end{equation*}
\end{definition}

\begin{proposition}\label{prop_mlim_implies_mbip}
Let $(X,X^+)$ be an ordered Banach space, and let $A:X^+ \rightarrow X^+$ be a nonlinear operator.
If the system \eqref{eq_monotone_system} has the MLIM property, then the operator $\id - A$ has the MBI property.
\end{proposition}

\begin{proof}
Assume that $(\id - A)(v) \leq w$ for some $v,w \in X^+$. We write this as $v \leq A(v) + w$. Hence, $x(\cdot) :\equiv v$ is a constant solution of \eqref{eq_monotone_system} corresponding to the constant input sequence $u(\cdot) :\equiv w$. By the MLIM property, there exists $\xi \in \Kinf$ (independent of $v,w$) so that for every $\eps>0$ there is $N$ with%
\begin{equation*}
  \|v\|_X = \|x(N)\|_X \leq \eps + \xi(\|u\|_{\infty}) = \eps + \xi(\|w\|_X).%
\end{equation*}
Since this holds for every $\eps>0$, we obtain $\|v\|_X \leq \xi(\|w\|_X)$, which completes the proof. %
\end{proof}

\emph{It is an open problem whether the MBI property is strictly weaker than the MLIM property or whether they are equivalent.} In the following proposition, though, we show that they are equivalent under certain assumptions on the operator $A$ or the cone $X^+$. Later, in Propositions~\ref{prop:eISS-criterion-linear-systems} and \ref{prop:join-morphism-eISS-criterion}, we show their equivalence for linear operators and for the gain operator $\Gamma_\otimes$ with linear gains, defined on $\ell_{\infty}(I)$.%

The cone $X^+$ is said to have the \emph{Levi property} if every decreasing sequence in $X^+$ is norm-convergent \cite[Def.~2.44(2)]{AlT07}. Typical examples are the standard cones in $L^p$-spaces for $p \in [1,\infty)$ and the standard cone in the space $c_0$ of real-valued sequences that converge to $0$. We note in passing that if the cone $X^+$ has the Levi property, then it is normal \cite[Theorem~2.45]{AlT07}.

\begin{proposition}\label{prop_compact_operators}
Let $(X,X^+)$ be an ordered Banach space with normal cone, and let $A:X^+ \rightarrow X^+$ be a nonlinear, continuous, and monotone operator. 
If the cone $X^+$ has the Levi property or if the operator $A$ is compact (i.e., it maps bounded sets onto precompact sets), then the following statements are equivalent:%
\begin{enumerate}
\item[(i)] The system \eqref{eq_monotone_system} satisfies the MLIM property.%
\item[(ii)] The operator $\id - A$ satisfies the MBI property.%
\end{enumerate}
\end{proposition}

\begin{proof}
In view of Proposition \ref{prop_mlim_implies_mbip}, it suffices to prove the implication ``(ii) $\Rightarrow$ (i)''. Hence, consider a constant input $u(\cdot) :\equiv w \in X^+$ and a decreasing sequence $x(\cdot)$ in $X^+$ such that $x(k+1) \leq A(x(k)) + w$ for all $k \in \Z_+$. As $A$ is monotone, the operator $\tilde{A}(x) := A(x) + w$, $\tilde{A}:X^+ \rightarrow X^+$, is monotone as well; if $A$ is compact, then so is $\tilde A$. Moreover,%
\begin{equation}\label{eq_newineq}
  x(k+1) \leq \tilde{A}(x(k)),\quad k \in \Z_+.%
\end{equation}
Now consider the sequence $y(k) := \tilde{A}(x(k))$, $k \in \Z_+$. As $x$ is decreasing, so is $y$. Next, we note that $y$ converges in the norm. Indeed, if the cone has the Levi property, this is clear. If the cone does not have the Levi property, then $A$, and thus $\tilde{A}$, is compact by assumption. So $y$ has a convergent subsequence. 
Since $y$ is decreasing and the cone is normal, it thus follows that $y$ converges itself.%

Let $y_* \in X^+$ denote the limit of the sequence $y$. Applying $\tilde{A}$ on both sides of \eqref{eq_newineq} yields%
\begin{equation*}
  y(k+1) \leq \tilde{A}(y(k)),\quad k \in \Z_+.%
\end{equation*}
Taking the limit for $k \rightarrow \infty$ and using the continuity of $A$ results in $y_* \leq \tilde{A}(y_*) = A(y_*) + w$. Since this can be written as $(\id - A)(y_*) \leq w$, the MBI property of $\id - A$ gives $\|y_*\| \leq \xi(\|w\|)$. As $X^+$ is a normal cone, there is $\delta>0$ such that for every $\varepsilon>0$ there is $k>0$ large enough, for which
\begin{equation*}
  \|x(k+1)\|_X \leq \delta\|\tilde{A}(x(k))\|_X \leq \eps + \delta\xi(\|w\|_X).%
\end{equation*}
This completes the proof. 
\end{proof}

\section{Small-gain theorems}\label{sec:Small-gain theorems}

\subsection{Small-gain theorems in semi-maximum formulation}

In this subsection, we prove small-gain theorems for UGS and ISS, both in the semi-maximum formulation. We start with UGS.%

\begin{theorem}[UGS small-gain theorem in semi-maximum formulation]
\label{thm_ugs_semimax_sg}
\index{ISS small-gain theorem!semi-maximum formulation}
\index{ISS small-gain theorem!trajectory-form}
Let $I$ be an arbitrary nonempty index set, $(X_i,\|\cdot\|_{X_i})$, $i\in I$, normed spaces, and $\Sigma_i = (X_i,PC_b(\R_+,X_{\neq i}) \tm \Uc,\bar{\phi}_i)$ forward complete control systems. Assume that the interconnection $\Sigma = (X,\Uc,\phi)$ of the systems $\Sigma_i$ is well-defined. Furthermore, let the following assumptions be satisfied:%
\begin{enumerate}
\item[(i)] Each system $\Sigma_i$ is UGS in the sense of Definition \ref{def_subsys_ugs_semimax} with $\sigma_i \in \K$ and nonlinear gains $\gamma_{ij},\gamma_i \in \K \cup \{0\}$.%
\item[(ii)] There exist $\sigma_{\max} \in \Kinf$ and $\gamma_{\max} \in \Kinf$ so that $\sigma_i \leq \sigma_{\max}$ and $\gamma_i \leq \gamma_{\max}$, pointwise for all $i \in I$.%
\item[(iii)] Assumption \ref{ass_gammamax_welldef} is satisfied for the operator $\Gamma_{\otimes}$ defined via the gains $\gamma_{ij}$ from (i) and $\id - \Gamma_{\otimes}$ has the MBI property.%
\end{enumerate}
Then $\Sigma$ is forward complete and UGS.%
\end{theorem}

\begin{proof}
Fix $(t,x,u) \in D_{\phi}$ and observe that%
\begin{equation*}
  \|\phi(t,x,u)\|_X = \sup_{i \in I} \|\phi_i(t,x,u)\|_{X_i} = \sup_{i \in I} \|\bar{\phi}_i(t,x_i,(\phi_{\neq i},u))\|_{X_i}.%
\end{equation*}
Abbreviating $\bar{\phi}_j(\cdot) = \bar{\phi}_j(\cdot,x_j,(\phi_{\neq j},u))$ and using assumption (i), we can estimate%
\begin{equation}\label{eq_ugs_firstest}
  \sup_{s \in [0,t]}\|\bar{\phi}_i(s,x_i,(\phi_{\neq i},u))\|_{X_i} \leq \sigma_i(\|x_i\|_{X_i}) + \sup_{j \in I}\gamma_{ij}(\|\bar{\phi}_j\|_{[0,t]}) + \gamma_i(\|u\|_{\Uc}).%
\end{equation}
From the inequalities (using continuity of $s \mapsto \phi(s,x,u)$)%
\begin{equation*}
  0 \leq \sup_{s\in[0,t]} \|\bar{\phi}_i(s,x_i,(\phi_{\neq i},u))\|_{X_i} \leq \sup_{s\in[0,t]} \|\phi(s,x,u)\|_X < \infty,\quad i \in I,%
\end{equation*}
it follows that%
\begin{equation*}
  \vec{\phi}_{\max}(t) := \Bigl( \sup_{s \in [0,t]}\|\bar{\phi}_i(s,x_i,(\phi_{\neq i},u))\|_{X_i} \Bigr)_{i \in I} \in \ell_{\infty}^+(I).%
\end{equation*}
From Assumption (ii), it follows that also the vectors $\vec{\sigma}(x) := (\sigma_i(\|x_i\|_{X_i}))_{i\in I}$ and $\vec{\gamma}(u) := ( \gamma_i(\|u\|_{\Uc}) )_{i \in I}$ are contained in $\ell_{\infty}^+(I)$. Hence, we can write the inequalities \eqref{eq_ugs_firstest} in vectorized form as%
\begin{equation*}
  (\id - \Gamma_{\otimes})(\vec{\phi}_{\max}(t)) \leq \vec{\sigma}(x) + \vec{\gamma}(u).%
\end{equation*}
By Assumption (iii), this yields for some $\xi\in\Kinf$, independent of $x,u$:%
\begin{equation*}
  \|\vec{\phi}_{\max}(t)\|_{\ell_{\infty}(I)} \leq \xi ( \|\vec{\sigma}(x) + \vec{\gamma}(u) \|_{\ell_{\infty}(I)} ) \leq \xi( \|\vec{\sigma}(x)\|_{\ell_{\infty}(I)} + \|\vec{\gamma}(u)\|_{\ell_{\infty}(I)} ).%
\end{equation*}
Since $\xi(a + b) \leq \max\{\xi(2a),\xi(2b)\} \leq \xi(2a) + \xi(2b)$ for all $a,b\geq 0$, this implies%
\begin{eqnarray*}
  \|\vec{\phi}_{\max}(t)\|_{\ell_{\infty}(I)} 
	&\leq& \xi(2\|\vec{\sigma}(x)\|_{\ell_{\infty}(I)}) + \xi(2\|\vec{\gamma}(u)\|_{\ell_{\infty}(I)})\\
	&\leq& \xi(2\sigma_{\max}(\|x\|_X)) + \xi(2\gamma_{\max}(\|u\|_{\Uc})),
\end{eqnarray*}
and we conclude that%
\begin{eqnarray*}
  \|\phi(t,x,u)\|_X \leq \|\vec{\phi}_{\max}(t)\|_{\ell_{\infty}(I)} \leq \xi(2\sigma_{\max}(\|x\|_X)) + \xi(2\gamma_{\max}(\|u\|_{\Uc})),%
\end{eqnarray*}
which is a UGS estimate with $\sigma(r) := \xi(2\sigma_{\max}(r))$, $\gamma(r) := \xi(2\gamma_{\max}(r))$ for $\Sigma$ for all $(t,x,u)\in D_\phi$. Since $\Sigma$ has the BIC property by assumption, it follows that $\Sigma$ is forward complete and UGS. %
\end{proof}

Now we are in a position to state the ISS small-gain theorem.%

\begin{theorem}[Nonlinear ISS small-gain theorem in semi-maximum formulation]
\label{thm_smallgain_iss_semimax}
\index{ISS small-gain theorem!semi-maximum formulation}
Let $I$ be an arbitrary nonempty index set, $(X_i,\|\cdot\|_{X_i})$, $i\in I$, normed spaces, and $\Sigma_i = (X_i,PC_b(\R_+,X_{\neq i}) \tm \Uc,\bar{\phi}_i)$ forward complete control systems. Assume that the interconnection $\Sigma = (X,\Uc,\phi)$ of the systems $\Sigma_i$ is well-defined. Furthermore, let the following assumptions be satisfied:%
\begin{enumerate}
\item[(i)] Each system $\Sigma_i$ is ISS in the sense of Definition \ref{def_subsys_iss_semimax} with $\beta_i \in \KL$ and nonlinear gains $\gamma_{ij},\gamma_i \in \K \cup \{0\}$.%
\item[(ii)] There are $\beta_{\max} \in \KL$ and $\gamma_{\max} \in \K$ so that $\beta_i \leq \beta_{\max}$ and $\gamma_i \leq \gamma_{\max}$ pointwise for all $i \in I$.%
\item[(iii)] Assumption \ref{ass_gammamax_welldef} holds and the discrete-time system%
\begin{equation}\label{eq_sg_system}
  w(k+1) \leq \Gamma_{\otimes}(w(k)) + v(k),%
\end{equation}
with $w(\cdot),v(\cdot)$ taking values in $\ell_{\infty}^+(I)$, has the MLIM property.%
\end{enumerate}
Then $\Sigma$ is ISS.
\end{theorem}

\begin{proof}
We show that $\Sigma$ is UGS and satisfies the bUAG property, which implies ISS by Theorem~\ref{thm:UAG_equals_ULIM_plus_LS}.

{\bf UGS}. This follows from Theorem \ref{thm_ugs_semimax_sg}. Indeed, the assumptions (i) and (ii) of Theorem \ref{thm_ugs_semimax_sg} are satisfied with $\sigma_i(\cdot) := \beta_i(\cdot,0) \in \K$ and the gains $\gamma_{ij},\gamma_i$ from the ISS estimates for $\Sigma_i$, $i\in I$. From Proposition \ref{prop_mlim_implies_mbip} and Assumption (iii) of this theorem, it follows that Assumption (iii) of Theorem \ref{thm_ugs_semimax_sg} is satisfied. Hence, $\Sigma$ is forward complete and UGS.%

{\bf bUAG}. As $\Sigma$ is the interconnection of the systems $\Sigma_i$ and since $\Sigma$ is forward complete, we have $\phi_i(t,x,u) = \bar{\phi}_i(t,x_i,(\phi_{\neq i},u))$ for all $(t,x,u) \in \R_+ \tm X \tm \Uc$ and $i \in I$, with the notation from Definition \ref{def_interconnection}.%

Pick any $r > 0$, any $u \in \clo{B_{r,\Uc}}$ and $x \in \clo{B_{r,X}}$. As $\Sigma$ is UGS, there are $\sigma^{\UGS},\gamma^{\UGS} \in \Kinf$ so that%
\begin{equation*}
  \|\phi(t,x,u)\|_X \leq \sigma^{\UGS}(r) + \gamma^{\UGS}(r) =: \mu(r),\quad t \geq 0.%
\end{equation*}
In view of the cocycle property, for all $i \in I$ and $t,\tau \geq 0$ we have%
\begin{align*}
  \phi_i(t + \tau,x,u) &= \bar{\phi}_i(t+\tau,x_i,(\phi_{\neq i},u)) \\
	                     &= \bar{\phi}_i(\tau,\bar{\phi}_i(t,x_i,(\phi_{\neq i},u)),(\phi_{\neq i}(\cdot+t),u(\cdot+t))).%
\end{align*}
Given $\eps>0$, choose $\tau^* = \tau^*(\eps,r) \geq 0$ such that $\beta_{\max}(\mu(r),\tau^*) \leq \eps$. Then%
\begin{align}\label{eq_firstest}
\begin{split}
  x \in \clo{B_{r,X}} &\ \ \wedge \ \  u \in \clo{B_{r,\Uc}} \ \ \wedge \ \  \tau \geq \tau^* \ \ \wedge \ \  t \geq 0 \quad \Rightarrow\quad\\
	&\qquad\Rightarrow \|\phi_i(t+\tau,x,u)\|_{X_i} \leq \beta_i(\|\bar{\phi}_i(t,x_i,(\phi_{\neq i},u))\|_{X_i},\tau) \\
	                                 &\qquad\qquad\qquad + \sup_{j \in I} \gamma_{ij}( \|\phi_j\|_{[t,t+\tau]} ) + \gamma_i(\|u(\cdot+t)\|_{\Uc}) \\
																	&\qquad\leq \beta_{\max}(\|\phi(t,x,u)\|_X,\tau^*) + \sup_{j \in I}\gamma_{ij}(\|\phi_j\|_{[t,\infty)}) + \gamma_i(\|u\|_{\Uc}) \\
																	&\qquad\leq \eps + \sup_{j \in I}\gamma_{ij}(\|\phi_j\|_{[t,\infty)}) + \gamma_i(\|u\|_{\Uc}).%
\end{split}
\end{align}
Now pick any $k \in \N$ and write%
\begin{equation*}
  B(r,k) := \clo{B_{r,X}} \tm \{ u \in \Uc : \|u\|_{\Uc} \in [2^{-k}r,2^{-k+1}r]\}.%
\end{equation*} 
Then, taking the supremum in the above inequality over all $(x,u) \in B(r,k)$, we obtain for all $i \in I$ and all $t \geq 0$ that%
\begin{equation*}
  \sup_{(x,u) \in B(r,k)}\|\phi_i(t+\tau^*,x,u)\|_{X_i} \leq \eps + \sup_{j \in I} \gamma_{ij}\Bigl( \sup_{(x,u) \in B(r,k)} \|\phi_j\|_{[t,\infty)} \Bigr) + \gamma_i(2^{-k+1}r).%
\end{equation*}
This implies for all $t \geq 0$ that%
\begin{align*}
  &\sup_{s \geq t + \tau^*}\sup_{(x,u) \in B(r,k)} \|\phi_i(s,x,u)\|_{X_i} \\
	&\quad \leq \eps + \sup_{j \in I} \gamma_{ij}\Bigl(\sup_{s \geq t} \sup_{(x,u) \in B(r,k)} \|\phi_j(s,x,u)\|_{X_j}\Bigr) + \gamma_i(2^{-k+1}r).%
\end{align*}
Now we define%
\begin{equation*}
  w_i(t,r,k) := \sup_{s \geq t} \sup_{(x,u) \in B(r,k)} \|\phi_i(s,x,u)\|_{X_i}%
\end{equation*}
and note that $w_i(t,r,k) \in [0,\mu(r)]$ for all $i \in I$ and $t \geq 0$. With this notation, we can rewrite the preceding inequality as%
\begin{equation*}
  w_i(t + \tau^*,r,k) \leq \eps + \sup_{j \in I} \gamma_{ij}(w_i(t,r,k)) + \gamma_i(2^{-k+1}r).%
\end{equation*}

Using vector notation $\vec{w}(t,r,k) := (w_i(t,r,k))_{i\in I}$ and $\vec{\gamma}(r) := (\gamma_i(r))_{i \in I}$, and denoting $\unit:=(1,1,\ldots)$, this can be written as%
\begin{equation*}
  \vec{w}(t+\tau^*,r,k) \leq \Gamma_{\otimes}(\vec{w}(t,r,k)) + \eps \unit + \vec{\gamma}(2^{-k+1}r).%
\end{equation*}
Observe that $\vec{w}(t,r,k) \in \ell_{\infty}^+(I)$, as the entries of the vector are uniformly bounded by $\mu(r)$, and $\vec{w}(t_2,r,k) \leq \vec{w}(t_1,r,k)$ for $t_2 \geq t_1$. Hence, $\vec{w}(l) := \vec{w}(l\tau^*,r,k)$, $l \in \Z_+$, is a monotone solution of \eqref{eq_sg_system} for the constant input $v(\cdot) \equiv \eps \unit + \vec{\gamma}(2^{-k+1}r)$. By assumption (iii) of the theorem, this implies the existence of a time $\tilde{\tau} = \tilde{\tau}(\eps,r,k)$ and a $\Kinf$-function $\xi$ such that%
\begin{align*}
  \|\vec{w}(\tilde{\tau},r,k)\|_{\ell_{\infty}(I)} &\leq \eps + \xi(\|\eps \unit + \vec{\gamma}(2^{-k+1}r)\|_{\ell_{\infty}(I)}) \\
	                                              &\leq \eps + \xi(\|\eps{\unit}\|_{\ell_{\infty}(I)} + \|\vec{\gamma}(2^{-k+1}r)\|_{\ell_{\infty}(I)}) \\
																								&\leq \eps + \xi(\eps + \gamma_{\max}(2^{-k+1}r)) \\
																								&\leq \eps + \xi(2\eps) + \xi(2\gamma_{\max}(2^{-k+1}r)).%
\end{align*}
By definition, this implies%
\begin{align*}
  &i \in I \ \ \wedge \ \  (x,u) \in B(r,k) \ \ \wedge \ \  t \geq \tilde{\tau}(\eps,r,k) \\
	&\Rightarrow \|\phi_i(t,x,u)\|_{X_i} \leq \eps + \xi(2\eps) + \xi(2\gamma_{\max}(2^{-k+1}r)).%
\end{align*}
Now define $k_0 = k_0(\eps,r)$ as the minimal $k \geq 1$ so that $\xi(2\gamma_{\max}(2^{1-k}r)) \leq \eps$ and let%
\begin{equation*}
  \hat{\tau}(\eps,r) := \max\{ \tilde{\tau}(\eps,r,k) : 1 \leq k \leq k_0(\eps,r) \}.%
\end{equation*}
Pick any $0 \neq u \in \clo{B_{r,\Uc}}$. Then there is $k \in \N$ with $\|u\|_{\Uc} \in (2^{-k}r,2^{-k+1}r]$. If $k \leq k_0$ (large input), then for $t \geq \hat{\tau}(\eps,r)$ we have%
\begin{align}\label{eq_349}
\begin{split}
  \|\phi(t,x,u)\|_X &\leq \eps + \xi(2\eps) + \xi(\gamma_{\max}(2^{-k+1}r)) \\
	&\leq \eps + \xi(2\eps) + \xi(2\gamma_{\max}(2\|u\|_{\Uc})).%
\end{split}
\end{align}
It remains to consider the case when $k > k_0$ (small input). For any $q \in [0,r]$, one can take the supremum in \eqref{eq_firstest} over $x \in \clo{B_{r,X}}$ and $u \in \clo{B_{q,\Uc}}$ to obtain%
\begin{align*}
  &\sup_{(x,u) \in \clo{B_{r,X}} \tm \clo{B_{q,\Uc}}}\|\phi_i(t+\tau,x,u)\|_{X_i} \\
	&\qquad \leq \eps + \sup_{j \in I}\gamma_{ij}\Bigl(\sup_{(x,u) \in \clo{B_{r,X}} \tm \clo{B_{q,\Uc}}}\|\phi_j\|_{[t,\infty)}\Bigr) + \gamma_i(q).%
\end{align*}
With $z_i(t,r,q) := \sup_{s \geq t}\sup_{(x,u) \in \clo{B_{r,X}} \tm \clo{B_{q,\Uc}}}\|\phi_i(s,x,u)\|_{X_i}$, analogous steps as above lead to the following: for every $\eps>0$, $r>0$ and $q \in [0,r]$ there is a time $\bar{\tau} = \bar{\tau}(\eps,r,q)$ such that%
\begin{align*}
  (x,u) \in \clo{B_{r,X}} \tm \clo{B_{q,\Uc}} \ \ \wedge \ \  t \geq \bar{\tau} \quad \Rightarrow \quad \|\phi(t,x,u)\|_X \leq \eps + \xi(2\eps) + \xi(2\gamma_{\max}(q)).%
\end{align*}
In particular, for $q_0 := 2^{-k_0(\eps,r)+1}$, we have%
\begin{equation}\label{eq_352}
  (x,u) \in \clo{B_{r,X}} \tm \clo{B_{q_0,\Uc}} \ \ \wedge \ \  t \geq \bar{\tau} \quad \Rightarrow \quad \|\phi(t,x,u)\|_X \leq 2\eps + \xi(2\eps),%
\end{equation}
since $\xi(2\gamma_{\max}(q_0)) = \xi(2\gamma_{\max}(2^{-k_0(\eps,r)+1})) \leq \eps$ by definition of $k_0$. Define $\tau(\eps,r) := \max\{\hat{\tau}(\eps,r),\bar{\tau}(\eps,r,q_0)\}$. Combining \eqref{eq_349} and \eqref{eq_352}, we obtain%
\begin{align*}
  &(x,u) \in \clo{B_{r,X}} \tm \clo{B_{r,\Uc}} \ \ \wedge \ \  t \geq \tau(\eps,r) \\
	&\qquad \Rightarrow \quad \|\phi(t,x,u)\|_X \leq 2\eps + \xi(2\eps) + \xi(2\gamma_{\max}(2\|u\|_{\Uc})).%
\end{align*}
As $r \mapsto \xi(2\gamma_{\max}(2r))$ is a $\Kinf$-function, we have proved that $\Sigma$ has the bUAG property which completes the proof. %
\end{proof}

For finite networks, Theorem~\ref{thm_smallgain_iss_semimax} was shown in \cite{Mir19b}. However, in the proof of the infinite-dimensional case, there are essential novelties because the trajectories of an infinite number of subsystems do not necessarily have a uniform speed of convergence. This also resulted in a strengthening of the employed small-gain condition.

In the special case when all interconnection gains $\gamma_{ij}$ are linear, the small-gain condition in our theorem can be formulated more directly in terms of the gains, as the following corollary shows.%

\begin{corollary}[Linear ISS small-gain theorem in semi-maximum formulation]
\label{cor:SGT-sup-linear-gains}
\index{ISS small-gain theorem!linear}
Given an interconnection $(\Sigma,\Uc,\phi)$ of systems $\Sigma_i$ as in Theorem \ref{thm_smallgain_iss_semimax}, additionally to the assumptions (i) and (ii) of this theorem, assume that all gains $\gamma_{ij}$ are linear functions (and hence can be identified with nonnegative real numbers), $\Gamma_{\otimes}$ is well-defined and the following condition holds:%
\begin{equation}\label{eq_sg_hsr_cond}
  \lim_{n \rightarrow \infty} \Bigl( \sup_{j_1,\ldots,j_{n+1}\in I} \gamma_{j_1j_2} \cdots \gamma_{j_{n}j_{n+1}}\Bigr)^{1/n} < 1.%
\end{equation}
Then $\Sigma$ is ISS.
\end{corollary}

\begin{proof}
We only need to show that Assumption (iii) of Theorem \ref{thm_smallgain_iss_semimax} is implied by \eqref{eq_sg_hsr_cond}. The linearity of the gains $\gamma_{ij}$ implies that the operator $\Gamma_{\otimes}$ is homogeneous of degree one and subadditive, see Remark \ref{rem_lg_sg_operator}. Then Proposition \ref{prop:eISS-criterion-homogeneous-systems} and Remark \ref{rem_homogen_spectral_radius} together show that  \eqref{eq_sg_hsr_cond} implies that the system%
\begin{equation*}
  w(k+1) \leq \Gamma_{\otimes}(w(k)) + v(k)%
\end{equation*}
is eISS (according to Definition~\ref{def:eISS-discrete-time-inequalities}), which easily implies the MLIM property for this system. %
\end{proof}

\subsection{Small-gain theorems in summation formulation}
\label{sec:Small-gain theorems in summation formulation}

Now we formulate the small-gain theorems for UGS and ISS in summation formulation.%

\begin{theorem}[UGS small-gain theorem in summation formulation]
\label{thm_ugs_summation_sg}
Let $I$ be a countable index set, $(X_i,\|\cdot\|_{X_i})$, $i\in I$, be normed spaces and $\Sigma_i = (X_i,PC_b(\R_+,X_{\neq i}) \tm \Uc,\bar{\phi}_i)$, $i\in I$ be forward complete control systems. Assume that the interconnection $\Sigma = (X,\Uc,\phi)$ of the systems $\Sigma_i$ is well-defined. Furthermore, let the following assumptions be satisfied:%
\begin{enumerate}
\item[(i)] Each system $\Sigma_i$ is UGS in the sense of Definition \ref{def_subsys_ugs_sum} (summation formulation) with $\sigma_i \in \K$ and nonlinear gains $\gamma_{ij},\gamma_i \in \K \cup \{0\}$.%
\item[(ii)] There exist $\sigma_{\max} \in \Kinf$ and $\gamma_{\max} \in \Kinf$ so that $\sigma_i \leq \sigma_{\max}$ and $\gamma_i \leq \gamma_{\max}$, pointwise for all $i \in I$.%
\item[(iii)]  Assumption \ref{ass_gammasum_welldef}  is satisfied for the operator $\Gamma_{\boxplus}$ defined via the gains $\gamma_{ij}$ from (i) and $\id - \Gamma_{\boxplus}$ has the MBI property.%
\end{enumerate}
Then $\Sigma$ is forward complete and UGS.%
\end{theorem}

\begin{proof}
The proof is exactly the same as for Theorem \ref{thm_ugs_semimax_sg}, with the operator $\Gamma_{\boxplus}$ in place of $\Gamma_{\otimes}$. %
\end{proof}

\begin{theorem}[Nonlinear ISS small-gain theorem in summation formulation]
\label{thm_smallgain_iss_summation}
\index{ISS small-gain theorem!summation formulation}
Let $I$ be a countable index set, $(X_i,\|\cdot\|_{X_i})$, $i\in I$ be normed spaces and $\Sigma_i = (X_i,PC_b(\R_+,X_{\neq i}) \tm \Uc,\bar{\phi}_i)$, $i\in I$ be forward complete control systems. Assume that the interconnection $\Sigma = (X,\Uc,\phi)$ of the systems $\Sigma_i$ is well-defined. Furthermore, let the following assumptions be satisfied:%
\begin{enumerate}
\item[(i)] Each system $\Sigma_i$ is ISS in the sense of Definition \ref{def_subsys_iss_sum} with $\beta_i \in \KL$ and nonlinear gains $\gamma_{ij},\gamma_i \in \K \cup \{0\}$.%
\item[(ii)] There are $\beta_{\max} \in \KL$ and $\gamma_{\max} \in \K$ so that $\beta_i \leq \beta_{\max}$ and $\gamma_i \leq \gamma_{\max}$, pointwise for all $i \in I$.%
\item[(iii)] Assumption \ref{ass_gammasum_welldef} holds and the discrete-time system%
\begin{equation}
\label{eq:Gamma-boxplus-discrete-time}
  w(k+1) \leq \Gamma_{\boxplus}(w(k)) + v(k),%
\end{equation}
with $w(\cdot),v(\cdot)$ taking values in $\ell_{\infty}^+(I)$, has the MLIM property.%
\end{enumerate}
Then $\Sigma$ is ISS.
\end{theorem}

\begin{proof}
The proof is almost completely the same as for Theorem \ref{thm_smallgain_iss_semimax}. The only difference is that instead of interchanging the order of two suprema $\sup_{s \geq t}$ and $\sup_{j \in I}$, we now have to use the estimate
\[
\sup_{s \geq t} \sum_{j \in I} \ldots \leq \sum_{j \in I} \sup_{s \geq t} \ldots,
\]
which is trivially satisfied.
\end{proof}

Again, we formulate a corollary for the case when all gains $\gamma_{ij}$ are linear.%

\begin{corollary}[Linear ISS small-gain theorem in summation formulation]
\label{cor:SGT-sum-linear-gains}
Consider an interconnection $(\Sigma,\Uc,\phi)$ of systems $\Sigma_i$ as in Theorem \ref{thm_smallgain_iss_summation}. Additionally to the assumptions (i) and (ii) of Theorem \ref{thm_smallgain_iss_summation}, assume that all gains $\gamma_{ij}$ are linear functions (and hence can be identified with nonnegative real numbers), the linear operator $\Gamma_{\boxplus}$ is well-defined (thus bounded) and satisfies the spectral radius condition $r(\Gamma_{\boxplus}) < 1$. Then $\Sigma$ is ISS.
\end{corollary}

\begin{proof}
By Proposition \ref{prop:eISS-criterion-linear-systems}, $r(\Gamma_{\boxplus}) < 1$ is equivalent to the MLIM property of the system
\eqref{eq:Gamma-boxplus-discrete-time}, hence to Assumption (iii) of Theorem \ref{thm_smallgain_iss_summation}. %
\end{proof}

\subsection{Example: a linear spatially invariant system}\label{examp: a linear spatially invariant system}

Let us analyze the stability of the spatially invariant infinite network
\begin{equation}\label{eq:linear-spatially-invariant-system}
  \dot{x}_i = ax_{i-1} - x_i + b x_{i+1} + u,\quad i\in\Z,%
\end{equation}
where $a,b>0$ and each $\Sigma_i$ is a scalar system with the state $x_i \in\R$, internal inputs $x_{i-1}$, $x_{i+1}$ and an external input $u$, belonging to the input space $\Uc:=L^\infty(\R_+,\R)$.

Following the general approach in Section \ref{sec:Infinite interconnections}, we define the state space for the interconnection of 
$(\Sigma_i)_{i\in\Z}$ as $X:=\ell_\infty(\Z)$.
Similarly to finite-dimensional ODEs, it is possible to introduce the concept of a mild (Carath\'eodory) solution for the equation 
\eqref{eq:linear-spatially-invariant-system}, for which we refer to \cite{KMS21}. As \eqref{eq:linear-spatially-invariant-system} is linear, it is easy to see that for each initial condition $x_0 \in X$ and for each input $u \in\Uc$ the corresponding mild solution is unique and exists on $\R_+$. We denote it by $\phi(\cdot,x_0,u)$. One can easily check that the triple $\Sigma:=(X,\Uc,\phi)$ defines a well-posed and forward complete interconnection in the sense of this chapter.%

Having a well-posed control system $\Sigma$, we proceed to its stability analysis.%

\begin{proposition}\label{prop:Stability-linear-systems} 
The coupled system \eqref{eq:linear-spatially-invariant-system} is ISS if and only if $a+b<1$.
\end{proposition}

\begin{proof}
\q{$\Rightarrow$}: For any $a,b>0$, the function $y: t \mapsto (e^{(a+b-1)t}x^*)_{i\in\Z}$ is a solution of \eqref{eq:linear-spatially-invariant-system} subject to an initial condition $(x^*)_{i\in\Z}$ and input $u\equiv 0$. This shows that $a+b \geq 1$ implies that the system \eqref{eq:linear-spatially-invariant-system} is not ISS.%

\q{$\Leftarrow$}: By variation of constants, we see that for any $i\in\Z$, treating $x_{i-1}, x_{i+1}$ as external inputs from $L^\infty(\R_+,\R)$, we have the following ISS estimate for the $x_i$-subsystem:%
\begin{align*}
  |x_i(t)| &= \Big|e^{-t}x_i(0) + \int_0^t e^{s-t}[a x_{i-1}(s) + b x_{i+1}(s) + u(s)] ds\Big|\\
  			   &\leq e^{-t}|x_i(0)| + a \|x_{i-1}\|_\infty  + b \|x_{i+1}\|_\infty + \|u\|_\infty,
\end{align*}
for any $t\geq 0$, $x_i(0)\in\R$ and all $x_{i-1}, x_{i+1},u \in L^\infty(\R_+,\R)$.%

This shows that the $x_i$-subsystem is ISS in summation formulation, and the corresponding gain operator is the linear operator $\Gamma:\ell_\infty^+(\Z) \to \ell_\infty^+(\Z)$, acting on $s=(s_i)_{i\in\Z}$ as $\Gamma(s) = (as_{i-1} + b s_{i+1})_{i\in\Z}$. It is easy to see that
\begin{align*}
\|\Gamma\| 
&:= \sup_{\|s\|_{\ell_\infty(\Z)}=1}\|\Gamma s\|_{\ell_\infty(\Z)}
=  \|\Gamma {\unit}\|_{\ell_\infty(\Z)}
= a+b <1,%
\end{align*}
and thus $r(\Gamma)<1$, and the network is ISS by Corollary~\ref{cor:SGT-sum-linear-gains}. 
\end{proof}

\subsection{Example: a nonlinear spatially invariant system}
\label{examp: a nonlinear spatially invariant system}

Consider the infinite interconnection (in the sense of the previous sections)%
\begin{equation}\label{eq:cubic-spatially-invariant-system}
  \dot{x}_i = - x_i^3 + \max\{ax_{i-1}^3,b x_{i+1}^3,u \},\quad i\in\Z,%
\end{equation}
where $a,b>0$. As in Section~\ref{examp: a linear spatially invariant system}, each $\Sigma_i$ is a scalar system with the state $x_i \in\R$, internal inputs $x_{i-1}$, $x_{i+1}$ and an external input $u$, belonging to the input space $\Uc:=L^\infty(\R_+,\R)$. Let the state space for the interconnection $\Sigma$ be $X:=\ell_\infty(\Z)$.%

First, we analyze the well-posedness of the interconnection \eqref{eq:cubic-spatially-invariant-system}. Define for $x = (x_i)_{i\in\Z} \in X$ and $v \in \R$%
\begin{equation*}
  f_i(x,v):= - x_i^3 + \max\{ax_{i-1}^3,b x_{i+1}^3,v \},\quad i \in\Z,%
\end{equation*}
as well as%
\begin{equation*}
  f(x,v):=(f_i(x,v))_{i\in\Z} \in \R^{\Z}.%
\end{equation*}
It holds that%
\begin{equation*}
  |f_i(x,v)| \leq  \|x\|_X^3 + {\max\{a,b\}}\max\{\|x\|^3_X,|v|\},%
\end{equation*}
and thus $f(x,v) \in X$ with $\|f(x,v)\|_X \leq  \|x\|_X^3 + {\max\{a,b\}} \max\{\|x\|^3_X,|v|\}$.%

Furthermore, $f$ is clearly continuous in the second argument. Let us show Lipschitz continuity of $f$ on bounded balls with respect to the first argument. For any $x = (x_i)_{i\in\Z} \in X$, $y = (y_i)_{i\in\Z} \in X$ and any $v \in \R$ we have%
\begin{align*}
\|f(x,u)-&f(y,u)\|_X 
= \sup_{i\in\Z}|f_i(x,u)-f_i(y,u)|\\
&= \sup_{i\in\Z}\big|- x_i^3 + \max\{ax_{i-1}^3,b x_{i+1}^3,v \}  + y_i^3 - \max\{ay_{i-1}^3,b y_{i+1}^3,v \}\big|\\
&\leq \sup_{i\in\Z}\big| x_i^3 - y_i^3\big|
+ \sup_{i\in\Z}\big|\max\{ax_{i-1}^3,b x_{i+1}^3,v \}  - \max\{ay_{i-1}^3,b y_{i+1}^3,v \}\big|.%
\end{align*}
By Birkhoff's inequality $|\max\{a_1,a_2,a_3\} - \max\{b_1,b_2,b_3\}|\leq \sum_{i=1}^3|a_i-b_i|$, which holds for all real $a_i,b_i$,  we obtain
\begin{align*}
\|f(x,u)&-f(y,u)\|_X 
\leq \sup_{i\in\Z}\big| x_i^3 - y_i^3\big|
+ a\sup_{i\in\Z}\big|x_{i-1}^3 -y_{i-1}^3\big|
+ b\sup_{i\in\Z}\big|x_{i+1}^3 -y_{i+1}^3\big|\\
&= (1+a+b)\sup_{i\in\Z}\big| x_i^3 - y_i^3\big|
\leq (1+a+b) \sup_{i\in\Z}\big| x_i - y_i\big| \sup_{i\in\Z}\big| x_i^2 +x_iy_i + y_i^2\big| \\
&\leq (1+a+b)\|x-y\|_X \big(\|x\|^2_X + \|x\|_X\|y\|_X + \|y\|^2_X\big),
\end{align*}
which shows Lipschitz continuity of $f$ with respect to the first argument on the bounded balls in $X$, uniformly with respect to the second argument.%

According to \cite[Theorem~2.4]{AuW96}\footnote{The cited result assumes a global Lipschitz condition and accordingly ensures forward completeness. However, via the retraction method, this result can easily be localized.}, this ensures that the Carath\'eodory solutions of \eqref{eq:cubic-spatially-invariant-system} exist locally, and are unique for any fixed initial condition $x_0\in X$ and external input $u\in\Uc$. We denote the corresponding maximal solution by $\phi(\cdot,x_0,u)$. One can easily check that the triple $\Sigma:=(X,\Uc,\phi)$ defines a well-posed interconnection in the sense of this chapter, and furthermore, $\Sigma$ has BIC property.

We proceed to the stability analysis:
\begin{proposition}
\label{prop:Stability-cubic-systems} 
The coupled system \eqref{eq:cubic-spatially-invariant-system} is ISS if and only if $\max\{a,b\}<1$.
\end{proposition}

\begin{proof}
\q{$\Rightarrow$}: For any $a,b>0$ consider the scalar equation %
\begin{equation*}
  \dot{z} = - (1-\max\{a,b\})z^3,%
\end{equation*}
subject to an initial condition $z(0)=x^*$. The function $t \mapsto (z(t))_{i\in\Z}$ is a solution of \eqref{eq:cubic-spatially-invariant-system} subject to an initial condition $(x^*)_{i\in\Z}$ and input $u\equiv 0$. This shows that for $\max\{a,b\} \geq 1$, the system \eqref{eq:cubic-spatially-invariant-system} is not ISS.

\q{$\Leftarrow$}:
Consider $x_{i-1}$, $x_{i+1}$ and $u$ as inputs to the $x_i$-subsystem of \eqref{eq:cubic-spatially-invariant-system}
and define $q:=\max\{ax_{i-1}^3,b x_{i+1}^3,u \}$.
The derivative of $|x_i(\cdot)|$ along the trajectory satisfies for almost all $t$ the following inequality:
\begin{equation*}
  \frac{d}{dt}|x_i(t)|\leq -|x_i(t)|^3 + q(t)   \leq -|x_i(t)|^3 + \|q\|_\infty.%
\end{equation*}
For any $\varepsilon>0$, if $\|q\|_\infty \leq \frac{1}{1+\varepsilon}|x_i(t)|^3$, we obtain%
\begin{equation*}
  \frac{d}{dt}|x_i(t)|\leq -\frac{\varepsilon}{1+\varepsilon}|x_i(t)|^3.
\end{equation*}
By the comparison principle (Proposition~\ref{prop:ComparisonPrinciple}), applied for $y(t):=|x_i(t)|$, we obtain that there is a certain $\beta\in\KL$
such that for all $t\geq 0$ it holds that
\begin{align*}
|x_i(t)| & \leq \beta(|x_i(0)|,t) +  \big((1+\varepsilon)\|q\|_\infty \big)^{1/3}\\
& = \beta(|x_i(0)|,t) +  \max\{a_1 \|x_{i-1}\|_\infty,b_1 \|x_{i+1}\|_\infty,(1+\varepsilon)^{1/3}\|u\|_\infty^{1/3} \}\\
& \leq \beta(|x_i(0)|,t) +  \max\{a_1 \|x_{i-1}\|_\infty,b_1 \|x_{i+1}\|_\infty\} +(1+\varepsilon)^{1/3}\|u\|_\infty^{1/3},
\end{align*}
where $a_1:=(1+\varepsilon)^{1/3}a^{1/3}$, $b_1:=(1+\varepsilon)^{1/3}b^{1/3}$.%

This shows that the $x_i$-subsystem is ISS in semi-maximum formulation with the corresponding homogeneous {of degree one} gain operator 
$\Gamma:\ell_\infty^+(\Z) \to \ell_\infty^+(\Z)$ given for all $s=(s_i)_{i\in\Z}$ by $\Gamma(s) = (\max\{a_1 s_{i-1}, b_1 s_{i+1}\})_{i\in\Z}$.%

The previous computations are valid for all $\varepsilon>0$. Now pick $\varepsilon>0$ such that $a_1<1$ and $b_1<1$, which is possible
as $a \in (0,1)$ and $b\in(0,1)$. The ISS of the network follows by Corollary~\ref{cor:SGT-sup-linear-gains}. %
\end{proof}

\section{Small-gain conditions}
\label{sec:Small-gain conditions}

Key assumptions in the ISS and UGS small-gain theorems are the monotone limit property, and monotone bounded invertibility property, respectively. In this section, we thoroughly investigate these properties. More precisely, in Section~\ref{sec:A uniform small-gain condition and the MBI property}, we characterize the MBI property in terms of the uniform small-gain condition, in Section~\ref{sec:Non-uniform small-gain conditions}, we relate the uniform small-gain condition to several types of nonuniform small-gain conditions which have already been exploited in the small-gain analysis of finite and infinite networks. In Section~\ref{sec:Finite-dimensional systems}, we derive new relationships between small-gain conditions in the finite-dimensional case.
{Finally, in Section~\ref{sec:Systems with linear gains}, we provide efficient criteria for the MLIM and the MBI property in the case of linear operators and operators of the form $\Gamma_\otimes$ induced by linear gains.}

\subsection{A uniform small-gain condition and the MBI property}
\label{sec:A uniform small-gain condition and the MBI property}

As we have seen in Section~\ref{sec:Small-gain theorems}, the monotone bounded invertibility is a crucial property for the small-gain analysis of finite and infinite networks. The next proposition yields small-gain type criteria for the MBI property. In the context of small-gain theorems in terms of trajectories, derived in this chapter, we are interested primarily in the case of $(X,X^+) = (\ell_\infty(I),\ell^+_\infty(I))$. However, we prove the results in a more general setting, which besides the mathematical appeal, also has important applications to Lyapunov-based small-gain theorems for infinite networks, where other choices for $X$ are useful, see, e.g., \cite{KMS21} where $X=\ell_p$ for finite $p\geq1$.%

\begin{proposition}
\label{prop:criteria-MBI-without-unit}
Let $(X,X^+)$ be an ordered Banach space with a generating cone $X^+$. For every nonlinear operator $A: X^+ \to X^+$, the following conditions are equivalent:
\begin{enumerate}[label=(\roman*)]
\item\label{itm:MBI-criterion-without-unit-1} $\id - A$ satisfies the MBI property.

\index{small-gain condition!uniform}
\item\label{itm:MBI-criterion-without-unit-2} The \emph{uniform small-gain condition} holds: There exists $\eta \in \Kinf$ such that%
\begin{equation}
\label{eq:uSGC-dist-form}
  \dist(A(x) - x,X^+) \geq \eta(\|x\|_X), \quad x \in X^+.%
\end{equation}
\end{enumerate}
\end{proposition}

\begin{proof}
	\ref{itm:MBI-criterion-without-unit-1} $\Rightarrow$ \ref{itm:MBI-criterion-without-unit-2}. Fix $x \in X^+$ and write $a := (A - \id)(x)$. Let $\varepsilon > 0$. We choose $z \in X^+$ such that $\|a-z\|_X \le \dist(a,X^+) + \varepsilon$ and we set $y := a-z$. If the constant $M > 0$ is chosen as in~\eqref{eq:bounded-decomposition}, we can decompose $y$ as $y = u-v$, where $u,v \in X^+$ and $\|u\|_X,\|v\|_X \le M \|y\|_X \le M\dist(a,X^+) + M\varepsilon$. Then we have
\begin{equation*}
  (\id - A)(x) = -a = -y-z = v - (u+z) \le v,
\end{equation*}
so it follows from 
{the MBI property of $\id - A$} that
\begin{equation*}
  \|x\|_X \le \xi(\|v\|_X) \le \xi\big(M \dist(a,X^+)+M\varepsilon\big).
\end{equation*}
Consequently,
\begin{equation*}
  \dist(a,X^+) \ge \frac{1}{M} \xi^{-1}(\|x\|_X) - \varepsilon.
\end{equation*}
Since $\varepsilon$ was arbitrary, this implies~ \ref{itm:MBI-criterion-without-unit-2} with $\eta := \frac{1}{M}\xi^{-1}$.

\ref{itm:MBI-criterion-without-unit-2} $\Rightarrow$ \ref{itm:MBI-criterion-without-unit-1}. Let $v,w \in X^+$ and $(\id-A)(v) \le w$. The vector $z := w + (A-\id)(v)$ is positive, so from  \ref{itm:MBI-criterion-without-unit-2} it follows that%
\begin{equation*}
  \eta(\|v\|_X) \leq \dist\big( (A-\id)(v), X^+ \big) \le \|(A-\id)(v) - z\|_X = \|-w\|_X = \|w\|_X.%
\end{equation*}
Hence, $\|v\|_X \le \eta^{-1}(\|w\|_X)$. %
\end{proof}

\begin{remark}
The uniform small-gain condition in {Proposition~\ref{prop:criteria-MBI-without-unit}\ref{itm:MBI-criterion-without-unit-2}} is a uniform version of the well-known small-gain condition, sometimes also called \emph{no-joint-increase condition}:%
\index{small-gain condition}
\index{no-joint-increase condition}
\begin{equation*}
  A(x) \not\geq x,\quad x \in X^+ \setminus \{0\}.%
\end{equation*}
Indeed, $A(x) \not\geq x$ is equivalent to $A(x) - x \not\geq 0$, which in turn is equivalent to $\dist(A(x) - x, X^+) > 0$.
\end{remark}

\begin{remark}
It is important to point out that the distance to the positive cone, which occurs in the uniform small-gain condition in Proposition~\ref{prop:criteria-MBI-without-unit}, can be explicitly computed on many concrete spaces. Indeed, many important real-valued sequence or function spaces such as $X = \ell_p$ or $X = L^p(\Omega,\mu)$ (for $p \in [1,\infty]$ and a measure space $(\Omega,\mu)$) are not only ordered Banach spaces but so-called \emph{Banach lattices}.%

Now, assume that $(X,X^+)$ is a Banach lattice and let $x \in X$. Then the vectors $x^+ := \frac{|x|+x}{2} \ge 0$ and $x^- := \frac{|x|-x}{2} \ge 0$ are called the \emph{positive} and \emph{negative part} of $x$, respectively; clearly, they satisfy $x^+ - x^- = x$ and $x^+ + x^- = |x|$. If $X$ is a concrete sequence or function space, then $x^-$ is $0$ at all points where $x$ is positive and equal to $-x$ at all points where $x$ is negative.%

In a Banach lattice $(X,X^+)$, we have the formula%
\begin{equation*}
  \dist(x,X^+) = \|x^-\|_X%
\end{equation*}
for each $x \in X$, as can easily be verified.
\end{remark}

If the cone of the ordered Banach space $(X,X^+)$ has nonempty interior, the uniform small-gain condition from Proposition~\ref{prop:criteria-MBI-without-unit} can also be expressed by a condition that involves a fixed interior point of $X^+$.%

\begin{proposition}\label{prop:criteria-MBI-with-unit}
Let $(X,X^+)$ be an ordered Banach space, assume that the cone $X^+$ has nonempty interior, and let $z$ be an interior point of $X^+$. For every nonlinear operator $A: X^+ \to X^+$, the following conditions are equivalent:%
\begin{enumerate}[label=(\roman*)]
\item\label{itm:MBI-criterion-with-unit-1} There is $\eta \in \Kinf$ such that%
\begin{equation}
\label{eq:uSGC-with-unit}
  A(x) \not\geq x - \eta(\|x\|_X)z,\quad  x \in X^+ \setminus \{0\}.%
\end{equation}
\item\label{itm:MBI-criterion-with-unit-2} The uniform small-gain condition from Proposition~\ref{prop:criteria-MBI-without-unit}\ref{itm:MBI-criterion-without-unit-2} holds.
\end{enumerate}
\end{proposition}

\begin{proof}
\ref{itm:MBI-criterion-with-unit-1} $\Rightarrow$ \ref{itm:MBI-criterion-with-unit-2}. Let \ref{itm:MBI-criterion-with-unit-1} hold with some $\eta\in\Kinf$. By Proposition~\ref{prop:order-units}, we can find a number $c>0$ such that for every $y \in X$ we have%
\begin{equation}\label{eq_oldlem3}
  \|y\|_X \le c \quad \Rightarrow \quad y \geq -z.%
\end{equation}
Assume towards a contradiction that~(ii) does not hold. Then \eqref{eq:uSGC-dist-form} fails, in particular, for the function $c\eta$. Thus, we can infer that there is $x \in X^+ \setminus \{0\}$ so that 
\begin{equation*}
  \dist\big( (A-\id)(x), X^+ \big) < c\eta(\|x\|_X).
\end{equation*}
Hence, there exists $y \in X^+$ such that
\begin{equation*}
  \left\| (A-\id)(x) - y \right\|_X \le c\eta(\|x\|_X).
\end{equation*}
Consequently, the vector $\frac{(A-\id)(x) - y}{\eta(\|x\|_X)}$ has norm at most $c$, so it follows from \eqref{eq_oldlem3} that $(A-\id)(x) - y \ge - \eta(\|x\|_X) z$. Thus,%
\begin{equation*}
  (A-\id)(x) \ge -\eta(\|x\|_X)z + y \ge -\eta(\|x\|_X)z,%
\end{equation*}
which shows that \eqref{eq:uSGC-with-unit} fails for the function $\eta$, a contradiction.%

\ref{itm:MBI-criterion-with-unit-2} $\Rightarrow$ \ref{itm:MBI-criterion-with-unit-1}. Let \ref{itm:MBI-criterion-with-unit-2} hold with a certain $\eta\in\Kinf$. We show that \eqref{eq:uSGC-with-unit} holds for the function $\frac{\eta}{2\|z\|_X}$ substituted for $\eta$. Assume towards a contradiction that \eqref{eq:uSGC-with-unit} fails for the function $\frac{\eta}{2\|z\|_X}$. Then there is $x \in X^+ \setminus \{0\}$ such that%
\begin{equation*}
  (A - \id)(x) + \frac{\eta(\|x\|_X)}{2\|z\|_X} z \geq 0.%
\end{equation*}
Hence, it follows that%
\begin{align*}
  \dist\big((A - \id)(x),X^+\big) &\le \Bigl\|(A-\id)(x) - \big((A - \id)(x) + \frac{\eta(\|x\|_X)}{2\|z\|_X} z\big) \Bigr\|_X = \frac{\eta(\|x\|_X)}{2},%
\end{align*}
which shows that~\eqref{eq:uSGC-dist-form} fails for the function $\eta$. %
\end{proof}

A typical example of an ordered Banach space whose cone has nonempty interior is $(X,X^+) = (\ell_{\infty}(I),\ell_{\infty}^+(I))$ for some index set $I$. For instance, the vector ${\unit}$ is an interior point of the positive cone in this space.%

\subsection{Nonuniform small-gain conditions}\label{sec:Non-uniform small-gain conditions}

In Propositions~\ref{prop:criteria-MBI-without-unit} and~\ref{prop:criteria-MBI-with-unit}, we characterized the MBI property in terms of the uniform small-gain condition. In this {subsection}, we recall several further small-gain conditions, which have been used in the literature for the small-gain analysis of finite and infinite networks \cite{DRW07,DRW10,DMS19a}, and relate them to the uniform small-gain condition.%

In this subsection, we always suppose that $(X,X^+) = (\ell_{\infty}(I),\ell_{\infty}^+(I))$ for some nonempty index set $I$ (which is precisely the space in which gain operators act).

\begin{definition}\label{def:SGC}
We say that a nonlinear operator $A:\ell_{\infty}^+(I) \to \ell_{\infty}^+(I)$ satisfies%
\begin{enumerate}[label=(\roman*)]
\item \label{NL-SGC-def-item1} the \emph{small-gain condition} if%
\begin{equation}\label{eq:SGC}
   A(x)\not\geq x,\quad x \in \ell_{\infty}^+(I)\setminus\{0\}.
\end{equation}

\index{small-gain condition!strong}
\item \label{NL-SGC-def-item2} the \emph{strong small-gain condition} if there exists $\rho\in\Kinf$ and a corresponding operator ${D_{\rho}}:\ell_{\infty}^+(I) \to \ell_{\infty}^+(I)$, defined for any $x\in \ell_{\infty}^+(I)$ by%
\begin{equation*}
  {D_{\rho}}(x) := \big((\id + \rho)(x_i)\big)_{i\in I},
\end{equation*}
such that%
\begin{equation}\label{eq:strong-SGC-nonlinear}
  {D_{\rho}}\circ A(x) \not\geq x,\quad x\in \ell_{\infty}^+(I)\setminus\{0\}.%
\end{equation}

\index{small-gain condition!robust}
\item 
\label{NL-SGC-def-item3} the \emph{robust small-gain condition} if there is $\omega\in\Kinf$ with $\omega<\id$ such that for all $i,j \in I$ the operator $A_{i,j}$ given by%
\begin{equation}
\label{eq:A-modified}
  A_{i,j}(x) := A(x) + \omega(x_j) e_i, \quad x \in \ell_{\infty}^+(I)%
\end{equation}
satisfies the small-gain condition \eqref{eq:SGC};
here, $e_i \in \ell_{\infty}(I)$ denotes the $i$-th canonical unit vector.

\index{small-gain condition!robust strong}
\item \label{NL-SGC-def-item4} the \emph{robust strong small-gain condition} if there are $\omega,\rho\in\Kinf$ with $\omega<\id$ such that for all $i,j \in I$ the operator
$A_{i,j}$ defined by \eqref{eq:A-modified} satisfies the strong small-gain condition \eqref{eq:strong-SGC-nonlinear} with the same $\rho$ for all $i,j$. 
\end{enumerate}
\end{definition}

The strong small-gain condition was introduced in \cite{DRW07}, where it was shown that if the gain operator satisfies the strong small-gain condition, then a finite network consisting of ISS systems (defined in a summation formulation) is ISS. The robust strong small-gain condition has been introduced in \cite{DMS19a} in the context of the Lyapunov-based small-gain analysis of infinite networks.%

\begin{remark}\label{rem:Small-gain conditions and cycles} 
For finite networks, so-called \emph{cyclic small-gain conditions} also play an important role, as they help to effectively check the small-gain condition \eqref{eq:SGC} in the case when $A = \Gamma_\otimes$. This is important for the small-gain theorems in the maximum formulation. See \cite{Mir19b} for more discussions on this topic. For infinite networks, the cyclic condition for $\Gamma_\otimes$ is implied by \eqref{eq:SGC}, see \cite[Lemma~4.1]{DMS19a}, but is far too weak for the small-gain analysis. For {max-linear} systems, Remark~\ref{rem_homogen_spectral_radius} and Corollary~\ref{cor:SGT-sup-linear-gains} are reminiscent of the cyclic small-gain conditions.
\end{remark}

Now we give a criterion for the robust strong small-gain condition.%

\begin{proposition}\label{prop:Criterion-robust-strong-SGC} 
A nonlinear operator $A:\ell_{\infty}^+(I) \to \ell_{\infty}^+(I)$ satisfies the robust strong small-gain condition if and only if there are $\omega,\eta\in\Kinf$ and an operator $\vec{\eta}:\ell_{\infty}^+(I) \to \ell_{\infty}^+(I)$, defined by%
\begin{equation}\label{eq:vec-eta}
  \vec{\eta}(x) := (\eta(x_i))_{i\in I}, \quad x\in \ell_{\infty}^+(I),%
\end{equation}
such that for all $k\in I$ it holds that%
\begin{equation}\label{eq:rsSGC-l-infty}
  A(x)\not\geq x - \vec{\eta}(x) - \omega(\|x\|_{\ell_{\infty}(I)}) e_k \mbox{\quad for all\ } x\in \ell_{\infty}^+(I)\setminus\{0\}.%
\end{equation}
\end{proposition}

\begin{proof}
\q{$\Rightarrow$}: Let the robust strong small-gain condition hold with corresponding {$\rho,\omega$ and $D_{\rho}$}. Then for any $x = (x_i)_{i\in I} \in \ell_{\infty}^+(I)\setminus\{0\}$ and any $j,k\in I$, it holds that%
\begin{equation}\label{eq:rsSGC-l-infty-componentwise-1}
  \exists i\in I:\quad \big[{D_{\rho}}\big(A(x) + \omega(x_j) e_k\big)\big]_i = (\id + \rho)\big([A(x) + \omega(x_j) e_k]_i\big) < x_i.%
\end{equation}
As $\rho\in\Kinf$, there is $\eta\in\Kinf$ such that $\id-\eta = (\id + \rho)^{-1}\in\Kinf$, which can be shown as in {Lemma~\ref{lem:uSGC-lemma:id-Kinf}(i)}. Thus, \eqref{eq:rsSGC-l-infty-componentwise-1} is equivalent to%
\begin{equation}\label{eq:rsSGC-l-infty-componentwise}
  \exists i\in I:\quad A(x)_i < x_i - \eta(x_i) - \big[\omega( x_j) e_k\big]_i.%
\end{equation}
As for each $x \in \ell_{\infty}^+(I)$ there is $j\in I$ such that $x_j \geq \frac{1}{2}\|x\|_{\ell_{\infty}(I)}$, the condition \eqref{eq:rsSGC-l-infty-componentwise} with this particular $j$ implies that%
\begin{align*}
\exists i\in I:\quad A(x)_i &< x_i - \eta(x_i) - \Big[\omega\big(\frac{1}{2} \|x\|_{\ell_{\infty}(I)}\big) e_k\Big]_i
													= \Big[x - \vec{\eta}(x) - \omega\big(\frac{1}{2} \|x\|_{\ell_{\infty}(I)}\big) e_k\Big]_i,%
\end{align*}
which is up to the constant the same as \eqref{eq:rsSGC-l-infty}.%

\q{$\Leftarrow$}: Let \eqref{eq:rsSGC-l-infty} hold with a certain $\eta_1\in\Kinf$ and a corresponding $\vec{\eta}_1$. By {Lemma~\ref{lem:uSGC-lemma3}}, one can choose $\eta\in\Kinf$, such that $\eta\leq \eta_1$ and $\id-\eta\in\Kinf$. Then \eqref{eq:rsSGC-l-infty} holds with this $\eta$ and a corresponding $\vec{\eta}$, i.e., for all $k\in I$ we have%
\begin{equation*}
  \exists i\in I:\quad A(x)_i < x_i - \eta(x_i) - \big[\omega( \|x\|_{\ell_{\infty}(I)}) e_k\big]_i.%
\end{equation*}
As $\|x\|_{\ell_{\infty}(I)}\geq x_j$ for any $j\in I$, this implies that for all $j,k\in I$ it holds that%
\begin{equation*}
  \exists i\in I:\quad A(x)_i < x_i - \eta(x_i) - \big[\omega( x_j) e_k\big]_i,%
\end{equation*}
and thus%
\begin{equation*}
  \exists i\in I:\quad \big[ A(x) + \omega( x_j) e_k\big]_i < (\id-\eta)(x_i).%
\end{equation*}
As $\eta\in\Kinf$ satisfies $\id-\eta\in\Kinf$, by {Lemma~\ref{lem:uSGC-lemma:id-Kinf}(i)} there is $\rho\in\Kinf$ such that $(\id-\eta)^{-1} = \id + \rho$, and thus for all $j,k\in I$ property \eqref{eq:rsSGC-l-infty-componentwise-1} holds, which shows that $A$ satisfies the robust strong small-gain condition. %
\end{proof}

Specialized to the strong small-gain condition, Proposition~\ref{prop:Criterion-robust-strong-SGC} reads as follows.%

\begin{corollary}\label{cor:Criterion-strong-SGC} 
A nonlinear operator $A:\ell_{\infty}^+(I) \to \ell_{\infty}^+(I)$ satisfies the strong small-gain condition if and only if there are $\eta\in\Kinf$ and an operator $\vec{\eta}:\ell_{\infty}^+(I) \to \ell_{\infty}^+(I)$, defined via \eqref{eq:vec-eta} such that%
\begin{equation*}
  A(x)\not\geq x - \vec{\eta}(x),\quad  x\in \ell_{\infty}^+(I)\setminus\{0\}.%
\end{equation*}
\end{corollary}

The next proposition shows that the uniform small-gain condition is at least not weaker than the robust strong small-gain condition.
\begin{proposition}\label{prop:uSGC-implies-sSGC} 
Let $A:\ell_{\infty}^+(I)\to \ell_{\infty}^+(I)$ be a nonlinear operator. If $A$ satisfies the uniform small-gain condition, then $A$ satisfies the robust strong small-gain condition.
\end{proposition}

\begin{proof}
As $A$ satisfies the uniform small-gain condition with $\eta$, from the proof of Proposition~\ref{prop:criteria-MBI-with-unit} with $z:= \unit$, we see that for all $x \in \ell_{\infty}^+(I) \setminus \{0\}$%
\begin{equation*}
  A(x) \not\geq x - \frac{1}{2\|{\unit}\|_{\ell_{\infty}(I)}}\eta(\|x\|_{\ell_{\infty}(I)}){\unit} 
	= x - \frac{1}{2}\eta(\|x\|_{\ell_{\infty}(I)}){\unit}.%
\end{equation*}
For any $x\in \ell_{\infty}^+(I)$ and any $k\in I$, it holds that%
\begin{align*}
\frac{1}{2}\eta(\|x\|_{\ell_{\infty}(I)}){\unit}
&=\frac{1}{4}\eta(\|x\|_{\ell_{\infty}(I)}){\unit} + \frac{1}{4}\eta(\|x\|_{\ell_{\infty}(I)}){\unit}\\
&\ge \frac{1}{4}\vec{\eta}(x) + \frac{1}{4}\eta(\|x\|_{\ell_{\infty}(I)})e_k,%
\end{align*}
and by Proposition~\ref{prop:Criterion-robust-strong-SGC}, $A$ satisfies the robust strong small-gain condition. %
\end{proof}

\subsection{The finite-dimensional case}
\label{sec:Finite-dimensional systems}

The case of a finite-dimensional $X$ is significant as it is a key to the stability analysis of finite networks.%

\begin{proposition}
\label{prop:small-gain-condition-n-dim-general}
Assume that $(X,X^+) = (\R^n,\R^n_+)$ for some $n\in\N$, where $\R^n$ is equipped with the maximum norm $\|\cdot\|$ and $\R^n_+$ denotes the standard positive cone in $\R^n$. Further, assume that the operator $A$ is continuous and monotone. Then the following statements are equivalent:%
\begin{enumerate}[label = (\roman*)]
	\item\label{itm:n-dim-small-gain-criterion-1} The system \eqref{eq_monotone_system} has the MLIM property.%
	\item\label{itm:n-dim-small-gain-criterion-2} The operator $\id - A$ has the MBI property.%
	\item\label{itm:n-dim-small-gain-criterion-3} The uniform small-gain condition holds: There is an $\eta \in \Kinf$ such that $\dist(A(x) - x,X^+) \geq \eta(\|x\|)$ for all $x \in X^+$.%
	\item\label{itm:n-dim-small-gain-criterion-4} There is an $\eta \in \Kinf$ such that%
\begin{equation*}
  A(x) \not\geq x - \eta(\|x\|){\unit}, \quad x \in X^+ \setminus \{0\}.%
\end{equation*}
\end{enumerate}

Additionally, if $A$ is either $\Gamma_\boxplus$ or $\Gamma_\otimes$, then the above conditions are equivalent to:

\begin{enumerate}[label = (\roman*), start = 5]
	\item\label{itm:n-dim-small-gain-criterion-5} $A$ satisfies the robust strong small-gain condition.
	\item\label{itm:n-dim-small-gain-criterion-6} $A$ satisfies the strong small-gain condition.
\end{enumerate}
\end{proposition}

\begin{proof}
\ref{itm:n-dim-small-gain-criterion-1} $\Rightarrow$ \ref{itm:n-dim-small-gain-criterion-2}. Follows from Proposition~\ref{prop_mlim_implies_mbip}.

\ref{itm:n-dim-small-gain-criterion-2} $\Iff$ \ref{itm:n-dim-small-gain-criterion-3} $\Iff$ \ref{itm:n-dim-small-gain-criterion-4}. 
Follows from Propositions~\ref{prop:criteria-MBI-with-unit},~\ref{prop:criteria-MBI-without-unit}.

\ref{itm:n-dim-small-gain-criterion-2} $\Rightarrow$ \ref{itm:n-dim-small-gain-criterion-1}. This follows from Proposition \ref{prop_compact_operators} since the cone $\R^n_+$ has the Levi property.

\ref{itm:n-dim-small-gain-criterion-4} $\Rightarrow$ \ref{itm:n-dim-small-gain-criterion-5}. Follows by Proposition~\ref{prop:uSGC-implies-sSGC}. 

\ref{itm:n-dim-small-gain-criterion-5} $\Rightarrow$ \ref{itm:n-dim-small-gain-criterion-6}. Clear.

\ref{itm:n-dim-small-gain-criterion-6} $\Rightarrow$ \ref{itm:n-dim-small-gain-criterion-2}. Follows by \cite[Theorem~6.1]{Rue10}. %
\end{proof}

\begin{remark}
\label{rem:MAFs} 
The class of operators for which the equivalence between \ref{itm:n-dim-small-gain-criterion-1}--\ref{itm:n-dim-small-gain-criterion-4} and \ref{itm:n-dim-small-gain-criterion-5}, \ref{itm:n-dim-small-gain-criterion-6} can be shown, can be made considerably larger using the monotone aggregation functions formalism, see \cite[Theorem~6.1]{Rue10}. However, the proof of this implication in \cite[Lemma~13]{DRW07} uses more structure of the gain operator than merely monotonicity. Thus, the question if this implication is valid for general monotone $A$ is still open.
\end{remark}

\section[Exponential ISS of homogeneous subadditive systems]{Exponential ISS of linear and homogeneous of degree one subadditive discrete-time systems}
\label{sec:Exponential ISS of discrete-time systems}

Here we characterize exponential ISS for homogeneous of degree one and subadditive operators.

\begin{definition}
\label{def:eISS-discrete-time-inequalities}
Let $(X,X^+)$ be an ordered Banach space. 
The system \eqref{eq_monotone_system} is called \emph{exponentially input-to-state stable (eISS)} if there are $M\geq 1$, $a\in(0,1)$ and $\gamma \in \Kinf$ such that for every $u \in \ell_{\infty}(\Z_+,X^+)$ and any solution $x(\cdot) = (x(k))_{k\in\Z_+}$ of \eqref{eq_monotone_system} it holds that%
\begin{equation}
\label{eq:eISS-inequalities}
  \|x(k)\|_X \leq M\|x(0)\|_X a^k + \gamma(\|u\|_{\infty}),\quad k \in \Z_+.%
\end{equation}
\end{definition}

\begin{proposition}\label{prop:eISS-criterion-homogeneous-systems}
Let $(X,X^+)$ be an ordered Banach space with a generating and normal cone $X^+$. Consider the system \eqref{eq_monotone_system} and assume that the operator $A:X^+\to X^+$ is monotone and satisfies the following properties:%
\begin{enumerate}[label = (\roman*)]
\index{operator!homogeneous of degree one}
\item\label{itm:eISS-criterion-homog-1} $A$ is homogeneous of degree one, i.e., $A(rx) = rA(x)$ for all $x\in X^+$ and $r \geq 0$.%
\index{operator!subadditive}
\item\label{itm:eISS-criterion-homog-2} $A$ is subadditive, i.e., $A(x + y) \leq A(x) + A(y)$ for all $x,y\in X^+$.%
\item\label{itm:eISS-criterion-homog-3}  $A$ satisfies%
\begin{equation*}
  C := \sup_{x\in X^+ ,\; \|x\|_X=1}\|A(x)\|_X < \infty.%
\end{equation*}
\end{enumerate}
Then $A$ is globally Lipschitz continuous, and the following statements are equivalent:
\begin{enumerate}[label = (\alph*)]
	\item\label{itm:eISS-criterion-equiv-1} The system \eqref{eq_monotone_system} is eISS. 

	\item\label{itm:eISS-criterion-equiv-2}  It holds that 
	\begin{equation}\label{eq_def_spectralradius}
  r(A) := \lim_{n \rightarrow \infty} \sup_{x \in X^+,\; \|x\|_X = 1}\|A^n(x)\|_X^{1/n} < 1.%
\end{equation}
\item\label{itm:eISS-criterion-equiv-3} There is a globally Lipschitz map $V:X^+\to\R_+$ and $L_1,L_2,\psi>0$, $\eta>1$, such that 
\begin{equation}
\label{eq:Sandwich}
  L_1\|x\|_X \leq V(x) \leq L_2 \|x\|_X,\quad x \in X^+,
\end{equation}
and for any $u \in \ell_{\infty}(\Z_+,X^+)$, and any solution of \eqref{eq_monotone_system} it holds that 
\begin{equation}
\label{eq:Dissipative-inequality}
  V(x(k+1)) \leq \eta^{-1} V(x(k)) + \psi \|u\|_{\infty},\quad k \geq 0.%
\end{equation}
\end{enumerate}
\end{proposition}

\begin{proof}
First, we show that $A$ is Lipschitz continuous.
Pick any $x,y \in X^+$. As $X^+$ is generating, there are $M>0$ (which does not depend on $x,y$) and $a,b\in X^+$ such that $x-y = a-b$ and 
$\|a\|_X\leq M\|x-y\|_X$, $\|b\|_X\leq M\|x-y\|_X$.
Hence, for all $x,y\in X^+$ we have 
\begin{align*}
A(x) - A(y)
 &=  A(x-y + y)  -A(y)  =  A(a-b + y)  -A(y)\\
 &\leq  A(a + y)  -A(y) \leq  A(a) + A(y)   - A(y) = A(a).
\label{eq:A(x)-A(y)-upper}
\end{align*}
Analogously, we obtain for all $x,y\in X^+$ that $A(x) - A(y) \geq - A(b)$. As $X^+$ is normal, due to \cite[Theorem~2.38]{AlT07}, there is $c>0$, depending only on $(X,X^+)$, such that 
\begin{align*}
\|A(x) - & A(y)\|_{X} 
\leq c\max\{\|A(a)\|_{X},\|A(b)\|_{X}\}\\
&\leq c\max\{\|a\|_X\|A(a/\|a\|_X)\|_{X},\|b\|_X\|A(b/\|b\|_X)\|_{X}\}
\leq cCM \|x-y\|_{X}.
\end{align*}

\ref{itm:eISS-criterion-equiv-1} $\Rightarrow$ \ref{itm:eISS-criterion-equiv-2}. If \eqref{eq_monotone_system} is eISS, then for $u\equiv 0$, any $x \in X^+$ and for the solution $x(k+1) = A(x(k))$ of  \eqref{eq_monotone_system}, the inequality \eqref{eq:eISS-inequalities} implies that $\|A^n(x)\|_X\leq Ma^n\|x\|_X$ for all $n\in\Z_+$.
Hence, $\sup_{x \in X^+,\; \|x\|_X = 1}\|A^n(x)\|_X^{1/n} \leq M^{1/n}a \to a$ as $n\to\infty$ and thus $r(A) \leq a<1$.%

\ref{itm:eISS-criterion-equiv-2} $\Rightarrow$ \ref{itm:eISS-criterion-equiv-3}. 
 From the assumptions \ref{itm:eISS-criterion-homog-1} and \ref{itm:eISS-criterion-homog-3} together it follows that%
\begin{equation}
\label{eq_abound}
  \|A(x)\|_X = \|x\|_X \|A(x/\|x\|_X)\|_X \leq C \|x\|_X,\quad x \in X^+ \setminus \{0\}.%
\end{equation}
Consider the sequence%
\begin{equation*}
  a_n := \sup_{x \in X^+,\; \|x\|_X = 1}\|A^n(x)\|_X,\quad n \in \Z_+.%
\end{equation*}
This sequence is submultiplicative, as for all $m,n\in\Z_+$ it holds that
\begin{align*}
  a_{n+m} &= \sup_{x \in X^+,\; \|x\|_X = 1}\|A^m(A^n(x))\|_X = \sup_{x \in X^+ ,\; \|x\|_X = 1} \|A^n(x)\|_X \Bigl\|A^m\Bigl(\frac{A^n(x)}{\|A^n(x)\|_X}\Bigr)\Bigr\|_X \\
	        &\leq \sup_{x \in X^+ ,\; \|x\|_X = 1} \|A^n(x)\|_X \cdot \sup_{x \in X^+ ,\; \|x\|_X=1} \|A^m(x)\|_X = a_n \cdot a_m.%
\end{align*}
By a submultiplicative version of Fekete's subadditive lemma, 
\[
\lim_{n\to\infty}a_n^{\frac{1}{n}} = \inf_{n\to\infty}a_n^{\frac{1}{n}}\leq a_1<\infty, 
\]
and thus the limit in \eqref{eq_def_spectralradius} exists.%

We fix $\eta > 1$ such that $\eta r(A) < 1$ and define a function $V:X^+ \rightarrow \R_+$ by%
\begin{equation*}
  V(x) := \sup_{n \in \Z_+} \eta^n \|A^n(x)\|_X,\quad x \in X^+.%
\end{equation*}
Setting $n:=0$ in the supremum, we see that $\|x\|_X\leq V(x)$ for all $x \in X^+$. Since $r(A) < \eta^{-1}$, there exists $N \in \N$ so that%
\begin{equation*}
  \sup_{x \in X^+ ,\; \|x\|_X = 1}\|A^n(x)\|_X \leq \eta^{-n},\quad n \geq N.%
\end{equation*}
By homogeneity of degree one of $A$, this implies%
\begin{equation*}
  \eta^n \|A^n(x)\|_X = \|x\|_X \eta^n \|A^n(\tfrac{x}{\|x\|_X})\|_X \leq \|x\|_X,\quad n \geq N,\ x \in X^+ \setminus \{0\}.%
\end{equation*}
By \eqref{eq_abound}, we have $\|A(x)\|_X \leq C\|x\|_X$ for all $x \in X^+$. Due to the homogeneity of $A$%
\begin{equation*}
  \|A^n(x)\|_X = \|A^{n-1}(x)\|_X \|A(\frac{A^{n-1}(x)}{\|A^{n-1}(x)\|_X})\|_X \leq C\|A^{n-1}(x)\|_X%
\end{equation*}
for all $x \in X^+$, and by induction $\|A^n(x)\|_X \leq C^n\|x\|_X$ for all $x \in X^+$.%

Since $\eta^0 \|A^0(x)\|_X=\|x\|_X$, with $\psi := \max_{0 \leq n < N} (\eta C)^n$ we have%
\begin{equation}\label{eq:LF-estimate-above}
  V(x) = \sup_{n \in \Z_+} \eta^n \|A^n(x)\|_X = \sup_{0\leq n < N} \eta^n \|A^n(x)\| \leq \psi \|x\|_X.%
\end{equation}
Also, observe that%
\begin{equation*}
  V(A(x)) = \sup_{n \in \Z_+} \eta^n \|A^{n+1}(x)\|_X = \eta^{-1} \sup_{n \in \Z_+} \eta^{n+1} \|A^{n+1}(x)\|_X \leq \eta^{-1} V(x).%
\end{equation*}
As $A$ is monotone and subadditive, it holds by induction for all $n\in\N$ that 
\begin{equation*}
  A^n(x+y) = A^{n-1}(A(x+y)) \leq A^{n-1}(A(x)+A(y)) \leq A^n(x) + A^n(y),%
\end{equation*}
that is, $A^n$ are subadditive as well.%

We can assume without loss of generality that the norm $\|\cdot\|_X$ is monotone, i.e., $0 \leq x \leq y$ implies $\|x\|_X \leq \|y\|_X$ for any $x,y\in X^+$. Otherwise, we choose an equivalent norm with this property and note that eISS in one norm implies eISS in any other equivalent norm and that the spectral radius does not depend on the choice of an equivalent norm.%

Together with the subadditivity of $A^n$, $n\in\N$, this implies for all $x,y \in X^+$ that%
\begin{align}\label{eq_c_subadditive}
\begin{split}
  V(x + y) &= \sup_{n \in \Z_+} \eta^n \| A^n(x + y) \|_X 
	\leq \sup_{n \in \Z_+} \eta^n \|A^n(x) + A^n(y)\|_X \\
	&\leq \sup_{n \in \Z_+} \eta^n (\|A^n(x)\|_X + \|A^n(y)\|_X) \leq V(x) + V(y),
\end{split}
\end{align}
and hence $V$ is subadditive as well. Now consider a sequence $x(\cdot)$ in $X^+$ such that%
\begin{equation}
\label{eq:Discrate-time-sys-aux}
  x(k+1) = A(x(k)) + u(k),\quad k \in \Z_+.%
\end{equation}
It then follows that%
\begin{align*}
  V(x(k+1)) &= V(A(x(k)) + u(k)) \leq V(A(x(k))) + V(u(k)) \leq \frac{1}{\eta} V(x(k)) + V(u(k)).%
\end{align*}
By \eqref{eq:LF-estimate-above}, we obtain%
\begin{equation*}
  V(x(k+1)) \leq \eta^{-1} V(x(k)) + \psi \|u\|_{\infty},\quad k \geq 0.%
\end{equation*}

As $V$ is homogeneous of degree one, subadditive, and monotone,
 $V$ is Lipschitz continuous by the argumentation at the beginning of the proof. 

\ref{itm:eISS-criterion-equiv-3} $\Rightarrow$ \ref{itm:eISS-criterion-equiv-1}. This standard argument it omitted.
%
\end{proof}

\section{Systems with linear gains}
\label{sec:Systems with linear gains}

Here we show that in the case of linear and sup-linear gain operators, the MBI and MLIM properties are equivalent and can be characterized via the spectral condition.

For linear systems, we obtain the following result, which we use to formulate an efficient small-gain theorem in summation formulation, see Corollary \ref{cor:SGT-sum-linear-gains}.
\begin{proposition}
\label{prop:eISS-criterion-linear-systems}
Let $(X,X^+)$ be an ordered Banach space with a generating and normal cone $X^+$. 
{Let the operator $A:X^+ \to X^+$ be the restriction to $X^+$ of a positive linear operator on $X$.} Then the following statements are equivalent:%
\begin{enumerate}[label = (\roman*)]
\item\label{itm:linear-LIM-criterion-1} The system \eqref{eq_monotone_system} is exponentially ISS.
\item\label{itm:linear-LIM-criterion-3} The system \eqref{eq_monotone_system} satisfies the MLIM property.%
\item\label{itm:linear-LIM-criterion-4} The operator $\id-A$ satisfies the MBI property.%
\item\label{itm:linear-LIM-criterion-5} The spectral radius of $A$ satisfies $r(A) < 1$.%
\end{enumerate}
\end{proposition}

\begin{proof}
The implication ``\ref{itm:linear-LIM-criterion-1} $\Rightarrow$ \ref{itm:linear-LIM-criterion-3}'' is trivial. 
By Proposition \ref{prop_mlim_implies_mbip}, \ref{itm:linear-LIM-criterion-3} implies \ref{itm:linear-LIM-criterion-4}. 

\ref{itm:linear-LIM-criterion-4} $\Rightarrow$ \ref{itm:linear-LIM-criterion-5}. It is easy to check that if $A$ is homogeneous of degree one  and $\id - A$ satisfies the MBI property with a certain $\xi\in\Kinf$, then $\id - A$ satisfies the MBI property with $r\mapsto \xi(1)r$ instead of $\xi$. 
The application of Theorem~\ref{thm:stability-for-pos-ops} shows \ref{itm:linear-LIM-criterion-5}.

\ref{itm:linear-LIM-criterion-5} $\Rightarrow$ \ref{itm:linear-LIM-criterion-1}. Follows from Proposition~\ref{prop:eISS-criterion-homogeneous-systems}.
\end{proof}

For sup-linear systems, MBI is again equivalent to eISS, and the following holds:
\begin{proposition}
\label{prop:join-morphism-eISS-criterion}
Assume that the gains $\gamma_{ij}$, $(i,j) \in I^2$, are all linear and that the associated gain operator $\Gamma_{\otimes}$ is well-defined. Then the following statements are equivalent:%
\begin{enumerate}[label = (\roman*)]
\item\label{itm:join-morphism-eISS-criterion-1} The operator $\id - \Gamma_{\otimes}$ satisfies the MBI property.%
\item\label{itm:join-morphism-eISS-criterion-2} There are $\lambda\in(0,1)$ and $s_0 \in \intt(\ell_{\infty}^+(I))$ such that 
\begin{eqnarray}
\Gamma_{\otimes}(s_0) \leq \lambda s_0.
\label{eq:Point-of-strict-decay}
\end{eqnarray}

\item\label{itm:join-morphism-eISS-criterion-3} The spectral radius of $\Gamma_{\otimes}:\ell_{\infty}^+(I) \to \ell_{\infty}^+(I)$ satisfies%
\begin{equation*}
  r(\Gamma_{\otimes}) = \lim_{n\rightarrow\infty} \sup_{s \in \ell_{\infty}^+(I) ,\; \|s\|_{\ell_{\infty}} = 1} \|\Gamma_{\otimes}^n(s)\|_{\ell_{\infty}(I)}^{1/n}
	= \lim_{n\rightarrow\infty} \|\Gamma_{\otimes}^n({\unit})\|_{\ell_{\infty}(I)}^{1/n} < 1.
\end{equation*}
\item\label{itm:join-morphism-eISS-criterion-4} The system \eqref{eq_monotone_system} with $A = \Gamma_{\otimes}$ is eISS.%
\item\label{itm:join-morphism-eISS-criterion-5} The system \eqref{eq_monotone_system} with $A = \Gamma_{\otimes}$ has the MLIM property.%
\end{enumerate}
\end{proposition}

\begin{proof}
By Proposition~\ref{prop:eISS-criterion-homogeneous-systems}, \ref{itm:join-morphism-eISS-criterion-3} is equivalent to \ref{itm:join-morphism-eISS-criterion-4}. Clearly, \ref{itm:join-morphism-eISS-criterion-4} implies \ref{itm:join-morphism-eISS-criterion-5}. By Proposition~\ref{prop_mlim_implies_mbip}, \ref{itm:join-morphism-eISS-criterion-5} implies \ref{itm:join-morphism-eISS-criterion-1}.%

\ref{itm:join-morphism-eISS-criterion-1} $\Rightarrow$ \ref{itm:join-morphism-eISS-criterion-2}. By Proposition \ref{prop:criteria-MBI-without-unit}, the MBI property of $\id - \Gamma_{\otimes}$ is equivalent to the uniform small-gain condition. Then Proposition \ref{prop:criteria-MBI-with-unit} shows that%
\begin{equation*}
  \Gamma_{\otimes}(s) \not\geq s - \eta(\|s\|_{\ell_{\infty}}){\unit},\quad s \in \ell_{\infty}^+(I) \setminus \{0\}%
\end{equation*}
for some $\eta \in \Kinf$. In particular,%
\begin{equation*}
  \Gamma_{\otimes}\Bigl(\frac{s}{\|s\|_{\ell_{\infty}}}\Bigr) \not\geq \frac{s}{\|s\|_{\ell_{\infty}}} - \eta(1){\unit},\quad s \in \ell_{\infty}^+(I) \setminus \{0\}.%
\end{equation*}
Multiplying this inequality by $\|s\|_{\ell_{\infty}}$, putting $\eta := \eta(1)$ and using the homogeneity of degree one of $\Gamma_{\otimes}$ yields%
\begin{equation*}
  \Gamma_{\otimes}(s) \not\geq s - \eta\|s\|_{\ell_{\infty}} {\unit},\quad s \in \ell_{\infty}^+(I) \setminus \{0\}.%
\end{equation*}
Then for any $s \in \ell_{\infty}^+(I)$ we have%
\begin{align*}
  (1 + \eps)\Gamma_{\otimes}(s) &\not\geq (1 + \eps)(s - \eta \|s\|_{\ell_{\infty}} {\unit}) = s + \eps s - (1 + \eps) \eta \|s\|_{\ell_{\infty}} {\unit}.
\end{align*}
As $s + \eps s - (1 + \eps) \eta \|s\|_{\ell_{\infty}} \unit \leq s - [(1 + \eps)\eta - \eps] \|s\|_{\ell_{\infty}}{\unit}$,
we have
\begin{align*}
  (1 + \eps)\Gamma_{\otimes}(s) &\not\geq s - [(1 + \eps)\eta - \eps] \|s\|_{\ell_{\infty}}{\unit}.
\end{align*}
Choosing $\varepsilon>0$ small enough, Proposition \ref{prop:uSGC-implies-sSGC} implies that $(1 + \eps)\Gamma_{\otimes}$ satisfies the robust strong small-gain condition. By Lemma \ref{lem_Q}, the operator%
\begin{equation*}
  Q^{\eps}(s) := \sup_{k \in \Z_+} (1 + \eps)^k \Gamma_{\otimes}^k(s),\quad s \in \ell_{\infty}^+(I)%
\end{equation*}
is well-defined and satisfies
\begin{equation*}
  \Gamma_{\otimes}(Q^{\eps}(s)) \leq \frac{1}{1 + \eps}Q^{\eps}(s),\quad s \in \ell_{\infty}^+(I).
\end{equation*}
In particular, this holds for $s = {\unit}$. Since $s_0 := Q^{\eps}({\unit}) \geq {\unit}$, we have $s_0 \in \intt(\ell_{\infty}^+(I))$. 

\ref{itm:join-morphism-eISS-criterion-2} $\Rightarrow$ \ref{itm:join-morphism-eISS-criterion-3}. By monotonicity and homogeneity of degree one of $\Gamma_{\otimes}$, we have%
\begin{equation*}
  \Gamma_{\otimes}^k(s_0) \leq \lambda^k s_0 \mbox{\quad for all\ } k \geq 1.%
\end{equation*}
There exists $n \in \N$ such that any $s \in \ell_{\infty}^+(I)$ with $\|s\|_{\ell_{\infty}} = 1$ satisfies $s \leq ns_0$. Hence,%
\begin{equation*}
  \Gamma_{\otimes}^k(s) \leq \Gamma_{\otimes}^k(ns_0) = n \Gamma_{\otimes}^k(s_0) \leq n \lambda^k s_0,\quad k \geq 1,\ \|s\|_{\ell_{\infty}} = 1.%
\end{equation*}
This implies $r(\Gamma_{\otimes}) \leq \lambda < 1$, which completes the proof. %
\end{proof}

\begin{remark}\label{rem:Paths-of-strict-decay} 
The special form of the operator $\Gamma_\otimes$ is used in Proposition~\ref{prop:join-morphism-eISS-criterion} only for the proof of 
the implication \ref{itm:join-morphism-eISS-criterion-1} $\Rightarrow$ \ref{itm:join-morphism-eISS-criterion-2}. The remaining implications are valid for considerably more general types of operators. Note that if $s_0$ is as in item \ref{itm:join-morphism-eISS-criterion-2}, then $ts_0$ also satisfies all conditions in item \ref{itm:join-morphism-eISS-criterion-2}, for any $t>0$. Thus, we can construct a path of strict decay $t\mapsto t s_0$ for the gain operator $\Gamma_\otimes$, which is an important ingredient for the proof of the Lyapunov-based ISS small-gain theorem, see \cite{DRW06b}.
\end{remark}

\section{Systems, governed by a max-form gain operator}
\label{sec:Systems, governed by a max-form gain operator}

Here, we study the properties of the operator $\Gamma_\otimes$ and its strong transitive closure. These results strengthen the corresponding results in \cite[Section~4]{DMS19a} and are motivated by them.
We use these results to characterize the MBI and MLIM properties for the operator $\Gamma_{\otimes}$ with linear gains in 
Proposition~\ref{prop:join-morphism-eISS-criterion}. However, the developments of this section are also useful for the construction of paths of strict decay for the nonlinear operator $\Gamma_{\otimes}$, which is essential for nonlinear Lyapunov-based small-gain theorems.

The powers of the operator $\Gamma_\otimes$ have a particularly simple representation:%

\begin{lemma}
\label{lem:Gamma_otimes_formula} 
For any $n\in\N$ and any $x \in \ell_\infty^+(I)$, it holds that
\begin{equation}\label{eq:Potenzen-join-morphism}
  \Gamma^n_\otimes(x) = \Big(\sup_{j_2,\ldots,j_{n+1}\in I}\gamma_{i j_2}\circ \cdots\circ \gamma_{j_{n}j_{n+1}} (x_{j_{n+1}})\Big)_{i\in I}.
\end{equation}
\end{lemma}

\begin{proof}
For $n=1$ the claim is clear. Let the claim hold for a certain $n\in\N$. 
Then
\begin{eqnarray*}
  \Gamma^{n+1}_\otimes(x) 
	&=& \Gamma_\otimes(\Gamma^{n}_\otimes(x)) 
	= \big(\sup_{j_2 \in I}\gamma_{ij_2} \circ \sup_{j_3,\ldots,j_{n+2}\in I}\gamma_{j_2 j_3}\circ \cdots\circ \gamma_{j_{n+1}j_{n+2}} (x_{j_{n+2}})\big)_{i \in I}.
\end{eqnarray*}
As $(\gamma_{ij})\subset\Kinf$, \quad
$  \Gamma^{n+1}_\otimes(x) 
	= \big(\sup\limits_{j_2,\ldots,j_{n+2}\in I}\vspace{-4mm}\gamma_{ij_2} \circ \gamma_{j_2 j_3}\circ \cdots\circ \gamma_{j_{n+1}j_{n+2}} (x_{j_{n+2}})\big)_{i \in I}.$
\end{proof}

\begin{remark}\label{rem_homogen_spectral_radius}
Let $\Gamma_{\otimes}$ be homogeneous, which is the case if and only if all the gains $\gamma_{ij}$ are linear. Then using monotonicity of $\Gamma_{\otimes}$, and the formula \eqref{eq:Potenzen-join-morphism}, we obtain the following expression for the spectral radius of $\Gamma_{\otimes}$:
\begin{eqnarray}
\label{eq:Spectral-radius-join-morphism}
  r(\Gamma_{\otimes}) 
	&=& \lim_{n \rightarrow \infty} \sup_{x \in \ell_\infty^+(I),\; \|x\|_{\ell_\infty(I)} = 1}\|\Gamma_{\otimes}^n(x)\|_{\ell_\infty(I)}^{1/n} = \lim_{n \rightarrow \infty} \|\Gamma_{\otimes}^n({\unit})\|_{\ell_\infty(I)}^{1/n}\nonumber\\	
       &=& \lim_{n \rightarrow \infty} \Bigl( \sup_{j_1,\ldots,j_{n+1}\in I} \gamma_{j_1j_2} \cdots \gamma_{j_{n}j_{n+1}}\Bigr)^{1/n}.%
\end{eqnarray}
\end{remark}

We need the following lemma, which we adopt from \cite[Lemma~4.1]{DMS19a}.
\begin{lemma}\label{lem_cycle_contr}
Assume that $\Gamma_\otimes$ satisfies the small-gain condition. Then all cycles built from the gains $\gamma_{ij}$ are contractions. That is,%
\begin{equation*}
  \gamma_{i_1i_2} \circ \cdots \circ \gamma_{i_{k-1}i_k}(r) < r%
\end{equation*}
for all $r>0$ if $i_1,\ldots,i_k$ is an arbitrary path with $i_1 = i_k$.
\end{lemma}

\begin{proof}
Assume that the assertion is not satisfied for some cycle:%
\begin{equation*}
  \gamma_{i_1i_2} \circ \cdots \circ \gamma_{i_{k-1}i_k}(r) \geq r,\quad i_1 = i_k,\ r > 0.%
\end{equation*}
Then define the vector
\begin{align*}
  s &:= re_{i_1} + e_{i_2}\gamma_{i_2i_3} \circ \cdots \circ \gamma_{i_{k-1}i_k}(r) \\
	 &\qquad + e_{i_3}\gamma_{i_3i_4} \circ \cdots \circ \gamma_{i_{k-1}i_k}(r) + \ldots + e_{i_{k-1}}\gamma_{i_{k-1}i_k}(r).%
\end{align*}
We check that $\Gamma_\otimes(s) \geq s$ to obtain a contradiction. Indeed, for any nonzero component of $s$ indexed by $i_{\nu}$ we have (recall that $i_k = i_1$)%
\begin{align*}
  \Gamma_{\otimes,i_{\nu}}(s) &= \sup_{j\in\N}\gamma_{i_{\nu}j}(s_j) \geq \gamma_{i_{\nu}i_{\nu+1}}(s_{i_{\nu+1}}) 
	= \gamma_{i_{\nu}i_{\nu+1}} \circ \cdots \circ \gamma_{i_{k-1}i_k}(r).%
\end{align*}
The last expression equals $s_{i_{\nu}}$ if $\nu \in \{2,\ldots,k-1\}$ and is greater than or equal to $r = s_{i_1}$ if $\nu = 1$. This implies $\Gamma_\otimes(s) \geq s$ in contradiction to the small-gain condition.
\end{proof}

The implication shown in Lemma~\ref{lem_cycle_contr} cannot be reversed in general. However, if the operator $\Gamma$ is block-diagonal with finite-dimensional components, then the implication can be ``upgraded'' to the following criterion.%

\begin{proposition}
\label{prop:SGC_block-diagonal-operators} 
Assume that the gain matrix $\Gamma =(\gamma_{ij})_{i,j\in\N}$, $\gamma_{ij}\in\Kinf\cup\{0\}$, is a block-diagonal matrix of the form 
$\Gamma = \diag(\Gamma_1, \Gamma_2, \ldots)$, where all $\Gamma_j:\R^{n_j}_+\to \R^{n_j}_+$ and $n_j <\infty$ for all $j\in\N$.
We write for the corresponding gain operators $\Gamma_\otimes = \diag(\Gamma_{1,\otimes}, \Gamma_{2,\otimes}, \ldots)$.
Then the following statements are equivalent:
\begin{itemize}
	\item[(i)]  $\Gamma_\otimes:\ell_\infty^+\to\ell_\infty^+$ satisfies the small-gain condition. 
	\item[(ii)] For all $i\in\N$, the operator $\Gamma_{i,\otimes}$ satisfies the small-gain condition.
	\item[(iii)] All cycles built from the gains of $\Gamma$ are contractions (as defined in Lemma~\ref{lem_cycle_contr}).
	\item[(iv)] All cycles built from the gains of $\Gamma_i$, $i\in\N$, are contractions.
\end{itemize}
\end{proposition}

\begin{proof}
(i) $\Rightarrow$ (ii). Assume that $\Gamma_\otimes$ satisfies the small-gain condition.
Pick any $i\in\N$ and consider the vector $s = (s_1,s_2,\ldots) \in\ell_\infty^+$ with $s_j=0$ for $j\neq i$ and $s_i\neq 0$.

For this operator, we have $(\Gamma_\otimes s)_i = \Gamma_{i,\otimes} s_i$, and $(\Gamma_\otimes s)_j=0$ for $j\neq i$.
As $\Gamma_\otimes s \not\geq s$ for all $s_i\neq 0$, it follows  that $\Gamma_{i,\otimes} s_i \not\geq s_i$, for all $s_i\neq 0$.
As $i\in\N$ has been chosen arbitrarily, the claim follows.

(ii) $\Rightarrow$ (i). Pick any $s\in\ell_\infty^+\backslash\{0\}$. We can write it as $s=(s_1,s_2,\ldots)$.
Pick any $i\in\N$ so that $s_i\neq 0$. As $s\neq 0$, such $i$ exists and it holds that $\Gamma_i s_i \not\geq s_i$.

As $\Gamma_\otimes s = (\Gamma_{1,\otimes}s_1,\Gamma_{2,\otimes} s_2, \ldots)$ and by definition of the order in $\ell_\infty$ 
it holds that $\Gamma_\otimes s \geq s$ $\Iff$ $\Gamma_{i,\otimes} s_i \geq s_i$ for all $i$, it follows that $\Gamma_\otimes s \not\geq s$.

(ii) $\Iff$ (iv). Well-known, see, e.g., \cite[p.~108]{DRW07}.

(iii) $\Iff$ (iv). This follows by noting that all the nonzero cycles of $\Gamma_\otimes$ are cycles of $\Gamma_{i,\otimes}$ for a certain $i\in\N$.
\end{proof}

We will use the following lemma, which is interesting in itself since it does not require the gains to be linear. It was shown in \cite[Lemma B.5]{MKG20}, and it is a strengthening of \cite[Lemma~4.3]{DMS19a}. This result is ultimately useful for the construction of so-called paths of strict decay (a.k.a.~$\Omega$-paths) in Lyapunov-based small-gain theorems, see \cite{KMZ22,DMS19a} for infinite networks and Chapter~\ref{chap:Interconnections} for finite networks.

\begin{lemma}
\label{lem_Q}
Assume that $\Gamma_{\otimes}:\ell_{\infty}^+(I) \rightarrow \ell_{\infty}^+(I)$, defined as in \eqref{eq:Gain-operator-semimax}, is well-defined, continuous and satisfies the robust small-gain condition. Then the operator%
\begin{equation*}
Q:\ell_{\infty}^+(I) \to \ell_{\infty}^+(I),\quad Q(s) := \sup_{k \in \Z_+} \Gamma_{\otimes}^k(s),\quad s \in \ell_{\infty}^+(I),%
\end{equation*}
where the supremum is taken componentwise, is well-defined and satisfies%
\begin{equation}
\label{eq:upper estimate on Q}
  s \leq Q(s) \leq \omega^{-1}(\|s\|_{\ell_{\infty}}){\unit},\quad s \in \ell_{\infty}^+(I),
\end{equation}
where $\omega\in\Kinf$, $\omega<\id$ stems from the robust small-gain condition.
In particular,
\begin{equation}
\label{eq:Sandwich-estimates-Q}
  \|s\|_{\ell_{\infty}} \leq\|Q(s)\|_{\ell_{\infty}} \leq \omega^{-1}(\|s\|_{\ell_{\infty}}),\quad s \in \ell_{\infty}^+(I).
\end{equation}
Furthermore,
\begin{equation}
\label{eq_Gamma_Q_ineq}
  \Gamma_{\otimes}(Q(s)) =Q(\Gamma_{\otimes}(s)) \leq Q(s),\quad s \in \ell_{\infty}^+(I).
\end{equation}
\end{lemma}

\begin{proof}
Considering $k=0$ in the definition of $Q$, it is clear that $Q(s) \geq s$ for all $s\in \ell_\infty^+(I)$.

Now assume that either $Q$ is not well-defined, or $Q(s) > \omega^{-1}(\|s\|_{\ell_{\infty}}){\unit}$ for a certain $s \in \ell_\infty^+(I)$.
In any case, this implies that there are $s \in \ell_\infty^+(I)$ and $i \in I$ such that $\sup_{k \in \Z_+} \big[\Gamma_{\otimes}^k(s)\big]_i > \omega^{-1}(\|s\|_{\ell_{\infty}})$.

By formula \eqref{eq:Potenzen-join-morphism}, we obtain that there exist $k \in \N$, indices $i,j \in \N$ and a path $j_1,\ldots,j_k$ such that%
\begin{equation}\label{eq_con_ass}
  \gamma_{ij_1} \circ \gamma_{j_1j_2} \circ \cdots \circ \gamma_{j_kj}(s_j) \geq \omega^{-1}(\|s\|_{\ell_{\infty}}).%
\end{equation}
For the given $i,j$, we define the operator%
\begin{equation*}
  \tilde{\Gamma}_{ji}(s)_l := \Bigl(\sup_{k\in\N}\left[\gamma_{lk}(s_k) + \delta_{jl}\delta_{ik} \omega(s_k)\right]\Bigr)_{l\in\N},\quad s \in \ell_{\infty}^+(I),%
\end{equation*}
where $\delta_{xy}$ is the Kronecker delta, and observe that%
\begin{align*}
  \tilde{\Gamma}_{ji}(s) &\leq \Bigl(\sup_{k\in\N}\gamma_{lk}(s_k) + \sup_{k\in\N} \delta_{jl}\delta_{ik} \omega(s_k)\Bigr)_{l\in\N} \\
	&= \Bigl(\sup_{k\in\N}\gamma_{lk}(s_k) + \delta_{jl} \omega(s_i)\Bigr)_{l\in\N} = \Gamma_{\otimes}(s) + \omega(s_i) e_j = \Gamma_{ji}(s).%
\end{align*}
Since $\Gamma_{ji}$ satisfies the small-gain condition by assumption, also $\tilde{\Gamma}_{ji}$ satisfies the small-gain condition $\tilde{\Gamma}_{ji}(s) \not\geq s$ for all $s \in \ell_{\infty}^+(I) \setminus \{0\}$. By Lemma~\ref{lem_cycle_contr}, all cycles built from the gains%
\begin{equation*}
  \tilde{\gamma}_{lk}(r) := \gamma_{lk}(r) + \delta_{jl}\delta_{ik} \omega(r),\quad (l,k) \in \N^2%
\end{equation*}
are contractions. In particular,%
\begin{equation*}
  \tilde{\gamma}_{ji} \circ \tilde{\gamma}_{ij_1} \circ \tilde{\gamma}_{j_1j_2} \circ \cdots \circ \tilde{\gamma}_{j_kj}(s_j) < s_j.%
\end{equation*}
With \eqref{eq_con_ass}, we thus obtain%
\begin{align*}
  s_j &> \tilde{\gamma}_{ji} \circ \tilde{\gamma}_{ij_1} \circ \tilde{\gamma}_{j_1j_2} \circ \cdots \circ \tilde{\gamma}_{j_kj}(s_j) \geq \tilde{\gamma}_{ji} \circ \gamma_{ij_1} \circ \gamma_{j_1j_2} \circ \cdots \circ \gamma_{j_kj}(s_j) \\
			&= ( \gamma_{ji} + \omega ) \circ \gamma_{ij_1} \circ \gamma_{j_1j_2} \circ \cdots \circ \gamma_{j_kj}(s_j) \geq \|s\|_{\ell_{\infty}} \geq s_j,%
\end{align*}
a contradiction. This shows \eqref{eq:upper estimate on Q}.

The sandwich estimate \eqref{eq:Sandwich-estimates-Q} follows from \eqref{eq:upper estimate on Q} by monotonicity of the norm in $\ell_\infty(I)$.%

Let us show \eqref{eq_Gamma_Q_ineq}. For every $i\in\N$, we have%
\begin{align*}
  [\Gamma_{\otimes}(Q(s))]_i &= \sup_{j\in\N} \gamma_{ij}( Q(s)_j ) = \sup_{j\in\N} \gamma_{ij}( \sup_{k\in\Z_+} \Gamma_{\otimes}^k(s)_j ) \\
	               &= \sup_{j\in\N}\sup_{k\in\Z_+} \gamma_{ij}(\Gamma_{\otimes}^k(s)_j) = \sup_{k\in\Z_+}\sup_{j\in\N} \gamma_{ij}(\Gamma_{\otimes}^k(s)_j) \\
								 &= \sup_{k\in\Z_+}\Gamma_{\otimes}(\Gamma_{\otimes}^k(s))_i = \sup_{k\in\Z_+}\Gamma_{\otimes}^{k+1}(s)_i \leq \sup_{k\in\Z_+} \Gamma_{\otimes}^k(s)_i = Q(s)_i,%
\end{align*}
and hence \eqref{eq_Gamma_Q_ineq} holds. 
Finally, for any $s \in \ell_\infty^+(I)$ we have
\begin{eqnarray*}
Q(\Gamma_{\otimes}(s))  
= \sup_{k \in \Z_+} \Gamma_{\otimes}^{k+1}(s)
= \sup_{k \in \Z_+} \Gamma_{\otimes} (\Gamma_{\otimes}^{k}(s))
= \Gamma_{\otimes} (\sup_{k \in \Z_+} (\Gamma_{\otimes}^{k}(s)))
= \Gamma_{\otimes} (Q(s)).
\end{eqnarray*}
\end{proof}

\begin{remark}
\label{rem:Kleene-star} 
\index{operator!Kleene star}
\index{strong transitive closure}
The operator $Q$ is of fundamental importance in max-algebra, and is sometimes called the strong transitive closure of $\Gamma_\otimes$, or Kleene star, see \cite[Section~1.6.2]{But10}, or just the closure of $\Gamma_\otimes$, see \cite[Section~1.4]{Rue17}.
\end{remark}

		\section{Concluding remarks}
		
		This chapter is based on \cite{MKG20}. For historical remarks and a brief overview of the small-gain theory for finite and infinite networks, we refer to Section~\ref{sec:Lyapunov-based-SGT-Concluding remarks}.


\cleardoublepage
\chapter[Lyapunov-based small-gain theorems]{Lyapunov-based small-gain theorems for finite interconnections}
\label{chap:Interconnections}

As in most cases, the ISS of nonlinear systems is verified by the construction of an appropriate ISS Lyapunov function,
a natural desire is to use in the formulation of the small-gain theorems the information about the ISS Lyapunov functions for subsystems.
In this chapter, we show that this is possible and allows for the construction of ISS Lyapunov functions for the overall system if the ISS Lyapunov functions for the subsystems are known.

\section[Interconnections of ISS systems]{Interconnections of input-to-state stable systems}
\label{GekoppelteISS_Systeme}

Consider an interconnected system of the form
\begin{equation}
\label{Kopplung_N_Systeme}
\left\{
\begin{array}{l}
\dot{x}_{i}=A_{i}x_{i}+f_{i}(x_{1},\ldots,x_{n},u),\quad x_{i}(t)\in X_{i}, \ u(t) \in U,\\
i=1,\ldots,n,
\end{array}
\right.
\end{equation}
where the state space of the $i$-th subsystem $X_{i}$ is a Banach space and $A_i$ is the generator of a $C_{0}$-semigroup on $X_{i}$, $i=1,\ldots, n$.
The input space we take as $\Uc:=PC_b(\R_+,U)$ for some Banach space of input values $U$. 
The state space of the system \eqref{Kopplung_N_Systeme} we define as
\[
X=X_{1}\times\ldots\times X_{n},
\]
which is a Banach space with the norm 
\[
\|\cdot\|_{X}:=\|\cdot\|_{X_{1}}+\ldots+\|\cdot\|_{X_{n}}.
\]
The internal input space for the $i$-th subsystem is
\[
X_{\neq i}:=X_1 \times \ldots \times X_{i-1} \times X_{i+1} \times \ldots \times X_n.
\]
The norm in $X_{\neq i}$ we introduce as
\[
\|\cdot\|_{X_{\neq i}}:=\|\cdot\|_{X_{1}}+\ldots + \|\cdot\|_{X_{i-1}} + \|\cdot\|_{X_{i+1}} + \ldots +\|\cdot\|_{X_{n}}.
\]
The elements of $X_{\neq i}$ we denote by $x_{\neq i}=(x_j)_{j\neq i}$.

The transition map of the $i$-th subsystem we denote by 
\[
\phi_i:\R_+ \times X_i \times (PC_b(\R_+, X_{\neq i}) \tm \Uc) \to X_i. 
\]
For $x=(x_i)_{i=1}^n$, $x_i \in X_i,\ i=1,\ldots, n$, define
\[
f(x,u):=(f_{1}(x,u)^T,{\ldots},f_{n}(x,u)^T)^T,
\]
and
\[
A:=\diag(A_{1},\ldots,A_{n}),
\quad
D(A):=D(A_{1})\times\ldots\times D(A_{n}).
\]
Clearly, $A$ is the generator of a $C_{0}$-semigroup on $X$.

We rewrite the system \eqref{Kopplung_N_Systeme} in the vector form:
\begin{equation}
\label{KopplungHauptSys}
\dot{x}=Ax+f(x,u).
\end{equation}

We assume that 
\begin{ass}
\label{ass:Assumption-network} 
The system \eqref{KopplungHauptSys} satisfies the following properties:
\begin{enumerate}[label=(\roman*)]
	\item $f$ is continuous on $X \tm U$.
	\item For each $i=1,\ldots,n$, $f_i$ is Lipschitz continuous w.r.t.\  $x_i$ on bounded subsets.
	\item Each subsystem is forward complete.
	\item (Well-posedness) For any initial condition $x \in X$ and any input $u\in\Uc$, there is a unique maximal solution of \eqref{KopplungHauptSys}.
	\item \eqref{KopplungHauptSys} satisfies the BIC property.
\end{enumerate}
\end{ass}

\section{Lyapunov-based small-gain theorems}
\label{sec:ISS-SGT-Lyapunov-form}

For the analysis of networks, we have to know the response of the subsystems (or of the corresponding Lyapunov functions) to the inputs from all other subsystems. To aggregate the individual inputs into the total input, we use the following formalism:
\begin{definition}
\label{def:Monotone aggregation functions}
\index{monotone aggregation function}
\index{function!monotone aggregation}
\index{MAF}
A function $\mu:\R^n_+ \to\R_+$ is called a \emph{monotone aggregation function (MAF)} if $\mu$ is continuous and satisfies the following properties:
\begin{enumerate}[label=(\roman*)]
	\item \emph{Positivity}: $\mu(v)\geq 0$, for all $v \in \R^n_+$.
	\item \emph{Strict increase}: if $v_i < w_i$ for all $i=1,\ldots,n$, then $\mu(v)<\mu(w)$.
	\item \emph{Unboundedness}: $\mu(v) \to\infty$ as $|v| \to\infty$.
	\item \emph{Sub-additivity}: $\mu(v+w) \leq \mu(v) + \mu(w)$, for all $x,y\in\R^n_+$.
\end{enumerate}
\end{definition}
Two important examples of monotone aggregation functions are 
\begin{eqnarray}
\mu_\boxplus: s = (s_i)_{i=1}^n \mapsto s_1+\ldots+s_n,
\label{eq:sum-MAF}
\end{eqnarray}
and 
\begin{eqnarray}
\mu_\otimes: s = (s_i)_{i=1}^n \mapsto \max\{s_1,\ldots,s_n\}.
\label{eq:max-MAF}
\end{eqnarray}

Now, in addition to Assumption~\ref{ass:Assumption-network}, let the following hold:
\begin{ass}
\label{ass_subsystem_iss}
For each $i=1,\ldots,n$, there exists a continuous function $V_i:X_i \rightarrow \R_+$, that satisfies the following properties:
\begin{enumerate}[label=(\roman*)]
\item There exist $\psi_{i1},\psi_{i2} \in \Kinf$ such that%
\begin{equation}
\label{eq_subsystem_iss_coerc}
  \psi_{i1}(\|x_i\|_{X_i}) \leq V_i(x_i) \leq \psi_{i2}(\|x_i\|_{X_i}),\quad x_i \in X_i.%
\end{equation}
\item There exists a MAF $\mu_i$, together with $\chi_{ij} \in \Kinf\cup\{0\}$, $j=1,\ldots,n$, where $\chi_{ii} = 0$, and $\chi_{iu} \in \K$, $\alpha_i \in \PD$ such that for all $x_{j}\in X_j$, $j=1,\ldots,n$, and all $u\in \Uc$ the following implication holds:%
\begin{align}
\label{eq_subsystem_orbitalder_est}
\begin{split}
  V_i(x_i) > \max\Bigl\{\mu_i \Big[\big(\chi_{ij}&(V_j(x_j)\big)\big)_{j=1}^n\Big],\chi_{iu}(\|u\|_{\Uc}) \Bigr\} \\
	&\Rightarrow \dot{V}_{i,u}(x_i) \leq - \alpha_i(V_i(x_i)).
\end{split}
\end{align}
\end{enumerate}
The functions $\chi_{ij}$ are called \emph{internal Lyapunov gains}, while the functions $\chi_{iu}$ are called \emph{external Lyapunov gains}.%
\end{ass}

The function $\mu_i$ combines (\q{aggregates}) the inputs from all $x_j$-subsystems into a total internal input to the $x_i$-subsystem from other elements of the network. 
If $\mu_i=\mu_\otimes$, then the total input to the $i$-th subsystem is the maximum of the inputs from the individual subsystems.
If $\mu_i=\mu_\boxplus$, see \eqref{eq:sum-MAF}, then the total input to the $i$-th subsystem is the sum of the inputs from the individual subsystems.
Allowing for such flexibility, we can obtain sharper conditions for the stability of the network in a wide variety of applications.

We collect the Lyapunov gains $\chi_{ij} \in \Kinf\cup \{0\}$ for $i,j\in\{1,2,\ldots,n\}$ 
into the matrix $\Gamma:=(\chi_{ij})_{i,j=1}^n$. 
We combine all MAFs $\mu_i$ into the vector of monotone aggregation functions $\mu:=(\mu_i)_{i=1}^n$.

\index{operator!gain}
Then the pair $(\Gamma,\mu)$ gives rise to the following monotone and continuous operator, called \emph{gain operator}
\begin{eqnarray}
\Gamma_\mu(s):=
\begin{pmatrix}
\mu_1\big(\chi_{11}(s_1),\ldots,\chi_{1n}(s_n)\big)\\
\vdots\\
\mu_n\big(\chi_{n1}(s_1),\ldots,\chi_{nn}(s_n)\big)
\end{pmatrix}
,\quad s = (s_i)_{i=1}^n \in \R^n_+.
\label{eq:Induced-operator-MAF}
\end{eqnarray}

If $\mu_i=\mu_\boxplus$ for all $i$, then $\Gamma_\mu$ corresponds to the operator 
\[
\Gamma_{\boxplus}:=\Big(\sum_{j=1}^n\chi_{ij}(s_j)\Big)_{i=1}^n,
\]
and if $\mu_i=\mu_\otimes$ for all $i$, then $\Gamma_\mu$ corresponds to the operator
\[
\Gamma_{\otimes}:=\Big(\max_{j=1}^n\chi_{ij}(s_j)\Big)_{i=1}^n.
\]
We will need the following technical lemmas that we state without the proof (see \cite[Section 3.5]{Mir23}).
\begin{lemma}
\label{LIM_MAX}
Let $(x^1_k)_{k=1}^{\infty}$, $\ldots$, $(x^m_k)_{k=1}^{\infty}$ be sequences of real numbers. 
Let the limit $\lim_{k \to \infty} \max_{i=1}^m x^i_k$ exist. 
Then 
\begin{equation}
\label{BeideRichtungen}
\lim_{k \to \infty} \max_{i=1}^m x^i_k = \max_{i=1}^m \Limsup_{k \to \infty}x^i_k,
\end{equation}
where $\Limsup\limits_{k \to \infty}x^i_k$ is the upper limit of the sequence $x^i_k$.
\end{lemma}

\begin{corollary}
Let $f_i:\R \to \R$ be defined and bounded in some neighborhood $D$ of $t=0$. The following holds
\begin{equation}
\label{BeideRichtungen_Cor}
\Limsup_{t \to 0} \max_{i=1}^m f_i(t)  = 
 \max_{i=1}^m \Limsup_{t \to 0} f_i(t) 
\end{equation}
\end{corollary}


To construct an ISS Lyapunov function for the whole interconnection, we will use the following notion:
\begin{definition} 
\index{$\Omega$-path}
\label{Omega-Path-Def}
Let $\Gamma:\R_{+}^{n} \rightarrow \R_{+}^{n}$ be a nonlinear operator.
A function $\sigma=(\sigma_{1},\dots,\sigma_{n})^{T}:\R_{+}^{n}\rightarrow\R_{+}^{n}$,
where $\sigma_{i}\in\K_{\infty}$, $i=1,\ldots,n$, is called a \emph{regular path of non-strict decay} with respect to the operator $\Gamma$, 
if it possesses the following properties:
\begin{enumerate}[label=(\roman*)]
\item For every compact interval $K \subset (0,\infty)$, there exist $0 < c \leq C < \infty$ such that for all $r_1,r_2 \in K$ and $i \in \N$%
\begin{equation}
\label{eq:Lipschitz-bounds-path-of-decay}
  c|r_1 - r_2| \leq |\sigma_i^{-1}(r_1) - \sigma_i^{-1}(r_2)| \leq C|r_1 - r_2|.%
\end{equation}
\item $\Gamma$ is nonincreasing on $\sigma(\R_+)$:
\begin{align}
\label{sigma cond 2-nonstrict}
\Gamma(\sigma(r))\leq\sigma(r)\quad \forall r\geq 0.
\end{align}
\end{enumerate}
If the map $\sigma$ satisfies also $\Gamma(\sigma(r)) \ll \sigma(r)$, $\forall r>0$, then $\sigma$ is called a 
\emph{regular path of strict decay}.
\end{definition}

The next theorem provides the construction of an ISS Lyapunov function for an interconnection of ISS subsystems, provided a path of non-strict decay for the gain operator is available.
\begin{theorem}[Lyapunov-based ISS small-gain theorem]
\index{small-gain theorem!in Lyapunov formulation!for ODE systems}
\label{Constr_of_LF_FinDim} 
Let Assumption~\ref{ass_subsystem_iss} hold. 
If there exists a regular path of non-strict decay $\sigma=(\sigma_{1},\ldots,\sigma_{n})^{T}$ corresponding to the operator $\Gamma_\mu$ defined by
\eqref{eq:Induced-operator-MAF}, then an ISS Lyapunov function for \eqref{KopplungHauptSys} can be constructed as
\begin{equation}
\label{NeuLyapFun}
V(x):=\max_{i=1}^n \sigma_{i}^{-1}(V_{i}(x_{i})),
\end{equation}
and the Lyapunov gain of the whole system is given by
\begin{equation}
\label{ChiDef}
\chi(r):=\max_{i=1}^n\sigma_{i}^{-1}(\chi_{i}(r)).
\end{equation}
\end{theorem}

\begin{proof} 
We divide the proof into several parts.

\textbf{(i). Partition of $X$.} 
Pick any $u \in \Uc$. To prove that $V$ is an ISS Lyapunov function, it is useful to divide its domain of definition into subsets on which $V$ takes a simpler form. Thus, for all $i\in \{1,\ldots,n\}$ define a set
\[
M_{i}=\left\{ x\in X:\sigma_{i}^{-1}(V_{i}(x_{i})) > \sigma_{j}^{-1}(V_{j}(x_{j}))\,\,\forall j=1,\ldots,n,\ j\neq i\right\}.
\]
From the continuity of $V_i$ and $\sigma_{i}^{-1}$, $i=1,\ldots,n$ it follows that all $M_i$ are open. Also, note that 
$X=\bigcup_{i=1}^{n}\clo{M_i}$ and for all $i,j$: $i \neq j$ holds $M_i \bigcap M_j= \emptyset$.

\textbf{(ii). Estimates for $\dot{V}_u(x)$ on $\bigcup_{i=1}^{n}{M}_{i}$.} 
Take some $i\in \{1,\ldots,n\}$ and pick any $x\in M_{i}$. Assume that $V(x)\geq\chi(\|u\|_{\Uc})$ holds, where $\chi$ is given by \eqref{ChiDef}. 
We obtain
\[
\sigma_{i}^{-1}(V_{i}(x_{i}))=V(x)\geq\chi(\|u\|_{\Uc})=\max_{j=1}^{n}\sigma_{j}^{-1}\circ\chi_{j}(\|u\|_{\Uc})\geq\sigma_{i}^{-1}(\chi_{i}(\|u\|_{\Uc})).
\]
Since $\sigma_{i}^{-1} \in \Kinf$, it holds that
\begin{eqnarray}
\label{eq:GainAbschaetzung}
V_{i}(x_{i}) \geq \chi_{i}(\|u\|_{\Uc}).
\end{eqnarray}
On the other hand, from the condition of non-strict decay \eqref{sigma cond 2-nonstrict}, and using the monotonicity of the MAF $\mu_i$, we obtain that
\begin{eqnarray*}
V_{i}(x_{i})=\sigma_{i}(V(x))    &\geq& \Big[\Gamma_\mu\Big(\sigma\big(V(x)\big)\Big)\Big]_i \\
& =&  \mu_i\Big[\Big(\chi_{ij}\big(\sigma_{j}\left(V\left(x\right)\right)\big)\Big)_{j=1}^{n}\Big]= \mu_i\Big[\Big(\chi_{ij}\big(\sigma_{j}\left(\sigma_{i}^{-1}\left(V_i\left(x_{i}\right)\right)\right)\big)\Big)_{j=1}^{n}\Big]\\
    &\ge&     \mu_i\Big[\Big(\chi_{ij}\big(\sigma_{j}\left(\sigma_{j}^{-1}\left(V_j\left(x_{j}\right)\right)\right)\big)\Big)_{j=1}^{n}\Big]
		=\mu_i\Big[\Big(\chi_{ij}\big(V_j\left(x_{j}\right)\big)\Big)_{j=1}^{n}\Big].
\end{eqnarray*}
Combining with \eqref{eq:GainAbschaetzung}, we obtain
\begin{eqnarray}
\label{eq:LyapImplikation}
V_{i}(x_{i})  \geq  \max\Big\{ \mu_i\Big[\Big(\chi_{ij}\big(V_j\left(x_{j}\right)\big)\Big)_{j=1}^{n}\Big],\chi_{i}(\|u\|_{\Uc})\Big\}.
\end{eqnarray}
Now, by \eqref{eq_subsystem_orbitalder_est} we have that 
\begin{eqnarray}
\dot{V}_{i,u}(x_i) \leq - \alpha_i(V_i(x_i)).
\label{eq:Decay-V_i_aux}
\end{eqnarray}

As $\sigma_{i}^{-1}(V_i(x_i))> \sigma_{j}^{-1}(V_j(x_j))$ for all $j\neq i$, and since the trajectories of all subsystems are continuous, 
there is a time $t>0$ such that for all $s \in [0,t]$ it holds that
$\sigma_{i}^{-1}\big(V_i(\phi_i(s,x,u))\big)> \sigma_{j}^{-1}\big(V_j(\phi_j(s,x,u))\big)$ for all $j\neq i$.
Thus, $V(\phi(s,x,u)) = \sigma_{i}^{-1}\big(V_i(\phi_i(s,x,u))\big)$ for $s \in[0,t]$.
Note also that $V_i(\phi_i(s,x,u)) \leq V_i(x_i)$ in view of \eqref{eq:Decay-V_i_aux}.

We have
\begin{eqnarray}
\label{eq:Difference-quotient-tmp}
\hspace{-8mm}\frac{1}{s}\Big(V(\phi(s,x,u)) - V(x)\Big) 
&=& \frac{1}{s}\Big(\sigma_{i}^{-1}\big(V_i(\phi_i(s,x,u))\big) - \sigma_{i}^{-1}\big(V_i(x_i)\big)\Big) \nonumber \\
&=& -\frac{1}{s}\Big|\sigma_{i}^{-1}\big(V_i(\phi_i(s,x,u))\big) - \sigma_{i}^{-1}\big(V_i(x_i)\big)\Big|.
\end{eqnarray}
Define (pointwise) $\sigma_{min}:=\min_{i=1}^n\sigma_i$ and $\sigma_{max}:=\max_{i=1}^n\sigma_i$.

Thanks to continuity of $s \mapsto \sigma_{i}^{-1}\big(V_i(\phi_i(s,x,u))\big)$, we choose $t$ small enough, such that for all $s \in[0,t]$ it holds that 
\[
\frac{1}{2}V_i(x_i) \leq V_i(\phi_i(s,x,u)) \leq V_i(x_i) = \sigma_i(V(x)) \leq \sigma_{max}(V(x)).
\]
Similarly, we estimate
\[
\frac{1}{2}V_i(x_i) = \frac{1}{2}\sigma_i(V(x)) \geq \frac{1}{2}\sigma_{min}(V(x)).
\]
Define for any $r\geq 0$ 
\[
K(r):=\Big[\frac{1}{2}\sigma_{min}(r), \sigma_{max}(r)\Big],
\]
and let $c = c(K(r)) > 0$ be the maximal constant such that%
\begin{equation*}
  |\sigma_i^{-1}(r_1) - \sigma_i^{-1}(r_2)| \geq c|r_1 - r_2|, \quad \forall r_1,r_2 \in K(r).%
\end{equation*}
For all $s\in(0,t)$, we obtain from \eqref{eq:Difference-quotient-tmp} that
\begin{align*}
\frac{1}{s}\Big(V(\phi(s,x,u)) - V(x)\Big)  
	&\leq -c\big(K(V(x))\big)\frac{1}{s}\Big|V_i\big(\phi_i(s,x,u)\big) - V_i(x_i)\Big|\\
	&= c(K(V(x)))\frac{1}{s}\big(V_i(\phi_i(s,x,u)) - V_i(x_i)\big).
\end{align*}
Taking the limit superior, we obtain
\begin{align*}
\dot{V}_u(x) &= \Limsup_{s\to+0} \frac{1}{s}\Big(V(\phi(s,x,u)) - V(x)\Big)  \\
	&\leq c(K(V(x)))\dot{V}_i(x_i)\\
	&\leq -c(K(V(x))) \alpha_i(V_i(x_i))\\
	&= -c(K(V(x))) \alpha_i(\sigma_i(V(x))).
\end{align*}
Define now 
\[
\hat{\alpha}(r):= c(K(r)) \min_{i=1}^n \alpha_i(\sigma_i(r)),\quad r\geq 0.
\]

\textbf{(iii). Lower bound for $\hat{\alpha}$.}  
It remains to lower bound $\hat{\alpha}$ by a positive definite function. For each $r>0$, define $K_2(r) := \bigcup_{q \in K(r)}\mathrm{str}(1,q)$, where $\mathrm{str}(1,q)$ equals $[1,q]$ for $q\geq 1$ and $[q,1]$ for $q<1$. Clearly, $K_2(r)$ is a compact subset of $(0,\infty)$. Further, we introduce%
\begin{equation*}
  \hat{\alpha}_2(r) := c(K_2(r)) \min_{i=1}^n \alpha_i(\sigma_i(r)),\quad \forall r > 0.%
\end{equation*}
As $K(r)\subset K_2(r)$ for any $r>0$, $\hat{\alpha}(r) \geq \hat{\alpha}_2(r)$ for all $r>0$. Furthermore, there is $r_{\min}$ such that $K_2(r_1) \supset K_2(r_2)$ for all $r_1,r_2 \in (0,r_{\min})$ with $r_1 < r_2$. This implies that $\hat{\alpha}_2$ is a nondecreasing positive function on $(0,r_{\min})$, and $\lim_{r\to -0}\hat{\alpha}_2(r) = 0$. Moreover, for all $r\in (0,r_{\min})$ we have%
\begin{equation*}
  \hat{\alpha}_2(r) = \frac{2}{r}\int_{r/2}^r \hat{\alpha}_2(r)\, ds \geq \frac{2}{r}\int_{r/2}^r \hat{\alpha}_2(s)\, ds,%
\end{equation*}
where $\hat{\alpha}_2$ is integrable on $(0,r_{\min})$ as it is monotone on this interval. Hence, $\hat{\alpha}_2$, and thus $\hat{\alpha}$, can be lower bounded by a continuous function on $[0,r_{\min}]$. Similarly, there is $r_{\max}$ such that $K_2(r_1)\supset K_2(r_2)$ for all $r_1,r_2 \in (r_{\max},\infty)$ with $r_1 > r_2$. This implies that $\hat{\alpha}_2$ is a 
nonincreasing positive function on $(r_{\max},\infty)$. Consequently, for $r\in (r_{\max},\infty)$, we have%
\begin{equation*}
  \hat{\alpha}_2(r) = \frac{1}{r}\int_{r}^{2r} \hat{\alpha}_2(r)\, ds \geq \frac{1}{r}\int_{r}^{2r} \hat{\alpha}_2(s)\, ds.%
\end{equation*}
Hence, $\hat{\alpha}_2$, and thus $\hat{\alpha}$, can be lower bounded by a continuous function on $[r_{\max},\infty)$. As $\hat{\alpha}$ assumes positive values and is bounded away from zero on every compact interval in $(0,\infty)$, $\hat{\alpha}$ can be lower bounded by a positive definite function, which we denote by $\alpha$.

Overall, for all $x\in \bigcup_{i=1}^{n}{M}_{i}$ it holds that
\begin{align*}
\dot{V}_u(x) \leq -\alpha(V(x)).
\end{align*}

\textbf{(iv). Estimates for $\dot{V}_u(x)$ on $X$.} 
Now let $x \notin \bigcup_{i=1}^{n}{M}_{i}$. From $X =\bigcup_{i=1}^{n}{\clo{M_i}}$ it follows that $x \in \bigcap_{i\in I(x)}\partial {M}_{i}$  for some index set $I(x)\subset\left\{ 1,\ldots,n\right\}$, $|I(x)| \geq 2$. Clearly, 
\begin{align*}
\bigcap_{i\in I(x)}\partial {M}_{i} = \big\{ x \in X: \forall i &\in I(x),\ \forall j \notin I(x)\  \sigma_{i}^{-1}(V_{i}(x_{i})) > \sigma_{j}^{-1}(V_{j}(x_{j}))\\
 &\forall i,j \in I(x) \ \sigma_{i}^{-1}(V_{i}(x_{i})) = \sigma_{j}^{-1}(V_{j}(x_{j})) \big\}. 
\end{align*}
Due to continuity of $\phi$, we have, there exists $t>0$, such that for all $s \in [0, t)$ it follows that 
\[
\phi(s,x,u) \in \left( \bigcap_{i\in I(x)} \partial {M}_{i} \right) \cup \left( \bigcup_{i\in I(x)}{M}_{i} \right).
\]

By definition of the Lie derivative, we obtain
\begin{eqnarray}
\label{HauptAbschaetzung}
\hspace{-8mm}\dot{V}_u(x) &=& \Limsup\limits _{t\rightarrow+0}\frac{1}{t}\Big(V(\phi(t,x,u))-V(x)\Big) \notag\\
&=&
\Limsup\limits _{t\rightarrow+0}\frac{1}{t} \left(\max_{i \in I(x)}\big\{\sigma_{i}^{-1}(V_{i}(\phi_i(t,x,u)))\}-\max_{i \in I(x)}\{\sigma_{i}^{-1}(V_{i}(x_{i}))\big\} \right).
\end{eqnarray}
From the definition of $I(x)$ it follows that
\[
\sigma_i^{-1}(V_{i}(x_{i}))=\sigma_j^{-1}(V_{j}(x_{j})) \quad \forall i,j\in I(x),
\]
and therefore the index $i$, on which the maximum of $\max_{i \in I(x)}\{\sigma_{i}^{-1}(V_{i}(x_{i}))\}$ is attained, may always be set equal to the index on which the maximum \linebreak
$\max_{i \in I(x)}\{\sigma_{i}^{-1}(V_{i}(\phi_i(t,x,u)))\}$ is reached. \\
We infer from the estimates \eqref{HauptAbschaetzung} that
\[
\dot{V}_u(x)=\Limsup\limits _{t\rightarrow+0}\max_{i \in I(x)}\Big\{\frac{1}{t}\left(\sigma_{i}^{-1}(V_{i}(\phi_i(t,x,u)))-\sigma_{i}^{-1}(V_{i}(x_{i}))\right) \Big\}.
\]
Using Lemma \ref{LIM_MAX} we obtain
\begin{eqnarray*}
\dot{V}_u(x) &=& \max_{i \in I(x)}\Big\{   \Limsup \limits _{t\rightarrow+0}\frac{1}{t}\left(\sigma_{i}^{-1}(V_{i}(\phi_i(t,x,u)))-\sigma_{i}^{-1}(V_{i}(x_{i}))\right)\Big\} \\
&\leq&   -\alpha(V(x)).
\end{eqnarray*}
Overall, we have that for all $x \in X$ and $u\in\Uc$ it holds that
\[
V(x)\geq\chi(\|u\|_{\Uc})   \qrq  \dot{V}_u(x)\leq-\alpha(V(x)),
\]
and the ISS Lyapunov function in implication form for the whole interconnection is constructed. ISS of the whole system follows by 
Theorem~\ref{LyapunovTheorem}.
\end{proof}

We have shown that the network of ISS systems is ISS, provided that there is a path of non-strict decay for the operator $\Gamma_\mu$. The existence of a path of strict (and thus also of a non-strict) decay (at least up to its regularity) is guaranteed provided that $\Gamma_\mu$ satisfies a so-called strong small-gain condition, see \cite[Appendix C, Section C.8]{Mir23}.
As a special case, consider the interconnection of 2 systems with the corresponding internal gains $\chi_{12}$, $\chi_{21}$, with $\mu_1=\mu_2=\id$.

\begin{corollary}[Coupling of 2 systems]
\index{small-gain theorem!in Lyapunov formulation!for ODE systems}
\label{cor:SGT-2-systems} 
Let $n=2$, and let Assumption~\ref{ass_subsystem_iss} hold with $\mu_1=\mu_2=\id$.
If 
\begin{eqnarray}
\chi_1 \circ \chi_2(r) < r,\quad r>0,
\label{eq:SGT-2-sys}
\end{eqnarray}
then  \eqref{KopplungHauptSys} is ISS.
\end{corollary}

\begin{proof}
In view of assumptions, the gain operator takes the form
\begin{eqnarray}
\Gamma\left(\begin{array}{c} s_1 \\ s_2\end{array}\right)
=
\left(\begin{array}{c}
\chi_{12}(s_2)\\
\chi_{21}(s_1)
\end{array}\right).
\end{eqnarray}
According to \cite[p. 108]{DRW07}, \cite{Mir23}, the condition \eqref{eq:SGT-2-sys} ensures the existence of a regular path of non-strict decay for our gain operator.
Application of Theorem~\ref{Constr_of_LF_FinDim} finishes the proof.
\end{proof}

\subsection{Example: coupling of linear systems}

As a particular example, consider the following system of interconnected linear reaction-diffusion equations
\begin{equation}
\label{GekoppelteLinSyst}
\left
\{
\begin {array} {l}
{x_{1,t}(z,t) = c_1 x_{1,zz}(z,t) + a_{12} x_2(z,t) , \quad z \in (0,d),\ t>0,}\\
{x_1(0,t) = x_1(d,t)=0;} \\
{x_{2,t}(z,t) = c_2 x_{2,zz}(z,t) + a_{21} x_1(z,t) , \quad z \in (0,d),\ t>0,}\\
{x_2(0,t) = x_2(d,t)=0.}
\end {array}
\right.
\end{equation}
Here $d$,  $c_1$, and $c_2$ are positive constants.

We choose the state space for subsystems as $X_1:=L^2(0,d)$, $X_2:=L^2(0,d)$.
The operators $A_i= c_i \frac{d^2}{dz^2}$ with $D(A_i)=H^1_0(0,d) \cap H^2(0,d)$, $i=1,2$, are the generators of analytic semigroups for the corresponding subsystems.

Both subsystems are ISS, moreover, $\sigma(A_i) = \{ - c_i \left( \frac{\pi n}{d} \right)^2 : n \in\N \}$, $i=1,2$.

Take $P_i = \frac{1}{2c_i} \left( \frac{d}{\pi} \right)^2 \id$, where $\id$ is the identity operator on $X_i$.
We have
\begin{eqnarray*}
\lel A_i x_i ,P_i x_i \rir + \lel P_i x_i,A_i x_i \rir &=& \frac{1}{c_i} \left(\frac{d}{\pi} \right)^2 \lel A_i x_i, x_i \rir \\
   &= &  \left( \frac{d}{\pi} \right)^2    \int_0^d x_{i,zz}(z) x_i(z) dz
		=  - \left( \frac{d}{\pi} \right)^2    \int_0^d  x_{i,z}^2(z) dz \\   &\leq & - \|x_i\|^2_{L^2(0,d)}.
\end{eqnarray*}
In the last estimate, we have used the Friedrichs' inequality \eqref{ineq:Friedrichs}.
The Lyapunov functions for subsystems are defined by
\[
V_i(x_i) = \lel P_i x_i,x_i \rir = \frac{1}{2c_i} \left( \frac{d}{\pi} \right)^2 \|x_i\|^2_{L^2(0,d)}, \quad x_i \in X_i.
\]
Using Cauchy-Schwarz inequality, we arrive at the following estimates for the derivatives:
\[
\dot{V}_1(x_1)  \leq -\|x_1\|^2_{L^2(0,d)} + \frac{1}{c_1}\left( \frac{d}{\pi} \right)^2 |a_{12}|  \|x_1\|_{L^2(0,d)}\|x_2\|_{L^2(0,d)},
\]
\[
\dot{V}_2(x_2) \leq -\|x_2\|^2_{L^2(0,d)} + \frac{1}{c_2}\left( \frac{d}{\pi} \right)^2 |a_{21}|  \|x_1\|_{L^2(0,d)}\|x_2\|_{L^2(0,d)}.
\]
We choose the gains in the following way:
\[
\gamma_{12}(r)=\frac{c_2}{c_1^3}\left(\frac{d}{\pi} \right)^4 \left|\frac{a_{12}}{1-\varepsilon}\right|^2 r,
\qquad
\gamma_{21}(r)=\frac{c_1}{c_2^3} \left( \frac{d}{\pi} \right)^4 \left|\frac{a_{21}}{1-\varepsilon}\right|^2 r,
\qquad 
r\geq 0.
\]
For all $x_1,x_2 \in L^2(0,d)$ it holds that 
\begin{eqnarray*}
V_1(x_1) \geq \gamma_{12} \circ V_2(x_2) &\qiq& V_1(x_1) \geq \gamma_{12}(1) V_2(x_2)  \\
 &\qiq& \sqrt{ \frac{c_1}{c_2} \gamma_{12}(1) }\|x_2\|_{L^2(0,d)} \leq \|x_1\|_{L^2(0,d)}  \\
& \Iff & \frac{1}{c_1}\left( \frac{d}{\pi} \right)^2 |a_{12}| \|x_2\|_{L^2(0,d)} \leq
(1-\eps)\|x_1\|_{L^2(0,d)}.
\end{eqnarray*}
In view of the above inequalities, for all $x_1,x_2\in L^2(0,d)$ we have the following implication:
\[
V_1(x_1) \geq \gamma_{12} \circ V_2(x_2) \qrq \dot{V}_1(x_1) \leq - \varepsilon \|x_1\|^2_{L^2(0,d)}.
\]
This shows that $V_1$ is an ISS Lyapunov function for the first subsystem of \eqref{GekoppelteLinSyst}.

Analogously, assuming that for some $x_1,x_2 \in L^2(0,d)$ the following inequality holds
\begin{eqnarray*}
V_2(x_2) \geq \gamma_{21} \circ V_1(x_1) \qiq \frac{1}{c_2}\left( \frac{d}{\pi} \right)^2 |a_{21}| \|x_1\|_{L^2(0,d)} \leq (1-\eps)\|x_2\|_{L^2(0,d)},
\end{eqnarray*}
we have that
\[
V_2(x_2) \geq \gamma_{21} \circ V_1(x_1) \qrq \dot{V}_2(x_2)  \leq  - \varepsilon \|x_2\|^2_{L^2(0,d)}.
\]
This shows that $V_2$ is an ISS Lyapunov function for the second subsystem of \eqref{GekoppelteLinSyst}.

Let us find conditions guaranteeing the validity of the small-gain condition:
\[
\gamma_{12} \circ \gamma_{21} < \id \quad  \Iff \quad
\frac{1}{c_1^2 c_2^2}\left( \frac{d}{\pi} \right)^8 \frac{\left|a_{12} a_{21} \right|^2}{(1-\varepsilon)^4} <1,
\]
where we can choose $\eps>0$ arbitrarily small. 

By Corollary~\ref{cor:SGT-2-systems}, if
\begin{equation}
\label{StabBedin_LinDiff}
|a_{12} a_{21}| <c_1c_2\left( \frac{\pi}{d} \right)^4
\end{equation}
is satisfied, then the whole system \eqref{GekoppelteLinSyst} is 0-UGAS.

\subsection{Example: coupling of nonlinear systems}

Consider now a coupling of two nonlinear parabolic systems:
\begin{equation}
\label{eq:Gekoppelte-parabolic-sys-ISS}
\left
\{
\begin {array} {l}
{x_{1,t}(z,t) =  q_1 x_{1,zz}(z,t) + x^2_2(z,t), \quad z \in (0,\pi),\ t>0,}\\
{x_1(0,t) = x_1(\pi,t)=0;} \\
{x_{2,t}(z,t) =  q_2 x_{2,zz}(z,t) + \sqrt{|x_1(z,t)|} , \quad z \in (0,\pi),\ t>0,}\\
{x_2(0,t) = x_2(\pi,t)=0.}
\end {array}
\right.
\end{equation}
Here $q_1,q_2>0$ are diffusion coefficients.

This system can model a chemical reaction network with two reagents whose densities at a point $z$, and at time $t$ are given by $x_i(z,t)$, $i=1,2$.

We assume that $x_1 \in L^2(0,\pi)=:X_1$ and $x_2 \in L^4(0,\pi):=X_2$. The state of the whole system \eqref{eq:Gekoppelte-parabolic-sys-ISS} is $X:=X_1 \times X_2$.

We choose the following Lyapunov functions for subsystems 1 and 2, respectively:
\[
V_1(x_1)=\int_0^\pi{x_1^2(z)dz} = \|x_1\|^2_{L^2(0,\pi)},
\]
\[
V_2(x_2)=\int_0^\pi{x_2^4(z)dz} = \|x_2\|^4_{L^4(0,\pi)}.
\]
Assume for a while that the maps $x_i$ are twice continuously differentiable.
Consider the Lie derivative of $V_1$:
\begin{eqnarray*}
\dot{V}_1(x_1) &=& 2 \int_0^\pi{x_1(z) \left( q_1 x_{1,zz}(z) + x^2_2(z) \right)dz} \\
&\leq&
-2q_1 \left\|x_{1,z} \right\|^2_{L^2(0,\pi)} + 2\|x_1\|_{L^2(0,\pi)} \|x_2\|^2_{L^4(0,\pi)}.
\end{eqnarray*}
In the last estimate, we used integration by parts for the first term and the Cauchy-Schwarz inequality for the second one.

Applying Friedrichs' inequality \eqref{ineq:Friedrichs} to the first term, we obtain the estimate
\begin{eqnarray*}
\dot{V}_1(x_1) & \leq & -2q_1 \|x_1\|^2_{L^2(0,\pi)} + 2\|x_1\|_{L^2(0,\pi)} \|x_2\|^2_{L^4(0,\pi)} \\
&= &
-2q_1  V_1(x_1) + 2 \sqrt{V_1(x_1)} \sqrt{V_2(x_2)}.
\end{eqnarray*}
Take
\[
\chi_{12}(r)= \frac{1}{a} r,\; \forall r>0,
\]
with an arbitrary $a>0$. We obtain for all $x_1 \in X_1$, $x_2\in X_2$ that
\begin{eqnarray}
V_1(x_1) \geq \chi_{12}(V_2(x_2)) \qrq \frac{d}{dt}V_1(x_1) \leq
 -2(q_1-a^{\frac{1}{2}})  V_1(x_1).
\label{eq:V1dot_ineq}
\end{eqnarray}
The derivation was made under the assumption that
$x_1,x_2$ are twice continuously differentiable functions. The above estimate holds due to the density argument for general $x_1 \in L^2(0,\pi)$.
To ensure that $V_1$ is an ISS Lyapunov function, the following condition must hold:
\begin{eqnarray}
2(q_1-a^{\frac{1}{2}})>0 \quad\Iff\quad  a<q_1^2.
\label{eq:V1dot_ineq_koef}
\end{eqnarray}

Now assuming that $x_1$ and $x_2$ are smooth enough, consider the Lie derivative of $V_2$:
\begin{align*}
\dot{V}_2(x_2) &= 4 \int_0^\pi x_2^3(z)\Big( q_2 x_{2,zz}(z) + \sqrt{|x_1(z)|}\Big)dz \\
&\leq - 3q_2  \int_0^\pi 4 x_{2,z}^2(z) x_{2}^2(z)  dz + 4 \int_0^\pi  x_2^3(z) |x_1(z)|^{\frac{1}{2}} dz\\
&= - 3q_2  \int_0^\pi \Big( \frac{\partial}{\partial z}(x^2_2) \Big)^2 dz + 4 \int_0^\pi  x_2^3(z) |x_1(z)|^{\frac{1}{2}} dz.
\end{align*}

Applying Friedrichs' inequality to the first term (note that $x_2^2 \in L^2(0,d)$) and the H\"older's inequality to the last one, we obtain
\[
\dot{V}_2(x_2) \leq - 3q_2  V_2(x_2) +4 (V_2(x_2))^{3/4}(V_1(x_1))^{1/4}.
\]
Let
\vspace{-2mm}
\[
\chi_{21}(r)= \frac{1}{b} r,\; \forall r>0,
\]
where $b>0$ is an arbitrary constant. The following implication holds:
\[
V_2(x_2) \geq \chi_{21}(V_1(x_1)) \quad \Rightarrow \quad
\dot{V}_2(x_2) \leq  -(3q_2 - 4b^{\frac{1}{4}}) V_2(x_2).
\]
To ensure that $V_2$ is an ISS Lyapunov function for the second subsystem, the following condition must hold:
\begin{eqnarray}
3q_2 - 4b^{\frac{1}{4}}>0 \quad\Iff\quad  b<\Big(\frac{3q_2}{4}\Big)^4.
\label{eq:V2dot_ineq_koef}
\end{eqnarray}
To ensure the stability of the interconnection, we apply the
small-gain condition
\vspace{-1mm}
\begin{eqnarray}
\label{H_Function}
\chi_{12} \circ \chi_{21} < \id    \quad \Iff \quad ab>1.
\end{eqnarray}
If we impose the condition
\begin{eqnarray}
q_1^2\Big(\frac{3q_2}{4}\Big)^4 >1,
\label{eq:Conditions_on_q1_and_q2}
\end{eqnarray}
then it is possible to find $a,b$ so that the conditions \eqref{eq:V1dot_ineq_koef}, \eqref{eq:V2dot_ineq_koef}
and small-gain condition \eqref{H_Function} will be satisfied.
Thus, if $q_1$ and $q_2$ satisfy \eqref{eq:Conditions_on_q1_and_q2}, the small-gain theorem 
(Corollary~\ref{cor:SGT-2-systems}) guarantees the UGAS property of the interconnection.

\section{Interconnections of integral ISS systems}
\label{sec:ISS-SGT-iISS-couplings}

Small-gain theory can also be developed to study couplings of strongly integral input-to-stable systems. 

In this section, we state an iISS small-gain theorem for the semilinear system 
\begin{align}
\label{eq:interconnection}
\begin{array}{l}
\dot{x}_i(t)=A_ix_i(t)+f_i(x_1,x_2,u) , \quad i=1,2, \\
x_i(t)\in X_i , \quad u\in \Uc , 
\end{array}
\end{align}
where $X_i$ is the state space of the $i$-th subsystem, and $A_i:D(A_i) \subset X_i \to X_i$ is the generator of a strongly continuous semigroup over $X_i$. 
Let $X:=X_1\times X_2$ which is the space of $x=(x_1, x_2)$, 
and the norm on $X$ is defined as 
$\|\cdot\|_X=\|\cdot\|_{X_1}+\|\cdot\|_{X_2}$. 
In this section, we assume that 
there exist continuous functions $V_i:X_i \to \R_+$, 
$\psi_{i1},\psi_{i2} \in \Kinf$, $\alpha_i \in \K$, $\sigma_i \in \K$  and $\kappa_i \in \K\cup\{0\}$ 
for $i=1,2$ such that 
\begin{equation}
\label{LyapFunk_1Eigi}
\psi_{i1}(\|x_i\|_{X_i}) \leq V_i(x_i) \leq \psi_{i2}(\|x_i\|_{X_i}) \quad \forall x_i \in X_i,
\end{equation}
and the system \eqref{eq:interconnection} satisfies 
\begin{equation}
\label{GainImplikationi}
\dot{V}_i(x_i) \leq -\alpha_i(\|x_i\|_{X_i})
 + \sigma_i(\|x_{3-i}\|_{X_{3-i}}) + \kappa_i(\|u(0)\|_U)
\end{equation}
for all $x_i\in X_i$, $x_{3-i}\in X_{3-i}$ and $u\in \Uc$, 
where the Lie derivative of $V_i$ corresponding to the inputs 
$u\in \Uc$ and $v\in PC_b( \R_+,X_{3-i})$ with $v(0)=x_{3-i}$ 
is defined by
\begin{equation}
\label{LyapAbleitungi}
\dot{V}_i(x_i)=\mathop{\overline{\lim}} \limits_{t \rightarrow +0} {\frac{1}{t}\Big(V_i\big(\phi_i(t,x_i,v,u)\big)-V_i(x_i)\Big) }.
\end{equation}

To present a small-gain criterion for the interconnected 
system \eqref{eq:interconnection} whose components are not 
necessarily ISS, we make use of 
a generalized expression of inverse mappings on the set of 
extended nonnegative numbers $\clo{\R_+}:=[0,\infty]$. 
For $\omega\in\K$, define the function $\omega^\ominus$: 
$\clo{\R_+}\to\clo{\R_+}$ as 
\[
\omega^\ominus(s) := \sup \{ v \in \R_+ : s \geq \omega(v) \}=
\begin{cases}
\omega^{-1}(s) & \text{, if } s\in \text{Im}\ \omega, \\ 
+\infty & \text{, otherwise. } 
\end{cases}
\]
The function $\omega\in\K$ is extended to 
$\omega$: $\clo{\R_+}\to\clo{\R_+}$ as 
$\omega(s) := \sup_{v\in\{y\in\R_+ \, : \, y \leq s\}} \omega(v)$. 
These notations are useful for presenting the following result succinctly. 

\begin{theorem}[Integral ISS small-gain theorem]
\label{theorem:intercon}
Suppose that 
\begin{align}
\label{eq:alpsig}
\displaystyle\lim_{s\rightarrow\infty}\!\alpha_i(s)=\infty
\ \mbox{or} \ 
\lim_{s\rightarrow\infty}\!\sigma_{3-i}(s)\kappa_i(1)<\infty
\end{align}
is satisfied for $i=1,2$. 
If there exists $c>1$ such that 
\begin{align}
\label{eq:sg}
\psi_{11}^{-1}\circ\psi_{12}\circ
\alpha_1^\ominus\circ c\sigma_1\circ
\psi_{21}^{-1}\circ\psi_{22}\circ
\alpha_2^\ominus\circ c\sigma_2(s)
\leq s
\end{align}
holds for all $s\in\R_+$, then the system \eqref{eq:interconnection} is iISS. 
Moreover, if additionally $\alpha_i\in\Kinf$ for $i=1,2$, then the system \eqref{eq:interconnection} is ISS. 
Furthermore, one can construct $\lambda_i\in\K$ (see \cite[Theorem 6]{MiI15b} for details) so that
\begin{align}
\label{eq:Vsum}
V(x)=\int_0^{V_1(x_1)}\lambda_1(s)ds + \int_0^{V_2(x_2)}\lambda_2(s)ds
\end{align}
is an iISS (ISS) Lyapunov function for \eqref{eq:interconnection}.
\end{theorem}

\begin{proof}
Omitted, see \cite[Theorem 6]{MiI15b}.
\end{proof}

\begin{remark}
\label{rem:Simplified_SGC} 
Condition \eqref{eq:sg} can be called an iISS small-gain condition. At first glance, it may seem a bit technical. However, it simplifies considerably if $V_i(x_i) = \psi_i(\|x_i\|_{X_i})$, for some $\psi_1,\psi_2\in\Kinf$. In this case $\psi_{11}=\psi_{12}=\psi_i$, $i=1,2$ and \eqref{eq:sg} takes the form
\begin{align}
\label{eq:sg_2}
\alpha_1^\ominus\circ c\sigma_1\circ
\alpha_2^\ominus\circ c\sigma_2(s)
\leq s.
\end{align}
The term $\alpha_i^\ominus\circ \sigma_i$ can be interpreted as a Lyapunov gain of the $i$-th subsystem.
This interpretation justifies the name \q{small-gain condition} for \eqref{eq:sg}.
\end{remark}

\begin{remark}
\label{rem:Decay-rates-in-iISS-formulation} 
In Theorem~\ref{theorem:intercon}, we required that the decay rates of the iISS Lyapunov functions $\alpha_i$ are $\K$-functions. 
It was shown in \cite[Theorem 1]{CAI14} that the existence of an iISS Lyapunov function with such a decay rate implies not only iISS but also strong iISS. This result can be transferred similarly to infinite-dimensional systems.
Thus, in Theorem~\ref{theorem:intercon}, we implicitly assume that the subsystems are not only iISS but strongly iISS.
\end{remark}

\subsection{Example: coupling of an ISS and iISS nonlinear reaction-diffusion systems}
\label{sec:Example-for-SGT}

Consider the coupled nonlinear reaction-diffusion system
\begin{equation}
\label{GekoppelteNonLinSyst}
\left
\{
\begin{array}{l} 
{x_{1,t}(z,t) = x_{1,zz}(z,t) + x_1(z,t)x^4_2(z,t),} \\[1.5ex]
{x_1(0,t) = x_1(\pi,t)=0;} \\[1ex]
{x_{2,t}(z,t) =  x_{2,zz}(z,t) + ax_2(z,t) - bx_2(z,t) \Big(x_{2,z}(z,t)\Big)^{\!\!2} + \Big( \frac{x^2_1(z,t)}{1+x^2_1(z,t)} \Big)^{\!\frac{1}{2}}}\hspace{-1ex}, \\
{x_2(0,t) = x_2(\pi,t)=0.} 
\end{array}
\right.	
\end{equation}
defined on the region $(z,t) \in (0,\pi) \times (0,\infty)$. 
To fully define the system, we should choose the state spaces of subsystems. We take $X_1:=L^2(0,\pi)$ for $x_1(\cdot,t)$ and $X_2:=H^1_0(0,\pi)$ for $x_2(\cdot,t)$. 
We divide the analysis into three parts. First, we prove that the 
$x_1$-subsystem is iISS based on the developments in 
Section \ref{sec:parablic}. 
Next, we prove that the $x_2$-subsystem is ISS using the result 
in Section \ref{sec:ISS_Parabolic}. 
In the last part, we exploit the small-gain theorem presented 
in Section \ref{GekoppelteISS_Systeme} 
to derive conditions guaranteeing UGAS of the overall system \eqref{GekoppelteNonLinSyst}.

\subsubsection{The first subsystem is iISS}\label{sssec:exsys1}

First, we invoke item \ref{item:iISS_boundary} of 
Theorem~\ref{theorem:parabplic_iISS} 
with $q=1$ for $X_1=L^2(0,\pi)$. 
Then $W(h)=h^2$ and due to $\dot{W}(h) \leq 2W(h) |y_2|^4$ 
and $y_1(0,t) = y_1(\pi,t)=0$ for all $t\in\R_+$, 
one can choose
\begin{eqnarray}
V_1(y_1):=\ln\Big(1+\|y_1\|^2_{L^{2}(0,\pi)} \Big)
\label{eq:V1}
\end{eqnarray}
as an iISS Lyapunov function for $x_1$-subsystem. Its Lie derivative according to \eqref{Zderivative} satisfies 
\begin{eqnarray}
\dot{V}_1(y_1) \leq  
- \frac{2\|y_1\|_{L^2(0,\pi)}^2}{1+\|y_1\|_{L^2(0,\pi)}^2} 
+ 2\|y_2\|_{L^\infty(0,\pi)}^4 . 
\label{eq:Z1dot}
\end{eqnarray}
Note that we have put $1=\frac{\pi^2}{\pi^2}$ instead of $\frac{\pi^2}{4\pi^2}$ in formula \eqref{Zderivative} because $y_1=0$ holds at both ends of the interval $[0,\pi]$, and thus the less conservative  Friedrichs' inequality instead of Poincare's inequality can be used in getting \eqref{Zderivative}. 
To replace $L^\infty(0,\pi)$ with $X_2=H^1_0(0,\pi)$ 
for the input space used in \eqref{eq:Z1dot}, 
we recall \eqref{eq:agmonfried0}, which results in 
\begin{eqnarray}
\dot{V}_1(y_1) \leq 
- \frac{2\|y_1\|_{L^2(0,\pi)}^2}{1+\|y_1\|_{L^2(0,\pi)}^2} 
+ 8  \|y_2\|^4_{H^1_{0}(0,\pi)}.
\label{eq:Z1dot_Final}
\end{eqnarray}
Thus, we arrive at \eqref{GainImplikationi} for $i=1$, and the  
$x_1$-subsystem is iISS with respect to the state space 
$X_1=L^2(0,\pi)$ and the space of input values $X_2=H^1_0(0,\pi)$.

\subsubsection{The second subsystem is ISS}\label{sssec:exsys2}

We invoke Theorem~\ref{theorem:parabolic_ISS_W_12q_norm} with $q=1$. 
To simplify notation, we denote 
\[
u_2(z):=\Big( \frac{y^2_1(z)}{1+y^2_1(z)} \Big)^{1/2}.
\] 
Following \eqref{eq:LF_W_12q_norm}, we consider the following Lyapunov function candidate:
\begin{eqnarray}
V_2(y_2) :=  \int_0^{\pi}{\Big(y_{2,z}(z)\Big)^2  dz} =  \|y_2\|^2_{H^1_0(0,\pi)} .
\label{eq:LF_H10_Example}
\end{eqnarray}
Notice that the $x_2$-subsystem is of the form \eqref{Nonlinear_Parabolic_W_12q} 
with $c=1$ and 
\[
f(y_2(z),y_{2,z}(z))=ay_2(z) - by_2(z) (y_{2,z}(z))^2.
\]
To arrive at \eqref{eq:Assumption_on_f}, we obtain by integration by parts with $y_2(0) = y_2(\pi)=0$ that
\begin{align*}
\int_0^L y_{2,zz}(z)\Big( ay_2(z) &- by_2(z) \Big(y_{2,z}(z)\Big)^2 \Big) dz
=  - a V_2(y_2) - b\hspace{-.5ex}\int_0^L\hspace{-.5ex}y_{2,zz}(z)y_2(z) \Big(y_{2,z}(z)\Big)^2 dz.
\end{align*}
Due to $y_2(0) = y_2(\pi)=0$, we have
\begin{align*}
\int_0^\pi{ y_{2,zz}(z)  y_2(z) \Big(y_{2,z}(z)\Big)^2 dz}
&=
- \int_0^\pi{ y_{2,z}(z) \Big( 2  y_2(z) y_{2,z}(z)y_{2,zz}(z) +  \Big(y_{2,z}(z)\Big)^3  \Big) dz} 
\nonumber \\
&=
- \int_0^\pi{  2  y_2(z) \Big(y_{2,z}(z)\Big)^2y_{2,zz}(z) dz}
- \int_0^\pi{ \Big(y_{2,z}(z)\Big)^4 dz}, 
\end{align*}
which implies that
\begin{eqnarray}
\int_0^\pi{ y_{2,zz}(z)  y_2(z) \Big(y_{2,z}(z)\Big)^2 dz}
=
- \frac{1}{3} \int_0^\pi{ \Big(y_{2,z}(z)\Big)^4 dz}.
\label{Useful_Expression}
\end{eqnarray}
Thus, we arrive at \eqref{eq:Assumption_on_f} with  
\begin{align}
\eta(s)=-as+\frac{b}{3}s^2 , \quad \forall s\in\R_+, 
\end{align}
which is convex if $b\geq 0$. We also obtain 
\begin{align}
\hat{\alpha}(s)=(1-a-\epsilon)s+\frac{b}{3\pi}s^2
\label{eq:exhatalpha}
\end{align}
for \eqref{eq:Assumption_on_a}. The inequality in 
\eqref{eq:Assumption_on_a} is achieved for $\epsilon=1-a>0$ if 
$a<1$. Hence, if $a<1$ and $b\geq 0$ hold, Theorem~\ref{theorem:parabolic_ISS_W_12q_norm} with $q=1$ 
proves that for $\omega\in(0,2(1-a)]$, the function 
$V_2$ satisfies (with $\epsilon = 1-a$)
\begin{align}
\dot{V}_2 (y_2) \leq -2\Big(1-a-\frac{\omega}{2}\Big)V_2(y_2)
-\frac{2b}{3\pi}V_2(y_2)^2
+ \frac{1}{\omega} \|u_2\|^{2}_{L^{2}(0,\pi)}, 
\label{eq:x2issumod}
\end{align}
as in \eqref{eq:W_12q_Final_q_equals_2}, 
and $V_2$ is an ISS Lyapunov function for the $x_2$-subsystem 
with respect to the state space 
$X_2=H^1_0(0,\pi)$ and the input space 
$U_2=L^{2}(0,\pi)$. 
Since $s \mapsto {s}/({1+s})$ is a concave function of 
$s\in\R_+$, Jensen's inequality for concave functions \eqref{ineq:Jensen-concave} yields

\begin{eqnarray*}
\int_0^\pi\hspace{-.7ex} |u_2(z)|^2 dz 
=\int_0^\pi\hspace{-.7ex}{ \frac{y^2_1(z)}{1+y^2_1(z)} dz} 
\leq \pi \frac{ (1/\pi)\|y_1\|^2_{L^{2}(0,\pi)}}{1+(1/\pi)\|y_1\|^2_{L^{2}(0,\pi)}} 
=
\frac{ \pi\|y_1\|^2_{L^{2}(0,\pi)}}{\pi+\|y_1\|^2_{L^{2}(0,\pi)}} 
\leq
\frac{ \pi\|y_1\|^2_{L^{2}(0,\pi)}}{1+\|y_1\|^2_{L^{2}(0,\pi)}} . 
\end{eqnarray*}
Using this property in \eqref{eq:x2issumod}, we have 
\begin{align}
\dot{V}_2(y_2) \leq 
-2\Big(1-a-\frac{\omega}{2}\Big)\|y_2\|^2_{H^1_0(0,\pi)}
-\frac{2b}{3\pi}\|y_2\|^4_{H^1_0(0,\pi)}
+ \frac{\pi}{\omega}\left( 
\frac{\|y_1\|^2_{L^{2}(0,\pi)}}{1+\|y_1\|^2_{L^{2}(0,\pi)}}\right) . 
\label{eq:x2iss}
\end{align}
Therefore, $V_2$ is an ISS Lyapunov function for the $x_2$-subsystem 
with respect to the state space 
$X_2=H^1_0(0,\pi)$ and the input space $X_1=L^2(0,\pi)$.

Although property \eqref{eq:x2iss} is satisfactory for 
establishing UGAS of the overall 
system \eqref{GekoppelteNonLinSyst}, it is
worth mentioning that 
we can obtain a different estimate for 
$\dot{V}_2$.

Using the same derivations as in Theorem~\ref{theorem:parabolic_ISS_W_12q_norm}, we can obtain from \eqref{BasicEstimate_For_Other_U} the estimate
\begin{eqnarray}
\dot{V}_2(y_2) \leq -2(1-a)V_2(y_2) - \frac{2b}{3} \int_0^\pi{ \Big(y_{2,z}(z)\Big)^4 dz} -2 \int_0^\pi{ y_{2,zz}(z)u_2(z) dz}.
\label{eq:ISS_2nd_subsystem_alter}
\end{eqnarray}
Now we are going to estimate the last term in \eqref{BasicEstimate_For_Other_U} as a sum of certain functions of the input and of the state in a way to compensate the resulting \q{state-space} term not by means of linear dynamics as we have done before, but with the help of nonlinear dynamics of the system. 
Integrating the last term in \eqref{eq:ISS_2nd_subsystem_alter} by parts, and using 
$u_2(0) = u_2 (\pi)=0$, which follows from \eqref{GekoppelteNonLinSyst}, we see that
\begin{eqnarray}
- \int_0^\pi{ y_{2,zz}(z)u_2(z) dz} &=&
\int_0^\pi{ y_{2,z}(z)\frac{\partial u_2}{\partial z}(z) dz}  \nonumber\\
&\leq&
\int_0^\pi{ \left|  y_{2,z}(z)\right|\left|\frac{\partial u_2}{\partial z}(z)\right| dz}  \nonumber\\
&\leq&
\frac{\omega}{4} \int_0^\pi{ \Big|y_{2,z}(z)\Big|^4 dz}
+ 
\frac{1}{\omega^{\frac{1}{3}}} \frac{3}{4}  \int_0^\pi{\Big|\frac{\partial u_2}{\partial z}(z)\Big|^{\frac{4}{3}} dz},
\label{eq:InputTerm_Estimate}
\end{eqnarray}
where the last estimate is Young's inequality with $p=4$ and $\omega >0$.
Substitution of \eqref{eq:InputTerm_Estimate} into \eqref{eq:ISS_2nd_subsystem_alter} leads to
\begin{eqnarray}
\phantom{eee}
\dot{V_2}(y_2) \leq -2(1-a)V_2(y_2)  - \Big(\frac{2b}{3}-\frac{\omega}{2}\Big) \int_0^\pi{ \Big(y_{2,z}(z)\Big)^4 dz}
+ 
\frac{1}{\omega^{\frac{1}{3}}}  \frac{3}{2} \|u_2\|^{\frac{4}{3}}_{W^{1,\frac{4}{3}}_0(0,\pi)}.
\label{Last_Estim_W43_example}
\end{eqnarray}
For $\omega < \frac{4b}{3}$ and $a \leq 1$ the above estimate proves ISS of the second subsystem for 
$y_2,u_2 \in C_0^{\infty}(0,\pi) $. The density argument shows that the above inequality is valid for all $y_2 \in X_2 = H^1_0(0,\pi)$ and all $u_2 \in W^{1,\frac{4}{3}}_0(0,\pi)$.

The estimates \eqref{eq:x2iss} and \eqref{Last_Estim_W43_example}
demonstrate that different choices of input spaces result in 
different properties of a single system, even if the state space 
remains the same. 
The dissipative estimate \eqref{eq:x2iss} also illustrates that ISS does not necessarily imply an exponential decay rate. 
In contrast, in \eqref{Last_Estim_W43_example}, the constructed function $V_2$ is guaranteed to 
exhibit an exponential or faster decay rate globally.

Finally, it is worth noting that if $a >1$,  then 
$x=0$ of the linearization 
of the $x_2$-subsystem for $y_1\equiv 0$ is not UGAS (see \cite[Theorem 5.1.3]{Hen81}). 
Thus, $x=0$ of the nonlinear $x_2$-subsystem for $y_1\equiv 0$ also cannot be 
UGAS if $a >1$.

\subsubsection{Interconnection is UGAS}

Now, we collect the findings of the previous subsections. 
Assume that $a<1$ and $b\geq 0$. 
For the space $X:=L^2(0,\pi)\times H^1_0(0,\pi)$, 
the Lyapunov functions defined as \eqref{eq:V1} and 
\eqref{eq:LF_H10_Example} for the two subsystems satisfy 
\eqref{LyapFunk_1Eigi} 
with the class $\K_\infty$-functions 
$\psi_{11}=\psi_{12}: s \mapsto \ln(1+s^2)$ and $\psi_{21}=\psi_{22}: s \mapsto s^2$. 
Due to \eqref{eq:Z1dot_Final} and \eqref{eq:x2iss}, 
we have \eqref{GainImplikationi} for 
\begin{align}
&
\alpha_1(s)=\frac{2s^2}{1+s^2} 
, \quad  
\sigma_1(s)=8s^4
, \quad  
\kappa_1(s)=0,
\\
&
\alpha_2(s)=2\left(1\!-\!a\!-\!\frac{\omega}{2}\right)\!s^2 + \frac{2b}{3\pi}\!s^4
, \quad
\sigma_2(s)=\frac{\pi}{\omega}\left(\!\frac{s^2}{1+s^2}\!\right), 
\quad
\kappa_2(s)=0, 
\end{align}
defined with $\omega\in(0,2(1-a)]$. 
For these functions, condition \eqref{eq:sg} holds 
for all $s\in\R_+$ if and only if 
\begin{align}
\frac{12c^2\pi^2}{b\omega}\left(\!\frac{s^2}{1+s^2}\!\right) 
\leq 
\frac{2s^2}{1+s^2} 
 \quad \forall s\in\R_+
\label{eq:exsg}
\end{align}
is satisfied. Thus, there exists $c>1$ such that \eqref{eq:sg} holds 
if and only if ${6\pi^2}/{b}<\omega$ holds. 
Combining this with $\omega\in(0,2(1-a)]$, $a<1$ and $b\geq 0$, 
Theorem \ref{theorem:intercon} establishes UGAS
of $x=0$ for the whole system \eqref{GekoppelteNonLinSyst} provided that 
\begin{align}
a+\frac{3\pi^2}{b}<1 , \quad  b\geq 0. 
\label{eq:expibassum}
\end{align}
Note that \eqref{eq:alpsig} is satisfied. 
Due to the boundary conditions of $y_2$, Friedrichs' inequality ensures 
$\|y_2(\cdot,t)\|_{L^2(0,\pi)}\leq \|y_2(\cdot,t)\|_{H^1_0(0,\pi)}$. 
Thus, the UGAS guarantees the existence of $\beta\in\KL$ 
such that 
\begin{align}
\left\| \phi(t,x,0) \right\|_{L^2(0,\pi)\times L^2(0,\pi)}
\leq \left\| \phi(t,x,0) \right\|_X 
\leq  \beta(\left\| x \right\|_{X},t) 
\label{eq:exxestimate}
\end{align}
holds for all $x \in X=L^2(0,\pi)\times H^1_0(0,\pi)$ and all $t\in\R_+$.

\section{Concluding remarks}
\label{sec:Lyapunov-based-SGT-Concluding remarks}

The analysis of interconnections plays a remarkable role in the system theory, as it allows for establishing the stability of a complex system by analyzing the stability properties of its less complex components and the type of interaction between them. 
Small-gain theory is a powerful tool to study feedback interconnections ubiquitous in the control literature. 
The first small-gain results appeared for linear systems characterized by input-output gains. See \cite{DeV09}, \cite[Chapter 5]{Kha02} for an overview of classical results of this kind.
In \cite{Hil91,MaH92}, the small-gain theorems were extended to nonlinear feedback systems within the input-output context.

\textbf{ISS small-gain theorems for finite networks of finite-dimensional systems.} 
Nonlinear small-gain theorems for general feedback interconnections of two ODE ISS systems were introduced in \cite{JMW96,JTP94}. These results have been extended to networks composed of $ n \geq 2 $ finite-dimensional ISS systems in \cite{DRW10,DRW06b,DRW07}. 
The proof of the small-gain theorem \cite[Theorem 5.3]{DRW10} is based on the concept of Clarke generalized gradient of finite-dimensional maps, which is less suited for the analysis of networks with infinite-dimensional components. 
An alternative proof of the small-gain theorem for a coupling of $n \in \N$ finite-dimensional systems was given in \cite{Mir23}, based on the notion of Dini derivatives. 
In Theorem~\ref{Constr_of_LF_FinDim}, we have transferred this proof to the infinite-dimensional setting. 
A special case of Theorem~\ref{Constr_of_LF_FinDim} was shown in \cite{DaM13}, which was one of the first small-gain results for networks with infinite-dimensional subsystems (for other results of a quite different kind, see \cite{KaJ11}).

A summary of iISS small-gain theorems for couplings of ODE systems, and the challenges encountered on this way (e.g.,\ the insufficiency of max-type Lyapunov functions, which we have used in Theorem~\ref{Constr_of_LF_FinDim}), have been excellently explained in \cite{Ito13}. 
In \cite{Ito06,ItJ09,IJD13,AnA07,KaJ12}, small-gain theorems for interconnections whose subsystems are not necessarily ISS have been developed.

We refer to the corresponding chapters of \cite{Mir23} for an overview of available small-gain results for finite and infinite networks of finite-dimensional systems. See also \cite{DES11} for a survey of small-gain results that appeared till 2011. Applications of small-gain theorems have been studied in \cite{LJH14}.

\textbf{Cascade interconnections.} 
Consider an interconnection of two evolution equations of the form \eqref{eq:interconnection}, where the right-hand side $f_2$ does not depend on $x_1$. In other words, the second subsystem does not depend on the dynamics of the first system, but it influences the dynamics of the first system.
Such interconnections are called \emph{cascade interconnections}. 

The fact that a \emph{cascade interconnection of two input-to-state stable ODE systems is itself ISS} was already shown in \cite{Son89}. It can be obtained as a consequence of small-gain theorems in trajectory formulations (specialized for finite networks) for cascade interconnections of an arbitrary finite number of systems of a rather general nature.
In contrast to this, a cascade interconnection of two iISS systems is iISS only under some additional conditions, see \cite{AAS02, ChA08}. If these conditions are not met, a cascade interconnection of two iISS systems is not necessarily iISS, as shown in \cite[Example 1]{AAS02}. The failure of this important property was one of the motivations to introduce the strong iISS property in \cite{CAI14}. In \cite{CAI14b}, it was shown that cascade interconnections of strongly iISS ODE systems are again strongly iISS.
We refer to \cite[Chapter 4]{Mir23} for the unified exposition of iISS theory, in particular for cascade interconnections of finite-dimensional iISS systems.

\textbf{Finite networks of infinite-dimensional systems.} 
Within the infinite-dimensional ISS theory, generalizations of the results in \cite{JTP94,JMW96,DRW07,DRW10} to couplings of finitely many infinite-dimensional systems have been proposed in \cite{BLJ18,DaM13,Mir19b,TWJ12}. We refer to \cite{Mir19b} for more details and references on small-gain results for finite couplings. For the case of trajectory-based ISS small-gain theorems for finite networks, the main difficulties in going from finite to infinite dimensions stem from the fact that the characterizations of ISS developed for ODE systems in \cite{SoW96} are no more valid for infinite-dimensional systems. As argued in \cite{Mir19b}, more general characterizations shown in \cite{MiW18b} have to be used, which requires major changes in the proof of the small-gain result. 
In \cite{MKG20}, the findings from \cite{Mir19b} have been extended to the case of infinite interconnections of infinite-dimensional systems. In Chapter~\ref{chap:Infinite interconnections: Non-Lyapunov methods}, we follow \cite{MKG20} very closely. 
Integral ISS small-gain theorems and the examples presented in Section~\ref{sec:ISS-SGT-iISS-couplings}, were shown in \cite{MiI15b}.

Small-gain theorems for finite networks have been applied to the stability analysis of coupled parabolic-hyperbolic partial differential equations (PDEs) in \cite{KaK18}. Small-gain based boundary feedback design for global exponential stabilization of 1-D semilinear parabolic PDEs have been proposed in \cite{KaK19b}.

\textbf{Lyapunov-based small-gain theorems for infinite networks.}	Recently, a significant research effort was dedicated to the extension of the Lyapunov-based small-gain theorems discussed in this chapter to interconnections of an infinite number of subsystems 
\cite{DMS19a,DaP20,KMS21,KMZ22,MNK21,NMK22}. Here we briefly discuss these results.

In~\cite{DaP20}, it is shown that a countably infinite network of continuous-time input-to-state stable systems is ISS, provided that the gain functions capturing the influence of the subsystems on each other are all less than the identity, which is a very conservative condition. In~\cite{DMS19a}, it was shown that classic max-form strong small-gain conditions developed for finite networks in \cite{DRW10} do not ensure the stability of infinite networks, even {in the linear case}. To address this issue, more restrictive robust strong small-gain conditions {have been} developed in~\cite{DMS19a}. Still, the ISS small-gain for infinite networks of ODE systems has been shown in \cite{DMS19a} under the quite strong restriction that there is a linear path of strict decay for the gain operator, which makes the result not fully nonlinear.%


In contrast, for networks consisting of exponentially ISS systems, possessing exponential ISS Lyapunov functions with linear gains, it was shown in \cite{KMS21} that if the spectral radius of the gain operator $r(\Gamma)$ is less than one, then the whole network is exponentially ISS and there is a coercive exponential ISS Lyapunov function for the whole network. This result is tight and provides a complete generalization of \cite[Proposition~3.3]{DIW11} from finite to infinite networks. 
In \cite{GlM21}, numerous criteria for the condition $r(\Gamma)<1$ for general positive linear bounded operators have been shown (see Theorem~\ref{thm:stability-for-pos-ops} for one of such results) that make it possible to check this condition in practice efficiently. 

The small-gain theorem from \cite{KMS21} has been extended in~\cite{NMK22} to ISS with respect to closed sets and applied to the stability analysis of infinite time-variant networks and to the design of distributed observers for infinite networks.

The first fully nonlinear Lyapunov-based ISS small-gain theorem for networks with nonlinear gains and sup-preserving gain operators has been developed in \cite{KMZ22}. It states that \emph{if there exists a (possibly nonlinear) path of strict decay for the operator $\Gamma_\otimes$ (and some uniformity conditions hold), then the whole network is ISS, and a corresponding ISS Lyapunov function for the network can be constructed}. This result partially extends the nonlinear Lyapunov-based small-gain theorem for finite networks in maximum formulation proved in \cite{DRW10} to the infinite-dimensional setting; it contains as a special case the Lyapunov-based small-gain theorem for infinite networks in \cite{DMS19a}, and partially recovers the main result in \cite{DaP20}. See \cite{KMZ22} for the detailed discussion.

Furthermore, by refining the technique for the construction of paths of strict decay introduced in \cite{DMS19a}, it is shown in \cite{KMZ22} that \emph{the robust small-gain condition of the gain operator together with the global attractivity of the discrete-time dynamical system, induced by the gain operator, ensure the existence of a path of strict decay for the gain operator}, and thus the existence of an ISS Lyapunov function for the interconnection. 

In \cite{MNK21}, an ISS small-gain theorem for systems with linear gains and homogeneous and subadditive gain operators (including linear and max-linear operators) has been shown, which is firmly based on the nonlinear small-gain theorems in \cite{KMZ22}.
For homogeneous gain operators, one can introduce the concept of spectral radius via a definition analogous to the well-known Gelfand's formula for the spectral radius of linear operators. In this case, the existence of a path of strict decay for the gain operator is equivalent to the existence of a point of strict decay, which in turn is equivalent to the spectral small-gain condition. These findings extend the small-gain criteria of the spectral small-gain conditions obtained in \cite{GlM21} to the case of homogeneous and subadditive monotone operators.


\cleardoublepage

\chapter*{Conclusion}
\label{chap:Conclusion}


\addcontentsline{toc}{chapter}{\numberline{}Conclusion}

In this work, we have discussed a broad palette of techniques for analyzing the input-to-state stability of distributed parameter systems stemming from nonlinear control, operator and semigroup theory, Lyapunov theory, nonlinear networks, and partial differential equations (PDEs).

We have studied coercive and non-coercive Lyapunov methods, representing an indispensable tool for the stability analysis of nonlinear control systems. We have proved small-gain theorems that are a powerful tool for stability analysis of finite and infinite networks with nonlinear stable components.
Last but not least, we discussed admissibility methods that are ultimately useful for analyzing general linear boundary control systems. 

Despite tremendous progress in the infinite-dimensional ISS theory, the ISS analysis, as well as robust control and observation of PDEs with boundary inputs, or more generally of boundary control systems, remains highly challenging.
Already for linear BCS, there are no general converse ISS Lyapunov theorems. This raises serious concerns about the applicability of Lyapunov methods to the ISS analysis of BCS. 
Furthermore, for some ISS boundary control systems, coercive ISS Lyapunov functions are known and pretty standard. Yet, coercive ISS Lyapunov functions have not been constructed for other systems, including the heat equation with Dirichlet boundary input.
One of the most important open problems in infinite-dimensional ISS theory is a precise understanding of the reasons for this. 
Having solved this problem, it will be of interest to analyze nonlinear problems.
Well-posedness theory for nonlinear BCS is still in its infancy, and few results in the ISS analysis have been obtained. In particular, no global ISS results are available for such systems as Burgers' equation.

The author is sure that the next years will bring further breakthroughs in the ISS theory of linear and nonlinear distributed parameter systems, 
making ISS the firm foundation for the control of infinite-dimensional systems. And he hopes that this work will be helpful in this way.

\begin{appendices}
\renewcommand{\thechapter}{\Alph{chapter}}%
\setcounter{chapter}{0}
\cleardoublepage
\chapter{Compendium of finite-dimensional ISS theory}
\label{chap:Finite-dim_ISS_Theory}

For a modern and comprehensive account of the ISS theory of ODE systems, we refer a reader to the monograph \cite{Mir23}.
Another recent introduction to the ISS theory with stress on interconnected systems is a book \cite{LJH14}.
Applications of the ISS theory to the design of adaptive stabilizing controls to
nonlinear control systems are considered in \cite{KKK95}.
Here we state some important results and relate them to the results shown in this work.

Classical ISS theory of ordinary differential equations (ODEs) deals with the systems of the form
\begin{eqnarray}
\dot{x}=f(x,u).
\label{xdot=f_xu}
\end{eqnarray}

\begin{ass}
\label{ass:Existence-Uniqueness} 
We suppose throughout this appendix that
$f$ satisfies Assumption~\ref{Assumption1}, which amounts to the fact that $f$ is continuous on $\R^n \tm\R^m$, and is Lipschitz continuous on bounded balls with respect to the first argument. 
\end{ass}

Let inputs $u$ belong to the space $L^{\infty}(\R_+,\R^m)$ of essentially bounded Lebesque measurable functions.

Standard properties of the solution maps of ODEs justify that \eqref{xdot=f_xu} satisfies all the properties of the control systems with $X:=\R^n$, $U:=\R^m$ and $\Uc:=L^{\infty}(\R_+,U)$ and with $\phi(\cdot,x,u)$ given by the maximal solution of \eqref{xdot=f_xu} with an initial condition $x \in \R^n$ and an input $u\in\Uc$.

The topology of finite-dimensional spaces implies the following important property due to \cite[Proposition 5.1]{LSW96}:
\begin{theorem} 
\label{thm:boundedness_reachability_sets}
Let Assumption~\ref{ass:Existence-Uniqueness} hold.
If \eqref{xdot=f_xu} is forward complete, then \eqref{xdot=f_xu} has BRS.
\end{theorem}

The following result is useful for getting tight superposition theorems for ODE systems, which is due to \cite[Proposition 13]{MiW18b} that, in turn, is based on \cite[Corollary III.3]{SoW96}.
\begin{proposition}
\label{prop:ULIM_equals_LIM_in_finite_dimensions}
Let Assumption~\ref{ass:Existence-Uniqueness} hold and \eqref{xdot=f_xu} be
forward complete. Then \eqref{xdot=f_xu} is LIM if and only if it is ULIM.
\end{proposition}

\begin{definition}
\label{prop:ISS-LF_for_Lipschitz_Continuous_functions}
A Lipschitz continuous on bounded balls function $V:\R^n \to \R_+$ is called an \textit{ISS Lyapunov function} for \eqref{xdot=f_xu} if there exist
$\psi_1,\psi_2 \in \Kinf$, $\alpha \in \Kinf$ and $\chi \in \K$
such that
\begin{equation}
\label{LyapFunk_1Eig-ODEs}
\psi_1(|x|) \leq V(x) \leq \psi_2(|x|), \quad x \in X
\end{equation}
holds, and for points $x \in \R^n$ of differentiability of $V$ and for any input value $v \in \R^m$ the following implication holds:
\begin{eqnarray}
V(x) \geq \chi(|v|)    \qrq \nabla V(x) \cdot f(x,v) \leq -\alpha(V(x)).
\label{ImplicationIneq}
\end{eqnarray}
\end{definition}

\begin{definition}
\label{def:LISS-LF-ODEs}
A Lipschitz continuous function $V:D \to \R_+$, $0 \in \intt(D) \subset \R^n$  is called a \textit{LISS Lyapunov function},  if there exist $r >0$, $\psi_1,\psi_2 \in \Kinf$, $\alpha \in \Kinf$ and $\xi \in \K$
such that $B_r \subset D$,
\begin{equation}
\label{eq:LyapFunk_1Eig_LISS-ODEs}
\psi_1(|x|) \leq V(x) \leq \psi_2(|x|), \quad  x \in B_r
\end{equation}
and
\begin{eqnarray}
|x| \leq r,\ \|u\|_{\infty} \leq r \qrq \dot{V}_u(x) \leq -\alpha(V(x)) + \xi (\|u\|_\infty).
\label{eq:LISS_DissipativeIneq}
\end{eqnarray}
\end{definition}

Next, we characterize local ISS.
\begin{theorem}[Characterization of LISS]
\label{thm:Characterization_LISS-ODEs}
Let Assumption~\ref{ass:Existence-Uniqueness} hold and let $f(0,0)=0$.

For the system \eqref{xdot=f_xu} the following properties are equivalent:
\begin{enumerate}
    \item[(i)] 0-UAS
    \item[(ii)] Existence of a Lipschitz continuous 0-UAS Lyapunov function
    \item[(iii)] Existence of a Lipschitz continuous LISS Lyapunov function
    \item[(iv)] LISS
\end{enumerate}
\end{theorem}

\begin{proof}
The proof of this result can be found in \cite{Mir23}. If $\Uc=PC_b(\R_+,\R^m)$, then we can derive this result from Theorem~\ref{Characterization_LISS}. 
Indeed, in finite-dimensional case \eqref{eq:Estimate_f_concerning_u}  holds with 
\[
\sigma(r):=\sup_{v\in \Uc: \|v\|_U \leq r} \sup_{x\in X: \|x\|_X \leq \rho} \|f(x,v)-f(x,0)\|_X+r,
\]
due to the compactness of closed balls in $\R^m$. 

Furthermore, Assumption~\ref{ass:Existence-Uniqueness} implies that also Assumption~\ref{Assumption1} holds. 
Hence all assumptions of Theorem~\ref{Characterization_LISS} are fulfilled, and it implies the claim.
\end{proof}

A counterpart of the ISS superposition Theorem~\ref{thm:UAG_equals_ULIM_plus_LS} for ODE systems is given by the following result shown in \cite{SoW96}, see also \cite{Mir23}.
\begin{theorem}[ISS superposition theorem]
\label{Characterizations_ODEs}
Let \eqref{xdot=f_xu} be forward complete, $f(0,0)=0$, and let Assumption~\ref{ass:Existence-Uniqueness} hold.
The following statements are equivalent:
\begin{multicols}{2}
\begin{itemize}
    \item[(i)] \eqref{xdot=f_xu} is ISS.
    \item[(ii)] \eqref{xdot=f_xu} is UAG.
    \item[(iii)] \eqref{xdot=f_xu} is AG and ULS.
    \item[(iv)] \label{cond:ULIM_ULS_BRS_is_ISS} \eqref{xdot=f_xu} is LIM and ULS.
    \item[(v)] \eqref{xdot=f_xu} is LIM and UGS.
    \item[(vi)] \eqref{xdot=f_xu} is AG and UGS.
    \item[(vii)] \eqref{xdot=f_xu} is LIM and 0-ULS.
    \item[(viii)] \eqref{xdot=f_xu} is LIM and 0-UGAS.
    \item[(ix)] \eqref{xdot=f_xu} is AG and 0-UGAS.
\end{itemize}
\end{multicols}
\end{theorem}

\begin{proof}
By Theorem~\ref{thm:boundedness_reachability_sets}, forward completeness implies the BRS property. 
By Lemma~\ref{lem:RobustEquilibriumPoint}, \eqref{xdot=f_xu} has the CEP property.
Proposition~\ref{prop:ULIM_equals_LIM_in_finite_dimensions} ensures that LIM is equivalent to ULIM.
Now, the equivalences (i)-(vi) are a consequence of Theorem~\ref{thm:UAG_equals_ULIM_plus_LS}.

In case of $\Uc=PC_b(\R_+,\R^m)$, the remaining implications can be obtained by additionally invoking Theorem~\ref{thm:MainResult_Characterization_ISS_EQ_Banach_Spaces}.
\end{proof}

The integral characterization of the ISS property is as follows:
\begin{corollary}
\label{cor:ni-ISS-and-ISS-ODEs} 
Let $f$ in \eqref{xdot=f_xu} be Lipschitz continuous in both variables. 

Assume that \eqref{xdot=f_xu} is forward complete. Then the following statements are equivalent:
\begin{itemize}
	\item[(i)] \eqref{xdot=f_xu} is ISS
	\item[(ii)] \eqref{xdot=f_xu} is integral-to-integral ISS
	\item[(iii)] \eqref{xdot=f_xu} is norm-to-integral ISS
\end{itemize}
\end{corollary}

\begin{proof}
(i) $\Iff$ (ii). Was shown in \cite[Theorem 1]{Son98}. Here Lipschitz continuity in both variables is used as the proof of this equivalence in 
\cite[Theorem 1]{Son98} is based on the smooth ISS converse Lyapunov theorem in \cite{SoW95}, which requires at least Lipschitz continuity in both variables (actually, it was assumed that $f$ is continuously differentiable in $(x,u)$).

(i) $\Iff$ (iii). It is well-known that under made assumptions on $f$ a forward complete system \eqref{xdot=f_xu} satisfies the BRS property (by Theorem~\ref{thm:boundedness_reachability_sets}) and CEP property (by Proposition~\ref{lem:RobustEquilibriumPoint}). The claim follows from Theorem~\ref{thm:ncISS_LF_sufficient_condition_NEW}.
\end{proof}

\cleardoublepage
\chapter[Comparison functions]{Comparison functions and comparison principles}
\label{chap:Comparison functions and principles}
                
In this appendix, we introduce the classes of comparison functions that are very useful for developing the stability theory and derive several important properties of these functions. 
Next, we define the Dini derivatives and introduce several comparison principles based on the formalism of comparison functions. 

We omit the proofs of most results, which can be found in \cite{Mir23}, together with some further results and numerous exercises. 

\section{Comparison functions}
\label{sec:Comparison_Functions}

Recall that $C(\R_+)$ denotes the vector space of continuous maps from $\R_+$ to $\R_+$.
In this section, we study the properties of the following classes of comparison functions:
\index{comparison!functions}
\index{class!$\K$}
\index{class!$\LL$}
\index{class!$\KL$}
\index{class!$\PD$}
\index{class!$\Kinf$}
\index{function!positive definite}
\begin{equation*}
\begin{array}{ll}
{\PD} &:= \left\{\gamma \in C(\R_+): \gamma(0)=0 \mbox{ and } \gamma(r)>0 \mbox{ for } r>0  \right\}, \\
{\K} &:= \left\{\gamma \in \PD : \gamma \mbox{ is strictly increasing}   \right\},\\
{\K_{\infty}}&:=\left\{\gamma\in\K: \gamma\mbox{ is unbounded}\right\},\\
{\LL}&:=\{\gamma\in C(\R_+): \gamma\mbox{ is strictly decreasing with } \lim\limits_{t\rightarrow\infty}\gamma(t)=0 \},\\
{\KL} &:= \left\{\beta \in C(\R_+\times\R_+,\R_+): \beta(\cdot,t)\in{\K}\ \ \forall t \geq 0, \  \beta(r,\cdot)\in {\LL}\ \ \forall r > 0\right\}.
\end{array}
\end{equation*}
Functions of class $\PD$ are also called \emph{positive definite functions}.
Functions of class $\K$ will be frequently called also class $\K$ functions or simply $\K$-functions. The same convention we use with other classes of comparison functions.

Various stability properties can be restated in terms of comparison functions, and proofs of many results in stability theory can be significantly simplified by using standard manipulations with comparison functions.

\subsection{Elementary properties of comparison functions}
\label{sec:Elementary properties of comparison functions}

We start with a list of simple properties of $\K$ and $\LL$ functions. Here we denote by $\circ$ a composition of the maps, that is $f\circ g(\cdot):=f(g(\cdot))$.
\begin{proposition}[Elementary properties of $\K$- and $\LL$-functions]
\label{prop:Elementary_props_comparison_functions}
The following properties hold:
\begin{enumerate}[label=(\roman*)]
    \item \label{itm:Prop-comparison-fun-item1} For all $f,g \in \K$ it follows that $f \circ g \in \K$.
    \item \label{itm:Prop-comparison-fun-item2} For any $f \in \Kinf$ there exists $f^{-1}$, which also belongs to $\Kinf$.
    \item \label{itm:Prop-comparison-fun-item3} For any $f \in \K$, $g \in \LL$ it holds that $f \circ g \in \LL$ and $g \circ f \in \LL$.
    \item \label{itm:Prop-comparison-fun-item4} $(\Kinf,\circ)$ is a group.
    \item \label{itm:Prop-comparison-fun-item5} For any $\gamma_1,\gamma_2\in\K$ it follows that $\gamma_+: s\mapsto \max\{\gamma_1(s),\gamma_2(s)\}$
    and $\gamma_-: s\mapsto \min\{\gamma_1(s),\gamma_2(s)\}$ also belong to class $\K$.
    \item \label{itm:Prop-comparison-fun-item6} For any $\beta_1,\beta_2\in\KL$ it follows that $\beta_+: (s,t)\mapsto \max\{\beta_1(s,t),\beta_2(s,t)\}$ and $\beta_-: (s,t)\mapsto \min\{\beta_1(s,t),\beta_2(s,t)\}$ also belong to class $\KL$.    
\end{enumerate}
\end{proposition}

\begin{remark}
\label{rem:Inverse_K_function}
If $f \in \K$, then $f^{-1}$ is defined on $\im(f) \subset \R_+$. We will slightly abuse the notation by setting $f(r):=+\infty$, for $r \in \R_+ \backslash \im(f)$.
\end{remark}

The elementary but useful result on $\K$-functions is
\index{inequality!weak triangle}
\begin{proposition}[Weak triangle inequality] For any $\gamma \in \K$, and any $\sigma \in \Kinf$ it holds that:
\begin{eqnarray}
\gamma(a+b) \leq \max\big\{\gamma(a+\sigma(a)), \gamma(b+\sigma^{-1}(b))\big\}.
\label{eq:WeakTriangle_Ver1}
\end{eqnarray}
\end{proposition}

In particular, setting $\sigma:=\id$ (identity operator on $\R_+$) in \eqref{eq:WeakTriangle_Ver1} we obtain for any $\gamma \in\K$ a simple inequality
\[
\gamma(a+b) \leq \max\big\{\gamma(2a), \gamma(2b)\big\}.
\]

Another useful inequality is
\index{$\Kinf$-inequality}
\begin{proposition}[$\Kinf$-inequality] 
\label{prop:Kinf-inequality} 
For all $a,b \geq 0$, for all $\alpha \in \Kinf$ and all $r>0$ it holds that 
\begin{eqnarray}
\label{ineq:Kinf}
ab \leq r a \alpha(a)  + b \alpha^{-1}(\tfrac{b}{r}).
\end{eqnarray}
\end{proposition}

We have also
\begin{proposition}
\label{prop:Kfun_lower_estimates} 
For any $\alpha\in\K$ and any $a,b\geq0$ it holds that 
\begin{eqnarray}
\alpha(a+b)\geq \frac{1}{2}\alpha(a) + \frac{1}{2}\alpha(b).
\label{eq:Kfun_lower_bounds}
\end{eqnarray} 
\end{proposition}

\begin{lemma}
\label{lem:uSGC-lemma:id-Kinf} 
For any $\eta\in\Kinf$ with $\id - \eta\in\Kinf$ it holds that:
\begin{enumerate}[label = (\roman*)]
	\item\label{itm:id-eta-itm1} There is $\rho\in\Kinf$ such that $(\id-\eta)^{-1} = \id + \rho$.
	\item\label{itm:id-eta-itm2}  There are $\eta_1,\eta_2\in\Kinf$ such that $\id-\eta_1,\id-\eta_2\in\Kinf$ and $\id-\eta = (\id-\eta_1) \circ (\id-\eta_2)$.
\end{enumerate}
\end{lemma}

Analogously, the following holds
\begin{lemma}
\label{lem:uSGC-lemma:id+Kinf} 
For any $\eta\in\Kinf$ it holds that
\begin{enumerate}[label=(\roman*)]
	\item\label{itm:id+eta-itm1} There is $\rho\in\Kinf$ such that $(\id+\eta)^{-1} = \id - \rho$.
	\item\label{itm:id+eta-itm2}  There are $\eta_1,\eta_2\in\Kinf$ such that $\id+\eta = (\id+\eta_1) \circ (\id+\eta_2)$.
\end{enumerate}
\end{lemma}

\subsection{Sontag's $\KL$-Lemma}
\label{sec:Sontags-KL-Lemma}

One of the most useful results concerning comparison functions is a so-called Sontag's $\KL$-Lemma. In particular, it is very helpful for the proof of converse Lyapunov theorems.
\begin{proposition}[Sontag's $\KL$-Lemma]
\index{Sontag's $\KL$-Lemma}
\label{Sontags_KL_Lemma}
For each $\beta \in \KL$ and any $\lambda>0$ there exist $\alpha_1,\alpha_2 \in \Kinf$: 
\begin{eqnarray}
\alpha_1(\beta(s,t)) \leq \alpha_2(s) e^{-\lambda t}, \quad t,s \in \R_+.
\label{eq:KL-Lemma_Estimate}
\end{eqnarray}
\end{proposition}

We will also need a variation of Proposition~\ref{Sontags_KL_Lemma}:
\begin{corollary}
\label{Sontags_KL_Lemma-2}
For each $\beta \in \KL$ and any $\lambda>0$ there exist $\tilde\alpha_1,\tilde\alpha_2 \in \Kinf$: 
\begin{eqnarray}
\beta(s,t) \leq \tilde\alpha_1 (\tilde\alpha_2(s) e^{-\lambda t}) \quad \forall s,t \geq 0.
\label{eq:KL-Lemma_Estimate-2}
\end{eqnarray}
\end{corollary}

\subsection{Positive definite functions}
\label{sec:Positive definite functions}

In this section, we characterize positive definiteness in terms of comparison functions.

\index{function!proper}
\begin{definition}
\label{def:Proper-function} 
Let $X,Y$ be metric spaces. 
\begin{enumerate}[label=(\roman*)]
    \item A map $f:X \to Y$ is called \emph{proper} if the preimage of any compact subset of $Y$ is compact in $X$.
    \item We say that a \emph{sequence $(x_i)_{i \geq 1}$ escapes to infinity} if for every compact set $S \subset X$ only finitely many points of $(x_i)_{i \geq 1}$ are contained in $S$.
\end{enumerate}
\end{definition}

Properness can be characterized as follows:
\begin{proposition}
\label{prop:proper-escaping-to-infinity}
A continuous map $f:X\to Y$ is proper $\qiq$ for every sequence of points $(x_i)_{i \geq 1} \subset X$ escaping to infinity the sequence $(f(x_i))_{i \geq 1}$ also escapes to infinity.
\end{proposition}

As a corollary of Proposition~\ref{prop:proper-escaping-to-infinity}, we obtain the criterion of properness for $V\in C(\R^n,\R_+)$.
\begin{corollary}
\label{cor:criterion-properness} 
Let $V\in C(\R^n,\R_+)$ be given. The following assertions are equivalent:
\begin{enumerate}[label=(\roman*)]
	\item\label{itm:criterion-properness-1} 
 $V$ is proper.
	\item\label{itm:criterion-properness-2} 
 $V$ is \emph{radially unbounded}, that is $V(x) \to \infty$ as $|x| \to \infty$.
	\item\label{itm:criterion-properness-3} 
 There is $\psi \in\Kinf$: $V(x)\geq \psi(|x|)$ for all $x \in\R^n$.	
\end{enumerate}
\end{corollary}

We extend the definition of positive definite functions to functions defined on $\R^n$.
\index{function!positive definite}
\begin{definition}
\label{def:Positive-definiteness-proper} 
A map $V\in C(\R^n,\R_+)$ is called \emph{positive definite} if $V(0)=0$ and $V(x)>0$ for any $x \neq 0$. 
\end{definition}

Positive definiteness can be characterized in terms of comparison functions.
\begin{proposition}
\label{prop:Positive-definiteness-criterion} 
A map $\rho\in C(\R^n,\R_+)$ is positive definite if and only if $\rho(0)=0$ and there are $\omega \in\Kinf$ and $\sigma \in \LL$ such that 
\begin{eqnarray}
\rho(x) \geq \omega(|x|)\sigma(|x|),\quad x\in\R^n.
\label{eq:Lower-estimate-pdf-functions}
\end{eqnarray}
\end{proposition}

We close the section with the criterion of positive definiteness and properness
\index{function!radially unbounded}
\begin{proposition}
\label{prop:criteria_positive_definiteness_proper}
$V\in C(\R^n, \R_+)$ is positive definite and proper $\Iff$ there are $\psi_1,\psi_2 \in \Kinf$ such that 
\begin{eqnarray}
\psi_1(|x|)\leq V(x) \leq \psi_2(|x|),\quad x\in\R^n.
\label{eq:Sandwich-ineq}
\end{eqnarray}
\end{proposition}

\subsection{Approximations, upper and lower bounds}
\label{sec:Upper and lower bounds}

It is ultimately helpful to have conditions that ensure the possibility of estimating a given function $z:\R_+\to\R_+$ from above/below by a $\K$-function. In addition, smooth approximations of $\K$-functions are often needed. In this section, we provide desirable conditions of this kind. 

The next result will be used frequently:
\begin{proposition}
\label{prop:upper_lower_estimates}
Let $z:\R_+\to\R_+$ be nondecreasing and continuous at zero function satisfying $z(0)=0$ and $z(r)>0$ for $r>0$. 
Then:
\begin{enumerate}[label=(\roman*)]
   \item\label{itm:prop-1.2.7-i} there exist $z_1,z_2\in\K \bigcap C^\infty\big((0,+\infty)\big)$ so that 
\begin{eqnarray}
z_1(r)\leq z(r)\leq z_2(r) \quad r\geq 0.
\label{eq:Comp_Fun_Estimates}
\end{eqnarray}
    \item\label{itm:prop-1.2.7-ii} if in addition $\lim_{r\to\infty}z(r) = \infty$, then $z_1,z_2$ can be chosen from class $\Kinf$.
    \item\label{itm:prop-1.2.7-iii} if $z\in C(\R_+,\R_+)$, then for every $\eps>0$ the functions $z_1,z_2$ can be chosen so that $\sup_{s\in\R_+}|z_2(s)-z(s)| <\eps$ and $\sup_{s\in\R_+}|z_1(s)-z(s)| <\eps$.
\end{enumerate}
\end{proposition}

Next, we state a counterpart of a part of Proposition~\ref{prop:upper_lower_estimates} for estimates of functions $\psi:\R_+\times\R_+\to\R_+$ by a $\KL$-function.
\begin{proposition}
\label{prop:Upper_estimate_as_KL_function}
Let $\psi:\R_+ \times \R_+ \to \R_+$ be any function nondecreasing and continuous at 0 w.r.t.\  the first argument, nonincreasing w.r.t.\  the second argument and so that $\lim_{t\to\infty}\psi(r,t)= 0$ for any $r\geq 0$. Let also $\psi(0,t)=0$ for any $t\geq 0$. Then there exists $\beta\in\KL$: 
\[
\psi(r,t)\leq \beta(r,t) \quad \forall r,t\geq \R_+.
\]
\end{proposition}

The following result will be exploited for the proof of the converse Lyapunov theorems
\begin{lemma}
\label{lem:KinfLipschitzLowerEstimate}
The following statements hold:

\begin{enumerate}[label = (\roman*)]
	\item\label{itm:KinfLipschitzLowerEstimate-itm1} For any $\alpha \in\PD$, and for any $L>0$ there is 
$\rho \in \PD$ so that $\rho(s) \leq \alpha(s)$ for all $s \in \R_+$ and $\rho$ is globally Lipschitz with 
a Lipschitz constant $L$.

One of such $\rho: \R_+ \to \R_+$ is given by
\begin{eqnarray}
\rho(r):=\inf_{y\geq 0} \big\{\alpha(y) + L|y-r|\big\}.
\label{eq:LowerEstimRho_def}
\end{eqnarray}

	\item\label{itm:KinfLipschitzLowerEstimate-itm2} If in \ref{itm:KinfLipschitzLowerEstimate-itm1}  $\alpha\in \K$, then $\rho$ given by \eqref{eq:LowerEstimRho_def} is a $\K$-function.
	
	\item\label{itm:KinfLipschitzLowerEstimate-itm3} If in \ref{itm:KinfLipschitzLowerEstimate-itm1}  $\alpha\in \Kinf$, then $\rho$ given by \eqref{eq:LowerEstimRho_def} is a $\Kinf$-function.
\end{enumerate}
\end{lemma}

As a corollary of Lemma~\ref{lem:KinfLipschitzLowerEstimate}, we obtain
\begin{lemma}
\label{lem:uSGC-lemma3} 
For any $\alpha\in\Kinf$ there is $\eta\in\Kinf$ such that $\eta(r)\leq\alpha(r)$ for all $r\geq 0$, and $\id-\eta\in\Kinf$.
\end{lemma}

We proceed with a useful result on functions of 2 arguments, which are $\K$-functions in both arguments.
\begin{lemma}
\label{lem:KK-upperbounds}
Let $\gamma \in C(\R_+\tm\R_+,\R_+)$ be such that 
$\gamma(r,\cdot)\in\K$ for each $r>0$ and $\gamma(\cdot,s)\in\K$ for each $s>0$.
Then there is $\sigma\in\Kinf$, such that 
\begin{eqnarray}
\gamma(r,s)\leq \sigma(r)\sigma(s),\quad r,s\geq 0.
\label{eq:Upper-bounds-KK-functions}
\end{eqnarray}
\end{lemma}

We close this section with a result for upper bounds on families of $\Kinf$-functions.
\begin{lemma}
\label{lem:gamma_M-decomposition}
Let $(\gamma_M)_{M\in\N}$ be a family of class $\Kinf$ functions. Then, there exists $\sigma\in\Kinf$ such that
\begin{eqnarray}
\gamma_M(r)\leq \sigma(M)\sigma(r),\quad M\in\N,\ r\geq 0.
\label{eq:gamma_M-decomposition}
\end{eqnarray}
\end{lemma}

\section{Marginal functions}
\label{sec:Marginal functions}

We will often deal with the functions defined as suprema or infima of continuous functions over certain regions. Such functions are occasionally called marginal functions. In this section, we state the basic results concerning such functions.
We start with
\begin{proposition}
\label{prop:LemmaMaxFunk}
Let $(S_X,\rho_X)$ and $(S_Y,\rho_Y)$ be metric spaces with the metrics $\rho_X$ and $\rho_Y$, respectively, and $X \subset S_X$ be a compact set, and $Y \subset S_Y$ be locally compact (that is, each point of $Y$ has a compact neighborhood). Assume that $f:X \times Y \rightarrow \R$ is continuous on $X \times Y$. Then the function $g:y\mapsto \max\limits_ {x\in X}f(x,y)$ is continuous on $Y$.
\end{proposition}

We have the following corollary:
\begin{corollary}
\label{cor:sup_over_spheres}
Let $g:\R^n \to \R^p$, $n,p\in\N$ be a continuous function. Then $v,w:\R_+\to\R_+$ defined by
\begin{equation*}
v(r):= \sup_{|x|=r} |g(x)|,\quad w(r):= \inf_{|x|=r} |g(x)|,
\end{equation*}
are continuous functions.
\end{corollary}

\begin{lemma}
\label{lem:Max_Cont_Fun}
Let $g:\R_+ \to \R^p$, $p\in\N$ be a continuous function.
Then $\psi_1,\psi_2:\R_+ \to \R_+$ defined for $r\geq0$ by $\psi_1(r):=\inf_{x \geq r} |g(x)|$ and $\psi_2(r):=\sup_{x \leq r} |g(x)|$
are continuous and nondecreasing functions.
\end{lemma}

Now we can state a corollary of Proposition~\ref{prop:LemmaMaxFunk} and Lemma~\ref{lem:Max_Cont_Fun}:
\begin{corollary}
\label{lem:Max_xt_fun}
Let $g:\R_+ \times \R_+\to \R$ be a continuous function.
Then $\psi_1,\psi_2:\R_+ \times\R_+ \to \R_+$ defined by $\psi_1(r,t):=\inf_{x \geq r} |g(x,t)|$ and $\psi_2(r,t):=\sup_{x \leq r} |g(x,t)|$
are continuous functions, which are nondecreasing w.r.t.\  the first argument.
If in addition $g(x,t)\geq 0$ for all $x,t\in\R_+\times\R_+$ and $g(r,\cdot) \in\LL$ for each $r>0$, then also $\psi_2(r,\cdot)\in\LL$ for all $r>0$.
\end{corollary}

\begin{lemma}
\label{Lemma_Sup_Inf}
Let $g:\R^n \to \R^p$, $n,p\in\N$ be a continuous function and $g(0)=0$. 
Then $\psi_1,\psi_2:\R_+ \to \R_+$ defined by $\psi_1(r):=\inf_{|x| \geq r} |g(x)|$ and $\psi_2(r):=\sup_{|x| \leq r} |g(x)|$
are well-defined, continuous, nondecreasing and $\psi_1(0)=\psi_2(0)=0$.
\end{lemma}

\begin{lemma}
\label{lem:Lower-bounds-of-PD-funtions-with-nonzero-limit} 
Let $\xi \in\PD$ and $\Liminf_{r\to\infty}\xi(r) >0$. Then there is a $\psi\in\K$: $\psi(r) \leq \xi(r)$ for $r\in\R_+$
and $\lim_{r\to\infty}\psi(r) = \Liminf_{r\to\infty}\xi(r)$.
In particular, if $\Liminf_{r\to\infty}\xi(r) = \infty$, then $\psi\in\Kinf$.
\end{lemma}

\section{Dini derivatives}
\label{sec:Dini derivatives}

As we know, absolutely continuous functions are differentiable almost everywhere.
However, sometimes one has to \q{differentiate} functions, which are merely continuous. 
To do this, we exploit Dini derivatives.
\index{derivative!Dini}
\index{Dini derivative}
For a continuous function $y:\R \to \R$, define 
the right upper Dini derivative and 
the right lower Dini derivative respectively by 
\[
D^+y(t):=\Limsup_{h \to +0}\frac{y(t+h)-y(t)}{h},
\quad  
D_+y(t):=\Liminf_{h \to +0}\frac{y(t+h)-y(t)}{h}.
\] 

We start with simple properties of Dini derivatives
\begin{lemma} 
\label{lem:Properties_Dini_Derivatives}
The following properties of Dini derivatives hold:
\begin{enumerate}[label=(\roman*)]
    \item\label{itm:Lem1.4.1-i} For all $f,g \in C(\R)$ it holds that $D^+(f + g) \leq D^+(f) + D^+(g)$. 
    \item\label{itm:Lem1.4.1-ii} For any $f \in C(\R)$ it holds that $D^+(-f) = -D_+(f)$.
\end{enumerate}
\end{lemma}

We proceed with the following generalization of the chain rule: 
\begin{lemma}
\label{Lemma2}
Let $g$ be a continuous function on a bounded interval $[a,b]$ and let $t$ be such that $f: g([a,b]) \to\R$ is continuously differentiable at $g(t)$ and
$\frac{df}{dr}(r)\Big|_{r=g(t)} > 0$.
Then
\[
D^{+}(f\circ g)(t) = \frac{df}{dr}(r)\Big|_{r=g(t)} D^{+}g(t).
\]
\end{lemma}

\begin{proof}
Consider a function $h: t\in\R \to \R$ that is continuous at $t=c$ 
and satisfies $h(c) > 0$. Let $k:\R \to \R$ be an arbitrary continuous function. 
Then the following holds (allowing the limits to be equal $\pm \infty$)
\[
\Limsup_{t \to c}\big(h(t)k(t)\big) = \lim_{t\to c}h(t) \cdot \Limsup_{t \to c}k(t).
\]
Then for $f$, $g$, and $t$ as in the statement of the lemma, the following holds
\begin{eqnarray*}
D^{+}(f\circ g)(t) &=& \Limsup_{h \to +0}\frac{f \circ g(t+h)- f \circ g(t)}{h} \\
                                 &=& \Limsup_{h \to +0}\frac{f ( g(t+h) - g(t) + g(t))- f \circ g(t)}{g(t+h) - g(t)} \cdot \frac{g(t+h) - g(t)}{h} \\
                                 &=& \lim_{h \to +0}\frac{f ( g(t+h) - g(t) + g(t))- f \circ g(t)}{g(t+h) - g(t)}\Limsup_{h \to +0}\frac{g(t+h) - g(t)}{h} \\
                                 &=&   \frac{df}{dr}(r)\Big|_{r=g(t)} D^{+}g(t).
\end{eqnarray*}
\end{proof}

We also state the following version of the fundamental theorem of calculus (which is a special case of \cite[Theorem 7.3, pp. 204--205]{Sak47}, \cite[Theorem 9]{HaT06}, see also \cite[Theorem 2.1]{Sza65}).
\begin{proposition}
\label{prop:fundamental theorem of calculus-Dini-derivative} 
Let $f: [a,b] \to \R$ be continuous, and $g: [a,b] \to \R$ be locally Lebesgue integrable and bounded on $[a,b]$.
Suppose that 
\begin{eqnarray}
D_+f(t) \leq g(t), \quad t\in[a,b].
\label{eq:Bound-on-Dini-derivative}
\end{eqnarray}
Then for all $t\in[a,b]$ we have
\begin{eqnarray}
f(t) - f(a) \leq \int_a^t g(s) ds.
\label{eq:Main-theorem-of-calculus}
\end{eqnarray}
\end{proposition}

We close the section with a lemma on derivatives of monotone functions.

\begin{lemma}
\label{lem:Integrals-of-monotone-functions}
Let $b:\R_+\to\R$ be a nonincreasing function. Then for each $t\in\R_+$ it holds that
\begin{eqnarray}
b(t) \geq D^+ \int_0^t b(s) ds \geq D_+ \int_0^t b(s) ds \geq \lim_{h\to +0}b(t+h).
\label{eq:Lower-Dini-Derivative}
\end{eqnarray}
\end{lemma}

\begin{proof}
Pick any $t\geq 0$. By definition of a Dini derivative, it holds that
\begin{eqnarray*}
D^+ \int_0^t b(s) ds
&=& \mathop{\overline{\lim}} \limits_{h \rightarrow +0} \frac{1}{h}\Big({\int_0^{t+h} b(s) ds - \int_0^t b(s) ds}\Big)\\
&=& \mathop{\overline{\lim}} \limits_{h \rightarrow +0} \frac{1}{h}{\int_t^{t+h} b(s) ds}\\
&\leq& \mathop{\overline{\lim}} \limits_{h \rightarrow +0} \frac{1}{h}{\int_t^{t+h} b(t) ds} = b(t).
\end{eqnarray*}
On the other hand, we have that
\begin{eqnarray*}
D_+ \int_0^t b(s) ds
&\geq& \mathop{\underline{\lim}} \limits_{h \rightarrow +0} \frac{1}{h}{\int_t^{t+h} b(t+h) ds}
= \mathop{\underline{\lim}} \limits_{h \rightarrow +0} b(t+h).
\end{eqnarray*}
The inequality $D^+ \int_0^t b(s) ds \geq D_+ \int_0^t b(s)ds$ is clear.
\end{proof}

\section{Comparison principles}
\label{sec:Comparison principles}

We want to state a comparison principle, giving a uniform estimate for solutions of certain differential inequalities.
This result can be understood as a nonlinear version of the Gronwall's inequality (Lemma~\ref{lem:Gronwall}).
\index{comparison!principle}
\begin{proposition}[Comparison principle]
\label{prop:ComparisonPrinciple}
For any $\alpha \in \PD$ there exists a $\beta \in \KL$ so that for any $y \in C(\R_+, \R_+)$ satisfying the differential inequality
\begin{eqnarray}
D^+y(t) \leq -\alpha(y(t)) \quad \forall t>0,
\label{eq:ComparisonPrinciple}
\end{eqnarray}
it holds that
\begin{eqnarray}
y(t) \leq \beta(y(0),t) \quad \forall t \geq 0.
\label{eq:ComparisonPrinciple_FinalEstimate}
\end{eqnarray}
\end{proposition}

We also state a comparison principle for systems with inputs:
\begin{proposition}
\label{prop:comparison-principle-with-inputs-max-form} 
For any $\rho\in\PD$ there is $\beta\in\KL$ such that
for any $\tilde{t}\in(0,+\infty]$, any $y\in C([0,\tilde{t}),\R)$ with $y(0)\geq0$ and any $v \in C([0,\tilde{t}),\R_+)$
if 
\begin{eqnarray}
D^+y(t) \leq -\rho(\max\{y(t) + v(t),0\})
\label{eq:comparison-principle-with-max-premise}
\end{eqnarray}
holds for all $t \in[0,\tilde{t})$, then denoting by $v_t$ the restriction of $v$ to $[0,t]$, we have 
\begin{eqnarray}
y(t) \leq \max\{\beta(y(0),t),\|v_t\|_\infty\},\quad t\in[0,\tilde{t}).
\label{eq:comparison-principle-with-max-endestimate}
\end{eqnarray}
\end{proposition}

The final result of this section is

\begin{proposition}
\label{prop:comparison-principle-with-inputs} 
For any $\rho\in\PD$ there is $\beta\in\KL$ such that
for any $\tilde{t}\in(0,+\infty]$, any $y\in C([0,\tilde{t}),\R_+)$ and any $v \in PC_b([0,\tilde{t}),\R_+)$ if 
\begin{eqnarray}
D^+y(t) \leq -\rho(y(t)) + v(t)
\label{eq:comparison-principle-premise}
\end{eqnarray}
holds for all $t \in[0,\tilde{t})$, then 
\begin{eqnarray}
y(t) \leq \beta(y(0),t) + \int_0^t2v(s)ds,\quad t\in[0,\tilde{t}).
\label{eq:comparison-principle-endestimate}
\end{eqnarray}
\end{proposition}

\begin{proof}
We assume that \eqref{eq:comparison-principle-premise} holds for all $t \in[0,\tilde{t})$.

By Lemma~\ref{lem:KinfLipschitzLowerEstimate}, there is a globally Lipschitz continuous $\rho_1 \in\PD$, such that $\rho\geq \rho_1$ pointwise.
Consider the following initial value problem:
\begin{eqnarray}
\dot{w}(t) = - \rho_1(w(t)) + v(t),\quad w(0)=y(0).
\label{eq:Auxiliary-problem}
\end{eqnarray}
As $\rho_1$ is globally Lipschitz, and $v$ is piecewise continuous, there is a unique piecewise continuously differentiable solution $w$ of \eqref{eq:Auxiliary-problem}.
Take some $\hat{t}\in (0,\tilde{t})$, such that $w$ is continuously differentiable on $(0,\hat{t})$.

We have
\begin{eqnarray}
D^+ y(t) - \dot{w}(t) \leq -\big(\rho_1(y(t))- \rho_1(w(t))\big),\quad \text{for all } t\in[0,\hat{t}).
\label{eq:comparison-principle-estimate}
\end{eqnarray} 
Assume that there is $\tau \in [0,\hat{t})$ such that $y(\tau)>w(\tau)$. Pick the maximal in $[0,\hat{t})$ interval $(s,s+\varepsilon)$ containing $\tau$ such that $y(t)>w(t)$ for all $t \in (s,s+\varepsilon)$. By maximality of the interval  $(s,s+\varepsilon)$, continuity of $y,w$ and as $y(0)=w(0)$, we have that $y(s)=w(s)$. 

As $\rho_1$ is globally Lipschitz, there is some $L>0$ such that for all $t\in(s,s+\varepsilon)$ we have
\begin{eqnarray*}
D^+ \big(y(t) - w(t)\big) \leq \big|\rho_1(y(t))- \rho_1(w(t))\big| \leq L |y(t)- w(t)| = L (y(t)- w(t)).
\end{eqnarray*} 
Define $z(t):=y(t) - w(t)$ on $(s,s+\varepsilon)$. 
By Proposition~\ref{prop:fundamental theorem of calculus-Dini-derivative}, this implies that
\[
z(t)-z(s) \leq L\int_s^tz(\tau)d\tau,\quad t\in[s,s+\varepsilon],
\]
which by Gr\"onwall's inequality leads to
\[
z(t) \leq z(s) e^{L(t-s)},\quad t\in[s,s+\varepsilon].
\]
As $z(s)=0$, we come to a contradiction to $z(t)>0$ for $t\in (s,s+\varepsilon)$.

From this fact and due to the continuity of $y,w$, one easily obtains that $0\leq y(t) \leq w(t)$ for all $t \in [0,\tilde{t})$.

Define for $t\in[0,\tilde{t})$ the following maps:
\[
v_1(t):=\int_0^t v(s) ds,\qquad w_1(t) = w(t) - v_1(t).
\]
The functions $w$, $v_1$ and $w_1$ are differentiable, and thus for $t\in[0,\tilde{t})$ it holds that 
\[
\dot{w}_1(t) = \dot{w}(t) - \dot{v}_1(t) = \dot{w}(t) - v(t) = - \rho_1(w(t)).
\]
Thus, for almost all $t \in [0,\tilde{t})$, we obtain
\[
\dot{w}_1(t) = - \rho_1(\max\{w_1(t)+v_1(t),0\}).
\] 
Proposition~\ref{prop:comparison-principle-with-inputs-max-form} ensures that there is $\beta\in\KL$ depending on $\rho_1$ solely, such that for all $t \in [0,\tilde{t})$ 
\[
w_1(t) \leq \max\{\beta(w_1(0),t),\|v_{1,t}\|_\infty\}.
\]
As $w_1(0)=w(0)=y(0)$, we infer for all $t \in [0,\tilde{t})$ that
\begin{eqnarray*}
y(t) \leq  w(t) = w_1(t) + v_1(t)
		 &\leq & \beta(y(0),t) +\|v_{1,t}\|_\infty + v_1(t)\\
		 &= & \beta(y(0),t) + 2\int_0^t v(s) ds.
\end{eqnarray*}
The result is shown.
\end{proof}

\section{Concluding remarks}

This Appendix is a shortened version of the appendix on comparison functions and comparison principles in \cite{Mir23}, where you can find the proofs of most results, some further results as well a number of exercises that extend the material. 
Some further properties of comparison functions can be found in a survey \cite{Kel14}, which is a nice complement to the material described in \cite{Mir23}. Some general ways for recasting $\eps-\delta$ definitions into comparison function formalism are analyzed in \cite{RaR15}. 

\textbf{Comparison functions.} 
Classes of comparison functions were introduced by Jose Massera in the 1950s and were widely used by Wolfgang Hahn in his book \cite{Hah67} to characterize stability notions of dynamical systems. The usefulness of these classes of functions in the ISS theory and other contexts made their use standard within the systems theoretic society. 

Lemma~\ref{lem:uSGC-lemma:id+Kinf} is due to \cite[Lemma~1.1.5, Lemma~1.1.3]{Rue07}. 
Proposition~\ref{Sontags_KL_Lemma} is due to Sontag, see \cite[Proposition 7]{Son98}. We employed the proof due to \cite[Lemma 3, p.326]{TeP00}.
Proposition \ref{prop:Positive-definiteness-criterion} was shown for scalar positive definite functions in \cite[Lemma IV.1]{ASW00} and extended to general positive definite functions in \cite[Lemma 18]{Kel14}. 
Proposition~\ref{prop:upper_lower_estimates} generalizes \cite[Lemma 2.5]{CLS98} (estimates from below/above by infinitely differentiable $\K$-functions) and \cite[Lemma 1.1.6]{Rue07} (approximation by smooth functions staying in $\eps$-neighborhood of a given function) and unifies these two features. 
Lemma~\ref{lem:KinfLipschitzLowerEstimate} is due to \cite[p.130]{KaJ11b}.
Lemma~\ref{lem:gamma_M-decomposition} was shown in \cite[Lemma 3.9]{ASW00b}.
Lemma~\ref{prop:LemmaMaxFunk} is well-known. A reader may consult \cite[Theorem 1.4.16]{AuF09} for more general results of this kind.

\textbf{Dini derivatives.} Dini has introduced the derivatives of continuous functions in \cite{Din1878}. For more on this theory we refer to \cite{HaT06}, \cite{CLS98}, \cite[Section 6]{BrL66}, \cite[Chapter 3]{KaK12}.

\textbf{Comparison principles.} 
For absolutely continuous $y$, Proposition~\ref{prop:ComparisonPrinciple} was proved in \cite[Lemma 4.4]{LSW96}, Proposition~\ref{prop:comparison-principle-with-inputs-max-form} was shown in \cite[Lemma IV.2]{ASW00}, and Proposition~\ref{prop:comparison-principle-with-inputs} was shown in \cite[Corollary IV.3]{ASW00}.
We have extended it to continuous $y$ by using the Dini derivatives formalism.
 In the time-delay realm, a notable result closely connected to comparison principles is a nonlinear Halanay's inequality \cite{Pep21}.

%
%
%
%
%

\cleardoublepage
\chapter{Ordered Banach spaces and positive operators}
\label{section:ordered-banach-spaces}

In this appendix, we recall some background knowledge on ordered Banach spaces, which is essential for us for the development of the
small-gain theorems for infinite networks.

\section{Ordered Banach spaces}

\begin{definition}
\label{def:Partial_order}
\index{partial order}
Let $X$ be any set. \emph{A partial order} on $X$ is a binary relation \q{$\leq$} that is:
\begin{itemize}
	\item[(i)] \emph{reflexive}: $x \leq x$ for all $x \in X$,
	\item[(ii)] \emph{transitive}: $x \leq y$ and $y \leq z$ implies $x \leq z$,
	\item[(iii)] \emph{antisymmetric}: $x \leq y$ and $y \leq x$ implies $x = y$.
\end{itemize}
\end{definition}

\begin{definition}
\label{def:Partial_pre-order}
\index{partial pre-order}
Let $X$ be any set. \emph{A partial pre-order} on $X$ is a binary relation \q{$\leq$} that is reflexive and transitive.
\end{definition}

Consider a vector space $X$ over $\R$. 
For two sets $A, B \subset X$ we define 
$A + B := \{ a + b : a \in A, b \in B \}$, $- A := \{- a : a \in A \}$, and $\R_+ \cdot A := \{ r \cdot a: a \in A, r \in \R_+\}.$
Given a topology, by $\intt A$, we denote the \emph{(topological) interior of $A$} and by $\clo{A}$ its \emph{closure}.

\begin{definition}
\label{def:Positive_cone}
\index{positive cone}
We say that $X^+\subset X$ is a \emph{(positive) cone} in a vector space $X$ if the following conditions hold:
\begin{itemize}
	\item[(i)]   $X^+ \cap (-X^+)=\{0\}$,
	\item[(ii)]  $\R_+ X^+ \subset X^+$,
	\item[(iii)] $X^+ + X^+ \subset X^+$.
\end{itemize}
$X^+$ is called a \emph{wedge} if $X^+$ satisfies the conditions (ii) and (iii).
\index{wedge}

A wedge (and, in particular, a cone) $X^+$ introduces the binary relation ``$\leq$'' on $X$ as $x \leq y\ \Iff \  y-x \in X^+$.
\end{definition}

\begin{proposition}
\label{prop:Wedge-generates-pre-order}
If $X^+$ is a wedge, then \q{$\leq$} is a partial pre-order relation; and it is a partial order if $X^+$ is a cone.
\end{proposition}

\begin{proof}
Let $X^+$ be a wedge. Then $0 \in X^+$ by axiom (ii) of Definition~\ref{def:Positive_cone}.
As $0 = x-x \in X^+$ for all $x \in X$, then $x \leq x$ holds for all $x \in X$, and \q{$\leq$} is reflexive.
Furthermore, if $x \leq y$ and $y\leq z$, then $y-x \in X^+$ and $z-y \in X^+$, and by axiom (iii) $(y-x)+(z-y) = z-x \in X^+$. Thus, $x \leq z$.

Let $X^+$ be a cone. Let $x,y \in X$ be such that $x\le y$ and $y \le x$. Then $y-x \in X^+$ and $y-x \in - X^+$, and as $X^+$ is a cone, it holds that $x=y$.
\end{proof}

\begin{definition}
\label{def:Ordered-set}
\index{OVS}
\index{ordered vector space}
A real vector space with a linear order is called an \emph{ordered vector space (OVS)}.
\end{definition}

\begin{definition}
\label{def:OBS}
\index{OBS}
\index{ordered Banach space}
Let $X$ be a Banach space ordered by a cone $X^+$. Then $(X,X^+)$ is called an \emph{ordered Banach space (OBS)}, if the cone $X^+$ is closed.
\end{definition}

\index{order interval}
\begin{definition}
\label{def:order interval} 
For any two vectors $x,z$ in an ordered Banach space, we call the set $[x,z] := \{y \in X: \, x \le y \le z\}$ the 
\emph{order interval} between $x$ and $z$; this order interval is nonempty if and only if $x \le z$.
\end{definition}

\section{Systematics of cones} 

\begin{definition}
\label{def:Total_cone}
\index{cone!total}
\index{cone!generating}
Let $(X,X^+)$ be an OBS. The cone $X^+$ is called
\begin{itemize}
	\item[(i)] \emph{total} (or \emph{spatial}) if $\clo{X^+ + (-X^+)} = X$.
	\item[(ii)] \emph{generating} (or \emph{reproducing}) if $X^+ + (-X^+) = X$ (that is, if $X^+$ spans $X$).
\end{itemize}
\end{definition}

In other words, the cone is generating if and only if each vector $x \in X$ can be decomposed as $x = y-z$ for two vectors $y,z \in X^+$. This decomposition can always be performed in a way that controls the norms of $y$ and $z$. 
Indeed, if the cone $X^+$ is generating, there exists $M > 0$ such that for each $x \in X$ there exist $y,z \in X^+$:
\begin{align}
	\label{eq:bounded-decomposition}
	x = y-z \qquad \text{and} \qquad \norm{y}, \norm{z} \le M \norm{x};
\end{align}
see for instance \cite[Theorem~2.37(1) and~(3)]{AlT07}. 

We adopt the following terminology from \cite[p.~627]{Ama76}.
\begin{definition}
\label{def:Norm_monotonicity}
Let $(X,X^+)$ be an ordered Banach space. The norm in $X$ is called:
\begin{itemize}
\index{norm!semimonotone}
	\item[(i)] \emph{semimonotone}, if there is $C>0$ so that
\begin{align}
	\label{eq:normality-constant}
	\norm{x} \le C \norm{y} \qquad \text{whenever } 0 \le x \le y.
\end{align}
\index{norm!monotone}
\index{cone!normal}
		\item[(ii)] \emph{monotone}, if for any $x,y \in X$: $0\leq x \leq y$ it follows that $\|x\|_X \leq \|y\|_X$.
\end{itemize}
If the norm in $X$ is semimonotone, the cone $X^+$ is called \emph{normal}.
\end{definition}
The cone is normal if and only if there exists an equivalent monotone norm on $X$ \cite[Theorem~2.38(1) and~(2)]{AlT07}.

\begin{definition}
\label{def:cone-with-a nonempty interior} 
\index{cone!with nonempty interior}
If $\intt(X^+)\neq\emptyset$, then we say that the \emph{cone $X^+$ has a nonempty interior}. 
\end{definition}

It is worthwhile to recall a few standard characterizations of interior points in the cone, summarized in \cite[Proposition~2.11]{GlW20}.
\begin{proposition} 
\label{prop:order-units}
	Let $(X,X^+)$ be an ordered Banach space and let $z \in X^+$. The following assertions are equivalent:
	\begin{enumerate}[label=(\roman*)]
		\item\label{prop:order-units:itm:interior-point} The vector $z$ is an element of the topological interior of $X^+$.
		
		\item\label{prop:order-units:itm:order-unit} The vector $z$ is an \emph{order unit}, i.e., for each $y \in X$ there exists $\varepsilon > 0$ such that $z \ge \varepsilon y$.
		
		\item\label{prop:order-units:itm:order-unit-variant} The cone $X^+$ is generating, and for each $y \in X^+$ there exists $\varepsilon > 0$ such that $z \ge \varepsilon y$.
		
		\item\label{prop:order-units:itm:principal-ideal} The so-called \emph{principal ideal} $\bigcup_{n \in \N} [-nz,nz]$ equals $X$.
		
		\item\label{prop:order-units:itm:ball} There exists $\varepsilon > 0$ such that $z \ge y$ for every $y \in X$ of norm $\norm{y} \le \varepsilon$.
	\end{enumerate}
\end{proposition}
In particular, if a cone in an ordered Banach space has a nonempty interior, then it is automatically generating. 
We point out that we do not a priori assume the cone $X^+$ to have nonempty interior in Proposition~\ref{prop:order-units}.

\begin{example}[Sequence spaces]
	\label{ex:sequence-spaces}
	\begin{enumerate}[label=(\roman*)]
		\item\label{ex:sequence-spaces:itm:ell_p} Let $X = \ell_p$ for $p \in [1,\infty]$, where for $p\in[1,+\infty]$
		\begin{align*}
			\ell_p := \{x = (x_n)_{n \in \N} \subset  \R^{\N}: \, \|x\|_{\ell_p} <\infty\},
		\end{align*}		
		and the norms in $\ell_p$ spaces are defined for $x \in\ell_p$ by $\|x\|_{\ell_p} := \Big(\sum_{n=1}^\infty |x_n|^p\Big)^{1/p}$ for $p<\infty$ and by $\|x\|_{\ell_\infty}:=\sup_{n=1}^\infty |x_n|$.
		
		We endow $\ell_p$ spaces with the cone
		\begin{align*}
			\ell_p^+ = \{x = (x_n)_{n \in \N} \in \ell_p: \, x_n \ge 0 \text{ for all } n \in \N\}.
		\end{align*}
		Then $(\ell_p,\ell_p^+)$ is an ordered Banach space, and in fact, even a Banach lattice (where the lattice operations can be performed entrywise); in particular, the cone $\ell_p^+$ is generating and normal. If $p = \infty$, then the cone has nonempty interior, whereas $\intt(\ell_p^+)=\emptyset$ for $p \in [1,\infty)$.
		
		\item\label{ex:sequence-spaces:itm:c-and-c_0} Let $X = c$ (the space of convergent real sequences) or $X = c_0$ (the space of real sequences that converge to $0$), endowed with the supremum norm and with the cone $X^+$ which is defined in the same way as in Example~\ref{ex:sequence-spaces}\ref{ex:sequence-spaces:itm:ell_p}. Then $(X,X^+)$ is an ordered Banach space and even a Banach lattice (with entrywise lattice operations). Thus, the cone $X^+$ is generating and normal. The interior of the cone is nonempty if $X = c$ and empty if $X = c_0$.
	\end{enumerate}
\end{example}

Several further examples of ordered Banach spaces are discussed, for instance, in \cite[Subsection~2.3]{ArN09} and \cite[Subsection~2.3]{GlW20}.

\section{Banach lattices}

\index{Banach lattice}

A particularly well-behaved class of ordered Banach spaces is the class of so-called \emph{Banach lattices} \cite{Sch74, Mey12}. An ordered Banach space $(X,X^+)$ is said to be \emph{lattice-ordered} if every two elements $x,y \in X$ have a smallest upper bound -- i.e., a so-called \emph{supremum} -- in $X$, which we then denote by $x \lor y$. This allows us to define the \emph{positive part}, the \emph{negative part} and the \emph{modulus} of a vector $x \in X$ by means of the formulae
\begin{align*}
	x^+  := x \lor 0, \qquad
	x^-  := (-x) \lor 0, \qquad
	\modulus{x}  := x \lor (-x).
\end{align*}
The cone in a lattice-ordered Banach space is always generating since we have $x = x^+ - x^-$ for each $x \in X$.

A \emph{Banach lattice} is a lattice-ordered Banach space that satisfies an additional compatibility assumption between norm and lattice operations, namely
\begin{align}
	\label{eq:compatibility-in-banach-lattices}
	\modulus{x} \le \modulus{y} \quad \Rightarrow \quad \norm{x} \le \norm{y}
\end{align}
for all $x,y \in X$. This implies, in particular, that the cone in a Banach lattice is always normal.

\section{The distance to the cone}

Recall that for a subset $S$ and a vector $x$ in a Banach space $X$, we denote by
\begin{align*}
	\dist(x,S) := \inf\left\{ \norm{x-y}: \, y \in S \right\}
\end{align*}
the distance from $x$ to $S$. If $(X,X^+)$ is an ordered Banach space, we will often need the distance of points to the positive cone $X^+$. Due to the specific properties of cones, the distance function $\dist(\cdot,X^+)$ quite nicely behaves; in particular, it is not difficult to see that we have for all $x,y \in X$ and all $\alpha \in [0,\infty)$
\begin{align*}
	& \dist(x+y,X^+) \le \dist(x,X^+) + \dist(y,X^+), \\
	\text{and} \qquad & \dist(\alpha x,X^+) = \alpha \dist(x,X^+).
\end{align*}

\section{Positive operators}

Let $X$ be an ordered vector space with a partial order $\leq$. We say that $x<y$ if $x\leq y$ and $x\neq y$.%

Having defined an order on a vector space, it is natural to consider maps on this set, which preserve the order in a certain sense.%

\begin{definition}\label{def:Increasing_maps}
Let $X,Y$ be ordered vector spaces with positive cones $X^+$ and $Y^+$, respectively. We say that the map $f:X\to Y$ is:
\index{function!monotone}
\index{function!positive}
\index{function!strictly positive}
\index{function!strongly positive}
\index{function!strictly monotone}
\begin{itemize}
	\item[(i)] \emph{monotone}, if $x\leq y$ implies $f(x)\leq f(y)$ for all $x,y\in X$.
	\item[(ii)] \emph{strictly monotone}, if $x < y$ implies $f(x)< f(y)$ for all $x,y\in X$.
	\item[(iii)] \emph{positive} if $f(X^+) \subset Y^+$.
	\item[(iv)] \emph{strictly positive} if $f(X^+\setminus\{0\}) \subset  Y^+\setminus\{0\}$.
	\item[(v)] If $(Y, Y^+)$ is an OBS and $Y^+$ has nonempty interior, then $f$ is called \emph{strongly positive} (in which case we write $f \gg 0$) if $f(X^+\setminus\{0\}) \subset  \intt Y^+$.
\end{itemize}
Note that if $f$ is linear, then positivity and monotonicity are equivalent properties. 
\end{definition}

We will be particularly interested in \emph{bounded} linear operators. Interestingly, if the cone $X^+$ is generating, then every positive linear operator $A: X \to Y$ is automatically bounded \cite[Theorem~2.32]{AlT07}.

\section{Duality of ordered Banach spaces} 
Let $(X, X^+)$ be an ordered Banach space. The subset
\begin{align*}
	(X')^+ := \{x' \in X': \, \langle x', x \rangle \ge 0 \text{ for all } x \in X^+\}
\end{align*}
of the dual space $X'$ is called the \emph{dual wedge} of $X^+$. The elements of $(X')^+$ are called the \emph{positive functionals} on $X$.
\index{dual wedge}

The dual wedge is also closed (even weak${}^*$-closed), convex and invariant with respect to multiplication by nonnegative scalars. The dual wedge $(X')^+$ is a cone -- i.e., its intersection with $-(X')^+$  is $\{0\}$ -- if and only if the cone $X^+$ is total in $X$.
In particular, if $(X,X^+)$ is an ordered Banach space with a total cone $X^+$, the dual space $(X', (X')^+)$ is an ordered Banach space too. The properties of $X^+$ and $(X')^+$ are closely related. 
The dual cone $(X')^+$ is generating if and only if $X^+$ is normal, see \cite[Theorem~4.5]{KLS89} or \cite[Corollary 2.43]{AlT07}.  Conversely, if $X$ is reflexive, the dual cone $(X')^+$ is normal if and only if $X^+$ is generating \cite[Theorem~4.6]{KLS89}, \cite[Corollary 2.43]{AlT07}.

\section{Complexifications} Since we are going to use the spectral theory of positive operators, a word on \emph{complexifications} is in order. The underlying scalar field of an ordered Banach space $(X,X^+)$ is real, but every real Banach space $X$ has a (in general, non-unique) \emph{complexification} $X_\C$, which is a Banach space over $\C$, and each bounded linear operator $A$ between real Banach spaces $X$ and $Y$ can be uniquely extended to a bounded $\C$-linear operator $A_\C: X_\C \to Y_\C$. Moreover, this extension always satisfies the norm estimate $\norm{A} \le \norm{A_\C} \le 2 \norm{A}$. Whenever we talk about spectral properties of $A$, we tacitly mean the corresponding spectral properties of the complex extension $A_\C$.
For an overview of complexifications of Banach spaces, we refer, for instance, to \cite{MST99} and \cite[Appendix~C]{Glu16}.

\section{Monotone control systems}

\begin{definition}
\index{control system!ordered}
We call a control system $\Sigma=(X,\Uc,\phi)$ (as defined in Definition~\ref{Steurungssystem}) \emph{ordered} if $X$ and $\Uc$ are normed vector spaces with a given order.
\end{definition}

Monotone control systems are an important in applications subclass of control system:
\begin{definition} 
\index{control system!monotone}
An ordered control system $\Sigma=(X,\Uc,\phi)$ is called \emph{monotone}, provided for all $t \geq 0$, all $x_1, x_2 \in X$ with $x_1 \leq x_2$ and all $u_1 \in\Uc$, $u_2 \in \Uc$ with $u_1 \leq u_2$ it holds that
$\phi(t,x_1,u_1)  \leq \phi(t,x_2,u_2)$.
\end{definition}

To treat situations when the monotonicity w.r.t. initial states is not available, the following definition is useful:
\begin{definition} 
\index{control system!monotone w.r.t.\  inputs}
An ordered control system $\Sigma=(X,\Uc,\phi)$ is called \emph{monotone w.r.t.\  inputs}, provided for all $t \geq 0$, all $x \in X$ and all $u_1,u_2 \in \Uc$ with $u_1 \leq u_2$ it holds that $ \phi(t,x,u_1)  \leq \phi(t,x,u_2)$.
\end{definition}

\section[Linear monotone discrete-time systems]{Linear monotone discrete-time systems and small-gain conditions}
\label{sec:Linear-discrete-time-systems}

\index{system!discrete-time}
Let $X$ be a normed vector space. Consider a map $A:X\to X$ and a corresponding induced \emph{discrete-time system}
\begin{eqnarray}
x(k+1) = A(x(k)), \quad k\in\N.
\label{eq:Gamma-discrete-time-system}
\end{eqnarray}
We denote the solution of this system at time $k\in\N$ subject to an initial condition $x\in X$ 
by $\phi(k,x) \in X$.

Let us define the stability notions for discrete-time systems.
\begin{definition}
\label{def:UGAS_discrete_time}
Discrete-time system \eqref{eq:Gamma-discrete-time-system} is called \emph{uniformly globally asymptotically stable (UGAS)}, if there is $\beta\in\KL$ so that
for each $k\geq 0$ and all $x\in X$ it holds that
\begin{eqnarray}
\|\phi(k,x)\|_X\leq \beta(k,\|x\|_X).
\label{eq:UGAS_discrete}
\end{eqnarray}
\end{definition}

It is also of virtue to define a stronger notion:
\begin{definition}
\label{def:UGES_discrete_time}
Discrete-time system \eqref{eq:Gamma-discrete-time-system} is called \emph{uniformly globally exponentially stable (UGES)},
if there are $a \in (0,1)$ and $M>0$ so that
for each $k\geq 0$ and all $x\in X$ it holds that
\begin{eqnarray}
\|\phi(k,x)\|_X\leq M a^{k}\|x\|_X.
\label{eq:UGES_discrete}
\end{eqnarray}
\end{definition}

In the rest of this section, we consider discrete-time systems \eqref{eq:Gamma-discrete-time-system} with $A\in L(X)$.
We start with the simple criterion of UGES for systems without any orderings.

\begin{proposition}
\label{prop:UGAS_discrete_systems_and_Small-gain_theorem_Operators_No_Order}
Let $X$ be a normed vector space and let $A \in L(X)$.
Then the following statements are equivalent:
\begin{enumerate}
\item[(i)]   $r(A) < 1$.
\item[(ii)]  $A^k \to 0$, for $k\to \infty$.
\item[(iii)]  \eqref{eq:Gamma-discrete-time-system} is UGES.
\item[(iv)]  \eqref{eq:Gamma-discrete-time-system} is UGAS.
\end{enumerate}
\end{proposition}

The positivity of the operator $A$ does not simplify the criteria for uniform exponential stability considered in Proposition~\ref{prop:UGAS_discrete_systems_and_Small-gain_theorem_Operators_No_Order}. On the other hand, though, positivity allows showing diverse characterizations of a very different nature, which are typically referred to as \emph{small-gain} type conditions.
The next result, shown in \cite[Theorem 3.3]{GlM21}, contains several such criteria. In the case if $X^+$ has a nonempty interior, or if $A$ is quasicompact, several additional criteria can be established, see \cite[Theorems 3.11, 3.14]{GlM21}.

\begin{theorem} 
	\label{thm:stability-for-pos-ops}
	Let $(X,X^+)$ be an ordered Banach space with generating and normal cone, and let $A \in \calL(X)$ be positive. Then the following assertions are equivalent:
	\begin{enumerate}[label=(\roman*)]
		\item\label{thm:stability-for-pos-ops:itm:stability} 
		\emph{Uniform exponential stability:} 
		The system~\eqref{eq:Gamma-discrete-time-system} satisfies the equivalent criteria of 
		Proposition~\ref{prop:UGAS_discrete_systems_and_Small-gain_theorem_Operators_No_Order}, 
		i.e., we have $r(A) < 1$.
		
		\item\label{thm:stability-for-pos-ops:itm:pos-resolvent} 
		\emph{Positivity of the resolvent at $1$:} 
		The operator $\id - A: X \to X$ is bijective and its inverse $(\id - A)^{-1}$ is positive.\footnote{Note that as $\id - A$ is invertible and bounded, then $(\id - A)^{-1}$ is closed, and since $A$ is surjective, $(\id - A)^{-1}$ is bounded by a closed graph theorem. Thus, $1 \in \rho(A)$ and $(\id - A)^{-1}$ is indeed a resolvent.}
		
		\item\label{thm:stability-for-pos-ops:itm:mbi} 
		\emph{Monotone bounded invertibility property:} 
		There exists a number $c \ge 0$ such that
		\begin{align*}
			(\id - A)x \le y \qquad \Rightarrow \qquad \norm{x} \le c \norm{y}
		\end{align*}
		holds for all $x,y \in X^+$.
				
		\item\label{thm:stability-for-pos-ops:itm:uniform-small-gain} 
		\emph{Uniform small-gain condition:} 
		There is a number $\eta > 0$ such that for each $x \in X^+$
		\begin{align*}
			\dist\big((A-\id)x,X^+\big) \ge \eta \norm{x}.
		\end{align*}

		\item\label{thm:stability-for-pos-ops:itm:perturbed-small-gain} 
		\emph{Robust small-gain condition:} 
		There exists a number $\varepsilon > 0$ such that
		\begin{align}
			\label{eq:Robust small-gain condition}
			(A+P)x \not\ge x
		\end{align}
		for every $0 \not= x \in X^+$ and for every positive operator $P \in \calL(X)$ of norm $\norm{P} \le \varepsilon$.
		
		\item\label{thm:stability-for-pos-ops:itm:perturbed-rank-1-small-gain} 
		\emph{Rank-$1$ robust small-gain condition:} 
		There exists a number $\varepsilon > 0$ such that 
		\begin{align*}
			(A+P)x \not\ge x
		\end{align*}
		for every $0 \not= x \in X^+$ and for every positive operator $P \in \calL(X)$ of rank $1$ and of norm $\norm{P} \le \varepsilon$.
	\end{enumerate}
\end{theorem}


\cleardoublepage
\chapter{Sobolev spaces and inequalities}
\label{chap:Function_Spaces_Inequalities}

The natural state spaces for the systems described by PDEs are $L^p$-spaces and Sobolev spaces.
In this chapter, we define these spaces and recall their main properties.
Special attention is devoted to the integral inequalities for the functions belonging to
$L^p$-spaces and Sobolev spaces. Mastering the use of these inequalities is crucial for the application of Lyapunov methods to PDE systems.

In this appendix, we denote by $G$ be an open bounded region in $\R^n$ (if not mentioned otherwise).
Under measurability, we always understand Lebesgue measurability.

\section{Spaces of integrable functions}
\label{sec:Function_spaces}

\begin{definition}
\label{def:Book_A.2.11}
\index{$L^p$}
Let $G \subset \R^n$ be an open nonempty set and let $p\geq1$ be a fixed real number.
Consider the set of (Lebesgue) measurable functions $x:G\to\R$ satisfying $\int\limits_G|x(z)|^p dz < \infty$ endowed with the map
\begin{gather*}
\|x\|:=\left(\int_G|x(z)|^p dz\right)^{1/p}.
\end{gather*}
This is a linear vector space with addition and scalar multiplication defined by:
\begin{gather*}
(x + y)(t) = x(t) + y(t),\qquad (\alpha x)(t) = \alpha x(t).
\end{gather*}
Clearly, $\|\lambda x\| = |\lambda|\|x\|$ for all scalar $\lambda$ and all $x$,and the triangle inequality follows by Minkowski's inequality (Proposition~\ref{prop:MinkowskiIneq}), which shows that $\|\cdot\|$ is a seminorm. However, $\|\cdot\|$ is not a norm on this space, as $\|x\| = 0$ only implies that $x(t) = 0$ almost everywhere.

Now consider the space of (equivalence) classes $[x]$
of functions that equal $x$ almost everywhere. These equivalence classes form a linear space and
\[
\|[x]\|_{L^p(G)} := \|x_1\|,\quad  \text{for any } x_1\in [x]
\]
defines a~norm; we call this normed vector space $L^p(G)$.

Usually, we write $x_1$ instead of $[x]$, where $x_1$ is any element of the equivalence class $[x]$.
\end{definition}

\begin{definition}\label{def:Book_A.2.12}
Let $G \subset \R^n$ be an open nonempty set. Consider all measurable functions $x: G \to \R$ satisfying the condition $\esssup_{z\in G}|x(z)|<\infty$.
Following Example~\ref{def:Book_A.2.11}, we construct equivalence classes $[x]$ containing functions that equal $x$ almost everywhere on $G$.
With the norm
\begin{gather*}
\|[x]\|_{\infty} := \esssup_{z\in G}|x_1(z)|, \qquad \text{for any } x_1\in [x],
\end{gather*}
this space is a~normed vector space, which we denote by $L^{\infty}(G)$.
Following the standard conventions, we usually write $x_1$ instead of $[x]$, where $x_1$ is any element of $[x]$.
\end{definition}

\begin{theorem}
\label{thm:Book_1.21}
If $p\in [1,\infty]$, and let $G \subset \R^n$ be an open nonempty set, then $L^p(G)$ is a~Banach space.

Furthermore, $L^2(G)$ is a Hilbert space with a scalar product given by
\[
\scalp{x}{y}_{L^2(G)}:=\int_G x(z)y(z)dz,\quad x,y\in L^2(G).
\]
\end{theorem}

\begin{proof}
See, e.g., \cite[p. 11, Theorem 1.21]{FHH11}, \cite[p. 12, Proposition 1.24]{FHH11}.
\end{proof}

For $k\in\N \cup \{+\infty\}$ we denote by $C_c^k(G)$  the space of $k$ times continuously differentiable functions 
$f:G \to \R$ with support compact in $G$.
The functions belonging to $C_c^{\infty}(G)$ are frequently called \emph{test functions}.

Test functions are dense in $L^p$-spaces for finite $p$.
\begin{proposition}
\label{prop:Properties_Lp_spaces}
Let $G \subset \R^n$ be an open nonempty set. 
$C_c^{\infty}(G)$ is dense in $L^p(G)$ (with respect to $L^p(G)$-norm) for $p\in[1,+\infty)$.
\end{proposition}

An important result in the theory of $L^p$-spaces is the scalar dominated convergence theorem, which we recall next:
\begin{theorem}[Scalar-valued dominated convergence theorem]
\label{thm:dominated convergence theorem}
Let $G \subset \R^n$ be an open nonempty set and let $(f_n)_{n\in\N}$ be a sequence of Lebesgue measurable scalar-valued functions on $G$, such that $f_n \to f$ a.e. in $G$, and there is $g \in L^1(G)$, such that 
\[
|f_n(z)| \leq g(z),\quad \text{ for a.e. } z \in G.
\] 
Then $f\in L^1(G)$ and 
\begin{eqnarray}
\lim_{n\to\infty}\int_G f_n(z)dz = \int_G f(z)dz.
\label{eq:Dominated-convergence-result}
\end{eqnarray}
\end{theorem}

\begin{proof}
See \cite[Theorem 5.6.1]{Ste05}.
\end{proof}

\section{Elementary inequalities}
\label{sec:Elementary_inequalities}

In this section, we collect several inequalities, which are instrumental for the stability analysis of infinite and finite-dimensional control systems and for verifying the important properties of function spaces.

\subsection{Inequalities in $\R^n$}
\label{sec:Elem_Rn_inequalities}

\begin{definition}
\label{def:convex-map} 
\index{function!convex}
A map $f:\R^n\to \R$ is called \emph{convex} provided that 
\[\
f(a x + (1-a)y) \leq af(x) + (1-a)f(y),
\]
for all $x,y \in\R^n$, and all $a \in[0,1]$.
\end{definition}

\begin{proposition}[Young's inequality]\label{prop:FHH11_1.11}
\index{inequality!Young's}
Let $p, q > 1$ be such that $\frac{1}{p}+\frac{1}{q}=1$.
Then for all $a, b \geq 0$ holds
\begin{eqnarray}
ab \leq \frac{a^p}{p} + \frac{b^q}{q}.
\label{eq:Youngs_Inequality_no_eps}
\end{eqnarray}
\end{proposition}

\begin{proof}
Since $x \mapsto e^x$ is a convex map, we have:
\[
ab = e^{\ln a  + \ln b}= e^{\frac{1}{p} \ln (a^p)  + \frac{1}{q}\ln (b^q)}\leq \frac{1}{p} a^p + \frac{1}{q} b^q.
\]
\end{proof}

\begin{proposition}[Young's inequality, general form]
\label{thm:Young}
For all $a,b \geq 0$ and all $\omega,p>0$ it holds
\begin{eqnarray}
\label{ineq:Young}
ab \leq \frac{\omega}{p}a^p + \frac{1}{\omega^{\frac{1}{p-1}}} \frac{p-1}{p}b^{\frac{p}{p-1}}.
\end{eqnarray}
\end{proposition}

\begin{proof}
Write $ab = \varepsilon^{\frac{1}{p}}a \cdot \frac{b}{\varepsilon^{\frac{1}{p}}}$, and apply for $q=\frac{p}{p-1}$ Young's inequality \eqref{eq:Youngs_Inequality_no_eps}.
\end{proof}

Often a special case of Proposition~\ref{thm:Young} (also called Cauchy's inequality with $\varepsilon$) is used:
\begin{proposition}[Young's inequality, special case]
\label{thm:Young_Simple}
For all $a,b \in\R$ and all $\varepsilon>0$ it holds
\begin{eqnarray}
\label{ineq:Young_Simple}
ab \leq \frac{\varepsilon}{2} a^2 + \frac{1}{2\varepsilon} b^2.
\end{eqnarray}
\end{proposition}

%

We close the section with (we adopt the formulation from \cite[Lemma 2.7, p.42]{Tes12}):
\begin{lemma}[Generalized Gronwall's inequality]
\label{lem:Gronwall} 
\index{inequality!generalized Gronwall's}
Let $\psi,\alpha,\beta \in C([0,T],\R)$ and let the inequality hold
	 \begin{equation}
			\psi(t) \le \alpha(t) + \int_{0}^{t}\beta(s)\psi(s)ds, \quad t \in [0,T],
	 \label{eq:2.34}
	 \end{equation}
	with $\alpha(t) \in \mathbb{R}$ and $\beta(t) \ge 0$, $t \in [0,T]$. Then:
\begin{enumerate}
	\item[(i)] $\psi(t) \le \alpha(t) + \int_{0}^{t} \alpha(s)\beta(s)\exp{\left(\int_{s}^{t}\beta(r)dr\right)}ds, \quad t \in [0,T].$
	\item[(ii)] If further $\alpha(s) \le \alpha(t)$ for $s \le t$, then it holds that 
	 \begin{equation*}
			\psi(t) \le \alpha(t)\exp{\left(\int_{0}^{t}\beta(s)ds\right)}, \quad t \in [0,T].
	 \end{equation*}
\end{enumerate}	
\end{lemma}

\subsection{Basic integral inequalities}
\label{sec:Elem_Integral_inequalities}

We collect here several basic inequalities in $L^p$ spaces.

\begin{proposition}[H\"older's inequality]
\label{prop:HoelderIneq}
\index{inequality!H\"older's (for integrals)}
Let $G$ be an open, bounded region in $\R^n$.
Assume $p,q \in [1, \infty]$ and $\frac{1}{p} + \frac{1}{q} = 1$. Then if $x \in L^p(G)$, $v \in L^q(G)$, then
\begin{eqnarray}
\label{HoelderIneq}
\int_G |x(z)v(z)| dz \leq \|x\|_{L^p(G)} \|v\|_{L^q(G)}.
\end{eqnarray}
\end{proposition}

\begin{proof}
Assume first that $\|x\|_{L^p(G)} = 1$ and $\|v\|_{L^q(G)} = 1$.
By Young's inequality \eqref{eq:Youngs_Inequality_no_eps} for $p, q$ as in the statement of the proposition, we have for a.e. $z \in G$
\begin{eqnarray}
|x(z)v(z)|  \leq \frac{1}{p}|x(z)|^p + \frac{1}{q} |v(z)|^q.
\label{eq:Hoelder_Inequality_HilfEstim}
\end{eqnarray}
The function $x v$ is measurable and \eqref{eq:Hoelder_Inequality_HilfEstim}
shows that its absolute value is upperbounded by an integrable function, thus $xv \in L^1(G)$ and
\begin{eqnarray*}
\int_G |x(z)v(z)| dz \leq \frac{1}{p}\int_G |x(z)|^p dz + \frac{1}{q}\int_G |v(z)|^q dz = \frac{1}{p} + \frac{1}{q} = 1
= \|x\|_{L^p(G)} \|v\|_{L^q(G)}.
\end{eqnarray*}
For general $x \in L^p(G)$, $v \in L^q(G)$, the result follows by homogeneity.
\end{proof}

Iterating H\"older's inequality, we obtain the following one:
\begin{proposition}[Generalized H\"{o}lder's inequality]
\label{prop:Generalized-Hoelder-integral}
\index{inequality!generalized H\"older's (for integrals)}
Let $G$ be an open, bounded region in $\R^n$.
Let $p_i \in[1,+\infty]$, $i=1,\ldots,m$ be such that $\frac{1}{p_1}+\ldots+\frac{1}{p_m}=1$ and let
$n\in\mathbb N$. 

Then for all $x_k \in L^{p_k}(G)$, $k=1,\ldots,m$, we have
\begin{gather}
\label{eq:Generalized-Hoelder-integral}
\int_G|x_1(z)\cdots x_m(z)|dz \leq \prod_{k=1}^m \|x_k\|_{L^{p_k}(G)}.
\end{gather}
\end{proposition}

\begin{proposition}[Cauchy-Bunyakovsky-Schwarz inequality in $L^2(G)$]
\label{prop:Cauchy-Schwarz}
\index{inequality!Cauchy-Schwarz}
For all $x \in L^2(G)$, $v \in L^2(G)$ it holds that
\begin{eqnarray}
\label{CauchySchwarzIneq}
\int_G |x(z)v(z)| dz \leq \|x\|_{L^2(G)} \|v\|_{L^2(G)}.
\end{eqnarray}
\end{proposition}

\begin{proof}
The result is a special case of H\"older's inequality \eqref{HoelderIneq} for $p=q=2$, as well as a special case of a general
Cauchy-Bunyakovsky-Schwarz' inequality in Hilbert spaces.
\end{proof}

The following result constitutes the triangle inequality in the $L^p(G)$-space.
\begin{proposition}[Minkowski inequality]
\label{prop:MinkowskiIneq}
\index{inequality!Minkowski (for integrals)}
Assume that $p,q \in [1, \infty]$ and $\frac{1}{p} + \frac{1}{q} = 1$. Then if $x,v \in L^p(G)$, then
\begin{eqnarray}
\label{eq:MinkowkiIneq_integral_version}
\|x+v\|_{L^p(G)}\leq \|x\|_{L^p(G)} + \|v\|_{L^p(G)}.
\end{eqnarray}
\end{proposition}

\begin{proof}
See \cite[p. 707]{Eva10}.
\end{proof}

An important integral inequality for convex functions is
\begin{proposition}[Jensen's inequality]
\label{thm:Jensen-general}
\index{inequality!Jensen's (for integrals)}
Let $G$ be open and bounded, and $f:\R^m \to \R$ be convex. Let also $x: G \to\R^m$ be summable. Then
\begin{eqnarray}
\label{ineq:Jensen-general}
 \frac{1}{|G|}\int_G f(x(z)) dz \geq f\Big(\frac{1}{|G|}\int_G x(z)dz\Big).
\end{eqnarray}
\end{proposition}

\begin{proof}
See \cite[p. 705]{Eva10}.
\end{proof}

In the scalar case, Jensen's inequality simplifies to
\begin{proposition}[Scalar Jensen's inequality]
\label{thm:Jensen}
For any convex $f:\R \to \R$ and any summable $x:[0,L] \to \R$, it holds that
\begin{eqnarray}
\label{ineq:Jensen}
 \int_0^L f(x(z)) dz \geq L f\Big( \frac{1}{L} \int_0^L x(z)dz\Big).
\end{eqnarray}
Similarly, for any concave $f:\R \to \R$ and any summable $x:[0,L] \to \R$, Jensen's inequality holds with the reverse inequality sign:
\begin{eqnarray}
\label{ineq:Jensen-concave}
 \int_0^L f(x(z)) dz \leq L f\Big( \frac{1}{L} \int_0^L x(z)dz\Big).
\end{eqnarray}
\end{proposition}

\begin{proof}
For convex $f$, the claim follows from the general Jensen's inequality. To see the second part, note that $f$ is concave iff $-f$ is convex. Applying Jensen's inequality to $-f$, we obtain the result for concave $f$.
\end{proof}


\section{Weak derivatives}
\label{sec:Weak derivatives}

Take any $x \in C^1(G)$, and any $\phi \in C_c^{\infty}(G)$, that is $\phi$ is infinitely times continuously differentiable with a compact in $G$ support. As $\phi$ is compactly supported in $G$, for each $i=1,\ldots,n$ the integration by parts formula shows that 
\begin{eqnarray}
\int_G x(z) \phi_{z_i}(z) dz = - \int_G x_{z_i}(z) \phi(z) dz.
\label{eq:weak-derivative-motivation}
\end{eqnarray}

Take any $k\in\N$, any $x \in C^k(G)$ and any multiindex $\alpha:=(\alpha_1,\ldots,\alpha_n)$ of order $|\alpha|:=\alpha_1+\ldots+\alpha_n = k$. We use the following notation for derivatives:
\[
D^\alpha \phi:= \frac{\partial^{\alpha_1}}{\partial z_1^{\alpha_1}}\ldots\frac{\partial^{\alpha_n}}{\partial z_n^{\alpha_n}}\phi.
\]

Then similarly as above it holds for all $\phi \in C_c^{\infty}(G)$ that 
\begin{eqnarray}
\int_G x(z) D^{\alpha} \phi(z) dz = (-1)^{|\alpha|} \int_G D^{\alpha}x(z) \phi(z) dz.
\label{eq:weak-derivative-motivation-multi-dim}
\end{eqnarray}

The expression in the left-hand side of \eqref{eq:weak-derivative-motivation} and 
\eqref{eq:weak-derivative-motivation-multi-dim} makes sense already for $x \in L^{1}_{\loc}(G)$. This motivates us to relax the differentiation concept as follows:

\begin{definition}
\label{def:definition-weak-derivative} 
\index{derivative!weak}
Suppose $x,v \in L^{1}_{\loc}(G)$, and $\alpha$ is a multi-index. We say that $v$ is the \emph{$\alpha^{th}$-weak partial derivative} of $x$, which we express by 
\begin{eqnarray*}
D^\alpha x = v,
\end{eqnarray*}
if for all test functions $\phi \in C_c^{\infty}(G)$ it holds that
\begin{eqnarray}
\int_G x(z) D^{\alpha} \phi(z) dz = (-1)^{|\alpha|} \int_G v(z) \phi(z) dz.
\label{eq:weak-derivative-multi-dim}
\end{eqnarray}
\end{definition} 

The existence of a weak derivative is a global property. We do not speak about \q{weak derivative of a function at a given point of $G$}.

The next lemma shows that the weak derivative is a well-defined concept if we identify the functions that differ on a set of measure zero. 
\begin{lemma}
\label{lem:Weak-derivatives-are-well-defined} 
Suppose $x\in L^{1}_{\loc}(G)$, and $\alpha$ is a multi-index.
If there exists a $\alpha^{th}$-weak partial derivative of $x$, then it is uniquely defined up to a set of measure 0.
\end{lemma}

\begin{proof}
Assume that there are $v_1,v_2$ such that 
\begin{eqnarray*}
\int_G x(z) D^{\alpha} \phi(z) dz = (-1)^{|\alpha|} \int_G v_1(z) \phi(z) dz = (-1)^{|\alpha|} \int_G v_2(z) \phi(z) dz,
\end{eqnarray*}
for all $\phi \in C_c^{\infty}(G)$. 	
Then for all $\phi \in C_c^{\infty}(G)$ we have also
\begin{eqnarray*}
\int_G (v_1(z) - v_2(z)) \phi(z) dz = 0.
\end{eqnarray*}
This implies that $v_1-v_2=0$ a.e.
\end{proof}

Consider an example of a function, which is not differentiable in the classical sense at a point in the interval, yet it does possess a weak derivative on the whole interval.
\begin{example}
\label{examp:weak-derivative-exercise}
Let $G=(0,2)\subset \R$, and 
\begin{eqnarray}
x(z) :=
\begin{cases}
z, \quad  &z\in (0,1),\\
1, \quad  &z \in [1,2).
\end{cases}
\end{eqnarray}
Let us show that the weak derivative of $x$ equals
\begin{eqnarray*}
v(z) :=
\begin{cases}
1, \quad  &z\in (0,1),\\
0, \quad  &z \in [1,2).
\end{cases}
\end{eqnarray*}
Indeed, we have for any $\phi \in C_c^{\infty}(G)$ that
\begin{eqnarray*}
\int_0^2 x(z) \phi'(z) dz =
\int_0^1 z \phi'(z) dz + \int_1^2 \phi'(z) dz
&=&  z\phi(z)|_0^1 - \int_0^1 \phi(z) dz - \phi(1)\\
&=&  - \int_0^2 v(z) \phi(z) dz.
\end{eqnarray*}
\end{example}

At the same time, adding a discontinuity to the previous example prevents a function from being weakly differentiable at the interval.
\begin{example}
\label{examp:weak-derivative-exercise-2}
Let $G=(0,2)\subset \R$, and 
\begin{eqnarray}
x(z) :=
\begin{cases}
z, \quad  &z\in (0,1],\\
2, \quad  &z \in (1,2).
\end{cases}
\end{eqnarray}
Let us show that the weak derivative of $x$ does not exist.

Assume that there is $v = D^1 x \in L^{1}_{\loc}(G)$. Then for any $\phi \in C_c^{\infty}(G)$ it holds that
\begin{eqnarray*}
-\int_0^2 v(z)\phi(z) dz 
&=& \int_0^2 x(z) \phi'(z) dz \\
&=& \int_0^1 z \phi'(z) dz + 2\int_1^2 \phi'(z) dz\\
&=&  z\phi(z)|_0^1 - \int_0^1 \phi(z) dz - 2\phi(1)
= - \int_0^1 \phi(z) dz - \phi(1).
\end{eqnarray*}
Take a sequence $(\phi_m)$ of test functions, such that 
$\phi_m(1)=1$ for all $m$, $\phi_m(z) \in[0,1]$ for all $z\in (0,2)$ and all $m\in\N$.
Assume also that $\supp(\phi_m)\subset [1-\frac{1}{m},1+\frac{1}{m}]$ for all $m\in\N$.

As $v\in L^1_{\loc}(G)$, we have for all $m\in\N$ that
\begin{eqnarray*}
1 = \phi_m(1) 
&=& \int_0^2 v(z)\phi_m(z) dz  - \int_0^1 \phi_m(z) dz\\
&\le& \int_0^2 |v(z)\phi_m(z)| dz  - \int_0^1 \phi_m(z) dz.
\end{eqnarray*}
Now for a. e. $z\in (0,2)$ it holds that $|v(z)\phi_m(z)| \to 0$ as $m\to\infty$, and 
\[
|v(z)\phi_m(z)| \leq |v(z)|,
\]
with $v \in L^1_{\loc}(G)$.

Dominated convergence theorem shows that
\begin{eqnarray*}
1 = \phi_m(1) 
&\le& \lim_{m\to\infty} \int_0^2 |v(z)\phi_m(z)| dz  - \lim_{m\to\infty}  \int_0^1 \phi_m(z) dz\\
&=& \int_0^2 \lim_{m\to\infty} |v(z)\phi_m(z)| dz \\
&=& 0,
\end{eqnarray*}
a contradiction. 
\end{example}

\section{H\"older spaces}
\label{sec:Hoelder spaces}

In this section, let $G\subset \mathbb{R}^n$ be open. We define the space $C^k(\clo{G})$ as the space of all functions $x \in C^k(G)$ such that $D^\alpha x$ is uniformly continuous on bounded subsets of $G$ for all multiindices $\alpha:|\alpha|\leq k$.

In particular, if $x\in C^k(\clo{G})$, then for each multiindex $\alpha$ such that $|\alpha|\leq k$, $D^\alpha x$ can be continuously extended to $\clo{G}$.

\begin{definition}
\label{def:Hoelder-continuity} 
Let $\gamma \in (0,1]$.
A~map $x : G \to \mathbb{R}$ is called \emph{H\"older continuous with exponent $\gamma$}, if there is $C>0$ such that 
\begin{eqnarray}
\label{eq:Eva10_2-Hoelder}
|x(z) - x(y)|\leq  C|z - y|^{\gamma},\qquad z,y\in G.
\end{eqnarray}
\end{definition}

\begin{definition}
Let $\gamma \in (0,1]$.

\item[(i)] For a bounded and continuous $x : G \to \mathbb{R}$, we denote its standard sup-norm by
\begin{eqnarray*}
\|x\|_{C(\clo{G})} :=\sup_{z\in G}|x(z)|.
\end{eqnarray*}
\item[(ii)] The \emph{$\gamma^{\text{th}}$-H\"older seminorm} of $x : G \to \mathbb{R}$ is defined by
\begin{eqnarray*}
[x]_{C^{0,\gamma}(\clo{G})} :=\sup_{\stackrel{z,y\in G}{z\ne y}}\left\{ \frac{|x(z)-x(y)|}{|z-y|^{\gamma}}\right\},
\end{eqnarray*}
and the \emph{$\gamma^{\text{th}}$-H\"older norm} is
\begin{eqnarray*}
\|x\|_{C^{0,\gamma}(\clo{G})} := \|x\|_{C(\clo{G})} + [x]_{C^{0,\gamma}(\clo{G})}.
\end{eqnarray*}
\end{definition}

\begin{definition}
\index{space!H\"older}
The \emph{H\"older space} $C^{k,\gamma}(\clo{G})$ consists of all maps $x\in C^k(\clo{G})$ with a finite norm\begin{eqnarray}
\label{eq:Eva10_HoelderNorm}
\|x\|_{C^{k,\gamma}(\clo{G})}
:= \sum_{|\alpha|\leq  k} \|D^{\alpha}x\|_{C(\clo{G})}
 + \sum_{|\alpha| = k} \big[D^{\alpha}x\big]_{C^{0,\gamma}(\clo{G})}.
\end{eqnarray}
\end{definition}

In other words, the space $C^{k,\gamma}(\clo{G})$ is comprised of $k$-times continuously
differentiable functions $x$ whose $k^{\text{th}}$-partial derivatives are bounded and
H\"older continuous with exponent $\gamma$.

Note that H\"older continuous functions are uniformly continuous, and thus they can be uniquely extended to $\clo{G}$, with the same H\"older (semi)norms. 
This is why we define (following \cite[Section 5.1, p. 254]{Eva10}) the seminorms and norms in H\"older spaces on $\clo{G}$. 

\begin{theorem}
\label{thm:Hoelder spaces are Banach spaces}
The space $C^{k,\gamma}(\clo{G})$ is a~Banach space.
\end{theorem}

\begin{proof}
See \cite[Theorem 1, p. 255]{Eva10}.
\end{proof}

\section{Sobolev spaces}
\label{sec:Sobolev_spaces}

Let $G,V \subset \R^n$ be open sets. We say that \emph{$V$ is compactly contained in $G$}, which we denote by $V \subset\subset G$, if $V \subset \clo{V} \subset G$, with compact $\clo{V}$.

The spaces of functions possessing weak derivatives of various orders, which lie in $L^p$ spaces will be very important for us in the sequel.  
\begin{definition}
\label{def:Sobolev space} 
\index{space!Sobolev}
Let $p \in[1,\infty]$ and $k\in\N$. 
The \emph{Sobolev space} $W^{k,p}(G)$ consists of functions $x \in L^p(G)$, such that for each multiindex $\alpha$ of order  $\leq k$, the weak derivative $D^\alpha x$ exists and belongs to $L^p(G)$.

For $p<\infty$ the norm in $W^{k,p}(G)$ is defined by
\begin{equation}
\label{SobNorm-glossary}
\|x\|_{W^{k,p}(G	)} := \Big( \sum_{|\alpha| \leq k}\int_G{\left| D^\alpha x (z)\right|^p dz} \Big)^{\frac{1}{p}}.
\end{equation}
The norm in $W^{k,\infty}(G)$ we define by
\begin{equation}
\label{eq:Wk-infty-norm}
\|x\|_{W^{k,\infty}(G)} := \sum_{|\alpha| \leq k}\left\| D^\alpha x\right\|_{L^\infty(G)}.
\end{equation}

By $W^{k,p}_0(G)$ we denote the closure of $C_c^\infty(G)$ in the norm of $W^{k,p}(G)$.

$W^{k,p}_{\loc}(G)$ is the space of all $x:G\to\R$, such that for any $V\subset\subset G$ the map $x$ belongs to 
$W^{k,p}(V)$.

Furthermore, if $p=2$ we further define $H^k(G):=W^{k,2}(G)$,  $H^k_0(G):=W^{k,2}_0(G)$.
\end{definition}

Several basic properties of weak derivatives are collected in the next proposition:
\begin{theorem}[Properties of weak derivatives]
\label{thm:Book_p}
Assume $x,v\in W^{k,p}(G)$, $|\alpha| \leq  k$. Then
\begin{itemize}
\item[(i)] $D^{\alpha}x\in W^{k-|\alpha|,p}(G)$
and $D^{\beta}(D^{\alpha}x) = D^{\alpha}(D^{\beta}x) = D^{\alpha+\beta}x$
for all
multiindices $\alpha$, $\beta$ with $|\alpha| + |\beta| \leq  k$.
\item[(ii)] For each $\lambda,\mu\in \mathbb{R}$,
$\lambda x + \mu v\in W^{k,p}(G)$
and $D^{\alpha}(\lambda x + \mu v)
= \lambda D^{\alpha}x +\mu D^{\alpha}v$,
$|\alpha| \leq k$.
\item[(iii)] If $V$ is an open subset of $G$, then $x\in W^{k,p}(V)$.
\item[(iv)] Leibniz's formula: If $\zeta \in C_c^{\infty}(G)$,
then $\zeta x\in W^{k,p}(G)$ and
\begin{eqnarray}\label{eq:Book_7}
D^{\alpha}(\zeta x) = \sum_{\beta \leq \alpha}
\begin{pmatrix}
\alpha \\ \beta
\end{pmatrix}
D^{\beta}\zeta D^{\alpha-\beta}x,
\end{eqnarray}
where $\begin{pmatrix}
\alpha \\ \beta
\end{pmatrix}=\frac{\alpha!}{\beta!(\alpha-\beta)!}$, and the factorial is defined for multiindices as $\alpha!:=(|\alpha|)!$.
\end{itemize}
\end{theorem}

\begin{proof}
See \cite[Theorem 1, p. 261]{Eva10}.
\end{proof}

The Sobolev spaces have a very nice mathematical structure.
\begin{theorem}
\label{thm:Sobolev spaces are Banach spaces} 
For each $k\in\N$ and each $p\in[1,+\infty]$ the Sobolev space $W^{k,p}(G)$ is a Banach space.
\end{theorem}

\begin{proof}
See \cite[Theorem 2, p. 262]{Eva10}.
\end{proof}

The next result gives a useful 
characterization of the Sobolev spaces on 1-dimensional intervals.
\begin{proposition}
\label{prop:Sobolev-spaces-in-1-d} 
Let $n=1$ and $p\in[1,+\infty)$. The Sobolev space $W^{1,p}(0,1)$ coincides (if we identify the functions which are equal up to the set with the measure zero) with the space of absolutely continuous functions on $(0,1)$, whose classical derivative exists a.e. on $(0,1)$ and belongs to $L^p(0,1)$.
\end{proposition}

\begin{proof}
See \cite[Exercise 4, p. 306]{Eva10}.
\end{proof}

Elements of Sobolev spaces can be approximated by smooth functions.

For a sequence $(x_k)$, and $x \in W^{k,p}(G)$, we will say that
\begin{itemize}
	\item  $x_k \to x$ in $W^{k,p}(G)$ if $\|x_k -x\|_{W^{k,p}(G)}\to 0$ as $k\to\infty$. 
	\item  $x_k \to x$ in $W^{k,p}_{\loc}(G)$ if for any $V\subset\subset G$ it holds that $\|x_k -x\|_{W^{k,p}(V)}\to 0$ as $k\to\infty$. 
\end{itemize}

One says that $\partial G \in C^k$, for $k \in \N$, if for every $z \in \partial G$ there is a ball $B_\delta(z)$ such that there is a $k$ times continuously differentiable function $\phi$ and an index $i$, such that for all $y \in \partial G \cap B_\delta(z)$ it holds that $y_i = \phi(y_1,\ldots,y_{i-1},y_{i+1},\ldots,y_n)$.

We have the following result:
\begin{theorem}[Approximation by smooth functions]
\label{thm:Book_local}
Assume $x\in W^{k,p}(G)$ for some $1 \leq p <\infty$. Then
\begin{itemize}
\item[(i)] There is $(x_m) \subset C^\infty_c(G)$, such that $x_k \to x$ in $W^{k,p}_{\loc}(G)$.
\item[(ii)] If $G$ is bounded, then there is $(x_m) \subset C^\infty(G)\cap W^{k,p}(G)$, such that $x_k \to x$ in $W^{k,p}(G)$.
\item[(ii)] If $G$ is bounded and $\partial G \in C^1$, then there is $(x_m) \subset C^\infty(\clo{G})$, such that $x_k \to x$ in $W^{k,p}(G)$.
\end{itemize}
\end{theorem}

\begin{proof}
See \cite[Theorems 1, 2, 3 in Section 5.3]{Eva10}.
\end{proof}

Next, we discuss how to define the boundary values of a  function $x\in W^{1,p}(G)$, assuming that $\partial G$ is $C^1$.
This is not straightforward, as $\partial G$ has $n$-dimensional Lebesgue measure zero, and $x$ is defined only up to the set of measure zero.
The notion of a~\emph{trace operator} resolves this problem.
\begin{theorem}[Trace Theorem]
\label{thm:Book_Trace}
\index{operator!trace}
Assume $G$ is bounded and $\partial G$ is $C^1$. Let further $1\leq p <\infty$.
Then there exists a bounded linear operator
\begin{eqnarray*}
T \colon W^{1,p}(G) \to L^p(\partial G)
\end{eqnarray*}
such that
\begin{itemize}
\item[(i)] $\left. Tx = x \right|_{\partial G}$ if
$x\in W^{1,p}(G) \cap C(\clo{G})$.

\item[(ii)] There is $C>0$ depending only on $p$ and $G$, such that
\begin{eqnarray*}
\|Tx\|_{L^p(\partial G)} \leq C\|x\|_{W^{1,p}(G)},\quad x\in W^{1,p}(G).
\end{eqnarray*}

\end{itemize}
\end{theorem}

\begin{proof}
See \cite[Theorem 1 in Section 5.5]{Eva10}.
\end{proof}

\begin{definition}
We call $Tx$ the \emph{trace} of $x$ on $\partial G$.
\end{definition}

Next, we characterize the space $W_0^{1,p}(G)$ in terms of the trace operator.
\begin{theorem}[Characterization of $W_0^{1,p}(G)$]
Assume $G$ is bounded and $\partial G$ is $C^1$.
Suppose furthermore that $x\in W^{1,p}(G)$.
Then
\begin{eqnarray}\label{eq:Book_Trace4}
x\in W_0^{1,p}(G)
\quad\text{if and only if}\quad
Tx =0
\quad\text{on}\quad
\partial G.
\end{eqnarray}
\end{theorem}

\begin{proof}
See \cite[Theorem 2 in Section 5.5]{Eva10}.
\end{proof}

\begin{definition}
\index{Sobolev conjugate}
If $1\leq p < n$, the \emph{Sobolev conjugate} of $p$ is
\begin{eqnarray}
\label{eq:Eva10_Ch-5_Sec-6_8}
p^{\ast}:= \frac{np}{n-p}.
\end{eqnarray}
Note that
\begin{eqnarray}
\label{eq:Eva10_Ch-5_Sec-6_8a}
\frac{1}{p^{\ast}}=\frac{1}{p}-\frac{1}{n},
\qquad \text{and}
\qquad 
p^{\ast}>p.
\end{eqnarray}
\end{definition}

We start by showing that if $x\in W^{1,p}(G)$, $1\leq p <n$, then $x\in L^{p^{\ast}}(G)$. As the Sobolev conjugate $p^\ast$ is larger than $p$, this means that the existence of weak derivatives for a function in $L^p(G)$ space already implies that it has better integrability properties.
\begin{theorem}[Sobolev embedding theorem for $W^{1,p}$, $1\leq p <n$]
\label{thm:Sobolev-embeddings-W1p-p-less-n}
Let $G$ be a~bounded, open subset of $\mathbb{R}^n$,
and suppose $\partial G$ is $C^1$.
Assume $1\leq p < n$,
and $x\in W^{1,p}(G)$.
Then $x\in L^{p^{\ast}}(G)$, with the estimate
\begin{eqnarray}
\label{eq:Eva10_Ch-5_Sec-6_15}
\|x\|_{L^{p^{\ast}}(G)} \leq C\|x\|_{W^{1,p}(G)},
\end{eqnarray}
the constant $C$ depending only on $p$, $n$, and $G$.
\end{theorem}

\begin{proof}
See \cite[Theorem 2 in Section 5.6]{Eva10}.
\end{proof}

For functions in $W_0^{1,p}(G)$ a stronger result can be achieved:
\begin{theorem}[Estimates for $W_0^{1,p}$, $1\leq p < n$]
\label{thm:Poincare-inequality-multidim}
Assume $G$ is a~bounded,
open subset of $\mathbb{R}^n$.
Suppose $x \in W_0^{1,p}(G)$ for some $1\leq p < n$.
Then we have the estimate
\begin{eqnarray*}
\|x\|_{L^q(G)} \leq C\|Dx\|_{L^p(G)}
\end{eqnarray*}
for each $q \in [1, p^{\ast}]$,
the constant $C$ depending only on $p$, $q$, $n$ and $G$.

In particular, we obtain for all $1\leq p \leq \infty$ the so-called \emph{Friedrichs' inequality} (sometimes called \emph{Poincar\'e's inequality}):
\begin{eqnarray}
\label{eq:Poincare-multidimensional}
\|x\|_{L^p(G)} \leq C \|Dx\|_{L^p(G)},\quad x \in W_0^{1,p}(G).
\end{eqnarray}
\end{theorem}

\begin{proof}
See \cite[Theorem 3 in Section 5.6]{Eva10}.
\end{proof}

\begin{remark}
\label{rem:Equivalent-norms-on-Sobolev-spaces-with-0} 
In view of Friedrichs' inequality, on $W_0^{1,p}(G)$ the norm $\|Dx\|_{L^p(G)}$
is equivalent to the Sobolev norm $x \mapsto \|x\|_{W^{1,p}(G)}$, if $G$ is bounded.
\end{remark}

Having analyzed Sobolev spaces $W^{1,p}(G)$ with $p<n$, we proceed to the analysis of such spaces with $n<p\le\infty$.

\begin{definition}
We say $x^{\ast}$ is a~\emph{version} of a~given function $x$ provided
\begin{eqnarray*}
x = x^{\ast} \text{ a.e.}
\end{eqnarray*}
\end{definition}

We will show that if $x\in W^{1,p}(G)$, with $p>n$, then $x$ is H\"older continuous if we redefine it on a~set of measure zero.
\begin{theorem}[Sobolev embedding theorem for $W^{1,p}$, $n < p \leq \infty$]
\label{thm5}
Let $G$ be a~bounded, open subset of $\mathbb{R}^n$,
and suppose $\partial G$ is $C^1$.
Assume $n < p \leq \infty$ and $x\in W^{1,p}(G)$.
Then $x$ has a~version $x^{\ast} \in C^{0,\gamma}(\clo{G})$,
for $\gamma = 1-\frac{n}{p}$, with the estimate
\begin{eqnarray*}
\|x^{\ast}\|_{C^{0,\gamma}(\clo{G})} \leq C\|x\|_{W^{1,p}(G)}.
\end{eqnarray*}
The constant $C$ depends only on $p$, $n$ and $G$.
\end{theorem}

\begin{proof}
See \cite[Theorem 5 in Section 5.6]{Eva10}.
\end{proof}

In view of Theorem~\ref{thm5}, we will henceforth always identify a~function
$x \in W^{1,p}(G)$ ($p > n$) with its continuous version.

\begin{remark}
\label{rem:The borderline case} 
Finally, assume that $p = n$. 
Thanks to Theorem~\ref{thm:Sobolev-embeddings-W1p-p-less-n} and as
$p^{\ast} = \frac{np}{n-p} \to +\infty$ whenever $p \to n$, one might
conjecture that $x\in W^{1,p}(G)$ implies that $x\in L^{\infty}(G)$.
However, this does not hold in general for $n > 1$.
For example, if $G = B^0(0,1)$,
the function $x: z \mapsto \log\log \left(1+\frac{1}{|z|}\right)$
belongs to
$W^{1,n}(G)$ but it is not in $L^{\infty} (G)$, see \cite[Exercise 14, p. 307]{Eva10}.
\end{remark}


\section{Integral inequalities for 1-dimensional domains}
\label{sec:Elem_Integral_inequalities-1d}

This section states several integral inequalities on $(0,L)$. 
In what follows, we denote by $\scalp{\cdot}{\cdot}$ the classical scalar product in $L^2(0,L)$. 
First, we show the well-known  Friedrichs' inequality (see, e.g., \cite[Appendix A]{Mur03}).
\begin{proposition}[Friedrichs' inequality]
\label{theorem:Friedrichs}
\index{inequality!Friedrichs'}
For all $x \in H^1_0(0,L)$ it holds that
\begin{align}
\label{ineq:Friedrichs}
\frac{L^2}{\pi^2} \int_0^L{\left( \frac{d x(z)}{d z} \right)^2 dz} \geq \int_0^L {x^2(z)}dz.
\end{align}
\end{proposition} 

\begin{proof}
Define $\phi^{<i>}(z):=\sqrt{\frac{2}{L}}\sin\frac{i\pi z}{L}$, $i\in\N$, $z\in (0,L)$. It is well-known, that $(\phi^{<i>})_{i\in\N}$ is an orthonormal family of functions in $L^2(0,L)$.

Pick any $x \in C^\infty_c(0,L)$. Then $x$ be extended by symmetry to a periodic odd function on $\R$, and thus it can be represented as a Fourier series of the form
\[
x(z) = \sum_{i=1}^\infty c_i \phi^{<i>}(z),\quad z \in (0,L),
\]
where $c_i := \scalp{x}{\phi^{<i>}}$, $i\in\N$. 

We have that
\begin{eqnarray*}
\int_0^L x^2_{z}(z)dz &=& \scalp{\sum_{i=1}^\infty c_i \phi^{<i>}_{z}}{\sum_{j=1}^\infty c_j \phi^{<j>}_{z}}\\
&=& \sum_{i=1}^\infty \sum_{j=1}^\infty c_i c_j \scalp{ \phi^{<i>}_{z}}{\phi^{<j>}_{z}}\\
&=& \sum_{i=1}^\infty \sum_{j=1}^\infty c_i c_j \int_0^1 \phi^{<i>}_{z}(z) \phi^{<j>}_{z}(z) dz.
\end{eqnarray*}
Using the integration by parts, we see that
\begin{eqnarray*}
\int_0^L x^2_{z}(z)dz &=& -\sum_{i=1}^\infty \sum_{j=1}^\infty c_i c_j \int_0^1 \phi^{<i>}(z) \phi^{<j>}_{zz}(z) dz\\
&=& \sum_{i=1}^\infty \sum_{j=1}^\infty c_i c_j \Big(\frac{j\pi}{L}\Big)^2 \int_0^1 \phi^{<i>}(z) \phi^{<j>}(z) dz.
\end{eqnarray*}
As $\{\phi^{<j>}:\ j\in\N\}$ is an orthonormal family in $L^2(0,1)$, we proceed to
\begin{eqnarray*}
\int_0^L x^2_{z}(z)dz&=& \sum_{i=1}^\infty c^2_i \Big(\frac{i\pi}{L}\Big)^2\\
&\geq& \Big(\frac{\pi}{L}\Big)^2 \sum_{i=1}^\infty c^2_i\\
&=& \Big(\frac{\pi}{L}\Big)^2 \int_0^L x^2(z)dz.
\end{eqnarray*}
As $ H^1_0(0,L)$ is a closure of $C^\infty_c(0,L)$ in $H^1(0,L)$-norm, holds the last inequality also for all $x\in H^1_0(0,L)$.
\end{proof}
The constant $\frac{L^2}{\pi^2}$ in \eqref{ineq:Friedrichs} is optimal. 
The essence of the above proof of the Friedrichs' inequality is the decomposition of the functions in an orthogonal sequence in the space $L^2(0,L)$. One can see that the functions $\phi^{<i>}$ are eigenfunctions of the Dirichlet Laplacian operator.
Using this idea, one can extend the above Friedrichs' inequality to multidimensional bounded regions. 
In this case, one has to use the series expansion of the smooth functions with respect to the eigenfunctions of the multidimensional Dirichlet Laplace operator.

Friedrichs' inequality is essential for studying parabolic problems with homogeneous Dirichlet boundary conditions.
If Dirichlet condition is given only on one side of the boundary, Friedrichs' inequality cannot be applied. Instead, Poincare's inequality plays a prominent role.

\begin{proposition}[Poincare's inequality]
\label{theorem:Wirtinger}
\index{inequality!Poincare's}
For every $x \in H^1(0,L)$ with either $x(0)=0$ or $x(L)=0$, it holds that
\begin{align}
\label{Wirtinger_Variation_Ineq}
\frac{4 L^2}{\pi^2} \int_0^L{\left( \frac{dx(z)}{d z} \right)^2 dz} \geq \int_0^L {x^2(z)}dz
\end{align}
\end{proposition}

\begin{proof}
See \cite[Section 7.6, Nr 256]{HLP52}.
\end{proof}

A notable inequality, that relates $L^\infty$-norm to $H^1$-norm, is:
\begin{proposition}[Agmon's inequality]
\label{thm:Agmon}
\index{inequality!Agmon's}
For all $x \in H^1(0,L)$ it holds that
\begin{align}
\label{ineq:Agmon}
\|x\|^2_{L^{\infty}(0,L)} \leq |x(0)|^2 + 2\|x\|_{L^{2}(0,L)}\Big\|\frac{dx}{dz}\Big\|_{L^{2}(0,L)}.
\end{align}
\end{proposition}

\begin{proof}
See \cite[Lemma 2.4., p. 20]{KrS08}.
\end{proof}

We state also a useful inequality from \cite[Lemma 1]{ZhZ18}:
\begin{proposition}
\label{lem:ZhZ17_Lem1}
Let $x\in C^1([a,b],\R)$. Then for each $c\in[a,b]$ it holds that
\begin{eqnarray}
x^2(c) \leq \frac{2}{b-a}\|x\|^2_{L^2(0,1)} + (b-a)\|x_z\|^2_{L^2(0,1)}
\label{eq:ZhZ17_Lem1}
\end{eqnarray}
\end{proposition}

\section{An inequality for the Kuramoto-Sivashinskiy equation}
\label{sec:An inequality for the Kuramoto-Sivashinskiy equation}

For a given $\lambda\in\R$, consider an eigenvalue problem
\begin{subequations}
\begin{eqnarray}
x_{zzzz} + \lambda x_{zz} = \sigma x, \quad z\in(0,1),\\
x(0) = x(1) = x_z(0) = x_z(1) = 0.
\end{eqnarray}
\label{eq:KS_Eigenvalue_Problem}
\end{subequations}
One can show, see \cite{LiK01}, that there are countable sequences of eigenfunctions $\{\phi_i:i\in\N\}$ and corresponding eigenvalues $(\sigma_i)_{i\in\N}$, such that $(\phi_i,\sigma_i)$ solve the above eigenvalue problem, and $\{\phi_i:i\in\N\}$ form an orthonormal basis of $L^2(0,1)$. It holds also that $\lim_{n\to\infty}\sigma_n = \infty$.

 We define 
\begin{eqnarray}
\sigma(\lambda):=\min_{n\in\N}\sigma_n.
\label{eq:min-sigman-KS-equation}
\end{eqnarray}

We have the following result, see \cite[Lemma 2.1]{LiK01}:
\begin{lemma}
\label{lem:LowerBound_of_the_Spectrum}
The function $\sigma(\cdot)$ is strictly decreasing on $\R$ and
\[
\sigma(4\pi^2)=0.
\]
\end{lemma}

Arguing similarly to Proposition~\ref{theorem:Friedrichs}, we can show that
\cite[Lemma 3.1]{LiK01}:
\begin{proposition}
\label{prop:KS-inequality}
Let $x\in H^2_0(0,1)$. Then
\begin{eqnarray}
\int_0^1 x^2_{zz}(z)dz - \lambda\int_0^1 x^2_{z}(z)dz\geq \sigma(\lambda)\int_0^1 x^2(z)dz.
\label{eq:KS_LiK01_inequality}
\end{eqnarray}
\end{proposition}

\begin{proof}
Pick any $x \in C^\infty_c(0,1)$. We develop it into a series with respect to the orthonormal basis $\{\phi_i:i\in\N\}$ of $L^2(0,1)$:
\[
x(z) = \sum_{i=1}^\infty c_i \phi^{<i>}(z),\quad z \in (0,1),
\]
where $c_i := \scalp{x}{\phi^{<i>}}$, $i\in\N$. 

We have that
\begin{eqnarray*}
\int_0^1 x^2_{zz}(z)dz &=& \scalp{\sum_{i=1}^\infty c_i \phi^{<i>}_{zz}}{\sum_{j=1}^\infty c_j \phi^{<j>}_{zz}}\\
&=& \sum_{i=1}^\infty \sum_{j=1}^\infty c_i c_j \scalp{ \phi^{<i>}_{zz}}{\phi^{<j>}_{zz}}\\
&=& \sum_{i=1}^\infty \sum_{j=1}^\infty c_i c_j \int_0^1 \phi^{<i>}_{zz}(z) \phi^{<j>}_{zz}(z) dz.
\end{eqnarray*}
Using the integration by parts twice, we see that
\begin{eqnarray}
\int_0^1 x^2_{zz}(z)dz &=& \sum_{i=1}^\infty \sum_{j=1}^\infty c_i c_j \int_0^1 \phi^{<i>}(z) \phi^{<j>}_{zzzz}(z) dz.
\label{eq:KS_Dist_Inp_HelpIneq_Tmp_1}
\end{eqnarray}
Analogously, we obtain that
\begin{eqnarray}
\int_0^1 x^2_{z}(z)dz &=& \sum_{i=1}^\infty \sum_{j=1}^\infty c_i c_j \int_0^1 \phi^{<i>}_{z}(z) \phi^{<j>}_{z}(z) dz\nonumber\\
&=& -\sum_{i=1}^\infty \sum_{j=1}^\infty c_i c_j \int_0^1 \phi^{<i>}(z) \phi^{<j>}_{zz}(z) dz.
\label{eq:KS_Dist_Inp_HelpIneq_Tmp_2}
\end{eqnarray}
Using \eqref{eq:KS_Dist_Inp_HelpIneq_Tmp_1}, \eqref{eq:KS_Dist_Inp_HelpIneq_Tmp_2}, we have that
\begin{eqnarray*}
\int_0^1 x^2_{zz}(z)dz - \lambda\int_0^1 x^2_{z}(z)dz =
\sum_{i=1}^\infty \sum_{j=1}^\infty c_i c_j \int_0^1 \phi^{<i>}(z) \big(\phi^{<j>}_{zzzz}(z) + \lambda \phi^{<j>}_{zz}(z)\big) dz.
\end{eqnarray*}
Since $\phi^{<j>}$ satisfies the eigenvalue problem \eqref{eq:KS_Eigenvalue_Problem}, it holds also that
\begin{eqnarray*}
\int_0^1 x^2_{zz}(z)dz - \lambda\int_0^1 x^2_{z}(z)dz = \sum_{i=1}^\infty \sum_{j=1}^\infty c_i c_j \int_0^1 \phi^{<i>}(z) \sigma_j \phi^{<j>}(z) dz.
\end{eqnarray*}
Since $\{\phi^{<j>}:\ j\in\N\}$ is an orthonormal basis of $L^2(0,1)$, we proceed to
\begin{eqnarray*}
\int_0^1 x^2_{zz}(z)dz - \lambda\int_0^1 x^2_{z}(z)dz &=& \sum_{i=1}^\infty c^2_i \sigma_i\\
&\geq& \sigma(\lambda) \sum_{i=1}^\infty c^2_i\\
&=& \sigma(\lambda) \int_0^1 x^2(z)dz.
\end{eqnarray*}
Here the last equality follows from the definition of $c_i$. 

As $ H^2_0(0,1)$ is a closure of $C^\infty_c(0,L)$ in $H^2(0,1)$-norm, the last inequality holds also for all $x\in H^2_0(0,1)$.
\end{proof}

\section{Concluding remarks}

In Section~\ref{sec:Sobolev_spaces}, we follow \cite[Chapter 5]{Eva10} in our treatment of weak derivatives and Sobolev spaces.
The proof of Proposition~\ref{theorem:Friedrichs}, achieving the optimal constant, can be found in many sources, e.g., in \cite[Appendix A]{Mur03}.
Section~\ref{sec:An inequality for the Kuramoto-Sivashinskiy equation} is motivated by \cite{LiK01}.


\end{appendices}

\addcontentsline{toc}{chapter}{\numberline{}Bibliography}
\bibliographystyle{abbrv}
\bibliography{Mir_LitList_NoMir,MyPublications}

\begin{thebibliography}{100}

\bibitem{AlT07}
C.~D. Aliprantis and R.~Tourky.
\newblock {\em Cones and Duality}.
\newblock American Mathematical Society, 2007.

\bibitem{Ama76}
H.~Amann.
\newblock Fixed point equations and nonlinear eigenvalue problems in ordered
  {B}anach spaces.
\newblock {\em SIAM Review}, 18(4):620--709, 1976.

\bibitem{Ang99}
D.~Angeli.
\newblock Intrinsic robustness of global asymptotic stability.
\newblock {\em Systems \& Control Letters}, 38(4-5):297--307, 1999.

\bibitem{Ang09}
D.~Angeli.
\newblock Further results on incremental input-to-state stability.
\newblock {\em IEEE Transactions on Automatic Control}, 54(6):1386--1391, 2009.

\bibitem{AnA07}
D.~Angeli and A.~Astolfi.
\newblock A tight small gain theorem for not necessarily {ISS} systems.
\newblock {\em Systems \& Control Letters}, 56:87--91, 2007.

\bibitem{AIS04b}
D.~Angeli, B.~Ingalls, E.~D. Sontag, and Y.~Wang.
\newblock Separation principles for input-output and integral-input-to-state
  stability.
\newblock {\em SIAM Journal on Control and Optimization}, 43(1):256--276, 2004.

\bibitem{ASW00b}
D.~Angeli, E.~Sontag, and Y.~Wang.
\newblock Further equivalences and semiglobal versions of integral input to
  state stability.
\newblock {\em Dynamics and Control}, 10(2):127--149, 2000.

\bibitem{ASW00}
D.~Angeli, E.~D. Sontag, and Y.~Wang.
\newblock A characterization of integral input-to-state stability.
\newblock {\em IEEE Transactions on Automatic Control}, 45(6):1082--1097, 2000.

\bibitem{AAS02}
M.~Arcak, D.~Angeli, and E.~Sontag.
\newblock A unifying integral {ISS} framework for stability of nonlinear
  cascades.
\newblock {\em SIAM Journal on Control and Optimization}, 40(6):1888--1904,
  2002.

\bibitem{ArK01}
M.~Arcak and P.~Kokotovi{\'c}.
\newblock Nonlinear observers: {A} circle criterion design and robustness
  analysis.
\newblock {\em Automatica}, 37(12):1923--1930, 2001.

\bibitem{ABH11}
W.~Arendt, C.~J. Batty, M.~Hieber, and F.~Neubrander.
\newblock {\em Vector-Valued {L}aplace Transforms and {C}auchy Problems}.
\newblock Springer Science \& Business Media, 2011.

\bibitem{ArN09}
W.~{Arendt} and R.~{Nittka}.
\newblock Equivalent complete norms and positivity.
\newblock {\em Archiv der Mathematik}, 92(5):414--427, 2009.

\bibitem{AuF09}
J.~Aubin and H.~Frankowska.
\newblock {\em Set-Valued Analysis}.
\newblock Birkh{\"a}user Boston, 2009.

\bibitem{Aug16}
B.~Augner.
\newblock {\em Stabilisation of infinite-dimensional port-{H}amiltonian systems
  via dissipative boundary feedback}.
\newblock PhD thesis, Bergische Universit{\"a}t Wuppertal, 2016.

\bibitem{AuW96}
B.~Aulbach and T.~Wanner.
\newblock Integral manifolds for {C}arath\'eodory type differential equations
  in {B}anach spaces.
\newblock In {\em Six Lectures on Dynamical Systems}, pages 45--119. World
  Scientific, River Edge, 1996.

\bibitem{AvL98}
G.~Avalos and I.~Lasiecka.
\newblock The strong stability of a semigroup arising from a coupled
  hyperbolic/parabolic system.
\newblock {\em Semigroup Forum}, 57(2):278--292, 1998.

\bibitem{BaK00}
A.~Balogh and M.~Krstic.
\newblock Boundary control of the {K}orteweg-de {V}ries-{B}urgers equation:
  {F}urther results on stabilization and well-posedness, with numerical
  demonstration.
\newblock {\em IEEE Transactions on Automatic Control}, 45(9):1739--1745, 2000.

\bibitem{BPD02}
B.~Bamieh, F.~Paganini, and M.~A. Dahleh.
\newblock Distributed control of spatially invariant systems.
\newblock {\em IEEE Transactions on Automatic Control}, 47(7):1091--1107, 2002.

\bibitem{BaV05}
B.~Bamieh and P.~G. Voulgaris.
\newblock A convex characterization of distributed control problems in
  spatially invariant systems with communication constraints.
\newblock {\em Systems \& Control Letters}, 54(6):575--583, 2005.

\bibitem{BLJ18}
A.~Bao, T.~Liu, Z.-P. Jiang, and L.~Zhang.
\newblock A nonlinear small-gain theorem for large-scale infinite-dimensional
  systems.
\newblock {\em Journal of Systems Science and Complexity}, 31(1):188--199,
  2018.

\bibitem{Bar10}
V.~Barbu.
\newblock {\em Nonlinear Differential Equations of Monotone Types in Banach
  Spaces}.
\newblock Springer Science \& Business Media, 2010.

\bibitem{BeJ17}
B.~Besselink and K.~H. Johansson.
\newblock String stability and a delay-based spacing policy for vehicle
  platoons subject to disturbances.
\newblock {\em IEEE Transactions on Automatic Control}, 62(9):4376--4391, 2017.

\bibitem{Bha66}
N.~P. Bhatia.
\newblock Weak attractors in dynamical systems.
\newblock {\em Boletin Sociedad Matematica Mexicana}, 11:56--64, 1966.

\bibitem{BhS02}
N.~P. Bhatia and G.~P. Szeg{\"o}.
\newblock {\em Stability Theory of Dynamical Systems}.
\newblock Springer Science \& Business Media, 2002.

\bibitem{BCC16}
B.~J. Bialy, I.~Chakraborty, S.~C. Cekic, and W.~E. Dixon.
\newblock Adaptive boundary control of store induced oscillations in a flexible
  aircraft wing.
\newblock {\em Automatica}, 70:230--238, 2016.

\bibitem{AWP12}
F.~Bribiesca~Argomedo, E.~Witrant, and C.~Prieur.
\newblock {$D^1$}-{I}nput-to-state stability of a time-varying nonhomogeneous
  diffusive equation subject to boundary disturbances.
\newblock In {\em Proc. of 2012 American Control Conference}, pages 2978--2983,
  2012.

\bibitem{BrL66}
A.~Bruckner and J.~Leonard.
\newblock Derivatives.
\newblock {\em The American Mathematical Monthly}, 73(4):24--56, 1966.

\bibitem{BDK74}
C.~Bruni, G.~Dipillo, and G.~Koch.
\newblock Bilinear systems: {A}n appealing class of ``nearly linear'' systems
  in theory and applications.
\newblock {\em IEEE Transactions on Automatic Control}, 19(4):334--348, 1974.

\bibitem{But10}
P.~Butkovi{\v{c}}.
\newblock {\em Max-Linear Systems: Theory and Algorithms}.
\newblock Springer Science \& Business Media, 2010.

\bibitem{CaT09}
C.~Cai and A.~Teel.
\newblock Characterizations of input-to-state stability for hybrid systems.
\newblock {\em Systems \& Control Letters}, 58(1):47--53, 2009.

\bibitem{CTG07}
C.~Cai, A.~R. Teel, and R.~Goebel.
\newblock Smooth {L}yapunov functions for hybrid systems--part {I}: {E}xistence
  is equivalent to robustness.
\newblock {\em IEEE Transactions on Automatic Control}, 52(7):1264--1277, 2007.

\bibitem{CTG08}
C.~Cai, A.~R. Teel, and R.~Goebel.
\newblock Smooth {L}yapunov functions for hybrid systems part {II}:
  ({P}re){A}symptotically stable compact sets.
\newblock {\em IEEE Transactions on Automatic Control}, 53(3):734--748, 2008.

\bibitem{CaH98}
T.~Cazenave and A.~Haraux.
\newblock {\em An Introduction to Semilinear Evolution Equations}.
\newblock Oxford University Press, New York, 1998.

\bibitem{ChA08}
A.~Chaillet and D.~Angeli.
\newblock Integral input to state stable systems in cascade.
\newblock {\em Systems \& Control Letters}, 57(7):519--527, 2008.

\bibitem{CAI14}
A.~Chaillet, D.~Angeli, and H.~Ito.
\newblock Combining i{ISS} and {ISS} with respect to small inputs: {T}he strong
  i{ISS} property.
\newblock {\em IEEE Transactions on Automatic Control}, 59(9):2518--2524, 2014.

\bibitem{CAI14b}
A.~Chaillet, D.~Angeli, and H.~Ito.
\newblock Strong i{ISS} is preserved under cascade interconnection.
\newblock {\em Automatica}, 50(9):2424--2427, 2014.

\bibitem{CKP22}
A.~Chaillet, I.~Karafyllis, P.~Pepe, and Y.~Wang.
\newblock The {ISS} framework for time-delay systems: a survey.
\newblock {\em arXiv preprint arXiv:2206.06167}, 2022.

\bibitem{ChZ09}
W.-H. Chen and W.~X. Zheng.
\newblock Input-to-state stability and integral input-to-state stability of
  nonlinear impulsive systems with delays.
\newblock {\em Automatica}, 45(6):1481--1488, 2009.

\bibitem{CLS98}
F.~H. Clarke, Y.~Ledyaev, and R.~J. Stern.
\newblock Asymptotic stability and smooth {L}yapunov functions.
\newblock {\em Journal of Differential Equations}, 149(1):69--114, 1998.

\bibitem{Col51}
J.~D. Cole.
\newblock On a quasi-linear parabolic equation occurring in aerodynamics.
\newblock {\em Quarterly of Applied Mathematics}, 9(3):225--236, 1951.

\bibitem{CoL15}
J.-M. Coron and Q.~L{\"u}.
\newblock {F}redholm transform and local rapid stabilization for a
  {K}uramoto--{S}ivashinsky equation.
\newblock {\em Journal of Differential Equations}, 259(8):3683--3729, 2015.

\bibitem{CrP69}
M.~G. Crandall and A.~Pazy.
\newblock Semi-groups of nonlinear contractions and dissipative sets.
\newblock {\em Journal of Functional Analysis}, 3(3):376--418, 1969.

\bibitem{CIZ09}
R.~Curtain, O.~V. Iftime, and H.~Zwart.
\newblock System theoretic properties of a class of spatially invariant
  systems.
\newblock {\em Automatica}, 45(7):1619--1627, 2009.

\bibitem{CuZ20}
R.~Curtain and H.~Zwart.
\newblock {\em Introduction to Infinite-Dimensional Systems Theory: {A}
  State-Space Approach}.
\newblock Springer, 2020.

\bibitem{CuZ95}
R.~F. Curtain and H.~Zwart.
\newblock {\em An Introduction to Infinite-Dimensional Linear Systems Theory}.
\newblock Springer, New York, 1995.

\bibitem{DaK74}
J.~Daleckii and M.~Krein.
\newblock {\em Stability of Solutions of Differential Equations in {B}anach
  Space}.
\newblock American Mathematical Society, 1974.

\bibitem{DaD03}
R.~D'Andrea and G.~E. Dullerud.
\newblock Distributed control design for spatially interconnected systems.
\newblock {\em IEEE Transactions on Automatic Control}, 48(9):1478--1495, 2003.

\bibitem{DIW11}
S.~Dashkovskiy, H.~Ito, and F.~Wirth.
\newblock On a small gain theorem for {ISS} networks in dissipative {L}yapunov
  form.
\newblock {\em European Journal of Control}, 17(4):357--365, 2011.

\bibitem{DaK13}
S.~Dashkovskiy and M.~Kosmykov.
\newblock Input-to-state stability of interconnected hybrid systems.
\newblock {\em Automatica}, 49(4):1068--1074, 2013.

\bibitem{DKM12}
S.~Dashkovskiy, M.~Kosmykov, A.~Mironchenko, and L.~Naujok.
\newblock Stability of interconnected impulsive systems with and without time
  delays, using {L}yapunov methods.
\newblock {\em Nonlinear Analysis: Hybrid Systems}, 6(3):899--915, 2012.

\bibitem{DaM13}
S.~Dashkovskiy and A.~Mironchenko.
\newblock Input-to-state stability of infinite-dimensional control systems.
\newblock {\em Mathematics of Control, Signals, and Systems}, 25(1):1--35,
  2013.

\bibitem{DaM13b}
S.~Dashkovskiy and A.~Mironchenko.
\newblock Input-to-state stability of nonlinear impulsive systems.
\newblock {\em SIAM Journal on Control and Optimization}, 51(3):1962--1987,
  2013.

\bibitem{DMS19a}
S.~Dashkovskiy, A.~Mironchenko, J.~Schmid, and F.~Wirth.
\newblock Stability of infinitely many interconnected systems.
\newblock {\em IFAC-PapersOnLine}, 52(16):550--555, 2019.

\bibitem{DaP20}
S.~Dashkovskiy and S.~Pavlichkov.
\newblock Stability conditions for infinite networks of nonlinear systems and
  their application for stabilization.
\newblock {\em Automatica}, 112:108643, 2020.

\bibitem{DRW10}
S.~Dashkovskiy, B.~R\"{u}ffer, and F.~Wirth.
\newblock Small gain theorems for large scale systems and construction of {ISS}
  {L}yapunov functions.
\newblock {\em SIAM Journal on Control and Optimization}, 48(6):4089--4118,
  2010.

\bibitem{DRW06b}
S.~Dashkovskiy, B.~S. R\"{u}ffer, and F.~R. Wirth.
\newblock On the construction of {ISS} {L}yapunov functions for networks of
  {ISS} systems.
\newblock In {\em Proc. of 17th Int. Symposium on Mathematical Theory of
  Networks and Systems}, pages 77--82, 2006.

\bibitem{DRW07}
S.~Dashkovskiy, B.~S. R\"{u}ffer, and F.~R. Wirth.
\newblock An {ISS} small gain theorem for general networks.
\newblock {\em Mathematics of Control, Signals, and Systems}, 19(2):93--122,
  2007.

\bibitem{DES11}
S.~N. Dashkovskiy, D.~V. Efimov, and E.~D. Sontag.
\newblock Input to state stability and allied system properties.
\newblock {\em Automation and Remote Control}, 72(8):1579--1614, 2011.

\bibitem{DaR10}
S.~N. Dashkovskiy and B.~S. R\"{u}ffer.
\newblock Local {ISS} of large-scale interconnections and estimates for
  stability regions.
\newblock {\em Systems \& Control Letters}, 59:241--247, 2010.

\bibitem{DeV09}
C.~A. Desoer and M.~Vidyasagar.
\newblock {\em Feedback Systems: Input-Output Properties}.
\newblock SIAM, Philadelphia, PA, 2009.

\bibitem{DiK15}
M.~Diagne and M.~Krstic.
\newblock State-dependent input delay-compensated bang-bang control:
  {A}pplication to 3{D} printing based on screw-extruder.
\newblock In {\em Proc. of 2015 American Control Conference}, pages 5653--5658,
  2015.

\bibitem{Din1878}
U.~Dini.
\newblock {\em Fondamenti per la Teorica Delle Funzioni di Variabili Reali}.
\newblock T. Nistri e C., 1878.

\bibitem{EmT00}
Z.~Emirsajlow and S.~Townley.
\newblock From {PDE}s with boundary control to the abstract state equation with
  an unbounded input operator: {A} tutorial.
\newblock {\em European Journal of Control}, 6(1):27--49, 2000.

\bibitem{EMJ17}
T.~Endo, F.~Matsuno, and Y.~Jia.
\newblock Boundary cooperative control by flexible {T}imoshenko arms.
\newblock {\em Automatica}, 81:377--389, 2017.

\bibitem{EnN00}
K.-J. Engel and R.~Nagel.
\newblock {\em One-Parameter Semigroups for Linear Evolution Equations}.
\newblock Springer, New York, 2000.

\bibitem{Eva10}
L.~C. Evans.
\newblock {\em Partial Differential Equations}.
\newblock American Mathematical Society, 2010.

\bibitem{FHH11}
M.~Fabian, P.~Habala, P.~H{\'a}jek, V.~Montesinos, and V.~Zizler.
\newblock {\em Banach Space Theory: The Basis for Linear and Nonlinear
  Analysis}.
\newblock Springer Science \& Business Media, 2011.

\bibitem{Fat68}
H.~O. Fattorini.
\newblock Boundary control systems.
\newblock {\em SIAM Journal on Control}, 6(3):349--385, 1968.

\bibitem{Glu16}
J.~{G}l{\"u}ck.
\newblock {\em Invariant sets and long time behaviour of operator semigroups}.
\newblock PhD thesis, Universit{\"a}t Ulm, 2016.
\newblock DOI: 10.18725/OPARU-4238.

\bibitem{GlM21}
J.~Gl\"uck and A.~Mironchenko.
\newblock Stability criteria for positive linear discrete-time systems.
\newblock {\em Positivity}, 25(5):2029--2059, 2021.

\bibitem{GlW20}
J.~{Gl\"uck} and M.~R. {Weber}.
\newblock Almost interior points in ordered {B}anach spaces and the long-term
  behaviour of strongly positive operator semigroups.
\newblock {\em Studia Mathematica}, 254(3):237--263, 2020.

\bibitem{God75}
A.~Godunov.
\newblock {P}eano's theorem in {B}anach spaces.
\newblock {\em Functional Analysis and Its Applications}, 9(1):53--55, 1975.

\bibitem{Gra95}
P.~Grabowski.
\newblock Admissibility of observation functionals.
\newblock {\em International Journal of Control}, 62(5):1161--1173, 1995.

\bibitem{Gru02}
L.~Gr\"{u}ne.
\newblock Input-to-state dynamical stability and its {L}yapunov function
  characterization.
\newblock {\em IEEE Transactions on Automatic Control}, 47(9):1499--1504, 2002.

\bibitem{Haa06}
M.~Haase.
\newblock {\em The Functional Calculus for Sectorial Operators}.
\newblock Springer, 2006.

\bibitem{HaT06}
J.~W. Hagood and B.~S. Thomson.
\newblock Recovering a function from a {D}ini derivative.
\newblock {\em The American Mathematical Monthly}, 113(1):34--46, 2006.

\bibitem{Hah67}
W.~Hahn.
\newblock {\em Stability of Motion}.
\newblock Springer, New York, 1967.

\bibitem{HaJ10}
P.~H{\'a}jek and M.~Johanis.
\newblock On {P}eano's theorem in {B}anach spaces.
\newblock {\em Journal of Differential Equations}, 249(12):3342--3351, 2010.

\bibitem{HaW97}
S.~Hansen and G.~Weiss.
\newblock New results on the operator {C}arleson measure criterion.
\newblock {\em IMA Journal of Mathematical Control and Information},
  14(1):3--32, 1997.

\bibitem{HaS11}
F.~M. Hante and M.~Sigalotti.
\newblock Converse {L}yapunov theorems for switched systems in {B}anach and
  {H}ilbert spaces.
\newblock {\em SIAM Journal on Control and Optimization}, 49(2):752--770, 2011.

\bibitem{HLP52}
G.~H. Hardy, J.~E. Littlewood, and G.~P{\'o}lya.
\newblock {\em Inequalities}.
\newblock Cambridge university press, 1952.

\bibitem{HCZ19}
A.~Hastir, F.~Califano, and H.~Zwart.
\newblock Well-posedness of infinite-dimensional linear systems with nonlinear
  feedback.
\newblock {\em Systems \& Control Letters}, 128:19--25, 2019.

\bibitem{HMR15}
J.~Heigel, P.~Michaleris, and E.~Reutzel.
\newblock Thermo-mechanical model development and validation of directed energy
  deposition additive manufacturing of {T}i--6{A}l--4{V}.
\newblock {\em Additive Manufacturing}, 5:9--19, 2015.

\bibitem{Hen81}
D.~Henry.
\newblock {\em Geometric Theory of Semilinear Parabolic Equations}.
\newblock Springer, Berlin, 1981.

\bibitem{HLT08}
J.~P. Hespanha, D.~Liberzon, and A.~R. Teel.
\newblock Lyapunov conditions for input-to-state stability of impulsive
  systems.
\newblock {\em Automatica}, 44(11):2735--2744, 2008.

\bibitem{Hil91}
D.~J. Hill.
\newblock A generalization of the small-gain theorem for nonlinear feedback
  systems.
\newblock {\em Automatica}, 27:1043--1045, 1991.

\bibitem{HiP00}
E.~Hille and R.~S. Phillips.
\newblock {\em Functional Analysis and Semi-Groups}.
\newblock American Mathematical Society, Providence, R. I., 1974.

\bibitem{HJS22}
R.~Hosfeld, B.~Jacob, and F.~L. Schwenninger.
\newblock Integral input-to-state stability of unbounded bilinear control
  systems.
\newblock {\em Mathematics of Control, Signals, and Systems}, pages 1--23,
  2022.

\bibitem{HuP21}
K.~Huhtala and L.~Paunonen.
\newblock Approximate local output regulation for a class of nonlinear fluid
  flows.
\newblock {\em European Journal of Control}, 62:136--142, 2021.

\bibitem{HMM13}
D.~Hundertmark, M.~Meyries, L.~Machinek, and R.~Schnaubelt.
\newblock Operator semigroups and dispersive equations.
\newblock In {\em Proc. of 16th Internet Seminar on Evolution Equations}, 2013.

\bibitem{InS02}
B.~Ingalls and E.~D. Sontag.
\newblock A small-gain theorem with applications to input/output systems,
  incremental stability, detectability, and interconnections.
\newblock {\em Journal of the Franklin Institute}, 339:211--229, 2002.

\bibitem{Ito06}
H.~Ito.
\newblock State-dependent scaling problems and stability of interconnected
  {iISS} and {ISS} systems.
\newblock {\em IEEE Transactions on Automatic Control}, 51:1626--1643, 2006.

\bibitem{Ito13}
H.~Ito.
\newblock Utility of i{ISS} in composing {L}yapunov functions.
\newblock In {\em Proc. of 9th IFAC Symposium on Nonlinear Control Systems},
  pages 723--730, 2013.

\bibitem{ItJ09}
H.~Ito and Z.-P. Jiang.
\newblock Necessary and sufficient small gain conditions for integral
  input-to-state stable systems: {A} {L}yapunov perspective.
\newblock {\em IEEE Transactions on Automatic Control}, 54(10):2389--2404,
  2009.

\bibitem{IJD13}
H.~Ito, Z.-P. Jiang, S.~Dashkovskiy, and B.~R\"{u}ffer.
\newblock Robust stability of networks of i{ISS} systems: {C}onstruction of
  sum-type {L}yapunov functions.
\newblock {\em IEEE Transactions on Automatic Control}, 58(5):1192--1207, 2013.

\bibitem{JDP95}
B.~Jacob, V.~Dragan, and A.~J. Pritchard.
\newblock Infinite dimensional time varying systems with nonlinear output
  feedback.
\newblock {\em Integral Equations and Operator Theory}, 22(4):440--462, 1995.

\bibitem{JMP20}
B.~Jacob, A.~Mironchenko, J.~R. Partington, and F.~Wirth.
\newblock Noncoercive {L}yapunov functions for input-to-state stability of
  infinite-dimensional systems.
\newblock {\em SIAM Journal on Control and Optimization}, 58(5):2952--2978,
  2020.

\bibitem{JNP18}
B.~Jacob, R.~Nabiullin, J.~R. Partington, and F.~L. Schwenninger.
\newblock Infinite-dimensional input-to-state stability and {O}rlicz spaces.
\newblock {\em SIAM Journal on Control and Optimization}, 56(2):868--889, 2018.

\bibitem{JaP04}
B.~Jacob and J.~R. Partington.
\newblock Admissibility of control and observation operators for semigroups:
  {A} survey.
\newblock In {\em Current Trends in Operator Theory and its Applications},
  pages 199--221. Birkh\"{a}user Basel, 2004.

\bibitem{JSZ19}
B.~Jacob, F.~L. Schwenninger, and H.~Zwart.
\newblock On continuity of solutions for parabolic control systems and
  input-to-state stability.
\newblock {\em Journal of Differential Equations}, 266:6284--6306, 2019.

\bibitem{JaW15}
B.~Jacob and S.-A. Wegner.
\newblock Asymptotics of evolution equations beyond banach spaces.
\newblock In {\em Semigroup Forum}, volume~91, pages 347--377. Springer, 2015.

\bibitem{JaZ12}
B.~Jacob and H.~J. Zwart.
\newblock {\em Linear Port-{H}amiltonian Systems on Infinite-Dimensional
  Spaces}.
\newblock Springer, Basel, 2012.

\bibitem{JLR09}
B.~Jayawardhana, H.~Logemann, and E.~P. Ryan.
\newblock Input-to-state stability of differential inclusions with applications
  to hysteretic and quantized feedback systems.
\newblock {\em SIAM Journal on Control and Optimization}, 48(2):1031--1054,
  2009.

\bibitem{JLZ15}
W.~Jiang, B.~Liu, and Z.~Zhang.
\newblock Robust observability for regular linear systems under nonlinear
  perturbation.
\newblock {\em Electronic Journal of Differential Equations}, 218:1--14, 2015.

\bibitem{JMW96}
Z.-P. Jiang, I.~M.~Y. Mareels, and Y.~Wang.
\newblock A {L}yapunov formulation of the nonlinear small-gain theorem for
  interconnected {ISS} systems.
\newblock {\em Automatica}, 32(8):1211--1215, 1996.

\bibitem{JTP94}
Z.-P. Jiang, A.~R. Teel, and L.~Praly.
\newblock Small-gain theorem for {ISS} systems and applications.
\newblock {\em Mathematics of Control, Signals, and Systems}, 7(2):95--120,
  1994.

\bibitem{JiW01}
Z.-P. Jiang and Y.~Wang.
\newblock Input-to-state stability for discrete-time nonlinear systems.
\newblock {\em Automatica}, 37(6):857--869, 2001.

\bibitem{JoB05b}
M.~R. Jovanovi\'{c} and B.~Bamieh.
\newblock On the ill-posedness of certain vehicular platoon control problems.
\newblock {\em IEEE Transactions on Automatic Control}, 50(9):1307--1321, 2005.

\bibitem{KFA69}
R.~E. Kalman, P.~L. Falb, and M.~A. Arbib.
\newblock {\em Topics in Mathematical System Theory}.
\newblock McGraw-Hill New York, 1969.

\bibitem{KaF18b}
W.~Kang and E.~Fridman.
\newblock Distributed sampled-data control of {K}uramoto--{S}ivashinsky
  equation.
\newblock {\em Automatica}, 95:514--524, 2018.

\bibitem{KaK12}
R.~Kannan and C.~K. Krueger.
\newblock {\em Advanced Analysis: On the Real Line}.
\newblock Springer Science \& Business Media, 2012.

\bibitem{Kar07a}
I.~Karafyllis.
\newblock A system-theoretic framework for a wide class of systems {I}:
  {A}pplications to numerical analysis.
\newblock {\em Journal of Mathematical Analysis and Applications},
  328(2):876--899, 2007.

\bibitem{KaJ11b}
I.~Karafyllis and Z.-P. Jiang.
\newblock {\em Stability and Stabilization of Nonlinear Systems}.
\newblock Springer, London, 2011.

\bibitem{KaJ11}
I.~Karafyllis and Z.-P. Jiang.
\newblock A vector small-gain theorem for general non-linear control systems.
\newblock {\em IMA Journal of Mathematical Control and Information},
  28:309--344, 2011.

\bibitem{KaJ12}
I.~Karafyllis and Z.-P. Jiang.
\newblock A new small-gain theorem with an application to the stabilization of
  the chemostat.
\newblock {\em International Journal of Robust and Nonlinear Control},
  22(14):1602--1630, 2012.

\bibitem{KaK16b}
I.~Karafyllis and M.~Krstic.
\newblock {ISS} with respect to boundary disturbances for 1-{D} parabolic
  {PDE}s.
\newblock {\em IEEE Transactions on Automatic Control}, 61(12):3712--3724,
  2016.

\bibitem{KaK17a}
I.~Karafyllis and M.~Krstic.
\newblock {ISS} in different norms for {1-D} parabolic {PDE}s with boundary
  disturbances.
\newblock {\em SIAM Journal on Control and Optimization}, 55(3):1716--1751,
  2017.

\bibitem{KaK17b}
I.~Karafyllis and M.~Krstic.
\newblock Decay estimates for 1-{D} parabolic {PDE}s with boundary
  disturbances.
\newblock {\em ESAIM Control, Optimisation and Calculus of Variations},
  24(4):1511--1540, 2018.

\bibitem{KaK18}
I.~Karafyllis and M.~Krstic.
\newblock Small-gain stability analysis of certain hyperbolic-parabolic {PDE}
  loops.
\newblock {\em Systems \& Control Letters}, 118:52--61, 2018.

\bibitem{KaK19}
I.~Karafyllis and M.~Krstic.
\newblock {\em Input-to-State Stability for {PDE}s}.
\newblock Springer, Cham, 2019.

\bibitem{KaK19b}
I.~Karafyllis and M.~Krstic.
\newblock Small-gain-based boundary feedback design for global exponential
  stabilization of one-dimensional semilinear parabolic {PDE}s.
\newblock {\em SIAM Journal on Control and Optimization}, 57(3):2016--2036,
  2019.

\bibitem{KPC22}
I.~Karafyllis, P.~Pepe, A.~Chaillet, and Y.~Wang.
\newblock Is global asymptotic stability necessarily uniform for time-delay
  systems?
\newblock {\em arXiv preprint arXiv:2202.11298}, 2022.

\bibitem{KPJ08}
I.~Karafyllis, P.~Pepe, and Z.-P. Jiang.
\newblock Input-to-output stability for systems described by retarded
  functional differential equations.
\newblock {\em European Journal of Control}, 14(6):539--555, 2008.

\bibitem{KPZ86}
M.~Kardar, G.~Parisi, and Y.-C. Zhang.
\newblock Dynamic scaling of growing interfaces.
\newblock {\em Physical Review Letters}, 56(9):889, 1986.

\bibitem{KMS21}
C.~Kawan, A.~Mironchenko, A.~Swikir, N.~Noroozi, and M.~Zamani.
\newblock A {L}yapunov-based small-gain theorem for infinite networks.
\newblock {\em IEEE Transactions on Automatic Control}, 66(12):5830--5844,
  2021.

\bibitem{KMZ22}
C.~Kawan, A.~Mironchenko, and M.~Zamani.
\newblock A {L}yapunov-based {ISS} small-gain theorem for infinite networks of
  nonlinear systems.
\newblock {\em Appeared online at IEEE Transactions on Automatic Control,
  https://ieeexplore.ieee.org/document/9813386}, 2022.

\bibitem{Kel14}
C.~M. Kellett.
\newblock A compendium of comparison function results.
\newblock {\em Mathematics of Control, Signals, and Systems}, 26(3):339--374,
  2014.

\bibitem{Kha02}
H.~K. Khalil.
\newblock {\em Nonlinear Systems}.
\newblock Prentice Hall, New York, 2002.

\bibitem{Kha03}
A.~Y. Khapalov.
\newblock Controllability of the semilinear parabolic equation governed by a
  multiplicative control in the reaction term: {A} qualitative approach.
\newblock {\em SIAM Journal on Control and Optimization}, 41(6):1886--1900,
  2003.

\bibitem{KLS89}
M.~A. {Krasnosel'skii}, E.~A. {Lifshits}, and A.~V. {Sobolev}.
\newblock {\em Positive Linear Systems. -- {T}he Method of Positive Operators}.
\newblock Berlin: Heldermann-Verlag, 1989.

\bibitem{KSW01}
M.~Krichman, E.~D. Sontag, and Y.~Wang.
\newblock Input-output-to-state stability.
\newblock {\em SIAM Journal on Control and Optimization}, 39(6):1874--1928,
  2001.

\bibitem{KKK95}
M.~Krstic, I.~Kanellakopoulos, and P.~V. Kokotovic.
\newblock {\em Nonlinear and Adaptive Control Design}.
\newblock Wiley, 1995.

\bibitem{KrS08}
M.~Krstic and A.~Smyshlyaev.
\newblock {\em Boundary Control of {PDE}s: {A} Course on Backstepping Designs}.
\newblock Society for Industrial and Applied Mathematics, Philadelphia, PA,
  2008.

\bibitem{Kur78}
Y.~Kuramoto.
\newblock Diffusion-induced chaos in reaction systems.
\newblock {\em Progress of Theoretical Physics Supplement}, 64:346--367, 1978.

\bibitem{KuT75}
Y.~Kuramoto and T.~Tsuzuki.
\newblock On the formation of dissipative structures in reaction-diffusion
  systems: {R}eductive perturbation approach.
\newblock {\em Progress of Theoretical Physics}, 54(3):687--699, 1975.

\bibitem{LaN03}
D.~S. Laila and D.~Ne{\v{s}}i{\'{c}}.
\newblock Discrete-time {L}yapunov-based small-gain theorem for parameterized
  interconnected {ISS} systems.
\newblock {\em IEEE Transactions on Automatic Control}, 48(10):1783--1788,
  2003.

\bibitem{LLM90}
V.~Lakshmikantham, S.~Leela, and A.~A. Martynyuk.
\newblock {\em Practical Stability of Nonlinear Systems}.
\newblock World Scientific, 1990.

\bibitem{LGE00}
L.~Lef{\`e}vre, D.~Georges, Z.~Emirsjlow, and S.~Townley.
\newblock Discussion on: '{F}rom {PDE}s with boundary control to the abstract
  state equation with an unbounded input operator: {A} tutorial'.
\newblock {\em European Journal of Control}, 6(1):50--53, 2000.

\bibitem{LSZ20}
H.~{Lhachemi}, D.~{Saussie}, G.~{Zhu}, and R.~{Shorten}.
\newblock Input-to-state stability of a clamped-free damped string in the
  presence of distributed and boundary disturbances.
\newblock {\em IEEE Transactions on Automatic Control}, 65(3):1248--1255, 2020.

\bibitem{LhS19}
H.~Lhachemi and R.~Shorten.
\newblock {ISS} property with respect to boundary disturbances for a class of
  {R}iesz-spectral boundary control systems.
\newblock {\em Automatica}, 109:108504, 2019.

\bibitem{LBG15}
H.~Li, R.~Baier, L.~Gr{\"u}ne, S.~Hafstein, and F.~Wirth.
\newblock Computation of local {ISS} {L}yapunov functions with low gains via
  linear programming.
\newblock Available at https://hal.inria.fr/hal-01101284, 2015.

\bibitem{LNT14}
D.~Liberzon, D.~Nesic, and A.~R. Teel.
\newblock Lyapunov-based small-gain theorems for hybrid systems.
\newblock {\em IEEE Transactions on Automatic Control}, 59(6):1395--1410, 2014.

\bibitem{LSW96}
Y.~Lin, E.~D. Sontag, and Y.~Wang.
\newblock A smooth converse {L}yapunov theorem for robust stability.
\newblock {\em SIAM Journal on Control and Optimization}, 34(1):124--160, 1996.

\bibitem{Lit89}
W.~Littman.
\newblock A generalization of a theorem of {D}atko and {P}azy.
\newblock In {\em Advances in Computing and Control}, pages 318--323. Springer,
  1989.

\bibitem{LLX11}
J.~Liu, X.~Liu, and W.-C. Xie.
\newblock Input-to-state stability of impulsive and switching hybrid systems
  with time-delay.
\newblock {\em Automatica}, 47(5):899--908, 2011.

\bibitem{LJH12}
T.~Liu, Z.-P. Jiang, and D.~J. Hill.
\newblock Lyapunov formulation of the {ISS} cyclic-small-gain theorem for
  hybrid dynamical networks.
\newblock {\em Nonlinear Analysis: Hybrid Systems}, 6(4):988--1001, 2012.

\bibitem{LJH14}
T.~Liu, Z.-P. Jiang, and D.~J. Hill.
\newblock {\em Nonlinear Control of Dynamic Networks}.
\newblock CRC Press, 2014.

\bibitem{LiK01}
W.-J. Liu and M.~Krsti\'{c}.
\newblock Stability enhancement by boundary control in the
  {K}uramoto-{S}ivashinsky equation.
\newblock {\em Nonlinear Analysis: Theory, Methods and Applications},
  43(4):485--507, 2001.

\bibitem{Lun12}
A.~Lunardi.
\newblock {\em Analytic Semigroups and Optimal Regularity in Parabolic
  Problems}.
\newblock Springer Science \& Business Media, 2012.

\bibitem{MaH92}
I.~M.~Y. Mareels and D.~J. Hill.
\newblock Monotone stability of nonlinear feedback systems.
\newblock {\em Journal of Mathematical Systems, Estimation, and Control},
  2:275--291, 1992.

\bibitem{MaP11}
F.~Mazenc and C.~Prieur.
\newblock Strict {L}yapunov functions for semilinear parabolic partial
  differential equations.
\newblock {\em Mathematical Control and Related Fields}, 1(2):231--250, 2011.

\bibitem{Mey12}
P.~Meyer-Nieberg.
\newblock {\em Banach Lattices}.
\newblock Springer Science \& Business Media, 2012.

\bibitem{Mir16}
A.~Mironchenko.
\newblock Local input-to-state stability: Characterizations and
  counterexamples.
\newblock {\em Systems \& Control Letters}, 87:23--28, 2016.

\bibitem{Mir17a}
A.~Mironchenko.
\newblock Uniform weak attractivity and criteria for practical global
  asymptotic stability.
\newblock {\em Systems \& Control Letters}, 105:92--99, 2017.

\bibitem{Mir19a}
A.~Mironchenko.
\newblock Criteria for input-to-state practical stability.
\newblock {\em IEEE Transactions on Automatic Control}, 64(1):298--304, 2019.

\bibitem{Mir20}
A.~Mironchenko.
\newblock Lyapunov functions for input-to-state stability of
  infinite-dimensional systems with integrable inputs.
\newblock In {\em IFAC-PapersOnLine}, volume~53, pages 5336--5341, 2020.

\bibitem{Mir21}
A.~Mironchenko.
\newblock Non-uniform {ISS} small-gain theorem for infinite networks.
\newblock {\em {IMA} Journal of Mathematical Control and Information},
  38(4):1029--1045, 2021.

\bibitem{Mir19b}
A.~Mironchenko.
\newblock Small gain theorems for general networks of heterogeneous
  infinite-dimensional systems.
\newblock {\em SIAM Journal on Control and Optimization}, 59(2):1393--1419,
  2021.

\bibitem{Mir22c}
A.~Mironchenko.
\newblock Well-posedness and properties of the flow for semilinear boundary
  control systems.
\newblock {\em Submitted}, 2022.

\bibitem{Mir23}
A.~Mironchenko.
\newblock {\em Input-to-State Stability}.
\newblock Accepted to Springer Communications and Control Engineering series,
  2023.

\bibitem{MiI15b}
A.~Mironchenko and H.~Ito.
\newblock Construction of {L}yapunov functions for interconnected parabolic
  systems: An i{ISS} approach.
\newblock {\em SIAM Journal on Control and Optimization}, 53(6):3364--3382,
  2015.

\bibitem{MiI16}
A.~Mironchenko and H.~Ito.
\newblock Characterizations of integral input-to-state stability for bilinear
  systems in infinite dimensions.
\newblock {\em Mathematical Control and Related Fields}, 6(3):447--466, 2016.

\bibitem{MKK19}
A.~Mironchenko, I.~Karafyllis, and M.~Krstic.
\newblock Monotonicity methods for input-to-state stability of nonlinear
  parabolic {PDE}s with boundary disturbances.
\newblock {\em SIAM Journal on Control and Optimization}, 57(1):510--532, 2019.

\bibitem{MKG20}
A.~Mironchenko, C.~Kawan, and J.~Gl\"uck.
\newblock Nonlinear small-gain theorems for input-to-state stability of
  infinite interconnections.
\newblock {\em Mathematics of Control, Signals, and Systems}, 33:573--615,
  2021.

\bibitem{MNK21}
A.~Mironchenko, N.~Noroozi, C.~Kawan, and M.~Zamani.
\newblock {ISS} small-gain criteria for infinite networks with linear gain
  functions.
\newblock {\em Systems \& Control Letters}, 157:105051, 2021.

\bibitem{MiP20}
A.~Mironchenko and C.~Prieur.
\newblock Input-to-state stability of infinite-dimensional systems: Recent
  results and open questions.
\newblock {\em SIAM Review}, 62(3):529--614, 2020.

\bibitem{MPW21}
A.~Mironchenko, C.~Prieur, and F.~Wirth.
\newblock Local stabilization of an unstable parabolic equation via saturated
  controls.
\newblock {\em IEEE Transactions on Automatic Control}, 66(5):2162--2176, 2021.

\bibitem{MiW15}
A.~Mironchenko and F.~Wirth.
\newblock A note on input-to-state stability of linear and bilinear
  infinite-dimensional systems.
\newblock In {\em Proc. of 54th IEEE Conference on Decision and Control}, pages
  495--500, 2015.

\bibitem{MiW17e}
A.~Mironchenko and F.~Wirth.
\newblock Input-to-state stability of time-delay systems: Criteria and open
  problems.
\newblock In {\em Proc. of 56th IEEE Conference on Decision and Control}, pages
  3719--3724. IEEE, 2017.

\bibitem{MiW18b}
A.~Mironchenko and F.~Wirth.
\newblock Characterizations of input-to-state stability for
  infinite-dimensional systems.
\newblock {\em IEEE Transactions on Automatic Control}, 63(6):1602--1617, 2018.

\bibitem{MiW17c}
A.~Mironchenko and F.~Wirth.
\newblock Lyapunov characterization of input-to-state stability for semilinear
  control systems over {B}anach spaces.
\newblock {\em Systems \& Control Letters}, 119:64--70, 2018.

\bibitem{MiW19b}
A.~Mironchenko and F.~Wirth.
\newblock Existence of non-coercive {L}yapunov functions is equivalent to
  integral uniform global asymptotic stability.
\newblock {\em Mathematics of Control, Signals, and Systems}, 31(4):1--26,
  2019.

\bibitem{MiW19a}
A.~Mironchenko and F.~Wirth.
\newblock Non-coercive {L}yapunov functions for infinite-dimensional systems.
\newblock {\em Journal of Differential Equations}, 105:7038--7072, 2019.

\bibitem{MYL18}
A.~Mironchenko, G.~Yang, and D.~Liberzon.
\newblock {L}yapunov small-gain theorems for networks of not necessarily {ISS}
  hybrid systems.
\newblock {\em Automatica}, 88:10--20, 2018.

\bibitem{MST99}
G.~A. {M}u{\~n}oz, Y.~{S}arantopoulos, and A.~{T}onge.
\newblock Complexifications of real {B}anach spaces, polynomials and
  multilinear maps.
\newblock {\em Studia Mathematica}, 134(1):1--33, 1999.

\bibitem{Mur03}
J.~D. Murray.
\newblock {\em Mathematical Biology. {II} Spatial Models and Biomedical
  Applications.}
\newblock Springer, New York, third edition, 2003.

\bibitem{NaS18}
R.~Nabiullin and F.~L. Schwenninger.
\newblock Strong input-to-state stability for infinite-dimensional linear
  systems.
\newblock {\em Mathematics of Control, Signals, and Systems}, 30(1):4, 2018.

\bibitem{NaB16}
V.~Natarajan and J.~Bentsman.
\newblock Approximate local output regulation for nonlinear distributed
  parameter systems.
\newblock {\em Mathematics of Control, Signals, and Systems}, 28(3):1--44,
  2016.

\bibitem{NZW19}
V.~Natarajan, H.-C. Zhou, G.~Weiss, and E.~Fridman.
\newblock Exact controllability of a class of nonlinear distributed parameter
  systems using back-and-forth iterations.
\newblock {\em International Journal of Control}, 92(1):145--162, 2019.

\bibitem{NST85}
B.~Nicolaenko, B.~Scheurer, and R.~Temam.
\newblock Some global dynamical properties of the {K}uramoto-{S}ivashinsky
  equations: Nonlinear stability and attractors.
\newblock {\em Physica D: Nonlinear Phenomena}, 16(2):155--183, 1985.

\bibitem{NKK15}
N.~Noroozi, A.~Khayatian, and H.~R. Karimi.
\newblock Semiglobal practical integral input-to-state stability for a family
  of parameterized discrete-time interconnected systems with application to
  sampled-data control systems.
\newblock {\em Nonlinear Analysis: Hybrid Systems}, 17:10--24, 2015.

\bibitem{NMK22}
N.~Noroozi, A.~Mironchenko, C.~Kawan, and M.~Zamani.
\newblock A small-gain theorem for set stability of infinite networks:
  {D}istributed observation and {ISS} for time-varying networks.
\newblock {\em European Journal of Control}, 67:100634, 2022.

\bibitem{Oos00}
J.~Oostveen.
\newblock {\em Strongly Stabilizable Distributed Parameter Systems}.
\newblock Society for Industrial and Applied Mathematics, Philadelphia, PA,
  2000.

\bibitem{PGC13}
A.~A. Paranjape, J.~Guan, S.~J. Chung, and M.~Krstic.
\newblock {PDE} boundary control for flexible articulated wings on a robotic
  aircraft.
\newblock {\em IEEE Transactions on Robotics}, 29(3):625--640, 2013.

\bibitem{PaY08}
P.~M. Pardalos and V.~Yatsenko.
\newblock {\em Optimization and Control of Bilinear Systems}.
\newblock Springer, Boston, MA, 2008.

\bibitem{Paz83}
A.~Pazy.
\newblock {\em Semigroups of Linear Operators and Applications to Partial
  Differential Equations}.
\newblock Springer, New York, 1983.

\bibitem{Pep21}
P.~Pepe.
\newblock A nonlinear version of {H}alanay's inequality for the uniform
  convergence to the origin.
\newblock {\em Mathematical Control \& Related Fields}, 2021.

\bibitem{PTS16}
C.~Prieur, S.~Tarbouriech, and J.~M.~G. da~Silva.
\newblock Wave equation with cone-bounded control laws.
\newblock {\em IEEE Transactions on Automatic Control}, 61(11):3452--3463,
  2016.

\bibitem{RZG17}
H.~Ramirez, H.~Zwart, and Y.~Le~Gorrec.
\newblock Stabilization of infinite dimensional port-{H}amiltonian systems by
  nonlinear dynamic boundary control.
\newblock {\em Automatica}, 85:61--69, 2017.

\bibitem{RaR15}
J.~Rawlings and M.~Risbeck.
\newblock On the equivalence between statements with epsilon-delta and
  {K}-functions.
\newblock Technical report, 2015.

\bibitem{RoF10}
H.~Royden and P.~Fitzpatrick.
\newblock {\em Real Analysis}.
\newblock Prentice Hall, 2010.

\bibitem{Rue07}
B.~R{\"{u}}ffer.
\newblock {\em Monotone dynamical systems, graphs, and stability of large-scale
  interconnected systems}.
\newblock PhD thesis, Fachbereich 3 (Mathematik \& Informatik) der
  Universit{\"{a}}t Bremen, 2007.

\bibitem{Rue10}
B.~S. R{\"u}ffer.
\newblock Monotone inequalities, dynamical systems, and paths in the positive
  orthant of {E}uclidean $n$-space.
\newblock {\em Positivity}, 14(2):257--283, 2010.

\bibitem{Rue17}
B.~S. R{\"u}ffer.
\newblock Nonlinear left and right eigenvectors for max-preserving maps.
\newblock In {\em Positive Systems}, volume 471 of {\em {L}ecture {N}otes in
  {C}ontrol and {I}nformation {S}ciences}, pages 227--237. Springer, Cham,
  2017.

\bibitem{Sac87}
P.~L. Sachdev.
\newblock {\em Nonlinear Diffusive Waves}.
\newblock Cambridge University Press, 1987.

\bibitem{Sak47}
S.~Saks.
\newblock {\em Theory of the Integral}.
\newblock Courier Corporation, 1947.

\bibitem{Sal87}
D.~Salamon.
\newblock Infinite-dimensional linear systems with unbounded control and
  observation: {A} functional analytic approach.
\newblock {\em Transactions of the American Mathematical Society},
  300(2):383--431, 1987.

\bibitem{Sch74}
H.~H. {Schaefer}.
\newblock {\em Banach Lattices and Positive Operators}.
\newblock Springer, Berlin, 1974.

\bibitem{Sch22}
J.~Schmid.
\newblock Well-posedness and stability of non-autonomous semilinear
  input-output systems.
\newblock {\em Evolution Equations and Control Theory}, 2022.

\bibitem{ScZ21}
J.~Schmid and H.~Zwart.
\newblock Stabilization of port-{H}amiltonian systems by nonlinear boundary
  control in the presence of disturbances.
\newblock {\em ESAIM: Control, Optimisation and Calculus of Variations}, 27,
  2021.

\bibitem{Sch20}
F.~L. Schwenninger.
\newblock Input-to-state stability for parabolic boundary control: {L}inear and
  semilinear systems.
\newblock In {\em Control Theory of Infinite-Dimensional Systems}, pages
  83--116. Springer, 2020.

\bibitem{ShL12}
Y.~Sharon and D.~Liberzon.
\newblock Input to state stabilizing controller for systems with coarse
  quantization.
\newblock {\em IEEE Transactions on Automatic Control}, 57(4):830--844, 2012.

\bibitem{SWT22}
S.~Singh, G.~Weiss, and M.~Tucsnak.
\newblock A class of incrementally scattering-passive nonlinear systems.
\newblock {\em Automatica}, 142:110369, 2022.

\bibitem{Siv77}
G.~Sivashinsky.
\newblock Nonlinear analysis of hydrodynamic instability in laminar
  flames--{I}. {D}erivation of basic equations.
\newblock {\em Acta Astronautica}, 4:1177--1206, 1977.

\bibitem{Son89}
E.~D. Sontag.
\newblock Smooth stabilization implies coprime factorization.
\newblock {\em IEEE Transactions on Automatic Control}, 34(4):435--443, 1989.

\bibitem{Son98}
E.~D. Sontag.
\newblock Comments on integral variants of {ISS}.
\newblock {\em Systems \& Control Letters}, 34(1-2):93--100, 1998.

\bibitem{Son08}
E.~D. Sontag.
\newblock Input to state stability: Basic concepts and results.
\newblock In {\em Nonlinear and Optimal Control Theory}, chapter~3, pages
  163--220. Springer, Heidelberg, 2008.

\bibitem{SoK03}
E.~D. Sontag and M.~Krichman.
\newblock An example of a {GAS} system which can be destabilized by an
  integrable perturbation.
\newblock {\em IEEE Transactions on Automatic Control}, 48(6):1046--1049, 2003.

\bibitem{SoW95}
E.~D. Sontag and Y.~Wang.
\newblock On characterizations of the input-to-state stability property.
\newblock {\em Systems \& Control Letters}, 24(5):351--359, 1995.

\bibitem{SoW96}
E.~D. Sontag and Y.~Wang.
\newblock New characterizations of input-to-state stability.
\newblock {\em IEEE Transactions on Automatic Control}, 41(9):1283--1294, 1996.

\bibitem{SoW97}
E.~D. Sontag and Y.~Wang.
\newblock Output-to-state stability and detectability of nonlinear systems.
\newblock {\em Systems \& Control Letters}, 29(5):279--290, 1997.

\bibitem{Sta05}
O.~Staffans.
\newblock {\em Well-posed Linear Systems}.
\newblock Cambridge University Press, 2005.

\bibitem{Ste05}
S.~Sternberg.
\newblock Theory of functions of a real variable.
\newblock {\em Harvard University, E-print disponible}, 2005.

\bibitem{Sza65}
J.~Szarski.
\newblock {\em Differential Inequalities}.
\newblock Polish Sci. Publ. PWN, Warszawa, Poland, 1965.

\bibitem{TMP18}
A.~Tanwani, S.~Marx, and C.~Prieur.
\newblock Local input-to-state stabilization of 1-{D} linear reaction-diffusion
  equation with bounded feedback.
\newblock In {\em Proc. of 23rd International Symposium on Mathematical Theory
  of Networks and Systems}, pages 576--581, 2018.

\bibitem{TPT16}
A.~Tanwani, C.~Prieur, and S.~Tarbouriech.
\newblock Input-to-state stabilization in $h^{1}$-norm for boundary controlled
  linear hyperbolic {PDE}s with application to quantized control.
\newblock In {\em Proc. of 55th IEEE Conference on Decision and Control}, pages
  3112--3117, 2016.

\bibitem{TPT18}
A.~Tanwani, C.~Prieur, and S.~Tarbouriech.
\newblock Stabilization of linear hyperbolic systems of balance laws with
  measurement errors.
\newblock In {\em Control Subject to Computational and Communication
  Constraints}, pages 357--374. Springer, 2018.

\bibitem{Tee98}
A.~R. Teel.
\newblock Connections between {R}azumikhin-type theorems and the {ISS}
  nonlinear small-gain theorem.
\newblock {\em IEEE Transactions on Automatic Control}, 43(7):960--964, 1998.

\bibitem{TeP00}
A.~R. Teel and L.~Praly.
\newblock A smooth {L}yapunov function from a class-estimate involving two
  positive semidefinite functions.
\newblock {\em ESAIM: Control, Optimisation and Calculus of Variations},
  5:313--367, 2000.

\bibitem{Tes12}
G.~Teschl.
\newblock {\em Ordinary Differential Equations and Dynamical Systems}.
\newblock American Mathematical Society, 2012.

\bibitem{TWJ12}
S.~Tiwari, Y.~Wang, and Z.~P. Jiang.
\newblock Nonlinear small-gain theorems for large-scale time-delay systems.
\newblock {\em Dynamics of Continuous, Discrete and Impulsive Systems Series A:
  Mathematical Analysis}, 19(1):27--63, 2012.

\bibitem{TuW09}
M.~Tucsnak and G.~Weiss.
\newblock {\em Observation and Control for Operator Semigroups}.
\newblock Birkh{\"{a}}user Basel, 2009.

\bibitem{TuW14}
M.~Tucsnak and G.~Weiss.
\newblock Well-posed systems -- the {LTI} case and beyond.
\newblock {\em Automatica}, 50(7):1757--1779, 2014.

\bibitem{ScJ14}
A.~van~der Schaft and D.~Jeltsema.
\newblock Port-{H}amiltonian systems theory: {A}n introductory overview.
\newblock {\em Foundations and Trends in Systems and Control}, 1(2-3):173--378,
  2014.

\bibitem{Nee96}
J.~van Neerven.
\newblock {\em The Asymptotic Behaviour of Semigroups of Linear Operators}.
\newblock Birkh{\"{a}}user Basel, 1996.

\bibitem{Vaz07}
J.~L. V{\'a}zquez.
\newblock {\em The Porous Medium Equation: Mathematical Theory}.
\newblock Oxford University Press, 2007.

\bibitem{VCL07}
L.~Vu, D.~Chatterjee, and D.~Liberzon.
\newblock Input-to-state stability of switched systems and switching adaptive
  control.
\newblock {\em Automatica}, 43(4):639--646, 2007.

\bibitem{Wei89b}
G.~Weiss.
\newblock Admissibility of unbounded control operators.
\newblock {\em SIAM Journal on Control and Optimization}, 27(3):527--545, 1989.

\bibitem{Wil72}
J.~C. Willems.
\newblock Dissipative dynamical systems part {I}: {G}eneral theory.
\newblock {\em Archive for Rational Mechanics and Analysis}, 45(5):321--351,
  1972.

\bibitem{Wil72b}
J.~C. Willems.
\newblock Dissipative dynamical systems part {II}: {L}inear systems with
  quadratic supply rates.
\newblock {\em Archive for Rational Mechanics and Analysis}, 45(5):352--393,
  1972.

\bibitem{YSH13}
S.~Yang, B.~Shi, and S.~Hao.
\newblock Input-to-state stability for discrete-time nonlinear impulsive
  systems with delays.
\newblock {\em International Journal of Robust and Nonlinear Control},
  23(4):400--418, 2013.

\bibitem{ZhX13}
C.-R. Zhao and X.-J. Xie.
\newblock Output feedback stabilization using small-gain method and
  reduced-order observer for stochastic nonlinear systems.
\newblock {\em IEEE Transactions on Automatic Control}, 58(2):523--529, 2013.

\bibitem{ZhZ19b}
J.~Zheng and G.~Zhu.
\newblock A {D}e {G}iorgi iteration-based approach for the establishment of
  {ISS} properties for {B}urgers' equation with boundary and in-domain
  disturbances.
\newblock {\em IEEE Transactions on Automatic Control}, 64(8):3476--3483, 2018.

\bibitem{ZhZ18}
J.~Zheng and G.~Zhu.
\newblock Input-to-state stability with respect to boundary disturbances for a
  class of semi-linear parabolic equations.
\newblock {\em Automatica}, 97:271--277, 2018.

\bibitem{ZhZ20b}
J.~Zheng and G.~Zhu.
\newblock A weak maximum principle-based approach for input-to-state stability
  analysis of nonlinear parabolic {PDE}s with boundary disturbances.
\newblock {\em Mathematics of Control, Signals, and Systems}, 32:157--176,
  2020.

\end{thebibliography}

\cleardoublepage
\pdfbookmark[0]{\indexname}{idx}
\printindex

%

\end{document}